\documentclass{book}

\usepackage{pstricks}
\usepackage{pst-plot}
\usepackage{times}
\usepackage{amsthm}
\usepackage{amsmath}
\usepackage{ifthen}
\usepackage{amssymb}
\usepackage{oldgerm}
\usepackage{euscript}
\usepackage{imakeidx}
\usepackage{delarray}
\usepackage{graphicx}
\usepackage{graphics}
\usepackage{subcaption}
\usepackage{enumerate}
\usepackage{calc}
\usepackage{kbordermatrix}
\usepackage{yfonts}
\usepackage[geometry]{ifsym}
\usepackage{tikz}
\usetikzlibrary{matrix, arrows, decorations.pathmorphing, calc}
\usepackage[titles]{tocloft}
\makeindex

\theoremstyle{definition}
\newtheorem{definition}{Definition} [section]
\newtheorem{example}{Example} [section]

\theoremstyle{remark}

\theoremstyle{plain}
\newtheorem{theorem}{Theorem}[section]

\newtheorem{corollary}[theorem]{Corollary}

\makeatletter
\def\@startsection#1#2#3#4#5#6{%
  \if@noskipsec \leavevmode \fi
  \par
  \@tempskipa #4\relax
  \@afterindenttrue
  \ifdim \@tempskipa <\z@
    \@tempskipa -\@tempskipa %\if@alwaysindent\else\@afterindentfalse\fi
  \fi
  \if@nobreak
    \everypar{}%
  \else
    \addpenalty\@secpenalty\addvspace\@tempskipa
  \fi
  \@ifstar
    {\@ssect{#3}{#4}{#5}{#6}}%
    {\@dblarg{\@sect{#1}{#2}{#3}{#4}{#5}{#6}}}}

\@addtoreset{figure}{section}

\@addtoreset{table}{section}

\makeatother

\begin{document}

\begin{titlepage}
	\begin{center}
		\vspace*{1cm}
		\Huge
		\textbf{Algebraic Topology for Data Scientists}
		
		\vspace{0.5cm}
		
		Michael S. Postol. Ph.D.
		
		\vspace{0.5cm}
		
		\Large
		
		\today
		
		\vspace{0.5cm}
		
		The MITRE Corporation
		
		\vspace{0.5cm}
		
		mpostol@mitre.org
		
		\vspace{0.5cm}
		
		All rights reserved.
		
		\vspace{2.0cm}
		
		The author's affiliation with The MITRE Corporation is provided for identification purposes only and is not intended to convey or imply MITRE's concurrence with, or support for, the positions, opinions, or viewpoints expressed by the author.
	\end{center}
\end{titlepage}

 \newpage
\chapter*{Acknowledgements}

I would like to thank the following people for their help in putting together this book in its current form:

First of all I would like to thank my long time collaborator, Bob Simon, for leading the development of our time series analysis tool (TSAT). This grew out of my interest in the time series discretization algorithm developed by Jessica Lin and Eamonn Keogh. Jessica Lin, gave us a lot of advice on how to use SAX for our novel cybersecurity applications. Candace Diaz has an algebraic topology background and we have had many discussions on current topics in topological data analysis. She was the person who introduced me to the Gidea-Katz work which we applied to extending SAX to the problem of classifying multivariate time series. Drew Wicke did a large portion of the software development part of TSAT and helped in the experiments described in our paper.

In addition to the people listed above, I would also like to thank Liz Munch, Rocio Gonzalez-Diaz, Andrew Blumberg, Emilie Purvine, and Tegan Emerson for helpful conversations. I would like to thank Marian Gidea, Yuri Katz, Eamonn Keogh, Jessica Lin, Ian Witten, Dmitry Cousin, Henry Adams, Liz Munch, Robert Ghrist, Vin de Silva, Greg Friedman, and Frances Sergereart for the use of illustrations from their works. Mary Ann Wymore helped me with the proper formatting and with legal expertise. I apologize for anyone I may have left out from the above list.

Finally, I would like to thank my wife Nadine for her help in pulling together the permissions. I am grateful to her for her support for this effort and all her support in general over the last 26 years.

\tableofcontents
\pagebreak
\listoftables
\pagebreak
\listoffigures
\pagebreak

\setcounter{secnumdepth}{10}

%\abstract{This paper gives the background needed to understand the emerging field of topological data analysis. After providing the background in point-set topology, abstract algebra, and algebraic topology needed for this subject, it describes persistent homology and some useful technique used to apply it to data science. It then describes some more advanced topics such as cohomology, homotopy, obrstruction theory and Sttenrod squares and explores some ideas for the application of these ideas in having a more complete understanding of data.}
\setcounter{page}{1}

\chapter{Introduction}

This textbook will give a thorough introduction to topological data analysis (TDA), the application of algebraic topology to data science. Algebraic topology was the topic of my Ph.D. dissertation, written 30 years ago at the University of  Minnesota \cite{Pos1, Pos2}. At the time, the subject was rather theoretical so I have always been intrigued by the idea of having a practical application. Today, there seems to be a growing interest in using topology for data science and machine learning. At the International Conference on Machine Learning Applications held in December 2019, there was a special session in topological data analysis that was filled beyond capacity for the entire session.

Algebraic topology is traditionally a very specialized field of math and most mathematicians have never been exposed to it, let alone data scientists, computer scientists, and analysts. I have three goals in writing this book. The first is to bring people up to speed who are missing a lot of the necessary background. I will describe the topics in point-set topology, abstract algebra, and homology theory needed for a good understanding of TDA. The second is to explain TDA and some current applications and techniques. Finally, I would like to answer some questions about more advanced topics such as cohomology, homotopy, obstruction theory, and Steenrod squares, and what they can tell us about data. It is hoped that readers will have the tools to start to think about these topics and where they might fit in.

One important question is whether TDA gives any information that more conventional techniques don't. Is it just a solution looking for a problem? What are the advantages? I strongly believe that TDA is a good tool to be used in conjunction with other machine learning methods. I am hoping to make that case in this discussion.

So what is the advantage of algebraic topology? First of all it is a classifier. The idea is to classify geometric objects. Roughly, two objects or {\it spaces} (more on that term in the next chapter) are equivalent, if it is possible to deform one into the other without tearing or gluing. You may have heard that a topologist has a hard time distinguishing between a donut and a coffee cup. In reality, though, it would be a hollow donut and a covered coffee cup with a hollow handle. (See figure 1.0.1) Another way to picture it, would be a hollow donut vs. a sphere with a hollow handle. 

\begin{figure}[ht]
\begin{center}
  \scalebox{0.4}{\includegraphics{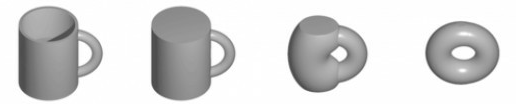}}
\caption{
\rm 
Coffee Cup Turning into a Donut \cite{MSX1}
}
\end{center}
\end{figure}

Can you tell the difference between a line and a plane? One way is to remove a point. The line is now in two pieces but the plane stays connected. But what about a plane and three dimensional space? They are actually quite different, but how do you know? If you remove a point from either, they stay in one piece. When we get to homology, I will show you a very easy way to tell them apart.

Now look at the data set plotted in Figure 1.0.2. It looks like a ring of points, but how does a computer know that? As far as the computer knows, it is just a bunch of disconnected points. As you will see, persistent homology, will easily determine its shape.

\begin{figure}[ht]
\begin{center}
  \scalebox{0.4}{\includegraphics{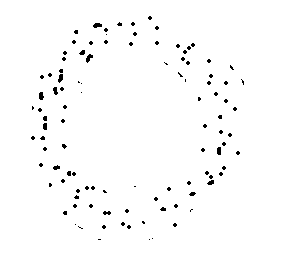}}
\caption{
\rm 
Ring of Data Points
}
\end{center}
\end{figure}

Another advantage is that topology is very bad at telling things apart. Think of all the coffee cups that have been mistakenly eaten. In algebraic topology, a circle, a square, and a triangle are pretty much the same thing. What this does for you is to assure you that the differences it detects are pretty significant and it should reduce your false positive rate.

Topological data analysis is easier than ever before to experiment with as there is now a choice of fast and easily accessible software.

Finally, you will learn a whole new collection of math puns to impress your friends. For example, what do you get when you cross an elephant with an ant? The trivial elephant bundle over the ant. If you want to understand that one, keep reading.

In chapters 2 and 3, I will give the necessary background on point-set topology and abstract algebra. If you are already fluent in these, you can skim through these chapters or skip them entirely. Chapter 4 will deal with homology theory in the more traditional way. Although many papers just jump to persistent homology  with half page descriptions, some exposure to more conventional homology theory will put you in a good position both to understand TDA and possibly to extend it. In chapter 5, I will present persistent homology and some of its visualizations such as bar codes, persistence diagrams, and persistence landscapes. I will also describe applications to data in the form of sublevel sets, graphs, and time series. This will include some of my own work on multivariate time series and the application to the classification of internet of things data. In chapter 6, I will discuss some applications such as Q-analysis, sensor coverage, Ayasdi's "mapper" algorithm, simplicial sets, and the dimesnionality reduction algorithm, UMAP. In chapter 7, I will discuss some ideas I have had in the past but never finished on the use of traditional homology. These include simplicial complexes derived from graphs and their potential application to change detection in a time series of graphs. Also, I will describe an idea for using homology for market basket analysis. How do you distinguish between a grocery store shopper who is buying for a large family and one who is single?

In the last part, I will discuss an idea I have for analyzing a cloud of data points using obstruction theory. (Don't worry if you don't know what all or even most the words in this paragraph mean. I will get you there.) Obstruction theory deals with the "extension problem."  If you have a map from a space A to a space Y, and A is a subspace of X, can you extend the map continuously to all of X while keeping its values on A the same as before. Obstruction theory combines cohomology and homotopy theory. So chapter 8 will deal with cohomology. Cohomology's strength is the idea of a "cup product" giving cohomology the structure of a ring. (See chapter 3 for rings.) Chapter 9 will deal with homotopy. Homotopy theory is much easier to describe than homology but computations are extremely difficult. The advantage is that when you can do the computations, you get a lot more information than homology theory gives you. In Chapter 10, I will describe obstruction theory and my ideas for how to apply it to data science. Chapter 11 will deal with Steenrod squares and reduced powers which give cohomology the structure of an algebra. (Again see chapter 3). Finally, Chapter 12 will describe homotopy groups of spheres. This is a problem that has not been completely solved but I feel that the early work that was done in the 1950's and 1960's holds the key to the potential application of obstruction theory to data science.

This is quite an ambitious list, but I hope to make it understandable even if you don't yet have all of the background. The idea is to get you used to the language and the concepts. TDA is a field with a lot of untapped potential so I would like to start you thinking about all of the possibilities for its future.

One final note. {\it Topology} is not the same as {\it topography}\index{topography}. Topography looks at features of a geographic area while topology deals with the classification of geometric objects. And {\it network topology} is an example of a common feature in math: word reuse.

\section{Note on the Use of This Book.}

The last four chapters of this book contain some very specialized subjects in algebraic topology that even most math Ph.D.'s have never seen. Reading the first 5 chapters will give you a strong background and allow you to read most of the current papers in the subject. Chapters 6 and 7 provide some interesting special topics, and cohomology which is covered in Chapter 8, will allow you to understand RIPSER, a very fast state of the art software package which will allow you to do most computations in topological data analysis.

The last four chapters are what make this book unique. In his 1956 American Mathematical Society Colloquium Lecture \cite{Ste1}, Norman Steenrod (my mathematical great-grandfather and someone you will soon hear more about) described his motivation for his construction of Steenrod squares (see Chapter 11). The idea is to encode more of the geometry of a shape into the algebra by introducing additional structure. This allows shapes to be classified when they can't be classified by easier methods. 

It seems to me that difficult classification problems in data science should work the same way. The additional structure provided by homotopy groups, obstructions, and Steenrod squares should allow for the solution of harder classification problems.

There are two issues that then need to be addressed. Are these methods tractable and are they useful? I believe that the answer to the first question is that homotopy groups, obstructions, and Steenrod squares are tractable for a range of problems. The question of their usefulness has yet to be answered. I feel. though, that it is important to lay out the issues and let the reader decide if there are problems that could be addressed using these methods. There is already some work in this direction and it could be an interesting research area for the future.

\chapter{Point-Set Topology Background}

What is space? According to the Hitchhiker's Guide to the Galaxy \cite{Ada}, "Space is big. You just won't believe how vastly, hugely, mind-bogglingly big it is. I mean, you may think it's a long way down the road to the chemist's, but that's just peanuts compared to space." 

Actually, there are really only 2 types of spaces: {\it topological spaces} and {\it vector spaces}. Any others are just special cases of these. The space we live in is really both, so it is a topological vector space. I will deal with topological spaces in this chapter and postpone vector spaces to chapter 3.

Point set topology is concerned with the properties of topological spaces, sets of points with a special class of subsets called {\it open sets}. Topological spaces correspond to geometric shapes. Topology describes the properties of geometric shapes we are used to and gives all sorts of examples of strange shapes where things break down. Fortunately in algebraic topology, almost everything we deal with is nicely behaved, and in TDA even more so. A big advantage in TDA is that there are only a finite number of data points, greatly simplifying the theory. By the time someone collects infinitely many samples, I will have plenty of time to revise this book. 

Rather than write a text for the usual semester length course, I will just focus on defining the terms and giving some basic properties. One very popular text is Munkres \cite{Mun2}. I will follow the book by Willard \cite{Wil}, the book that I used in college, and most of the material in this chapter is taken from there. In section 2.1, I will describe the basics of sets and functions, mainly so we can agree on notation. I will also briefly discuss, commutative diagrams. At the end, I will define one of the most familiar type of topological space, metric spaces. Section 2.2  will formally define topological spaces. In section 2.3, I will discuss continuous functions which have a much easier description than when you first met them in calculus. Section 2.4 will cover subspaces, product spaces and quotient spaces, all of which will have analogues in abstract algebra as you will see in Chapter 3. Section 2.5 will present the hierarchy of spaces defined by the separation axioms. In Section 2.6, I will do a little more set theory and discuss the difference between countably infinite and uncountable infinite sets. I will conclude with compactness in Section 2.7 and connectedness in Section 2.8, two ideas that will be very important in algebraic topology.

\section{Sets, Functions, and Metric Spaces}
\subsection{Sets}
A {\it set}\index{set} is just a collection of objects called {\it elements}. In point-set topology, the elements are called {\it points}. If $a$ is an element of set $A$ we write $a\in A$. Otherwise, we write $a\notin A$. If $A$ and $B$ are sets, we say that $A$ is {\it contained} in $B$ or $A\subset B$ if every element of $A$ is also an element of $B$. In this case, $A$ is a {\it subset} of $B$.\index{subset}  If  $A\subset B$ and $B\subset A$ then A and B have exactly the same elements and $A=B$. 

Write $A-B$ for the set of all elements in A that are not in B. If  $B\subset A$, then $A-B$ is the {\it complement}\index{complement} of $B$ in $A$. 

Figure 2.1.1 illustrates some of these concepts. Here we show the union of overlapping sets $A$ and $B$ and have a set $C$ such that $C\subset A$. A picture like this is called a {\it Venn Diagram.}

\begin{figure}[ht]
\begin{center}
  \scalebox{0.4}{\includegraphics{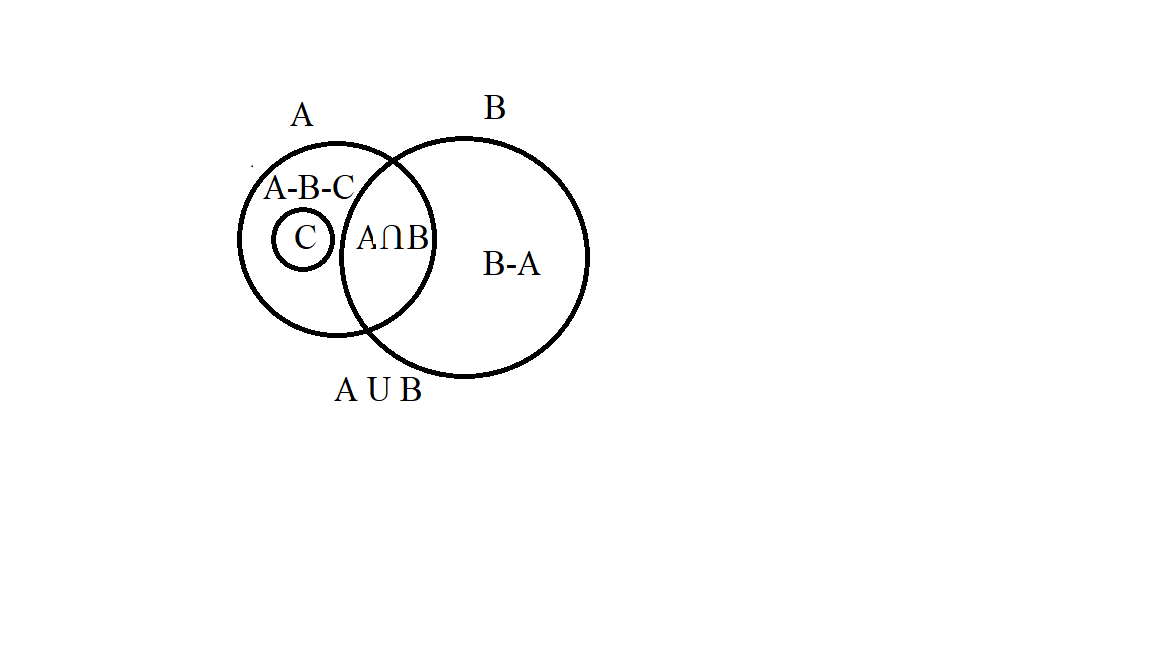}}
\caption{
\rm 
Union, Intersection and Set Difference
}
\end{center}
\end{figure}

If $A$ consists of the elements 1, 2, x, cat, dog, and horse, we write $A=\{1, 2, x, cat, dog, horse\}$. A set with no elements is the {\it empty set}\index{empty set} and is written $\emptyset$ or $\{\}$.
 
A set $A$ can be an element of another set $B$ or even itself.  $A\in B$ does not mean the same thing as  $A\subset  B$ as in the latter case the elements of $A$ are all elements of $B$ while in the former case $A$ itself is an element of $B$. This leads to a question. Is there a set of all sets? If there is, this leads to {\it Russell's Paradox}\index{Russell's Paradox}. If there was a set of all sets, it would have to include the set of all sets that don't contain themselves as an element, ie $Q=\{A|A\notin A\}$. (The vertical bar means "such that".)

The famous analogy is the town with a male barber who shaves every man who doesn't shave himself. In that case, who shaves the barber? A more descriptive example is based on an analogy given in \cite{Smu}. Suppose there is an island where every inhabitant has a club named after them and every inhabitant is in at least one club. Not everyone, though is in the club named after them. Also, every club is named after someone. Call a person {\it sociable} if they are a member of the club named after them and {\it unsociable} if they are not. Could there be a club made up of all of the unsociable people? Suppose there were. It would have to be named after someone. If it was the Michael Postol club, I could not be a member as I would be sociable. Then I would not be in the club named after me so I would be unsociable so then I would have to be in the club. 

For this reason we don't talk about the set of all sets. Instead, the collection of sets form a {\it category}. Category theory deals with collections of things like sets, topological spaces, groups, etc and special maps between them. The proper place to discuss this subject is in the next chapter but category theory and algebraic topology are very closely linked. For now, think about the different types of sets we will learn about (especially in chapter 3) and the things they have in common.

If $A$ and $B$ are sets, then the {\it union}\index{union} $A\cup B$ is the set of elements that are in $A$ or $B$ or both. The  {\it intersection}\index{intersection} $A\cap B$ is the set of elements that are in both $A$ and $B$. If $A\cap B=\emptyset$, we say that $A$ and $B$ are {\it disjoint}.

For example, if $A=\{1,3,5\}$ and $B=\{5,7,9\}$ then $A\cup B=\{1,3,5,7,9\}$ and $A\cap B=\{5\}$.

Similar definitions apply for more than two and even infinitely many sets where the union denotes elements that are in at least one of the sets and the intersection is the set of elements that are in all of the sets.

Here are some sets that will arise frequently:

\begin{align}
R &: \text{The set of real numbers.}\\
R^n&:  \text{The set of ordered }n\text{-tuples of real numbers.}\\
Z&:  \text{The set of integers }\{\cdots,-2,-1,0,1,2,\cdots\}\\
Z_n&:\text{The set }\{0,1,2,\cdots,n-1\} \text{ of integers modulo }n.\\
Q&: \text{The set of rational numbers}.
\end{align}

\subsection{Functions and Relations}

A method of building a new set from an old one is to form the {\it (Cartesian) product}. 

\begin{definition}
If $X_1$ and $X_2$ are two sets then the {\it Cartesian product} \index{Cartesian product} (or simply the {\it product}) of  $X_1$ and $X_2$ is the set $X_1 \times X_2$ is the set of all ordered pairs $(x_1,x_2)$ such that $x_1 \in X_1$ and $x_2\in X_2$. 
\end{definition}

When the sets are geometric shapes we can think of the Cartesian product as the shape resulting from putting a copy of one set at every point of the other set. For example, $R^n \times R=R^{n+1}$. As another example, let $I$ be the interval $[0,1]$ and $S^n$ \index{$S^n$}be the hollow $n-$sphere defined to be the set of all points $(X_1,\cdots, X_{n+1})$ such that $X_1^2+\cdots+X_{n+1}^2=1$. So $S^1$ is a circle and $S^2$ is the surface of a sphere.Then $S^1 \times I$ is a cylinder and $S^1 \times S^1$ is a torus (a hollow donut). A Moebius strip looks like a cylinder held up close but has a twist in it. Later we will look at Cartesian products with a "twist". These are called {\it fiber bundles}. See Figure 2.1.2 for some pictures.

\begin{figure}[ht]
	\begin{subfigure}{.5\textwidth}
		\centering
		\includegraphics[width=.4\linewidth]{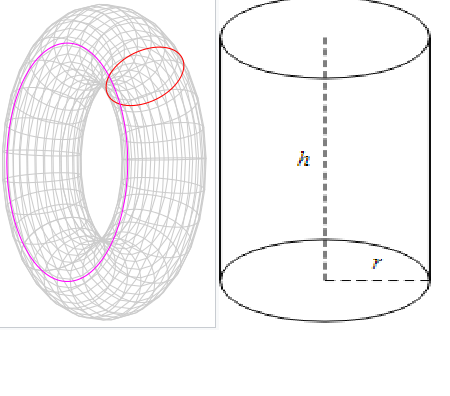}
		\caption{Torus and Cylinder}
		\label{fig:sub1}
	\end{subfigure}%
	\begin{subfigure}{.5\textwidth}
		\centering
		\includegraphics[width=.4\linewidth]{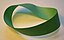}
		\caption{Moebius Strip}
		\label{fig:sub2}
	\end{subfigure}
	
	\caption{
		\rm 
		Torus \cite{WikT}, Cylinder\cite{WikC}, and Moebius Strip\cite{WikMB}
	}
	
\end{figure}

\begin{definition}
A {\it function}\index{function} $f$ from a set $X$ to a set $Y$ written $f: X\rightarrow Y$ is a subset of $X\times Y$ with the property that for each $x\in X$ there is one and only one $y\in Y$ such that $(x,y)\in f$. In this case, we write $f(x)=y.$
\end{definition}

In this above definition, $X$ is the {\it domain}\index{function: domain} of $f$ and $Y$ is the  {\it range}\index{function: range} of $f$. The {\it image} \index{function: image} $f(X)$ of $f$ is the set defined to be the set of $y\in Y$ such that $f(x)=y$ for some $x \in X$. If $X=Y$ there is a special function $f:X\rightarrow X$ called the {\it identity function}\index{identity function}. In this case, $f(x)=x$ for all $x\in X$. 

If  $f:X\rightarrow Y$ and  $g:Y\rightarrow Z$ then there is a function $gf: X\rightarrow Z$  called the {\it compositio}n of $g$ with $f$ and defined by $gf(x)=g(f(x))$. We can also write it this way.
\begin{tikzpicture}
  \matrix (m) [matrix of math nodes,row sep=3em,column sep=4em,minimum width=2em]
  {
   X & Y & Z\\};
  \path[-stealth]
    (m-1-1) edge node [above] {$f$} (m-1-2)
    (m-1-2)  edge node [above] {$g$} (m-1-3);
\end{tikzpicture}

 For example, if $f(x)=x^2+1$ and $g(x)=\sin(x)$, then $gf(x)=\sin(x^2+1)$.

If $f:X\rightarrow Y$ , we would like a way to undo $f$. To do this, we need two things to happen. First of all, for every $y\in Y$, we need at least one $x\in X$ for it to map to. The second condition is that for every $y\in Y$, we need at most one $x\in X$ for it to map to. The first condition is called onto and the second is called one-to-one.

\begin{definition}
A function $f$ is {\it one-to-one}\index{one-to-one} or {\it injective}\index{injective} or an {\it injection}\index{injection} if for $x_1, x_2\in X$, if $x_1\neq x_2$ implies $f(x_1)\neq f( x_2)$. A function $f$ is {\it onto}\index{onto} or {\it surjective}\index{surjective} or a {\it surjection}\index{surjection} if for any $y\in Y$, there is an $x\in X$ such that $f(x)=y$. In this case, $f(X)=Y$, i.e. the range of $f$ is all of $Y$.
\end{definition}

\begin{example}
The function $f(x)=x^2$ is not one-to-one since $f(-1)=1=f(1)$. It is also not onto if $Y=R$, since there is no real value of $x$ with $x^2=-1$. But $g(x)=x^3$ is both one-to-one and onto. 
\end{example}

\begin{definition}
A function $f:X\rightarrow Y$  has a {\it (two-sided) inverse} $g:Y\rightarrow X$ if $gf=1_X$ and $fg=1_Y$ where $1_X$ is the identity function on $X$. In this case, we write $g=f^{-1}(x)$. So $f^{-1}$ undoes $f$. By our remarks above, a function has an inverse if and only if it is one-to-one and onto.
\end{definition}

It is possible for a function to be one-to-one, but not onto. Then it has an inverse function whose domain is the original function's range. An example of this is the function $f(x)=e^x.$ This is one-to-one, but its range is $R^+$, the postitive real numbers as opposed to all of $R$. This function has an inverse $f^{-1}(x)=ln(x)$ , which is only defined on the range $R^+$ of $f$. 

\begin{definition}
Let $X\subset Y$. Then $f: X\rightarrow Y$ is an {\it inclusion}\index{inclusion} if $f(x)=x$ for all $x\in X$. Then $f$ is obviously one-to-one, but $f$ is only onto if $X=Y$.
\end{definition}

Closely related to functions are {\it relations}. A {\it relation}\index{relation} is any subset of $A\times A$. So a function $F: A\rightarrow A$ is a special type of relation in which the subset can not contain both $(a_1,b_1)$ and $(a_1,b_2)$ for $a_1, b_1, b_2\in A$ as this is saying that $f(a_1)=b_1$ and $f(a_1)=b_2$. There are two types of relations that will be important in what follows. These are {\it equivalence relations} and {\it partial orders}.

\begin{definition}
We will define a particular relation between elements in $A$ writing $a\sim b$ if $(a,b)$ is in the relation. Then $\sim$ is an  {\it equvalence relation}\index{equvialence relation} if the following three properties hold for all $a, b, c\in A$:\begin{enumerate}
\item Reflexive: $a\sim a$
\item Symmetric: If $a\sim b$ then $b\sim a$
\item Transitive: If $a\sim b$ and $b\sim c$ then $a\sim c$. 
\end{enumerate}
In this case, if $a\sim b$, we say that $a$ is equivalent to $b$.
\end{definition}

Equality between real numbers is an example of an equivalence relation. Congruence of triangles is another. We will meet many more.

\begin{definition}
We will define a particular relation between elements in $A$ writing $a\leq b$ if $(a,b)$ is in the relation. Then $\leq$ is an  {\it partial order}\index{partial order} if the following three properties hold for all $a, b, c\in A$:\begin{enumerate}
\item Reflexive: $a\leq a$
\item Antisymmetric: If $a\leq b$ and $b\leq a$, then $a=b$
\item Transitive: If $a\leq b$ and $b\leq c$ then $a\leq c$. 
\end{enumerate}
\end{definition}

Real numbers obviously form a partial order where $\leq$ has the usual meaning. Note that we didn't require that $a\leq b$ or $b\leq a$ for any $a, b\in A$. If this holds, then we have a {\it linear order}. The real numbers have a linear order. 

An example of a partial order that is not a linear order is the following: Consider subsets of a set $X$. We say that for $A, B\subset X$, we have $A\leq B$ if $A\subset B$. You can easily check that all of the properties hold, but it is possible that neither $A\subset B$ nor $B\subset A$ is true. In  this case, we say that $A$ and $B$ are {\it incomparable}.

To conclude this section, I will introduce the idea of a {\it commutative diagram}\index{commutative diagram}. If you are an algebraic topologist, you learn to love them. 

Consider the following diagram:

\begin{tikzpicture}
  \matrix (m) [matrix of math nodes,row sep=3em,column sep=4em,minimum width=2em]
  {
   A & B & C \\
     D & E & F \\};
  \path[-stealth]
    (m-1-1) edge node [left] {$f_1$} (m-2-1)
 (m-1-2) edge node [left] {$f_2$} (m-2-2)
 (m-1-3) edge node [left] {$f_3$} (m-2-3)
  (m-1-1) edge node [above] {$g_1$} (m-1-2)
  (m-1-2) edge node [above] {$g_2$} (m-1-3)
  (m-2-1) edge node [above] {$g_3$} (m-2-2)
  (m-2-2) edge node [above] {$g_4$} (m-2-3)
           (m-1-2) edge node [above] {$h$} (m-2-3) ;
\end{tikzpicture}

The vertices represent sets and the arrows are functions between them. As before, following two arrows in succession represents function compostiton. The idea is that no matter what path you take, as long as you are following arrows, you get to the same place. So let $a\in A$. Then we can get to $F$ in several ways, and $g_4g_3f_1(a)$, $g_4f_2g_1(a)$, $f_3g_2g_1(a)$, and $hg_1(a)$ are all the same element of $F$.

\subsection{Metric Spaces}

In the next section I will finally introduce topological spaces, but here we will start with the most famous special case. A metric space is basically a space with the notion of a distance between two points. 

\begin{definition}
A {\it metric space}\index{metric space} is a set $M$ together with a function $\rho: M\times M\rightarrow R$ such that we have the following for $x,y,z \in M$:\begin{enumerate}
\item $\rho(x,y)\geq 0$ and $\rho(x,y)=0$ if and only if $x=y$.
\item $\rho(x,y)=\rho(y,x)$
\item Triangle Inequality: $\rho(x,y)+\rho(y,z)\geq \rho(x,z)$. 
\end{enumerate}
\end{definition}

So the first property says that the difference between two points is a positive number unless they are the same point in which it is 0. The second property is that the distance from $x$ to $y$ is the same as the distance from $y$ to $x$. The triangle inequality says that the sum of the lengths of two sides of a triangle is greater than or equal to the length of the third side. They are equal if the three points lie on the same line. An equivalent statement is that the shortest distance between two pints is a straight line.

The function $\rho$ is called a distance or equivalently, a metric.

\begin{example}
$M$ is the real numbers where $\rho(x,y)=|x-y|$ for $x,y\in R$.
\end{example}

\begin{example}
$M$ is $R^n$. For $x=(x_1, x_2, \cdots, x_n), y=(y_1, y_2, \cdots, y_n)\in R^n$, we have $\rho(x,y)=\sqrt{\sum_{i=1}^n{(x_i-y_i)^2}}$. This is the usual {\it Euclidean} metric.
\end{example}

\begin{example}
$M$ is $R^2$. For $x=(x_1, x_2), y=(y_1, y_2)$, we have $\rho(x,y)=|x_1-y_1|+|x_2-y_2|$. This is called the {\it taxi cab} or  {\it Manhattan} distance.
\end{example}

\begin{example}
Let $L^p$  for $p$ a positive integer be the set of real valued functions $f$ of a single real variable such that $\int_R |f|^p<\infty$.Then for $f, g\in L^p$ we have the {\it $L^p$ distance} $\rho(f,g)=(\int_R |f-g|^p)^\frac{1}{p}$.
\end{example}

The last two examples will serve as opposite extremes when we discuss topological spaces.

\begin{example}
The {\it discrete metric space}: Let $M$ be any set. For $x, y\in M$ we let $\rho(x,y)=0$ if $x=y$, and let $\rho(x,y)=1$ if $x\neq y$.
\end{example}

You should check that this is an actual metric space. If you lived in this world, everyone would be one unit tall, your commute to work would be one unit, the distance to the sun (assuming you had one) would be one unit, etc. Sort of convenient but as we will see, everything would fall apart and there would be no path from anywhere to anywhere else.

The opposite extreme is the {\it trivial pseudometric space} in which $\rho(x,y)=0$ for any $x$ and $y$. This is a pseudometric rather than a metric as property one is violated but properties 2 and 3 still hold.

Finally, given a metric, we have the idea of a ball. This is the set of points whose distance from a center point is  less than or equal to a particular value called the {\it radius}. 

\begin{definition}
Let $M$ be a metric space and $x\in M$. The {\it open  ball of radius $\delta$ centered at x} written $B_\delta(x)$ is the set of $y\in M$ such that $\rho(x,y)<\delta$. If instead, we have $\rho(x,y)\leq \delta$, the set is called a {\it closed ball} and we will write it as $\bar{B}_\delta(x)$.
\end{definition}

Note that if $M=R$, then $B_\delta(x)$ is the open interval $(x-\delta, x+\delta)$. If $M=R^2$, then $B_\delta(x)$ is the circle centered at $x$ of radius $\delta$ (without the boundary). 

The open balls in a metric space will be an example of the more general concept of a base for a topology. We will explore this further in the next section.

\section{Topological Spaces}

We are now ready to define a topological space. This is a set with a special class of subsets called {\it open sets.}

\begin{definition}
A {\it topology}\index{topology} on a set $X$ is a collection $\tau$ of subsets of $X$ called the {\it open sets}\index{open set} satisfying:\begin{enumerate}
\item Any union (finite or infinite) of subsets of $\tau$ belong to $\tau$.
\item Any {\it finite} intersection of subsets of $\tau$ belong to $\tau$.
\item $\emptyset$ and $X$ belong to $\tau$.
\end{enumerate}
The set $X$ together with the topology $\tau$ is called a {\it topological space}\index{topological space}. In algebraic topology, the word {\it space} will always mean a topological space unless it is explicitly stated otherwise.
\end{definition}

\begin{definition}
Given two topologies $\tau_1$ and $\tau_2$ on the same set $X$, we say that $\tau_1$ is {\it weaker (smaller, coarser)} than $\tau_2$ or equivalently that $\tau_2$ is {\it stronger (larger, finer)} than $\tau_1$ if $\tau_1\subset\tau_2$.
\end{definition}

\begin{example}
Given a nonempty set $X$, the weakest topology is the {\it trivial topology}\index{topology!trivial} in which $\tau_1={X, \emptyset}$. At the other extreme is the {\it discrete topology}\index{topology!discrete} in which $\tau_2$ consists of every subset of $X$. Obviously  $\tau_1\subset\tau_2$. In fact, the discrete topology is the strongest possible topology on $X$. You should check that both of these satisfy Definition 2.2.1.
\end{example}

\begin{example}
Let $M$ be a metric space. Then the {\it metric topology} or {\it usual topology} consists of sets $G$ which are defined to be open if and only if for any $x\in G$, there exists a real number $\delta>0$ such that $B_\delta(x)\subset G$. 
\end{example}

In other words, G is open if for any point in $G$, we can draw an open ball around it that stays in $G$. For example, the open interval $(0, 1)$ is open. Let $x=.1$. Then if we let $\delta=.05$, then $x-\delta=.05$, and $x+\delta=.15$, so the open ball is $(.05, .15)$ and that is contained in $(0,1)$. It should be clear that we can do something similar no matter how close we get to 0 or 1. So $(0,1)$ is open.

We can also see why the union of open sets can be infinite but we needed to specify a finite union. Consider the collection of open intervals $(-\frac{1}{n},\frac{1}{n})$, for $n=1, 2, \cdots$. These are all open, but $$\bigcap_{n=1}^\infty(-\frac{1}{n},\frac{1}{n})=\{0\}$$
The set $\{0\}$ is not open since any open interval containing $0$ will contain points other than $0$.

\begin{definition}
If $X$ is a topological space, then a subset $F$ is {\it closed}\index{closed set} if the complement of $F$ in $X$ is open.
\end{definition}

For the rest of this section, you should convince yourself of the truth of the results by drawing pictures. The first one is analogous to the definition of open set and can be shown by looking at complements.

\begin{theorem}
Let $X$ be a topological space and let \textfrak{F} be the collection of closed subsets of $X$. Then the following hold:\begin{enumerate}
\item Any intersection (finite or infinite) of subsets of \textfrak{F} belong to \textfrak{F}.
\item Any {\it finite} union of subsets of\textfrak{F} belong to \textfrak{F}.
\item $\emptyset$ and $X$ belong to\textfrak{ F}.
\end{enumerate}
\end{theorem}

On $R$, closed intervals are closed sets as their complements are open. For intervals, we can make an open interval closed by adding in its boundary. This motivates the following:

\begin{definition}
If $X$ is a topological space, and $E\subset X$, then the {\it closure}\index{closure} of $E$ written $\bar{E}$ is the intersection of all closed subsets of $X$ which contain $E$. By the previous theorem, $\bar{E}$ must also be closed.
\end{definition}

\begin{theorem}
If $A\subset B$ then $\bar{A}\subset\bar{B}$.
\end{theorem}

\begin{theorem}
Closure of subsets in a topological space have the following properties:\begin{enumerate}
\item $E\subset\bar{E}$
\item $\bar{\bar{E}}=\bar{E}$
\item $\overline{A\cup B}=\bar{A}\cup\bar{B}$
\item $\bar{\emptyset}=\emptyset$
\item $E$ is closed in $X$ if and only if $\bar{E}=E$.
\end{enumerate}
\end{theorem}

Note that the closure of the intersection is not necessarily the intersection of the closure. Let $A$ be the rationals in $R$ and $B$ be the irrationals in $R$. It turns out that $R=\bar{A}=\bar{B},$ so $\bar{A}\cap\bar{B}=R$, but $A\cap B=\emptyset$, so $\overline{A\cap B}=\bar{\emptyset}=\emptyset$.

You can think of closure as putting in the rest of the boundary of a set. The opposite of that is the interior operation.

\begin{definition}
If $X$ is a topological space, and $E\subset X$, then the {\it interior}\index{interior} of $E$ written $E^\circ$ is the union of all open subsets of $X$ which are contained in $E$. By the definition of an open set, $E^\circ$ must also be open.
\end{definition}

\begin{theorem}
If $A\subset B$ then $A^\circ\subset B^\circ$.
\end{theorem}

\begin{theorem}
Interiors of subsets in a topological space have the following properties:\begin{enumerate}
\item $E\supset E^\circ$
\item $E^{\circ\circ}= E^\circ$
\item $(A\cap B)^\circ=A^\circ\cap B^\circ$
\item $X^\circ=X$
\item $G$ is open in $X$ if and only if $G=G^\circ$.
\end{enumerate}
\end{theorem}

Again let $A$ be the rationals and $B$ be the irrationals. Then $A\cup B=R$, and $R^\circ=R$. But $A^\circ=B^\circ=\emptyset$. So  $(A\cup B)^\circ\neq A^\circ\cup B^\circ$

I alluded to a boundary that is the difference between an open and a closed set. I will now make this precise:

\begin{definition}
If $X$ is a topological space, and $E\subset X$, then the {\it boundary}\index{boundary} or {\it frontier}\index{frontier} of $E$ written $Bd(E)$ is the set $Bd(E)=\bar{E}\cap\overline{(X-E)}.$
\end{definition}

\begin{theorem}
For any subset $E$ of a topological space $X$:\begin{enumerate}
\item $\bar{E}=E\cup Bd(E)$
\item $E^\circ=E-Bd(E)$
\item $X=E^\circ\cup Bd(E)\cup (X-E)^\circ$
\end{enumerate}
\end{theorem}

To help you picture these results, consider Figure 2.2.1. The set $E$ consists of the small region including the solid but not the broken line. The closure $\bar{E}$ includes the broken  line as well. The interior $E^\circ$ is everything on the inside but does not include even the solid line. The frontier or boundary $Bd(E)$ is the broken and solid lines together. You should convince yourself of the last few theorems at least in this case.

\begin{figure}[ht]
\begin{center}
  \scalebox{0.4}{\includegraphics{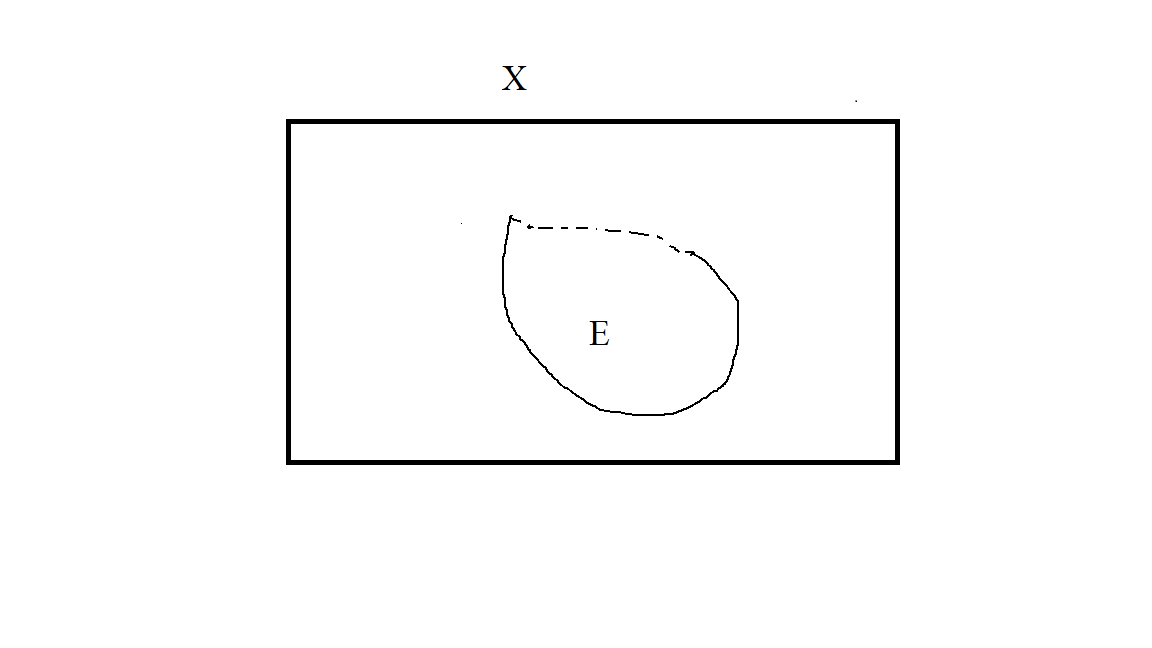}}
\caption{
\rm 
Closure, Interior, and Frontier(Boundary)
}
\end{center}
\end{figure}

In the remainder of this section, I will briefly discuss neighborhoods and bases. These generalize the open balls in a metric space.

 \begin{definition}
If $X$ is a topological space, and $x\in X$, then a {\it neighborhood}\index{neighborhood} of $x$ is a set $U$ which contains an open set $V$ containing $x$. So $U$ is a neighborhood if and only if $x\in U^\circ$. The collection $\textfrak{U}_x$ of all neighborhoods of $x$ is called a  {\it neighborhood system}\index{neighborhood system} at $x$.
\end{definition}

A neighborhood system can be used to define a topology on $X$.

\begin{theorem}
The neighborhood system $\textfrak{U}_x$ at $x$ in a topological space $X$ has the following properties:\begin{enumerate}
\item If $U\in\textfrak{U}_x$, then $x\in U$.
\item  If $U, V\in\textfrak{U}_x,$ then $U\cap V\in\textfrak{U}_x.$
\item  If $U\in\textfrak{U}_x$, then there is a $V\in\textfrak{U}_x$ such that $U\in\textfrak{U}_y$ for each $y\in V.$
\item If $U\in\textfrak{U}_x$ and $U\subset V$ then $V\in\textfrak{U}_x$.
\item A subset $G$ of $X$ is open if and only if $G$ contains a neighborhood of each of its points.
\end{enumerate}
\end{theorem}

Since any set containing a neighborhood is also a neighborhood, we will limit ourselves to a smaller collection called a {\it neighborhood base}.

 \begin{definition}
If $X$ is a topological space, and $x\in X$, then a {\it neighborhood base}\index{neighborhood base} at $x$ is a subcollection $\textfrak{B}_x$ taken from the neighborhood system  $\textfrak{U}_x$ having the property that each $U\in\textfrak{U}_x$ contains some $V\in\textfrak{B}_x$. An element  of a neighborhood base is a {\it basic neighborhood}\index{basic neighbohood}
\end{definition}

\begin{example}
For any topological space, the open neighborhoods of $x$ form a neighborhood base since if $U$ is a neighborhood of $x$, $U^\circ$ must also be one.
\end{example}

\begin{example}
For $x$ in a metric space, the open balls $B_\delta(x)$ form a neighborhood base at $x$.
\end{example}

\begin{example}
If $X$ has the discrete topology, the set $\{x\}$ is a neighborhood base at $x$ since every subset of $X$ containing $x$ is open and thus a neighborhood of $x$. If $X$ has the trivial topology, the only neighborhood base at any point is $\{X\}$ as $X$ is the only nonempty open set and thus the only neighborhood of any point $X$. 
\end{example}

\begin{theorem}
The neighborhood base $\textfrak{B}_x$ at $x$ in a topological space $X$ has the following properties:\begin{enumerate}
\item If $V\in\textfrak{B}_x$, then $x\in V$.
\item  If $U, V\in\textfrak{B}_x,$ then there is a $W\in\textfrak{B}_x$ such that $W\subset U\cap V$.
\item If $V\in\textfrak{B}_x$ then there is some $V_0\in\textfrak{B}_x$ such that if $y\in V_0$, then there is some $W\in\textfrak{B}_y$ with $W\subset V$.
\item A subset $G$ of $X$ is open if and only if $G$ contains a basic neighborhood of each of its points.
\end{enumerate}
\end{theorem}

By property 4, a subset $G$ of a metric space $M$ is open if and only if for any $x\in G$, there is a $\delta>0$ such that the open ball $B_\delta(x)\subset G$. This the usual definition of an open set that you may have seen.

The next theorem characterizes some of the concepts we have seen in terms of neighborhood bases.

\begin{theorem}
Let $X$ be a topological space and suppose we have fixed a neighbohood base at each point $x\in X$. Then we have:\begin{enumerate}
\item $G\subset X$ is open if and only if $G$ contains a basic neighborhood of each of its points.
\item  $F\subset X$ is closed if and only if each point $x\notin F$ has a basic neighborhood disjoint from $F$.
\item $\bar{E}$ is the set of points $x$ in $X$ such that each basic neighborhood of $x$ meets $E$. 
\item $E^\circ$ is the set of points $x$ in $X$ such that some basic neighborhood of $x$ is contained in $E$.
\item $Bd(E)$ is the set of points $x$ in $X$ such that every basic neighborhood of $x$ meets both $E$ and $X-E$. 
\end{enumerate}
\end{theorem}

A more global version of the above is a base for a topology.

\begin{definition}
If $X$ is a topological space with topology $\tau$, a {\it base}\index{topology!base} for $\tau$ is a collection $\textfrak{B}\subset \tau$, such that $\tau$ is produced by taking all possible unions of elements in  $\textfrak{B}.$
\end{definition}

\begin{example}
For a metric space $X$,  the open balls $B_\delta(x)$ for $\delta>0$ and $x\in X$ form a base for the usual topology..
\end{example}

The connection with neighborhood bases is that for each $x\in X$, the sets in \textfrak{B} which contain $x$ form a neighbohood base at $x$.

Finally, you may see a topology defined in terms of a subbase.

\begin{definition}
If $X$ is a topological space with topology $\tau$, a {\it subbase}\index{topology!subbase} for $\tau$ is a collection $\textfrak{C}\subset \tau$, such that the set of all possible intersections of elements in  $\textfrak{C}$ is a base for $\tau$.
\end{definition}

\section{Continuous Functions}

As we will see in Chapter 3 as well, for each type set we describe, we will look at a privileged function that preserves the key defining property of the set. In abstract algebra, this function that preserves the algebraic structure is called a {\it homomorphism}. (Lots more on those later.) For a topological space, the special functions are called {\it continuous.}

If you remember continuous functions from calculus but hated epsilons and deltas, I am going to make you very happy. Continuous functions simply preserve the idea of open sets, although in the reverse direction.

\begin{definition}
Let $X$ and $Y$ be topological spaces and $f: X\rightarrow Y$. Then $f$ is {\it continuous}\index{continuosus function} at $x_0\in X$ if for each neighborhood V of $f(x_0)$ in Y, there is a neighborhood U of $x_0$ in X such that $f(U)\subset V$. We say that $f$ is continuous on $x$ if it is continuous at every point of $X$. 
\end{definition}

By what we have said about the metric space topology, you should convince yourself that the definition of continuity that you learned in calculus is equivalent to our definition here.

The next theorem gives a much more useful definition of continuity. It is the one we tend to use in practice.

\begin{theorem}
Let $X$ and $Y$ be topological spaces and $f: X\rightarrow Y$. Then the following are equivalent:\begin{enumerate}
\item $f$ is continuous.
\item  for each open set $G$ in $Y$, $f^{-1}(G)$ is open in $X$.
\item  for each closed set $F$ in $Y$, $f^{-1}(F)$ is closed in $X$
\item For each $E\subset X$, $f(\bar{E})\subset\overline{f(E)}$. 
\end{enumerate}
\end{theorem}

The characterization of a continuous function as a function having the property that the inverse image of an open set is open immediately implies the following: 

\begin{theorem}
If $X$, $Y$, and $Z$ are topological spaces, and $f:X\rightarrow Y$ and $g: Y\rightarrow Z$ are continuous, then the composition $gf: X\rightarrow Z$ is also continuous.
\end{theorem}

From now on, the terms function or map between topological spaces will always be assumed to be continuous unless otherwise specified. 

Note that for a continuous function, we have not claimed that the image of an open set is always open. For example, let $f$ be a real valued function of one variable defined by $f(x)=0$. It is continuous as the inverse image of any open set containing $0$ is all of $R$, which is also open. The inverse image of an open set not containing $0$ is $\emptyset$ which is also open. But as an example $f(0,1)=\{0\}$ so $f$ takes the open set $(0,1)$ to the closed set$\{0\}$. The maps which preserve open sets in both directions are called {\it homeomorphisms}. 

\begin{definition}
Let $X$ and $Y$ be topological spaces and $f: X\rightarrow Y$. Then $f$ is a {\it homeomorphism}\index{homeomorphism} if $f$ is continuous, one-to-one, and onto and if $f^{-1}$ is continuous. In this case we say that $X$ and $Y$ are {\it homeomorphic}. 
\end{definition}

\begin{theorem}
Let $X$ and $Y$ be topological spaces and $f: X\rightarrow Y$. Then the following are equivalent:\begin{enumerate}
\item $f$ is a homeomorphism.
\item  A set $G$ is open in $X$ if and only if $f(G)$ is open in $Y$
\item  A set $F$ is closed in $X$ if and only if $f(F)$ is closed in $Y$
\item For each $E\subset X$, $f(\bar{E})=\overline{f(E)}$. 
\end{enumerate}
\end{theorem}

If $X$ and $Y$ are homeomorphic, then from the point of view of topology, they are basically the same object. Note that they may not seem exactly the same. For example, a circle, a square and a triangle are all homeomorphic. But they can be continuously deformed into each other without tearing or gluing anything. (Those terms will be more precise when we discuss connectedness.) The equivalent term for the algebraic objects we will discuss in chapter 3 is {\it isomorphic}. The relation of being homeomorphic is an equivalence relation. The main goal of algebraic topology is to classify topological spaces by whether they are homeomorphic to each other. In general, this is a very hard problem.

\section{Subspaces, Product Spaces, and Quotient Spaces}

Here are some new topological spaces you can create from old ones. These all have analogues in abstract algebra.

\subsection{Subspaces}
We get this one almost for free. When is a subset of a topological space another topological space? Always. You just need to define an open set in the right way.

\begin{definition}
Let $X$ be a topological space with topology $\tau,$ and let $A\subset X$. The collection $\tau^\prime=\{G\cap A|G\in\tau\}$ is a topology for $A$ called the {\it relative topology}\index{relative topology} for $A$. In this case, we call the subset $A$ a {\it subspace}\index{subspace} of $X$.
\end{definition}

So basically, we get our open sets in $A$ by taking the open sets of $X$ and intersecting them with $A$.

\begin{theorem}
If $A$ is a subspace of a topological space $X$ then:\begin{enumerate}
\item $H\subset A$ is open in $A$ if and only if $H=G\cap A$ where $G$ is open in $X$. 
\item $F\subset A$ is closed in $A$ if and only if $F=K\cap A$ where $K$ is closed in $X$. 
\item for $E\subset A$, the closure of $E$ in $A$ is the intersection of $A$ with the closure of $E$ in $X$.
\item If $x\in A$, then $V$ is a neighborhood of $x$ in $A$ if and only if $V=U\cap A$ where $U$ is a neighborhood of $x$ in $X$.
\item If $x\in A$, and $\textfrak{B}_x$ is a neighborhood base at $x$ in $X$, then $\{B\cap A|B\in \textfrak{B}_x\}$ is a neighborhood base at $x$ in $A$.
\item If $\textfrak{B}$ is a base for $X$, then $\{B\cap A|B\in \textfrak{B}\}$ is a base for $A$.
\end{enumerate}
\end{theorem}

Note that we haven't talked about interiors or frontiers. The next example shows why.

\begin{example}
Let $X$ be $R^2$ with the usual topology and let $A=E$ be the $x$-axis. Then the interior in $A$ of $E$ is $A$ while the interior in $X$ of $E$ is $\emptyset$, and $A\neq\emptyset\cap A$. It is always true, though that the interior of $E$ in $A$ is contained in the intersection of $A$ and the interior of $E$ in $X$.
In the same case, the frontier of $E$ in $A$ is $\emptyset$ while the frontier of $E$ in $X$ is $A$ and $A\neq\emptyset\cap A$. It is always true, though that the frontier of $E$ in $A$ is contained in the intersection of $A$ and the frontier of $E$ in $X$.
\end{example}

Finally, we mention that when restricted to subspaces, continuous functions stay continuous. 

\begin{definition}
If $f: X\rightarrow Y$ and $A\subset X$, then the {\it restriction} of $f$ to $A$ written $f|A$ is the map of $A$ into $Y$ defined by $(f|A)(a)=f(a)$ for every $a\in A$. 
\end{definition}

\begin{theorem}
If $A\subset X$ and $f:X\rightarrow Y$ is continuous, then $(f|A): A\rightarrow Y$ is continuous. 
\end{theorem}

\begin{theorem}
If $X=A\cup B$, where $A$ and $B$ are both open or both closed in $X$ and if  $f:X\rightarrow Y$ is a function such that both $(f|A)$ and $(f|B)$ are continuous, then $f$ is continuous.
\end{theorem}

\subsection{Product Spaces}

We discussed the cartesian product of two sets in Section 2.1.2. We repeat the definition here for reference.

\begin{definition}
If $X_1$ and $X_2$ are two sets then the {\it Cartesian product} \index{Cartesian product} (or simply the {\it product}) of  $X_1$ and $X_2$ is the set $X_1 \times X_2$ is the set of all ordered pairs $(x_1,x_2)$ such that $x_1 \in X_1$ and $x_2\in X_2$. 
\end{definition}

We will now loosen up the definition a bit. First of all, we will allow as many of these sets as we want, even infinitely many. For the product of $n$ sets, the product is $$\prod_{i=1}^n X_i=\{(x_1,\cdots,x_n)\},$$ where $x_i\in X_i$ for $1\leq i\leq n.$ To define the product of infinitely many sets, we use the idea of an {\it index set}\index{index set}. An index set is a set $A$ of objects where we have one term for each object in $A$. So here is the most general definition we will use:

\begin{definition}
Let $X_\alpha$ be a set for each $\alpha\in A$. The {\it Cartesian product} \index{Cartesian product} (or simply the {\it product}) of  the $X_\alpha$ is the set $$\prod_{\alpha\in A} X_\alpha=\{x: A\rightarrow\bigcup_{\alpha\in A} X_\alpha|x(\alpha)\in X_\alpha {\rm\hspace{2.0 pt}  for\hspace{2.0 pt} each\hspace{2.0 pt}} \alpha\in A\}.$$
\end{definition}

So basically, we define a point in the product as a function which picks out one element from each of the sets. In order to do this we need the {\it Axiom of Choice}. I won't go into a long discussion of this but informally, suppose you had infinitely many shoes. You could easily pick out one from each pair by just taking all of the left ones. But suppose you had infinitely many pairs of socks. It is not obvious that you could pick out one from each pair as this would keep you busy all morning and every other morning until the end of time. So the Axiom of Choice says that you can pick out one element from each of an infinite collection of sets. If you believe it or if you don't, all of the other axioms of set theory still hold.

The one case where we may be interested in an infinite Cartesian product is when talking about spaces whose points are actually functions. Consider a real valued function of one variable defined on all of $R$. Any such function is a point in the space $R^R$. This is the product of the reals taken infinitely many times, one for each real number. For example, take the function $f(x)=x^2$. If we think of a coordinate for every real number, then the coordinate corresponding to 3 will be $3^2=9$. 

It will be helpful in what follows to picture products as finite to understand the results in this section but I will do this in more generality so you can understand the proofs in the paper on template functions that I will discuss in Chapter 6.

Since we are working with topological spaces, we want to know what open sets look like on a product of spaces. To get there we need a couple more terms. The value $x(\alpha)\in X_\alpha$ is usually denoted $x_\alpha$ and called the $\alpha$th coordinate of $X$. The space $X_\alpha$ is the $\alpha$th factor space. 

The map $\pi_\beta:\prod X_\alpha\rightarrow X_\beta$ defined by $\pi_\beta(x)=x_\beta$ is called the $\beta$th {\it projection map}.

As an example, consider $R^3$ which is the product of the reals with itself three times. The index set $A$ is simply $\{1,2,3\}$. Let $x=(3,-4.5, 100)$ be a point in $R^3$. The three coordinates are $x_1=3$, $x_2=-4.5$, and $x_3=100$. The projection map $\pi_2$ takes a point in $R^3$ to the second coordiante so $\pi_2(3,-4.5, 100)=-4.5.$

Now we are ready to define a topology on a prodcut space.

\begin{definition}
The {\it Tychonoff topology} or {\it product topology}\index{Tychonoff topology}\index{product topology} on $\prod X_\alpha$ is obtained by taking as a base for the open sets the sets of the form $\prod U_\alpha$ where\begin{enumerate}
\item $U_\alpha$ is open in $X_\alpha$ for each $\alpha$ in $A$.
\item For all but finitely many coordinates, $U_\alpha=X_\alpha.$
\end{enumerate}
\end{definition}

In a finite product, we can dispense with the second condition, so open sets are unions of sets of the form $\prod_{i=1}^n U_\alpha$ where  $U_\alpha$ is open in $X_\alpha.$

\begin{definition}
If $X$ and $Y$ are topological spaces, and $f:X\rightarrow Y$, then $f$ is an {\it open map} (\it{closed map})\index{open map}\index{closed map} if for each open (closed) set $A$ in $X$, $f(A)$ is open (closed) in $Y$.
\end{definition}

An open map is almost the opposite of a continuous map. For an open map, the image of an open set is open while for a continuous map, the inverse image of an open set is open. If $F$ is one-to-one, onto, and open then $f^{-1}$ is continuous. So a one-to-one, onto, and continuous map is a homeomorphism. The same facts hold for closed maps.

\begin{theorem}
The $\beta$th projection map $\pi_\beta:\prod X_\alpha\rightarrow X_\beta$ is continuous and open but need not be closed.
\end{theorem}

As an example, consider the graph of $y=\frac{1}{x}$ as a subset of $R^2$, Then $\pi_1$ is the projection onto the x axis. If $G=\{(x,\frac{1}{x})\}$ is this graph, then $G$ is closed in $R^2$, but $\pi_1(G)=(-\infty,0)\cup(0,\infty)$, which is an open set. So in this case $\pi_1:R^2\rightarrow R$ is not closed.

Finally, we have the following:

\begin{theorem}
A map $f:Y\rightarrow\prod X_\alpha$ is continuous if and only if $\pi_\alpha\circ f$ is continuous for each $\alpha\in A$.
\end{theorem}

In abstract algebra, the analogue of a Cartesian product is a {\it direct product}.

\subsection{Quotient Spaces}
Here is where things will start to get sticky. By that, I mean that topologists love to glue things together. Twists and turns that would shred ordinary paper to pieces are totally legal in topology. In this section I will show you some examples and then define a quotient space (the result of gluing) in general. We will want to know how to define open sets and learn what happens to continuous functions. Also note, that there are analogous quotient objects in algebra. Especially quotient groups will play a vital role in algebraic topology.

\begin{figure}[ht]
\begin{center}
  \scalebox{0.4}{\includegraphics{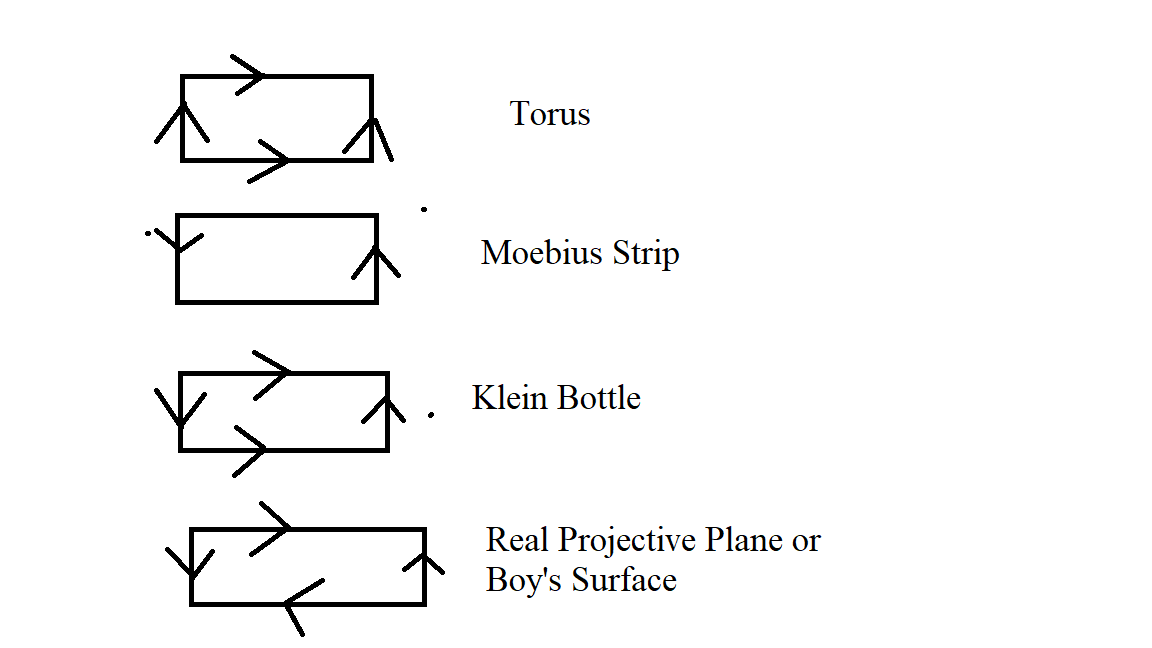}}
\caption{
\rm 
Folded Rectangles
}
\end{center}
\end{figure}

Here is an experiment you can try at home. Take a rectangular piece of paper and fold it in the four different ways shown in Figure 2.4.1. In the first example, we glue along both edges and get our old friend the torus. (See Figure 2.1.2.) In the second example, glue the short edges with a twist and leave the long ones alone. This gives a Moebius strip (also shown in Figure 2.1.2.). In the third case. glue the short ends with a twist and the long ones without one. This gives a {\it Klein Bottle}\index{Klein bottle} shown in Figure 2.4.2. If you tried to do this, you would probably have a pretty hard time. A Klein bottle only fits in 4 dimensions, but in 3 dimensions, we can do an {\it immersion}. This means that it may intersect itself. So you have to imagine that the narrow tube is not hitting a wall as it comes up on the right. If you walked through a Klein bottle and started on the inside, you would end up on the outside and vice versa. 

\begin{figure}[ht]
\begin{center}
  \scalebox{0.4}{\includegraphics{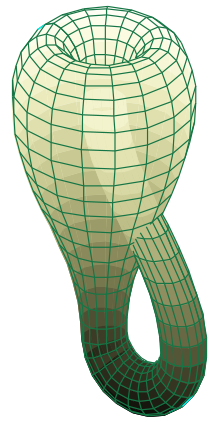}}
\caption{
\rm 
Klein Bottle \cite{WikIm}
}
\end{center}
\end{figure}

This brings up two excuses for not doing your math homework:\begin{enumerate}
\item I put it in a Klein bottle and it fell out.
\item I put it in a Klein bottle and a four dimensional dog ate it.
\end{enumerate}

In the last case we have the {\it real projective plane}\index{real projective plane}. Here we twist both the long and short edges. Again, we have a four dimensional object and we picture a three dimensional immersion called {\it Boy's surface}\index{Boy's surface} in Figure 2.4.3. Another way to make the real projective plane is as follows: Take a (hollow) sphere. Two points on a sphere are {\it antipodal} if the lie on a line through the center. So for example, the North Pole and the South Pole are antipodal. Now glue together every pair of antipodal points. Before you get to the Equator, you are fine. For example, just glue the Southern hemisphere to the Northern one. Once you start trying to stick the points on the Equator together, you will have the same trouble (as you would expect.) Projective spaces provide important examples in algebraic topology

\begin{figure}[ht]
\begin{center}
  \scalebox{0.4}{\includegraphics{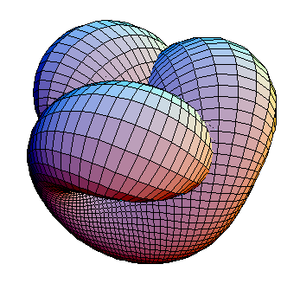}}
\caption{
\rm 
Boy's Surface \cite{WikIm}
}
\end{center}
\end{figure}

One final example will be important later. Let $X$ be a topological space and $I=[0,1]$. The {it cone}\index{cone} $CX$ of $X$ is the result of taking $X\times I$ and identifying the points $(x,1)$ to a single point. If $X$ is a circle, then $CX$ is an actual cone. Another way of thinking of a cone is to take every point in $X$ and connect it with a straight line to a fixed point outside of $X$. Suppose instead we have two fixed outside points and we connect every point of $X$ to each of them with a straight line. (It is easier to picture if you think of a point above and a point below.) This is  $X\times I$ where the points $(x, 0)$ are all identified and the points $(x,1)$ are all identified. This is called the {\it suspension}\index{suspension} and is written $SX$ or $\Sigma X$. Assuming everything is made of silly putty you can think of a cone as plugging up holes and a suspension as raising the dimension of a sphere by one. So you can flatten a cone to be a filled in circle or deform a suspension of a circle to be a sphere. Figure 2.4.4 helps illustrate this where the top or bottom half is a cone and the entire object is a suspension of a circle. (Two ice cream cones stuck together at their open ends).

\begin{figure}[ht]
\begin{center}
  \scalebox{0.4}{\includegraphics{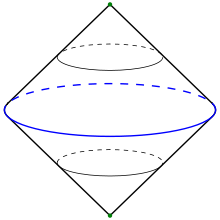}}
\caption{
\rm 
Suspension of a Circle \cite{WikSu}
}
\end{center}
\end{figure}

In practice we do the gluing by taking the original space to the new one with a function that is onto but not one-to-one as points we want to identify will be sent to the same place as each other. We define the topology on the new space in a way which makes this map continuous. The new space will be called the {\it quotient space}. 

\begin{definition}
If $X$ is topological space, $Y$ is a set, and $g:X\rightarrow Y$ is an onto mapping, then $G$ is defined to be open in $Y$ if and only if $g^{-1}(G)$ is open in $X$. The topology on $Y$ induced by $g$ and defined in this manner is called the {\it quotient topology}\index{quotient topology} and $Y$ is called a {\it quotient space}\index{quotient space} of $X$. The map $g$ is called a {\it quotient map}\index{quotient map}.
\end{definition}

Now suppose $Y$ was already a topological space. Are open sets in $Y$ the same as those defined in the quotient topology? If $g$ is continuous, then an open set in the original topology has to be open in the quotient topology by definition. To go the other way, we need $g$ to be open or closed. 

\section{Separation Axioms}

The separation algorithms form a hierarchy among topological spaces. They are classified by how well open sets can separate points in these spaces. The classification is denoted by the names $T_0, T_1, T_2, T_3,$ and $T_4$. In this section, I will define these types and show the relationships between them. I will also give some examples of weird spaces that fulfill some but not all of the conditions. This will help make sure you understand the concepts. Also, imagining weird spaces is part of the fun of point-set topology or point-set pathology as it is nicknamed. If you want to see more of these, I refer you to \cite{StSe}. The material in this section is taken from there as well as from \cite{Wil}. Note that in algebraic topology, we almost always deal with spaces that are $T_2$. In data science applications, the spaces we deal with will always be $T_4$.

\begin{definition}
A topological space $X$ is $T_0$\index{$T_0$}  if whenever $x$ and $y$ are distinct points in $X$, there is an open set containing one and not the other.
\end{definition}

\begin{example}
A space with at least two points having the trivial topology is not even $T_0$. The only nonempty open set is the entire space which has to contain both of these points.
\end{example}

\begin{definition}
A topological space $X$ is $T_1$\index{$T_1$}  if whenever $x$ and $y$ are distinct points in $X$, there is an open set containing {\it each} one and not the other.
\end{definition}

\begin{example}
This example will illustrate the difference between $T_0$ and $T_1$. Let $X=\{a,b\}$ be a space consisting of two points and suppose that the opens sets are $\tau=\{X,\emptyset,\{a\}\}$. (Check that this satisfies the definition of a topology.) Then $X$ is $T_0$ since there is an open set containing $a$ and not $b$, namely $\{a\}$. But  $\{b\}$ is not open (note that it is closed though being the complement of $\{a\}$.) The smallest open set containing $b$ is all of $X$ which also contains $a$. So $X$ is $T_0$ but not $T_1$
\end{example}

It should be obvious that any $T_1$ space is $T_0$. Here is an interesting fact about $T_1$ spaces. I will provide the proof which is short and instructive.

\begin{theorem}
The following are equivalent for a topological space $X$:\begin{enumerate}
\item $X$ is $T_1$
\item Each one point set in $X$ is closed.
\item Each subset of $X$ is the intersection of the open sets containing it.
\end{enumerate}
\end{theorem}

{\bf Proof:}

$(1)\Rightarrow (2)$: If $X$ is $T_1$ and $x\in X$ then each $y\neq x$ is contained in an open set disjoint from $x$, so $X-\{x\}$ is a union of open sets and thus open. So $x$ is closed.

$(2)\Rightarrow (3)$: If $A\subset X$ then $A$ is the intersection of all sets of the form $X-\{x\}$ for $x\notin A$. Each of these is open since one point sets are closed.

$(3)\Rightarrow (1)$: If (3) holds then $\{x\}$ is the intersection of all open sets containing it, so for any $y\neq x$, there is an open set containing $x$ and not $y$. $\blacksquare$

The next level is $T_2$. Another name for a $T_2$ space is a {\it Hausdorff space}. I will feel free to use the two terms interchangeably. In algebraic topology, we work almost exclusively with Hausdorff spaces. They were named after Felix Hausdorff, whose 1914 book {\em Set Theory} was the first textbook on point-set topology. The original version is in German, but see \cite{Haus} for an English translation.

\begin{definition}
A topological space $X$ is $T_2$ or {\it Hausdorff} \index{$T_2$}\index{Hausdorff}  if whenever $x$ and $y$ are distinct points in $X$, there are disjoint open subsets $U$ and $V$ of $X$ such that $x\in U$ and $y\in V$. 
\end{definition}

\begin{example}
Any metric space $M$ is Hausdorff. Let $x$ and $y$ be in $X$ and let $\rho(x,y)=\delta$. (Recall that  $\rho(x,y)$ is the distance between $x$ and $y$. Then letting $U(x,\epsilon)$ be an open ball centered at $x$ of radius $\epsilon$, we have that $U(x,\delta/2)$ and $U(y,\delta/2)$ are disjoint open sets containing $x$ and $y$ respectively.
\end{example}

It should be obvious that a Hausdorff space is $T_1$, but there are $T_1$ spaces that are not Hausdorff as the next example shows.

\begin{example}
Let $X$ be an infinite set with the {\it cofinite topology}\index{topology!cofinite}. In this topology, open sets are those whose complements are finite, i.e. the complements of open sets have finitely many points. We also add in $X$  and $\emptyset$. (Check that this satisfies the definition of a topology.) This means that the closed sets are the finite subsets of $X$, with the exception of $X$ itself.  Then one point sets are closed so by Theorem 2.5.1, $X$ is $T_1$. But $X$ is not $T_2$. Now it is a fact that the complement of an intersection of two sets is a union of their complements. (Try to prove this or convince yourself using a Venn diagram.) So if $U$ and $V$ are open sets and $U\cap V=\emptyset$, then the union of their complements is all of $X$. But this is impossible since $X-U$ and $X-V$ are both finite, but $X$ is infinite. Thus, $X$ is not $T_2$.
\end{example}

\begin{theorem}
\begin{enumerate}
\item Every subspace of a $T_2$ space is $T_2$.
\item A nonempty product space is $T_2$ if and only if each factor is $T_2$.
\item Quotients of $T_2$ spaces need not be $T_2$.
\end{enumerate}
\end{theorem}

\begin{example}
Things are even worse then part 3 of the previous theorem. A continuous open image of a Hausdorff space need not be Hausdorff. Let $A$ consist of the two horizontal lines in $R^2$, $y=0$ and $y=1$. Let $f:A\rightarrow B$ be a quotient map which identifies $(x,0)$ and $(x,1)$ for $x\neq 0$. So every point on the lower line except for $(0,0)$  is glued to the point directly above it. Then the space $B$ is one line except that the points $(0,0)$ and $(0,1)$ are still separated. But any open neighborhood containing one of these points has to intersect any open neighborhood containing the other. So $B$ is not Hausdorff.
\end{example}

For our next trick, we will try to replace points by arbitrary closed sets and separate them.

\begin{definition}
A topological space $X$ is {\it regular}\index{regular}  if whenever $A$ is closed in $X$, and $x\notin A$, there are disjoint open subsets $U$ and $V$ of $X$ such that $x\in U$ and $A\subset V$. A regulat $T_1$ space is called a $T_3$ \index{$T_3$}space.
\end{definition}

Note that we have not entirely preserved our hierarchy. A regular space is not necessarily Hausdorff. For example, let $X$ have the trivial topology. Let $A=\emptyset$ and $x\in X$. Then since $A$ is empty, $x\notin A$. Let $U=X$ and $V=A$. Then $U$ and $V$ are both open, $x\in U$, $A\subset V$, and $U\cap V=\emptyset$. So $X$ is regular. (Note that we have not insisted that $A$ is a proper subset of $V$ and in general we won't as we describe separation axioms.) $X$ is obviously not Hausdorff as the only open set which contains $x$ or $y$ in $X$ is all of $X$. So a regular space need not be Hausdorff.

To fix this, we defined a $T_3$ \index{$T_3$}space to be a regular $T_1$ space. Then since we saw that in a $T_1$ space, one point sets are closed, we automatically get that a $T_3$ set is Hausdorff ($T_2$). So it is important to remember that $T_2$ and {\it Hausdorff} mean the same thing, but $T_3$ and {\it regular} do not.

\begin{example}
Not every $T_2$ space is $T_3$. Let $X$ be the real line with neighborhoods of nonzero points as in the usual topology, but neighborhoods of $0$ are in the form $U-A$ where $U$ is a neighborhood of 0 in the usual topology and $$A=\{1/n|n=1, 2, \cdots\}.$$ Then $X$ is obviously Hausdorff (draw small enough balls around $x$ and $y$). But $X$ is not $T_3$ since $A$ is closed in $X$, but can't be separated from $0$ by disjoint open sets. 
\end{example}

The next two theorems state some properties of regular spaces.

\begin{theorem}
The following are equivalent for a topological space $X$:
\begin{enumerate}
\item $X$ is regular.
\item If $U$ is open in $X$ and $x\in U$, then there is an open set $V$ containing $x$ such that $\overline{V}\subset U$. 
\item Each $x\in X$ has a neighborhood base consisting of closed sets. 
\end{enumerate}
\end{theorem}

\begin{theorem}   
\begin{enumerate}
\item Every subspace of a regular ($T_3$) space is regular ($T_3)$.
\item A nonempty product space is regular  ($T_3$) if and only if each factor space is regular ($T_3)$.
\item Quotients of $T_3$ spaces need not be regular.
\end{enumerate}
\end{theorem}

At the next level, we can start to talk about separation with the use of real valued functions. This will eventually lead to our first result on extending continuous functions. I will get much more into this when I talk about obstruction theory. 

\begin{definition}
A topological space $X$ is {\it completely regular}\index{completely regular}  if whenever $A$ is closed in $X$, and $x\notin A$, there is a continuous function $f: X\rightarrow I$ such that $f(x)=0$ and $f(A)=1$. (Recall that $I=[0,1]$ is the closed unit interval.) Equivalently, we can find a continuous function  $f: X\rightarrow R$ such that $f(x)=a$ and $f(A)=b$, wher $a$ and $b$ are real numbers with $a\neq b$. Any such function $f$ is said to {\it separate} $A$ and $x$. A completely regular $T_1$ space is called a {\it Tychonoff space}\index{Tychonoff space}.
\end{definition}

Note that completely regular spaces are regular. If $A$ is closed, $x\notin A$, and  $f: X\rightarrow I$ such that $f(x)=0$ and $f(A)=1$, then $f^{-1}([0,\frac{1}{2}))$ and $f^{-1}((\frac{1}{2},1])$ are disjoint open sets containing $x$ and $A$ respectively. 

Since the symbol $T_4$ is used for normal spaces, which we will cover next, some books refer to Tychonoff spaces as $T_{3\frac{1}{2}}$\index{$T_{3\frac{1}{2}}$} spaces. 

\begin{example}
Every metric space is Tychonoff. Let $X$ be a metric space with distance function $\rho$, $A$ be closed in $X$, and $x\notin A$. For $y\in X$, let $f(y)=\rho(y,A)=\min_{a\in A}\rho(y,a)$. Then $f$ is a continuous function from $A$ to $R$, $f(A)=0$, and $f(x)>0.$ Thus, $X$ is Tychonoff.
\end{example}

\begin{theorem}
\begin{enumerate}
\item Every subspace of a completely regular (Tychonoff) space is completely regular (Tychonoff).
\item A nonempty product space is completely regular  (Tychonoff) if and only if each factor space is completely regular (Tychonoff).
\item Quotients of Tychonoff spaces need not be regular or even $T_2$.
\end{enumerate}
\end{theorem}

Here is a result that is useful in homotopy theory,

\begin{theorem}
A topological space is a Tychonoff space if and only if  it is homeomorphic to a subspace of some cube. (Here a cube\index{cube} is a product of copies of I.)
\end{theorem}

I will stop at the fourth level of our heirarchy, normal spaces. These are the logical next step but much less well behaved. Even subspaces of normal spaces are not normal. I will state some results but avoid the proofs and most of the counterexamples which will take us too far afield. I encourage you to look at \cite{Wil}, \cite{StSe}, and \cite{Mun2}  if you are curious.

\begin{definition}
A topological space $X$ is {\it normal}\index{normal space}  if whenever $A$ and $B$ are disjoint closed sets in $X$, there are disjoint open subsets $U$ and $V$ of $X$ such that $A\subset U$ and $B\subset V$. A normal $T_1$ space is called a $T_4$  \index{$T_4$}space. 
\end{definition}

Here are two relatively easy examples.

\begin{example}
A normal space need not be regular. Let $X$ be the real line and let the open sets be the sets of the form $(a, \infty)$ for $a\in R$. Then $X$ is normal as no two nonempty closed sets are disjoint. But X is not regular since the point $\{1\}$ can not be separated from the closed set $(-\infty, 0]$ by disjoint open sets. 
\end{example}

\begin{example}
If $X$ is a metric space, then $X$ is normal. Let $\rho$ be the distance function and let $A$ and $B$ be two disjoint closed sets in $X$. Since we already know that a metric space is Hausdorff, for each $x\in A$, we can pick $\delta_x>0$ such that the open ball $U(x,\delta_x)$ centered at $x$ of radius $\delta_x$ does not meet $B$. Also, for each $y\in B$, we can pick $\epsilon_y>0$ such that $U(x,\epsilon_y)$  does not meet $A$. Let $$U=\bigcup_{x\in A} U(x, \delta_x/3), V=\bigcup_{y\in B} U(y, \epsilon_y/3).$$ Then $U$ and $V$ are open sets in $X$ containing $A$ and $B$ respectively. To see that $U$ and $V$ are disjoint, suppose $z\in U\cap V$. Then for $x\in A$, and $y\in B$, $\rho(x,z)<\delta_x/3$ and $\rho(y,z)<\epsilon_y/3$, Without loss of generality we can assume that $\delta_x=\max\{\delta_x,\epsilon_y\}$. Then by the triangle inequality, $\rho(x,y)<\delta_x/3+\epsilon_y/3<\delta_x$. But then $y\in U(x,\delta_x)$ which is impossible. Thus $U$ and $V$ must be disjoint so $X$ is normal.
\end{example}

\begin{theorem}
\begin{enumerate}
\item Closed subspaces of normal spaces are normal.
\item Product spaces of normal spaces need not be normal.
\item The closed continuous image of a normal ($T_4$) space is normal ($T_4$).
\end{enumerate}
\end{theorem}

Our main interest in normal spaces is Tietze's Extension Theorem, our first example of the solution to an extension problem. The proof involves Urysohn's Lemma which is interesting in its own right.

\begin{theorem} {\bf Urysohn's Lemma}
A space $X$ is normal if and only if whenever $A$ and $B$ are disjoint closed sets in $X$, there is a continuous function $f: X\rightarrow [0,1]$ with $f(A)=0$ and $f(B)=1$. An immediate consequence is that every $T_4$ space is Tychonoff.
\end{theorem}

Let $X$ and $Y$ be topological spaces, and $A\subset X$. Let $f:A\rightarrow Y$ be continuous. The extension problem asks if we can extend $f$ to a continuous function $g:X\rightarrow Y$ whose restriction to $A$, $g|A: A\rightarrow Y$ is equal to $f$. One of the main goals of homotopy theory is to answer this question. The next theorem provides an answer in the special case where $A$ is closed in $X$ and $Y=R$. 

\begin{theorem} {\bf Tietze's Extension Theorem}
A space $X$ is normal if and only if whenever $A$ is a closed subset of $X$, and $f:A\rightarrow R$ is continuous, there is a continuous function $g: X\rightarrow R$ such that $g|A=f$.
\end{theorem}

\section{Two Ways to be Infinite}
In this section, we will discuss different sizes of infinity. The argument came from the work of Georg Cantor and was originally pretty controversial. It is pretty much accepted now and resolves a big paradox in probability theory. This will be a digression back to set theory, but it will be needed in the next section as well as many other areas of mathematics.

\begin{definition}
Let $A$ be a set. Then the {\it cardinality}\index{cardinality} of $A$, written $|A|$ is the number of elements in $A$.
\end{definition}

So if $X=\{{\rm cat, dog, horse}\}$, then $|X|=3$. So two sets with the same cardinality have the same number of elements. We can show that two sets have the same cardinality by exhibiting a {\it one-to-one correspondence.} For every element of set $A$ we pair it with exactly one element of set $B$. We will do the same thing whether these sets are finite or infinite. For infinite sets you just need to keep this in mind even though things will go against your intuition if you have never seen these arguments before.

\begin{definition}
Let $A$ be a set. Then $A$ is {\it countable}\index{countable} if there is a one-to-one correspondence between the elements of $A$ and the positive integers. If $A$ is countable, we write the cardinality of $A$ as $|A|=\aleph_0$.\index{$\aleph_0$}
\end{definition}

So a set is countable if you can count it: 1, 2, 3, 4, $\cdots$ 

\begin{example}
The positive integers are obviously countable. So are the integers. To see this we write them as 0, -1, 1, -2, 2, -2, 3, $\cdots$. Then we have the correspondence $1\leftrightarrow 0$,  $2\leftrightarrow -1$,$3\leftrightarrow 1$,$4\leftrightarrow -2$,$5\leftrightarrow 2$, $\cdots$ You might think that there are more integers than positive integers, but we can match them up. 
\end{example}

An illustration of this is in Hilbert's paradox of the Grand Hotel. There have been several versions of this, but see \cite{Gam} for one of them. This will be my own version. Suppose there is a hotel with infinitely many rooms. There is a single hallway and the rooms are numbered, 1, 2, 3, $\cdots$. Suppose you come to the hotel and find out that all of the rooms are full. There is an easy solution. Just have everyone move down one room and you can get room 1. The next day, turns out to be the first day of IICI, the Infinite International Convention on Inifinte International Conventions on Infinite International Conventions on .... The conference is very popular and there are countably infinite many attendees. There is still no problem, have everyone move to the room that is double their current number. Now there are countably many empty rooms. 

It turns out that there are also countably many rational numbers. To see this, consider the following array. I have represented every rational number at least once with a lot of duplicates. Now just count them following the arrows: 

\begin{tikzpicture}
  \matrix (m) [matrix of math nodes,row sep=3em,column sep=4em,minimum width=2em]
  {
   0 & -1 & 1 & -2 & 2 &-3 & 3 & \cdots\\
     0/2 & -1/2 & 1/2 & -2/2 & 2/2 &-3/2 & 3/2 &\cdots\\
 0/3 & -1/3 & 1/3 & -2/3 & 2/3 &-3/3 & 3/3 &  \cdots\\
0/4 & \cdots & \cdots & \cdots & \cdots &\cdots & \cdots&  \cdots\\};
  \path[-stealth]
    (m-1-1) edge (m-2-1)
 (m-2-1) edge (m-1-2)
 (m-1-2) edge (m-1-3)
 (m-1-3) edge (m-2-2)
 (m-2-2) edge (m-3-1)
(m-3-1) edge (m-4-1)
 (m-4-1) edge (m-3-2)
 (m-3-2) edge (m-2-3)
 (m-2-3) edge (m-1-4)
(m-1-4) edge (m-1-5)
(m-1-5) edge (m-2-4)
(m-2-4) edge (m-3-3)
;
\end{tikzpicture}

By continuing the zig-zag pattern and letting each row $i$ have $i$ in the denominator, we can be sure that we count all of the rationals at least once. 

If we replace row $i$ in the diagram with the set $X_i=\{x_{i1},x_{i2}, x_{i3},\cdots\}$ the same argument shows that $\cup_{i=1}^{\infty}X_i$ is countable. So we have shown the following:

\begin{theorem}
A countable union of countable sets is countable.
\end{theorem}

Now I will claim that the real numbers are {\it uncountable}, in other words, that there are somehow more of them. So what would that actually mean for infinite sets? It has to do with matching. 

Suppose you bought your silverware as a set and you started with the same number of forks and spoons. But you suspect that there is a small black hole in your kitchen and that your spoons are being pulled into a parallel universe. You have a big dinner party planned for that night. (This is post-COVID.) When you set the table, you find that there are places with forks but no spoons. Since you ran out of spoons but had forks left over, you can conclude that there are more forks than spoons. So our plan is to match real numbers up with positive integers and still have some left over. 

Construct a sequence $x_1, x_2, x_3, \cdots$ where each $x_i\in[0,1]$ is written as a binary expansion. Then let $x$ be the number where the $i$-th bit is complemented from $x_i$. Then $x$ differs from every $x_i$ on the list and can't appear on it. So any countable sequence of numbers in $[0,1]$ has to leave something out, like the fork with no spoon. So the numbers in $[0,1]$ and hence the reals are not countable.

So how many real numbers are there? To find out, we use the concept of a {\it power set}\index{power set}.

\begin{definition}
Let $A$ be a set. Then the {\it power set} of $A$ denoted $P(A)$ is the set of all subsets of $A$.
\end{definition}

\begin{example}
Let $A=\{a,b,c\}.$ Listing all of the subsets of $A$ we get:\begin{enumerate}
\item $\emptyset$
\item $\{a\}$
\item $\{b\}$
\item $\{c\}$
\item $\{a,b\}$
\item $\{a,c\}$
\item $\{b,c\}$
\item $\{a,b,c\}=A$
\end{enumerate}
So $|P(A)|=8.$
\end{example}

In general, if $|A|=n$, then $P(A)=2^n$. To see this, if $A=\{a_1,a_2,\cdots,a_n\}$, then for each $i$ where $1\leq i\leq n$, we can include $x_i$ in our subset or exclude it. So this makes the number of subsets $2^n$. 

Now looking at real numbers in $[0,1]$ we write the binary expansion as $.x_1x_2x_3\cdots$. We are in the situation above if 0 means exclude a member of a set and 1 means to include it. Since the binary digits map to the positive integers, they form a set of cardinality $\aleph_0$, so the cardinality of their power set is $2^{\aleph_0}$. We sometimes write $2^{\aleph_0}=\textfrak{C}$, where $\textfrak{C}$ stands for {\it continuum}. The {\it continuum hypothesis}\index{$\textfrak{C}$}\index{continuum hypothesis} states that there is no set with cardinality $c$ such that $\aleph_0<c<2^{\aleph_0}$. It has never been proven or disproven and whether it is true or false is independent of the other axioms of set theory. 

A related concept which appears in many topology examples is the set of {\it ordinals} $\Omega.$\index{ordinals} Note that this means something different than the usual definition of ordinal numbers (first, second, third, etc.) $\Omega$ is an uncountable linearly ordered set with a largest element $\omega_1$. (A set is $A$ is linearly ordered if it is partially ordered and has the property that for $a, b\in A$, either $a\leq b$ or $b\leq a$.) You can think of $\Omega$ as the set starting with positive integers 1, 2, 3, 4, $\cdots$. There is a smallest ordinal larger than all of these which we call $\omega_0$, the first infinite ordinal. Then we have $\omega_0+1$, $\omega_0+2$, $\omega_0+3$, $\cdots$ until we reach $2\omega_0$.Eventually we get to $\omega_0^2$, $\omega_0^3$, $\cdots$ until we finally get to the first uncountable ordinal $\omega_1$ and stop. We define $\Omega_0=\Omega-\omega_1$ and call it the set of {\it countable ordinals}.

Now we can build one of my favorite weird spaces. First attach {0} to $\Omega$ and insert a copy of $(0,1)$ between any successive pair of ordinals. This makes a very long number line called the {\it long line}\index{long line}. I will follow \cite{StSe} and let the long line $L$ only include the countable ordinals. If we extend all the way to $\omega_1$ we have the {\it extended} long line $L^*$. 

So next time you are frustrated by standing in a long line, be grateful it isn't as long as {\bf the} long line.

\section{Compactness}

If you ever took a class in Real Analysis/Advanced Calculus, you probably learned something about compactness. I found it a little confusing the first time, so hopefully I can make it clearer. In $R^n$, compactness is equivalent to being closed and bounded (pairs of points are closer than some fixed number.) Since we live in $R^3$, this is the case in our world. This is convenient as if a compact car was not closed and bounded the edges would be kind of fuzzy and it would be very hard to fit into a parking space. And if you had a compact disc (anyone remember what those are?) it would be hard to find something that could play them. But unless you are in $R^n$ with the usual metric space topology, this is not the case. Later we will see a space that is closed and bounded but not compact. 

Since we now understand what a countable set is, I will start with a little more generality. 

\begin{definition}
A subset $D$ of a set $X$  is {\it dense}\index{dense} in $X$ if $\overline{D}=X$. An equivalent definition is that every open set in $X$ contains an element of $D$.
\end{definition}

\begin{definition}
A topological space $X$ is {\it separable}\index{separable} if it has a countable dense subset.
\end{definition}

\begin{example}
It is obvious that every open interval in $R$ contains a rational. Since these form a base for the topology on $R$, we know that every open set in $R$ contains a rational. So the rationals are dense in $R$. Since we know that the rationals are countable, $R$ is separable.
\end{example}

\begin{theorem}
\begin{enumerate}
\item A continuous image of a separable space is separable. 
\item Subspaces of a separable space need not be separable. But an open subspace of a separable set is separable.
\item A product of separable spaces is separable provided each factor is separable, Hausdorff, and contains 2 points and that there are no more than $2^{\aleph_0}$ of them.
\end{enumerate}
\end{theorem}

\begin{definition}
An {\it open cover}\index{open cover} of a topological space $X$ is a collection of subsets of $X$ such that $X$ is contained in their union. We say that this collection of open sets {\it covers} $X$. A {\it subcover} is a subcollection of a cover which still covers $X$.
\end{definition}

\begin{definition}
A topological space $X$ is {\it Lindel\"{o}f}\index{Lindel\"{o}f space} if every open cover of $X$ has a countable subcover. A topological space $X$ is {\it compact}\index{compact space} if every open cover of $X$ has a {\bf finite} subcover.
\end{definition}

Be careful here. I didn't say for example that a space is compact if there is some open cover with a finite subcover. That is always true. Just cover $X$ with itself. But if I give you {\bf any} open cover, you need to pick out finitely many which still cover $X$.

We will start with Lindel\"{o}f spaces. 

\begin{example}
The real line $R$ with the usual topology is Lindel\"{o}f. Let $\{U_\alpha\}$ be an open cover of $R$. For each rational, choose an element of $\{U_\alpha\}$ containing it. Then there are countably many of these (at most) and they must still cover $R$ since the rationals are dense in $R$. So every open cover has a countable subcover.
\end{example}

\begin{theorem}
\begin{enumerate}
\item A continuous image of a Lindel\"{o}f space is Lindel\"{o}f. 
\item Closed subspaces of  Lindel\"{o}f spaces are Lindel\"{o}f. Arbitrary subspaces of  Lindel\"{o}f spaces need not be Lindel\"{o}f.
\item Products of even two Lindel\"{o}f spaces need not be Lindel\"{o}f. 
\end{enumerate}
\end{theorem}

I will prove 1. and the first part of 2.

{\bf Proof:}\newline
1) Let $f: X\rightarrow Y$ be continuous and onto and $X$ is  Lindel\"{o}f. Let $\{U_\alpha\}$ be an open cover of $Y$. Then $\{f^{-1}(U_\alpha)\}$ is an open cover of $X$. Since $X$ is Lindel\"{o}f, we can choose a countable subcollection $\{f^{-1}(U_1), f^{-1}(U_2),\cdots\}$ which covers $X$. Then $\{U_1, U_2,\cdots\}$ is a countable subcolection which covers $Y$, and $Y$ is  Lindel\"{o}f.$\blacksquare$

2) Let $F$ be closed in a Lindel\"{o}f space $X$. If $\{U_\alpha\}$ is an open cover of $F$, we have for each $\alpha$ an open subset $V_\alpha$ of $X$ such that  $U_\alpha=F\cap V_\alpha$. Then $X-F$ together with $\{V_\alpha\}$ form an open cover of $X$. Since $X$ is Lindel\"{o}f,there is a countable subcover $\{X-F, V_1, V_2,\cdots\}$. Then the corresponding $\{U_1, U_2,\cdots\}$ form a countable subcover of $F$.$\blacksquare$

\begin{example}
As an example of a subset of a  Lindel\"{o}f space that is not  Lindel\"{o}f, we let $\Omega$ be the set of ordinals, and $\Omega_0=\Omega-\omega_1$, where $\omega_1$ is the first uncountable ordinal. We provide $\Omega$ with the order topology in which a basic neighborhood of $\alpha\in\Omega_0$ is of the form $(\alpha_1,\alpha_2)$ where $\alpha_1<\alpha<\alpha_2$. A basic neighborhood of $\omega_1$ is of the form $(\gamma, \omega_1]$. Then $\Omega$ is a  Lindel\"{o}f space. To see this, given any open cover of $\Omega$ one of the elements, say $U$,  contains $\omega_1$. Then $U$ contains an interval  $(\gamma, \omega_1]$ for some $\gamma$. But then we only need to cover $[1,\gamma]$ which is a countable set, so we have a countable subcover and $\Omega$ is a  Lindel\"{o}f space.

But the subspace $\Omega_0$ is not  Lindel\"{o}f. Let $U_\alpha=[1,\alpha)$ for $\alpha\in\Omega_0$. Then the set of $U_\alpha$ is an open cover of $\Omega_0$. But there is no countable subcover. If there were, let $\{U_{\alpha1}, U_{\alpha2},\cdots\}$ be this subcover. Then we would have that the least upper bound of $\{\alpha1,\alpha2,\cdots\}$ is $\omega_1$. But any element of $\Omega$ larger than any number in this set would have a countable number of predecessors and would need to be strictly less than $\omega_1$. So we have a contradiction and $\Omega_0$ is not  Lindel\"{o}f.  
\end{example}

Now we move on to the much more important compact spaces. It should be obvious that a compact space is always Lindel\"{o}f. The first example shows that not every  Lindel\"{o}f space is compact. 

\begin{example}
The real line $R$ with the usual topology is not compact.Consider the open intervals $\{(-n,n)\}$ , where $n=1, 2, 3, \cdots$. Then these intervals cover $R$ but there is no finite subcover. If there were, then there would be a largest interval $(-N, N)$ and it would not include any $x\in R$ with $|x|>N$.
\end{example}

\begin{theorem}
\begin{enumerate}
\item A continuous image of a compact space is compact. 
\item Closed subspaces of  a compact space are compact.
\end{enumerate}
\end{theorem}

We get these for free from Theorem 2.7.2 by replacing "countable" with "finite" in the proofs. 

\begin{theorem}
A compact subset of a Hausdorff space is closed. 
\end{theorem}

\begin{example}
To see a compact subset of a space which is not closed, consider the two point space $X=\{a, b\}$ where $\{a\}$ is open but $\{b\}$ is not. Then  $\{a\}$ is not closed in $X$ but it is compact as there are only finitely many open sets in $X$. Of course, $X$ is not Hausdorff.
\end{example}

The next result which you may know from real analysis gives us lots of examples of compact sets:

\begin{theorem}
{\bf Heine-Borel Theorem} A subset of $R^n$ (with the usual topology) is compact if and only if it is closed and bounded. 
\end{theorem}

So in particular, the unit interval $I$ is compact as is the n-cube $I^n$. 

\begin{example}
To see a closed and bounded set which is not compact, let $X$ be the discrete metric space with infinitely many points. Recall that for any $x\neq y$, $\rho(x,y)=1$. Then $X$ is bounded as the distance between any two points is bounded by 1 and $X$ is closed as it is equal to a closed ball of radius 1. But $X$ is not compact. Take as the open cover, open balls of radius $\frac{1}{2}$ centered on each point of $X$. Then each of these sets consist of a single point and no subset of them (finite or not) can cover $X$. So $X$ is closed and bounded but not compact. 
\end{example}

Finally, products of compact spaces are particularly nice.

\begin{theorem}
{\bf Tychonoff Theorem} A product of compact spaces is compact.
\end{theorem}

\section{Connectedness}

{\it Everything grows together, because you're all one piece.}--Fred Rogers, Mister Rogers Neighborhood.

This brings up some serious questions. What does it mean to be one piece. Are you really in one piece? What about the cells, molecules, atoms, subatomic particles, etc. Fortunately, for topological spaces, it is easy to define being in one piece.

While we are on the subject of Mister Rogers, what do you call an open subset of the complex plane containing the number $i$? The Neighborhood of Make Believe.

Actually, connectedness and its variants form one of the most important concepts in algebraic topology. Here is the definition:

\begin{definition}
An topological space $X$ is {\it disconnected}\index{connected} if there are nonempty disjoint open subsets $H$ and $K$ of $X$ such that $X=H\cup K$. We say that $H$ and $K$ disconnect $X$. If these subsets do not exist, we say that $X$ is {\it connected}.
\end{definition}

We could replace open with closed in this definition. So a space $X$ is connected if the only subsets that are both open and closed are $X$ and $\emptyset.$

\begin{definition}
If $x\in X$, the largest connected subset of $X$ containing $x$ is called the {\it (connected) component}\index{component} of $x$. We will denote it as $C_x$. It is the union of all connected subsets of $x$ containing $X$.
\end{definition}

So you can think of components as the individual pieces of $X$. 

\begin{example}
Let $X\subset R$ where $X=[0,1]\cup (4,5)$. Then $[0,1]$ and $(4,5)$ are the two components of $X$.
\end{example}

If $x, y$ are two distinct points in $X$, then either $C_x=C_y$ or $C_x\cap C_Y=\emptyset$. If this was not the case, then $C_x\cup C_Y$ would be a connected set containing both $x$ and $y$ which is larger than both $C_x$ and $C_y$ which is impossible. So components of $X$ form a partition of $X$ into maximal connected subsets. (A partition means a collection of disjoint subsets whose union is all of $X$.) 

Here are some examples of connected and disconnected spaces.

\begin{example}
The Sorgenfrey line $E$ is the real line with an alternate topology. Here the basic neighborhoods of $x\in E$ are of the form $[x,z)$, where $z>x$. Then $E$ is disconnected. If we let $H=\cup_{x<0}[x,0)$ and $K=\cup_{z>0}[0,z)$ then $H$ and $K$ are both open, $H\cap K=\emptyset$, and $H\cup K=E$. So $E$ is disconnected.
\end{example}

\begin{example}
The discrete metric space is disconnected if it has at least 2 points. Since each point is an open set, any 2 subsets of the space form a disjoint union of open sets.
\end{example}

\begin{example}
The unit interval $I=[0,1]$ is connected. Suppose it was disconnected by $H$ and $K$. Let $1\in H$. Then $H$ contains a neighborhood of $1$, so let $c=\sup K$. (The term $\sup$ means least upper bound or the smallest number that is greater than or equal to everything in $K$. The greatest lower bound is denoted $\inf$.) We know that $c\neq 1$ since $H$ contains a neighborhood of $1$. Since $c$ belongs to $H$ or $K$ and both of these are open, some neighborhood of $c$ is contained in $H$ or $K$. But any neighborhood of $c$ must contain points of $K$ to the left and points of $H$ to the right. This is a contradiction so $I$ is connected. By the same argument, so is any closed interval in $R$.
\end{example}

Connected sets are preserved under continuous functions.

\begin{theorem}
The continuous image of a connected set is connected. 
\end{theorem}

{\bf Proof:} Let $X$ be connected and $f$ be a continuous function of $X$ onto $Y$.Then if $H$ and $K$ disconnect $Y$, $f^{-1}(H)$ and $f^{-1}(K)$ disconnect $X$ which is a contradiction. so $Y$ is connected. $\blacksquare$.

Now I will show you something really cool. Do you remember the Intermediate Value Theorem from calculus. Your book probably didn't even try to prove it. With topology, the proof is two lines.

\begin{corollary}
{\bf Intermediate Value Theorem:} Let $f$ be a real valued function of a single real variable which is continuous on $[a,b]$. (Without loss of generality we can assume $f(b)\geq f(a)$ or we can make the obvious modification.) Then for every $y$ such that $f(a)\leq y\leq f(b)$, there is an $x\in[a,b]$ such that $f(x)=y$. 
\end{corollary}

{\bf Proof:} Since $[a,b]$ is connected, $f([a,b])$ is also connected by the previous theorem. So if $f(a)\leq y\leq f(b)$ but $y$ is not in the image of $f$ then the image of $f$ is disconnected which would be a contradiction. $\blacksquare$.

The next set of results deal with subsets of connected spaces. These need not be connected as $X=[0,.1]\cup (.3,.5)$ is a subset of the connected set $I$ which is not connected.

\begin{definition}
Sets $H$ and $K$ are {\it mutually separated}\index{mutually separated} in $X$ if $$H\cap\overline{K}=K\cap\overline{H}=\emptyset.$$
\end{definition}

\begin{theorem}
A subspace $E$ of $X$ is connected if and only if there are no nonempty mutually separated sets $H$ and $K$ in $X$ with $E=H\cup K$.
\end{theorem}

\begin{corollary}
If $H$ and $K$ are mutually separated in $X$ and $E$ is a connected subset of $H\cup K$ then $E\subset H$ or $E\subset K$.
\end{corollary}

We now have some more ways of proving a space is connected.

\begin{theorem}
\begin{enumerate}
\item If $X=\cup X_\alpha$ where each $X_\alpha$ is connected and $\cap X_\alpha\neq\emptyset$, then $X$ is connected.
\item If each pair of points in $X$ lies in some connected subset of $X$ then $X$ is connected.
\item If $X=\cup_{i=1}^n X_i$ where each $X_i$ is connected and $X_i\cap X_{i+1}\neq\emptyset$ for all $i$, then $X$ is connected.
\end{enumerate}
\end{theorem}

\begin{theorem}
If $E$ is a connected subset of $X$ and $E\subset A\subset\overline{E}$ then $A$ is connected.
\end{theorem}

{\bf Proof:} It is enough to show that $\overline{E}$ is connected since if $E\subset A\subset\overline{E}$ then $A$ is equal to the closure in $A$ of $E$, so we can replace the closure in $X$ of $E$ (i.e. $\overline{E}$) by the closure in $A$ of $E$. Suppose $\overline{E}=H\cup K$  where $H$ and $K$ are disjoint nonempty open sets in  $\overline{E}$. Then $E=(H\cap E)\cup(K\cap E)$ and so the union of disjoint nonempty open sets in $E$. So if $\overline{E}$ is disconnected, so is $E$. $\blacksquare$.

\begin{example}
The real line $R$ is connected. To see this, we know that $R=\cup_{n=1}^\infty[-n, n]$, and each interval $[-n, n]$ is homeomorphic to $I$ and thus connected. Since their intersection is nonempty (it contains 0), part 1 of Theorem 2.8.5 shows that $R$ is connected.
\end{example}

\begin{example}
$R^n$ is connected by the same theorem. $R^n$ is the union of all straight lines through its origin, and each of these is homeomorphic to $R$ which is connected by the last example. Also, the origin lies in their intersection. So $R^n$ is connected.
\end{example}

\begin{example}
$R^2$ and $R^1$ are not homeomorphic. To see this, remove a point from both of them. Then $R^2$ is still connected, but $R^1$ becomes disconnected. When we do homology theory in Chapter 4, I will show you why $R^n$ and $R^m$ are not homeomorphic if $n\neq m$. We don't yet have the tools for it now, but it is an easy consequence of homology theory.
\end{example}

Connectivity behaves well with product spaces due to continuity of projection maps.

\begin{theorem}
A nonempty product space is connected if and only if each factor space is connected. 
\end{theorem}

This follows from the continuity of projection maps and the fact that the continuous image of a connected space is connected.

Finally, we have the following result on components of a space.

\begin{theorem}
The components of $X$ are closed sets.
\end{theorem}

{\bf Proof:} If $C_x$ is the component of $x$ in $X$, then $\overline{C_x}$ is a connected set containing $x$. So $\overline{C_x}\subset C_x$ and $C_x$ is closed.$\blacksquare$.

Now that we have defined connectedness, you may wonder if you can always get there from here. Is there a "path" connecting any two points? In homotopy theory you learn to love paths and especially "loops" which are paths that start and end at the same place. We will now define the notion of {\it pathwise connectedness}.

\begin{definition}
A space $X$ is {\it pathwise connected}\index{pathwise connected} if for any two points $x,y\in X$, there is a continuous function $f: I\rightarrow X$ such that $f(0)=x$ and $f(1)=y$. The function as well as its range $f(I)$ is called a {\it path}\index{path} from $x$ to $y$ with {\it initial point}\index{initial point} $x$ and {\it terminal point}\index{terminal point} $y$. A continuous function $f: I\rightarrow X$ such that $f(0)=f(1)=x$ is called a {\it loop}\index{loop}.
\end{definition}

You will frequently see the term {\it arcwise connected}\index{arcwise connected} which means that for a pair of points we can find a path is a homeomorphism onto its image. (So it doesn't intersect itself.) It turns out that a Hausdorff space is arcwise connected if and only if it is pathwise connected. (See \cite{Wil} for a proof.) Since we won't need spaces that are not Hausdorff in data science applications, we can treat these terms as interchangeable when you see them in books and papers. I will use the term pathwise connected in the rest of this paper.

\begin{theorem}
Every pathwise connected space is connected.
\end{theorem}

{\bf Proof:} If $H$ and $K$ disconnect the pathwise connected space $X$, let $f:I\rightarrow X$ be a path between $x\in H$ and $y\in K$. Then $f^{-1}(H)$ and $f^{-1}(H)$ disconnect $I$ which is a contradiction.$\blacksquare$.

\begin{example} 
Not every connected space is path connected.  Figure 2.8.1 shows the {\it topologist's sine curve}\index{topologist's sine curve} which is the union of the graph of $\sin(\frac{1}{x})$ for $x>0$ and the $x$-axis for $x\leq 0$. The right hand side, oscillates between -1 and 1 and the oscillation gets faster and faster as you approach 0.  This set is connected as any neighborhood of 0 contains points on both pieces of the graph, so it cannot be disconnected by open sets. But there is no path from $(0,0)$ to $(x, \sin(\frac{1}{x})$ for any $x>0$. So the topologist's sine curve is connected but not path connected.
\end{example}

\begin{example} 
In the discrete metric space, every point is both open and closed, so every point is its own connected component. Such a space is called {\it totally disconnected}\index{totally disconnected}. In this space, there is no path between any pair of points as this would disconnect the unit interval $I$. So if you lived in this space, even though you would have a short commute to work (only 1 unit), you would have no way to get there. And you wouldn't be all one piece, but infinitely many. 
\end{example}

\begin{figure}[ht]
\begin{center}
  \scalebox{1.4}{\includegraphics{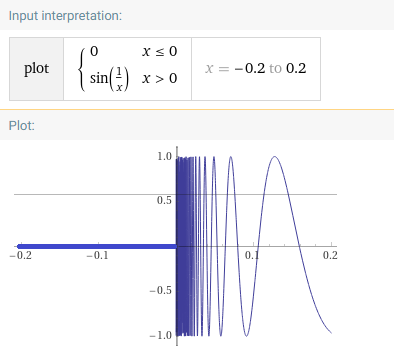}}
\caption{
\rm 
Topologist's Sine Curve
}
\end{center}
\end{figure}

Now you should know enough topology to handle algebraic topology. Next we will focus on the algebra.

\chapter{Abstract Algebra Background}

Like point-set topology, abstract algebra is a huge area. I took an entire year of it in college and an additional year in graduate school. Fortunately, the amount you need for algebraic topology is not nearly as extensive. So in this chapter, I will present the amount you need to have a basic understanding of homology theory. We can add a little more as we need it when we go into more advanced topics.

In topology, we give sets additional structure by defining a distinguished class of subsets, the open sets. In algebra, by contrast, the structure relates to operations where a subset of elements point to another element in a set. In a binary operation we take two elements and produce a third. A {\it group} is a set in which there is a single operation called addition or multiplication. It doesn't really matter which one as the operation may not resemble what we do to numbers. But we want to write the operations in a group in a familiar form, so for example, if we are adding, the identity element (which leaves the other fixed) is written as a 0, and in multiplication, it is written as a 1. The elements could be numbers but could also be rocks, flowers, or people as long as you can multiply or add them in a way which follows a fixed set of rules. So in Section 3.1, I will discuss the basics of groups. In algebraic topology, groups have an especially nice structure which is very easy to visualize.

Section 3.2 will deal with the idea of an {\it exact sequence}. This is a chain of groups and maps between them with some very nice properties. You can often compute the structure of groups in an exact sequence if you know the structure of their neighbors. In algebraic topology this is one of your most effective tools. 

In Section 3.3, we will introduce {\it rings} which have two binary operations, both addition and multiplication. Now there is a big difference between them. For example, if $r$ is and element of a ring $R$, then $0+r=r$, but $0r=0.$ What ties them together is a {\it distributive law}. A very special type of ring is a {\it field}. The real, rational, and complex numbers are all fields. An important example of a ring that is not a field is the integers.

In Section 3.4, we introduce a new operation called scalar multiplication. A vector space is a group under addition with a way to multiply elements by elements of an external field called the field of {\it scalars}. Vector spaces should be familiar to you if you have ever taken a class in linear algebra. If we use a ring instead of a field, we get a {\it module}. So a vector space is just a special case of a module. Therefore, a vector space on the moon is a lunar module. If you have addition, multiplication, and scalar multiplication, you get an algebra. Square matrices are an example.

So now you are going to object and say, "None of this is abstract enough for me." Well you don't have to worry because we will finish this chapter off with a section on {\it category theory}. As I pointed out when we were discussing set theory, Russell's paradox doesn't allow for the set of all sets. So we can talk about the category of sets. Categories consist of {\it objects} and special maps between them called {\it morphisms}. Examples are the categories of sets, topological spaces, groups, and rings. Maps between categories that take objects to objects and morphisms to morphisms are called {\it functors}. Homology is a functor that takes a geometric object (actually a topological space) and replaces it by a group in each dimension up to the highest dimension of the space. (I will make this more precise in Chapter 4, but we can always think of these spaces as subsets of some $R^n$ so $n$ is their dimension.)

There are lots of textbooks on abstract algebra, but my personal favorite to start with is Herstein's {\em Topics in Algebra} \cite{Her1}. This book seems to be out of print and hard to find now. It has been replaced by the posthumously published {\em Abstract Algebra, 3rd Edition} \cite{Her2}, which covers somewhat less material. Most of Sections 3.1, 3.2, and 3.4 will come from  \cite{Her1}, with some additional material from Munkres \cite{ Mun1}. Section 3.2 is entirely from \cite{Mun1}. Another one of my favorites is Jacobson's {\em Basic Algebra}. The only problem is that the name is a misnomer and it is almost embarrassing to walk around with it as people will think I never learned my high school algebra. In actuality, though, if you know all the material in this two volume set, you will be well prepared for any Ph.D. qualifying exam in algebra. There is a good short chapter on category theory in the second volume which I will use in Section 3.5. If you want to learn more, the classic book by one of the inventors of category theory is MacLane's {\em Categories for the Working Mathematician} \cite{MacL1}. 

\section{Groups}

Since as mentioned in the introduction, a group's operation can either be multiplication or addition, I will define it in both cases.

\begin{definition}
A nonempty set $G$ is a {\it group}\index{group} if there is a binary operation {\it sum} which is denoted by $+$  or {\it product} which is denoted by $\cdot$ which takes two elements of $G$ and produces a third and obeys the rules listed in Table 3.1.1. (Here we write $ab$ for $a\cdot b$.)
\end{definition}

\begin{table} 
\begin{center}
\begin{tabular}{|c|c|}
\hline
Rule  & Addition\\
\hline\hline
Closure & $a, b\in G$ implies that $a+b\in G$\\
\hline
Associative Law & $a+(b+c)=(a+b)+c$ \\
\hline
Identity & There  exists an element $0\in G$ such that for every $a\in G$, $0+a=a+0=a$ \\
\hline
Inverse & For every $a\in G$ there exists $-a\in G$ such that $a+(-a)=(-a)+a=0$  \\
\hline\hline
Rule  & Multiplication\\
\hline\hline
Closure & $a, b\in G$ implies that $ab\in G$\\
\hline
Associative Law & $a(bc)=(ab)c$\\
\hline
Identity & There  exists an element $1\in G$ such that for every $a\in G$, $1a=a1=a$\\
\hline
Inverse  & For every $a\in G$ there exists $a^{-1}\in G$ such that $aa^{-1}=a^{-1}a=1$ \\
\hline
\end{tabular}
\caption{Group Laws.}
\end{center}
\end{table}

Arranging the rules in a table shows the explicit parallels betwen groups under addition and multiplication. They are the same thing other than a difference in notation. 

Note that there is something missing from the rules you learned about numbers in high school algebra. There is no commutative law. This law only holds for a special type of group.

\begin{definition}
A group $G$ is called an {\it abelian group}\index{abelian group} if the commutative law holds. It is written $a+b=b+a$ for addition or $ab=ba$ for multiplication.A group which is not abelian is a {\it non-abelian group.}
\end{definition}

What is purple and commutes? An abelian grape.

A group can be either finite or infinite. If $G$ is finite, the number of elements in $G$ is called the {\it order} of $G$ and is written as $|G|$. 

\begin{example} 
The trivial group having only one element is $\{1\}$ written multiplicatively or $\{0\}$ written additively.We will write these simply as 1 and 0.
\end{example}

\begin{example} 
Let $G=Z$, the set of integers. Then $G$ is a group under addition. In fact, it is an abelian group. But $Z$ is not a group under multiplication. No integer other than $1$ and $-1$ have inverses that are also integers. The nonzero integers are an example of a {\it monoid}\index{monoid}. You may remember a monoid as the one-eyed giant that Odysseus met. (Or was that a Cyclops?) Actually a monoid is a similar to a group except that it might not have inverses. As we will see, a ring is an abelian group under addition and its nonzero elements form a monoid under multiplication.
\end{example}

\begin{example} 
Let $G=R$, the set of real numbers. Then $R$ is a group under addition and the nonzero elements of $R$ form a group under multiplication. The same holds if $R$ is replaced by the set $Q$ of rationals.The irrationals are not a group under addition and the nonzero irrationals are not a group under multiplication as they are missing both $0$ and $1$.
\end{example}

\begin{example} 
An example of a finite group is the integers mod $n$ written as $Z_n$. $Z_n=\{0, 1, 2, \cdots, n-1\}$ which is a group under addition. Addition is done in the usual way, but $(n-1)+1=0$. Think of a clock with a $0$ where the 12 usually goes. Then for example, $11+3=2$. 
\end{example}

\begin{example} 
An example of a non-abelian group is the set of invertible $n\times n$ matrices with real entries for a fixed $n$. The operation is matrix multiplication which we know is not commutative. The set of all $n\times n$ matrices with real entries is a group under addition.
\end{example}

\begin{example} 
A more exotic yet important group is the {\it permutation group}\index{permutation group} $S_n$. (Don't confuse this with the sphere $S^n$.) This group consists of all possible permutations of $n$ objects where permutations are multiplied by taking them in succession. $|S_n|=n!$ as can be seen from the fact that a permutaion can put the first object in $n$ positions, the second in $n-1$ positions, etc. Also, this group is non-abelian for $n\geq 3$. In fact, $S_3$ is the smallest non-abelian group.

As an example, suppose you are the zookeeper in a small zoo which has only 3 animals: a lion, an elephant, and a giraffe. Your boss likes to constantly move animals around to confuse any potential visitors. So he gives you a card every day with instructions on moving animals. The cards are written in {\it cycle} notation. If two objects are enclosed in the same set of parentheses, you move each object to the one on its right and circle around if necessary. So for example, if the cages are in a row in the order [Lion, Elephant, Giraffe] or $[L,E,G]$ and the card says $(E) (L,G)$ this means to leave the elephant in its place and switch the lion and the giraffe resulting in the new configuration $[G, E, L]$. We usually omit a cycle with one element so we write the previous permutation as $(L,G)$. If we have the cycle $(L)$ or equivalently either of the other two letters, this means leave everything alone. We multiply permutations by doing then in succession so $(E,L)(L,E,G)$ means first switch the elephant and the lion, then lion to elephant, elephant to giraffe, and giraffe to lion. So $[L,E,G]$ becomes $[E, L, G]$ and then $[L, G, E]$. Note that reversing the order is $(L,E,G)(E,L)$ so we do the circular switch first and then switch the elephant and the lion. So $[L, E, G]$ becomes, $[G,L, E]$ and then $[G,E,L].$ So reversing the order gives a different configuration, and $S_3$ is non-abelian. To show it is an actual group, we can undo a cycle with two elements by doing it again. So applying $(L,E)$ twice switches the lion and the elephant and then switches them back. We undo $(L,E, G)$ with $(G, E, L)$.
\end{example}

Note that $Z_n$ is a group formed by taking an element and all multiples until one of them becomes 0. An equivalent multiplicative group is the group $G=\{1, a, a^2, a^3, \cdots, a^{n-1}\}$ where $a(a^{n-1})=1.$ This type of group is generated by a single element and is called a {\it cyclic group}\index{cyclic group} of order $n$. The integers form an {\it infinite cyclic group} under addition. An equivalent multiplicative group is a group of the form $\{g^i\}$ where $i$ is any integer and $g^0=1$. I will make precise what I mean by {\it equivalent} a little later.

What do you call 5 friends who all like to ride bicycles together? A cyclic group of order 5.

Recall that topological spaces had subspaces, product spaces, and quotient spaces. I will now describe the analogues for groups. These will be critical in algebraic topology. I will start with subgroups. Unlike the case of topological spaces, not every subset of a group is a subgroups. We need this subset to be a group itself.

\begin{definition}
A nonempty subset $H$ of a group $G$ is a {\it subgroup}\index{subgroup} of $G$ if it is itself a group.  
\end{definition}

$H$ inherits the binary operation from $G$ so it is additive if $G$ is additive and multiplicative if $G$ is multiplicative.

I am going to use multiplicative notation for the next part of the discussion. Later I will assume we are dealing with abelian groups which are traditionally represented as additive. Homology and cohomology only involve abelian groups. Non-abelian groups can come up in homotopy theory, though, so I will have more to say about them when we get to Chapter 9. I will always specify when the groups being discussed are restricted to being abelian.

So what does $H$ need to be a group itself? It will automatically have the associative law. So we need closure, the identity, and inverses. We can drop the requirement that the identity be in $H$ if we have closure and the existence of inverses since $h\in H$ implies $h^{-1}\in H$, so by closure, $h\cdot h^{-1}=1\in H$. So we have shown:

\begin{theorem}
A nonempty subset $H$ of a group $G$ is a subgroup of $G$ if and only if $a, b\in H$ implies $ab\in H$ and $a\in H$ implies $a^{-1}\in H$.  
\end{theorem}

\begin{example} 
If $Z$ is the integers under addition, then let $3Z$ be the multiples of 3. Then since the sum of multiples of 3 is always a multiple of 3, and $a$ is a multiple of 3 if and only if $-a$ is a multiple of $3$, we have that $3Z$ is a subgroup of $Z$.
\end{example}

\begin{example} 
If $Z$ is the integers under addition, then the even numbers are equal to $2Z$ so they form a subgroup by the argument in the last example. The odd numbers, do not form a subgroup. For example, 3 and 5 are odd, but 3+5=8 which is an even number, so closure fails. And the odd numbers do not contain the identity element 0.
\end{example}

\begin{example} 
The integers form a subgroup of the reals under addition.
\end{example}

\begin{example} 
Let $G$ be the nonzero real numbers under multiplication. Then the nonzero rationals form a subgroup as do the nonzero positive rationals.
\end{example}

If our group is finite, there is an important relationship between the order of $G$ and the order of its subgroups. We start with the idea of a {\it coset}.

\begin{definition}
If $H$ is a subgroup of $G$, and $a\in G$ then we define the {\it right coset}\index{coset} $Ha=\{ha|h\in H\}$. The left coset aH has an analogous definition. If $G$ is additive, cosets are of the form $H+a$.
\end{definition}

\begin{theorem}
Any two right cosets of $H$ in $G$ are either identical or disjoint. Since $H=H1$, $H$ is itself a coset and the theorem implies that the cosets partition $G$.
\end{theorem}

{\bf Proof:} Suppose $Ha$ and $Hb$ have a common element $c$. Then $c=h_1a=h_2b$ for some $h_1, h_2\in H$. So $h_1^{-1}h_1a=h_1^{-1}h_2b$. So $a=h_1^{-1}h_2b$ which implies that $a\in Hb$. This implies $Ha\subset Hb$. A similar argument shows $Hb\subset Ha$. So $Ha=Hb$. $\blacksquare$

Note that the above argument shows that $Ha=Hb$ if $ab^{-1}\in H$.

Now every coset is the same size as $H$. This is true since for any two distinct elements of $h_1, h_2 \in H$, $h_1a$ and $h_2a$ are distinct as we can see by multiplying both by $a^{-1}$ on the right. So $G$ is a union of disjoint sets of size $|H|$. So we have proves the very important result: 

\begin{theorem}
{\bf Langrange Theorem} If $G$ is a finite group, then $H$ is a divisor of $G$.
\end{theorem}

So for a group of prime order, the only subgroups are $G$ itself and the trivial subgroup $1$.

\begin{definition}
If $H$ is a subgroup of $G$, then the {\it index} of $H$ in $G$ written $[G:H]$ is the number of right cosets of $H$ in $G$.
\end{definition}

\begin{theorem}
If $H$ is a subgroup of $G$, then $[G:H]=|G|/|H|$.
\end{theorem}

\begin{example} 
Let $G=Z_{10}$ and $H=\{0,5\}$. Then $|G|=10$ and $|H|=2$. $H$ is a subgroup of $G$ since in $Z_{10}$, $5+5=0$. The cosets of $H$ are $H+0=\{0,5\}$,  $H+1=\{1,6\}$,  $H+2=\{2,7\}$,  $H+3=\{3,8\}$,  $H+4=\{4,9\}$.These exhaust $Z_{10}$ and there are 5 of them, so $[G:H]=5=10/2=|G|/|H|$ as we expect. Note that for an additive group, $H+a=H+b$ if $a-b\in H$. For example $H+2=\{2,7\}$, and $H+7=\{7,2\}=H+2$. We check that $7-2=5\in H$.
\end{example}

\begin{example} 
An infinite group can have a subgroup with finite index. An easy example is letting $G$ be the integers and $H$ the even integers under addition. $H+1$ is the set of odd integers, and $G=H\cup (H+1)$. So we have that the index is $[G:H]=2.$
\end{example}

We will talk about the order of an element of a group.

\begin{definition}
If $G$ is a group and $g\in G$, then the {\it order} of $g$, written $o(g)$ is the smallest power $m$ of $g$ such that $g^m=1$. This is finite for every element if $g$ is finite. the analogue for an additive group is $mg=0$.
\end{definition}

If $g\in G$ has finite order $m$, then $\{1, g, g^2, \cdots, g^{m-1}\}$ is a subgroup of $G$. So we have:

\begin{theorem}
If $G$ is a finite group and $g\in G$ has order $m$, then $m$ divides $|G|$.
\end{theorem}

If a group is abelian, left and right cosets coincide. In general, though, we want to only look at subgroups where the cosets behave nicely. Such a subgroup is called {\it normal}.

\begin{definition}
A subgroup $N$ of a group $G$ is called {\it normal}\index{normal subgroup} if for any $g\in G$, $n\in N$, we have $gng^{-1}\in N$. 
\end{definition}

Be careful. We are not saying that $gng^{-1}=n$. That would happen if the group was abelian. But we do have the following:

\begin{theorem}
$N$ is a normal subgroup of $G$ if and only if $gNg^{-1}=N$ for all $g\in G$.
\end{theorem}

This says that for $n\in N$, $gng^{-1}=n_1$ for some $n_1\in N$. For a normal subgroup, a right coset $Ng$ is the left coset $gN$ for the same $g$. This gives us a way to multiply right cosets. For a normal subgroup $N$ of $G$ and $g_1, g_2\in G$, $Ng_1Ng_2=N(g_1N)g_2=N(Ng_1)g_2=Ng_1g_2.$ With this multiplication, $N=N1$ the identity, and the inverse of $Ng$ equal to $Ng^{-1}$, the right cosets form a group. 

\begin{definition}
If $N$ is a normal subgroup of $G$, the set of right cosets with the multiplication written above is called a {\it quotient group}\index{quotient group} and written $G/N$. By our comments above, $|G/N|=|G|/|N|$
\end{definition}

I can't overemphasize how important quotient groups are in algebraic topology. Homology is defined as a quotient group, and every TDA paper I have ever seen assumes you already know what a quotient group is. 

If $G$ is written additively, then the cosets are of the form $N+a$. $N+0=N$ is the identity, $(N+a)+(N+b)=N+(a+b)$, and the inverse of $N+a$ is $N+(-a)=N-a$. 

If $G$ is abelian, then every subgroup is normal but this is not true if $G$ is non-abelian. 

\begin{example} 
Let $G=Z_{10}$ and $H=\{0,5\}$. We listed the cosets of $H$ in $G$ in Example 3.1.11 above. The quotient group $G/H$ is formed by adding the cosets $(H+a)+(H+b)=H+(a+b)\mod 10$. (Recall $a\mod b$ is the remainder you get when dividing $a$ by $b$.)
\end{example}

\begin{example} 
Let $G=Z$ under addition, and $H=nZ$ which is the set of multiples of $n$. We saw that $H$ is a subgroup of $G$. The cosets of $H$ are $H+0=H, H+1, H+2, \cdots, H+(n-1)$. These exhuast $Z$ since a multiple of $n$ plus $n$ is still a multiple of $n$, so $H+n=H$. We can write out $G/H$ as $\{0,1,2,\cdots,n-1\}$ so we see that $Z/nZ$ is just $Z_n$. That is why you will often see $Z_n$ written as $Z/nZ$ or even $Z/n$ in books. I will stick to the $Z_n$ notation.
\end{example}

\begin{definition}
A group $G$ is {\it simple}\index{simple group} if it has no normal subgroups except $G$ itself and the trivial group $1$.
\end{definition}

An abelian group can only be simple if it is a cyclic group of order $p$ where $p$ is a prime. Non-abelian groups can be simple and very complicated. There was a big project in hte 1960s-1980s to classify all finite simple groups, and the proof is about 10,000 pages long. There are a few families and 26 sporadic ones. The largest of these is the Fischer-Griess Monster Group, also known as the ''Friendly Giant" or simply the "Monster."  Its order is 808,017,424,794,512,875,886,459,904,961,710,757,005,754,368,000,000,000.  If you are curious about where it comes from see \cite{Gri} for the original 100 page construction or \cite{Con} for the simpler 30 page version.

Now that we have seen subgroups and quotient groups, you may wonder if there is an analogue to Cartesian products. There are a few constructions that are useful but the one we will be concerned with is a {\it direct product}\index{direct product}.

\begin{definition}
Let $G$ and $H$ be groups. Then the direct product of $G$ and $H$ denoted $G\times H$ has an underlying set which is the cartesian product of the underlying sets of $G$ and $H$ and $(g_1, h_1)(g_2, h_2)=(g_1g_2, h_1h_2)$. If $G$ and $H$ are written additively, we write $G\oplus H$ and use the term {\it direct sum}\index{direct sum}.
\end{definition}

We can think of $G$ and $H$ as subgroups of $G\times H$ by writing $G=\{(g,1)\}$ for $g\in G$ and $H=\{(1,h)\}$ for $h\in H$. Then we have the following properties:
\begin{enumerate}
\item $G\cap H=1$
\item Every element of $G\times H$ can be written uniquely as the product of an element of $G$ and an element of $H$.
\item Every element of $G$ commutes with every element of $H$.
\item Both $G$ and $H$ are normal subgroups of $G\times H$.
\item The order of $G\times H$ is the product of the order of $G$ and the order of $H$.
\end{enumerate}

\begin{example} 
There can only be one group of order 2, order 3, and order 5. These are the cyclic groups $Z_2$, $Z_3$, and $Z_5$ respectively. We know that there are no others as these are of prime order so the can not have any nontrivial subgroups.
\end{example}

\begin{example} 
There are two distinct groups of order 4. The group $Z_4$ is the familiar cyclic group. We also have the group $Z_2\oplus Z_2$. This group contains the identity $(0,0)$ and three elements of order two: $(1,0)$, $(0,1)$, and $(1,1)$ have order 2, since in $Z_2$, $1+1=0.$This group is called the {\it Klein four group.}
\end{example}

\begin{example} 
For order 6, there are 3 groups. They are $Z_6$, $Z_3\oplus Z_2$ and the non-abelian permutation group $S_3$ that we saw earlier. As we will see later, though, $Z_6$, $Z_3\oplus Z_2$  can be thought of as the same group.
\end{example}

Now we want to know what it means for two groups to be the same. The analogue of a continuous function from topology is a {\it homomorphism}. This preserves the algebraic structure, this time in the forward rather than the reverse direction. The analogue of a homeomorphism is an {\it isomorphism}. Two groups that are isomorphic are basically the same group. I will now make these concepts precise.

\begin{definition}
A mapping $f$ from a group $G$ to a group $H$ is a {\it homomorphism}\index{homomorphism} if $f(g_1g_2)$=$f(g_1)f(g_2)$ for all $g_1, g_2\in G$. (In additive notation, $f(g_1+g_2)$=$f(g_1)+f(g_2)$.) A homomorphism $f$ is a {\it monomorphism}\index{monomorphism} if it is one to one. It is a  {\it epimorphism}\index{epimorphism} if it is onto. If the homomorphism $f$ is both one to one and onto (implying the existence of an inverse), then $f$ is called an {\it isomorphism}\index{isomorphism} and we say that $G$ and $H$ are {\it isomorphic}. 
\end{definition}

Two groups that are isomorphic are pretty much indistinguishable in the way that two topological spaces are indistinguishable if they are homeomorphic.

\begin{example} 
One easy example of a homomorphism is defined as  f(g)=1 for all $g\in G$. In additive notation, $f(g)=0$, and this is called the {\it zero homomorphism}. We see it often in algebraic topology. 
\end{example}

\begin{example} 
If $f: G\rightarrow G$, then $f$ is the {\it identity homomorphism} if $f(g)=g$ for $g\in G$.
\end{example}

\begin{example} 
Let $G$ be the group of real numbers under addition and $H$ the group of nonzero real numbers under multiplication. Let $f: G\rightarrow H$, be defined by $f(x)=2^x$. This is a homomorphism since $f(a+b)=2^{a+b}=2^a2^b=f(a)f(b).$
\end{example}

\begin{example} 
 Let $f: Z\rightarrow Z$, be defined by $f(x)=2x$. This is a homomorphism since $f(a+b)=2(a+b)=2a+2b=f(a)f(b).$
\end{example}

\begin{example} 
 Let $G=\{0,1\}=Z_2$ under addition and $H=\{1,-1\}$ under the usual multiplication. Let $f: G\rightarrow H$ with $f(0)=1$ and $f(1)=-1$. This is one to one and onto and you can check that f is a homomorphism. For example, $f(1+1)=f(0)=1$ and $f(1)f(1)=(-1)(-1)=1.$ Checking the other possibilities gives that $f(a+b)=f(a)f(b)$. So f is a homomorphism and we have defined it to be one to one and onto. So $G$ and $H$ are isomorphic.
\end{example}

\begin{example}
As mentioned above, $Z_6$ is isomorphic to $Z_3\oplus Z_2$. To see this, let $g\in Z_6$. Let $f(g)=(a,b)\in Z_3\oplus Z_2$, where $g=3b+a$. (I.e. $b$ and $a$ arer the quotient and remainder respectively when $g$ is divided by $3$. We have $f^{-1}((a,b))=3b+a$. It is easily checked that $f(g+h)=f(g)+f(h)$ and that $f^{-1}$  is actually the inverse of $f$. So $Z_6$ is isomorphic to $Z_3\oplus Z_2$.
\end{example}

\begin{definition}
Let $f: G\rightarrow H$ be a homomorphism. Then the {\it kernel}\index{kernel} $K$ of $f$ is $\{g\in G| f(g)=1\}$. (In additive notation $f(g)=0$.) The {\it image} $f(G)$ of $f$ is $\{h\in H| h=f(g){\rm \hspace{.02 in} for\hspace{.02 in}some\hspace{.02 in}} g\in G\}$. 
\end{definition}

\begin{theorem}
Let $f: G\rightarrow H$ be a homomorphism. Then $f(1_G)=1_H$ and $f(g^{-1})=f(g)^{-1}$, where $1_G$ and $1_H$ are the identity elements of $G$ and $H$ respectively.
\end{theorem}

{\bf Proof:} Let $g\in G$. Then $f(g)=f(g1_G)=f(g)f(1_G)$. Multiplying both sides by $f(g)^{-1}$ on the left gives $1_H=f(1_G)$. So $1_H=f(1_G)=f(gg^{-1})=f(g)f(g^{-1})$ So multiplying on the left by $f(g)^{-1}$ on both sides gives $f(g)^{-1}=f(g^{-1})$. $\blacksquare$

Using the above result it is easy to show the following: 

\begin{theorem}
Let $f: G\rightarrow H$ be a homomorphism. Then the kernel of $f$ is a normal subgroup of $G$.
\end{theorem}

Try this one yourself using the formulas.

\begin{theorem}
Let $f: G\rightarrow H$ be a homomorphism. Then $f$ is one to one if and only if the kernel of $f$ is 1.
\end{theorem}

{\bf Proof:} Supose the kernel of $f$ is 1. Let $g_1, g_2\in G$, and suppose $f(g_1)=f(g_2)$. Then $1=f(g_1)f(g_1)^{-1}=f(g_1)f(g_2)^{-1}=f(g_1g_2^{-1}).$. Since the kernel of $f$ is 1, $g_1g_2^{-1}=1$, so $g_1=g_2$ and $f$ is one to one. If $f$ is one to one and $g$ is in the kernel of $f$, then  $1=f(g)=f(1)$ so we must have $g=1$. $\blacksquare$

The next theorem is very important. See \cite{Her1} for the proof.

\begin{theorem}
Let $f: G\rightarrow H$ be a homomorphism with kernel $K$. Then $G/K\approx H$. (Here $\approx$ denotes an isomorphism.)
\end{theorem}

Now let $G$ and $H$ be abelian. We will always write abelian groups additively. So the above result says that $f$ is one to one if the kernel $K=0.$ Also, recall that every subgroup of an abelian group is normal. So the image $f(G)$ is a normal subgroup of $H$ and can form the quotient group $H/f(G)$. We call this group the {\it cokernel}\index{cokernel} of $f$.If $f$ is onto, then $f(G)=H$, so $H/f(G)=H/H=0$. So we have:

\begin{theorem}
Let $f: G\rightarrow H$ be a homomorphism of {\bf abelian} groups. Then $f$ is onto if and only if the cokernel of $f$ is 0. $f$ is an isomorphism if the kernel and the cokernal of $f$ are both 0.
\end{theorem}

In both homology and cohomology theory, all of the groups we care about have an especially nice form. To understand this, we will need a couple more definitions.

\begin{definition}
Let $G$ be a group. A set $S$ of {\it generators}\index{generators} for $G$ is a set of elements of $G$ such that every element of $G$ can be expressed as a product of finitely many of the elements of $S$ and their inverses. If $G$ is generated by a finite set $S$, then $G$ is {\it finitely generated}\index{finitely generated}. If $S$ consists of a single element $g$ then $G$ is a cyclic group generated by $g$.
\end{definition}

\begin{definition}
Let $G$ be a group. A set $R$ of relations is a set of strings $g_1^{r_1}g_2^{r_2}\cdots g_n^{r_n}$ of generators $g_i$ and integers $r_i$ whose product is 1. Every group has the relation $gg^{-1}=1$ for every generator $g$, but there may be others as well. A group with no other relations is called a {\it free group}\index{free group}. A group with only finitely many relations is called {\it finitely presented}\index{finitely presented}.
\end{definition}

\begin{example} 
 Let $G$ be the free group with two generators $x$ and $y$. Then elements of $G$ are the element 1 along with strings of $x$, $y$, $x^{-1}$ and $y^{-1}$. Some examples are $x^8$, $x^6y^3$, and $x^3({y^{-1}})^5y^3x^8$.
\end{example}

In homology and cohomology, every group we will ever deal with is finitely generated and abelian. In homotopy theory, there is one exception called the {\it fundamental group} that could be non-abelian. We will deal with it in Chapter 9. 

What's purple, commutes, and is worshiped by a limited number of people? A finitely venerated abelian grape.

What's purple, commutes, and given gifts by a limited number of people? A finitely presented abelian grape.

\begin{definition}
Let $G$ be a finitely generated group.Then $G$ is a {\it free abelian group}\index{free abelian group}, if $G=Z\oplus Z\oplus\cdots\oplus Z$. The number of summands of $Z$ is called the {\it rank} of $G$. In algebraic topology, the rank is also called the {\it betti number.}\index{betti number}
\end{definition}

Warning: Do not confuse a free group and a free abelian group. A free group is about as non-abelian as you can get.

A free abelian group of rank $n$ has a {\it basis} $\{g_1,\cdots, g_n\}$ which is defined to be a set of elements of $G$ such that every element of $G$ can be written uniquely as a sum $\sum_{i=1}^n c_ig_i$ where $c_i\in Z$ for $1\leq i\leq n.$

One way of constructing a free abelian group from a finite set of generators is called a {\it formal sum}\index{formal sum}. In this case, we start with a generating set $S=\{g_1,\cdots, g_n\}$ and a function $\phi: S\rightarrow Z$, so the elements of our group are of the form $\sum_{i=1}^n\phi(g_i)g_i$

\begin{example} 
Suppose you go to a fruit stand that has apples, bananas, and cherries. You buy 2 apples, 4 banananas and 12 cherries. This represents the element $2a+4b+12c$, where $\{a, b, c\}$ is the generating set. You could also have an element like $-6a$ where you sell them 6 apples. You can add and subtract peieces of fruit of the same kind, but for different kinds, we just put a plus or minus sign between them.
\end{example}

When we build homology groups, we will start with the group of {\it chains} which will be a set of formal sums. We will discuss these in Chapter 4.

Recall that for an additive group, the order of an element $g$ is $m$ if $m$ is the smallest positive integer such that $mg=0$. The set of all elements of $G$ of finite order is called the {\it torsion subgroup}\index{torsion subgroup}. If $T=0$, we say that $G$ is {\it torsion free}\index{torsion free}.

We conclude this section with an important result which is the basis of homology theory. The good news is that it means that homology groups take a very simple form. See \cite{Mun1} for a proof.

\begin{theorem}
{\bf The Fundamental Theorem of Finitely Generated Abelian Groups:} Let $G$ be a finitely generated abelian group. Then we can write $G$ uniquely as $$G\equiv(Z\bigoplus\cdots\bigoplus Z)\bigoplus Z_{t_1}\bigoplus\cdots\bigoplus Z_{t_k}$$ with $t_i>1$ and $t_1|t_2|\cdots|t_k$ where $a|b$ means that $a$ divides $b$. The number of summands $Z$ is called the {\it betti number} $\beta$ of $G$. The {\it torsion subgroup} is $Z_{t_1}\bigoplus\cdots\bigoplus Z_{t_k}$.
\end{theorem}

\begin{example} 
$G\equiv(Z\bigoplus Z\bigoplus Z)\bigoplus Z_3\bigoplus Z_6 \bigoplus Z_{12}$. Then $g$ is the set of $6-$tuples $(a, b, c, d, e, f, g)$ where addition is coordinatewise and $a, b, c\in Z$, $d\in Z_3$, $e\in Z_6$, $f\in Z_{12}$. A typical element could be $(1, 5, -3, 2, 4, 8)$. Note that $3|6|12$.
\end{example}

Note that for $m$ and $n$ relatively prime, $Z_m\bigoplus Z_n\equiv Z_{mn}$ (prove this yourself), so any finite cyclic group can be written as a direct sum of cyclic groups whose orders are powers of primes. So we can also write $G$ uniquely as $$G\equiv(Z\bigoplus\cdots\bigoplus Z)\bigoplus Z_{a_1}\bigoplus\cdots\bigoplus Z_{a_p}$$ where the $a_i$ are powers of distinct primes.

So that was a rather long section, but it barely scratched the surface of group theory. It is pretty much most of what you need to know for algebraic topology. The next section will be a first taste of a branch of abstract algebra called {\it homological algebra} which is specialized for algebraic topology while being a subject in its own right. You will see a lot more of it in Chapter 4. 

\section{Exact Sequences}

Exact sequences are a key tool in computation in homology, cohomolgy, and homotopy theory. This topic really belongs in Chapter 4, but I am including it here as it illustrates some of the group theory concepts we have just learned.

In this section, all groups will be abelian and written additively, but there is no reason we can't have an exact sequence of non-abelian groups.

\begin{definition}
Consider a sequence (finite or infinite) of abelian groups and homorphisms.

\begin{tikzpicture}
  \matrix (m) [matrix of math nodes,row sep=3em,column sep=4em,minimum width=2em]
  {
   \cdots & A_1 & A_2 & A_3 & \cdots\\};
  \path[-stealth]
    (m-1-1) edge[dashed]  (m-1-2)
    (m-1-2)  edge node [above] {$\phi_1$} (m-1-3)
(m-1-3)  edge node [above] {$\phi_2$} (m-1-4)
(m-1-4) edge[dashed]  (m-1-5);

\end{tikzpicture}

The sequence is {\it exact} at $A_2$ if image $\phi_1$=kernel $\phi_2$. If this happens at every group (except the left and right ends if they exist), then the sequence is an {\it exact sequence}\index{exact sequence}.
\end{definition}

We will use some of what we have learned in the last section about homomorphisms to prove some facts about exact sequences.

\begin{theorem}
The sequence

\begin{tikzpicture}
  \matrix (m) [matrix of math nodes,row sep=3em,column sep=4em,minimum width=2em]
  {
A_1 & A_2 & 0\\};
  \path[-stealth]
    (m-1-1)  edge node [above] {$\phi$} (m-1-2)
(m-1-2)  edge  (m-1-3);
\end{tikzpicture} is exact if and only if $\phi$ is an epimorphism.
\end{theorem}

{\bf Proof:} Since the second map is the zero homomorphism, its kernel is all of $A_2$. So the image of $\phi$ is also all of $A_2$,and $\phi$ is an epimorphism. If $\phi$ is an epimorphism, then the image of $\phi$ is all of $A_2$, so the second map has all of $A_2$ as its kernel and is the zero homomorphism. $\blacksquare$.

\begin{theorem}
The sequence

\begin{tikzpicture}
  \matrix (m) [matrix of math nodes,row sep=3em,column sep=4em,minimum width=2em]
  {
0 & A_1 & A_2 \\};
  \path[-stealth]
    (m-1-1)  edge (m-1-2)
(m-1-2)  edge  node [above] {$\phi$} (m-1-3);
\end{tikzpicture} is exact if and only if $\phi$ is an monomorphism.
\end{theorem}

{\bf Proof:} Since the first map is the zero homomorphism, its image is 0. So the kernel of $\phi$ is 0, and $\phi$ is a monomporphism by Theorem 3.1.9.  If $\phi$ is a monomorphism, then the kernel of $\phi$ is 0 so the first map is the zero homomorphism. $\blacksquare$.

You should be very happy whenever you see this as a piece of an exact sequence:

\begin{tikzpicture}
  \matrix (m) [matrix of math nodes,row sep=3em,column sep=4em,minimum width=2em]
  {
   0 & A_1 & A_2 & 0\\};
  \path[-stealth]
    (m-1-1) edge  (m-1-2)
    (m-1-2)  edge node [above] {$\phi$} (m-1-3)
(m-1-3)  edge (m-1-4);

\end{tikzpicture}

By the last two theorems, the map $\phi$ is both a monomorphism and an epimorphism. So $\phi$ is an isomorphism and $A_1\equiv A_2$, and if you already know one of the groups, you also know the other. This is a very common situation in algebraic topology and this gives exact sequences a lot of their power.

I will give you two more results to try to prove on your own. We will need two new definitions.

\begin{definition}
Let $G$ be an abelian group with subgroup $H$. Then the homomorphism $i: H\rightarrow G$ is an {\it inclusion} if $i(h)=h$ for $h\in H$. (Recall the analogous definition for topological spaces.)  The homomorphism $p: G\rightarrow G/H$ with $p(g)=g+H$ is called a {it projection}. For nonabelian groups we modify the definition of a projection by having $H$ be a {\it normal} subgroup of $G$ and $p(g)=Hg=gH$. (In some mathematical word reuse, if $G=H\oplus K$, the maps of $(h, k)$ to either the first or the second coorcinate are also called projections.)
\end{definition}

\begin{definition}
An exact sequence of the form 

\begin{tikzpicture}
  \matrix (m) [matrix of math nodes,row sep=3em,column sep=4em,minimum width=2em]
  {
  0 & A_1 & A_2 & A_3 & 0\\};
  \path[-stealth]
    (m-1-1) edge  (m-1-2)
    (m-1-2)  edge node [above] {$\phi$} (m-1-3)
(m-1-3)  edge node [above] {$\psi$} (m-1-4)
(m-1-4) edge  (m-1-5);

\end{tikzpicture}

is called a {\it short exact sequence}\index{short exact sequence}.
\end{definition}

Short exact sequences are also very common. Try to prove the next two results with the facts we have learned.

\begin{theorem}
For the short exact sequence defined above, $\psi$ induces an isomorphism between $A_2/\phi(A_1)$ and $A_3$. Conversely, if $\psi: A\rightarrow B$ is an epimorphism with kernel $K$, then the sequence 

\begin{tikzpicture}
  \matrix (m) [matrix of math nodes,row sep=3em,column sep=4em,minimum width=2em]
  {
  0 & K & A & B & 0\\};
  \path[-stealth]
    (m-1-1) edge  (m-1-2)
    (m-1-2)  edge node [above] {$i$} (m-1-3)
(m-1-3)  edge node [above] {$\psi$} (m-1-4)
(m-1-4) edge  (m-1-5);

\end{tikzpicture}

is exact, where $i$ is inclusion.
\end{theorem}

\begin{theorem}
Suppose we have an exact sequence 

\begin{tikzpicture}
  \matrix (m) [matrix of math nodes,row sep=3em,column sep=4em,minimum width=2em]
  {
A_1 & A_2 & A_3 & A_4\\};
  \path[-stealth]
    (m-1-1)  edge node [above] {$\alpha$}  (m-1-2)
    (m-1-2)  edge node [above] {$\phi$} (m-1-3)
(m-1-3)  edge node [above] {$\beta$} (m-1-4);

\end{tikzpicture}

then the following are equivalent:
\begin{enumerate}
\item $\alpha$ is an epimorphism.
\item $\beta$ is a monomorphism.
\item $\phi$ is the zero homomorphism.
\end{enumerate}
\end{theorem}

We can get even more information from a short exact sequence if it is {\it split}.

\begin{definition}
Consider a short exact sequence

\begin{tikzpicture}
  \matrix (m) [matrix of math nodes,row sep=3em,column sep=4em,minimum width=2em]
  {
  0 & A_1 & A_2 & A_3 & 0.\\};
  \path[-stealth]
    (m-1-1) edge  (m-1-2)
    (m-1-2)  edge node [above] {$\phi$} (m-1-3)
(m-1-3)  edge node [above] {$\psi$} (m-1-4)
(m-1-4) edge  (m-1-5);

\end{tikzpicture}

The sequence is {\it split}\index{split exact sequence} if the group $\phi(A_1)$ is a direct summand of $A_2$.
\end{definition}

In this case the sequence is of the form 

\begin{tikzpicture}
  \matrix (m) [matrix of math nodes,row sep=3em,column sep=4em,minimum width=2em]
  {
  0 & A_1 & \phi(A_1)\bigoplus B & A_3 & 0.\\};
  \path[-stealth]
    (m-1-1) edge  (m-1-2)
    (m-1-2)  edge node [above] {$\phi$} (m-1-3)
(m-1-3)  edge node [above] {$\psi$} (m-1-4)
(m-1-4) edge  (m-1-5);

\end{tikzpicture}

where $\phi$ is an isomorphism of $A_1$ and $\phi(A_1)$, and $\psi$ is an isomorphism of $B$ and $A_3$.

Here are two more facts about split exact sequences.

\begin{theorem}
Suppose we have an exact sequence 

\begin{tikzpicture}
  \matrix (m) [matrix of math nodes,row sep=3em,column sep=4em,minimum width=2em]
  {
  0 & A_1 & A_2 & A_3 & 0.\\};
  \path[-stealth]
    (m-1-1) edge  (m-1-2)
    (m-1-2)  edge node [above] {$\phi$} (m-1-3)
(m-1-3)  edge node [above] {$\psi$} (m-1-4)
(m-1-4) edge  (m-1-5);

\end{tikzpicture}

then the following are equivalent:
\begin{enumerate}
\item The sequence splits.
\item There is a map $p: A_2\rightarrow A_1$ such that $p\circ\phi=i_{A_1}$, where $i_A$ is the identity map on an abelian group $A$.
\item There is a map $j: A_3\rightarrow A_2$ such that $\psi\circ j=i_{A_3}$
\end{enumerate}
\end{theorem}

The picture for the last two statements is 

\begin{tikzpicture}
  \matrix (m) [matrix of math nodes,row sep=3em,column sep=4em,minimum width=2em]
  {
  0 & A_1 & A_2 & A_3 & 0.\\};
  \path[-stealth]
    (m-1-1) edge  (m-1-2)
    ([yshift = 2pt]m-1-2.east)  edge node [above] {$\phi$} ([yshift = 2pt]m-1-3.west)
 ([yshift = -2pt]m-1-3.west)  edge node [below] {$p$} ([yshift = -2pt]m-1-2.east)
([yshift = 2pt]m-1-3.east)  edge node [above] {$\psi$} ([yshift = 2pt]m-1-4.west)
([yshift = -2pt]m-1-4.west)  edge node [below] {$j$} ([yshift = -2pt]m-1-3.east)
(m-1-4) edge  (m-1-5);

\end{tikzpicture}

{\bf Proof:}

To show that (1) implies (2) and (3) consider the sequence

\begin{tikzpicture}
  \matrix (m) [matrix of math nodes,row sep=3em,column sep=4em,minimum width=2em]
  {
  0 & A_1 & A_1\bigoplus A_3 & A_3 & 0.\\};
  \path[-stealth]
    (m-1-1) edge  (m-1-2)
    (m-1-2)  edge node [above] {$i$} (m-1-3)
(m-1-3)  edge node [above] {$\pi$} (m-1-4)
(m-1-4) edge  (m-1-5);

\end{tikzpicture}

Then just let $p: A_1\oplus A_3\rightarrow A_3$ be projection and $j: A_3\rightarrow  A_1\oplus A_3$ be inclusion. Then you can easily check that (2) and (3) hold.

To show that (2) implies (1). We will show that $A_2=\phi(A_1)\oplus(\ker p)$. For $x\in A_2$, we can write $x=\phi(p(x))+(x-\phi(p(x))).$ The first term is in $\phi(A_1)$ and the second is in $\ker p$, since $p(x-\phi(p(x)))=p(x)-p(\phi(p(x))=p(x)-p(x)=0.$ Secondly, if $x\in\phi(A_1)\cap(\ker p)$ then $x=\phi(y)$ for some $y\in A_1$, so $p(x)=p(\phi(y))=y$. But since $x\in (\ker p)$, $p(x)=0$, so $y=0$, which implies that $x=\phi(y)=0$. Thus $A_2=\phi(A_1)\oplus(\ker p)$.

To show that (3) implies (1). We will show that $A_2=(\ker\psi)\oplus j(A_3)$. Since $(\ker\psi)=\phi(A_1)$, we will be done. First, for $x\in A_2$, we write $x=(x-j(\psi(x)))+j(\psi(x))$.The first term is in $(\ker\psi)$ since $\psi(x-j(\psi(x)))=\psi(x)-\psi(j(\psi(x)))=\psi(x)-\psi(x)=0.$ The second term is in $j(A_3).$ Secondly, If $x\in(\ker\psi)\cap j(A_3)$, then $x=j(z)$ for some $z\in A_3$, so $\psi(x)=\psi(j(z))=z$. Since $x\in(\ker \psi)$, $\psi(x)=0$, so $z=0$, and $x=j(z)=0.$ Thus, $A_2=(\ker\psi)\oplus j(A_3)$. $\blacksquare$

\begin{corollary}
Suppose we have an exact sequence 

\begin{tikzpicture}
  \matrix (m) [matrix of math nodes,row sep=3em,column sep=4em,minimum width=2em]
  {
  0 & A_1 & A_2 & A_3 & 0.\\};
  \path[-stealth]
    (m-1-1) edge  (m-1-2)
    (m-1-2)  edge node [above] {$\phi$} (m-1-3)
(m-1-3)  edge node [above] {$\psi$} (m-1-4)
(m-1-4) edge  (m-1-5);

\end{tikzpicture}

Then if $A_3$ is free abelian, the sequence splits;
\end{corollary}

{\bf Proof:} Choose a basis for $A_3$ and for each basis element $e$, define $j(e)$ to be any element of the nonempty set $\psi^{-1}(e)$. The set is nonempty since $\psi$ is onto.

We will see more about exact sequences in Chapter 4, but this is a good test of your understanding of homomorphisms and their properties. The full power of exact sequences will wait until we discuss homology theory and the important Eilenberg-Steenrod axioms.

\section{Rings and Fields}

The purpose of this section and the next will not be to teach you an entire course on abstract algebra or linear algebra. Instead, the idea is to make sure that you know all of the terms that we will use later.

As I mentioned earlier, a group has a single binary operation, either addition or multiplication. A ring has both with each one playing a separate role. They are tied together with a distributive law. There are also some variants on the detinition I will give. First of all, there is an element that plays the role of the number 1 which is called a {\it unit element}\index{unit element}. I will never use a ring with out a unit element, although \cite{Her1} allows for one. Since the multiplicative identiy can be written as $I$ instead of 1, we can refer to a ring without a unit element as a {\it rng}\index{rng}. Also, most rings have the associative property for multiplication, but there are exceptions. One example is the {\it octonions}. They are a generalization of the {\it quaternions} which I will soon define, but multiplication is not associative. They play a role in the homotopy groups of spheres and obstruction theory by extension, but there are easier problems to solve before we can think about that. So unless otherwise specified, I will always use associative rings. Here is a definition.

\begin{definition}
A nonempty set $R$ is called an {\it (associative) ring}\index{ring}\index{associative ring} if it has two binary operations called addition and multiplication, it is an abelian group under addition, and it has the following additional properties:\begin{enumerate}
\item If $a, b\in R$, then $ab\in R$. (closure for multiplication)
\item There exists an element $1\in R$ called the unit element such that for $r\in R$, $1r=r1=r.$. (identity for multiplication)
\item For $a, b, c\in R$, $a(bc)=(ab)c.$ (associative law for multiplication)
\item For $a, b, c\in R$, $a(b+c)=ab+ac$ and $(b+c)a=ba+ca$ (distributive laws).
\end{enumerate}
\end{definition}

\begin{definition}
A {\it commutative ring}\index{commutative ring} is a ring for which multiplication is commutative, i.e. for $a,b\in R$, $ab=ba$. (Remember that for a ring, addition is always commutative.) A ring that is not a commutative ring is a {\it noncommutative  ring}\index{noncommutative ring}. 
\end{definition}

\begin{example} 
The integers with usual addition and multiplication form a commutative ring. So do the rationals and the reals. The irrationals are not a ring at all as they don't contain either 0 or 1.
\end{example}

\begin{example} 
The complex numbers consist of numbers of the form $a+bi$, where $a, b$ are real numbers and $i^2=-1$. If $a+bi$ and $c+di$ are two complex numbers, then we define $(a+bi)+(c+di)=(a+c)+(b+d)i$, and $(a+bi)(c+di)=(ac-bd)+(ad+bc)i.$ The complex numbers form a commutative ring with these definitions of addition and multiplication. The identies for addition and multiplication are $0=0+0i$ and $1=1+0i$ respectively. 
\end{example}

\begin{example} 
Our first example of a noncommutative ring is the ring of {\it quaternions}\index{quaternion}. The quaternions are numbers of the form $a+bi+cj+dk$ where $a,b, c, $ and $d$ are real numbers and $i^2=j^2=k^2=-1$. We also have $ij=k, jk=i, ki=j,  ji=-k, kj=-i,$ and $ik=-j$. So the quaternions are not commutative. One way to remember the correct sign is through ths diagram:

\begin{tikzpicture}
  \matrix (m) [matrix of math nodes,row sep=3em,column sep=4em,minimum width=2em]
  {
   i &  & j \\
     & k &  \\};
  \path[-stealth]
(m-1-1) edge[bend left] (m-1-3)
(m-1-3) edge [bend left] (m-2-2)
(m-2-2) edge [bend left] (m-1-1)
;
\end{tikzpicture}

Following the arrows gives a positive sign (eg. $ij=k$), while opposing them gives a negative sign (eg. $ji=-k$.) We have $$(a+bi+cj+dk)+(e+fi+gj+hk)=(a+e)+(b+f)i+(c+g)j+(d+h)k.$$ For multiplication, we have $$(a+bi+cj+dk)(e+fi+gj+hk)=(ae-bf-cg-dh)+(af+be+ch-dg)i+(ag+ce+df-bh)j+(ah+de+bg-cf)k.$$ If we reverse the order of the multiplication, then we switch the signs on the last two terms of the coefficients of $i$, $j$, and $k$.
\end{example}

\begin{example} 
The $n\times n$ matrices with real entries form a noncommutative ring ring under the usual addition and multiplication of square matrices.
\end{example}

We now look at some particular types of rings. 

\begin{definition}
If $R$ is a commutative ring. then a nonzero element $a\in R$ is a {\it zero-divisor}\index{zero-divisor} if there exists a $b\in R$ such that $b\neq 0$ and $ab=0$. A commutative ring is called an {\it integral domain}\index{integral domain} if it has no zero divisors.
\end{definition}

\begin{example} 
The integers form an integral domain.
\end{example}

\begin{example} 
Let $R=Z_p$ where $p$ is prime. If $a, b\in Z_p$ and $ab=0$ then $p$ must divide $ab$. Since $p$ is prime, it must divide $a$ or $b$. But both $a$ and $b$ are less than $p$, so either $a=0$ or $p=0$. So $Z_p$ is an integral domain.
\end{example}

\begin{example} 
Let $R=Z_m$ where $m$ is not a prime. Then $Z_m$ is not an integral domain as $m=ab$ for some $a<m$ and $b<m$. For example, if $m=15$ then $3\cdot 5=0$ in $Z_{15}$, so 3 and 5 are both zero-divisors.
\end{example}

The other thing we didn't assume about rings was an inverse for multiplication. 

\begin{definition}
A ring is a {\it division ring}\index{division ring} if all of its nonzero elements have multiplicative inverses. 
\end{definition}

\begin{definition}
A commutative division ring is called a {\it field}\index{field}. 
\end{definition}

\begin{example} 
The integers are a commutative ring but not a division ring and thus not a field. For example, there is no integer we can multiply by 2 to get 1. The rationals, reals, and complex numbers are all fields.
\end{example}

Here is an easy result for rings that we will need. I will prove the first statement. try to prove the rest.

\begin{theorem}
Let $R$ be a ring and $a, b\in R$.
\begin{enumerate}
\item a0=0a=0
\item a(-b)=(-a)b=-(ab)
\item (-a)(-b)=ab
\item (-1)a=(-1)
\item (-1)(-1)=1
\end{enumerate}
\end{theorem}

{\bf Proof of statement 1:} $a0=a(0+0)=a0+a0$. So $a0=0$ since an additive inverse is unique. Similarly, $0a=(0+0)a=0a+0a,$ so $0a=0$. $\blacksquare$

Now a division ring can not have zero divisors. To see this, suppose $D$ is a division ring, and $a, b$ are nonzero elements of $D$ with $ab=0$. Then $a$ has a multiplicative inverse $a^{-1}$ and $0=ab=a^{-1}(ab)=(a^{-1}a)b=1b=b,$ so $b=0$ which is a contradiction. So a division ring is an integral domain. This shows that $Z_m$ for $m$ not prime can not be a field. 

We would like to show that $Z_p$ for $p$ a prime is a field. To do this, we use the famous {\it Pigeonhole Principle}\index{Pigeonhole Principle}. 

The image you should have is a wall of mailboxes. These resemble a design of dovecote where the birds nest in an array of containers. (Doves and pigeons are really the same family of birds.) So suppose there are 100 mailboxes and the mail carrier has 101 letters, Then if she puts one letter in each box, there will still be one left over. So at least one lucky person got at least two letters. This obvious situation is stated below:

\begin{theorem}
{\bf Pigeonhole Principle:} If $n$ objects are distributed over $m$ places, and if $n>m$, then some place receives at least two objects.
\end{theorem}

\begin{example} 
Here is a more interesting example. Suppose that the average person has 100,000 hairs on their head. We know that Max, has 500,000 hairs on his head and that is more than anyone else in the U.S. Then there are two people in the U.S. with exactly the same number of hairs on their heads. Suppose we have boxes numbered 0-500,000. We also have index cards with the names of all 300 million Americans. For each person, we put their card in the box corresponding to the number of hairs on their head. Then we have more cards than boxes so at least two cards must go in the same box and we are done.
\end{example}

\begin{theorem}
A finite integral domain is a field.
\end{theorem}

{\bf Proof:} Let $D$ be a finite integral domain. (All our rings are assumed to have a unit element.) Let $r_1, r_2, \cdots, r_n$ be all of the elements in $D$. Suppose $a\in D$ is nonzero. Since integral domains are assumed to be commutative we are done if we can produce a multiplicative inverse for $a$. Consider the elements  $r_1a, r_2a, \cdots, r_na$. They are all distinct, since if $r_ia=r_ja$, then $0=r_ia-r_ja=(r_i-r_j)a$. Since $D$ is an integral domain and $a\neq 0$, we must have $r_i-r_j=0$, so $r_i=r_j$. So $r_1a, r_2a, \cdots, r_na$ must include all elements of $D$ by the Pigeonhole Principle as there are $n$ of these and $D$ has n elements. So one of these elements, say $r_ka$ must equal 1 and $r_k=a^{-1}$. $\blacksquare$

Now we know that $Z_p$ for $p$ a prime is a field by the previous result. It turns out that any finite field must have $p^n$ elements for $p$ a prime. See Chapter 7 in \cite{Her1} for details.

As in the case of groups, the special maps between rings are called {\it homomorphisms}. In this case, they are required to preserve both addition and multiplication.

\begin{definition}
Let $R$ and $S$ be rings. A map $\phi: R\rightarrow S$ is a {\it homomorphism} if 
\begin{enumerate}
\item $\phi(a+b)=\phi(a)+\phi(b)$
\item $\phi(ab)=\phi(a)\phi(b)$
\end{enumerate}
If a homomorphism of rings is one to one and onto, it is an {\it isomorphism}.
\end{definition}

Now we start to run into a situation that does not quite match groups. First of all, how do we define the kernel of a homomorphism? As a ring is an abelian group under addtion, we will let the kernel of $\phi$ be the set $\{r\in R|\phi(r)=0\}$. So the kernel can easily be shown to be a subgroup of the additive group of $R$, and we don't have to worry about whether it is a normal subgroup since the additive group of $R$ is asssumed to be abelian. But is the kernel a subring, i.e. is it a ring itself? Not quite. For one thing, it may not contain the unit element 1. Also, while $\phi(0)=0$, $\phi(1)$ may not equal the unit element of $S$.

\begin{example} 
Let $R=Z_{10}$ and $S=Z_{20}$. Let $\phi(r)=2r$, for $r\in Z_{10}$. Then $\phi(0)=0,$ but $\phi(1)=2$. Also the kernel of $\phi$ is $\{0\}$.
\end{example}

Here is a way the kernel of a homomorphism interacts with multiplication:

\begin{theorem}
Let $R$ and $S$ be rings. Let $\phi: R\rightarrow S$ be a {\it homomorphism} and let $K=\ker\phi$. Then if $u\in K$, then $ru$ and $ur$ are both in $K$ for all $r\in R$.
\end{theorem}

{\bf Proof:} Let $r$ be in $R$ and $u\in K$. Then $\phi(ru)=\phi(r)\phi(u)=\phi(r)0=0$, so $ru\in K$. $ur\in K$ by a similar argument. $\blacksquare$

So we would like to define a subset of a ring that would act in a way that is similar to the kernel of a ring homomorphism in order to define quotient rings which would be analogus to quotient groups. What would be the ideal subset to define? Hopefully in wishing for such an object, we are not being too idealistic. How about an {\it ideal}?

\begin{definition}
A nonempty subset $U$ of a ring $R$ is a (two sided) {\it ideal}\index{ideal} of $R$ if it is a subgroup of $R$ under addition and for $u\in U$, we have $ur\in U$ and $ru\in U$ for all $r\in R$. 
\end{definition}

\begin{example} 
The prime example you should keep in mind is letting $R=Z$ and letting $U$ be the even integers. Then any integer times an even number is even. Multiples of $m$ for $m\in Z$ also works. Odd numbers are not an ideal as an even number times an odd number is even. They are not even a subgroup as 0 is not odd. They are a coset $1+U$ of $U$ where $U$ is the even numbers. 
\end{example}

\begin{definition}
Let $R$ be a ring and $U$ an ideal of $R$. The quotient ring $R/U$ is the set of all cosets $r+U$ of $U$, where $r\in R$. Addition is defined by $(r_1+U)+(r_2+U)=(r_1+r_2)+U$, and multiplication is defined by $(r_1+U)(r_2+U)=(r_1r_2)+U$. 
\end{definition}

My only comment here will be about multiplication. Suppose $a_1+U=a_2+U$, i.e. they represent the same coset. Also, let  $b_1+U=b_2+U$. Then for multiplication to make sense, we need $a_1b_1+U=a_2b_2+U$. From our discussion of cosets of groups, we know that if $a_1+U=a_2+U$, then $a_1-a_2\in U$, so $a_1=a_2+u_1$ for some $u_1\in U$. Similarly, $b_1=b_2+u_2$ for some $u_2\in U$.So $a_1b_1=(a_2+u_1)(b_2+u_2)=a_2b_2+u_1b_2+a_2u_2+u_1u_2.$ Since $U$ is an ideal, the last 3 terms are all in $U$. Thus, $a_1b_1+U=a_2b_2+U$ and we are done. You should be able to convince yourself that the cosets form a ring with these operations.

\begin{theorem}
Let $R$ and $S$ be rings. Let $\phi: R\rightarrow S$ be onto and let $K=\ker\phi$. We now know that $K$ is an ideal in $R$, and $S\equiv R/K$.
\end{theorem}

As a final comment, in cohomology, we use what is known as a {\it graded ring}. We assign a {\it dimension} to every element and let $r_1r_2$ have dimension $p+q$ if $r_1$ has dimension $p$ and $r_2$ has dimension $q$. You already know about at least one graded ring. The polynomials in $x$ with real coefficients is of this form where we use the degree of the polynomial as its dimension. The polynomial rings are the prototype for cohomology rings. We will meet them in Chapter 8.

So that is more than enough ring theory to do algebraic topology. The main purpose of this section was to define terms and show the analogues with groups.

\section{Vector Spaces, Modules, and Algebras}

In this section, I will just review some basic concepts of vector spaces which are analogous to what we saw for groups and rings as well as some new properties. This should be a review for most readers. There are a lot of good books on linear algebra if you have gaps. I will take most of the material from \cite{Her1}. A classic old book that is very comprehensive is \cite{HoKu}. The material on modules may be new. It is also taken from \cite{Her1}. It turns out that a vector space is a special type of module but so is an abelian group. There is a structure theorem for finitely generated modules that closely parallels the one for abelian groups. Finally, the Steenrod Squares I will talk about in Chapter 11 form an algebra. Algebras have addition, multiplication, and scalar multiplication. They come in all sorts of strange shapes so I will content myself with defining them in general and talk about specific ones in their proper place.

\begin{definition}
A nonempty set $V$ is a {\it vector space}\index{vector space} over a field $F$ if $V$ is an abelian group under addition and if for every $a\in F$ and $v\in V$, there is an element $av\in V$ such that for all $a, b\in F$ and $v, w\in V$, :
\begin{enumerate}
\item $a(v+w)=av+aw.$
\item $(a+b)v=av+bv.$
\item $a(bv)=(ab)v$
\item $1v=v$, where 1 is the unit element of $F$.
\end{enumerate}
The elements of a vector space are called {\it vectors} and the elements of the asscoaiated field are called {\it scalars}. The product of a scalar and a vector is called {\it scalar multiplication}.
\end{definition}

\begin{example} 
The most common vector space we will see is $R^n$, the set of ordered $n$-tuples of real numbers. Here, $F=R$. For $v=(v_1,\cdots,v_n)$ and  $w=(w_1,\cdots,w_n)$, we have $v+w=(v_1+w_1,\cdots,v_n+w_n)$. For a real number $a$, $av=(av_1,\cdots,av_n)$.
\end{example}

For vector spaces, we have the ideas of subspaces, homomorphisms, and quotient spaces. We define those next.

\begin{definition}
A nonempty subset $W$ of $V$ is a {\it subspace}\index{subspace} of $V$ is $W$ is itself a vector space. Equivalently, $W$ is a subspace of $V$ is given $w_1, w_2\in W$ and $a, b\in F$, $aw_1+bw_2\in W$.
\end{definition}

As you probably can guess now, a homomorphism of vector spaces preserves addition and scalar multiplication. A vector space homomorphism has a special name that is much more common: a {\it linear transformation}\index{linear transformation}.

\begin{definition}
Let $V$ and $W$ be vector spaces. A function $T: V\rightarrow W$ is called a {\it linear transformation} if for $v_1, v_2\in V$ and $a\in F$, $T(v_1+v_2)=T(v_1)+T(v_2)$ and $T(av_1)=aT(v_1)$. 
\end{definition}

\begin{example} 
Let $f: R\rightarrow R$. Then if $f(x)=ax+b$, the graph of $f$ is a line. Is $f$ a linear transformation? We have $f(x_1+x_2)=a(x_1+x_2)+b=ax_1+ax_2+b$, while $f(x_1)+f(x_2)=ax_1+b+ax_2+b=ax_1+ax_2+2b$. So $f(x_1)+f(x_2)=f(x_1+x_2)$ if and only if $b=0.$ Also $f(cx)=acx+b$ for $c\in R$. while $cf(x)=acx+cb$. Since $c$ is arbitrary, we must again have $b=0$. So a line is the graph of a linear transformation if and only if it passes through the origin. Otherwise, it is the graph of an {\it affine transformation}.
\end{example}

As practice, try to prove this next theorem:

\begin{theorem}
If $V$ is a vector space over $F$. Let $0_V$ be the zero element of $V$ and $0_F$ be the zero element of $F$. Then 
\begin{enumerate}
\item $a0_V=0_V$, for $a\in F$.
\item $0_Fv=0_V$, for $v\in V$.
\item $(-a)v=-(av)$ for $a\in F$, $v\in V$. 
\item If $v\neq 0_V$, then $av=0_V$ implies $a=0_F$.
\end{enumerate}
\end{theorem}

Now we would like to construct quotient spaces. Let $W$ be s subspace of $V$. Remenber that $V$ and $W$ are abelian groups so we can certainly form a quotient group $V/W$ consisting of cosets $v+W$ for $v\in V$. We would like $V/W$ to be a vector space as well. The obvious scalar multiplication is for $a\in F, a(v+W)=av+W$. This is fine as long as it is true that if $v$ and $v'$ genertate the same coset (i.e. $v+W=v'+W$), $av+W=av'+W$. This is where we need the previous theorem. If $v+W=v'+W$, $(v-v')\in W$, and since W is a subspace, $a(v-v')\in W$. Now using part 3 of the previous theorem, $-v'=-(1v')=(-1)v'$, so $a(v-v')=av+a(-1)v'=av-av'\in W$. So $av+W=av'+W$ and scalar multiplication is well defined. So we have constructed a quotient space $V/W$. 

I will just mention that vector spaces have direct sums with the meaning you would expect. $V=U\oplus W$ if every element $v\in V$ can be written uniquely in the form $v==u+w$, where $u\in U$, and $w\in W$. 

We will finish off vector spaces with the notion of its {\it dimension}. Then you will know how to answer the next person who asks, "Isn't the fourth dimension time?". To get there we need a couple more definitions.

\begin{definition}
Let $V$ be a vector space over $F$. If $v_1,\cdots, v_n\in V$, then any element of the form $a_1v_1+\cdots +a_nv_n$, where the $a_i\in F$ is called a {\it linear combination}\index{linear combination} of  $v_1,\cdots, v_n$ over $F$.
\end{definition}

\begin{definition}
If $S$ is a nonempty subset of vector space $V$ over $F$, the {\it linear span}\index{linear span} $L(S)$ of $S$ is the set of all linear combinations of finite subsets of $S$. 
\end{definition}

\begin{definition}
If $V$ is a vector space over $F$, and $v_1,\cdots, v_n\in V$, we say that they are {\it linearly dependent}\index{linear independence} if for some $a_1,\cdots, a_n\in F$ which are not all zero, $a_1v_1+\cdots+a_nv_n=0$. If they are not linearly dependent, then they are linearly independent. 
\end{definition}

\begin{definition}
If $V$ is a vector space over $F$, a linearly independent set $v_1,v_2, \cdots\in V$  is a {\it basis} of $V$ if $V=L(S)$. The number of elements in the basis (which can be finite or infinite) is called the {\it dimension} of $V$. For $V$ finite dimensional it can be shown that any 2 bases have the same number of elements, so the definition makes sense.
\end{definition}

\begin{example} 
$R^n$ is an $n$-dimensional vector space. For example, let $n=2$. Then $(1,0)$ and $(0,1)$ form a basis. Any element $(a,b)\in R^2$ can be written as $a(1,0)+b(0,1)$, and they are linearly independent as  $a(1,0)+b(0,1)=(0,0)$ if and only if $a=b=0$. Another basis is $(2,1)$ and $(3,3)$. $a(2,1)+b(3,3)=(2a+3b, a+3b)$. If this is $(0,0)$, Then $2a+3b=0$, and $a+3b=0$. Solving for $a$ and $b$, we see that again $a=b=0$, so they are linearly independent. If $(c,d)$ is an arbitrary element of $R^2$, we have $c=2a+3b$, and $d=a+3b$. So $a=c-d$, and $d=a+3b=c-d+3b$, so $b=(2d-c)/3$. This shows that $(2,1)$ and $(3,3)$ span $R^2$, and they form a basis since they are linearly independent as well. 
\end{example}

Recall that in Chapter 2, we said that $R^n$ was a topological space as well. The usual topology on $R^n$ is that of a metric space. So $R^n$ is actually a topological vector space. Finite dimensional topological vector spaces are sort of boring, but there is an entire subject called {\it functional analysis}  that deals with infinite dimensional topological vector spaces. If you are curious, a great reference is {\em Functional Analysis} by Rudin \cite{Rud}.

We now move on to modules, a generalization of abelian groups. We won't have a lot of need for them,, but they come up frequently in the field of {\it homological algebra}. As this is closely related to algebraic topology, I will include a short description for completeness. A module is like a vector space but the scalars can be any kind of ring, not just a field. So any vector space is a module. That's why if you see a vector space on the moon it is a lunar module. 

\begin{definition}
Let $R$ be a ring. A nonempty set $M$ is an $R${\it -module}\index{module} if $M$ is an abelian group under addition and for every $r\in R$ and $m\in M$, there is an element $rm$ in $M$ such that for $a,b\in M$ and $r, s\in R$, 
\begin{enumerate}
\item $r(a+b)=ra+rb$
\item $r(sa)=(rs)a$
\item $(r+s)a=ra+sa$
\item If 1 is the unit element of $R$, then $1m=m$.
\end{enumerate}
\end{definition}

Now we need to be a little careful. If $R$ is a non-commutative ring, then what we have defined is a left module and we could define an analogous right module. But we won't have any reason to worry about that in this book.

\begin{example} 
Every abelian group is a module over the integers. If $g\in G$, and $k\in Z$, then $kg=g+\cdots+g$ ($k$ times). You can check that is satisfies all of the properties of a module.
\end{example}

\begin{example} 
Let $R$ be a commutative ring and $M$ be an ideal of $R$. Then we know that $rm$ and $mr$ are both in $M$ and are equal to each other.  So $M$ is a module over $R$.
\end{example}

\begin{example} 
$R$ is a module over itself. 
\end{example}

We can define direct sums of modules in the usual way. We will develop a structure theorem for modules that is analogous to the Fundamental Theorem of Finitely Generated Abelian Groups. We will need to make some restrictions on the ring of scalars first. For the rest of this section, $R$ will be a commutative ring with unit element 1. We will want $R$ to be a {\it Euclidean ring}\index{Euclidean ring} in which we have a process like long division.

 \begin{definition}
An integral domain $R$ is a {\it Euclidean ring} if for every $a\neq 0$ in $R$, there is a non-negative integer $d(a)$ called the {\it degree} of $a$ such that 
\begin{enumerate}
\item For all nonzero $a, b\in R, d(a)\leq d(ab)$.
\item For all nonzero $a, b\in R,$ there exists $t, r\in R$ such that $a=tb+r$, where either $r=0$ or $d(r)<d(b)$. 
\end{enumerate}
\end{definition}

We do not assign a value to $d(0)$. 

\begin{example} 
The integers form a Euclidean ring where for $a\neq 0, d(a)=|a|$. (The absolute value of $a$.) Then $|ab|\geq |a|$, where equality holds if $b=1$ or $b=-1$. For the second part, suppose for example $a=30,$ and $b=7$. Then if we divide 30/7, we get a quotient of 4 and a remainder of 2. So let $t=4$ and $r=2$. Then $30=4(7)+2$ with $|2|<|7|.$ So we are basically doing long division.
\end{example}

\begin{example} 
The other interesting example is to let $R$ be the ring of polynomials in $x$ with coefficients in a field. (Say with real coefficients.) Then $d(f)$ is the degree of $f$. You can check that this is also a Euclidean ring and you can perform long division with polynomials as well. 
\end{example}

Ideals of Euclidean rings have a special property. They are all {\it principal ideals}, which means that they are of the form $(a)=\{xa|x\in R\}$. In other words, they are generated by a single element. 

\begin{definition}
An integral domain $R$ is a {\it principal ideal domain}\index{principal ideal domain} if every ideal of $R$ is a principal ideal.
\end{definition}

\begin{theorem}
Every Euclidean ring is a principal ideal domain. 
\end{theorem}

{\bf Proof:}
Let $R$ be a Euclidean ring and $A$ an ideal of $R$. Assume that $A\neq 0$ or we could generate $A$ with the single element $0$. Let $a_0\in A$ be of minimal degree in $A$, i.e. there is no nonzero element with a strictly smaller degree. This is always possible as the degree is a non-negative integer. Let $a\in A$. Then there exists $t, r\in R$, such that $a=ta_0+r$ where $r=0$ or $d(r)<d(a_0)$. Since $A$ is an ideal of $R$, $ta_0\in A$, so $r=a-ta_0\in A$. If $r\neq 0$, then $d(r)<d(a_0)$. But this is impossible since we assumed $a_0$ was of minimal degree in $A$. So $r=0$ and $a=ta_0$. Since $a$ was an arbitrary element of $A$, $A=(a_0)$. So $A$ is a principal ideal domain. $\blacksquare$

\begin{definition}
An module over $R$ is {\it cyclic} if there is an element $m_0\in M$ such that every $m\in M$ is of the form $m=rm_0$ for some $r\in R$. $M$ is {\it finitely generated} if there is a finite set of elements $m_1, \cdots, m_n$ such that every $m\in M$ is of the form $r_1m_1+\cdots+r_nm_n$ for some $r_1, \cdots, r_m\in R$. 
\end{definition}

Now we can state the promised structure theorem:

\begin{theorem}
Let $M$ be a finitely generated module over a principal ideal domain $R$. Then there is a unique decreasing sequence of principal ideals $(d_1)\supseteq (d_2) \supseteq \cdots\supseteq (d_n)$ such that $$M\equiv R\bigoplus \cdots\bigoplus R \bigoplus R/(d_1)\bigoplus R/(d_2) \bigoplus \cdots \bigoplus R/(d_n).$$ The generators of the ideals are unique and for each $i$ for $1\leq i\leq n-1$, $d_i$ divides $d_{i+1}$.
\end{theorem}

Finally, we have an {\it algebra} which combines addition, multiplication, and scalar multiplication.

\begin{definition}
A ring $A$ is called an {\it algebra}\index{algebra (object)} over a field $F$ if $A$ is a vector space over $F$ such that for all $a, b\in A$, there is an element $ab\in A$ with the property that if $\alpha\in F$, $\alpha(ab)=(\alpha a)b=a(\alpha b)$.
\end{definition}

\begin{example} 
One example of an algebra is the square $n\times n$ matrices with real entries. Then we have the usual matrix addition and multiplication. Scalars are the real numbers. If $m$ is a matrix and $a$ is a real number, then $am$ is the matrix $m$ with all of its entries multiplied by $a$.
\end{example}

In algebraic topology, the Steenrod squares form an algebra. We will see then in chapter 11.

So you now should have seen a lot of patterns that repeat for groups, rings, vector spaces, etc. They each have special maps preserving their structure as well as subobjects, product objects, and quotient objects. Category theory will allow us to talk about these properties in general terms.

\section{Category Theory}

I know what you are thinking now. ''This chapter is not abstract enough for me. I need even more abstraction." Never fear, in this last section, we will talk about category theory. I will take the material from the second volume of Jacobson \cite{Jac}. For a very complete description, see MacLane \cite{MacL1}. A more recent and currently popular book is Riehl \cite{Rie}.

Category theory addresses some issues we have already seen. The first is Russell's paradox. Recall we have seen that the collection of all sets is not a set. Instead of the set of sets, we will talk about the {\it category} of sets. Also, I have shown you that there a lot of repeated patterns. Topological spaces, groups, rings, and vector spaces, all have special subsets, products, quotients, and special maps between them that preserve their main structure. Sometimes, it is convenient to just prove results about these objects in general using the properties they have in common.  A category will contain a type of set or {\it object} and a special map between objects called a {\it morphism}.  Finally, we will want nicely behaved maps between categories. These are called {\it functors}. As an example, homology theory takes a topological space and maps it to a collection of abelian groups: one for each dimension. This collection is also known as a {\it graded group.}

What do you call the students in an abstract algebra class a week after the final? A graded group. 

The amount of category theory used in algebraic topology has varied over time. There are some modern TDA papers that are filled with it. My personal preference is to go easy on it, but you can't do algebraic topology without at least some category theory. In this section, I will provide basic definitions and examples and we will come back to it later as needed.

\begin{definition}
A {\it category}\index{category} {\bf C} consists of a class {\it ob }{\bf C} of {\it objects}, and a set ${\it hom}_{\bf C}(A,B)$ of {\it morphisms}\index{morphism} for each ordered pair $(A,B)$ of objects in {\bf C}. (We drop the subscript {\bf C} if the context is obvious and just write ${\it hom}(A,B)$.) For each ordered triple of objects $(A, B, C)$, there is a map $(f, g)\rightarrow gf$ called {\it composition} taking ${\it hom}(A,B) \times {\it hom}(B,C)$ to ${\it hom}(A,C)$. The objects and morphisms satisfy the following conditions:
\begin{enumerate}
\item If $(A, B)\neq(C, D)$, then ${\it hom}(A, B)$ and ${\it hom}(C, D)$ are disjoint.
\item (Associativity) If $f\in{\it hom}(A, B)$, $g\in{\it hom}(B, C)$, and $h\in{\it hom}(C, D)$, then $(hg)f=h(gf)$.
\item (Unit) For every object $A$, we have a unique element $1_A\in{\it hom}(A, A)$, such that $f1_A=f$ for every $f\in {\it hom}(A, B)$, and $1_Ag=g$, for every $g\in {\it hom}(B, A)$.
\end{enumerate}
\end{definition}

If the objects in a category form an actual set, the category is called a {\it small category}.

I will now give some examples. Convince yourself that they satisfy the definition.

\begin{example} 
Our first example is the category of sets. The objects are sets and the morphisms are functions. 
\end{example}

\begin{example} 
The category of groups. The morphisms are group homomorphisms. 
\end{example}

\begin{example} 
The category of vector spaces over a given field $F$. The morphisms are linear transformations.
\end{example}

\begin{example} 
The category of topological spaces. The morphisms are continuous functions.
\end{example}

\begin{definition}
A category {\bf D} is a subcategory of {\bf C} if any object of {\bf D} is an object of {\bf C}, and for any pair $A, B$ of objects of {\bf D}, ${\it hom}_{\bf D}(A,B)\subset {\it hom}_{\bf C}(A,B)$.
\end{definition}

\begin{example} 
The category of abelian groups. Again. let the morphisms be group homomorphisms. Abelian groups are a subcategory of the category of groups. 
\end{example}

Categories have special morphisms called {\it isomorphisms}. If $f\in{\it hom}(A, B)$, then $f$ is an isomorphism if there exists $g\in{\it hom}(B , A)$ such that $fg=1_B$, and $gf=1_A$. If such a map $g$ exists, it is unique and we write $g=f^{-1}$. In that case, $(f^{-1})^{-1}=f$.  If $f$ and $h$ are isomorphisms, and $fh$ is defined, then $(fh)^{-1}=h^{-1}f^{-1}$. In the category of sets, isomorphisms are bijective sets, while in the category of groups, they are group isomorphisms.

Here are two ways to get new categories from old ones. 

\begin{definition}
Let {\bf C} be a category. We define the {\it dual category}\index{dual category} or {\it opposite category}\index{opposite category} ${\bf C}^{op}$ to have the same objects as {\bf C}, but  ${\it hom}_{\bf C^{op}}(A,B)= {\it hom}_{\bf C}(B,A)$. Commutative diagrams in the dual category correspond to diagrams in the original category, but the arrows are reversed. If $f: A\rightarrow B$ in {\bf C}, then $f: B\rightarrow A$ in ${\bf C}^{op}$ , and if 

\begin{tikzpicture}
  \matrix (m) [matrix of math nodes,row sep=3em,column sep=4em,minimum width=2em]
  {
   A & B\\
  & C\\
};
  \path[-stealth]
    (m-1-1) edge node [above] {$f$} (m-1-2)
    (m-1-1)  edge node [below] {$h$} (m-2-2)
    (m-1-2) edge node [right] {$g$} (m-2-2);
\end{tikzpicture}

is commutative in {\bf C}, then 

\begin{tikzpicture}
  \matrix (m) [matrix of math nodes,row sep=3em,column sep=4em,minimum width=2em]
  {
   A & B\\
  & C\\
};
  \path[-stealth]
    (m-1-2) edge node [above] {$f$} (m-1-1)
    (m-2-2)  edge node [below] {$h$} (m-1-1)
    (m-2-2) edge node [right] {$g$} (m-1-2);
\end{tikzpicture}

is commutative in ${\bf C}^{op}$ .

\end{definition}

In general, the term ''dual" or the prefix "co-" will mean to reverse arrows.

We can also talk about the {\it product} of two categories.

\begin{definition}
If {\bf C} and {\bf D} are two categories,then we can define the product category {\bf C} $\times $ {\bf D}. The objects are ordered pairs $(A, A')$, where $A$ is an object in {\bf C} and $A'$ is an object in {\bf D}. Morphisms are also ordered pairs, and if $A, B, C$ are objects of {\bf C},  $A', B', C'$ are objects of {\bf D}, $f\in{\it hom}_C(A, B)$, $g\in{\it hom}_C(B, C)$, $f'\in{\it hom}_D(A', B')$, $g'\in{\it hom}_D(B', C')$, then $$(g, g')(f, f')=(gf, g'f').$$
\end{definition}

The next concept is very important in algebraic topology. A {\it functor}\index{functor} is a map between categories. It takes objects to objects and morphisms to morphisms. There are two types of functors. Covariant functors keep morphisms going int the same direction and contravariant functors reverse the direction. we will meet them both in algebraic topology.

\begin{definition}
If {\bf C} and {\bf D} are categories, a {\it covariant functor}\index{covariant functor} $F$ from {\bf C} to {\bf D} consists of a map taking an object $A$ in {\bf C} to an  object $F(A)$ in {\bf D}, and for every pair of objects $A, B$ in {\bf C}, a map of ${\it hom}_C(A, B)$ into ${\it hom}_D(F(A), F(B))$,such that $F(1_A)=1_{F(A)}$, and if $gf$ is defined in {\bf C}, then $F(gf)=F(g)F(f)$. 
\end{definition}

\begin{definition}
If {\bf C} and {\bf D} are categories, a {\it contravariant functor}\index{contravariant functor} $F$ from {\bf C} to {\bf D} consists of a map taking an object $A$ in {\bf C} to an  object $F(A)$ in {\bf D}, and for every pair of objects $A, B$ in {\bf C}, a map of ${\it hom}_C(A, B)$ into ${\it hom}_D(F(B), F(A))$,such that $F(1_A)=1_{F(A)}$, and if $gf$ is defined in {\bf C}, then $F(gf)=F(f)F(g)$. 
\end{definition}

A contravariant functor can also be thought of as a covariant functor from ${\bf C}^{op}$ to {\bf D}.

I will now give some examples of functors. They are all covariant unless otherwise specified.

\begin{example} 
Let {\bf D} be a subcategory of {\bf C}. The {\it injection} functor takes an object in D to itself and maps a morphism in {\bf D} to the same morphism in {\bf C}. In the special case where {\bf D}={\bf C}, it is called the {\it identity} functor.
\end{example}

\begin{example} 
A {\it forgetful functor}\index{forgetful functor} ''forgets" some of the structure of a category. For example, let $F$ go from the category of groups to the category of sets where it takes a group to its underlying set. Another example takes rings to abelian groups and a ring homomorphism to just its additive property to become an additive abelian group homomorphism.
\end{example}

\begin{example} 
The {\it power functor} goes from the category of sets to itself and takes any set to its power set (i.e. the set of its subsets. See  Section 2.6.) If $f: A\rightarrow B$, this functor sends $f$ to $f_P: P(A)\rightarrow P(B)$ which sends a subset $A_1\subset A$ to the subset $f(A_1)\subset B$. The empty set is sent to itself.
\end{example}

\begin{example} 
The {\it projection functor} goes from the product category {\bf C} $\times $ {\bf D} to {\bf C}. It takes $A\times A'$ to $A$, and $(f,g)\in {\it hom}((A, B), (A', B'))$ into $f\in{\it hom}(A, B)$.
\end{example}

\begin{example} 
Homology, which we will define in Chapter 4 is a functor from the category of topological spaces to the category of graded groups. A continuous function corresponds to a group homomorphism in each direction.
\end{example}

\begin{example} 
Finally, I will give an example of a contravariant functor which will be useful in cohomology.  (We will define cohomology in Chapter 8.) Let $G$ be a fixed abelian group. Define a functor which takes an abelian group $A$ to the group of homomorphisms $Hom(A, G)$ from $A$ to $G$. A map from $A$ to $B$ corresponds to a map from $Hom(B, G)$ to $Hom(A, G)$. To see this, let $f: A\rightarrow B$, and $g\in Hom(B,G)$. Then $gf\in Hom(A, G)$, so we have a map $Hom(B,G)\rightarrow Hom(A,G)$ given by $g\rightarrow gf$. Note that this is a contravariant functor. If we take $A$ to $Hom(G, A)$ for fixed $G$ then this is a covariant functor.
\end{example}

Finally, I will define a map between two functors called a {\it natural transformation}\index{natural transformation}.

\begin{definition}
Let $F$ and $G$ be functors from {\bf C} to {\bf D}. We define a {\it natural transformation} $\eta$ from $F$ to $G$ to be a map that assigns to every object $A$ in {\bf C} a morphism $\eta_A\in {\it hom}_{\bf D}(FA, GA)$ such that for any objects $A, B$ of {\bf C}, and any $f\in {\it hom}_{\bf C}(A,B)$, the rectangle below is commutative:

\begin{tikzpicture}
  \matrix (m) [matrix of math nodes,row sep=3em,column sep=4em,minimum width=2em]
  {
   A & FA & GA\\
  B & FB & GB\\
};
  \path[-stealth]
    (m-1-1) edge node [left] {$f$} (m-2-1)
    (m-1-2)  edge node [left] {$F(f)$} (m-2-2)
   (m-1-3)  edge node [right] {$G(f)$} (m-2-3)
   (m-1-2)  edge node [above] {$\eta_A$} (m-1-3)
    (m-2-2) edge node [below] {$\eta_B$} (m-2-3);
\end{tikzpicture}

If every $\eta_A$ is an isomorphism, then $\eta$ is called a natural isomorphism.

\end{definition}

I will save examples of natural transformations until the next chapter. 

This gives you a taste of category theory. We will revisit it when necessary. 

Now we have all the background we need. The next chapter will start to describe algebraic topology.

\chapter{Traditional Homology Theory}

Now I will put the last two chapters together and talk about algebraic topology. Before I get into persistent homology, which is the basis of topological data analysis as it is currently practiced, I will give an introduction to homology using the traditional approach. Persistent homology will be an easy special case. If we want to get beyond it and start exploring the possible applications of cohomology and homotopy, we will need a better understanding of the subject as it was originally developed. Rather than make this chapter 100 pages long, I will use it as an outline and you can look in any of the textbooks that I will list for more information. One thing I will sacrifice, though, is giving you a full appreciation for how hard it is to prove that all of the machinery works the way it is supposed to. Like a lot of other subjects in math, you go through the pain of developing the tools and then never look back.

The first section of this chapter will be the longest and deal with {\it simplicial homology}. A topological space is represented by  a collection of vertices, edges, triangles, etc.  There is an equivalent combinatorial approach where we are looking at subsets of a collection of objects. The latter approach is the one that is more applicable to data science. As we mentioned at the end of the last chapter, homology is a functor which takes a topological space to a graded group, i.e. an abelian group in every dimension up to the dimension of the space. This means that we not only need to produce these groups, but for a continuous map between two different spaces, there has to be a corresponding group homomorphism in every dimension between those spaces' homology groups. Also, if the spaces are homeomorphic, we would like these homomorphisms to be actual isomorphisms. This is the property of {\it topological invariance}. We can do all of this, but proving it works turns out to pretty difficult. I won't go through it all, but you can find it in any of the references.Along the way, I will briefly discuss computation. There is a very ugly way to compute these groups directly using the definition. Exact sequences give a shortcut in a lot of special cases. But a computer can always compute homology groups. I will briefly outline the method in \cite{Mun1}. There are much better methods developed in the last 30 years and there is software that is publicly available, so I won't go into a lot of detail on this particular method. I will finish the section with some interesting applications to maps of spheres and fixed points. The former is the key to obstruction theory, which I will describe in more detail later. 

The next section is about the Eilenberg Steenrod Axioms. In \cite{ES}, Eilenberg and Steenrod took an axxiomatic approach to the subject. In other words, there are no pictures. A big contribution was the idea of a homology theory that would satisfy a set of axioms. The full list is a big help in computation and has interesting parallels in cohomology and homotopy. Given the notion of a homology theory, you would expect there to be more than one. I will briefly describe two others in the last section. Simplicial homology only works on spaces that can be represented by simplices or {\it triangulated}. Singular homology is more general and proving topological invariance is much easier than in simplicail homology. But computations are not at all practical. A good compromise is cellular homology. This has the advantage that your topological space does not have to plugged into a wall outlet to use it. It replaces simplicial complexes with the more flexible {\it CW complexes}. For any situation we would meet in data science, you can use any of these theories. In the case where they all apply, the homology groups are all the same.

The material in this chapter is entirely taken from Munkres book, {\em Elements of Algebraic Topology} \cite{Mun1}, which is my personal favorite. Here are some other references. As far as I know, the earliest book on algebraic topology was from Alexandrov and  Hopf in 1935 \cite{AlHo}. (We will hear more about Hopf later.) It is in German. The first American (and English language) textbook was written in 1942 by Solomon Lefschetz \cite{Lef}. Eilenberg and Steenrod's book \cite{ES} that I metioned above was the standard reference for many years. Another popular book is Spanier \cite{Spa} which first came out in 1966. Currently, the main competitor to Munkres is Hatcher \cite{Hat}. Hatcher's book seems to be the most popular these days, and it had a free online version for some time before it was published in print. While it covers more material than Munkres such as Steenrod squares and some homotopy theory, I find it to be a much harder book. There are books on homotopy that I like better and will mention them when we get to chapter 9.

\section{Simplicial Homology}

I will now describe simplicial homology. In TDA (as it stands now), it is the only homology theory we will ever need. 

\subsection{Simplicial Complexes}

The basic building block of homology theory is a {\it simplex}. This is a point in dimension 0, a line segment in dimension 1, a triangle in dimension 2, a tetrahedoron in dimension 3, etc. Note that for an $n$-dimensional simplex, there are $n+1$ vertices. Simplices are glued together to form {\it simplicial complexes}. Most shapes you have ever seen or thought about can be built in this way. As an example, Epcot's Spaceship Earth at Disney World is  is a copy of the sphere $S^2$ modeled as a simplicial complex. (See Figure 4.1.1) Here are the precise definitions.

\begin{figure}[ht]
\begin{center}
  \scalebox{0.4}{\includegraphics{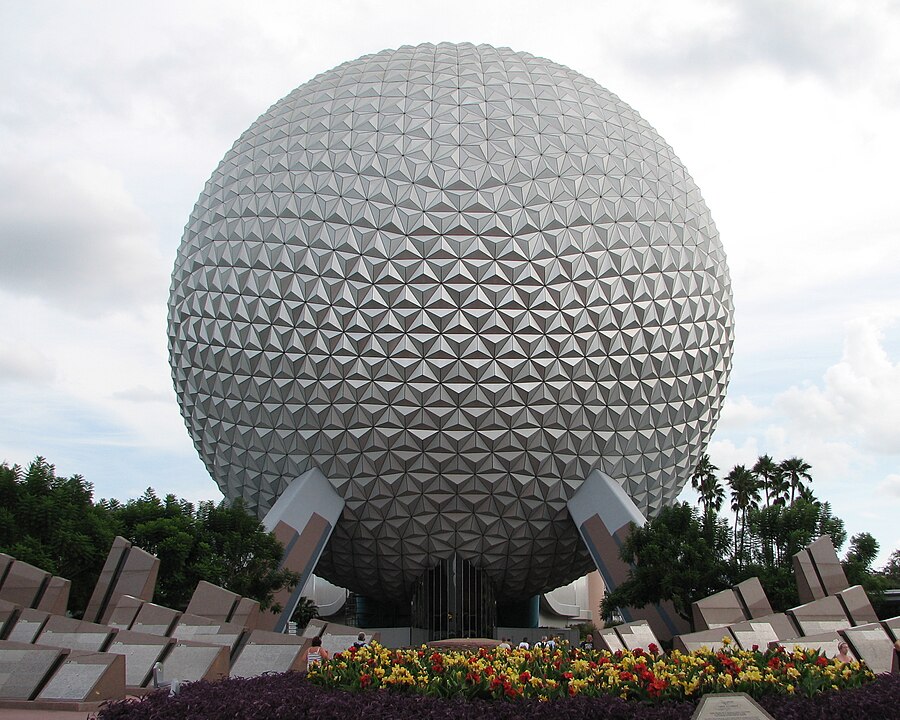}}
\caption{
\rm 
Spaceship Earth \cite{WikSE}
}
\end{center}
\end{figure}

\begin{definition}
A subset $X$ of $R^n$ is {\it convex}\index{convex set} if given any two points $x_1, x_2\in X$, the line segment connecting them is contained in $X$. (See Figure 4.1.2.) This line is of the form $x=tx_1+(1-t)x_2$ where $0\leq t\leq 1$. If $x_0, x_1,\cdots, x_k\in R^n$, the {\it convex hull}\index{convex hull} of these points is the smallest convex subset of $R^n$ containing them.
\end{definition}

\begin{figure}[ht]
\begin{center}
  \scalebox{0.4}{\includegraphics{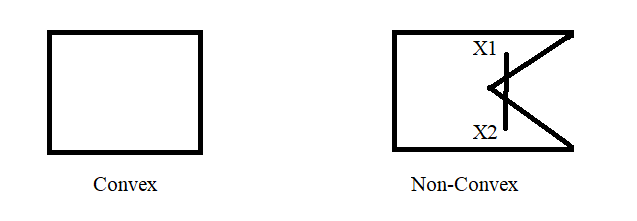}}
\caption{
\rm 
Convex vs Non-Convex Set
}
\end{center}
\end{figure}

\begin{definition}
A set of points $x_0, x_1,\cdots, x_k\in R^n$ is {\it geometrically independent}\index{geometrically independent} if for any real scalars $t_i$, $\sum_{i=0}^k t_i=0$ and $\sum_{i=0}^k t_ix_i=0$ imply $t_i=0$ for all $i$. This is the same as saying that any three of these points are not on the same line, any four are not on the same plane, etc. 
\end{definition}

\begin{definition}
Let  $x_0, x_1,\cdots, x_k\in R^n$ be geometrically independent. The $k${\it-simplex}\index{simplex} $\sigma$ spanned by  $x_0, x_1,\cdots, x_k$ is the convex hull of these points. In formulas, it is the set $$\{x|x= \sum_{i=0}^k t_ix_i\} \hspace{.1 in}{\rm where }\sum_{i=0}^k t_i=1  \hspace{.1 in}{\rm and  }  \hspace{.1 in}t_i\geq 0  \hspace{.1 in}{\rm   for \hspace{.1 in} all  } \hspace{.1 in} i.$$ The points,  $x_0, x_1,\cdots, x_k$ are the {\it vertices} of $\sigma$, and we say that $\sigma$ has dimension $k$. If $\tau$ is a simplex spanned by a subset of these vertices, then $\tau$ is a {\it face} of $\sigma$. The real numbers $\{t_i\}$ are called the {\it barycentric coordinates} of $\sigma$.
\end{definition}

Note that a $k$-simplex has $k+1$ vertices. For example, a triangle is a 2-simplex with three vertices, $x_0, x_1$, and $x_2$. Letting the $n$-ball $B^n$ be the interior of the the $(n-1)$-sphere $S^{n-1}$, an $n$-simplex is homeomorphic to $B^n$ and its boundary is homeomorphic to $S^{n-1}$. So a triangle (remember it is made of silly putty) is homeomorphic to $B^2$. Its boundary can be deformed into the circle, $S^1$.

We now stick simplices together into a {\it simplicial complex}.

\begin{definition}
A {\it simplicial complex}\index{simplicial complex} $K$ in $R^n$ is a collection of simplices such that
\begin{enumerate}
\item Every face of a simplex of $K$ is in $K$. 
\item The intersection of any two simplices of $K$ is a face of each of them.
\end{enumerate}
\end{definition}

Figure 4.1.3  shows some examples of good and bad simplicial complexes. The main thing is that the intersection of two faces must be an entire face. 

\begin{figure}[ht]
\begin{center}
  \scalebox{0.4}{\includegraphics{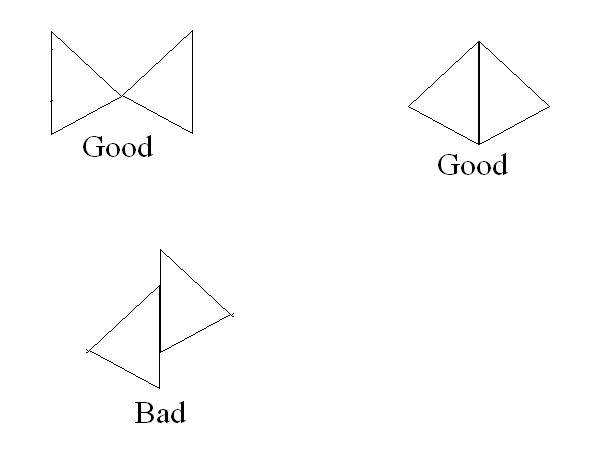}}
\caption{
\rm 
Good and Bad Simplical Complexes
}
\end{center}
\end{figure}

The next definition is very important.

\begin{definition}
Let $K$ be a simplicial complex and $K^p$ be the collection of all simplices of dimension at most $p$. Then $K^p$ is called the $p${\it -skeleton}\index{skeleton} of $K$. 
\end{definition}

To visualize where this term comes from, imagine Spaceship Earth with all its triangles missing. Then there would only be vertices and edges and you would see an outline of the sphere. That complex would be the 1-skeleton,  Spaceship Earth$^1$.

Now, we look at the space taken up by the complex $K$. 

\begin{definition}
Let $K$ be a simplicial complex. The union of all of the simplices of $K$ as a subset of $R^n$ is called the {\it polytope} or {\it underlying space} of$K$ and denoted by $|K|$.
\end{definition}

The polytope $|K|$ has some nice properties. First of all, any simplex is a closed bounded subset of $R^n$ so it is compact by the Heine-Borel Theorem. So for a finite complex $K$, $|K|$ is compact as it is a finite union of compact sets. It is easily proved that for a space $X$, $f: |K|\rightarrow X$ is continuous if and only if $f|\sigma$ is continuous for every simplex $\sigma$ in $K$. I will give the proof that $|K|$ is Hausdorff. First, we need to extend the definition of barycentric coordinates to a simplicial complex.

\begin{definition}
Let $K$ be a simplicial complex. If $x$ is a point of $|K|$, then $x$ is a vertex of $K$ or it is in the interior of exactly one simplex $\sigma$ of $K$. Suppose the vertices of $\sigma$ are $a_0, \cdots, a_k$. Then $$x=\sum_{i=0}^k t_ia_i,$$ where $t_i>0$ for each $i$ and $\sum t_i=1$. If $v$ is an arbitrary vertex of $K$ we define the barycentric coordinate $t_v(x)$ of $x$ as $t_v(x)=0$ if $v$ is not one of the vertices $a_i$ and $t_v(x)=t_i$, if $v=a_i$.
\end{definition}

\begin{theorem}
$|K|$ is Hausdorff.
\end{theorem}

{\bf Proof:} If $x\neq y$, then there is at least one vertex $v$ of $K$ with $t_v(x)\neq t_v(y)$. Without loss of generality, assume $t_v(x)<t_v(y)$, and choose a number $r$ with $t_v(x)<r<t_v(y)$. Then $\{z\in K|t_v(z)<r\}$ and $\{z\in K|t_v(z)>r\}$ are disjoint open sets containing $x$ and $y$ respectively. $\blacksquare$

Homology will be a functor that will take simplicial complexes to a group in each dimension. We will then need special maps between simplicial complexes. 

\begin{theorem}
Let $K$ and $L$ be complexes, and let $f: K^0\rightarrow L^0$ be a map. Suppose whenever the vertices $v_0, \cdots, v_n$ of $K$ span a simplex of $K$, $f(v_0), \cdots, f(v_n)$ span a simplex of $L$. Then $f$ can be extended to a continuous map $g: |K|\rightarrow |L|$ such that $$x=\sum_{i=0}^n t_iv_i$$ implies $$g(x)=\sum_{i=0}^n t_if(v_i).$$ We call $G$ the {\it simplicial map}\index{simplicial map} induced by $f$.
\end{theorem}

Note that we don't insist that the $f(v_i)$ are always distinct. If they are and $f$ is bijective, then $g$ turns out to be a homeomorphism. Also note that the composition of two simplical maps is a simplicial map.

We end this subsection with a different kind of simplicial complex. This version is much more useful in TDA applications.

\begin{definition}
An {\it abstract simplicial complex}\index{abstract simplicial complex} is a collection $\frak{S}$ of nonempty sets such that if $A$ is an element of $\frak{S}$, so is every nonempty subset of $A$. $A$ is called a {\it simplex} of $\frak{S}$. The {\it dimension} of $A$ is one less than the number of elements of $A$. A subset of $A$ is called a {\it face} of $A$.
\end{definition}

Note that in an abstract simplicial complex, any two simplices have to intersect in an edge. So the only property we need to worry about is that the complex contains all of the faces of any of its simplices. We can sometimes define an interesting abstract simplicial complex by defining {\it maximal simplices}, i.e. simplices that are not contained in any larger simplex. We then get the complex by including the maximal simplices and all of their faces.

\begin{example}
Here is an idea of mine. The Apriori Algorithm \cite{AIS} looks at frequent common subsets of a large collection of objects. For example, a supermarket might analyze the items that customers buy together. It determines rules like a customer that buys peanut butter and bread is also likely to buy jelly. For each customer, we can assign an abstract simplicial complex. The maximal simplices are the items bought in each shopping trip and then we make this into a simplicial complex by including all of the faces. My idea is to use the homology to cluster customers. For example, can you distinguish a customer who is single from one who is buying for a large family? I will discuss this more in Chapter 7.
\end{example}

\begin{example}
Abstract simplicial complexes are also closely related to {\it hypergraphs}. Suppose we have a set of objects called {\it vertices}. An (undirected) graph is the collection of vertices along with a set of pairs of vertices called {\it edges}. The two vertices are the {\it ends} of the edge and they can be the same if we allow loops. A hypergraph is a variant of a graph where the {\it hyperedges} can contain more than two vertices. For example, if the vertices are employees in an office, one person can send an email to five of the employees making a hyperedge of size 6. Taking all of the hyperedges in a hypergraph along with all of their subsets forms an abstract simplicial complex.
\end{example}

One more thing to note is that an abstract simplicial complex can always be represented as a simplicial complex as we first defined it. If there are $N+1$ vertices, we use the $N$ standard basis vectors along with the origin in $R^N$. (Recall that a standard basis vector has one coordinate a one and the rest zeros.) Label these points to match the vertices in your abstract complex, and for each simplex in the complex, form the simplex from its corresponding vertices in $R^N$. Thus we have a pairing between both kinds of simplicial complexes. 

In the next part, we will see how to turn complexes into groups.

\subsection{Homology Groups}

Now I will define homology groups. For now, I will talk about homology groups whose coefficients are integers. If you can understand these, going to other types of coefficients are not that hard. In persistent homology, integers are replaced with $Z_2$ coefficients. I will discuss the differences in the next subsection. There is a general formula for changing between coefficients but it is closely related to cohomology and will be discussed in Chapter 8 when we develop some additional tools. One issue that moving to $Z_2$ eliminates is the issue of {\it orientation} which I will describe next. 

\begin{definition}
Let $\sigma$ be a simplex with vertices $v_0, v_1, \cdots, v_p$. We choose an ordering of its vertices and say that two orderings are equivalent if they differ by an even permutation, i.e. an even number of swaps of pairs of them. The orderings then divide into two equivalence classes called {\it orientations}\index{orientation}. We choose a particular orientation and call $\sigma$ an {\it oriented simplex}, written $\sigma=[v_0, v_1, \cdots, v_p].$
\end{definition}

\begin{definition}
Let $K$ be a simplicial complex. A $p${\it -chain}\index{chain} on $K$ is a function $c$ from the oriented $p$-simplices of $K$ to the integers such that: 
\begin{enumerate}
\item $c(\sigma)=-c(\sigma')$ if $\sigma$ and $\sigma '$ are opposite orientations of the same simplex.
\item $c(\sigma)=0$ for all but finitely many oriented $p$-simplices $\sigma$.
\end{enumerate}
We can add $p$-chains by adding their values. They form an abelian group denoted $C_p(K)$. Letting the dimension {\it dim }$K$ be the dimension of the highest dimensional simplex in $K$, we let $C_p(K)=0$ if $p>${\it dim }$K$ or $p<0$.
\end{definition}

\begin{definition}
Let $\sigma$ be an oriented simplex. Then the {\it elementary chain} $c$ correspoding to $\sigma$ is defined by 
\begin{enumerate}
\item $c(\sigma)= 1$
\item $c(\sigma')=-1$ if $\sigma$ and $\sigma '$ are opposite orientations of the same simplex.
\item $c(\tau)=0$ for any other simplex $\tau$.
\end{enumerate}
\end{definition}

\begin{theorem}
$C_p(K)$ is a free abelian group. The basis is obtained by orienting each $p$-simplex and using the corresponding elementary chains.
\end{theorem}

\begin{definition}
We now define a special homomorphism $$\partial_p : C_p(k)\rightarrow C_{p-1}(K).$$ called the {it boundary operator}. From the previous theorem, we just need to define $\partial_p$ on the oriented simplices of $K$.  Let $\sigma=[v_0, v_1, \cdots, v_p]$ be an oriented simplex with $p>0$. Then $$\partial_p\sigma=\partial_p[v_0, v_1, \cdots, v_p]=\sum_{i=0}^p(-1)^i[v_0, \cdots,\hat{v}_i, \cdots, v_p]$$ where the symbol $\hat{v}_i$ means that the vertex $v_i$ is deleted. Since $C_p(K)=0$ for $p<0$, $\partial_p=0$ for $p\leq 0$.
\end{definition}

\begin{figure}[ht]
\begin{center}
  \scalebox{0.4}{\includegraphics{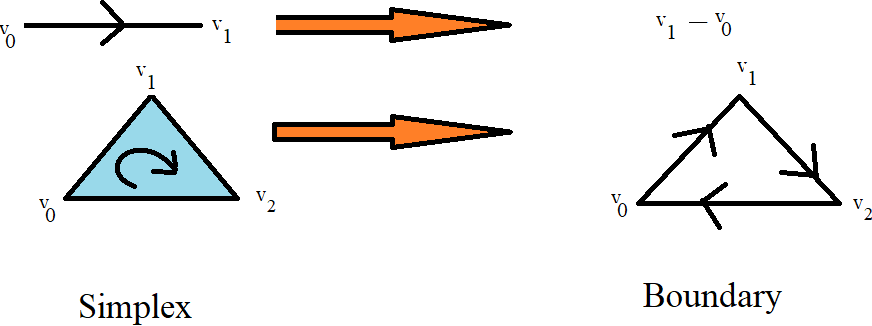}}
\caption{
\rm 
Boundary Examples
}
\end{center}
\end{figure}

\begin{example}
Figure 4.1.4 shows an example of the boundary of both a 1-simplex and a 2-simplex. For the 1-simplex $[v_0, v_1]$, $\partial_1[v_0, v_1]=v_1-v_0.$ so the boundary is just the end point minus the start point. 

For the 2-simplex, $\partial_2[v_0, v_1, v_2]=[v_1,v_2]-[v_0,v_2]+[v_0,v_1]$. Each i-simplex starts at the left vertex and moves towards the right one unless it is negative in the sum in which case we move in the opposite direction. The triangle $[v_0, v_1, v_2]$ is oriented clockwise as shown in the figure. So when taking the boundary, the inside is removed and we move along the edges as follows: Start at $v_1$ and move to $v_2$. Then traverse $[v_0,v_2]$ in the reverse direction so we move from $v_2$ to $v_0$. Then we move from $v_0$ to $v_1$. So we are moving on the boundary of the triangle in the same order as shown on the left hand side of the figure. 
\end{example}

Notice if we take the boundary $\partial_1([v_1,v_2]-[v_0,v_2]+[v_0,v_1])$, we get $v_2-v_1+v_0-v_2+v_1-v_0=0$. So the boundary of the boundary is 0. Does that always happen? 

\begin{theorem}
$\partial_{p-1}\circ\partial_p=0.$
\end{theorem}

Try to prove this yourself from the formula. It is messy but not too hard.

Now we can finally define a homology group.

\begin{definition}
The kernel of $\partial_p : C_p(k)\rightarrow C_{p-1}(K)$ is called the group of $p${\it -cycles} denoted $Z_p(K)$. The image of $\partial_{p+1} : C_{p+1}(k)\rightarrow C_p(K)$ is called the group of $p${\it -boundaries} and denoted $B_p(K)$. The previous theorem shows that $B_p(K)\subset Z_p(K)$. We can take the quotient group $H_p(K)=Z_p(K)/B_p(K)$ and we call it the $p${\it th homology group}\index{homology group} of $K$.
\end{definition}

What do elements of a homology group look like? We know every boundary is a cycle, but we are looking for cycles that are not boundaries. The next two examples from \cite{Mun1} are instructive. 

\begin{figure}[ht]
\begin{center}
  \scalebox{0.4}{\includegraphics{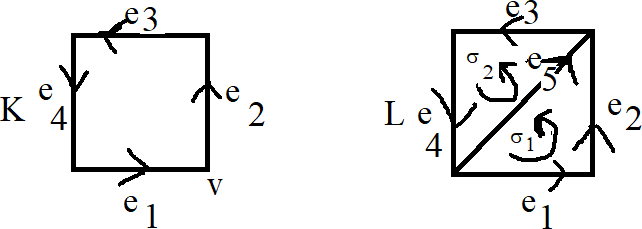}}
\caption{
\rm 
Homology Examples
}
\end{center}
\end{figure}

\begin{example}
Consider the complex $K$ of Figure 4.1.5. It is the boundary of a square with edges $e_1, e_2, e_3, e_4$. The group $C_1(K)$ has rank 4, and a general 1-chain is of the form $c=\sum_{i=1}^4n_ie_i$. The value on vertex $v$ (bottom right) of $\partial_1c$ is $n_1-n_2$. So for c to be a cycle, we must have $n_1=n_2$. Looking at the other vertices, we must have $n_1=n_2=n_3=n_4$. So $Z_1(K)\cong Z$ generated by the cycle $e_1+ e_2+ e_3+ e_4$. Since there are no 2-simplices, $B_1(K)=0$, so $H_1(K)=Z_1(K)=Z$. 
\end{example}

\begin{example}
The complex $L$ on the right side of Figure 4.1.5. is a filled in square with edges $e_1, e_2, e_3, e_4$. We also have an edge $e_5$ so a general 1-chain is of the form $c=\sum_{i=1}^5n_ie_i$. For $c$ to be a cycle we need $n_1=n_2$, $n_3=n_4$, and $n_5=n_3-n_2.$ As an example, the contributions on the top right vertex are $-n_3$, $+n_2$, $+n_5$, so we want $n_2+n_5-n_3=0$. We can arbitraily assign $n_2$ and $n_3$ and the other values are then determined. So $Z_1(L)$ is free abelian with rank 2. We get a basis by letting $n_2=1$ and $n_3=0$ making our generic 1-chain $e_1+e_2-e_5$ and letting $n_2=0$ and $n_3=1$ making our generic 1-chain $e_3+e_4+e_5$. The first is $\partial_2\sigma_1$ and the second is $\partial_2\sigma_2$. So $H_1(L)=Z_1(L)/B_1(L)=0$. The general 2-chain is $m_1\sigma_1+m_2\sigma_2$, which is a cycle only if $m_1=m_2=0$. So $H_2(L)=0$.
\end{example}

This would get pretty ugly if this was the only way to compute homology groups. Suppose you are a typical visitor to Disney World. You would probably want to calculate the homology groups of Spaceship Earth. Imagine spending your entire vacation labeling all of the vertices, edges, and triangles. But there is a lot of cancellation, and we can replace any chain with a chain whose difference is a boundary. (Remember for an abelian group $G$ and subgroup $H$, $g_1$ and $g_2$ are in the same coset of $H$ if $g_1-g_2\in H$.) In a homology group, two chains $c_1, c_2$ are {\it homologous} if $c_1-c_2=\partial d$ for some chain $d$. (Note we will drop the subscript on $\partial$ when it won't cause confusion.) So homologous chains represent the same element in the homology group. If $c=\partial d$ we say that $c$ {\it bounds} or is {\it homologous to zero.} So the way to get some homology groups more easliy is to subtract off pieces of chains that are also boudaries. See \cite{Mun1} for a lot of examples. 

The main thing I want you to take away from the first two examples is to think of what a cycle is and what is a cycle that is not a boundary. In dimension one, a cycle is a loop and if it is not the boundary of something, it corresponds to a hole. So identifying homologous chains as being in the same {\it class}, a class which is a generator in a homology group corresponds to a hole.

Here are some homology groups of interesting objects. See \cite{Mun1} or your favorite textbook for proofs. Basically, you eliminate pieces of chains that are also boundaries. 

\begin{theorem}
$H_p(S^n)=0$ for $p\geq 1$ and $p\neq n$, and $H_n(S^n)\cong Z$ for $n\geq 1$.
\end{theorem}

So now you don't have to compute the homology of Spaceship Earth as it is basically $S^2$.

\begin{theorem}
Let $T$ be a complex whose underlying space is a torus. Then $H_1(T)\cong Z\oplus Z$, and $H_2(T)\cong Z$.
\end{theorem}

\begin{theorem}
Let $K$ be a complex whose underlying space is a Klein bottle. Then $H_1(K)\cong Z\oplus Z_2$, and $H_2(K)\cong 0$.
\end{theorem}

\begin{theorem}
Let $P^2$ be a complex whose underlying space is a projective plane. Then $H_1(P^2)\cong Z_2$, and $H_2(P^2)\cong 0$.
\end{theorem}

You may notice that I haven't discussed zero dimensional homology. I will do that next.

\begin{theorem}
Let $K$ be a complex. Then the group $H_0(K)$ is free abelian. If $\{v_\alpha\}$ is a collection consisting of one vertex from each connected component of $|K|$, then the homology classes of the chains $v_\alpha$ form a basis for $H_0(K)$. In particular, if $|K|$ is connected, than $H_0(K)\cong Z$.
\end{theorem}

The idea behind the proof is to first show that if 2 vertices are in the same component of $|K|$, there is a path consisting of edges from one to the other. Fix a component of $|K|$ and let $v$ be the corresponding vertex in the collection. If $w$ is another vertex, let $[v, a_1], [a_1, a_2], \cdots, [a_{n-1}, w]$ be 1-simplices forming a path from $v$ to $w$. Then the 1-chain  $[v, a_1]+ [a_1, a_2]+ \cdots+[a_{n-1}, w]$ has boundary $w-v$. So $w$ is homologous to $v$ and every 0-chain in K is homologous to a linear combination of the $\{v_\alpha\}$. Also no chain of the form $\sum n_\alpha v_\alpha$ can be a boundary as any one simplex can only lie in a single component. This proves the theorem.

There is another common version of zero dimensional homology called {\it reduced homology}\index{reduced  homology}.

\begin{definition}
Let $\epsilon: C_0(K)\rightarrow Z$ be the surjective homomorphism defined by $\epsilon(v)=1$, for each vertex $v\in K$. If $c$ is a 0-chain, then $\epsilon(c)$ is the sum of the values of $c$ on the vertices of $K$. The map $\epsilon$ is called an {\it augmentation map}\index{augmentation map} for $C_0(k)$. Now if $d=[v_0, v_1]$ is a 1-simplex, then $\epsilon(\partial_1(d))=\epsilon(v_1-v_0)=1-1=0,$, so im $\partial_1\subset\ker \epsilon$. We define the {\it reduced homology group}\index{reduced homology} $\tilde{H}_0(K)=\ker \epsilon/$im $\partial_1$. If $p\geq1$, we define $\tilde{H}_p(K)\equiv H_p(K)$.
\end{definition}

\begin{theorem}
The group $\tilde{H}_0(K)$ is free abelian and $$\tilde{H}_0(K)\oplus Z\cong H_0(K).$$ So for $|K|$ connected, $\tilde{H}_0(K)=0$. We get a basis by fixing an index $\alpha_0$ and choosing a collection  $\{v_\alpha\}$  consisting of one vertex from each connected component of $|K|$. Then $\{v_\alpha-v_{\alpha_0}\}$ form a basis of $\tilde{H}_0(K)$.
\end{theorem}

\begin{definition}
A simplicial complex $K$ is {\it acyclic}\index{acyclic complex} if its reduced homology is zero in all dimensions. 
\end{definition}

Now recall the definition of a cone from Section 2.4.3. It is formed from a space $X$ by taking the cartesian product $X\times I$ and identifying the points $(x,1)$ to a single point. Another way to look at it is that you pick a point outside $X$ and connect every point in $X$ to that point with a straight line. What do you think this operation does to homology?

Imagine $X=S^1$ and take $Y=CX$. We have then turned a circle into an ice cream cone. We have also just plugged up a hole. Taking cones plugs up all of the holes. So homology becomes zero for $p\geq 1$ and every point of $X$ is connected to the external point, so in general, $CX$ is connected even if $X$ didn't start out that way. 

Suppose you start with a simplicial complex $K$. We form a cone $w*K$ in which for any simplex $\sigma=[a_0, \cdots, a_p]$ of $K$, we include the simplex $[w,\sigma]=[w,a_0, \cdots, a_p]$ and all of its faces. Note that any simplex of dimension $p>0$ can be formed by taking a cone of a simplex of dimension $p-1$.

\begin{theorem}
If $w*K$ is a cone, then $\tilde{H}_p(w*K)=0$ for all $p$.
\end{theorem}

So we now know the homology of a simplex of dimension $p>1$ as it is the cone of a simplex of dimension $p-1$. Also, the homology of a point $x$ has no simplices in dimensions higher than 0 and has 1 connected component, so $\tilde{H}_p(x)=0$ for all $p$.

As a final topic in this subsection, I will discuss the idea of {\it relative homology}\index{relative homology}. In this case we have a smaller complex $K_0\subset K$, and chains are represented as a quotient group.. This is surprisingly important for the theory as you will see when we discuss the Eilenberg Steenrod axioms. 

\begin{definition}
Let $K$ be a simplicial complex, and $K_0$ a subcomplex. Then a $p$-chain of $K_0$ can be thought of as a $p$-chain of $K$ whose value is set to zero on any simplices not in $K_0$. So $C_p(K_0)$ is a subgroup of $C_p(K)$ and we can define the group $C_p(K, K_0)$ of {\it relative} $p$-chains of $K$ modulo $K_0$ as the quotient group $$C_p(K, K_0)=C_p(K)/C_p(K_0).$$
\end{definition}

The group $C_p(K, K_0)$ is free abelian with basis $\sigma_i+C_p(K_0)$ where $\sigma_i$ runs over all simplices of $K$ that are not in $K_0$. 

We say that a chain is carried by a subcomplex if it is zero for every simplex not in the subcomplex. Now we need a boundary map $\partial: C_p(K, K_0)\rightarrow C_{p-1}(K, K_0)$. Since the restriction of the usual boundary $\partial$ to $C_p(K_0)$ takes a $p$-chain carried by $K_0$ to a $p-1$-chain carried by $K_0$. We can define a boundary $\partial: C_p(K, K_0)\rightarrow C_{p-1}(K, K_0)$ by $\partial(c_p+C_p(K_0))=\partial(c_p)+C_{p-1}(K_0)$ where the use of $\partial$ as the original or the relative boundary should be clear from the context. Then $\partial\circ\partial=0$. so we have \begin{enumerate}
\item $Z_p(K, K_0)=\ker\partial: C_p(K, K_0)\rightarrow C_{p-1}(K, K_0)$
\item $B_p(K, K_0)=$ im $\partial: C_{p+1}(K, K_0)\rightarrow C_p(K, K_0)$
\item $H_p(K, K_0)=Z_p(K, K_0)/B_p(K, K_0).$
\end{enumerate}

These groups are called the group of {\it relative} $p$-cycles, the group of {\it relative} $p$-boundaries, and the {\it relative homology group} in dimension $p$ respectively. Note that a relative $p$-chain $c_p+C_p(K_0)$ is a relative cycle if and only if $\partial c_p$ is carried by $K_0$ and a relative $p$-boundary if and only if there is a $p+1$-chain $d_{p+1}$ of $K$ such that $c_p-\partial d_{p+1}$ is carried by $K_0$.

We conclude with an important result which will also be included as one of the Eilenberg Steenrod axioms.

\begin{theorem}
{\bf Excision Theorem:} Let $K$ be a complex, and let $K_0$ be a subcomplex. Let $U$ be an open set contained in $|K_0|$, such that $|K|-U$ is the underlying space of a subcomplex $L$ of $K$. Let $L_0$ be the subcomplex of $K$ whose underlying space is $|K_0|-U$. Then we have an isomorphism $H_p(L, L_0)\cong H_p(K, K_0).$ 
\end{theorem}

{\bf Proof:} Consider the composite map $\phi$ $$C_p(L)\rightarrow C_p(K)\rightarrow C_p(K)/C_p(K_0).$$  This is inclusion followed by projection. It is onto since $C_p(K)/C_p(K_0)$ has as a basis all cosets $\sigma_i+C_p(K_0)$ with $\sigma_i$ not in $K_0$ and $L$ contains all of these simplices. The kernel is then $C_p(L_0)$. So the map induces an isomorphism $$C_p(L)/C_p(L_0)\cong C_P(K)/C_p(K_0),$$ for all $p$. Since the boundary operator is preserved by this isomorphism,  $H_p(L, L_0)\cong H_p(K, K_0).$ 

\begin{figure}[ht]
\begin{center}
  \scalebox{0.4}{\includegraphics{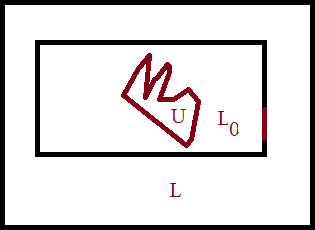}}
\caption{
\rm 
Excision Example
}
\end{center}
\end{figure}

Figure 4.1.6 shows an example. This is the problem of adversarial algebraic topology. You have a nice simplicial complex with a subcomplex and a bully comes and rips a hole in it. As long as the hole is an open set and contained within the subcomplex, you don't have a problem, and everything stays the same. We will see that this is not the case with homotopy groups and is one of the factors that makes them much harder to compute than homology groups.

\subsection{Homology with Arbitrary Coefficients}

This subsection will be pretty short but relevant to topological data analysis. In most applications, we use homology over $Z_2$. 

\begin{definition}
Let $K$ be a simplicial complex and $G$ an abelian group. Then a $p$-chain of $K$ with coefficients in $G$ is a function $c_p$ from the oriented $p$-simplices of $K$ to $G$ that vanishes on all but finitely many $p$-simplices and $c_p(\sigma)=-c_p(\sigma')$ if $\sigma$ and $\sigma'$ are opposite orientations of the same simplex. We write $C_p(K; G)$ where $C_p(K)$ is understood to equal $C_p(K: Z)$. \index{homology with arbitrary coefficients}
\end{definition}

The boundary operator $\partial: C_p(K; G)\rightarrow C_{p-1}(K: G)$. is defined by the formula $\partial(g\sigma)=g(\partial(\sigma))$ where $g\in G$. We can then define cycles, boundaries, and homology with coefficients in $G$ in the obvious way. Relative homology can also be defined.

When we let $G=Z_2$, we look at the simplices as being included or excluded from a chain. Also for $g\in Z_2$, $g=-g$ and $2g=0$. So the boundary map loses its minus signs and we no longer talk about orientation. This changes the calculations a bit for the homology of the objects we have seen. Persistent homology needs $G$ to be a field as we will see in the next chapter. 

If we know homology with integer coefficients. we can find it for any other coefficient group using the {\it Universal Coefficient Theorem}. I will defer this theorem until we talk about cohomology as it uses the machinery of homological algebra, especially tensor and torsion products. There is a Universal Coefficient Theorem for cohomology as well which allows us to compute cohomology groups if we know homology groups. This is good news but also bad news in the sense that cohomology groups do not give us any additional information. The power of cohomology comes from its ring structure using the operation of cup products. I will get to all of that in Chapter 8.

\subsection{Computability of Homology Groups}

In this subsection, I will briefly outline the automated method of computing homology groups described in Munkres\cite{Mun1}. This has some theoretical interest, but there is much better software out there these days. Remember that we are talking state of the art in 1984. Fast computation is not really my area, but I will talk about persistent homology software in the next chapter. Also, we can calculate homology groups with exact sequences in many interesting cases using old fashioned pencil and paper methods. 

We start with a theorem on subgroups of a free abelian group.

\begin{theorem}
Let $A$ be a free abelian group. Then any subgroup $B$ of $A$ is also free abelian. If $A$ is of finite rank $m$, then $B$ has rank $r\leq m$.
\end{theorem}

To compute homology groups, we will work with the matrix of the boundary map. (Remember that chain groups are always free abelian.) We will explicitly define a matrix of a map between two free abelian groups. 

\begin{definition}
Let $G$ and $G'$ be free abelian groups with bases $a_1, \cdots, a_n$ and $a'_1, \cdots, a'_m$ respectively. Then if $f: G\rightarrow G'$ is a homomorphism, then $$f(a_j)=\sum_{i=1}^m\lambda_{ij}a'_i$$ for unique integers $\lambda{ij}$. The matrix $(\lambda_{ij})$ is called the matrix of $f$ relative to the given bases for $G$ and $G'$.
\end{definition}

\begin{theorem}
Let $G$ and $G'$ be free abelian groups of ranks $n$ and $m$ respectively, and let $f: G\rightarrow G'$ be a homomorphism. Then there are bases for $G$ and $G'$ such that relative to these bases, the matrix has the form $$\left[\begin{array}{c c c| c c c}
b_1 &  & 0 & & &\\
 &\ddots & & & 0 &\\
0 & & b_k & & &\\
\hline
& & & & &\\
& 0 & & & 0 &\\
& & & & & 
\end{array}\right]$$
where $b\geq 1$ and $b_1|b_2|\cdots|b_k$. The matrix is called the {\it Smith normal form}\index{Smith normal form}.
\end{theorem}

I will give an outline of the procedure. See \cite{Mun1} or your favorite linear algebra book for details and examples. Also, Matlab has a built in function to do this calculation, and people have developed functions for python and Mathematica. 

Throughout the process we use the three {\it elementary row operations}:\begin{enumerate}
\item Exchange row $i$ and row $j$.
\item Multiply row $i$ by $-1$.
\item Replace row $i$ by (row $i$) +$q$(row $k$), where $q$ is an integer and $k\neq i$.
\end{enumerate}

There are also similar column operations. 

Now for a matrix $A=(a_{ij})$ let $\alpha(A)$ denote the smallest nonzero element of the absolute value $|a_{ij}|$ of the entires of $A$. We call $a_{ij}$ a {\it minimal entry} of $A$ if $|a_{ij}|=\alpha(A)$. The reduction proceeds in two steps. The first brings the matrix to a form where $\alpha(A)$ is as small as possible. The second reduces the dimensions of the matrix involved.

{\bf Step 1:} To decrease the value of $\alpha(A)$ we use the following fact: If the number $\alpha(A)$ fails to divide some entry of $A$ then it is possible to decrease the value of $\alpha(A)$ by applying elementary row and column operations to $A$. The converse is also true.

The idea is then to perform row and column operations until $\alpha(A)$ divides every entry of $A$. See \cite{Mun1} for explicit steps.

{\bf Step 2:} At this point the minimal nonzero element divides all other nonzero elements. Bring it to the top left corner and make it positive. Since it divides all entries in its row and column, we can apply row and column operations to make all of those entries zero.

Now repeat Steps 1 and 2 by ignoring the first row and column and working on the smaller matrix. 

{\bf Step 3:} The algorithm terminates when the smaller matrix is all zeros or disappears. At this point, the matrix is in Smith normal form and we have guaranteed that each diagonal element will divide all elements below it.

We will state a general theorem on {\it chain complexes} which we will now define. Simpliicial homology groups are a special case.

\begin{definition}
A {\it chain complex}\index{chain complex} $\mathcal{C}$ is a sequence $$
\begin{tikzpicture}
  \matrix (m) [matrix of math nodes,row sep=3em,column sep=4em,minimum width=2em]
  {
   \cdots & C_{p+1} & C_p & C_{p-1} & \cdots\\};
  \path[-stealth]
    (m-1-1) edge  (m-1-2)
    (m-1-2)  edge node [above] {$\partial_{p+1}$} (m-1-3)
(m-1-3)  edge node [above] {$\partial_p$} (m-1-4)
(m-1-4) edge  (m-1-5);

\end{tikzpicture}$$ of abelian groups $C_i$ and homomorphisms $\partial_i$ indexed by the integers such that $\partial_p\circ\partial_{p+1}=0$ for all $p$. The $p$th homology group of $C$ is defined by the equation $H_p(\mathcal{C})=\ker \partial_p/$im $\partial_{p+1}$.
\end{definition}

\begin{theorem}
Let $\{C_p, \partial_p\}$ be a chain complex such that each group $C_p$ is free and of finite rank. Then for each $p$, there are subgroups $U_p, V_p, W_p$ of $C_p$ such that $$C_p=U_p\bigoplus V_p\bigoplus W_p,$$ where $\partial_p(U_p)\subset W_{p-1}$ and $\partial_p(V_p)=\partial_p(W_p)=0.$ In addtion, there are bases for $U_p$ and $W_{p-1}$ relative to which $\partial_p: U_p\rightarrow W_{p-1}$ has a matrix of the form $$B=\left[\begin{matrix}
b_1 & & 0\\
& \ddots &\\
0 & & b_k
\end{matrix}\right],$$ where $b_i\geq 1$ and $b_1|b_2|\cdots|b_k.$
\end{theorem}

{\bf Outline of Proof \cite{Mun1}:}

Step 1:

 Let $Z_p$ be the group of $p$-cycles and $B_p$ be the group of $p$-boundaries. $W_p$ consists of all elements $c_p\in C_p$ such that $mc_p\in B_p$ for some non-zero integer $m$. Since $C_p$ is torsion free, $mc_p\neq 0$ and $mc_p=\partial_{p+1}d_{p+1}$ so $m\partial_p c_p=\partial mc_p=0$ implies that $\partial_p c_p=0$, so $W_p\subset Z_p$. $W_p$ is called the group of {\it weak boundaries}. Munkres shows that $W_p$ is a direct summand of $Z_p$ and lets $V_p$ be the group such that $Z_p=V_p\oplus W_p$.

Step 2:

Choose basis $e_1, \cdots, e_n$ for $C_p$ and $e'_1,\cdots,e'm$ for $C_{p-1}$ so that the matirx of $\partial _p: C_p\rightarrow C_{p-1}$ has the Smith normal form 

$$\kbordermatrix{
\mbox{} & e_1 & \cdots & e_k & \vrule & e_{k+1} & \cdots & e_n \\
e'_1 & b_1 &  & 0 & \vrule & & &\\
\cdots & &\ddots & & \vrule & & 0 &\\
e'_k & 0 & & b_k & \vrule & & &\\
\hline
e'_{k+1} & & & & \vrule & & &\\
\cdots & & 0 & & \vrule & & 0 &\\
e'_m & & & & \vrule & & & 
}$$ where $b_i\geq 1$ and $b_1|b_2|\cdots|b_k.$ Munkres shows that the following hold:
\begin{enumerate}
\item $e_{k+1},\cdots, e_n$ is a basis for $Z_p$. 
\item $e'_1,\cdots, e'_k$ is a basis for $W_{p-1}$
\item $b_1e'_1,\cdots, b_ke'_k$ is a basis for $B_{p-1}$
\end{enumerate}

Step 3: 

Now choose bases for $C_p$ and $C_{p-1}$ using Step 2. Define $U_p$ to be the group spanned by $e_1, \cdots, e_k$. then $C_p=U_p\oplus Z_p$. Using Step 1, choose $V_p$ so that $Z_p=V_p\oplus W_p$. Then we have a decompostition of $C_p$ such that $\partial_p(V_p)=\partial_p(W_p)=0$. The existence of the desired bases for $U_p$ and $W_{p-1}$ follow from Step 2. $\blacksquare$

Now we state the method of computing homology groups. 

\begin{theorem}
The homology groups of a finite complex $K$ are effectively computable.
\end{theorem}

{\bf Proof:}
By the preceeding theorem. we have a decomposition $$C_p(K)=U_p\bigoplus Z_p=U_p\bigoplus V_p\bigoplus W_p,$$ where $Z_p$ is the group of $p$-cycles and $W_p$ is the group of weak $p$-boundaries. Now $$H_p(K)=Z_p/B_p\cong V_p\bigoplus(W_p/B_p)\cong (Z_p/W_p)\bigoplus (W_p/B_p).$$ By the proof of the previous theorem, $Z_p/W_p$ is free and $W_p/B_p$ is the torsion subgroup. 

Orient the simplices of $K$ and choose a basis for the groups $C_p(K)$. The matrix of $\partial_p$ has entries in the set $\{-1, 0, 1\}$. Reducing it to Smith normal form and looking at Step 2 of the proof of the previous theorem we have these facts \begin{enumerate}
\item The rank of $Z_p$ equals the number of zero columns.
\item The rank of $W_{p-1}$ equals the number of nonzero rows.
\item There is an isomorphism $$W_{p-1}/B_{p-1}\cong Z_{b_1}\bigoplus  Z_{b_2}\bigoplus\cdots\bigoplus  Z_{b_k}.$$
\end{enumerate}

So we get the torsion coefficients of $K$ of dimension ${p-1}$ from the Smith normal form of the matrix of $\partial_p$. This normal form also gives the rank of $Z_p$ and the normal form of $\partial_{p+1}$ gives the rank of $W_p$. The difference of these numbers is the rank of $Z_p$ minus the rank of $W_p$ which is the betti number of $K$ in dimension $p$. $\blacksquare$

\begin{example}
Suppose we want to compute $H_1(K)$ and we have the following matrices already reduced to Smith normal form:$$\partial_2=\left[\begin{array}{c c c c| c c c}
1 & 0 & 0 & 0 & 0 & 0 & 0\\
 0 & 3 & 0 & 0 & 0 & 0 & 0\\
0 & 0 & 6 & 0 & 0 & 0 & 0\\
0 & 0 & 0 & 12 & 0 & 0 & 0\\
\hline
0 & 0 & 0 & 0 & 0 & 0 & 0\\
0 & 0 & 0 & 0 & 0 & 0 & 0\\
0 & 0 & 0 & 0 & 0 & 0 & 0\\
0 & 0 & 0 & 0 & 0 & 0 & 0
\end{array}\right],$$ and
$$\partial_1=\left[\begin{array}{c c| c c c c c c}
1 & 0 & 0 & 0 & 0 & 0 & 0 & 0\\
0 & 5 & 0 & 0 & 0 & 0 & 0 & 0\\
\hline
0 & 0 & 0 & 0 & 0 & 0 & 0 & 0\\
\end{array}\right].$$

Then we have the rank of $Z_1$ is 6. The rank of $W_1$ taken from $\partial_2$ is 4. So the betti number is 6-4=2. The torsion coefficients from the matrix for $\partial_2$ are 3, 6, and 12. so $$H_1(K)\cong Z\bigoplus Z\bigoplus Z_3 \bigoplus Z_6\bigoplus Z_{12}.$$
\end{example}

\subsection{Homomorphisms Induced by Simplicial Maps}

I promised you earlier to produce a functor from topological spaces to groups in each dimension. We have the groups now but we need to take maps between spaces and turn them into homomorphisms between the groups. We would like to do this for any continuous map but we are not quite there yet. Earlier, though we talked about simplicial maps. From now on we will talk about simplicial maps $f: K\rightarrow L$ and we mean that $f$ is a continuous map of $|K|$ into $|L|$ that maps each simplex of $K$ into a simplex of $L$ of the same or lower dimension. So $f$ maps each vertex of $K$ into a vertex or $L$, and this map equals the simplical map we defined in Theorem 4.1.2. 

\begin{definition}
Let $f: K\rightarrow L$ be a simplicial map. If $[v_0,\cdots, v_p]$ is s simplex of $K$ then the points $f(v_1),\cdots, f(v_p)$ span a simplex of $L$. Define a homomorphism $f_\#:C_p(K)\rightarrow C_p(L)$ by $f_\#([v_0,\cdots, v_p])=[f(v_1),\cdots, f(v_p)]$ if all of the vertices are distinct and $f_\#([v_0,\cdots, v_p])=0$ otherwise. $f_\#$ is called a {\it chain map}\index{chain map} induced by $f$.
\end{definition}

This important theorem can pe proved from the above definition. Care must be taken when the vertices of the image are not distinct. 

\begin{theorem}
The homomorphism $f_\#$ commutes with $\partial$. Therefore, $f_\#$ induces a homomorphism $f_*: H_p(k)\rightarrow H_p(L)$.The chain map $f_\#$ also commutes with the augmentation map $\epsilon$, so it induces a homomorphism $f_*$ of reduced homology groups.
\end{theorem}

\begin{theorem}
\begin{enumerate}
\item Let $i:K\rightarrow K$ be the identity simplicial map. Then $i_*: H_p(K)\rightarrow H_P(k)$ is the identity homomorphism for all $p$.
\item Let $f: K\rightarrow L$ and $g: L\rightarrow M$ be simplicial maps. Then $(gf)_*=g_*f_*.$ This makes homology a functor from topological spaces and simplicial maps to groups and homomorphisms.
\end{enumerate}
\end{theorem}

It turns out that more than one simplicial map can lead to the same homomorphism on homology groups. To see when this happens we need the following definition:

\begin{definition}
Let $f, g: K\rightarrow L$ be simplicial maps. Suppose that for each $p$ suppose there is a homomorphism $D: C_p(K)\rightarrow C_{p+1}(L)$ such that $$\partial D+D\partial=g_\#-f_\#.$$ Then $D$ is called a {\it chain homotopy}\index{chain homotopy} between $f_\#$ and $g_\#$.
\end{definition}

\begin{theorem}
If there is a chain homotopy between $f_\#$ and $g_\#$ then the induced homomorphisms $f_*$ and $g_*$ for both ordinary and reduced homology are equal.
\end{theorem}

{\bf Proof} Let $z$ be a $p$-cycle of $K$. Then $g_\#(z)-f_\#(z)=\partial D(z)+D\partial(z)=\partial D(z)+0.$ So  $f_\#(z)$ and $g_\#(x)$  differ by a boundary and are in the same homology class. Thus  $f_*=g_*$.  $\blacksquare$

We conclude with a condition that allows us to find a chain homotopy.

\begin{definition}
Let $f, g: K\rightarrow L$ be simplicial maps. These maps are {\it contiguous} if for each simplex $[v_0, \cdots, v_p]$ of $K$, the vertices $f(v_0), \cdots, f(v_p), g(v_0), \cdots, g(v_p)$ span a simplex of $L$. This simplex can have dimension anywhere between 0 and $2p+1$ depending on how many are distinct.
\end{definition}

The proof of this result is lengthy but not too difficult. See \cite{Mun1}.

\begin{theorem}
Let $f, g: K\rightarrow L$ be contiguous simplicial maps. Then there is a chain homotopy between $f_\#$ and $g_\#$, so $f_*$ and $g_*$ for both ordinary and reduced homology are equal.
\end{theorem}

\subsection{Topological Invariance}

One thing we would like is for homeomorphic topological spaces to have the same homotopy groups and for homeomorphisms to translate to isomorphisms on homology groups. This does turn out to be true, but unfortunately, having the same homology groups does not imply that two spaces are homemorphic. This issue led to the discovery of cohomology and Steenrod operations which provide additional algebraic structure.

For now, we know that simplicial complexes with invertible simplicial maps have the same homology groups. So we have some issues to consider. First of all, can any topological space be represented as a simplicial complex? The answer is no. For this to be true, the space needs to be {\it triangulable} which means it is homeomorphic to a simplicial complex. Almost all spaces you could name are triangulable and you will never have to deal with this issue in topological data analysis. In Section 4.3.2, I will describe the only counterexample I have ever seen.

The hard part is moving from simplicial maps to continuous maps.  Munkres \cite{Mun1} gives the entire proof of topological invariance (i. e. homeomorphic spaces lead to isomorphic homology groups) in sections 14-18. As it is long and complicated, I won't repeat it here but just describe some of the ideas behind it.

We will assume all of our spaces are the underlying spaces or {\it polytopes}\index{polytope} of simplicial complexes. 

\begin{definition}
If $K$ is a simplicial complex, then the {\it star} of a vertex $v\in K$ denoted $St(v)$\index{star of a vertex} is the union of the interiors of all simplices in $K$ that have $v$ as a vertex.
\end{definition}

\begin{definition}
Let $h: |K|\rightarrow |L|$ be a continuous map. If $f: K\rightarrow L$ is a simplicial map such that $h(St(v))\subset St(f(v))$ for each vertex $v\in K$, then $f$ is called a {\it simplicial approximation}\index{simplicial approximation} to $h$.
\end{definition}

The map $f$ is an approximation to $h$ in the sense that given $x\in K$, there is a simplex in $L$ that contains both $f(x)$ and $h(x)$. 

To find such an approximation, we subdivide $K$ into smaller simplices. 

\begin{definition}
Let $K$ be a complex. A complex $K'$ is a {\it subdivision} of $K$ if each simplex of $K'$ is contained in a simplex of $K$ and each simplex of $K$ equals the union of finitely many simplices of $K'$. 
\end{definition}

Now recall that a simplex $\sigma$ with vertices $x_0, \cdots, x_k$ is the set of all points $x=\sum_{i=0}^k t_ix_i$ where $\sum t_i=1$ and $t_i\geq 0$. The $t_i$ are called the {\it barycentric coordinates}\index{barycentric coordinates}. The place where they are all equal is $$\hat{\sigma}=\sum_{i=0}^k \frac{1}{k+1} v_i.$$ $\hat{\sigma}$ is called the {\it barycenter}\index{barycenter} of $\sigma$.

Now we divide up $K$ into a subdivision $sd(K)$ as follows. First, include all of the vertices in $sd(K)$ that were originally in $K$.  We divide up the $p$-skeletons for $p=1, 2, \cdots$. Suppose we already have subdivided the $p$-skeleton of $K$. Then for any $p+1$-simplex $\sigma$ of $K$, find its barycenter $\hat{\sigma}$. Then for every for every $p$-simplex $\tau\subset\sigma$ with vertices $v_0,\cdots,v_p$ in $sd(K)$ we include the $p+1$-simplex $[\hat{\sigma}, v_0,\cdots,v_p]$. When we are done, the new complex $sd(K)$ is called the {\it first barycentric subdivision}\index{barycentric subdivision} of $K$. We can do the process a second time creating the {\it second barycentric subdivision} $sd^2(K)$. Continuing in this way, we can make the simplices in $K$ as small as we want. Figure 4.1.7 demonstrates the process and shows the first 2 subdivisions of the complex on the left.

\begin{figure}[ht]
\begin{center}
  \scalebox{0.4}{\includegraphics{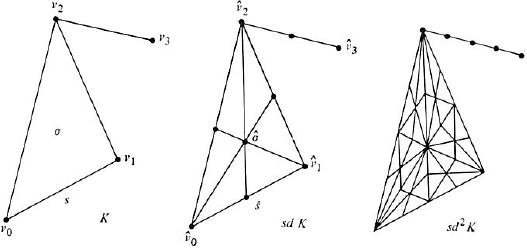}}
\caption{
\rm 
Barycentric Subdivision \cite{Mun1}
}
\end{center}
\end{figure}

The idea in finding a simplicial approximation to a continuous map is to subdivide $K$ until each star of a vertex in $K$  is in one of the sets $h^{-1}(St(w))$ where $w$ is a vertex of $L$. This gives the following: 

\begin{theorem}
{\bf The finite simplicial approximation theorem:} Let $K$ and $L$ be complexes where $K$ is finite. Given a continuous map $h: |K|\rightarrow |L|$ there is an integer $N$ such that $h$ has a simplicial approximation $f: sd^N(K)\rightarrow L$.
\end{theorem}

With this result, we are now ready to show topological invariance.

\begin{definition}
Let $K$ and $L$ be simplicial complexes and let $h: |K|\rightarrow |L|$ be continuous. Choose a subdivision $K'$ of $K$ such that $h$ has a simplicial approximation $f: K'\rightarrow L$. So $K'=sd^N(K)$ for some $N$ if $K$ is finite. Let $\frak{C}(K)$ be the chain complex consisting of the chain groups of $K$ and the usual boundaries. Let $\lambda: \frak{C}\rightarrow\frak{C'}$ be the chain map representing subdivision. (Munkres \cite{Mun1} Section 17 makes this precise.) Then the {\it homomorphism induced by} $h$ is $h_*: H_p(K)\rightarrow H_p(L)$ defined by $h_*=f_*\lambda_*$.
\end{definition}

Munkres shows that these maps have all the nice properties we had with simplicial maps. 

\begin{theorem}
The identity map $i: |K|\rightarrow |K|$ induces the identity homomorphism $i_*: H_p(K)\rightarrow H_p(K)$. If $h: |K|\rightarrow |L|$ and $k: |L|\rightarrow |M|$ are continuous, then $(kh)_*=k_*h_*$. The same holds for reduced homology.
\end{theorem}

\begin{theorem}
{\bf Topological invariance of homology groups:} If $h: |K|\rightarrow |L|$ is a homeomorphism then $h_*: H_p(K)\rightarrow H_p(L)$ is an isomorphism. The same holds for reduced homology.
\end{theorem}

This is immediate from the previous theorem by composing $h$ and $h^{-1}$ which exists since $h$ is a homeomorphism. I will remark that everything holds for relative homology as well.

We can even do better than this. I will show this result holds for the weaker condition of {\it homotopy equivalence}. I will discuss this next.

\begin{definition}
If $X$ and $Y$ are topological spaces, two continuous maps $h, k: X\rightarrow Y$ are {\it homotopic}\index{homotopic maps} if there is a continuous map $F: X\times I\rightarrow Y$ such that $F(X,0)=h(x)$ and $F(x,1)=k(x)$. The map $F$ is called a {\it homotopy} of $h$ to $k$. 
\end{definition}

We can think of $F$ as continuously deforming $h$ into $k$. The set of all maps from $X$ to $Y$ is partitioned into equivalence classes called {\it homotopy classes} where two maps are in the same class if they are homotopic to each other. 

As a preview of chapter 9, the elements of the homotopy group $\pi_n(X)$ will be homotopy classes of continuous maps from $S^n$ to $X$.

\begin{theorem}
 If $h, k: |K|\rightarrow |L|$ are homotopic, then $h_*, k_*: H_p(K)\rightarrow H_p(L)$ are equal. The same holds for reduced homology.
\end{theorem}

For relative homology, we modify the definition of homotopy to include maps of pairs. Let $h, k: (|K|, |K_0|)\rightarrow (|L|,|L_0|)$ be maps of pairs where $K_0$ is a subcomplex of $K$ and $L_0$ is a subcomplex of $L$. This means that they map $|K|$ to $L$ while mapping $K_0$ into $L_0$. Then they are homotopic if there is a homotopy $H: |K|\times I\rightarrow |L|$ of $h$ to $k$ such that $H$ maps $|K_0|\times I$ into $|L_0|$.

\begin{theorem}
 If $h, k$ are homotopic as maps of pairs of spaces, then $h_*=k_*$ as homomorphisms of relative homology groups.
\end{theorem}

\begin{definition}
Two spaces are {\it homotopy equivalent}\index{homotopy equivalent} or have the same {\it homotopy type}\index{homotopy type} if there are maps $f: X\rightarrow Y$ and $g: Y \rightarrow X$ such that $gf=i_X$ and $fg=i_Y$. The maps $f$ and $g$ are called {\it homotopy equivalences}\index{homotopy equivalence} and $g$ is a {\it homotopy inverse}\index{homotopy inverse} of $f$.
\end{definition}

Recall that a complex is acyclic if all reduced homology groups are zero. 

\begin{definition}
A space X is {\it contractible}\index{contractible}  if it has the homotopy type of a single point.
\end{definition}

Now you know what every word in my dissertation title \cite{Pos1} means.

\begin{example}
The unit ball is contractible. Let $F(x,t)=(1-t)x$, for $x\in B^n$. Then $F(x,0)=x$ and $F(x,1)=0$. Note that the unit ball and the origin are not homeomorphic. The constant map $f:B^n\rightarrow 0$ is  obviously not one-to-one and so not invertible.
\end{example}

So now we have a stronger result. If two spaces are not necessarily homeomorphic but only homotopy equivalent, they have isomorphic homology groups. This is immediate from Theorem 4.1.24 and formally stated here:

\begin{theorem}
  If $h: |K|\rightarrow |L|$ is a homotopy equivalence, then $h_*$ is an isomorphism on homology groups. In particular, if $|K|$ is contractible, then $K$ is acyclic.
\end{theorem}

\begin{definition}
Let $A\subset X$. A {\it retraction}\index{retraction} of $X$ onto $A$ is a continuous map $r: X\rightarrow A$ such that $r(a)=a$ for each $a\in A$. If there is a retraction of $X$ onto $A$, we say that $A$ is a {\it retract} of $X$. A {\it deformation retraction}\index{deformation retraction} of $X$ onto $A$ is a continuous map $F: X\times I\rightarrow X$ such that \begin{enumerate}
\item $F(x,0)=x$ for $x\in X$.
\item $F(x,1)\in A$ for $x\in X$.
\item $F(a,t)=a$ for $a\in A$.
\end{enumerate}
If such an $F$ exists, then $A$ is called a {\it deformation retract} of $X$.
\end{definition}

Every deformation retract is also a retract. To see this, let $r=F(x,1)$. Then $r: X\rightarrow A$ and $r(a)=a$ for $a\in A$. A retract is not necessarily a deformation retract as the next example shows. 

\begin{example}
Let $X=S^1$ (the unit circle) and $A=(0,1)$, i.e. $A$ is a single point. Letting $r(x)=(0,1)$ for $x\in S^1$. Then $r$ is obviously a retraction. But $A$ is not a deformation retract of $S^1$. To see this, let $X$ be any topological space, $A$ a deformation retract of $X$. If $r: X\rightarrow A$ is defined by $r=F(x,1)$ where $F$ is as above and $j: A\rightarrow X$ is inclusion, then $F$ is a homotopy between $jr$ and identity $i_X$ of $X$. So $r$ and $j$ are homotopy inverses of each other and the homology groups of $X$ and $A$ are isomorphic. In our case, $H_1(S^1)\cong Z$ and $H_1(A)=0$ for $A$ a single point. So $A$ is a retract of $X$ but not a deformation retract.
\end{example}

\begin{theorem}
 The unit sphere $S^{n-1}$ is a deformation retract of punctured Euclisean space $R^n-0$.
\end{theorem}

{\bf Proof:} Let $X=R^n-0$. Define $F: X\times I\rightarrow X$ by $$F(x,t)=(1-t)x+tx/||x||,$$ where $||x||$ is the standard Euclidean norm, i.e. if $x=(x_1, \cdots, x_n)$ then $||x||=\sqrt{x_1^2+\cdots+x_n^2}.$ Then $F$ moves $x$ along the line through the origin to it's intersection with the unit sphere, so it is a deformation retraction of $R^n-0$ onto  $S^{n-1}$.  $\blacksquare$

Now we can finally settle the question about Euclidean spaces of different dimension being homeomorphic.

\begin{theorem}
 The Euclidean spaces $R^m$ and $R^n$ are not homeomorphic if $m\neq n$.
\end{theorem}

{\bf Proof:} Assume $m, n>1$ as otherwise, we can remove a point and $R^1$ becomes disconnected. Remove a point from each of $R^m$ and $R^n$. Then the resulting spaces are homeomorphic to $R^m-0$ and $R^n-0$ respectively.  But these are deformation retracts of $S^{m-1}$ and $S^{n-1}$ respectively and so have isomorphic homology groups to them. But $H_{m-1}(S^{m-1})\cong Z$, while $H_{m-1}(S^{n-1})=0$. So their homology groups differ and they can not be homeomorphic. $\blacksquare$

Note that we could have made a similar argument by adding a point rather than removing one. Adding a point called the {\it point at infinity} to $R^n$ produces the sphere $S^n$. This is called a {\it one point compactification}. It is most easily visualized for $n=2$. 

\begin{figure}[ht]
\begin{center}
  \scalebox{0.4}{\includegraphics{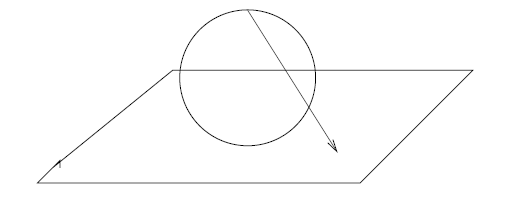}}
\caption{
\rm 
Stereographic Projection 
}
\end{center}
\end{figure}

The {\it stereographic projection}\index{stereographic projection} from $S^2$ to $R^2$ involves placing the South pole of $S^2$ at the origin of $R^2$. (See Figure 4.1.8.) Given any point $x\in S^2$ other than the North pole, we map it to the point of $R^2$ which is on the line containing the North pole and $x$. The South pole maps to the origin and the equator maps to the unit circle in $R^2$. It can be checked that this map is a homeomorphism between $S^2$ with the North pole removed and $R^2$. As you move closer to the north pole the line gets more and more horizontal and intersects  in a point of increasing magnitude so we can think of the projection of the North pole as infinity. With this projection and similar ones for any $R^m$, we can see that $R^n$ and $R^m$ are homeomorphic if and only if $S^n$ and $S^m$ have the same homology groups. As we saw above, this is only true if $m=n$.

\subsection{Application: Maps of Spheres and Fixed Point Theorems, and Euler Number}

At this point, I can give you an idea of how obstruction theory works. I will make this more precise when we see cohomology and homotopy. Suppose, we have a simplicial complex $K$ and $L$ is a subcomplex of $K$. Let $g: L\rightarrow Y$ and we want to extend $g$ to a function $f: K\rightarrow Y$ where $f$ restricted to $L$ is equal to $g$. We will suppose $Y=S^n$. Then the idea is to extend $F$ to the vertices of $K-L$, then the edges, then the 2-simplices, etc. We can always extend $g$ to the vertices and if $Y$ is path connected, we can extend to the edges. (This is fine if $Y$ is a sphere.) Suppose we have extended $g$ to the $n$-simplices. For an $n+1$-simplex, $f$ is defined on the boundary which we know is homeomorphic to $S^n$. If $Y=S^n$, then extending $f$ to the interior would amount to extending a map $h: S^n\rightarrow S^n$ to a map $k: B^{n+1}\rightarrow S^n$. Maps between spheres of the same dimension are classified by an integer called the {\it degree}. That will be the first topic in this section.

So if $Y=S^n$, what do we do if we are not at the step where we go from $n$-simplices to $n+1$-simplices?  In general, if $Y=S^n$ and we are going from $m$ to $m+1$, then $f$ restricted to the boundary of an $m+1$ simplex is a map from $S^m$ to $S^n$. The homotopy classes of these maps form the group $\pi_m(S^n)$. As we are about to see, we can extend to the interior of this simplex if $f$ is homotopic to a constant map which means that $f$ is the zero element in this group. If $n=m$, we are in the case described here. If $m<n$,  $\pi_m(S^n)=0,$ so the extension is always possible. If $m>n$, things get harder. The groups $\pi_m(S^n)$ have not all been calculated even for $n=2$, and unlike the case of homology, they are mostly nonzero. As we will see, though, there are many special cases where the problem is tractable and may be of some interest in data science. 

\begin{definition}
Let $n\geq 1$ and $f: S^n\rightarrow S^n$ be a continuous map. If $\alpha$ is one of the two generators of the group $H_n(S^n)\cong Z$, then $f_*(\alpha)=d\alpha$ for some integer $d$ as these are the only homomorphisms from $Z$ to itself. Since $f_*(-\alpha)=d(-\alpha)$ it doesn't matter which generator we pick. The number $d$ is called the {\it degree}\index{degree of a map} of $f$. 
\end{definition}

As I mentioned above, a constant map can always be extended. If $L$ is a subcomplex of $K$ and $g(l)=y_0$ for $l\in L$ and a fixed $y_0$ in $Y$, then this map is continuous and we can extend $g$ to the constant map $f(k)=y_0$ for $k\in K$. A constant map will have degree 0 as we will now show.

\begin{theorem}
Degree has the following properties:\begin{enumerate}
\item If $f$ is homotopic to $g$ then deg $f=$ deg $g$.
\item If $f$ extends to a continuous map $h: B^{n+1}\rightarrow S^n$ then deg $f$=0.
\item The identity map has degree 1.
\item $\deg(fg)=(\deg(f))(\deg(g))$.
\end{enumerate}
\end{theorem}

{\bf Proof:} Property 1 holds since homotopic maps give rise to the same homomorphisms in homology. Property 2 follows from the fact that $f_*: H_n(S^n)\rightarrow H_n(S^n)$ equals the composite $$
\begin{tikzpicture}
  \matrix (m) [matrix of math nodes,row sep=3em,column sep=4em,minimum width=2em]
  {
H_n(S^n) & H_n(B^{n+1}) & H_n(S^n),\\};
  \path[-stealth]
    (m-1-1)  edge node [above] {$j_*$} (m-1-2)
(m-1-2)  edge node [above] {$h_*$} (m-1-3);
\end{tikzpicture}$$ where $j$ is inclusion. Since $B^{n+1}$ is acyclic, the composite is the zero homomorphism. Properties 3 and 4 follow from Theorem 4.1.22. $\blacksquare$

\begin{theorem}
There is no retraction $r: B^{n+1}\rightarrow S^n$.
\end{theorem}

{\bf Proof:} The retraction $r$ would be an extension of the identity map $i: S^n\rightarrow S^n$. Since $i$ has degree $1\neq 0$, the extension can not exist. $\blacksquare$

We can now easily prove an interesting theorem about fixed points, ie points where $f(x)=x.$

\begin{theorem}
{\bf Brouwer Fixed-Point Theorem:} Every continuous map $\phi: B^n\rightarrow B^n$ has a fixed point.
\end{theorem}

{\bf Proof:} Suppose $\phi: B^n\rightarrow B^n$ has no fixed point. Then we can define a map $h: B^n\rightarrow S^{n-1}$ by the equation $$h(x)=\frac{x-\phi(x)}{||x-\phi(x)||},$$ since we know $x-\phi(x)\neq 0$. (Note that dividing by the norm makes the points have norm 1 so they lie on the surface of the sphere.) If $f: S^{n-1}\rightarrow S^{n-1}$ is the restriction of $h$ to $S^{n-1}$ then $f$ has degree 0 by property 2 of Theorem 4.1.29. 

But we can show that $f$ has degree 1. Let $H: S^{n-1}\times I\rightarrow S^{n-1}$ be a homotopy defined by $$H(u, t)=\frac{u-t\phi(u)}{||u-t\phi(u)||}.$$ Then the denominator is always nonzero since for $t=1, u\neq \phi(u)$ as $\phi$ does not have a fixed point, and for $t<1,$ we know $||u||=1$ and  $||t\phi(u)||=t||\phi(u)||\leq t<1.$ The map $H$ is a homotopy between $f$ and the identity of $S^{n-1}$ so $f$ has degree 1, which is a contradiction. Hence $\phi$ must have a fixed point.
$\blacksquare$

Next we will see a few more examples of degrees of certain maps and what can be done with them.

\begin{definition}
The {\it antipodal map} $a: S^n\rightarrow S^n$ is the map defined by the equation $a(x)=-x$ for $x\in S^n$. It takes a point to its antipodal point which lies on the other side of the sphere on a line from $x$ to the center of the sphere.
\end{definition}

\begin{theorem}
Let $n\geq 1.$ The degree of the antipodal map $a: S^n\rightarrow S^n$ is $(-1)^{n+1}$.
\end{theorem}

The idea of the proof is that if we reflect a coordinate in $R^{n+1}$ by multiplying it by $(-1)$, the reflection map has degree -1. Then just compose the $n+1$ reflection maps corresponding to each coordinate. See \cite{Mun1} for the full proof which has to take into account simplicial approximations. 

\begin{theorem}
If $h: S^n\rightarrow S^n$ has a degree different from $(-1)^{n+1}$ then $h$ has a fixed point.
\end{theorem}

{\bf Proof:} Suppose $h$ has no fixed point. We will prove $h$ is homotopic to $a$. Otherwise, since $h(x)$ and $x$ are not antipodal, we can form a homotopy which moves $h(x)$ to $-x$ along the shorter great circle path. The equation is $$H(x,t)=\frac{(1-t)h(x)+t(-x)}{||(1-t)h(x)+t(-x)||}.$$ So we are done as long as the denominator is never 0. Suppose we had $(1-t)h(x)=tx$ for some $x$ and $t$. Then taking norms of both sides, we would get $1-t=t=1/2.$ So $h(x)=x$ which is a contradiction since $h$ has no fixed point. So $h$ is homotopic to $a$ and has degree $(-1)^{n+1}$.
$\blacksquare$

\begin{theorem}
If $h: S^n\rightarrow S^n$ has a degree different from 1, then $h$ carries some point $x$ to its antipode $-x$. 
\end{theorem}

{\bf Proof:} In this case, if $a$ is the antipodal map, then $ah$ has degree different from $(-1)^{n+1}$, so $ah$ has a fixed point $x$. Then $a(h(x))=x$, so $h(x)=-x$. 
$\blacksquare$

An immediate consequence is a famous theorem about tangent vector fields. These fields associate a tangent vector to each point on a sphere. A tangent vector at a point $x$ must be perpendicular to the line from the center of the sphere to $x$ so the inner product with $x$ must be zero. 

\begin{theorem}
$S^n$ has a non-zero tanngent vector field if and only if $n$ is odd. 
\end{theorem}

This shows that since the surface of your head has the homotopy type of $S^2$, you can't brush your hair so that every hair is lying perfectly flat.

{\bf Proof:} If $n$ is odd, let $n=2k-1$. Then for $x$ in $S^n$ let $v(x)=v(x_1, \cdots, x_{2k})=(-x_2, x_1, -x_4, x_3, \cdots, -x_{2k}, x_{2k-1}$. Then the inner product $\langle x, v(x)\rangle=x_1(-x_2)+x_2x_1+\cdots+x_{2k-1}(-x_{2k})+x_{2k}x_{2k-1}=0$. So $v(x)$ is perpendicular to $x$ and the $\{v(x)\}$ form a nonzero tangent vector field to $S^n$.

If the  $\{v(x)\}$ form a nonzero tangent vector field to $S^n$, then for each $x$, $h(x)=v(x)/||v(x)||$ is a map of $S^n$ into $S^n$. Since $h(x)$ is perpendicular to $x$, $h(x)$ can not equal $x$ or $-x$. So the degree of $h$ is $(-1)^{n+1}$ and the degree of $h$ is 1. So $n$ must be odd.
$\blacksquare$

Munkres introduces the {\it Euler number} at this point as part of a discussion of more general fixed point theorems. I won't get into those here but you should know what an Euler number is as it is a useful topological invariant.

\begin{definition}
Let $K$ be a finite complex. Then the {\it Euler number}\index{Euler number} or {\it Euler characteristic}\index{Euler characteristic} $\chi(K)$\index{$\chi(K)$} is defined by the equation $$\chi(K)=\sum(-1)^p {\rm rank}(C_p(K)).$$
\end{definition}

For example, if $K$ is 2-dimensional then $\chi(K)=\# {\rm vertices}-\#{\rm edges}+\#{\rm faces}.$ Munkres \cite{Mun1} proves the following: 

\begin{theorem}
Let $K$ be a finite complex. Let $\beta_p={\rm rank}(H_p(K)/T_p(K))$, i.e. the betti number of $K$ in dimension $p$. Then $$\chi(K)=\sum_p(-1)^p\beta_p.$$
\end{theorem}

We note that for surfaces, there is a closely related idea called the {\it genus}\index{genus}. The relationship is $\chi(X)=2-2g$ where $g$ is the genus of $X$. if $X=S^2$ then $\beta_0=\beta_2=1$ and $\beta_1=0$. So $\chi(X)=2$ and $g=0$. For the torus $T$, $\beta_0=\beta_2=1$ and $\beta_1=2$. So $\chi(T)=0$ and $g$=1. Generally, $g$ is equal to the number of handles we attach to a sphere. So a donut and a coffee cup both have genus one and Euler number 0, and that is the reason a topologist can't tell them apart. 

\section{Eilenberg Steenrod Axioms}

In their textbook \cite{ES}, Eilenberg and Steenrod took an axiomatic approach to homology theory. They started with a series of axioms that a homology theory needed to satisfy and showed that simplicial complexes and simplicial maps satisfied them. In this section, I will list the axioms. First, though, there are two preliminaries that will help to unsderstand them. The first, is the long exact homology sequence of a pair. A related idea is the Mayer-Vietoris sequence which computes the homology of the union of two complexes using another long exact sequence. Finally, we will state the axioms and see how they relate to simplicial homology. 

\subsection{Long Exact Sequences and Zig-Zagging}

Recall that in Section 3.2, we discussed exact sequences of abelian groups in which the kernel of each map is the image of the previous map. A sequence of 5 groups whose first and last ones  are 0 is called a {\it short exact sequence}. A sequence indexed by the integers is called a {\it long exact sequence}\index{long exact sequence}. An important long exact sequence involves the homology of a pair. I will state the version for simplicial homology. One of the Eilenberg-Steenrod axioms will generalize this. Later, we will see variants for cohomology and homotopy.

\begin{theorem}
{\bf The long exact homology sequence of a pair:} Let $K$ be a complex, $K_0$ a subcomplex. Then there is a long exact sequence $$\begin{tikzpicture}
  \matrix (m) [matrix of math nodes,row sep=3em,column sep=4em,minimum width=2em]
  {
\cdots & H_p(K_0) & H_p(K) & H_p(K, K_0) & H_{p-1}(K_0) & \cdots,\\};
  \path[-stealth]
    (m-1-1)  edge  (m-1-2)
(m-1-2)  edge node [above] {$i_*$} (m-1-3)
(m-1-3)  edge node [above] {$j_*$} (m-1-4)
(m-1-4)  edge node [above] {$\partial_*$} (m-1-5)
    (m-1-5)  edge  (m-1-6);
\end{tikzpicture}$$  where $i: K_0\rightarrow K$ and $j: (K,\emptyset)\rightarrow (K, K_0)$ are inclusions.
\end{theorem}

The third map $\partial_*$ is called the {\it homology boundary homomorphism}. Given a cycle $z$ in $C_p(K,K_0)$, it is in $d+C_p(K_0)$ where $\partial d$ is carried by $K_0$. So if $\{\}$ means homology class then define $\partial_*\{z\}=\{\partial d\}$.

This result will be an immediate consequence of the more general "zig-zag lemma" which I will now describe. It featured in the only movie I have ever seen depicting algebraic topology: a Dutch movie called "Antonia's Line". Her granddaughter, a math professor, was writing it on the board. Definitely the most important part.

\begin{definition}
Let $\frak{C}, \frak{D}$, and $\frak{E}$ be chain complexes. Let 0 be the trivial chain complex whose groups vanish in every dimension. Let $\phi: \frak{C}\rightarrow\frak{D}$ and $\psi: \frak{D}\rightarrow\frak{E}$ be chain maps. (Remember this means they commute with boundaries.) We say the sequence $$\begin{tikzpicture}
  \matrix (m) [matrix of math nodes,row sep=3em,column sep=4em,minimum width=2em]
  {
  0 & \frak{C} & \frak{D} & \frak{E} & 0.\\};
  \path[-stealth]
    (m-1-1) edge  (m-1-2)
    (m-1-2)  edge node [above] {$\phi$} (m-1-3)
(m-1-3)  edge node [above] {$\psi$} (m-1-4)
(m-1-4) edge  (m-1-5);

\end{tikzpicture}$$ is a {\it short exact sequence of chain complexes} if $$\begin{tikzpicture}
  \matrix (m) [matrix of math nodes,row sep=3em,column sep=4em,minimum width=2em]
  {
  0 & C_p & D_p & E_p & 0.\\};
  \path[-stealth]
    (m-1-1) edge  (m-1-2)
    (m-1-2)  edge node [above] {$\phi$} (m-1-3)
(m-1-3)  edge node [above] {$\psi$} (m-1-4)
(m-1-4) edge  (m-1-5);

\end{tikzpicture}$$ is an exact sequence of groups for every $p$.
\end{definition}

In the previous theorem, we had the exact sequence of complexes $$\begin{tikzpicture}
  \matrix (m) [matrix of math nodes,row sep=3em,column sep=4em,minimum width=2em]
  {
  0 & \frak{C}(K_0) & \frak{C}(K) & \frak{C}(K, K_0) & 0.\\};
  \path[-stealth]
    (m-1-1) edge  (m-1-2)
    (m-1-2)  edge node [above] {$i$} (m-1-3)
(m-1-3)  edge node [above] {$\pi$} (m-1-4)
(m-1-4) edge  (m-1-5);

\end{tikzpicture}$$ Here $i$ is inclusion and $\pi$ is projection. The sequence is exact because $C_p(K, K_0)=C_p(K)/C_p(K_0)$.

\begin{theorem}
{\bf Zig-zag Lemma:} If $\frak{C}, \frak{D}$, and $\frak{E}$ are chain complexes such that the sequence $$\begin{tikzpicture}
  \matrix (m) [matrix of math nodes,row sep=3em,column sep=4em,minimum width=2em]
  {
  0 & \frak{C} & \frak{D} & \frak{E} & 0.\\};
  \path[-stealth]
    (m-1-1) edge  (m-1-2)
    (m-1-2)  edge node [above] {$\phi$} (m-1-3)
(m-1-3)  edge node [above] {$\psi$} (m-1-4)
(m-1-4) edge  (m-1-5);

\end{tikzpicture}$$ is exact. Then there is a long exact homology sequence $$\begin{tikzpicture}
  \matrix (m) [matrix of math nodes,row sep=3em,column sep=4em,minimum width=2em]
  {
\cdots & H_p(\frak{C}) & H_p(\frak{D}) & H_p(\frak{E}) & H_{p-1}(\frak{C}) & \cdots,\\};
  \path[-stealth]
    (m-1-1)  edge  (m-1-2)
(m-1-2)  edge node [above] {$\phi_*$} (m-1-3)
(m-1-3)  edge node [above] {$\psi_*$} (m-1-4)
(m-1-4)  edge node [above] {$\partial_*$} (m-1-5)
    (m-1-5)  edge  (m-1-6);
\end{tikzpicture}$$ where $\partial_*$ is induced by the boundary operator in $\frak{D}$. 
\end{theorem}

\begin{figure}[ht]
\begin{center}
  \scalebox{0.4}{\includegraphics{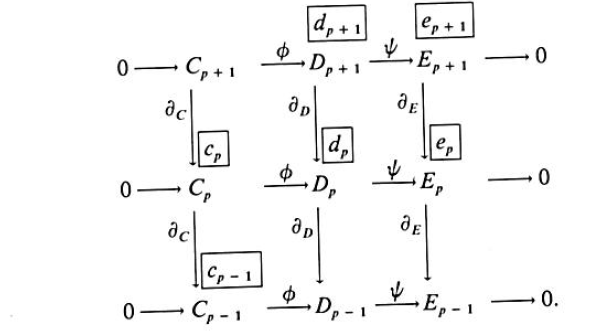}}
\caption{
\rm 
Zig-zag Lemma \cite{Mun1}
}
\end{center}
\end{figure}

The proof is long but not that hard. It is an example of what is called "diagram chasing". Refer to Figure 4.2.1 from Munkres \cite{Mun1}. You can find the full proof there. I will list the steps and write out step 1. The steps are:\begin{enumerate}
\item Define $\partial_*$.
\item Show that $\partial_*$ is a well defined homomorphism, i.e. two cycles in $E_p$ that are homologous are sent to  homologous elements of $C_{p-1}$ by $\partial_p$.
\item Show exactness at $H_p(\frak{D})$.
\item Show exactness at $H_p(\frak{E})$.
\item Show exactness at $H_{p-1}(\frak{C})$.
\end{enumerate}

As an example, to define $\partial_*$: Let $e_p\in E_p$ be a cycle. Since $\psi$ is onto by exactness, choose $d_p\in D_p$ so that $\psi(d_p)=e_p.$ The element $\partial_D(d_p)\in D_{p-1}$ is in $\ker\psi$ since $\psi(\partial_D(d_p))=\partial_E(\psi(d_p))=\partial_E(e_p)=0,$ since we assumed that $e_p$ was a cycle. So there is an element $c_{p-1}\in C_{p-1}$ such that $\phi(c_{p-1})=\partial_D(d_p)$ since $\ker\psi={\rm im }\phi$. This element is unique since $\phi$ is injective. Also $c_{p-1}$ is a cycle since $\phi(\partial_C(c_{p-1}))=\partial_D(\phi(c_{p-1}))=\partial_D\partial_D(d_p)=0.$ So $\partial_C(c_{p-1})=0$ since $\phi$ is injective.So if $\{\}$ denotes homology let $\partial_*\{e_p\}=\{c_{p-1}\}$.

So this gives you a taste of how the rest of the proof goes. Look at Figure 4.2.1 for help.

\subsection{Mayer-Vietoris Sequences}

In this very short section, I will give another useful exact sequence for computing homology.

\begin{theorem}
Let $K$ be a complex and let $K_0$ and $K_1$ be subcomplexes such that $K=K_0\cup K_1$. Let $A=K_0\cap K_1$. Then there is an exact sequence  $$\begin{tikzpicture}
  \matrix (m) [matrix of math nodes,row sep=3em,column sep=4em,minimum width=2em]
  {
\cdots & H_p(A) & H_p(K_0)\bigoplus H_p(K_1) & H_p(K) & H_{p-1}(A) & \cdots,\\};
  \path[-stealth]
    (m-1-1)  edge  (m-1-2)
(m-1-2)  edge  (m-1-3)
(m-1-3)  edge (m-1-4)
(m-1-4)  edge (m-1-5)
    (m-1-5)  edge  (m-1-6);
\end{tikzpicture}$$ called the {\bf Mayer-Vietoris sequence}\index{Mayer-Vietoris sequence} of $(K_0, K_1)$. There is a similar sequence in reduced homology if $A\neq\emptyset.$
\end{theorem}

{\bf Proof:} I will give the proof for ordinary homology. See \cite{Mun1} for the modifications in the reduced homology case.  

We want to construct a short exact sequence of chain complexes 
$$\begin{tikzpicture}
  \matrix (m) [matrix of math nodes,row sep=3em,column sep=4em,minimum width=2em]
  {
  0 & \frak{C}(A) & \frak{C}(K_0)\bigoplus\frak{C}(K_1) & \frak{C}(K) & 0.\\};
  \path[-stealth]
    (m-1-1) edge  (m-1-2)
    (m-1-2)  edge node [above] {$\phi$} (m-1-3)
(m-1-3)  edge node [above] {$\psi$} (m-1-4)
(m-1-4) edge  (m-1-5);

\end{tikzpicture}$$ and apply the zig-zag lemma. 

First, we need to define the chain complex in the middle. It's chain group in dimesnsion $p$ is $C_p(K_0)\bigoplus C_p(K_1)$. The boundary operator is defined by $\partial'(d,e)=(\partial_0d, \partial_1e)$, where $\partial_0$ and $\partial_1$ are the boundary operators for $ \frak{C}(K_0)$ and $ \frak{C}(K_1)$ respectively.

Next, we need to define the chain maps $\phi$ and $\psi$. Consider the commutative diagram where all maps are inclusions: 

$$\begin{tikzpicture}
  \matrix (m) [matrix of math nodes,row sep=3em,column sep=4em,minimum width=2em]
  {
  & K_0 & \\
A & & K.\\
& K_1 & \\};
  \path[-stealth]
    (m-2-1) edge node [above] {$i$} (m-1-2)
(m-2-1) edge node [below] {$j$} (m-3-2)
(m-2-1) edge node [above] {$m$} (m-2-3)
(m-1-2) edge node [above] {$k$} (m-2-3)
(m-3-2) edge node [below] {$l$} (m-2-3)
;

\end{tikzpicture}$$

Define homomomorphisms $\phi(c)=(i_\#(c), -j_\#(c))$ and $\psi(d, e)=k_\#(d)+l_\#(e).$ These are both chain maps. 

To check exactness, $\phi$ is injective since $i_\#$ and $j_\#$ are inclusion of chains. To see that $\psi$ is surjective, let $d\in C_p(K)$, write $d$ as a sum of oriented simplices and let $(d_0)$ be those carried by $K_0$. Then $d-d_0$ is carried by $K_1$ and $\psi(d_0, d-d_0)=d$.

Now check exactness at the middle term. Note that $\psi\phi(c)=m_\#(c)-m_\#(c)=0.$ Conversely, if $\psi(d, e)=0$, then $d=-e$ as chains in $K$. Since $d$ is carried by $K_0$ and $e$ is carried by $K_1$, they must both be carried by $K_0\cap K_1=A$. So $(d,e)=(d, -d)=\phi(d)$ and we are done.

The homology of the middle chain complex in dimension $p$ is $$\frac{\ker\partial'}{{\rm im}\partial'}=\frac{\ker\partial_0\bigoplus\ker\partial_1}{{\rm im} \partial_0\bigoplus{\rm im}\partial_1}\cong H_p(K_0)\bigoplus H_p(K_1).$$ The result now follows from the zig-zag lemma. $\blacksquare$

As an application, recall that the {\it cone of a space} is the result of taking an external point and including all of the lines from that point to a point of the original sapce. For a simplicial complex $K$, if $w\notin K$, then if $v_0,\cdots, v_p$ is a $p$-simplex of $K$, we add the $p+1$-simplex $[w, v_0,\cdots, v_p]$ to the cone. Munkres writes the cone as $w*K$. Cones basically plug up holes and are contractible which means they are also acyclic.

\begin{definition}
Let $K$ be a complex and $w_0*K$ and $w_1*K$ be two cones of $K$ whose polytopes intersect in $|K|$. Then $S(K)=(w_0*K)\cup(w_1*K)$ is called the {\it suspension}\index{suspension} of $K$. 
\end{definition}

I briefly mentioned suspension in Section 2.4 and Figure 2.4.4. The definition here is equivalent to the one we gave there but now stated in terms of simplicial complexes. I also mentioned that suspension raises the dimension of a sphere by one. (Hold two ice cream cones together with their wide ends touching.) Then the homology in dimension $p$ of the suspension of a sphere is equal to the homology in dimension $p-1$ of the original sphere. This is true even for other spaces as we now show.

\begin{theorem}
Let $K$ be a complex. Then for all $p$ there is an isomorphism $$\tilde{H}_p(S(K))\rightarrow\tilde{H}_{p-1}(K).$$
\end{theorem}

{\bf Proof:} Let $K_0=(w_0*K)$ and $K_1=(w_1*K)$. Then $S(K)=K_0\cup K_1$ and $K_0\cap K_1=K$. Using the Mayer-Vietoris sequence, we get $$\tilde{H}_p(K_0)\bigoplus\tilde{H}_p(K_1)\rightarrow\tilde{H}_p(S(K))\rightarrow\tilde{H}_{p-1}(K)\rightarrow\tilde{H}_{p-1}(K_0)\bigoplus\tilde{H}_{p-1}(K_1).$$ Since cones are acyclic, the terms at the end are zero so the middle map s an isomorphism. $\blacksquare$

I will end by mentioning that one interesting fact about the Austrian mathematician Leopold Vietoris was his unusually long lifespan. He died in 2002, just under two months before his 111th birthday.

\subsection{The Axiom List or What is a Homology Theory?}

As I mentioned earlier, in their textbook \cite{ES}, Eilenberg and Steenrod defined a homology theory as one satisfying a particular set of axioms. This implies that there could be multiple homology theories. That is in fact true. I will list the axioms here and spend the rest of the chapter briefly discussing two important alternatives to the simplicial homology theory we have already described. Homology theories work with special pairs of topological spaces and subspaces called {\it admissible pairs}.

\begin{definition}
Let $\mathcal{A}$ be a class of pairs $(X, A)$ of topological spaces such that:\begin{enumerate}
\item If $(X, A)$ belongs to $\mathcal{A}$, then so do $(X, X), (X, \emptyset), (A, A),$ and $(A, \emptyset)$
\item If $(X, A)$ belongs to $\mathcal{A}$, then so does $(X\times I, A\times I)$.
\item There is a one point set $P$ such that $(P, \emptyset)$ is in $\mathcal{A}$.
\end{enumerate}
Then $\mathcal{A}$ is called an {\it admissible class of spaces} for a homology theory.
\end{definition}

We now state the axioms themselves.\index{Eilenberg-Steenrod Axioms for Homology.}

\begin{definition}
If $\mathcal{A}$ is a set of admissible pairs, a homology theory on $\mathcal{A}$ consists of three functions:\begin{enumerate}
\item A function $H_p$ defined for each integer $p$ and each pair $(X, A)$ in $\mathcal{A}$ whose value is an abelian group.
\item A function for each integer $p$ that assigns to a continous map $h: (X, A)\rightarrow (Y, B)$ (recall that this means $h(A)\subset B$), a homomorphism $$h_*: H_p(X, A)\rightarrow H_p(Y, B).$$
\item A function for each integer $p$ that assigns to each pair $(X, A)\in\mathcal{A}$ a homomorphism $$\partial_*: H_p(X, A)\rightarrow H_{p-1}(A),$$ where $A$ denotes the pair $(A, \emptyset)$. 
\end{enumerate}
\end{definition}

These functions satisfy the following axioms where all pairs of spaces are in $\mathcal{A}$:
\begin{itemize}
\item[$-$]{\bf Axiom 1:} If $i$ is the identity, then $i_*$ is the identity.
\item[$-$]{\bf Axiom 2:} $(kh)_*=k_*h_*.$
\item[$-$]{\bf Axiom 3:} If $f: (X, A)\rightarrow (Y, B)$, then the following diagram commutes:$$\begin{tikzpicture}
  \matrix (m) [matrix of math nodes,row sep=3em,column sep=4em,minimum width=2em]
  {
H_p(X, A) & H_p(Y, B) \\
H_{p-1}(A) & H_{p-1}(B) \\};
  \path[-stealth]
    (m-1-1) edge node [above] {$f_*$} (m-1-2)
(m-2-1) edge node [above] {$(f|A)_*$} (m-2-2)
(m-1-1) edge node [right] {$\partial_*$} (m-2-1)
(m-1-2) edge node [right] {$\partial_*$} (m-2-2)
;

\end{tikzpicture}$$
\item[$-$]{\bf Axiom 4:} (Exactness Axiom) The sequence $$\begin{tikzpicture}
  \matrix (m) [matrix of math nodes,row sep=3em,column sep=4em,minimum width=2em]
  {
\cdots & H_p(A) & H_p(X) &  H_p(X, A) & H_{p-1}(A) &\cdots \\};
  \path[-stealth]
    (m-1-1) edge (m-1-2)
(m-1-2) edge node [above] {$i_*$} (m-1-3)
(m-1-3) edge node [above] {$\pi_*$} (m-1-4)
(m-1-4) edge node [above] {$\partial_*$} (m-1-5)
  (m-1-5) edge (m-1-6)
;

\end{tikzpicture}$$ is exact where $i: A\rightarrow X$ and $\pi: X\rightarrow (X,A)$ are inclusion maps.
\item[$-$]{\bf Axiom 5:} (Homotopy Axiom) If $h, k: (X, A)\rightarrow(Y,B)$ are homotopic then $h_*=k_*$.
\item[$-$]{\bf Axiom 6:} (Excision Axiom) Given $(X, A)$, let $U$ be an open subset of $X$ such that $\bar{U}\subset A^\circ$ (i.e. the closure of $U$ is contained in the interior of $A)$. Then if $(X-U, A-U)$ is admissible, then inclusion induces an isomorphism $$H_p(X-U, A-U)\cong H_p(X, A).$$
\item[$-$]{\bf Axiom 7:} (Dimension Axiom) If $P$ is a one point space then $H_p(P)=0$ for $p\neq 0$, and $H_0(P)\cong Z$. If we are using homology with coefficients in an arbitrary abelian group $G$, then $H_0(P)\cong G$.
\item[$-$]{\bf Axiom 8:} (Axiom of Compact Support) If $\alpha\in H_p(X, A)$, there is an admissible pair $(X_0, A_0)$ with $X_0\subset X$ and $A_0\subset A$ compact such that $\alpha$ is in the image of the homomorphism $H_p(X_0, A_0)\rightarrow H_p(X, A)$ induced by inclusion.
\end{itemize}

Note that there are variant homology theories where Axiom 7 may fail and a point may have nonzero homology in dimensions other than zero. Examples of this are $K$-theory and cobordism. Neither of them will be covered in this book. Such a theory is called an {\it extraordinary homology theory}. One book that briefly discusses both cobordism and K-theory is \cite{May1}.

A you might expect, simplicial homology is an actual homology theory, but there are some subtleties in showing this. I will refer you to Munkres \cite{Mun1} for most of the details but we should at least define what type of pair corresponds to an admissible pair in simplicial homology theory. What we need is called a {\it triangulable pair.}

\begin{definition}
Let $A$ be a subspace of $X$. A {\it triangulation} of the pair $(X, A)$ is a complex $K$ and a subcomplex $K_0$ along with a homeomorphism $h: (|K|, |K_0|)\rightarrow (X, A)$. If such a triangulation exists, then $(X, A)$ is a {\it triangulable pair}. If $A$ is empty, then $X$ is a {\it triangulable space}.\index{triangulation}\index{triangulable pair}\index{triangulable space}
\end{definition}

\section{Singular and Cellular Homology}

In this final section of the chapter I will briefly describe two other important homology theories. Singular homology has a much nicer theory, but it is hard to visualize or actually compute anything. Unlike simplicial homology, though, it is very easy to prove topological invariance and the only condition on $(X, A)$ is that $A$ is a subspace of $X$. Cellular homology is a good middle ground. It is more flexible than simplicial homology and simple shapes are easier to construct. Best of all, you don't have to worry about tripping over cords when you are using it. Finally, I will use cellular homology to compute the homology of the {\it projective spaces}. These will turn out to be very important for the construction of Steenrod squares.

\subsection{Singular Homology}

Start with the space $R^\infty$ which has coordinates indexed by the positive integers and has only finitely many non-zero components. Let $\Delta_p$ be the $p$-simplex in $R^\infty$ with vertices $e_0, \cdots, e_p$ such that $e_0$ is the origin and $e_i$ has a one in the $i$-th coordinate and zeros elsewhere. Then $\Delta_p$ is called the {\it standard p-simplex}.\index{standard $p$-simplex} If $X$ is a topological space, then a {\it singular p-simplex}\index{singular $p$-simplex} is a continuous map $T: \Delta_p\rightarrow X$. Note that $T$ has no other restrictions and can even be a constant map. 

The free abelian group generated by all singular $p$-simplices is called the {\it singular chain group}\index{singular chain group} of $X$ in dimension $p$ and denoted $S_p(X)$. You can already tell why this would not be practical for computation. In most cases, this group would be enormous. 

To define faces of a simplex we need a special map called the {\it linear singular simplex}\index{linear singular simplex}. If $a_0, \cdots. a_p$ are points in $R^\infty$ which are not necessarily independent, then we have a map $\frak{l}:\Delta_p\rightarrow R^\infty$ that takes $e_i$ to $a_i$ for $i=0, 1, \cdots, p$. It is defined by $$\frak{l}(x_1, \cdots, x_p, 0, \cdots)=a_0+\sum_{i=1}^p x_i(a_i-a_0).$$ We denote it by $\frak{l}(x_1, \cdots, x_p)$. Then $\frak{l}(e_1, \cdots, e_p)$ is inclusion of $\Delta_p$ into $R^\infty$. If we leave one coordinate out, we have a map $\frak{l}(e_1, \cdots, \hat{e}_i, \cdots, e_p)$ which takes $\Delta_{p-1}$ onto the face $e_0e_1\cdots e_{i-1}e_{i+1}\cdots e_p$ of $\Delta_p$. 

Now the $i$th face of the simplex $T:\Delta_p\rightarrow X$ is $T\circ\frak{l}(e_1, \cdots, \hat{e}_i, \cdots, e_p)$. So we can define a boundary homomorphism $\partial: S_p(X)\rightarrow S_{p-1}(X)$ by $$\partial T=\sum_{i=0}^p (-1)^i T\circ\frak{l}(e_1, \cdots, \hat{e}_i, \cdots, e_p).$$

If $f: X\rightarrow Y$ is a continuous map, then define $f_\#(T): S_p(X)\rightarrow S_p(Y)$ by $f_\#(T)=f\circ T.$ 

\begin{theorem}
The homomorphism $f_\#$ commutes with $\partial$ and $\partial\partial=0.$
\end{theorem}

We now have a chain complex with groups $S_p(X)$ and a boundary homomorphism $\partial$ so we can define homology groups in the usual way. 

We can also define reduced homology by letting the augmentation map $\epsilon: S_0(X)\rightarrow Z$ be defined by $\epsilon(T)=1$ for any 0-simplex $T$. Now if T is a one simplex, then $\epsilon(\partial T)$=$\epsilon(T(\frak{l}(e_0))-T(\frak{l}(e_1)))=1-1=0$. So we can also define reduced homology. 

\begin{theorem}
If $i: X\rightarrow X$ is the identity, then $i_*: H_p(X)\rightarrow H_p(X)$ is the identity. If $f:X\rightarrow Y$ and $g:Y\rightarrow Z$ then $gf_*=g_*f_*$. The same holds in reduced homology.
\end{theorem}

{\bf Proof:} Both equations actually hold for chains since $i_\#(T)=iT=T$ and $(gf)_\#(T)=(gf)T=g(fT)=g_\#(f_\#(T)).$ Now the result follows since everything commutes with $\partial$.

The preceding result immediately implies topological invariance. 

\begin{theorem}
If $h: X\rightarrow Y$ is a homeomorphism, then $h_*$ is an isomorphism.
\end{theorem}

The point is that we have proved topological invariance with much less work than in the simplicial case. 

I will leave out the rest of the details, but we can define relative homology in the same way as before. If $X$ is a topological space and $A$ is a subspace, then $(X, A)$ is automatically an admissible pair, and the Eilenberg-Steenrod Axioms all hold, so singular homology is a true homology theory. Most importantly for a space with the homotopy type of a simplicial complex, its singular and simplicial homology groups are the same. See \cite{Mun1} for all of the details.

\subsection{CW Complexes}

I will finish my discussion of homology with a discussion of cellular homology. In many cases, this makes spaces much easier to describe than in simplicial theory and aids in computations when you don't have the use of a computer. Spaces are represented as {\it CW complexes}. Instead of simplices, CW complexes are made up of {\it cells} which are homeomorphic to open balls. So a 0-cell is a point, a 1-cell is a line segment, a 2-cell is a solid circle, and a 3-cell is a solid ball. Now lets go back to Disneyworld. (Sounds good as I am writing this in January.) Remember the representation of Spaceship Earth as a simplicial complex. There are 11,520 triangles making up this sphere, although it is actually missing some due to supports and doors but we won't worry about that. Imagine if you wanted to compute the simplicial homology by hand. It would be a lot of work and most of them would cancel out anyway, so it would be a big waste of time. We can form $S^2$  as a CW complex with a 2-cell and glue is entire boundary to a 0-cell. Now the complex has only 2 pieces. (Note there are other ways to produce $S^2$ so we need to be a little careful as we will see later.) In this subsection and the next, I will make this definition formal and show you how to compute homology.  

First we need to say what we mean by gluing spaces together.

\begin{definition}
Let $X$ and $Y$ be disjoint topological spaces. Let $A$ be a closed subset of $X$ and let $f:A\rightarrow Y$ be continuous. Then let $X\cup Y$ be the topological sum (i,e, a disjoint union). Form a quotient space by identifying each set $\{y\}\cup f^{-1}(y)$ for $y\in Y$ to a point. The set will consist of these identified points along with te points $\{x\}$ of $X-A$. We denote this quotient space by $X\cup_f Y$\index{$X\cup_f Y$} and call it the {\it adjunction space determined by } $f$.\index{adjunction space}.
\end{definition}

\begin{definition}
A space is a {\it cell}\index{cell} of dimension $m$ or an $m$-{\it cell} if it is homeomorphic to the ball $B^m$. It is an {\it open cell} of dimension $m$ if it is homeomorphic to the interior of $B^m$. 
\end{definition}

\begin{definition}
A {\it CW complex}\index{CW complex} is a space $X$ and a collection of disjoint open cells $e_\alpha$ whose union is $X$ such that \begin{enumerate}
\item $X$ is Hausdorff.
\item For each open $m$-cell $e_\alpha$ of the collection, there exists a continuous map $f_\alpha: B^m\rightarrow X$ that maps the interior of $B^m$ homeomorphically onto $e_\alpha$ and carries the boundary of $B^m$ into a finite union of open cells, each of dimension less than $m$. 
\item A set $A$ is closed in $X$ if $A\cap\overline{e}_\alpha$ is closed in $\overline{e}_\alpha$ for each $\alpha$. 
\end{enumerate}
\end{definition}

CW complexes were discovered by J. H. C. Whitehead, but they are not named after him. The second condition is called "closure finiteness", and the third condition epresses the fact that $X$ has what Whitehead called the "weak topology" with respect to the $\{\overline{e}_\alpha\}$. The letters C and W come from these conditions. 

To match the notation in Munkres, we will use the notation $\dot{e}_\alpha$ for $\overline{e}_\alpha-e_\alpha.$ Then the definition implies that $f_\alpha$ carries $B^m$ onto $\overline{e}_\alpha$, and the boundary of $B^m$ onto $\dot{e}_\alpha$. 

The map $f_\alpha$ is called a {\it characteristic map} for the open cell $e_\alpha$. We use $X$ both for the CW complex and its underlyng space.

A {\it finite CW complex} is one which has finitely many cells. Such a complex is always compact. Also, in homotopy theory, a lot of ugliness is avoided for spaces having the homotopy type of a finite CW complex. Any space that we will ever see in data science will be of this type.

\begin{theorem}
Let $X$ be a CW complex with open cells $e_\alpha$. A function $f: X\rightarrow Y$ is continuous if and only if $f|\overline{e}_\alpha$ is continuous for each $\alpha$.  A function $F: X\times I\rightarrow Y$ is continuous if and only if $f|(\overline{e}_\alpha\times I)$ is continuous for each $\alpha$. 
\end{theorem}

Now we give some examples. 

\begin{example}
Recall that a torus can be created by folding a rectangle. We can represent a Torus as a CW complex consisting of one 2-cell (the image of the interior of the rectangle), two 1-cells (the images of the two short ends glued together and the two long ends glued together), and one 0-cell (the image of the four vertices all glued together). The short ends become the circle through the inside of the torus and the two long ends form the circle around the circumference of the torus. 
\end{example}

A CW complex has a notion of a $p$-skeleton analogous to that of a simplicial complex. 

\begin{definition}
If $X$ is a CW complex, then the $p$-{\it skeleton} $X^p$ is the subspace of $X$ consisting of open cells of dimension at most $p$.The {\it dimension} of $X$ is the dimension of the highest dimesnional cell in $X$.
\end{definition}

\begin{example}
We saw earlier that we can construct $S^2$ as a CW complex consisting of a circle (a 2-cell) whose boundary is glued to a single point (a 0-cell). So this complex has 2 cells. Another representation of $S^2$ has four cells. Take a point and attach a one-cell to form a circle. Now letting the circle be the equator, attach two 2-cells for the Northern and Southern hemispheres respectively. This complex is also of homotopy type $S^2$. Let $A$ be the first type of complex and $B$ the second type. Then $A$ and $B$ are homotopy equivalent. But their one skeletons are not homotopy equivalent. $A^1$ consists of a single point, while $B^1$ has the homotopy type of the circle, $S^1$.
\end{example}

As promised, I will now give an example of a CW complex which is non-triangulable. Recall that for a triangulable space $X$, there is a simplicial complex $K$ and a homeomorphism $h: |K|\rightarrow X$. The example as Munkres presents it is a little confusing. What I thinks he means is that he constructs a finite CW complex which can not be triangulated by a {\bf finite} simplicial complex. Assuming that is the case, I can't find an example of a space that can't be triangulated at all. We need one more definition:

\begin{definition}
If a CW complex $X$ is triangulated by a complex $K$ in such a way that each skeleton $X^p$ of $X$ is triangulated by a subcomplex of $K$ of dimension at most $p$, then we say that $X$ is a {\it triangulable CW complex.}
\end{definition}

Now consider the weird space pictured in Figure 4.3.1. I wouldn't try to build this one at home,

\begin{figure}[ht]
\begin{center}
  \scalebox{0.4}{\includegraphics{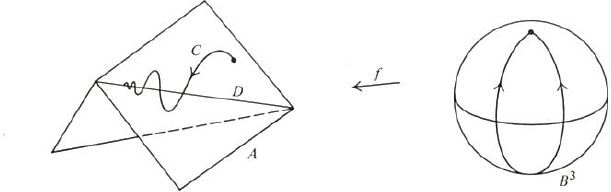}}
\caption{
\rm 
Non-triangulable CW Complex \cite{Mun1}
}
\end{center}
\end{figure}

$A$ is the subspace of $R^3$ consisting of a square and a triangle one of whose edges coincides with the diagonal $D$ of the square. So $A$ is the space of a complex consisting of 3 triangles with an edge in common. Now let $C$ be a 1-cell in the square intersecting the diagonal in an infinite disconnected set. (Letting $D$ be a piece of the positive x axis starting from the origin and let $C$ be the graph of the function $y=x\sin(1/x)$ will accomplish this.) Now take a ball $B^3$ and attach it to the curve $C$ by gluing every arc from the South pole to the North pole of $B^3$ to the curve $C$. The resulting adjunction space $X$ is a finite CW complex with cells the open simplices of $A$ in dimensions 0, 1, and 2, and the 3-cell $e_3$ which is the interior of $B^3$.

Suppose $h: |K|\rightarrow X$ is a triangulation, where $K$ is a finite simplicial complex. We can write $X$ as the disjoint union $$X=(A-C)\cup C\cup e_3.$$ Munkres uses arguments from the theory of {\it local homology groups} (see \cite{Mun1} section 35) to argue that for a simplex $\sigma$ in $K$, the image of the interior of $\sigma$ can not intersect both $(A-C)$ and $e_3$, so $h(\sigma)$ lies in $A$ or $\overline{e_3}$.  So if $h$ triangulates $A$ and $\overline{e_3}$, it must triangulate $A\cap\overline{e_3}=C$. A similar argument shows that $D$ is triangulated by $h$ and so is $C\cap D$. But the latter is an infinite set of disconnected points and can't be triangulated by a finite complex, so we have a contradiction. 

In the above I left out some details that won't show up in later sections. See \cite{Mun1} for the details. Fortunately, there is no data science application I can imagine where we will ever see a space like the one just described. We all always assume we can represent our spaces as finite simplicial complexes. 

I will conclude this subsection by stating two helpful results that will help you picture where CW complexes come from. Also, recall from Chapter 2 that a topological space $X$ is normal if whenever $A$ and $B$ are disjoint closed sets in $X$, there are disjoint open subsets $U$ and $V$ of $X$ such that $A\subset U$ and $B\subset V$. In his discussion of adjucntion spaces, Munkres shows that the adjunction of two normal spaces is normal and the normality is closed under {\it coherent union} in which the intersection of a closed subset of the union of a collection of spaces and each of these spaces individually is closed. The two theorems below then imply that CW complexes are always normal.

\begin{theorem}
\begin{enumerate}
\item Let $X$ be a CW complex of dimension $p$. Then $X$ is homeomorphic to an adjunction space formed from $X^{p-1}$ and a topological sum $\sum B_\alpha$ of closed $p$-balls by means of a continuous map $g: \sum Bd B_\alpha\rightarrow X^{p-1}$. 
\item If Y is a CW complex of dimension at most $p-1$, $\sum B_\alpha$ is a topological sum of closed $p$-balls, and if $g:\sum  Bd B_\alpha\rightarrow Y$ is a continuous map, then the adjunction space $X$ formed from $Y$ and $\sum B_\alpha$ by means of $g$ is a CW complex and Y is its $p-1$-skeleton.
\end{enumerate}
\end{theorem}

\begin{theorem}
\begin{enumerate}
\item Let $X$ be a CW complex. Then $X^p$ is a closed subspace of $X^{p+1}$ for each $p$, and $X$ is the coherent union of the spaces $X^0\subset X^1\subset X^2\subset\cdots$. 
\item Suppose that $X_p$ is a CW-complex for each $p$ and $X_p$ is the $p$-skeleton of $X_{p+1}$ for each $p$. If $X$ is the coheret union of the spaces $X_p$, then $X$ is a CW complex having $X_p$ as its $p$-skeleton.
\end{enumerate}
\end{theorem}

\subsection{Homology of CW Complexes}

It now remains to see how to compute the homology of a CW complex. To do this, we have to start with a definition of chains and boundaries. It might seem logical to have chains be formal sums of cells analogous to how simplicial chains were defined. The problem is that we can't really define boundaries in the same way. We will define cellular chains a little differently and immediately show that they are the same thing for the undelrying space of a simplicial complex. 

Let $X$ denote a CW complex with open cells $e_\alpha$ and characteristic maps $f_\alpha$. The symbol $H_p$ can designate singular homology, but since we will never worry about any type of space other than a triangulable finite CW complex, we can just as easily let it denote simplicial homology. Recall the exact homology sequence of a pair:  $$\begin{tikzpicture}
  \matrix (m) [matrix of math nodes,row sep=3em,column sep=4em,minimum width=2em]
  {
\cdots & H_p(A) & H_p(X) &  H_p(X, A) & H_{p-1}(A) &\cdots \\};
  \path[-stealth]
    (m-1-1) edge (m-1-2)
(m-1-2) edge node [above] {$i_*$} (m-1-3)
(m-1-3) edge node [above] {$j_*$} (m-1-4)
(m-1-4) edge node [above] {$\partial_*$} (m-1-5)
  (m-1-5) edge (m-1-6)
;

\end{tikzpicture}$$ where $i: A\rightarrow X$ and $j: X\rightarrow (X,A)$ are inclusion maps.

\begin{definition}
If $X$ is a CW complex, let $$D_p(X)=H_p(X^p, X^{p-1}).$$ Let $\partial: D_p(x)\rightarrow D_{p-1}(X)$ be the composite $$\begin{tikzpicture}
  \matrix (m) [matrix of math nodes,row sep=3em,column sep=4em,minimum width=2em]
  {
 H_p(X^p, X^{p-1}) & H_{p-1}(X^{p-1}) & H_{p-1}(X^{p-1}, X^{p-2}),\\};
  \path[-stealth]

(m-1-1) edge node [above] {$\partial_*$} (m-1-2)
(m-1-2) edge node [above] {$j_*$} (m-1-3)

;

\end{tikzpicture}$$ where $j$ is inclusion. The fact that $\partial^2=0$ follows from the fact that $$\begin{tikzpicture}
  \matrix (m) [matrix of math nodes,row sep=3em,column sep=4em,minimum width=2em]
  {
 H_{p-1}(X^{p-1}) &  H_{p-1}(X^{p-1}, X^{p-2}) &  H_{p-2}(X^{p-2})\\};
  \path[-stealth]

(m-1-1) edge node [above] {$j_*$} (m-1-2)
(m-1-2) edge node [above] {$\partial_*$} (m-1-3)

;

\end{tikzpicture}$$ is exact as it is part of the long exact sequence of the pair $(X^{p-1}, X^{p-2})$. The chain complex $\frak{D}(X)=\{D_p(x), \partial\}$ is called the {\it cellular chain complex}\index{cellular chain complex} of $X$.
\end{definition}

\begin{example}
To help picture this definition, we would like to see how it compares to the simplicial chain complex when the CW complex $X$ is the underlying space of a simplicial complex $K$. Letting $X^p$ be the cellular $p$-skeleton of $X$ and $K^{(p)}$ be the simplicial $p$-skeleton of $K$, assume that the open cells of $X$ coincide with the open simplices of $K$. We want to compute $H_p(X^p, X^{p-1})$. The simplicial chain group $C_i(K^{(p)}, K^{(p-1)})$ is equal to 0 for $i\neq p$ and $C_p(K^{(p)})=C_p(K)$. Thus, $$H_p(X^p, X^{p-1})=H_p(K^{(p)}, K^{(p-1)})=C_p(K).$$ Also, the boundary operator is the ordinary simplicial boundary. This shows that in this case we can use the CW complex to compute the homology of $X$.
\end{example}

Munkres next shows that the cellular chain complex of a CW complex $X$ behaves in a lot of ways like a simplicial complex. He shows that $D_p(X)$ is free abelian with a basis consisting of the {\it oriented} $p$-cells. (I will define what we mean by oriented cells next.) Also, he shows that $\frak{D}(X)$ can be used to compute the singular and thus the simplicial homology of $X$. Rather than reproduce his rather long proof, I will refer you to \cite{Mun1}, Section 39.

To understand some examples, we will introduce some new definitions. For each open $p$-cell $e_\alpha$ of the CW complex $X$, the group $H_p(\overline{e}_\alpha, \dot{e}_\alpha)$ is infinite cyclic (i.e. isomorphic to the group $Z$ of integers.) This is the case since we are taking a closed ball and identifying its boundary to a point. For example, $e_2$ is a circle and identifying the boundary of a closed circle to a point produces $S^2$. We know that $H_2(S^2)=Z$. Now the group $Z$ has two generators, $1$ and $-1$. These generators will be called the two {\it orientations} of the cell $e_\alpha$. The cell $e_p$ together with an orientation is called an {\it oriented p-cell}\index{oriented p-cell}. 

Assuming that $X$ is triangulable as we always can in data science applications, let $K$ be the complex that triangulates $X$. Since $X^p$ and $X^{p-1}$ are subcomplexes of $K$, any open $p$-cell $e_\alpha$ is a union of open simplices of $K$, so $\overline{e}_\alpha$ is  is the polytope of a subcomplex of $K$. The group $H_p(\overline{e}_\alpha, \dot{e}_\alpha)$ is the group of $p$-chains carried by  
$\overline{e}_\alpha$ whose boundaries are carried by $\dot{e}_\alpha$. This group is $Z$ and either generator of the group is called a {\it fundamental cycle}\index{fundamental cycle} for $(\overline{e}_\alpha, \dot{e}_\alpha)$.

The cellular chain group $D_p(X)$ is the group of all simplicial $p$-chains of $X$ carried by $X^p$ whose boundaries are carried by $X^{p-1}$. Any such $p$-chain can be written uniquely as a finite linear combination of fundamental cycles for those pairs $(\overline{e}_\alpha, \dot{e}_\alpha)$ for whiich dim $e_\alpha=p.$

\begin{figure}[ht]
\begin{center}
  \scalebox{0.4}{\includegraphics{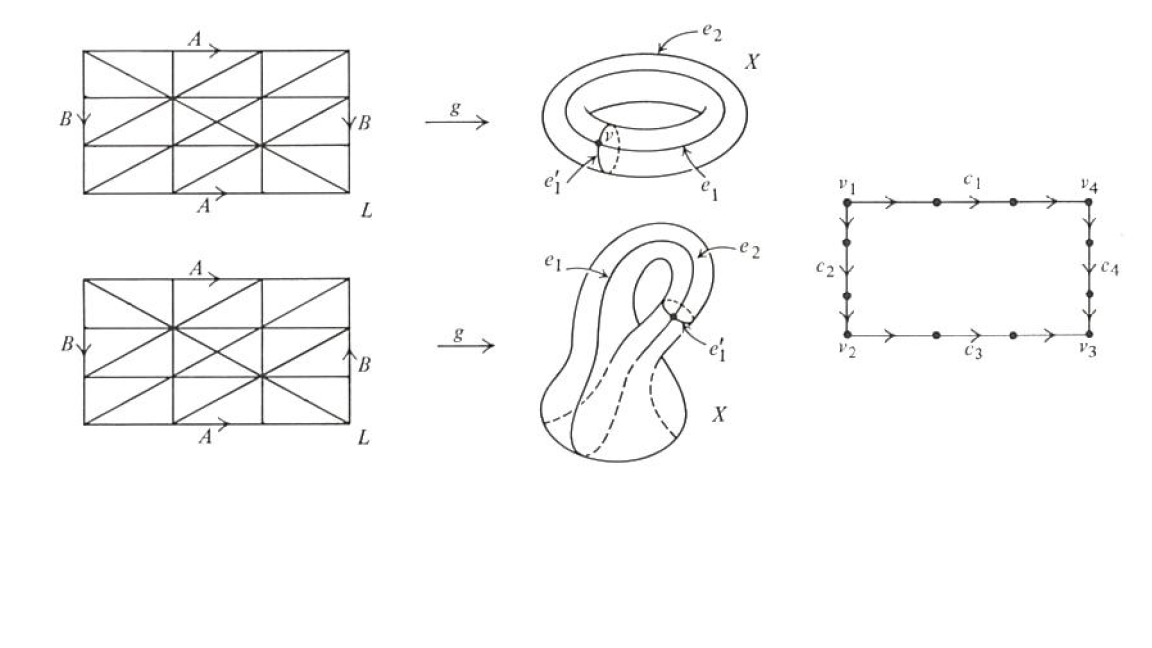}}
\caption{
\rm 
Torus and Klein Bottle as CW Complexes. \cite{Mun1}
}
\end{center}
\end{figure}

\begin{example}
Now we will compute the homology of the torus and the  Klein bottle. You can follow along on Figure 4.3.2.

We let $X$ be the torus or Klein bottle produced by folding the rectangle $L$ as shown in the figure. In either case, $X$ is a triangulable CW complex having an open 2-cell $e_2$, two open 1-cells $e_1$ and $e'_1$ which are the images of $A$ and $B$ respectively, and one 0-cell $e_0$. So $D_2(X)\cong Z$, $D_1(X)\cong Z\oplus Z$, and $D_0(X)\cong Z$. To compute homology, we need generators for these chain groups. Letting $d$ be the the sum of all of the 2-simplices of $L$ oriented counterclockwise, we get a cycle of $(L, Bd L)$. In fact, $d$ is a fundamental cycle for $(L, Bd L)$. By Munkres' arguments about the equivalence of cellular and simplicial homology, $\gamma=g_\#(d)$ where $g$ is the identification map as shown in the figure, is a fundamental cycle for $(\overline{e}_2, \dot{e}_2)$

Now look at the rectangle at the right side of the figure. $w_1=g_\#(c_1)$ is a fundamental cycle for $(\overline{e}_1, \dot{e}_1)$ as is $g_\#(c_3)$. Also, $z_1=g_\#(c_2)$ is a fundamental cycle for $(\overline{e}'_1, \dot{e}'_1)$ as is $g_\#(c_4)$. In the complex $L$, we have $$\partial c_1=v_4-v_1, \hspace{.5 in} \partial c_2=v_2-v_1,$$ $$\partial d=-c_1+c_2+c_3-c_4.$$

The difference between the torus and the Klein bottle is how $L$ is folded, so we take it into account when we apply $g$. Now $\partial w_1=g_\#(\partial c_1)=0$ and $\partial z_1=g_\#(\partial c_2)=0$ for both the torus ad the Klein bottle. In the case of the torus, $\partial(\gamma)=g_\#(\partial d)=0$ since $g_\#(c_1)=g_\#(c_3)$ and $g_\#(c_2)=g_\#(c_4)$. In the case of the Klein bottle, $\partial(\gamma)=g_\#(\partial d)=2g_\#(c_2)=2z_1$,  since $g_\#(c_1)=g_\#(c_3)$ and $g_\#(c_2)=-g_\#(c_4)$. 

Now it is easy to compute homology. If $X$ is the torus, $\gamma$ is a 2-cycle and there are no 2-boundaries, so $H_2(X)\cong Z$. The 1-chains of $X$ are generated by $w_1$ and $z_1$. Both of these are 1-cycles and there are no 1-boundaries so $H_1(X)\cong Z\oplus Z$. Since the torus is connected, $H_0(X)\cong Z$. 

Letting $X$ be the Klein bottle, $\gamma$ generates the 2-chains and it is not a cycle, so $H_2(X)=0$. The 1-chains of $X$ are generated by $w_1$ and $z_1$. Both of these are 1-cycles and the boundaries are generated by $\partial\gamma=2z_1$, so  $H_1(X)\cong Z\oplus Z_2$. Since the Klein bottle is connected, $H_0(X)\cong Z$. 

These computations agree with what was stated for simplicial homology.
\end{example}

\begin{example}
Let $S^n$ be an $n$-sphere with $n>1$. $S^n$ has the homotopy type of a CW complex with one $n$-cell and one 0-cell. Then the cellular chain complex is infinite cyclic in dimensions $n$ and 0 and 0 otherwise. So $H_i(S^n)\cong Z$ if $i=0, n$ and $H_i(S^n)=0$ otherwise. For $i=1$, the one cell is a cycle as its endpoints are identified. So $H_i(S^1)\cong Z$ if $i=0, 1$ and $H_i(S^1)=0$ otherwise.
\end{example}

\subsection{Projective Spaces}

Munkres \cite{Mun1} includes this section as an application of cellular homology. I am interested in it both for that reason and for its role in the construction of Steenrod squares which I will discuss in Chapter 11. This section comes entirely from \cite{Mun1}.

\begin{definition}
Starting with the $n$-sphere $S^n$, identify all pairs of antipodal points $x$ and $-x$. The resulting quotient space is called the {\it real projective n-space} and denoted $P^n$.\index{real projective space}\index{$P^n$} 
\end{definition}

Represent $R^n$ as the set of sequences of real numbers $\{x_1, x_2, \cdots\}$ such that $x_i=0$ for $i>n$. Then we have $R^n\subset R^{n+1}$. Letting $S^n$ be the set  $\{x_1, x_2, \cdots\}\subset R^{n+1}$ such that $x_1^2+x_2^2+\cdots+x_{n+1}^2=1$, we see that $S^{n-1}\subset S^n$. Now $S^{n-1}$ is just the equator of $S^n$ and if $x\in S^{n-1}$ then $-x\in S^{n-1}$. So we have that $P^{n-1}\subset P^n$. In fact, $P^{n-1}$ is a closed subspace of $P^n$.

\begin{theorem}
The space $P^n$ is a CW complex having one cell in each dimension $0\leq j\leq n$. Its $j$-skeleton is $P^j$.
\end{theorem}

The idea of the proof is that we proceed by induction on $n$. $P^0$ identifies the two points of $S^0$, so it consists of a single point. Supposing the theorem is true for $P^{n-1}$. Any pair of points that are off of the equator of $S^n$ have at least one point in the open Northern hemisphere $E_+^n$. Now  $E_+^n$ is an open $n$-cell and the quotient map takes its boundary $S^{n-1}$ onto $P^{n-1}$ which is assumed to be a CW-complex with one cell in each dimension up to $n-1$. So the quotient map restricted to $E_+^n$ becomes a characteristic map for an $n$-cell. and we are done.

Note that $P^1$ has a 0-cell and an open 1-cell so it is homeomorphic to $S^1$.

\begin{definition}
$P^0\subset P^1\subset P^2\cdots$ is an increasing sequence of spaces. Let $P^\infty$ be the {\it coherent union}\index{coherent union}\index{$P^\infty$} of these spaces. (A {\it coherent union} of a sequence of spaces $\{X_\alpha\}$ is the union of these spaces in which a subset $A$ is defined to be open if $A\cap X_\alpha$ is open for every $\alpha$.) $P^\infty$ is called infinite dimensional projective space. By the previous theorem, it is a CW complex with one $j$-cell for all $j\geq 0$. Its $n$-skeleton is $P^n$.
\end{definition}

We can also construct an analogous {\it complex projective space}. Start the way we did with real numbers defining the $n$-dimensional complex vector space $C^n$ to be the space of sequences of complex numbers  $\{z_1, z_2, \cdots\}$ such that $z_i=0$ for $i>n$. Then we have $C^n\subset C^{n+1}$. There is a homeomorphism $\rho: C^{n+1}\rightarrow R^{2n+2}$ called the {\it realification operator} defined by $\rho(z_1, z_2, \cdots)=(Re(z_1), Im(z_1), Re(z_2), Im(z_2), \cdots),$ where $Re(z_j)$ and $Im(z_j)$ are the real and imaginary parts of $z_j$ respectively. Let $\overline{z}_j$ be the complex conjugate of $z_j$ and $|z_j|$ be the norm of $z_j$. If $z_j=a+bi$ then $\overline{z}_j=a-bi$. $|z_j|=(a^2+b^2) =(a+bi)(a-bi)=z_j\overline{z}_j$. So for $z=(z_1, z_2,\cdots, z_{n+1})$, $$|z|=\sqrt{\sum_{j=1}^{n+1}z_j\overline{z}_j}.$$ The subspace of $C^{n+1}$ consisting of all points of norm 1 is called the complex $n$-sphere. It corresponds to $S^{2n+1}$ under the operator $\rho$  so we will denote it by $S^{2n+1}$ both as a subset of $C^{n+1}$ and of $R^{2n+2}$. 

\begin{definition}
In the complex $n$-sphere $S^{2n+1}$, define an equivalence relation $\sim$  by $$(z_1, \cdots, z_{n+1}, 0, \cdots)\sim(\lambda z_1, \cdots, \lambda z_{n+1}, 0, \cdots),$$ where $\lambda$ is a complex number such that $|\lambda|=1$. The resulting quotient space is called {\it complex projective n-space}\index{complex projective $n$-space} and is denoted $CP^n$.\index{$CP^n$}
\end{definition}

Analogous to the real case, $C^n\subset C^{n+1}$ for all $n$, so $S^{2n-1}\subset S^{2n+1}$. Passing to quotient spaces we get $CP^{n-1}\subset CP^n$. 

\begin{theorem}
The space $CP^n$ is a CW complex having one cell in each dimension $2j$ for $0\leq j\leq n$. Its dimension is $2n$ and its $2j$-skeleton is $CP^j$.
\end{theorem}

The proof is similar to the proof in the real case but now $CP^n-CP^{n-1}$ is an open $2n$-cell. Note that I have glossed over some of the point set topology that shows that $P^n$ and $CP^n$ satisfy all of the requirements for a CW complex. See Munkres \cite{Mun1} for more details. 

One useful fact is that $CP^0$ is a single point and $CP^1$ is formed by attaching an open 2-cell. This makes $CP^1$ homeomorphic to $S^2$.

\begin{definition}
Since $CP^0\subset CP^1\subset CP^2\cdots$ is an increasing sequence of spaces, let $CP^\infty$ be the coherent union of these spaces. $CP^\infty$ is called infinite dimensional complex projective space. By the previous theorem, it is a CW complex with one cell in each non-negative even dimension.
\end{definition}

Now we will compute the homology of these spaces. The easier case is the complex projective spaces. Since there is one cell in every even dimension, $D_i(CP^n)\cong Z$ if $i$ is even and  $D_i(CP^n)=0$ otherwise. So in even dimensions every chain is a cycle and no chain bounds. This gives the following result:

\begin{theorem}
$H_i(CP^n)\cong Z$ if $i$ is even and $0\leq i\leq 2n$, and $H_i(CP^n)=0$ otherwise. The group $H_i(CP^\infty)\cong Z$ if $i$ is even and $i\geq 0$, and $H_i(CP^\infty)=0$ otherwise.
\end{theorem}

The case for $P^n$ is harder. We know that for $P^n$, the cellular chain groups are infinite cyclic for $0\leq k\leq n$.  So we need to compute the boundary operators. Since the open k-cell $e_k$ of $P^n$ equals $P^k-P^{k-1}$ and $\dot{e}_k=P^{k-1}$, $D_k(P^n)=H_k(P^k, P^{k-1})$. So we need to compute the boundary $\partial_*: H_{k+1}(P^{k+1}, P^k)\rightarrow H_k(P^k, P^{k-1}).$

\begin{theorem}
Let $p: S^n\rightarrow P^n$ be the quotient map where $n\geq 1$, Let $j: P^n\rightarrow (P^n, P^{n-1})$ be inclusion. The composite homomorphism $$\begin{tikzpicture}
  \matrix (m) [matrix of math nodes,row sep=3em,column sep=4em,minimum width=2em]
  {
 H_n(S^n) &  H_n(P^n) &  H_n(P^n, P^{n-1})\\};
  \path[-stealth]

(m-1-1) edge node [above] {$p_*$} (m-1-2)
(m-1-2) edge node [above] {$j_*$} (m-1-3)

;

\end{tikzpicture}$$ is zero if $n$ is even and multiplication by 2 if $n$ is odd.
\end{theorem}

{\bf Proof:} Munkres provides a proof on the chain level and one on the homology level. I will give the simpler chain level proof and refer you to \cite{Mun1} for the homology level proof. 

A theorem of \cite{Mun1} says that we can triangulate $S^n$ so that the antipodal map $a: S^n\rightarrow S^n$ is simplicial and $P^n$ can be triangulated so that $p: S^n\rightarrow P^n$ is simplicial. Then we can use simplicial homology. Let $c_n$ be a fundamental cycle for $(E_+^n, S^{n-1})$. Then $p_\#(c_n)$ is  a fundamental cycle for $(P^n, P^{n-1})$. (See \cite{Mun1} Theorem 39.1.)

Let $\gamma_n=c_n+(-1)^{n-1}a_\#(c_n)$ be a chain of $S^n$. The antipodal map $a: S^{n-1}\rightarrow S^{n-1}$ has degree $(-1)^n$, so $a_\#$ equals multiplication by $(-1)^n$ on $C_{n-1}(S^{n-1})$. So $$\partial\gamma_n=\partial c_n+(-1)^{n-1}a_\#(\partial c_n)=\partial c_n+(-1)^{2n-1}\partial c_n=\partial c_n-\partial c_n=0,$$ so $\gamma_n$ is a cycle. It is a multiple of no other cycle of $S^n$ since its restriction to $E_+^n$ is $c_n$ which is a fundamental cycle for $(E_+^n, S^{n-1})$. So $\gamma_n$ is a fundamental cycle for $S^n$. 

Finally, $p_\#(\gamma_n)=p_\#(c_n+(-1)^{n-1}a_\#(c_n))$. Since $p\circ a=p$, $p_\#(\gamma_n)=[1+(-1)^{n-1}]p_\#(c_n)$. This is 0 for $n$ even and multiplication by 2 for $n$ odd. The theorem then follows from the fact that $\gamma_n$ is a fundamental cycle for $S^n$ and $p_\#(c_n)$ is a fundamental cycle for $(P^n, P^{n-1})$. $\blacksquare$

Now we can compute the cellular boundary maps.

\begin{theorem}
The map $\partial_*: H_{n+1}(P^{n+1}, P^n)\rightarrow H_n(P^n, P^{n-1})$ is zero if $n$ is even and muliplies a generator by 2 if $n$ is odd.
\end{theorem}

{\bf Proof:} The map $p': (E_+^{n+1}, S^n)\rightarrow (P^{n+1}, P^n)$ is a characteristic map for the open $n+1$-cell of $P^{n+1}$ so it induces an isomorphicm in homology. Consider the following commutative diagram:

$$\begin{tikzpicture}
  \matrix (m) [matrix of math nodes,row sep=3em,column sep=4em,minimum width=2em]
  {
H_{n+1}(E_+^{n+1}, S^n) & H_n(S^n)&\\
H_{n+1}(P^{n+1}, P^n) &  H_n(P^n) &  H_n(P^n, P^{n-1})\\};
  \path[-stealth]

(m-1-1) edge node [above] {$\partial_*$} (m-1-2)
(m-1-1) edge node [below] {$\cong$} (m-1-2)
(m-1-1) edge node [left] {$\cong$} (m-2-1)
(m-1-1) edge node [right] {$p'^*$} (m-2-1)
(m-1-2) edge node [right] {$p^*$} (m-2-2)
(m-2-1) edge node [above] {$\partial_*$} (m-2-2)
(m-2-2) edge node [above] {$j_*$} (m-2-3)

;

\end{tikzpicture}$$ 

The map $\partial_*$ at the top is an isomorphism by the long exact reduced homology sequence of $(E_+^{n+1}, S^n)$, since $E_+^{n+1}$ is acyclic. By the preceding theorem, $(j\circ p)_*$ is zero if $n$ is even and multiplication by 2 if $n$ is odd, so we are done. $\blacksquare$

So all chains in odd dimensions are cycles and even multiples of generators bound, while in even dimensions, there are no cycles. So we have the following:

\begin{theorem}
The homology of projective space is as follows:
$$\tilde{H}_i(P^{2n+1})\cong
\left
\{
\begin{array}{ll}
Z/2 & \mbox{if }i \mbox{ is odd and } 0<i<2n+1,\\
Z & \mbox{if }i=2n+1,\\
0 & \mbox{otherwise.}
\end{array}
\right.
$$
$$\tilde{H}_i(P^{2n})\cong
\left
\{
\begin{array}{ll}
Z/2 & \mbox{if }i \mbox{ is odd and } 0<i<2n,\\
0 & \mbox{otherwise.}
\end{array}
\right.
$$

$$\tilde{H}_i(P^\infty)\cong
\left
\{
\begin{array}{ll}
Z/2 & \mbox{if }i \mbox{ is odd and } 0<i,\\
0 & \mbox{otherwise.}
\end{array}
\right.
$$

\end{theorem}

We are now ready to apply what we have learned to topological data analysis.

\chapter{Persistent Homology}

Recall from my introduction to the book that we want to get an idea of the {\it shape} of our data. Recall Figure 1.0.2 reproduced here as Figure 5.0.1:

\begin{figure}[ht]
\begin{center}
  \scalebox{0.4}{\includegraphics{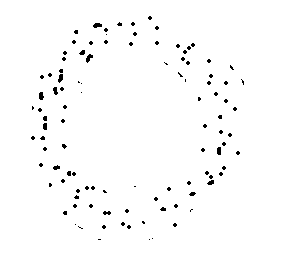}}
\caption{
\rm 
Ring of Data Points
}
\end{center}
\end{figure}

From the point of view of algebraic topology, this figure is pretty boring. It consists of N disconnected points. The 0-th Betti number is $N$ and all of the other homology grups are zero. But your eyes think this should be a circle. So how do we know? And how do we convince a computer of that fact? 

The answer is that we grow the points into circles of increasing radius.As we do, the circles will start to overlap and create holes. As they grow, the holes will be plugged up. There will be a lot of bubbles along the circumference and a big hole in the middle. Keeping track of when the holes form (ie at what radius) and when they are plugged up, we can fully describe the shape. The big hole is the most important feature and it {\it persists} the longest, so we decide it is the most important feature. In the language of the last chapter, each component is a 0-homology class and each hole is a 1-homology class. Once we find these classes, we can visualize them. The three standard methods are {\it bar codes, persistence diagrams}, and {\it persistence landscapes}. These three visualizations display the same information, the birth and death times for each homology class. 

In general, we will replace the points with balls of growing radius To compute homology at each stage, we will need a simplicial complex. Section 5.1 shows how to do this with a point cloud and formally defines persistent homology. As the balls grow, we get a simpliciall complex which is growing as simplices are added. The series of complexes is called a {\it filtration}. Section 5.2 discusses computation issues and explains why persistent homology is generally done with coefficients in a field as opposed to the integral coefficients used in Chapter 4. In Section 5.3, we will see the visualizations described above. 

To do machine learning tasks such as clustering, we would like to have some notion of distance between persistence diagrams. The standard measures are {\it bottleneck distance} and {\it Wasserstein distance}. Persistence landscapes have their own distance measure which I used in my own time series work. All of these measures will be described in Section 5.4. 

As I mentioned, the simplicial complexes in persistent homology grow monitonically forming a filtration. Can we still do something if it this is not the case? Section 5.5 will deal with this case which is called {\it zig-zag persistence}. 

 There are other ways in which persistence diagrams can arise besides point clouds. One example is sublevel sets of a {\it height} function assigned to a data set. An important example of this is grayscale images. These are treated in Section 5.6. Graphs can also lead to persistence diagrams if they are weighted. I will show how this works in Section 5.7. In addition, I will show how several different simplicial complexes can be produced from graphs and describe a quick and easy distance measure between graphs. 

In Section 5.8, I will talk about time series and how they relate to TDA. A time series is a set of equally spaced measurements over time. In Section 5.9, I will describe the SAX method of converting a time series into a string of symbols. This allows for an easy method to detect anomalous substrings as well as classify time series using substrings that are common in one class and rare in the other. Using the a method developed by Gidea and Katz \cite{GK} for financial applications, TDA can be used to convert a multivariate time series into a univariate one. I will conclude the section with my own use of these techniques to classify internet of things devices.

Is there a good general method of extracting features from persistence diagrams? I will describe two very promising techniques. The first is the technique of {\it persistence images} developed by Adams, et. al.\cite{AEKNPSCHMZ}. Another technique is to use {\it template functions}. These are a family of real valued functions highlighting different regions of a diagram. It was developed by Perea, Munch, and Khasawneh \cite{PMK}. I will conclude with a list of software which can be used for TDA applications.

Here are some references on TDA. A good survey paper is \cite{PXL}. I can also recommend four books. Edelsbrunner and Harer \cite{EH} was written by two of the founders of the subject. Other books include \cite{Gh1} and \cite{Zo1}. The recent book \cite{RB} is concerned specifically with applications in biology with a heavy emphasis on virus evolution and cancer genomics. It has a great condensed introduction to persistent homology, and I will rely on it pretty heavily in this chapter.

\section{Definition of Persistent Homology}

We will start with turning a set of discrete points into a simplicial complex. The idea is to consider that these points are just a sample of points from the true shape. We would like to somehow reconstruct the actual shape and then take the homology. But we can't really be sure of the true scale. Is the shape folded up? Which holes are actually there and which are just artifacts of our sampling. The idea of persistent homology is to take several scales at once. We do this by growing each point to a ball such that every ball has the same radius. We will freeze the process for a set of radii, build a complex, and then compute its homology.

We would like to build an abstract simplicial complex whose homology reflects that of the shape produced by taking the union of balls around points. We start with a complex that has been around for a long time. The definitions and results in this section are taken from \cite{RB} unless otherwise specfied. 

\begin{definition}
For $x\in R^n$ and $\epsilon>0$, let $B_\epsilon(x)$ be the open ball of radius $\epsilon$ centered at $x$. The {\it \u{C}ech complex}\index{Cech complex} $C_\epsilon(X)$ for a finite set of points $X\subset R^n$ is an abstract simplicial complex with the  points of $X$ as vertices and a $k$-simplex $[v_0, v_1, \cdots, v_k]$ for $\{v_0, v_1, \cdots, v_k\}\subset X$ if $$\bigcap_i B_\epsilon(v_i)\neq 0.$$
\end{definition}

It turns out that the underlying space of the \u{C}ech complex is actually homeomorphic to the union of balls centered at the points of $X$. (See the references cited in \cite{RB}.) The problem is that for high dimensions, it is not easy to tell when the intersections of the balls are non-empty. An easier task would be to look at whether the balls intersect pairwise. In this case, we would only need to know pairwise distances between the points. We would then include a simplex if the distances between each pair of vertices is less than or equal to than twice the radius of the balls. The resulting complex is the Vietoris-Rips complex. Note that we can now define a simplex in any metric space with distance $d$.

\begin{definition}
Let $X$ be a finite metric space with distance $d$ and fix $\epsilon>0$. The {\it Vietoris-Rips complex}\index{Vietoris-Rips complex} $VR_\epsilon(X)$ is an abstract simplicial complex with the  points of $X$ as vertices and a $k$-simplex $[v_0, v_1, \cdots, v_k]$ whenever $d(v_i, v_j)\leq 2\epsilon$ for all $0\leq i, j\leq k$. 
\end{definition}

The Vietoris-Rips complex is easier to compute than the \u{C}ech complex but is not always the same as the following example from \cite{RB} shows.

\begin{example}
Let $X=\{(0, 0), (1,0), (\frac{1}{2}, (\frac{\sqrt{3}}{2})\}\subset R^2$. These points are the vertices of an equilateral triangle with side 1. If $\epsilon>\frac{1}{2}$, then any pair of balls of radius epsilon intersect so the Vietoris-Rips and \u{C}ech complexes have the same 1-skeletons. But the centroid of the triangle which is equidistant from all three vertices is a distance $\frac{\sqrt{3}}{3}\approx .577$ from each vertex. So for $\frac{1}{2}< \epsilon< \frac{\sqrt{3}}{3}$, the pairs of balls interesect but the intersection of all three is empty. So the  Vietoris-Rips complex has a 2-simplex but the \u{C}ech complex does not. As the union of the balls has a hole, the 1-homology of the union has Betti-number of 1 which agrees with the \u{C}ech complex but not the Vietoris-Rips complex which is acyclic.
\end{example}

It turns out that there is a relationship between the Vietoris-Rips complex and the \u{C}ech complex. I will state it here without proof:

\begin{theorem}
If $X$ is a finite set of points in $R^n$ and $\epsilon>0$ then we have the inclusions $$C_\epsilon(X)\subset VR_\epsilon(X)\subset C_{2\epsilon}(X).$$
\end{theorem}

One thing to notice is that for both the Vietoris-Rips complex and the \u{C}ech complex, increasing the scale parameter $\epsilon$ increases the number of simplices in the complex. This is a consequence of the fact that increasing the radius of the balls centered at each point makes them more likely to intersect. This fact is the motivation for my interest in the connection between obstruction theory and TDA. Obstruction theory looks at extending a map from $L$ to $Y$ continuously to a map from $K$ to $Y$ if $L$ is a subcomplex of $K$. Picking the proper map may give interesting information about the data points. I will discuss this subject further in later chapters after I define cohomology and homotopy groups as they are needed to provide the criterion for extending maps in this way.

The increasing sequence of complexes is called a {\it filtration}\index{filtration}. The idea will be to choose several values of epsilon and compute the resulting homology in dimension $k$ of the homology of the Vietoris-Rips complex, $VR_\epsilon(X)$. As $\epsilon$ grows, we keep track of each $k$-dimensional homology class and determine at what value of $\epsilon$ the class is born and at what values it dies. The idea is that classes with a longer lifetime are the most important. The lifetime of a homology class is its {\it persistence} and this process computes the {\it persistent homology}\index{persistent homology} of a point cloud. Figure 5.1.1 shows the resulting complexes in an example point cloud as $\epsilon$ grows.

\begin{figure}[ht]
\begin{center}
  \scalebox{0.4}{\includegraphics{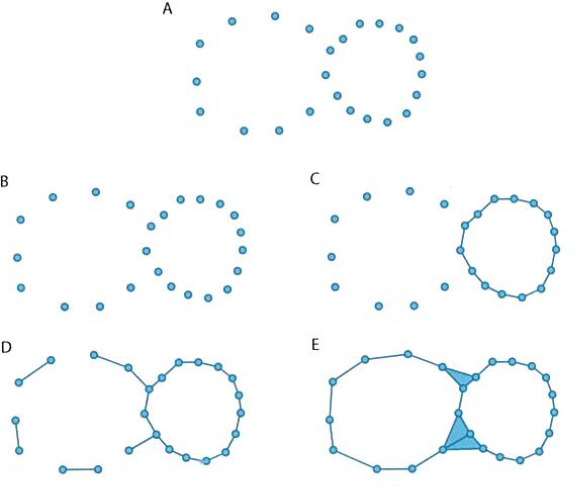}}
\caption{
\rm 
Vietoris-Rips Complexes of a Point Cloud with $\epsilon$ Growing \cite{RB}
}
\end{center}
\end{figure}

This leads to some issues. First of all, what coefficients should we use for our homology groups. The most common is to use $Z_2$. This avoids the problem of orientation as in $Z_2$, $1=-1$. But you are throwing away information. Why not use integer coefficients as we do in Chapter 4? It turns out, we need to have coefficients in a {\it field} for persistent homology and we will explore the reasons in the next section. Also, what values of $\epsilon$ should we use? We would like to capture the changes and with finitely many points, simplices are added only finitely many times leading to only finitely many changes in the homology groups. So equal spaceing may work but you need the intervals to be adjusted to how far apart your points are. Finally, what values of $k$ should we use to compute $H_k(VR_\epsilon(X)$. It is best to start with small values as these are easiest to interpret. $H_0$ represents the number of components, $H_1$ represents the number of holes, and $H_2$ represents the number of voids (think of the spherical holes in a block of Swiss cheese. In Section 5.3, we will see how to visualize these classes. 

\section{Computational Issues}

In this section, I will discuss the paper of Zomorodian and Carlsson \cite{ZC}, which is one of the earliest papers on computing persistent homology. I won't copy the entire algorithm. You can look there, but keep in mind that there are much faster algorithms that have been developed since then and there is a lot of open source software that will do these computations for you. Our interest will be in providing a more algebraic description of persistent homology and answering the question of why we need our coefficients to be in a field. 

First lets review some ring theory. Recall that a Euclidean ring allows for long division. Also, recall Theorem 3.4.2 which states that any Eucllidean ring is a principal ideal domain (PID). You may wonder if the opposite is true. Herstein \cite{Her1} says that the answer is no, and gives the paper by Motzkin \cite{Mot} as a reference. That paper is lacking in details and somewhat difficult to understand. A much better paper is that of Wilson \cite{Wil2}. Letting $R$ be the subring of the complex numbers of the form $$R=\{a+b\frac{1+\sqrt{-19}}{2}| a, b\in Z\},$$ he shows that $R$ is a PID which is not a Euclidean ring.

Recall that we also mentioned that the polynomial ring $F[x]$ is a Euclidean ring if $F$ is a field. I will outline the proof from Herstein \cite{Her1}.

First of all, an element of a Euclidean ring has a {\it degree}. For a polynomial, this is the usual degree. A generic element of $F[X]$ is $f(x)=a_0+a_1x+a_2x^2+\cdots+a_nx^n$, where the coefficients $a_0, a_1, \cdots, a_n\in F.$ In this case, the degree of $f$ is deg $f=n$.

Now let $f(x)=a_0+a_1x+a_2x^2+\cdots+a_nx^n$ and $g(x)=b_0+b_1x+b_2x^2+\cdots+b_mx^m$.Assume $a_n, b_m\neq 0$.  Multiplying polynomials in the usual way gives that the highest degree term is $a_nb_mx^{n+m}$. So we get the folllowing:

\begin{theorem}
If $f(x)$ and $g(x)$ are nonzero elements of $F[x]$, then deg $(f(x)g(x))=$ deg $f(x)+$ deg $g(x)$.
\end{theorem}

\begin{theorem}
If $f(x)$ and $g(x)$ are nonzero elements of $F[x]$, then deg $f(x)\leq$ deg $(f(x)g(x))$.
\end{theorem}

{\bf Proof:} This follows from Theorem 5.2.1 and the fact that for $g(x)\neq 0$, deg $g(x)\geq 0. \blacksquare$

\begin{theorem}
$F[x]$ is an integral domain.
\end{theorem}

This immediate from the previous result.

We will now be able to show that $F[x]$ is a Euclidean ring if we can prove the division algorithm.

\begin{theorem}
If $f(x)$ and $g(x)$ are nonzero elements of $F[x]$, then there exist polynomials $t(x)$ and $r(x)$ in $F[x]$ such that $f(x)=t(x)g(x)+r(x)$ where $r(x)=0$ or deg $r(x)<$ deg $g(x)$.
\end{theorem}

{\bf Proof:} If the degree of $f(x)$ is less than the degree of $g(x)$ we can set $t(x)=0$ and $r(x)=f(x)$. So assume that deg $f(x)\geq$ deg $g(x)$;.

Let $f(x)=a_0+a_1x+a_2x^2+\cdots+a_mx^m$ and $g(x)=b_0+b_1x+b_2x^2+\cdots+b_nx^n$ with $a_m, b_n\neq 0$, and $m\geq n$. We can proceed by induction on the degree of $f(x)$. If $f$ has degree 0, then so does $g$. So let $f(x)=a_0$ and $g(x)=b_0$. Then since $F$ is a field and $a_0$ and $b_0$ are nonzero, we have $f(x)=(a_0/b_0)g(x).$ Here division means multiplying by $b_0^{-1}$ which we can do since $F$ is a field and all nonzero elements have a multiplicative inverse.

If the degree of $f(x)$ is greater than zero, let $f_1(x)=f(x)-(a_m/b_n)x^{m-n}g(x)$, Then deg $f_1(x)\leq m-1$, so by induction on the degree of $f(x)$, we can write $f_1(x)=t_1(x)g(x)+r(x)$ where $r(x)=0$ or deg $r(x)<$ deg $g(x)$. So $f(x)-(a_m/b_n)x^{m-n}g(x)=t_1(x)g(x)+r(x)$ or $f(x)=(a_m/b_n)x^{m-n}g(x)+t_1(x)g(x)+r(x)=t(x)g(x)+r(x),$ where $t(x)=(a_m/b_n)x^{m-n}+t_1(x).$ So we are done. $\blacksquare$

\begin{theorem}
If $F$ is a field, then the polynomial ring $F[x]$ is a Euclidean ring. Consequently, it is also a principal ideal domain. 
\end{theorem}

Now letting $Z$ be the integers, we see the problem with this proof in the case of $Z[x]$. Unless the coefficient of the highest degree term of $g(x)$ is either 1 or -1, we can't divide by it. In fact, $Z[x]$ is not even a pricipal ideal domain. Let the ideal $U$ be generated by 2 and $x$. A principal ideal would have to be generated by an element of minimal degree that is of the form $2a+bx$ where $a$ and $b$ are integers. So the only possible generators would be 2 or -2. But $x$ is not a multiple of $2$ since $2$ has no inverse in $Z$. We would have the same problem with $-2$. So $Z[x]$ is not even a PID, despite the integers being a PID. 

Now we will show that the chain complexes we will use for persistent homology will be modules over the polynomial ring $R[t]$. This will be a PID if $R$ is a field. In that case, we will make use of the structure theorem of modules over a PID (Theorem 3.4.3).

Now I will summarize the description of persistence from \cite{ZC}. Assume we have a filtration of complexes. The $i$-th complex will be denoted $K^i$. Note that this section will not have any reference to cohomology in which superscripts are generally used. Here the superscript will be an index for a filtration. For the Vietoris-Rips complex, we could think of the superscript indexing the radius of the balls centered at points in a point cloud. For $K^i$, $C_k^i$, $Z_k^i$, $B_k^i$, and $H_k^i$,  are the groups of $k$-chains, $k$-cycles, $k$-boundaries, and $k$-th homology group respectively. The {\it p-persistent k-th homology group of} $K^i$ is $$H_k^{i,p}=Z_k^i/(B_k^{i+p}\bigcap Z_k^i).$$

Note that both groups in the denominator are subgroups of $C_k^{i+p}$ as the complexes are getting bigger as the superscript increases so the chain group $C_k^i$ is a subgroup of $C_k^j$ for $j>i$. So the intersection of the groups in the denominator is also a group and a subgroup of the numerator.

Now recall from Theoremm 3.4.3 that a finitely generated module over a principal ideal domain R is the direct sum of copies of $R$ (the free part) and summands of the form $R/(d_i)$ (the torsion part). The {\it Betti number} $\beta$ is the {\it rank} or number of summands of $R$ in the free part.

\begin{definition}
The {\it p-persistent k-th Betti number} of $K^i$ is the rank $\beta_k^{i,p}$ of the free subgroup of $H_k^{i,p}$. We can also define persistent homology groups using the injection $\eta_k^{i,p}: H_k^i\rightarrow H_k^{i+p}$ coming from the map on chains. We can define $H_k^{i,p}$ as the image of $\eta_k^{i,p}$.
\end{definition}

We are still left with the problem of finding compatible bases for $H_k^i$ and $H_k^{i+p}$. To get around this problem, \cite{ZC} combines all of the complexes of the filtration into a single algebraic structure. They then establish a correspondence that gives a simple description when the coefficients are in a field. Persistent homology is represented as the homology of a graded module over a polynomial ring. 

\begin{definition}
A {\it persistence complex}\index{persistence complex} $\mathcal{C}$ is a family of chain complexes $\{C_*^i\}_{i\geq 0}$ over $R$ together with chain maps $f^i: C_*^i\rightarrow C_*^{i+1}$ so that we have the following diagram:

$$\begin{tikzpicture}
  \matrix (m) [matrix of math nodes,row sep=3em,column sep=4em,minimum width=2em]
  {
 C_*^0 &  C_*^1 &  C_*^2 & \cdots\\};
  \path[-stealth]

(m-1-1) edge node [above] {$f^0$} (m-1-2)
(m-1-2) edge node [above] {$f^1$} (m-1-3)
(m-1-3) edge node [above] {$f^2$} (m-1-4)

;

\end{tikzpicture}$$

\end{definition}

The filtered complex $K$ with inclusion maps for the simplices is a persistence complex. Here is a part of the complex. The filtration indices increase as you move right and the dimension goes down as you move down. $\partial_k$ are the usual boundary operators. 

$$\begin{tikzpicture}
  \matrix (m) [matrix of math nodes,row sep=3em,column sep=4em,minimum width=2em]
  {
\hspace{1pt} & \hspace{1pt} & \hspace{1pt} & \\
 C_2^0 &  C_2^1 &  C_2^2 & \cdots\\
 C_1^0 &  C_1^1 &  C_1^2 & \cdots\\
 C_0^0 &  C_0^1 &  C_0^2 & \cdots\\};
  \path[-stealth]

(m-2-1) edge node [above] {$f^0$} (m-2-2)
(m-2-2) edge node [above] {$f^1$} (m-2-3)
(m-2-3) edge node [above] {$f^2$} (m-2-4)

(m-3-1) edge node [above] {$f^0$} (m-3-2)
(m-3-2) edge node [above] {$f^1$} (m-3-3)
(m-3-3) edge node [above] {$f^2$} (m-3-4)

(m-4-1) edge node [above] {$f^0$} (m-4-2)
(m-4-2) edge node [above] {$f^1$} (m-4-3)
(m-4-3) edge node [above] {$f^2$} (m-4-4)

(m-1-1) edge node [left] {$\partial_3$} (m-2-1)
(m-1-2) edge node [left] {$\partial_3$} (m-2-2)
(m-1-3) edge node [left] {$\partial_3$} (m-2-3)

(m-2-1) edge node [left] {$\partial_2$} (m-3-1)
(m-2-2) edge node [left] {$\partial_2$} (m-3-2)
(m-2-3) edge node [left] {$\partial_2$} (m-3-3)

(m-3-1) edge node [left] {$\partial_1$} (m-4-1)
(m-3-2) edge node [left] {$\partial_1$} (m-4-2)
(m-3-3) edge node [left] {$\partial_1$} (m-4-3)

;

\end{tikzpicture}$$

\begin{definition}
A {\it persistence module} $\mathcal{M}$ is a family of $R$--modules $M^i$ together with homomorphisms $\phi^i: M^i\rightarrow M^{i+1}$. 
\end{definition}

\begin{definition}
A peristence complex (resp. persistence module) is of {\it finite type} if each component complex (module) is a finitely generated $R$-module, and if the maps $f^i$ (resp. $\phi^i$) are isomorphisms for $i\geq m$ for some $m$.
\end{definition}

For a Vietoris-Rips complex $K$ produced by finitely many points, the associated persistence complex and persistence module (the homology of the complex) are of finite type. If $m$ is the largest distance between any pair of points, the complex $VR_m(K)$ consists of one simplex whose vertices are all of the points and it no longer changes when $\epsilon>m$.  

Now we need to make the algebra a little more precise. A {\it graded ring} $R$ will be a direct sum of abelian groups $R_i$ such that if $a\in R_i$ and $b\in R_j$ then $ab\in R_{i+j}$. Elements in a single $R_i$ are called {\it homogeneous}. The polynomial ring $R[t]$ is a graded ring with standard grading $R_i=Rt^i$ for all nonnegative integers $i$. $at^i$ and $bt^j$ for $a, b\in R$ are homogeneous elements while there sum is not. The product $abt^{i+j}$ has degree $i+j$ as required. We define a graded module $M$ over a graded ring $R$ in a similar way. $M=\oplus M_i$ where the sum is over all integers, and for $r\in R_i$ and $m\in M_j$, we have $rm\in M_{i+j}$. A graded ring (resp. module) is non-negatively graded if $R_i=0$ (resp. $M_i=0$) for $i<0$. 

In this notation, we write the structure theorem in the following form: 

\begin{theorem}
If $D$ is a PID, then every finitely generated $D$-module can be uniquely decomposed as $$D\cong D^\beta\oplus (\bigoplus_{i=1}^m D/d_iD),$$ for $d_i\in D, \beta\in Z,$ and $d_i$ divides $d_{i+1}$ for all $i$. A graded module $M$ over a graded PID $D$ decomposes uniquely as $$M\cong (\bigoplus_{1=1}^n \Sigma^{\alpha_i}D)\oplus(\bigoplus_{j=1}^m \Sigma^{\gamma_j} D/d_jD),$$ where $d_j\in D$ are homogeneous elements so that $d_j$ divides $d_{j+1}$, $\alpha_j, \gamma_j\in Z$, and $\Sigma^\alpha$ denotes and $\alpha$-shift upward in grading.
\end{theorem}

For a persistence module $\mathcal{M}=\{M^i, \phi^i\}_{i\geq 0}$ over a ring $R$, give $R[t]$ the standard grading and define a graded module over $R[t]$ by $$\alpha(\mathcal{M})=\bigoplus_{i=0}^\infty M^i,$$ where the the action of $r$ is just the sum of the actions on each component and $$t(m_0, m_1, m_2, \cdots)=(0, \phi^0(m_0), \phi^1(m_1), \phi^2(m_2), \cdots).$$ So $t$ shifts elements of the module up a degree.What this means in our case is that if a simplex is added at time $i$, $t$ shifts it to time $i+1$ but preserves the memory that it had been added earlier by applying $\phi^i$. 

So $\alpha$ provides a correspondence between persistence modules of finite type over $R$ and finitely generated non-negatively graded modules over $R[t]$. So we want $R$ to be a field so that $R[t]$ becomes a PID and we can use the structure theorem. So we now have that our persistence module is of the form $$(\bigoplus_{1=1}^n \Sigma^{\alpha_i}R[t])\oplus(\bigoplus_{j=1}^m \Sigma^{\gamma_j} R[t]/(t^{n_j})),$$ since $R[t]$ is a PID and all principal ideals are of the form $(t^n)$.

Now we get to the key point. We want a more meaningful description of the possible graded modules over $R[t]$. Define a $\mathcal{P}$-{\it interval} is an ordered pair $(i, j)$ with $0\leq i<j$, where $i\in Z$ and $j\in Z\cup\{+\infty\}$. Let $Q(i, j)$ for a $\mathcal{P}$-interval $(i,j)$ to be $\Sigma^i R[t]/(t^{j-i})$, and $Q(i, +\infty)=\Sigma^iR[t]$. For a set of $\mathcal{P}$-intervals $\mathcal{S}=\{(i_1, j_1), (i_2, j_2),\cdots,( i_n, j_n)\}$ , define $$Q(\mathcal{S})=\bigoplus_{k=1}^n Q(i_k, j_k).$$ So we have a one-to one correspondence between finite sets of  $\mathcal{P}$-intervals and finitely generated modules over $R[t]$ where $R$ is a field. 

So we know now why we wanted the coefficients to be in a field. Usually the field will be $Z_2$ since we don't have to worry about orientations. But what do the $\mathcal{P}$-intervals represent? Suppose we have the interval  $(i, j)$. This corresponds to the summand $\Sigma^i R[t]/(t^{j-i})$. If the module represents the homology, the shift $\Sigma^i$ means that the homology class in this summand doesn't show up until time $i$, while the quotient $(t^{j-i})$ means that it disappears by time $j$. The time $i$ is called the {\it birth} time of the class, and the time $j$ is the {\it death} time. The correspondence and the structure theorem imply that we can characterize the persistent homology by keeping track of all of the birth and death times of the homology classes in $H_k(X)$ for each $k$. This will be done using the visualizations in the next section.

The rest of \cite{ZC} provides an early algorithm that modifies the algorithm we saw in Section 4.1.4 for the persistent case where the coefficients are in a field. There is also a modification for coefficients in a PID such as the integers but only for a pair of times $i$ and $i+p$. Knowing the birth and death times of the classes doesn't uniquely determine the persistence module in this case, but we can compute regular homology at fixed times over the integers. As we will see in Chapter 8, we can use the {\it Universal Coefficient Theorem} to obtain the homology in any other abelian group from the homology with integer coefficients. We will describe that process when we cover some homological algebra. Meanwhile, in the rest of this chapter, we will always work in a field and use $Z_2$ unless otherwise specified.

I will not say anything more about the algorithm in \cite{ZC}. It has historical interest, but has since been replaced by much more efficient algorithms.

\section{Bar Codes, Persistence Diagrams, and Persistence Landscapes}

Now we know that we can characterize persistent homology by birth and death times of the homology classes in each dimension. The next step is to produce a visualization. The three most popular are bar codes, persistence diagrams, and persistence landscapes. 

These visualizations are best explained by using pictures. The first picture I will show is from the survey paper by Chazal and Michel \cite{CM}. 

\begin{figure}[ht]
\begin{center}
  \scalebox{0.4}{\includegraphics{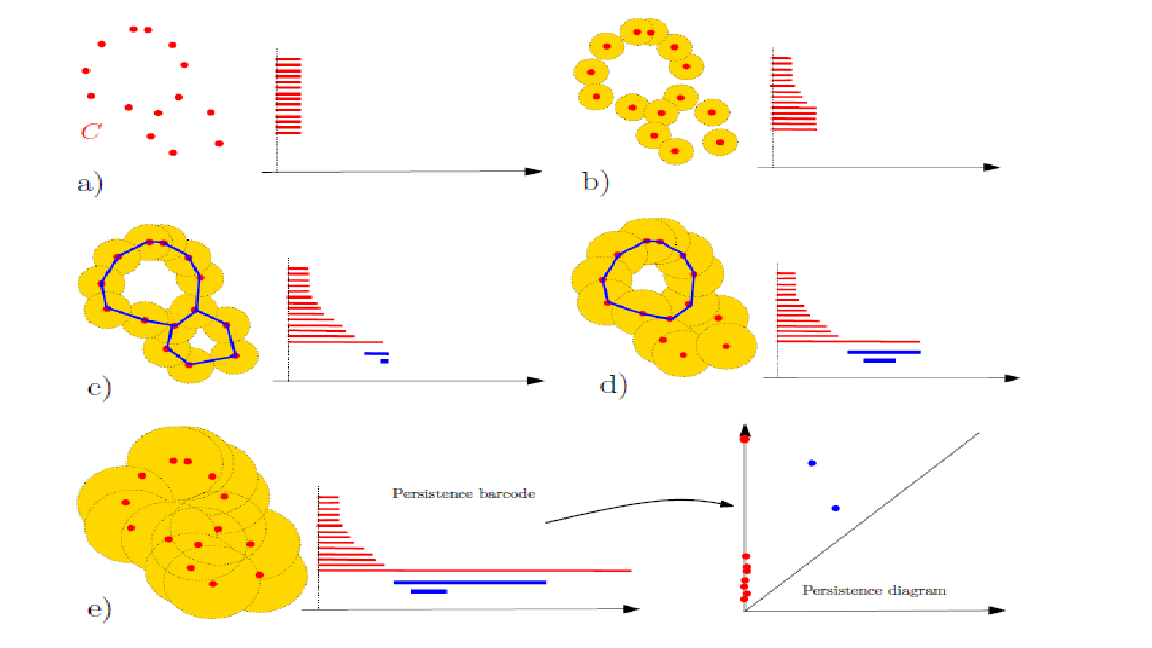}}
\caption{
\rm 
Bar Codes and Persitence Diagrams \cite{CM}
}
\end{center}
\end{figure}

{\it Bar codes}\index{bar code} were the first visualization tried. The $x$-axis represents the filtration index or time (here the radius of the balls.) The $y$-axis doesn't have a meaning. There is just a bar for every homology class. In Figure 5.3.1, the 0 and 1-dimensional homology are combined on one chart. Usually, we make a separate chart for each dimension, but this diagram is simple enough to combine them. The red bars are in dimension 0 and the blue ones are dimension 1. If you can't see colors on your copy, there are two blue bars at the bottom of charts c, d, and e. 

In picture a, we have a set of disconnected dots. The bars should probably have length 0 at this stage, but they are longer so you can see them. In b, many of them have merges so many of the 0-dimensional bars have died by then. In picture c, we have two complete loops so we get the two blue bars on the bottom. From picture b, we see that the top loop was born before the bottom one, so the corresponding bar starts further to the left. Also, we now have a single component, so there is only one red bar that is still growing. In picture d, the bottom hole has closed up so the bottom bar has ended. There is now one red bar corresponding to the one component and one blue bar corresponding to the one hole. In picture e, all of the holes have closed and the homology stops evolving. We now have the full bar code picture with one surviving red bar and all of the others have ended.

Bar codes present some problems though. The main one is an easy distance measure between bar codes for different point clouds. Also, they can be a little clumsy if there are a lot of homology classes. A more elegant visualization is a {\it persistence diagram}\index{persistence diagram}. At the bottom right of Figure 5.3.1, the final bar chart is transformed into a persistence diagram. This is more of a conventional scatter plot with the $x$-axis representing birth times and the $y$-axis representing death times. Again, these times correspond to the radius of the balls when the class is formed and when it dies. For each class, we plot a single point whose $x$ coordinate is its birth time and whose $y$ coordinate is its death time. In Figure 5.3.1, the 0-dimensional and 1-dimensional classes are again plotted together. All of the 0-dimensional classes sit on the $y$-axis as they were all born at time 0 but die at different times. The two dots towards the right are the 1-dimensional classes. Note that the higher one represents the top loop as it is born earlier(so it is further left) and dies later (so it is higher) than the bottom loop.

Note that the persistence diagram always sits above the diagonal line $y=x$. After all, a class can't die before it is born. Also, the set of points plotted is known as a {\it multiset} as the same point can appear multiple times. There is no reason that two or more classes can't both be born and die at the same time as each other. 

{\it Persistence landscapes}\index{persistence landscape} were developed by Bubenik \cite{Bu1} who felt they wee more suitable for statistical analysis. The reason has to do with distance measures so I will defer that to Section 5.5. Like bar codes and presistence diagrams, we start with a set of birth time-death time intervals. Rather than form a scatterplot, for a persistence landscape, we construct a function for each of these intervals. The function is shaped like a triangle and its base and height are proportional to its lifetime. Suppose homology class $i$ is born at time $b_i$ and dies at time $d_i$. Then we define a piecewise linear function $$f_{(b_i, d_i)}=
\left
\{
\begin{array}{ll}
x-b_i & \mbox{if }b_i<x\leq \frac{b_i+d_i}{2},\\
-x+d_i & \mbox{if }\frac{b_i+d_i}{2}<x<d_i,\\
0 & \mbox{otherwise.}
\end{array}
\right.$$

These functions are triangles with base from $b_i$ to $d_i$ and height $\frac{d_i-b_i}{2}$. So they are wider and taller the longer their lifetime. 

\begin{figure}[ht]
\begin{center}
  \scalebox{0.4}{\includegraphics{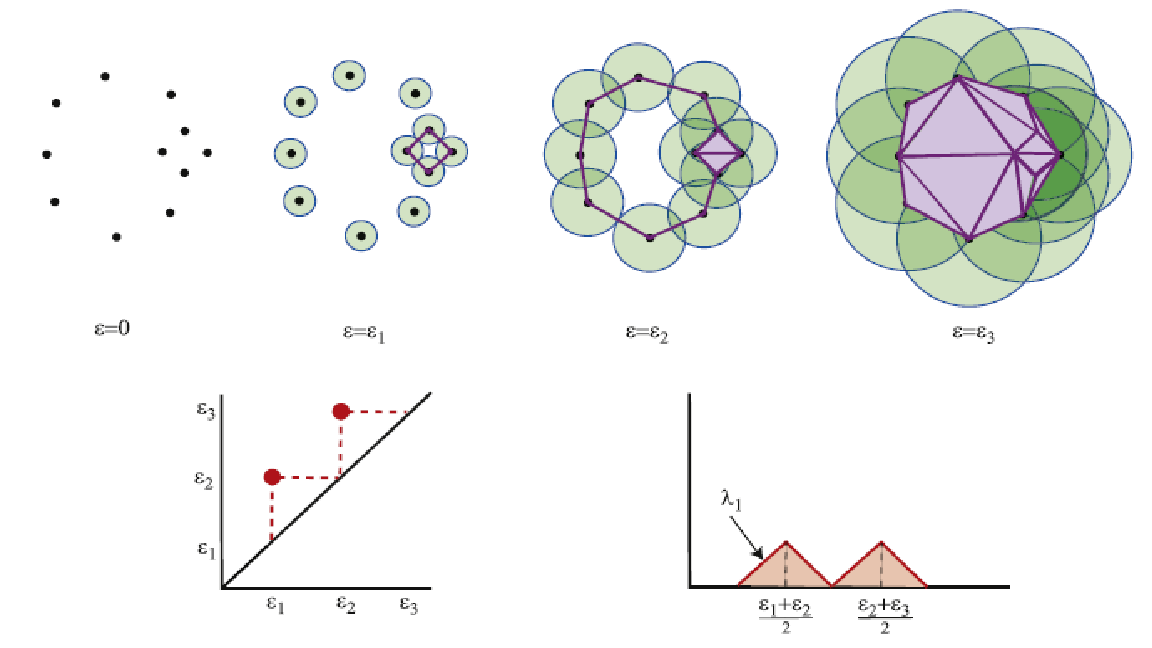}}
\caption{
\rm 
Persitence Diagrams vs. Persistence Landscapes \cite{GK}
}
\end{center}
\end{figure}

Figure 5.3.2 shows an example from \cite{GK}. In this case we are only plotting the 1-dimesnional homology. Let $\epsilon$ be the radius of the circles. The small hole along the right edge is born at $\epsilon=\epsilon_1$ and dies at $\epsilon=\epsilon_2$. The large hole in the middle is born at $\epsilon=\epsilon_2$ and dies at $\epsilon=\epsilon_3$. The persistence diagram at the bottom left includes points at $(\epsilon_1, \epsilon_2)$ and $(\epsilon_2, \epsilon_3)$. In the persistence landscape at the bottom right, we have 2 triangular functions. The one on the left has a base from $\epsilon_1$ to $\epsilon_2$, while the one on the right has a base from $\epsilon_2$ to $\epsilon_3$. The heights are $\frac{\epsilon_2-\epsilon_1}{2}$ and $\frac{\epsilon_3-\epsilon_2}{2}$ respectively. 

In a real life case there will be lots of triangles and they will overlap. This is a foreshadowing of what that $\lambda_1$ is doing in the picture.

I will end this section on a cliffhanger and finish this story in Section 5.5 after one short diversion.

\section{Zig-Zag Persistence}

In our previous discussion we were dealing with a {\it filtration} in which are complexes are growing monotonically. If this is not the case, there is an alternate notion of {\it zig-zag persistence}\index{zig-zag persistence} \cite{CDM}. The material here is taken from \cite{RB}. 

Suppose we take multiple sets of samples $X_i$ form a metric space $X$. We can form a sequence $$\begin{tikzpicture}
  \matrix (m) [matrix of math nodes,row sep=3em,column sep=4em,minimum width=2em]
  {
X_1 & X_1\cup X_2 & X_2 &  X_2\cup X_3 & X_3 & \cdots\\};
  \path[-stealth]

(m-1-1) edge  (m-1-2)
(m-1-3) edge  (m-1-2)
(m-1-3) edge  (m-1-4)
(m-1-5) edge  (m-1-4)
(m-1-5) edge  (m-1-6)

;

\end{tikzpicture}$$ where the maps are all inclusions. This gives rise to the corresponding diagram $$\begin{tikzpicture}
  \matrix (m) [matrix of math nodes,row sep=3em,column sep=4em,minimum width=2em]
  {
H_k(VR_\epsilon(X_1)) & H_k(VR_\epsilon(X_1\cup X_2)) &  H_k(VR_\epsilon(X_2)) & \\
& & H_k(VR_\epsilon(X_2\cup X_3)) & \\
& & H_k(VR_\epsilon(X_3)) & \cdots\\};
  \path[-stealth]

(m-1-1) edge  (m-1-2)
(m-1-3) edge  (m-1-2)
(m-1-3) edge (m-2-3)
(m-3-3) edge  (m-2-3)
(m-3-3) edge  (m-3-4)

;

\end{tikzpicture}$$

We now define zigzag diagrams more generally. They can have other shapes besides the one just shown.Keep in mind that homology is taken with coefficients in a field $F$ so the homology groups are actually vector spaces over $F$. 

\begin{definition}
A {\it zigzag diagram} or {\it zigzag module}\index{zigzag diagram}\index{zigzag module} of shape $S$ is a sequence of linear transformations between $F$ vector spaces: $$\begin{tikzpicture}
  \matrix (m) [matrix of math nodes,row sep=3em,column sep=4em,minimum width=2em]
  {
X_1 & X_2 & \cdots & X_k\\};
  \path[-stealth]

(m-1-1) edge[-]  node[above]{$f_1$}(m-1-2)
(m-1-2) edge[-]  node[above]{$f_2$} (m-1-3)
(m-1-3) edge[-]  node[above]{$f_{n-1}$} (m-1-4)

;

\end{tikzpicture}$$ where each map $f_i$ has its direction specified by the $i$-th letter in $S$. Each letter is $R$ or $L$.
\end{definition}

\begin{example}
If $S$ is a string of $R's$ or $L's$ then we have an ordinary filtration.
\end{example}

\begin{example}
If $S=RRL$ then a zigzag diagram of shape $S$ is $$\begin{tikzpicture}
  \matrix (m) [matrix of math nodes,row sep=3em,column sep=4em,minimum width=2em]
  {
X_1 & X_2 & X_3 & X_4\\};
  \path[-stealth]

(m-1-1) edge (m-1-2)
(m-1-2) edge  (m-1-3)
(m-1-4) edge (m-1-3)

;

\end{tikzpicture}$$
\end{example}

\begin{example}
If $S=RLRLRL$ then a zigzag diagram of shape $S$ is $$\begin{tikzpicture}
  \matrix (m) [matrix of math nodes,row sep=3em,column sep=4em,minimum width=2em]
  {
X_1 & X_2 & X_3 & X_4 & X_5 & X_6 & X_7\\};
  \path[-stealth]

(m-1-1) edge (m-1-2)
(m-1-3) edge  (m-1-2)
(m-1-3) edge (m-1-4)
(m-1-5) edge (m-1-4)
(m-1-5) edge (m-1-6)
(m-1-7) edge (m-1-6)

;

\end{tikzpicture}$$
\end{example}

In the case of filtrations, we would look for features that survived a chain of inclusion maps. For a zig-zag filtration, we look for consistency across parts of the diagram. For example, if a zigzag diagram has shape RL, $$\begin{tikzpicture}
  \matrix (m) [matrix of math nodes,row sep=3em,column sep=4em,minimum width=2em]
  {
X_1 & X_2 & X_3,\\};
  \path[-stealth]

(m-1-1) edge node [above] {$f_1$} (m-1-2)
(m-1-3) edge  node [above] {$f_2$} (m-1-2)

;

\end{tikzpicture}$$ we look for features $x_1\in X_1$, $x_2\in X_2$, and $x_3\in X_3$ such that $f_1(x_1)=x_2=f_2(x_3).$

In the next part, we use the term {\it zigzag module} to emphasize the algebraic structure.

\begin{definition}
A {\it zigzag submodule}\index{zigzag submodule} $N$ of a zigzag module $M$ of shape $S$ is a zigzag module of shape $S$ such that each space $N_i$ in $N$ a subspace of the space $M_i$ in $M$ and the maps in $N$ are the restrictions of those in $M$.
\end{definition}

\begin{example}
Let $F=R$, and suppose we have the zigzag module $$\begin{tikzpicture}
  \matrix (m) [matrix of math nodes,row sep=3em,column sep=4em,minimum width=2em]
  {
R & R^2 & R,\\};
  \path[-stealth]

(m-1-1) edge node [above] {$f_1$} (m-1-2)
(m-1-3) edge  node [above] {$f_2$} (m-1-2)

;

\end{tikzpicture}$$ where for $x\in R$, $f_1(x)=(x,0)$ and $f_2(x)=(0, x).$ Then we have a zigzag submodule $$\begin{tikzpicture}
  \matrix (m) [matrix of math nodes,row sep=3em,column sep=4em,minimum width=2em]
  {
R & R & R\\};
  \path[-stealth]

(m-1-1) edge node [above] {$g_1$} (m-1-2)
(m-1-3) edge  node [above] {$g_2$} (m-1-2)

;

\end{tikzpicture}$$ where the $R$ in the middle is the first coordinate of $R^2$ so that for $x\in R, g_1(x)=x$, and $g_2(x)=0$.
\end{example}

A zigzag module is {\it decomposable} if it  can be written as the direct sum of non-trivial submodules. Otherwise it is {\it indecomposable.}

\begin{theorem}
Any zigzag module of shape $S$ can be written as a direct sum of indecomposables in a way that is unique up to permutation. 
\end{theorem}

\begin{definition}
An {\it interval zigzag module} of shape $S$ is a zigzag module $$\begin{tikzpicture}
  \matrix (m) [matrix of math nodes,row sep=3em,column sep=4em,minimum width=2em]
  {
X_1 & X_2 & \cdots & X_k\\};
  \path[-stealth]

(m-1-1) edge[-]  node[above]{$f_1$}(m-1-2)
(m-1-2) edge[-]  node[above]{$f_2$} (m-1-3)
(m-1-3) edge[-]  node[above]{$f_{n-1}$} (m-1-4)

;

\end{tikzpicture}$$ where for fixed $a\leq b$, $X_i=F$ for $1\leq a\leq i\leq b\leq k$, and $X_i=0$ otherwise.
\end{definition}

\begin{theorem}
The indecomposable zigzag modules are all in the form of interval zigzag modules. 
\end{theorem}

This theorem allows us to produce barcodes in this case by writing one bar extending from $a_i$ to $b_i$ for each indecomposable zigzag submodule whose nonzero terms extend from $a_i$ to $b_i$.

Here is one type of zigzag module we can build in a metric space. Let $X=\{x_1, x_2, x_3, \cdots, x_n\}$ be a finite metric space where an ordering of the points has been chosen, and for $1\leq k\leq n$, $X_k=\{x_1, x_2, x_3, \cdots, x_k\}$ consists of the first $k$ points of the ordering. 

\begin{definition}
Let $A$ and $B$ be non-empty subsets of a metric space $X$ with distance function $d$. The {\it Hausdorff distance}\index{Hausdorff distance}  between $A$ and $B$ is $$d_H(A, B)=\max(\sup_{a\in A}\inf_{b\in B} d(a,b), \sup_{b\in B}\inf_{a\in A} d(a,b)).$$
\end{definition}

So for sets $A$ and $B$, for each element of $A$, find the closest element of $B$ and find the distance. Then take the largest of these as we run over the elements of $A$. Then do the same thing with the roles of $A$ and $B$ reversed. Finally take the largest of those two values. 

Now for $1\leq k\leq n$, let $\epsilon_k=d_H(X_k, X)$. Note that $\epsilon_k\geq\epsilon_{k+1}$ since $X_{k+1}$ includes an additional point so it will always be at least as close to the full set $X$ as $X_k$ in Hausdorff distance.

\begin{definition}
Choose real numbers $\alpha>\beta>0$. The {\it Rips zigzag}\index{Rips zigzag} consists of the zigzag module specified by the diagram of simplicial complexes $$\begin{tikzpicture}
  \matrix (m) [matrix of math nodes,row sep=3em,column sep=4em,minimum width=2em]
  {
\cdots & VR_{\beta\epsilon_{i-1}}(X_{i-1}) & VR_{\alpha\epsilon_{i-1}}(X_i) &  VR_{\beta\epsilon_i}(X_i) & \\
& & & VR_{\alpha\epsilon_i}(X_{i+1}) & \\
& & & VR_{\beta\epsilon_{i+1}}(X_{i+1}) & \cdots\\};
  \path[-stealth]

(m-1-1) edge  (m-1-2)
(m-1-2) edge  (m-1-3)
(m-1-4) edge  (m-1-3)
(m-1-4) edge (m-2-4)
(m-3-4) edge  (m-2-4)
(m-3-4) edge  (m-3-5)

;

\end{tikzpicture}$$
\end{definition} 

The $\alpha$ and $\beta$ control the size of the complexes and were introduced for computational efficiency.

This is all I will say about zigzag persistence. See \cite{RB}, Section 5.4.3 for a use in studying dementia in AIDS patients.

\section{Distance Measures}

Now that we can produce a bar code or persistence diagram from a set of data points, we would like to define a distance measure between these diagrams.This will allow us to perform machine learning tasks such as clustering and classification. What we need from these measures is that they have the property of {\it stability}. In other words, we want small changes in our data to result in small changes in the distance between our persistence diagrams. An important aspect of TDA is that there are measures that behave in this way. The material in this section is taken from \cite{CM} and \cite{RB}.

Recall that a persistence diagram consists of a multiset of points on a graph each of which represent a homology class. The $x$-coordinate is the birth time and the $y$-coordinate is the death time. We also will include the diagonal and consider each point on it to have infinite multiplicity. 

\begin{definition}
Let $a=(a_1, \cdots, a_n)$ and $b=(b_1, \cdots, b_n)$ be points in $R^n$. Then the $L^p$ {\it distance} between $a$ and $b$ for $p$ a positive integer is $$d_p(a, b)=(\sum_{i=1}^n |a_i-b_i|^p)^{\frac{1}{p}}.$$ The $L^\infty$ distance is $$d_\infty(a, b)=\max_{1\leq i\leq n}|a_i-b_i|.$$
\end{definition}

Now suppose we have two diagrams $dgm_1$ and $dgm_2$. A matching $m$ between $dgm_1$ and $dgm_2$ is a bijective map between their points. A point of $dgm_1$ can be mapped to a point of $dgm_2$ or it can be mapped to the closest point on the diagonal. If $p=(b_1, d_1)$ and $q=(b_2, d_2)$ then $d_\infty(p, q)=\max(|b_1-b_2|, |d_1-d_2|)$. If $p=(b, d)$ and we want to match $p$ to a point on the diagonal, we map it to $(\frac{b+d}{2}, \frac{b+d}{2})$ and the $L^\infty$ distance is $\frac{|d-b|}{2}$. We want to find the matching that will minimize the worst distance between any pair of points in the matching. 

\begin{definition}
The {\it bottleneck distance}\index{bottleneck distance} between two diagrams $dgm_1$ and $dgm_2$ is $$d_B(dgm_1, dgm_2)=\inf_m\max_{(p, q)\in m}d_\infty(p, q).$$ Here the infimum is taken over all possible matchings where p is in $dgm_1$ or on the diagonal and $q$ is in $dgm_2$ or on the diagonal. 
\end{definition}

\begin{figure}[ht]
\begin{center}
  \scalebox{0.4}{\includegraphics{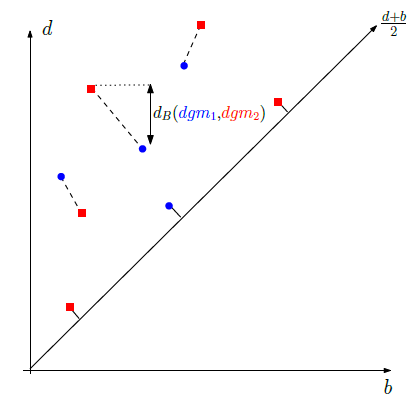}}
\caption{
\rm 
Example Matching of Two Persistence Diagrams \cite{CM}
}
\end{center}
\end{figure}

Figure 5.5.1 shows an example of matching 2 persistence diagrams. The blue circles belong to one diagram and the red squares belong to another. Note that some points from each diagram are mapped to the nearest point on the diagonal. 

The bottleneck distance is entirely determined by the furthest pair of points. To give more influence to the other points, we can use the {\it Wasserstein} distance instead. 

\begin{definition}
The {\it p-Wasserstein distance}\index{Wasserstein distance} between two diagrams $dgm_1$ and $dgm_2$ is $$d_{W_p}(dgm_1, dgm_2)=(\inf_m\sum_{(p, q)\in m}d_\infty(p, q)^p)^\frac{1}{p}.$$ As before, the infimum is taken over all possible matchings where p is in $dgm_1$ or on the diagonal and $q$ is in $dgm_2$ or on the diagonal. 
\end{definition}

Now recall form the last section that the Hausdorff distance $d_H$ is a nearness measure on subsets of a metric space. What if we want to generalize this when the metric spaces may not be the same. 

\begin{definition}
A function $f: X\rightarrow Z$ is an {\it isometry}\index{isometry} if $d_X(x_1, y_1)=d_Z(f(x_1), f(y_1))$. If $f$ is one-to-one then $f$ is an {\it isometric embedding} of $X$ into $Z$. 
\end{definition}

\begin{definition}
If $X_1$ and $X_2$ are compact metric spaces and can be isometrically embedded into $Z$ then we define the {\it Gromov-Hausdorff distance}\index{Gromov-Hausdorff distance} as $$d_{GH}(X_1, X_2)=\inf d_H(f_1(X_1), f_2(X_2))$$ where the infimum is taken over all possible isometric embeddings $f_i:X_i\rightarrow Z$, for $i=1, 2$, and $d_H$ denotes Hausdorff measure over subsets of $Z$. 
\end{definition}

The following stability theorem is from \cite{RB} and credited by them to \cite{CEH}. See \cite{CM} and its references for several other variants.

\begin{theorem}
Let $X$ and $Y$ be finite metric spaces. Let $PH_k(VR(X))$ be the $k$-dimensional persistence diagram of the Vietoris-Rips complex of $X$ and let $PH_k(VR(Y))$ be defined similarly. Then for all $k>0$, $$d_B(PH_k(VR(X)), PH_k(VR(Y))\leq d_{GH}(X, Y).$$ 
\end{theorem}

So a small change in a point cloud results in a small change in bottleneck distance between the respective persistence diagrams. Similar theorems hold for Wasserstein distance and persistence diagrams built in the other ways I will describe in the next few sections. 

You may wonder if there is an efficient way to calculate these distances. The answer is yes, but I won't go into them. The persistent homology packages that I will list at the end of this chapter can all do them for you. I will refer you to their descriptions and references if you are interested in efficient computation methods.

There has been some efforts to develop a statistical theory for persistent homology. For example, if you have a point cloud, how much of the true geometry have you covered? You really only have a sample. Is there a statistical theory complete with hypothesis testing, confidence intervals, etc.? This appears to be a pretty undeveloped area. If you are interested, see Chapter 3 of \cite{RB} or the relevant sections of \cite{CM}.

What I will say is that persistence landscapes \index{persistence landscape} provide the framework for some familiar statistical results such as the law of large numbers and the central limit theorem. The key is that we have the notion of a {\it mean} persistence landscape. In addition, we have a metric which makes landscapes into a {\it Banach space}, which I will now define.

\begin{definition} 
Let $X$ be a vector space over the real or complex numbers. Then $X$ is a {\it normed space}\index{normed space} if for every vector $x\in X$ we define a real number $||x||$ called the {\it norm} of $x$ with the following properties:
\begin{enumerate}
\item For every vector $x$, $||x||\geq 0$, and $||x||=0$ if and only if $x=0$.
\item For every vector $x$ and scalar $\alpha$, $||\alpha x||=|\alpha|||x||$.
\item {\bf Triangle Inequality:} For $x, y\in X$, $||x+y||\leq ||x||+||y||$.
\end{enumerate}
\end{definition}

Any normed space is automatically a metric space where we define for $x, y\in X$, the distance $d(x, y)=||x-y||.$

\begin{example}
$R^n$ is a normed space. If $x=(x_1, \cdots, x_n)\in R^n$, then define $$||x||=\sqrt{x_1^2+x_2^2+\cdots+x_n^2}.$$
\end{example}

To define a Banach space, we need the notion of completeness.

\begin{definition} 
Let $X$ be a metric space (or more generally any topological space). A {\it sequence}\index{sequence} of points in $X$ is the image of a function $f: Z^+\rightarrow X$ where $Z^+$ is the set of positive integers. We write a sequence as $\{x_1, x_2, \cdots\}$, where $x_i=f(i)$ for all $i\in Z^+$.
\end{definition}

\begin{definition} 
Let $X$ be a metric space. A sequence of points $\{x_1, x_2, \cdots\}$ in $X$ {\it converges} to a point $y\in X$ if given $\epsilon>0$, there exists an integer $N$ such that $d(y, x_i)<\epsilon$ for $i>N$. The point $y$ is called the {\it limit} of the sequence.
\end{definition}

\begin{definition} 
Let $X$ be a metric space. A sequence of points $\{x_1, x_2, \cdots\}$ in $X$ is a {\it Cauchy sequence}\index{Cauchy sequence} if given $\epsilon>0$, there exists an integer $N$ such that $d(x_i, x_j)<\epsilon$ for $i, j>N$. 
\end{definition}

\begin{definition} 
Let $X$ be a metric space.Then $X$ is {\it complete}\index{complete metric space} if every Cauchy sequence of points in $X$ converges to a point in $X$. If $X$ is a complete metric space which is also a normed space, then $X$ is a {\it Banach space}\index{Banach space}.
\end{definition}

\begin{example}
$R$ is a Banach space  and so is every closed subspace. The open interval $(0,1)$ is not a Banach space. Let $x_i=\frac{1}{i+1}$ for every positive integer $i$. This is a Cauchy sequence but its limit is $0\notin(0,1)$.
\end{example}

Now recall that to create a persistence landscape, we started with a persistence diagram and for each non-diagonal point $(b_i, d_i)$ we constructed a triangular shaped function $$f_{(b_i, d_i)}=
\left
\{
\begin{array}{ll}
x-b_i & \mbox{if }b_i<x\leq \frac{b_i+d_i}{2},\\
-x+d_i & \mbox{if }\frac{b_i+d_i}{2}<x<d_i,\\
0 & \mbox{otherwise.}
\end{array}
\right.$$

We can safely assume that the persistence diagram has finitely many off-diagonal points. Recall that we have a triangle for each of these and in practice there is a lot of overlap. For every positive integer $k$, we have a function $\lambda_k$ defined by $$\lambda_k(x)=kmax\{f_{(b_i, d_i)}(x)|(b_i, d_i)\in P\},$$ where $P$ is our persistence diagram and $kmax$ takes the $k$-th largest value of the overlapping triangles. Set $\lambda_k(x)=0$ if there is no $k$-th largest value.

Now let $\lambda=\{\lambda_k\}_{k\in Z^+}$. This sequence of functions is what we will call a persistence landscape. Persistence landscapes form a vector subspace of the Banach space $L^p(Z^+\times R)$ consisting of sequences of functions of this form whose $L^p$ norm defined below is finite. For the vector space structure, if $\lambda$ and $\eta$ are two landscapes and $c\in R$, then we define $(\lambda+\eta)_k$ to be $\lambda_k$+$\eta_k$ and $(c\lambda)_k=c\lambda_k$ for all $k$. It is a Banach space with the norm $$||\lambda||_p=(\sum_{k=1}^\infty ||\lambda_k||_p^p)^\frac{1}{p},$$ where $||f||_p$ denotes the $L^p$ norm $(\int_R|f|^p)^{1/p}$.

Now we can talk about the mean of persistence landscapes. For example, the mean of $\eta$ and $\lambda$ is $(\eta+\lambda)/2$. This doesn't necessarily correspond to a persistence landscape but it does allow us to derive a law of large numbers and a central limit theorem for persistence landscapes allowing for a statistical analysis that is more natural than if we worked with persistence diagrams. See Bubenik \cite{Bu1} for more details.

My own interest in persistence landscapes and especially the work of Gidea and Katz \cite{GK} is the application to time series analysis. I will go into more details in Sections 5.8 and 5.9.

So now we can review the main steps for analyzing a point cloud: 

Step 1: Look at the points. If they are low dimensional we may be able to see some features we want to capture. Otherwise we can use a standard dimensionality reduction technique like prinicpal componenet analysis or multidimensional scaling.

Step 2: Compute homology at different values of an increasing radius $\epsilon_i$. Do this for some differing low dimensions.

Step 3: Look at a persistence diagram. Points that are further from the main diagonal are more "persistent" and may be the most interesting. As an example, we have points in the shape of a ring. The hole in the middle will last the longest. 

Step 4: Use bottleneck or Wasserstein distance to perform clustering and classification tasks. 

I will describe some additional ideas in the rest of this chapter and the next two. 

So we are left with the question: Is TDA good for anything? It sees to be mostly potential now but in the internet of things classification work that I will describe in Section 5.9, it worked better than other methods we tried for noisy long term data. In \cite{RB} there are a lot of biology examples. For example, they describe work on the evolution of viruses. If two strains of a virus infect a single cell, they can combine into a new strain creating a cycle in their family tree. The size of the cycle corresponds to the persistence of a one dimensional homology class and can be informative in determining a virus's ancestry. There is also a lot of discussion on the use of TDA in cancer research. They particularly cite the success of the related technique, Mapper, in the discovery of a previously unknown line of breast cancer cell. I will discuss that result in Section 6.3.  Gidea and Katz \cite{GK} have a financial application in mind. I will say more about that in Section 5.9. There are a lot more applications in the literature, but I will leave it to subject matter experts to decide how useful they are. 

What I can say for sure is that there is a lot of interest in the subject and it will be easy to conduct experiments as there is free and open source software that can do all of the steps for you. It will be interesting to see where all of this leads in the future.

\section{Sublevel Sets}

In various problems in topology and especially in the field of {\it differential geometry} we deal with curved versions of $R^n$ called {\it manifolds}. Think of the surface of the Earth that looks flat when you are on it but you can see its spherical structure as you move back. We know that gravity bends space and time so differential geometry is the main tool in general relativity. Here is a precise definition.

\begin{definition}
A topological space $M$ is an $n${\it -manifold}\index{manifold} if every point $x\in M$ has an open neighborhood homeomorphic to $R^n$. 
\end{definition}

A curve is a 1-manifold and a surface is a 2-manifold.

Given a manifold $M$, we may want to look at the inverse images of a {height function} $h: M\rightarrow R$. The study of these functions is called {\it Morse theory}\index{Morse theory}. A full description can be found in \cite{Mil1}. A source of persistent homology diagrams is the sublevel sets of this function, ie. $f^{-1}((-\infty, a])$ for $a\in R$. This is best illustrated by some examples.

\begin{figure}[ht]
\begin{center}
  \scalebox{0.4}{\includegraphics{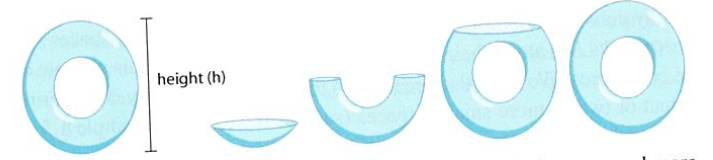}}
\caption{
\rm 
Sublevel Sets of the Height Function of a Torus \cite{RB}
}
\end{center}
\end{figure}

\begin{example}
Figure 5.6.1 \cite{RB} shows the sublevel sets of the height function of a torus that is standing on its side. The four examples on the right side of the figure represent a disk, a (bent) cylinder, a torus with a disk removed, and then a full torus. 
\end{example}

\begin{figure}[ht]
\begin{center}
  \scalebox{0.4}{\includegraphics{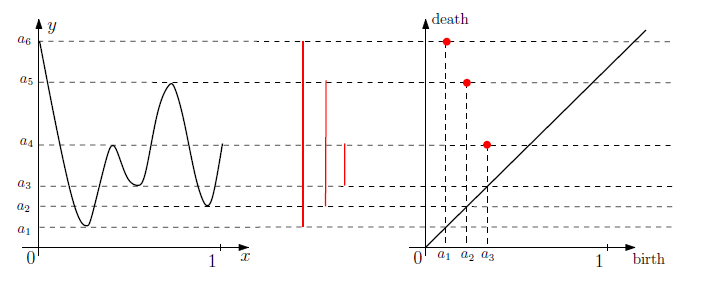}}
\caption{
\rm 
Barcode and Persistence Diagram of a Function $f: [0, 1]\rightarrow R$ \cite{CM}.
}
\end{center}
\end{figure}

\begin{example}
Figure 5.6.2 from \cite{CM} shows a real valued function a single variable. Looking at the sublevel sets, we start at $y=0$ and move upward, looking at the sets $f^{-1}((-\infty, y])$ for increasing values of $y$. The resulting sets are empty until we hit the minimum at $a_1$. Now we have an interval on the $x$-axis, so we have created a new 0-dimensional homology class. A second one is created at $a_2$ which is another minimum. When we hit a function maximum, a class is destroyed by the merging of two intervals.There are merges at $a_4$ and $a_5$. The resulting barcodes and persistence diagrams are in the center and right of the figure.
\end{example}

\begin{figure}[ht]
\begin{center}
  \scalebox{0.4}{\includegraphics{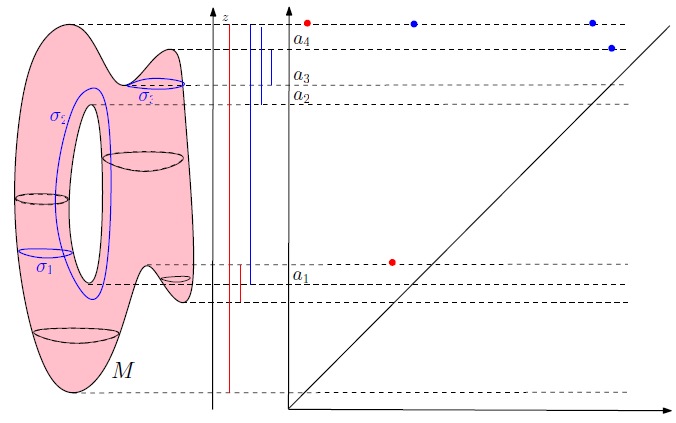}}
\caption{
\rm 
Barcode and Persistence Diagram of a Height function of a surface in $R^3$ \cite{CM}.
}
\end{center}
\end{figure}

\begin{example}
Figure 5.6.3 also comes from \cite{CM}. In this case we have a surface in $R^3$ and the sublevel sets have both 0-dimensional homology (left 2 bars) and 1-dimensional homology (right 3 bars). The corresponding persistence diagram is on the right. If you can see the colors, the red bars and dots are 0-dimensional classes and the blue ones are 1-dimensional homology. Note that this surface is homeomorphic to a torus.
\end{example}

\begin{example}
Suppose we have a surface which is the graph of a function of 2 variables. The sublevel sets form the {\it lawnmower} complex. I once used this as an idea for detecting a large change in a set of geographic points. Here is a potential use case.

The Island of the Lost Elephants is a little known island off the coast of Alaska originally settled by a herd of elephants with very bad navigation skills who had been looking for a shortcut to Antarctica. The island is rectangular and is divided into clearly marked squares of equal size. Every year, the island is visited by a team from the North American Society for the Counting of Elephants (NASCE). They count the number of elephants in each square and plot a three dimensional surface whose height is the number of elephants in each square. The motivation is to find out whether the location of the elephants has significantly changed. Computing a persistence diagram for each year based on the sublevel sets, we can then plot bottleneck or Wasserstein distances between them. If we see a spike or significant level change we can issue an ECLA (elephant change of location alert). This is a better solution than using Euclidean distance because an elephant moving into a neighboring square should not be as significant as a longer move. 
\end{example}

\begin{example}
Suppose we are interested in a set of grayscale images. Each pixel represents a number between 0 and 255. This is similar to the previous case in that we have a surface which is the graph of a function of two variables (pixel value vs  row and column). We can then compute the associated persistence diagram and look at the distance measures. Again, this is preferebale to Euclidean distance as it will be less sensitive to small differences  in the positions of the objects in the image. For a color image, there are three values between 0 and 255 for each pixel. We can simply combine them through an operation such as adding them or taking their average value.,
\end{example}

\begin{example}
The Vietoris-Rips complex of a point cloud is a special case of a sublevel set. Our height function is then the Euclidean distance of each point of $R^n$ to the nearest point of our finite data set where our data points are embedded in $R^n$.
\end{example}

In Chapter 6, I will discuss the data visualization method, Mapper. Mapper also involves looking at the inverse images of a height function and is a variation on the theme of this section. See the description there for details.

\section{Graphs}
Another structure that can be analyzed with persistent homology is a graph. Graphs can also lead to more traditional simplicial complexes. I will start with a quick review of graph terminology.

\subsection{Review of Graph Terminology}

\begin{figure}[ht]
\begin{center}
  \scalebox{0.4}{\includegraphics{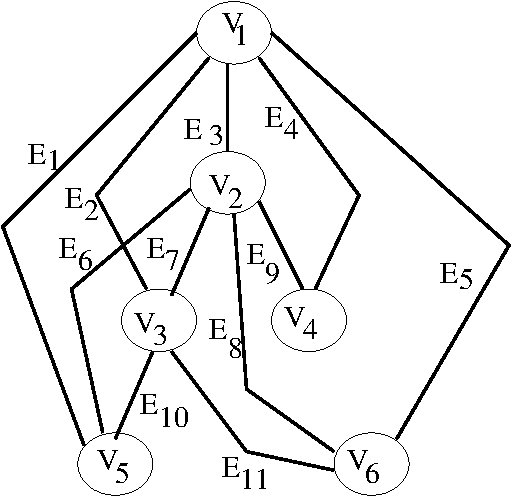}}
\caption{
\rm
Example of a Graph
}
\end{center}
\end{figure}

\begin{definition}
A {\it graph} \index{graph}$G$ is a nonempty set $V(G)$\index{V(G)@$V(G)$} of elements called {\it vertices} \index{vertices}or {\it nodes} \index{nodes}together with a set $E(G)$ \index{E(G)@$E(G)$}of elements called {\it edges} \index{edges} such that for each edge we associate a pair of vertices called its {\it ends}\index{edges!ends}. An edge with identical ends is called a {\it loop}\index{loop}. The 2 ends of an edge are {\it adjacent}\index{vertices!adjacent} or {\it joined} \index{vertices!joined by an edge} by the edge, and we refer to adjacent vertices as {\it neighbors}.\index{neighbors} We say that the edge is {\it incident} \index{incident}to each of its ends.
\end{definition}

Figure 5.7.1  represents a graph with 6 vertices and 11 edges. The edge $E_1$, for example, joins the 2 vertices $V_1$ and $V_5$ so that $V_1$ and $V_5$ are adjacent to each other or $V_1$ is a neighbor of $V_5$ and vice versa. $E_1$ is incident to $V_1$ and $V_5$.

\begin{definition}
A graph $G$ is {\it complete} \index{graph!complete}if every pair of distinct vertices is joined by an edge. It is {\it bipartite} \index{graph!bipartite}if the vertices can be partitioned into 2 sets $X$ and $Y$ such that every edge of $G$ has one end in $X$ and one end in $Y$.
\end{definition}

It is sometimes convenient to represent a graph in terms of a matrix. We define 2 types of matrices which contain all of the information represented in a graph.

\begin{definition}
Given a graph $G$, the {\it adjacency matrix} $A(G)$ \index{adjacency matrix}is defined to be the matrix $(a_{ij})$ whose rows and columns are indexed by the vertices of $G$ and such that $a_{ij}$ is the number of edges joining vertices $i$ and $j$. (Note that this matrix is always symmetric since if an edge joins $i$ and $j$ it also joins $j$ and $i$. The {\it incidence matrix} \index{incidence matrix}$B(G)$ is defined to be the matrix $(b_{ve})$ whose rows are indexed by the vertices of $G$ and whose columns are indexed by the edges of $G$. We let  $b_{ve}=2$ if edge $e$ is a loop at $v$, $b_{ve}=1$ if $e$ joins $v$ to another vertex, and $b_{ve}=0$ otherwise.
\end{definition}

For the graph in the figure, taking the vertices and edges in numerical order we have that the adjacency matrix is 
\[
A(G)=\begin{bmatrix}
0 & 1 & 1 & 1 & 1 & 1\\
1 & 0 & 1 & 1 & 1 & 1\\
1 & 1 & 0 & 0 & 1 & 1\\
1 & 1 & 0 & 0 & 0 & 0\\
1 & 1 & 1 & 0 & 0 & 0\\
1 & 1 & 1 & 0 & 0 & 0
\end{bmatrix}
\]

and the incidence matrix is 

\[
\setcounter{MaxMatrixCols}{20}
B(G)=\begin{bmatrix}
1 & 1 & 1 & 1 & 1 & 0 & 0 & 0 & 0 & 0 & 0\\
0 & 0 & 1 & 0 & 0 & 1 & 1 & 1 & 1 & 0 & 0\\
0 & 1 & 0 & 0 & 0 & 0 & 1 & 0 & 0 & 1 & 1\\
0 & 0 & 0 & 1 & 0 & 0 & 0 & 0 & 1 & 0 & 0\\
1 & 0 & 0 & 0 & 0 & 1 & 0 & 0 & 0 & 1 & 0\\
0 & 0 & 0 & 0 & 1 & 0 & 0 & 1 & 0 & 0 & 1
\end{bmatrix}.
\]

Unless otherwise specified we will always be working with {\it simple} \index{graph!simple}graphs, i.e. graphs with no loops or multiple edges.

\begin{definition}
In a simple graph, the {\it degree} \index{vertices!degree} of a vertex is the number of its neighbors. A graph is {\it d-regular} \index{graph!d-regular} if every vertex has degree $d$, and it is {\it regular} \index{graph!regular}if it is $d$-regular for some $d$.
\end{definition}

For the graph in the figure, $V_1$ and $V_2$ have degree 5, $V_3$ has degree 4, $V_5$ and $V_6$ have degree 3, and $V_4$ has degree 2. Note that in the adjacency matrix, the degree of a vertex is the weight (i.e. number of ones) of its corresponding row and corresponding column. It is also equal to the weight of the corresponding row in the incidence matrix.

\begin{definition}
Let $G$ and $H$ be graphs. If $V(H)\subseteq V(G)$, $E(H)\subseteq E(G)$ and every edge in $H$ has the same pair of ends as it has in $G$, we say that $H$ is a {\it subgraph} \index{subgraph}of $G$ or $G$ is a {\it supergraph} \index{supergraph}of $H$. We say that $H$ {\it spans} \index{graph!spanned by a subgraph}$G$ if $V(H)=V(G)$.
\end{definition}

For example, if $H$ is the subgraph of the graph $G$ in our figure such that 
\[E(H)=\{ E_3, E_7, E_9, E_{10}, E_{11}\},\]
all of the vertices of $G$ are included in the endpoints of these edges, so $H$ is a proper subgraph which spans $G$.

\begin{definition}
A {\it walk} \index{walk}in a graph $G$ is a sequence $W=v_0e_1v_1e_2\ldots e_kv_k$, where the $v_i$ are vertices of $G$, the $e_i$ are edges of $G$, and for $1\leq i\leq k$, $v_{i-1}$ and $v_i$ are the ends of $e_i$. The walk is {\it open} \index{walk!open}if $v_0\neq v_k$ and closed if $v_0=v_k$. The {\it length} \index{walk!length}of $W$ is the number of its edges. The walk is a {\it trail} \index{trail}if the edges are all distinct and a {\it path} \index{path}if the vertices are all distinct. A closed trail of positive length whose vertices (apart from its ends $v_0$ and $v_k$) are distinct is called a {\it cycle}\index{cycle}.
\end{definition}

For example, in our graph $G$, $V_1E_1V_5E_{10}V_3E_2V_1$ is a cycle. $G$ contains many examples of walks, paths, trails, and cycles. The reader should look for some examples.

\begin{definition}
A graph is {\it connected} \index{graph!connected}if any 2 of its vertices are connected by a path. A {\it component} \index{graph!component of}of a graph is a connected subgraph which is maximal in the sense that it is not properly contained in any larger connected subgraph.
\end{definition}

\begin{definition}
A graph is {\it acyclic} \index{graph!acyclic}or a {\it forest} \index{forest}if it contains no cycle. A {\it tree} \index{tree}is an acyclic connected graph. A vertex in a tree of degree one is called a {\it leaf}. \index{leaf}A {\it spanning tree} \index{spanning tree}of a graph $G$ is a subgraph $H$ which is a tree and spans $G$.
\end{definition}

We saw that the subgraph $H$ with edge set 
\[E(H)=\{ E_3, E_7, E_9, E_{10}, E_{11}\}\]
spans our graph $G$. Figure 5.7.2 shows this subgraph.

\begin{figure}[ht]
\begin{center}
  \scalebox{0.4}{\includegraphics{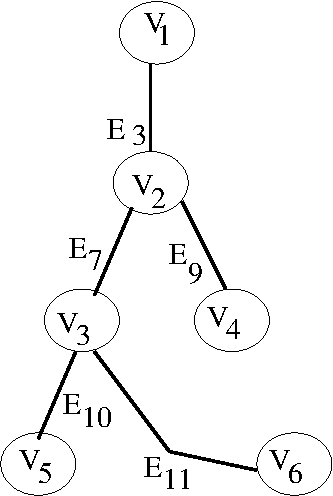}}
\caption{
\rm
Spanning Tree
}
\end{center}
\end{figure}
 
It is now easily seen that this graph is a tree, so it is a spanning tree. There are a number of others. The leaves of this tree are $V_1, V_4, V_5, V_6$. 

\begin{theorem}
If $H$ is a spanning tree for $G$ then the number of edges of $H$ is $|E(H)|=|V(G)|-1=|V(H)|-1.$
\end{theorem}   

So far, I have only been discussing {\it undirected graphs} \index{graph!undirected}in which the edge does not have a preferred direction. I now give the following definition:

\begin{definition}
A {\it directed graph} \index{graph!directed}or {\it digraph} \index{digraph}$D$ is a graph $G$ in which each edge is assigned a direction with the starting end called its {\it tail} \index{tail}and the finishing end called its {\it head}\index{head}. 
\end{definition}

Figure 5.7.3 shows a directed graph with the same edge set and vertex set as the graph in Figure 5.1.

\begin{figure}[ht]
\begin{center}
  \scalebox{0.4}{\includegraphics{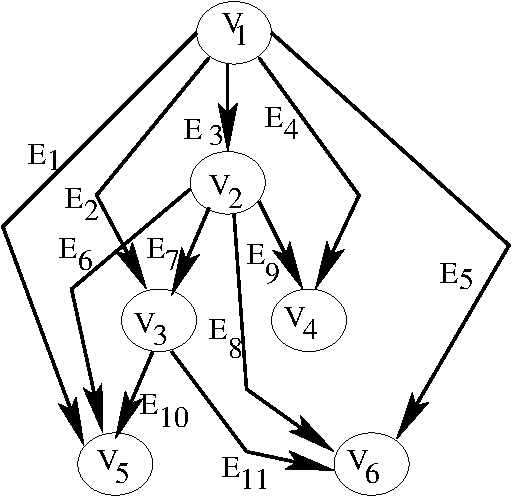}}
\caption{
\rm
Directed Graph
}
\end{center}
\end{figure}

\begin{definition}
The {\it associated digraph} $D(G)$ \index{digraph!associated}of an undirected graph $G$ is the digraph obtained from $G$ by replacing each edge of $G$ by 2 edges pointing in opposite directions. The {\it underlying graph} \index{graph!underlying}$G(D)$ of a digraph $D$ is the graph obtained from $D$ by ignoring the directions of its edges.
\end{definition}

Note that $G(D(G))\neq G$ and $D(G(D))\neq D$. Figure 5.7.1 is the underlying graph of Figure 5.7.3, but Figure 5.7.3 is not the associated digraph for Figure 5.7.1.

I conclude this section with some analogous definitions for digraphs to those given above for undirected graphs.

\begin{definition}
Given a digraph $D$, the {\it adjacency matrix} $A(D)$ is defined to be the matrix $(a_{ij})$ whose rows and columns are indexed by the vertices of $D$ and such that $a_{ij}$ is the number of edges with tail $i$ and head $j$. (Note that this matrix is not generally symmetric.)
\end{definition}

For the graph in Figure 5.7.3, taking the vertices and edges in numerical order we have that the adjacency matrix is 
\[
A(G)=\begin{bmatrix}
0 & 1 & 1 & 1 & 1 & 1\\
0 & 0 & 1 & 1 & 1 & 1\\
0 & 0 & 0 & 0 & 1 & 1\\
0 & 0 & 0 & 0 & 0 & 0\\
0 & 0 & 0 & 0 & 0 & 0\\
0 & 0 & 0 & 0 & 0 & 0
\end{bmatrix}
\]

\begin{definition}
In a digraph $D$, the {\it outdegree} \index{vertices!outdegree} of a vertex $v$ is number of edges of $D$ whose tail is $v$. The {\it indegree} \index{vertices!indegree}of a $v$ is number of edges of $D$ whose head is $v$. The sum of the outdegree and the indegree of $v$ is the {\it total degree} of \index{vertices!total degree}$v$. This is equal to the degree of $v$ in the underlying graph $G(D)$. A {\it sink} \index{sink}is a vertex of outdegree 0, and a {\it source} \index{source}is a vertex of indegree 0.
\end{definition}

The terms {\it sink} and {\it source} are actually derived from fluid mechanics. Fluids are sucked into sinks and spray out of sources. The digraph in Figure 5.7.3 has $V_1$ as its only source and sinks $V_4, V_5,$ and $V_6$.

The indegree and outdegree of each vertex can be read off from the adjacency matrix. The outdegree is the weight of the corresponding row, and the indegree is the weight of the corresponding column. For example, $V_2$ has indegree 1 and outdegree 4.

\begin{definition}
A {\it directed walk} \index{walk!directed}in a digraph $D$ is a sequence $W=v_0e_1v_1e_2\ldots e_kv_k$, where the $v_i$ are vertices of $G$, the $e_i$ are edges of $G$, and for $1\leq i\leq k$, $v_{i-1}$ is the tail of $e_i$, and $v_i$ is the head of $e_i$. Directed trails, paths and cycles are defined in an analogous way.
\end{definition}

Note that the digraph in Figure 5.7.3 has no directed cycles.

Finally, I will define a weighted graph. 

\begin{definition}
A {\it weighted graph} $G$ has a weighting function $w: E(G)\rightarrow R$ which assigns a real number called the {\it weight}\index{graph!weighted} to each edge of the graph. 
\end{definition}

Note that both undirected and directed graphs can be weighted. In either case we modify the defnition of the adjecency matrix by replacing each one with the weight of the corresponding edge. 

\subsection{Graph Distance Measures}
In areas such as network change detection or graph classification, we would like to have a good way to measure a {\it distance} between two graphs. A good reference on this subject with lots of interesting ideas is the book of Wallis, et. al. \cite{WBDK}. If two graphs have a similar set of nodes a quick and easy metric to compute is {\it graph edit distance}\index{graph edit distance}. It has served me well on many network change detection problems.

\begin{definition}
Let $G_1$ and $G_2$ be two graphs. Suppose for $i=1, 2$, $|V_i$ and $E_i$ represent the number of vertices and edges respectively in graph $G_i$, and $|V_1\cap V_2|$ and $|E_1\cap E_2|$ are the number of vertices and edges respectively that the two graphs have in common. Then the {\it (unweighted) graph edit distance} is $$d(G_1, G_2)=|V_1|+|V_2|-2|V_1\cap V_2|+|E_1|+|E_2|-2|E_1\cap E_2|.$$ If the graphs are weighted we can add an additional term. Let $\{e_j\}$ be the set of common edges and let $\beta_{j1}$ and $\beta_{j2}$ be the respective weights of edge $e_j$ in graphs $G_1$ and $G_2$. Then the weighted edit distance is  $$d(G_1, G_2)=|V_1|+|V_2|-2|V_1\cap V_2|+|E_1|+|E_2|-2|E_1\cap E_2|+\alpha\sum_j\frac{|\beta_{j1}-\beta_{j2}|}{\max(|\beta_{j1}, \beta_{j2})},$$ where $\alpha$ is a scaling factor which is usually set to one.
\end{definition}

Note that if $G_1=G_2$ then $|V_1|=|V_2|=|V_1\cap V_2|$ and $|E_1|=|E_2|=|E_1\cap E_2|$, so  $$d(G_1, G_2)=|V_1|+|V_2|-2|V_1\cap V_2|+|E_1|+|E_2|-2|E_1\cap E_2|=0.$$

\begin{figure}[ht]
\begin{center}
  \scalebox{0.4}{\includegraphics{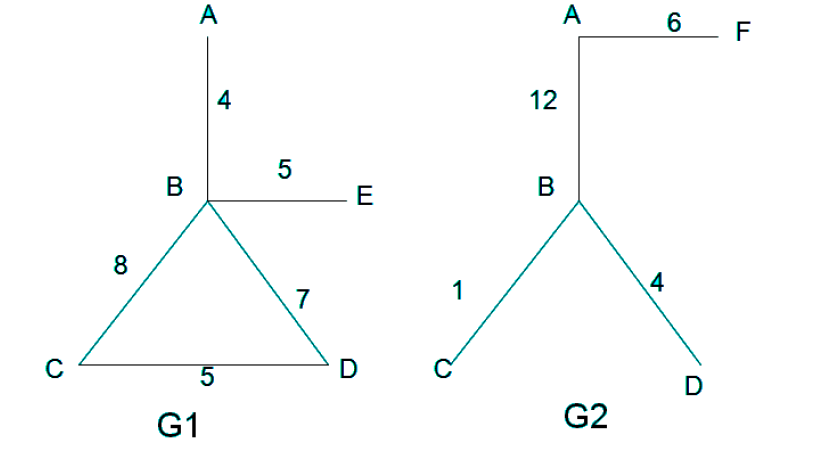}}
\caption{
\rm
Graph Edit Distance Example
}
\end{center}
\end{figure}

\begin{example}
Consider the two graphs in Figure 5.7.4. There are 5 vertices in each of the two graphs and they have 4 of them in common. $G_1$ has 5 edges and $G_2$ has 4 edges and the two graphs have 3 edges in common. The unweighted graph edit distance is $$d(G_1, G_2)=5+5-2(4)+5+4-2(3)=5.$$ If we look at weights there are three common edges: $AB$, $BC$, and $BD$. For $AB$, we have $|4-12|=8$ and $\max(4, 12)=12.$  For $BC$, we have $|8-1|=7$ and $\max(8, 1)=8$. For $BD$, we have $|7-4|=3$, and $\max(7,3)=7$. So to compute the weighted graph edit distance for $\alpha=1$, we add the unweighted distance of 5 to the terms in the summation giving $$d(G_1, G_2)=5+\frac{8}{12}+\frac{7}{8}+\frac{3}{7}=5+\frac{112+147+72}{168}\approx 6.97.$$
\end{example}

In my own work I have gotten better results by dropping the weight factor even if the graph is weighted. 

\subsection{Simplicial Complexes from Graphs}

Graphs give rise to several natural simplicial complexes. I will start with those. Then I will talk about how they give rise to persistence diagrams.

Recall that in an abstract simplicial complex we start with a collection of objects called vertices. A simplex is a subset of the vertices that will be included in the complex. If the subset is of order $n$ when we have an $n-1$-simplex. A {\it  maximal simplex} is a simplex which is not a subset of any other simplex. So in order to define an abstract simplicial complex, it suffices to define its vertices and maximal simplices. The complex then consists of the maximal simplices and all of their subsets. 

Before I list these complexes, I will give one more definition that will be helpful in understanding the third example.

\begin{definition} 
A {\it partially ordered set}\index{partially ordered set} or {\it poset}\index{poset} is a set $X$ with a binary relation $\leq$ called a {\it partial order} with the following three properties:\begin{enumerate}
\item {\bf Reflexivity:} If $a\in X$, then $a\leq a$.
\item {\bf Antisymmetry:} If $a, b\in X$ then $a\leq b$ and $b\leq a$ implies $a=b$.
\item {\bf Transitivity:} If $a, b, c\in X$, then $a\leq b$ and $b\leq c$ implies $a\leq c$.
\end{enumerate}
Note that in a partially ordered set it is possible that for $a, b\in X$ that $a\leq b$ and $b\leq a$ are both false. In this case we say that $a$ and $b$ are {\it incomparable}. 
\end{definition}

An easy example of a partially ordered set is a collection of sets where $A\leq B$ if $A\subseteq B$.

Now here are five simplicial complexes you can build out of a graph:

\begin{enumerate}
\item The graph itself. This is a complex of dimension one. We can think of any higher dimensional complex as strategically throwing information away. This is useful if we are trying to detect only major changes in a graph.
\item {\bf Complete Subgraph Complex:} The vertices are the graph edges and the maximal simplices are the edges sets forming a complete subgraph. Note that a complete subgraph is also called a {\it clique}.
\item{\bf Neighborhood Complex:} The vertices are the graph vertices and the maximal simplices consist each vertex along with their neighbors.
\item{\bf Poset Complex:} The vertices are the graph vertices and the maximal simplices are the maximal directed paths. Note this is defined only in a directed graph. We can also do something similar for any finte partially ordered set.
\item{\bf Matching Complex:} The vertices are the graph edges and the simplices are the edge sets forming a {\it matching}. A matching is a set of edges in which no two are adjacent.
\item{\bf Broken Circuit Complex:} The vertices are the graph edges and the maximal simplices are the are maximal connected edge sets containing no circuit. (These are the spanning trees for a connected graph.) Such edge sets form the independent sets of a matroid (a generalization of linear independence) and have particularly simple homology. For details see \cite{Bjo1}.
\end{enumerate}

Here are some potential features for a machine learning algorithm:
\begin{enumerate}
\item Dimension of the complex.
\item Number of simplices in each dimension.
\item Euler characteristic. 
\item Betti numbers of homology groups
\end{enumerate}

Finally, suppose we have a weighted graph. We can make a Vietoris-Rips complex and the corresponding persistence diagram. For each $\epsilon$, we include a set of vertices in a simplex if for each pair of vertices, there is an path between them of weight less than or equal to $\epsilon$. (The weight of a path is the sum of the weight of its edges.) Then we can make bar codes, persistence diagrams, and persistence landscapes just as in the case of point clouds.

\section{Time Series}

A {\it time series}\index{time series}  $\{x_1, x_2, x_3, \cdots\}$ is a sequence of measurements at successive time intervals. For example, a time series could represent the price of a stock on each day or the temperature every hour. We can also have a multidimensional times series where we take more than one measurement simultaneously. For example. you might keep track of your height and weight at noon Greenwich Mean Time on your birthday each year. This is a two dimensional time series.

There is a highly developed theory of time series in economic modeling. If we expect the measurement for each time period to be a linear combination of the measurements for the  last $p$ times for some $p$ (with some random noise added) we have an {\it autoregressive (AR) model}. A {\it moving average (MA) model} assumes the value for a time period is a linear combination of some average value and $q$ white noise terms for some $q$. An ARMA model is the sum of an AR and a MA model. 

There are also more complicated models such as ARIMA in which first (or higher order) differences are involved, and ARCH in which the variance is allowed to vary. What all of these have in common is that our goal is to determine the unknown parameters and be able to make predictions. 

I won't get into details here as time series modeling is a whole subject of its own. The  classic book on this type of modeling is the book by Box, Jenkins, and Reinsel \cite{BJR}. A more modern book is Hamilton \cite{Ham1}. Two books with an economics focus are \cite{End1} and \cite{Tsay}. Finally, Cryer and Chan \cite{CC} is a full description of the time series analysis (TSA) package in R.

I have found that traditional time series modeling is not readily applicable to cyber defense applications. Economists are trying to model what is normal, but I was more interested in anomalies. Also, in the applications I looked at, there was not a good justification for thinking that the quantities I was measuring were going to be linear combinations of the measurements at previous times. 

Another issue is finding a distance measure between two time series. The obvious one is Euclidean distance in which you add up the squared distances of the values of the the two time series at each time period. But this does not take into account that the series might just be shifts of each other or differ in scaling. A solution to this is {\it dynamic time warping}\index{dynamic time warping} or {\it DTW}. In DTW, we match indices from the first series to indices from the second. An index in one series is allowed to match more than one index in the other and indices can be skipped as long as we satisfy the following condition: If index $i$ is matched to index $f(i)$ then $i>j$ implies $f(i)\geq f(j)$. We then take take the match in which $\sum_i |i-f(i)|$ is minimized. See \cite{BC} for details of how this is implemented. This approach can be rather slow, though, even with more modern variants that are faster.

Another issue involves finding anomalies in time series. It is not hard to find spikes or unusal values. One way to do it is to take a sliding window of $n$ samples $\{x_i, x_{i+1}, \cdots, x_{i+n-1}\}$. Let $\mu$ be the mean of these samples and $\sigma$ the standard deviation. Then we say that $x_{i+n}$ is an anomalous value if $\frac{|x_{i+n}-\mu|}{\sigma}$ is above some threshold. But what if we instead want to find an anomalous substring? For example, an EKG has obvious spikes but we don't want to say that a flat line would be normal since we have eliminated the spikes. 

A good solution to all of these problems was to use the time series discretization method, SAX \cite{LKWL}. I will tell the story of how I used SAX to build a time series analysis tool that performs classification and anomaly detection in the next section. But what does any of this have to do with topological data analysis? TDA was the key to dealing with multivariate time series. The result was a successful effort to classify internet of things devices based on network traffic \cite{DPSW}.

\section{SAX and Multivariate Time Series}

While trying to find a different approach to time series analysis, I found an interesting paper. It was a Master's Thesis in Statistics by Caroline Kleist of Humboldt Universit\"{a}t in Berlin \cite{Kle1} that provides a survey of data mining techniques when the data points are time series. The paper had multiple references to Jessica Lin at George Mason University who had formulated a time series analysis method called SAX in which the series is represented by a string of symbols from a small alphabet. It turns out that one of my previous offices was already working with George Mason University and I started a project with Lin's colleague Robert Simon to apply SAX to cyber defense problems. Our internet of things paper \cite{DPSW} was the result of this project. We also developed a software tool called TSAT. TSAT is a JAVA program which runs on either Windows or Linux.

In this section, I will break down the steps of our analysis. In Section 5.9.1, I will introduce SAX and describe its advantages. In Section 5.9.2, I will discuss the SEQUITUR algorithm which combines symbols into {\it grammar rules} and gives SAX the ability to detect anomalies. Next, I will describe the Representative Pattern Mining or RPM algorithm which we use to classify time series. In Section 5.9.4, I wil give the method of Gidea and Katz for converting a multivariate time series into a univariate one using topological data analysis. Then I will describe our internet of things results. Finally, in Section 5.9.6 I will describe some work of Elizabeth Munch at Michigan State University on TDA and time series. It is hoped that this work will help us improve our current algorithms.

\subsection{SAX}
In this section, I will briefly outline the SAX algorithm. For more information, see \cite{LKWL} and its references.

SAX stands for Symbolic Aggregate approXimation. (The X comes from the middle of the last word.) The algorithm replaces a sequence of real valued measurements with a string of symbols from an alphabet size which is typically 3-12 symbols. Doing this allows for extremely fast dstance computations as there are a small number of pairs of symbols and their distances can be kept in a look up table. This greatly speeds up both Euclidean distance computation and dynamic time warping. In addition, SAX has the following applications: \begin{enumerate}
\item Finding {\it motifs} or common patterns within a series.
\item Finding {\it surprising} subsequences representing anomalies. 
\item Finding {\it contrast sets}. Given two sets of time series, are there features which distinguish one set from the other?
\item The RPM (Representative Pattern Mining) Algorithm uses SAX to find strings which distinguish one set of time series from another set. 
\end{enumerate}

The first step in the process is optional and referred to as {\it Piecewise Aggregate Approximation} or {\it PAA}. PAA reduces the number of terms in a time series by dividing it into $w$ equal sized {\it frames.} If $n$ is the length of the time series, we will assume in our example that $w$ to be a divisor of $n$. For each frame, we assign the value to the mean value of time steps in that frame. This converts the series into a series of length $n/w$. For example, if $n=128$, and $w=16$, we get the time series of length 8 shown in Figure 5.9.1.

\begin{figure}[ht]
\begin{center}
  \scalebox{0.4}{\includegraphics{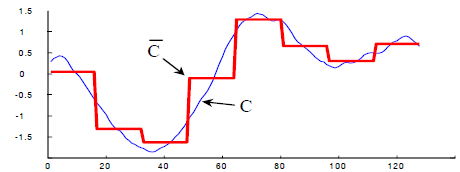}}
\caption{
\rm
Piecewise Aggregate Approximation (PAA) Example \cite{LKWL}
}
\end{center}
\end{figure}

The next set is to assign a symbol from our alphabet. To do this, let $\mu$ be the mean and $\sigma$ be the standard deviation of the time series resulting from PAA. For each $x_i$ in this series, we {\it normalize} it by replacing it with $\frac{x_i-\mu}{\sigma}$. Then pretend that the values come form a standard normal distribution. This is not necessarily true, but we at least want the symbols to be roughly equally likely. Divide the range of the standard normal into $\alpha$ parts of equal probability where $\alpha$ is the size of the alphabet we want to use. 

\begin{figure}[ht]
\begin{center}
  \scalebox{0.4}{\includegraphics{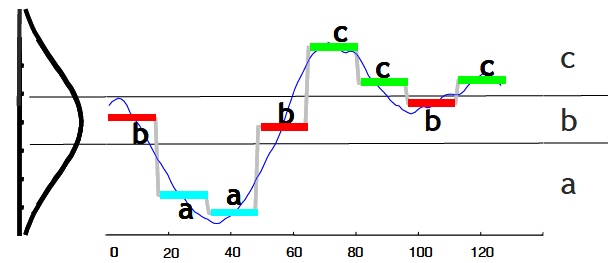}}
\caption{
\rm
SAX Discretization Example: $\alpha=3$ \cite{LKWL}
}
\end{center}
\end{figure}

In Figure 5.9.2 we have an example where the alphabet is of size 3 and consists of the symbols $a$, $b$, and $c$. For the standard normal distribution, $$P(x<-.43)=P(-.43<x<.43)=P(x>.43)=\frac{1}{3}.$$ So for each normalized value $x$ in our series, we replace it with $a$ if $x<-.43$, with $c$ if $x>.43$, and $b$ otherwise. We have now converted our time series into the string of symbols, $baabccbc$.

 To detect anomalies and classify series, the symbols are converted into rules using a {\it context free grammar}.  The rough idea is as follows:

Suppose we receive the string $abcabcaaaabcabc$. We see that $abc$ is common, so let $R_1=abc$ be a grammar rule. We now have $R_1R_1aaaR_1R_1$. Now let $R_2=R_1R_1$, giving $R_2aaaR_2$. Any substrings left over at the end are anomalous. In addition, we can classify 2 sets of time series by finding grammar rules which are common in one set but rare in another.

I did cheat a little though. We want to have a systematic algorithm to assign a set of rules. That is our next subject.

\subsection{SEQUITUR} 
In SAX, grammar rules are assigned using the {\it SEQUITUR algorithm}\index{SEQUITUR algorithm} of Nevill-Manning and Witten. (For computer scientists, SEQUITUR is an algorithm which computes a {\it context free grammar}.) I will outline their algorithm here. See their paper \cite{NW} for details.

SEQUITUR infers a hierarchical structure from a sequence of discrete symbols. The rules are new symbols that replaces strings of existing symbols. Rules then include other rules leading to the hierarchy. We are going to impose two constraints:\begin{enumerate}
\item {\bf Digram uniqueness:} No pair of adjacent symbols appears more than once in the grammar.
\item {\bf Rule utility:} Every rule is used more than once.
\end{enumerate}

The next four examples refer to Figure 5.9.3 taken from \cite{NW}.

\begin{figure}[ht]
\begin{center}
  \scalebox{0.4}{\includegraphics{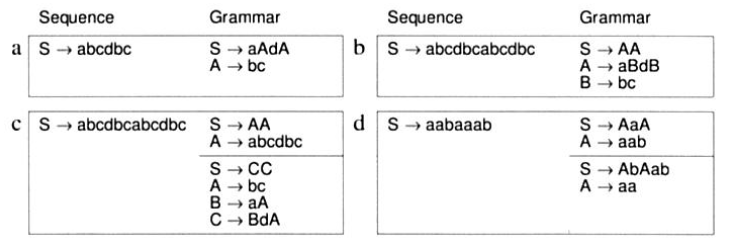}}
\caption{
\rm
Grammar Rule Examples \cite{NW}
}
\end{center}
\end{figure}

\begin{example}
The first example introduces their terminology. In box a, the original rule is the sequence $S\rightarrow abcdbc$. Since $bc$ repeats we replace it with the rule $A\rightarrow bc$. Now we get the sequence $S\rightarrow aAdA$. 
\end{example}

\begin{example}
Box b shows how a rule can be reused in longer rules. We start with the sequence $S\rightarrow abcdbcabcdbc$. If we impose the rule $B\rightarrow bc$ and then $A\rightarrow aBdB$, we get the sequence $S\rightarrow AA$. This gives a rule within a rule showing a heirarchical structure. 
\end{example}

\begin{example}
So far we haven't thought about our constraints. In box c, we have two possible substitutions for the original sequence $S$. In the top one, digram uniqueness is violated as we have $bc$ appearing twice in rule $A$, so we have introduced redundancy. In the bottom one, rule utility is violated as $B$ is used only once. We would have a more concise grammar if it was removed. 
\end{example}

\begin{example}
Box d shows that both constrains can be satisfied but the rules are not necessarily unique. Both the top and bottom are acceptable and satisfy both constraints.  
\end{example}

To show how the algorithm proceeds and enforces the constraints, Neville-Manning and Witten use a table which I have reproduced as Figure 5.9.4. 

\begin{figure}[ht]
\begin{center}
  \scalebox{0.4}{\includegraphics{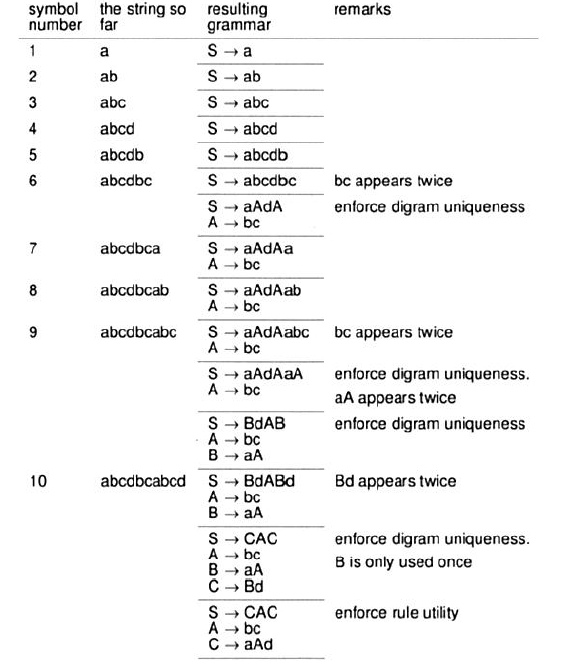}}
\caption{
\rm
SEQUITUR Example \cite{NW}
}
\end{center}
\end{figure}

SEQUITUR reads a sequence from left to right. If a digram appears more than once, a rule is created, imposing digram uniqueness. If a rule is only used once, it is eliminated, imposing rule utility. Reading down the table, the first interesting event is when symbol 6 is added. Now we have $bc$ appearing twice, so we create the rule $A\rightarrow bc$ and make the substitution in our sequence. At symbol 9,we have $bc$ appearing again, so we replace it with $A$, getting the sequence $aAdAaA$. But now $aA$ appears twice, so we create the new rule $B\rightarrow aA$. This gives the sequence $BdAB$.

Now at symbol 10, we get $BdABd$. Now $Bd$ appears twice, so we create $C\rightarrow Bd$. So the sequence is $CAC$. Now we see that the only use for $B$ was to define $C$, so we are violating rule utility. Our sequence doesn't change if we eliminate rule $B$ and just let $C\rightarrow aAd$. That is how we get longer rules. 

The remainder of \cite{NW} shows that the algorithm is linear in both time in space. The interest for SAX is that the rules represent common substrings in a sequence. Any symbols that are not covered by these rules are by definition anomalous. We can also classify time series by looking at rules that are common for one type of series and rare for the other. That is the subject of our next section.

\subsection{Representative Pattern Mining}

Representative Pattern Mining \cite{WLSOGBCF} is the algorithm that TSAT uses to classify time series. The idea is to find {\it representative patterns}, i.e. subsequences that are common in one class of time series and rare in the other. I will briefly summarize the steps and refer you to the paper for details. 

The algorithm consists of two steps: Learning  representative patterns during a training stage and then using them for classification. The training stage consists of preprocessing the data, generating representative pattern candidates, and then selecting the most representative candidates. Then classification is done by letting each pattern be a feature and computing a distance from a time series we want to classify to each pattern. Support Vector Machines are used in \cite{WLSOGBCF} as a classifier, but in TSAT, we used random forests instead.

The preprocessing step in the training phase consists of turning the time series into symbols using SAX and then using SEQUITUR to generate grammar rules. The algorithm keeps track of the position of the rules in the sequence. The algorithm eliminates patterns that are too similar to each other through a clustering process and eventually chooses a set that will be the best discriminator of the two classes. 

SAX extracts subsequences using a sliding window, and a PAA frame size and alphabet size must be chosen. In TSAT, a user can choose these parameters or the algorithm can choose them for the user by finding the parameter set that will do the best job of discriminating the training data. See \cite{WLSOGBCF} for details.

\subsection{Converting Multivariate Time Series to the Univariate Case or Predicting Economic Collapse}

So what does any of this have to do with algebraic topology? The answer comes when we try to use SAX to classify multivariate time series.

One approach might be to take each coordinate individually. But each coordinate could individually seem normal but taken together could represent an anomaly. For example, if you represented a person's height and weight over time, it might not be unusual that they weigh 100 pounds or are 6 feet tall but being both would be pretty unusual. TDA allows us to do the conversion without decoupling coordinates.

We used the method of Gidea and Katz \cite{GK}. Suppose our time series consisted of $m$ measurements at each time period. Then letting $X=\{x_1, x_2, \cdots\}$, each $x_i$ represents a point in $R^m$. Choose a window size $w$. Then slide the window across the series. At position $k$, our window is the set $\{x_k, x_{k+1}, \cdots, x_{k+w-1}\}$. We then have a cloud of $w$ points in $R^m$ and we can compute the persistence landscape of these points. If we step by one time unit, we now have a cloud of points $\{x_{k+1}, x_{k+2}, \cdots, x_{k+w}\}$. Now letting $\lambda_1$ and $\lambda_2$ be the two landscapes, we compute the distance $||\lambda_1-\lambda_2||_p$ as defined in Section 5.5, where we generally let $p=1$ or $p=2$. This gives a real number for the difference between each window setting, so we have a univariate time series and we can use SAX and RPM for classification as before.

Gidea and Katz used this technique for forecasting economic downturns. They started with the daily time series of four US stock market indices: S\&P 500, DJIA, NASDAQ, and Russell 2000 from December 23,1987 to December 8, 2016 (7301 trading days). This is a series of points in $R^4$. They used a sliding window of size $w=50$ and $w=100$. Two interesting dates were chosen, the dot com  crash on March 10, 2000 and the Lehman bankruptcy on September 15, 2008. Figure 5.9.5 shows in top to bottom order the normalized S\&P series, the normalized $L^1$ norm of the persistence landscapes calculated with a sliding window of 50 days, and the volatility index (VIX) which measures the volatility of the S\&P index over the next 30 days and is often used for forecasting. The left side represents the 1000 trading days before the dot com crash and the right side represent the 1000 trading days before the Lehman bankruptcy. Looking at the middle row, we see the dramatic behavior of the persistence landscapes in the days before the two crashes.

\begin{figure}[ht]
\begin{center}
  \scalebox{0.4}{\includegraphics{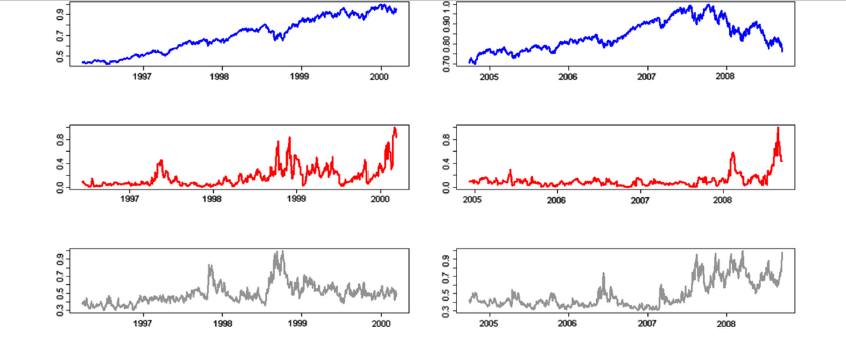}}
\caption{
\rm
Behavior of Three Indices Immediately Before Two Economic Crises. \cite{GK}
}
\end{center}
\end{figure}

I will briefly mention another Gidea and Katz paper that includes three additional authors \cite{GGKRS}. This paper analyzed four cryptocurrencies: Bitcoin, Etherium, Litecoin, and Ripple. The problem was to find a warning of collapse of the value of these currencies. The method was related to the previous paper with a few differences. For each of these currencies, start with a univariate time series representing the log-returns. The input for TDA is a coordinate embedding $d=4$ and a sliding window of size $w=50$. So if $\{x_i\}$ is the original time series, let $z_0=(x_0, x_1, x_2, x_3), z_1=(x_1, x_2, x_3, x_4), \cdots$. Then our sliding windows are of the form $w_i=\{z_i, z_{i+1}, \cdots, z_{i+49}\}$. We consider each window to be a point cloud, compute it's persistence landscape, and then compute the $L^1$ norm $$||\lambda||_1=(\sum_{k=1}^\infty ||\lambda_k||_1),$$ where $\lambda_k$ is as defined in Section 5.5. 

Now perform a $k$-means clustering to find unusual dates in the series. The input is a series of points $(x, y, z)\in R^3$ where $x$ is the log of the price of the asset, $y$ is the log-return of the asset, and $z$ is the $L^1$ norm produced by TDA. For each of these currencies, there were clusters consisting of points representing dates leading up to a crash. See \cite{GGKRS} for more details.

With these techniques, we were in a good position to use multivariate time series to classify internet of things devices.

\subsection{Classification of Internet of Things Data}

The material in this section comes from my own paper \cite{DPSW}. The idea was to classify internet of things (IoT) devices using data that would still be visible under encryption such as statistics on packet size and interarrival time. The experiments were performed on a proprietary testbed consisting of 183 devices produced by a variety of manufacturers. The devices were placed in multiple rooms and corridors designed to emulate how they would be used in commercial and residential environments. 

All network traffic over a 9 month period was collected and each recieved packet was time-stamped. From the logs, we could look at source and destination MAC and IP addresses,transport layer ports, and the device type for each payload. This gave us a labeled data set and for each conversation we could compute the total number of bytes sent in a time window and the mean and variance of interarrival time. We experimented with both single and multi-attribute time series. For the univariate series we processed the data using SAX, and used Representative Pattern Matching to perform classification. For the multivariate series, we used both observed and derived attributes. We used the TDA method of Gidea and Katz to transform a multivariate series into a univariate one and then used SAX and RPM for classification. For the paper, we used the Dionysus package \cite{Mor1}, but later sped up the computations dramatically with the Ripser package \cite{Bau1}. Ripser is now included in TSAT. As an alternate algorithm for classifying multivariate time series is WEASEL+MUSE \cite{SL} which looks at patterns in each coordinate separately and includes the coordinate as part of their label. 

We looked at three types of IoT devices:\begin{enumerate}
\item Cameras including audio speakers and characterized by high-volume burst data.
\item Sensors such as environmental sensors.
\item Multi-purpose devices such as tablets or certain cameras that allowed for two-way and streaming audio,
\end{enumerate}

The collection was noisy and incomplete. There were interruptions due to system maintenance and equipment upgrades and much higher traffic on weekdays during work hours than at night, on weekends, and during holiday breaks. For single attribute time series, the best results were for a day of training data and a day of test data. The attribute that worked best was total number of bytes followed by packet interarrival time.

For the multivariate time series, WEASEL+MUSE worked as well as TDA for small amounts of data. When we tried increasing the size of the test data, TDA worked much better than other methods. Table 5.9.1 shows the results of TDA on one month of training and 8 months of testing data. Considering the amount of noise, these results are quite good. They are evaluated using the F1 score and the Matthews Correlation Coefficient (MCC). Letting TP, FP, TN, and FN be true positives, false postives, true negatives, and false negatives respectively, we have $$F1=\frac{2TP}{2TP+FP+FN}$$ and $$MCC=\frac{TP\times TN-FP\times FN}{\sqrt{(TP+FP)(TP+FN)(TN+FP)(TN+FN)}}.$$ Note that the F1 score ranges from 0 to 1 while MCC ranges from -1 to 1.

\begin{table} 
\begin{center}
\begin{tabular}{|c|c|c|}
\hline
Device Type  & F1 & MCC\\
\hline\hline
Cameras & .687 & .535\\
\hline
Sensors & .835 & .765\\
\hline
Multi-purpose & .8 & .69\\
\hline
Weighted & .77 & .66\\
\hline
\end{tabular}
\caption{Results for TDA Multivariate Case.}
\end{center}
\end{table}

To conclude, while univariate time series and more traditional methods worked well for shorter time periods, TDA outperformed the other methods for a noisy multi-month period. See \cite{DPSW} for more details.

\subsection{More on Time Series and Topological Data Analysis}

In this section, I will discuss two additional papers that built on the work of the last section. Also, there are two remaining practical questions:
\begin{enumerate}
\item How do we pick the size of the window which determines the size of the point clouds used to compute persistence landscapes?
\item In our work we slid the window by one time step at a time. Is this always the right thing to do? How do we know?
\end{enumerate}

These are still open problems but I will conclude by discussing a paper \cite{MMK} that was suggested to me by Elizabeth Munch of Michigan State University to begin to answer these questions. 

First, though, I will summarize a paper \cite{CICC} which cites my IoT work and presents an alternative topology based algorithm.The paper classifies IoT devices using the single attribute of interarrival time. The advantage, though, is that they can work in smaller time windows and they create a simpler simplical complex.

\begin{figure}[ht]
\begin{center}
  \scalebox{0.4}{\includegraphics{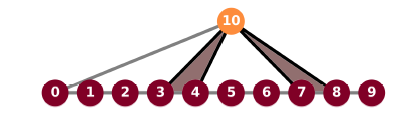}}
\caption{
\rm
Clique Complex of Window $\omega_{10}(0)$ After First 4 Nonzero Edges are Added \cite{CICC}
}
\end{center}
\end{figure}

The complex is best described using a picture. Figure 5.9.6 is taken from \cite{CICC}. We are given a time series of packets $\{p_0, p_1, \cdots, p_n\}$. If packet $p_i$ arrives at time $T(p_i)$, we define the interarrival time $\Delta_T(p_i)$ to be $T(p_i)-T(p_{i-1}).$ Now choose a window size $k$ and we classify windows $\omega_k(i)=(p_i, p_{i+1}, \cdots, p_{i+k-1})$ as to the type of IoT device that produced them. In this paper, they chose $k=25$.

Now build a filtration of complexes as follows. Figure 5.9.6 shows a window of size 10 with starting point 0 denoted $\omega_k(i)$ with $k=10$, and $i=0$. For the window $\omega_k(i)$, we draw a graph with nodes $i, i+1, \cdots, i+k-1$ along the bottom and node $i+k$ above them. For the bottom nodes, connect node $j$ to node $j+1$. Then we connect node $i$ to node $i+k$. This is our initial graph. Now for nodes $i+1$ to node $i+k-1$, add the edge from node $j$ to node $i+k$ in increasing order of interarrival time $T(p_j)$. When two consecutive bottom nodes are connected to the top node, we add the 2-simplex bounded by those edges and the edge on the bottom between them. In the figure, we have added the edges connecting 10 to 3, 4, 7, and 8, and this also adds the 2-simlplices $[3, 4, 10]$ and $[7, 8, 10]$. Adding edges and corresponding 2-simplices in succession creates our filtration. We only look at 1-dimensional homology classes. The resulting diagram is converted to a persistence image. I will explain how this is done in the next section, but for now think of it as a $128\times 128$ array of numbers. The images are then classified using a convolutional neural net. 

According to \cite{CICC}, the technique gave an accuracy of 95.34\%, a recall $(TP/TP+FN)$ of 95.27\%, and a precicison $(TP/TP+FP)$ of $95.46\%$ where TP, FP, and FN are true positive, false positives, and false negatives respectively.

These results came from a single experiment and would not have necessarily done as well in the case we looked at, but it still provides a quick and easy method to try.

Neither this method or the method of Gidea and Katz that we used has a good way of choosing window size. The paper of Myers, Munch, and Khasawneh \cite{MMK}  addreses this issue but only for the univariate case. I will describe it next.

For a time series we need to choose the {\it time lag} $\tau$ and the embedding dimension $d$. The parameter $\tau$ tells us when to take our measurements, and we may not be able to control it if our time series is discrete. Once $\tau$ is chosen, we can think of $d$ as a window size, where our window is $\{x_t, x_{t+\tau}, x_{t+2\tau},\cdots, x_{t+(d-1)\tau}\}$. Instead of a point cloud, though, we use these windows to construct a type of graph called an {\it ordinal partition graph}\index{ordinal partition graph}. I will define this graph first and then discuss how we choose $\tau$ and $d$.

Suppose for example that $d=3$. So each window is of the form $\{x_t, x_{t+\tau}, x_{t+2\tau}\}$. Let 6 consecutive terms of our time series be $\{58, 35, 73, 54, 40, 46\}$. The vertices of the graph will be permutations in the group $S_d$ (recall Example 3.1.6) and we want to find the permutation $\pi\in S_d$ such that $x_{\pi(1)}\leq\cdots\leq x_{\pi(d)}$. Recall that the order of $S_d$ is $d!$. In our case I will list the 6 elements of $S_3$ (in no particular order) as follows:\begin{align*}
\pi_1&=(1)\\
\pi_2&=(12)\\
\pi_3&=(13)\\
\pi_4&=(23)\\
\pi_5&=(123)\\
\pi_6&=(321)
\end{align*}

The first window is $\{58, 35, 73\}$, and to put these in ascending order we need to switch the first 2 elements while keeping the last one fixed. So we use $\pi_2$. Listing the first four windows and the corresponding permutations, we get 

\begin{align*}
 \{58, 35, 73\}&\rightarrow\pi_2\\
 \{35, 73, 54\}&\rightarrow\pi_4\\
 \{73, 54, 40\}&\rightarrow\pi_3\\
 \{54, 40, 46\}&\rightarrow\pi_6
\end{align*}

Now we form a graph with the permutations as vertices and a directed edge from $\pi_i$ to $\pi_j$ if there are consecutive windows corresponding to $\pi_i$ and $\pi_j$. The edges are weighted with weights corresponding to the number of times this transition appears in the graph. Figure 5.9.7 shows the graph corresponding to our example. In practice, since there are $d!$ vertices, we don't include the ones that never appear in our time series.

 \begin{figure}[ht]
\begin{center}
  \scalebox{0.4}{\includegraphics{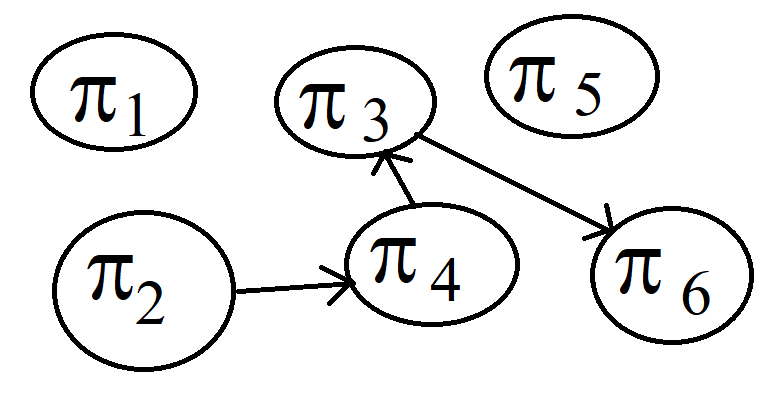}}
\caption{
\rm
Ordinal Partition Graph
}
\end{center}
\end{figure}

Once they build the ordinal partition graph, Myers, et. al. build a filtration using the shortest distance between every pair of vertices. This can be computed using the all\_pairs\_shortest\_path\_length function from the python NetworkX package. The one-dimensional persistence diagram is then computed using the python wrapper Scikit-TDA for the software package Ripser.

It remains to see how they found the time delay $\tau$ and the embedding dimension $d$. First fix $d=3$, and plot the {\it permutation entropy} $$H(d)=-\sum p(\pi_i)\log_2p(\pi_i),$$ where $p(\pi_i)$ is the probability of permutation $\pi_i$ for a range of values of $\tau$. Choose $\tau$ to be the first prominent peak of $H(d)$ in the $H$ as a function of $\tau$ curve. (Note that changing $\tau$ changes the probabilities of each permutation.) It was shown in \cite{RMW}, that the first peak of $H(d)$ is independent of $d$ for $d\in \{3, 4, 5, 6, 7, 8\}$. 

Once we choose $\tau$, define the permutation entropy per symbol to be $$h'(d)=\frac{1}{d-1}H(d),$$ where we make $d$ a free parameter that we want to determine. We then plot $h'(d)$ for $d$ ranging from 3 to 8, and choose the value of $d$ to that maximizes $h'(d)$. 

This leaves the question of whether we could do something similar for multivariate series. Here $d$ is actually the window size. In the example of my internet of things work, $\tau$ is already fixed. One idea might be to take the largest value of $d$ for each of the coordinates and have that be the window size. That would need to be the subject of future work. 

I will conclude this section with some ways of getting a single number summary or {\it score} of a persistence diagram to help with classification. These scores are all listed in \cite{MMK}. If a homology class $x$ is born at time $b$ and dies at time $d$, then the persistence is ${\it pers}(x)=d-b$.  \begin{enumerate}
\item {\bf Maximum Persistence:} This is the persistence of the class with the largest peristence. For diagram $D$, denote this as ${\it maxpers}(D)$.
\item {\bf Periodicity Score:} Let $G'$ be a cycle graph with $n$ vertices in which the nodes and edges form one big cycle. Suppose all edges have weight one and the weight is just the shortest path. Then there is  a single cycle that is born at time 1 and dies at time $\lceil\frac{n}{3}\rceil$ where $\lceil x\rceil$ is the smallest integer that is greater than or equal to $x$. The corresponding diagram $D'$ has $${\it maxpers}(D')=\lceil\frac{n}{3}\rceil-1.$$ If we call this quantity $L_n$, we can compare it to another unweighted graph $G$. (In our case we can use the unweighted ordinal partition graph.) If $G$ has persistence diagram $D$, then define the {\it network periodicity score} $$P(D)=1-\frac{{\it maxpers}(D)}{L_n}.$$ The score ranges from 0 to 1 with $P(D)=0$ if $G$ is a cycle graph. 
\item {\bf The ratio of the number of homology classes to the graph order:} This is $$M(D)=\frac{|D|}{|V|}.$$ The number is not reallly useful in dimension 0 as the number of 0 dimensional classes is $n-1$ for a graph with $n$ vertices. For a higher dimensional diagram, we expect a periodic time series to have a smaller number of classes than a chaotic one so this number should be smaller in this case.
\item {\bf Normalized persistent entropy:} This function is calculated from the lifetimes of the classes in the diagram, and is defined as $$E(D)=-\sum_{x\in D}\frac{{\it pers}(x)}{\frak{L}(D)}\log_2\left( \frac{{\it pers}(x)}{\frak{L}(D)}\right),$$ where $\frak{L}(D)=\sum_{x\in D}{\it pers}(x)$ is the sum of the lifetimes of the points in the diagram. Since it is hard to compare this quantity to diagrams with different numbers of points, we normalize $E$ as $$E'(D)=\frac{E(D)}{\log_2\left(\frak{L}(D)\right)}.$$
\end{enumerate}

Speaking of periodicity, I will mention the paper of Perea and Harer \cite{PH} which uses persistence diagrams to learn about the periodicity of a time series. 

One last comment on that. A quick and easy thing to do is to plot the first derivative of the function against the second derivative. A sine or cosine function plots to a perfect circle. You can tell if a periodic function starts to drift by connecting the points in this plot and looking at the path. (You can do this in R for example.) If the path starts to deviate form your circle, something is happening and those times may be anomalous.

The scores mentioned in \cite{MMK} can be used as machine learning features. I will discuss two more interesting ways to derive features from a persistence diagram in the next section.

\section{Persistence Images and Template Functions} 

Recall that in the last section we saw that \cite{CICC} classified IoT devices using {\it persistence images}. Now I will explain where these images come from.

Persistence images were proposed in the paper of Adams, et. al. \cite{AEKNPSCHMZ}. 

We start form a persistence diagram $B$ in birth-death coordinates. We transform this diagram to birth-lifetime coordinates using $T(x, y)=(x, y-x).$ Let $g_u(x, y)$ be a normalized symmetric Gaussian probability distribution over $R^2$ with mean $u=(u_x, u_y)$ and variance $\sigma^2$ defined as $$g_u(x, y)=\frac{1}{2\pi\sigma^2}e^{-[(x-u_x)^2+(y-u_y)^2]/2\sigma^2}.$$  Fix a nonnegative weighting function $f: R^2\rightarrow R$ that is zero on the horizontal axis, continuous, and piecewise differentiable. We then transform the persistence diagram into a real valued function over the plane. 

\begin{definition}
For a persistence diagram $B$, the corresponding {\it persistence surface}\index{persistence surface} $\rho_B: R^2\rightarrow R$ is the function $$\rho_B(z)=\sum_{u\in T(B)}f(u)g_u(z).$$
\end{definition}

To get our image, we fix a grid of $n$ pixels in the plane and integrate $\rho_B$ over each one to get an $n$-long vector of real numbers. 

\begin{definition}
For a persistence diagram $B$, the corresponding {\it persistence image}\index{persistence image} is the collection of pixels $$I(\rho_B)_p=\iint_p\rho_Bdydx.$$
\end{definition}

Persistence images can also combine persistence diagrams representing homology in different dimensions into a single object by concatenating.

Note that we need to make three choices when generating a persistence image:\begin{enumerate}
\item The resolution of the grid.
\item The distribution function $g_u$ with its mean and standard deviation.
\item The weighting function $f$.
\end{enumerate}

The resolution of the image corresponds to overlaying a grid on the persistence diagram (in birth-lifetime coordinates). It does not have a huge effect on classification accuracy,

The distribution used in the paper centers a Gaussian at each point. Letting the points in the diagram be the means, only the variance need be chosen. This is an open problem but changes in variance do not have a big effect.

The weighting function $f$ from the brth-lifetime plane to $R$ needs to be zero on the horizontal axis (the analogue of the diagonal in birth-death coordinates), be continuous, and be piecewise differentiable. The weighting function used in the paper was piecewise linear and was chosen to give the most weight to classes of maximum persistence. 

Let $b>0$ and define $w_b: R\rightarrow R$ as 
$$w_b(t)=
\left
\{
\begin{array}{ll}
0 & \mbox{if }t\leq 0,\\
t/b & \mbox{if }0<t<b, \mbox{ and}\\
1  & \mbox{if }t\geq b.
\end{array}
\right.
$$

Then set $f(x,y)=w_b(y)$, where $b$ is the persistence value of the longest lasting class in the persistence diagram.

Adams, et. al. show that small changes in Wasserstein distance lead to small changes in persistence images. 

Figure 5.10.1 shows the pipeline for persistent images from data to the persistence diagram to the transformed diagram and finally to the surface and image at different grid scales. 

 \begin{figure}[ht]
\begin{center}
  \scalebox{1.5}{\includegraphics{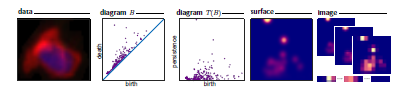}}
\caption{
\rm
Pipeline For Persistent Images \cite{AEKNPSCHMZ}. 
}
\end{center}
\end{figure}

I will comment that the result is more of a numerical matrix than an image. The picture is actually a Matlab heat map. If we rescale to be in the range 0-255, we can make a grayscale image or overlay three of them to get a color image. The numbers can be flattened into a one dimensional vector and sent to your favorite classifier or kept in two-dimensions and sent to a convolutional neural net. Better still, you can keep the Matlab image and sell it to your local modern art museum. (Now that should prove that TDA is practical.) 

For the last major topic in this chapter, I will summarize the paper by Perea, Munch, and Khasawneh on {\it Tempate Functions} \cite{PMK}. The paper is very long and is heavily mathematical, using both point-set topology and {\it functional analysis}. The latter can be thought of of infinite dimensional linear algebra. Finite dimensional spaces are sort of boring from the point of view of topology. All 5-dimensional spaces, for example, are pretty much the same (i.e. they are all homeomorphic to each other). Infinite dimensional spaces can vary a lot more. The typical example is function spaces. For example, integrable periodic real valued functions can be represented by a {\it Fourier series} which is the sum of an infinite number of sine and cosine functions of varying frequencies. We can think of these functions as an infinite vector space basis for the set of periodic functions. The paper cites the book by Conway \cite{Con1} as a reference for functional analysis, but my personal favorite is Rudin, \cite{Rud}.

Like persistent images, the main purpose of template functions is to turn a persistence diagram into a vector of real numbers that can be sent to a classifier. The approach is to build a family of continuous and compactly supported (i.e. only nonzero on a compact set) functions from persistence diagrams to the real numbers. They build two families of these functions: {\it tent functions} which emphasize local contributions of points in a persistence diagram, and {\it interpolating polynomials} which capture global pairwise interactions. Template fucntions provide high accuracy rates, and Perea, et. al. found that in most cases, interpolating polynomials did better than tent functions. Both types of functions can be computed using the open source {\it Teaspoon} package.

I wil now summarize the theory presented in \cite{PMK}. See there for details and especially for the proofs. Recall that a $d$-dimensional persistence diagram is a set $S$ of points in $R^2$ where each point $(x, y)$ represents a homology class of dimension $d$ that is born at time $x$ and dies at time $y$. Here we have $y>x\geq 0$, so we are only looking at points above the diagonal.. We assume the dimension is fixed and we will no longer refer to it. As it is possible for more than one class to be born and die at the same time, the points have multiplicity a nonnegative integer. So let $\mu: S\rightarrow N$ be the {\it multiplicity} where $N$ represents the natural numbers, i.e. the positive integers. We can write $D=(S, \mu)$ for a diagram, and we will often use $D$ and $S$ interchangeably where $S$ is understood to be the set of points in $D$. (So for example, we can talk about subsets of $D$ when we mean subsets of $S$.) 

\begin{definition}
Given $D=(S, \mu)$ and $U\subset R^2$, the {\it multiplicity} of $D$ in $U$ is $$Mult(D, U)=\sum_{x\in S\cap U} \mu(x)$$ if this quantity is finite, and $\infty$ otherwise. 
\end{definition}

Here is some more notation used throughout the paper. The {\it diagonal} is denoted $\Delta=\{(x, x)\}\in R^2$. The {\it wedge} $W$ is the region above and not including the diagonal. If $z=(x, y)\in W$, then {\it pers}$(z)=y-x$. The portion of $W$ with persistence greater than $\epsilon$ is denoted $W^\epsilon$. If we want to include the lower boundary, we write $\overline{W^\epsilon}$.

\begin{definition}
The {\it space of persistence diagrams} denoted $\mathcal{D}$ is the collection of diagrams $D=(S, \mu)$ where:\begin{enumerate}
\item $S\subset W$ called the {\it underlying set of D}, has the property that $Mult(D, W^\epsilon)$ is finite for any $\epsilon>0$.
\item $\mu$ is a function from $S$ to the positive integers. In particular, $\mu(x)$ is the multiplicity of $x\in S$. 
\end{enumerate}
The {\it space of finite persistence diagrams} is $\mathcal{D}^0=\{(S, \mu)\in \mathcal{D}| S \mbox{ is finite}\}$.
\end{definition}

We can describe the multiplicity of multiple persistence diagrams in a subset $U\subset R^2$. by adding their individual multiplicities. 

The next step is to characterize compactness in the space $\mathcal{D}$. We will need this to look at the topology of the space of real valued continuous functions on this space. The space $\mathcal{D}$ is made into a metric space using the bottleneck distance $d_B$. (See Section 5.5.)

\begin{definition}
A subspace of a topological space is {\it relatively compact}\index{relatively compact} if its closure is compact.
\end{definition}

The main theorem is as follows:

\begin{theorem}
A set $\mathcal{S}$ of persistence diagrams in $\mathcal{D}$ is relatively compact if and only if: \begin{enumerate}
\item It is bounded.
\item It is off-diagonal birth bounded (ODBB).
\item It is uniformly off-diagonally finite (UODF). 
\end{enumerate}
\end{theorem}

I will now define each of these conditions. To help picture this, consider Figure 5.10.2 from \cite{PMK}.

\begin{figure}[ht]
\begin{center}
  \scalebox{0.8}{\includegraphics{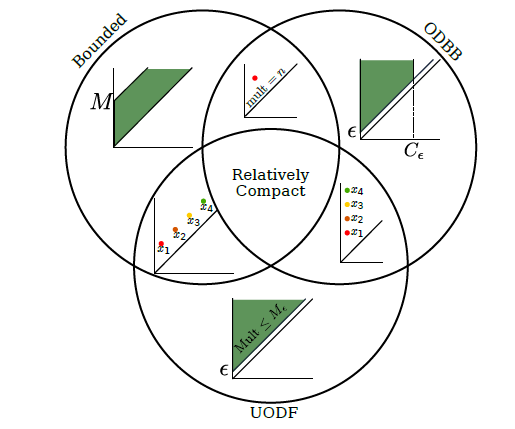}}
\caption{
\rm
Three Criterion for Relative Compactness on Sets of Persistence Diagrams \cite{PMK}. 
}
\end{center}
\end{figure}

\begin{definition}
A subspace of a topological space is {\it bounded}\index{bounded} if it is contained in an open ball of finite radius.
\end{definition}

Note that with bottleneck distance, unmatched points are matched against points on the diagonal. In \cite{PMK}, our definition of Bottleneck distance is modified so that the distance between an unmatched point and the diagonal is $1/2$ the perpendicular distance. This is then half of the points persistence value. So if $\emptyset$ is the empty diagram, $d_B(D, \emptyset)=\frac{1}{2}\max\{pers(x)|x\in S\}$. So if $B_C(D)=\{D'\in\mathcal{D}|d_B(D, D')<C\}$ denotes the ball of radius $C>0$ centered at diagram $D$, then if $\mathcal{S}\subset\mathcal{D}$ is bounded, there exists a $C>0$ with $\mathcal{S}\subset B_C(\emptyset).$ This is the situation pictured in the upper left circle of Figure 5.10.2.

\begin{definition}
A set  $\mathcal{S}\subset\mathcal{D}$ is {\it off-diagonal birth bounded (ODBB)} if for every $\epsilon>0$ there exists a constant $C_\epsilon\geq 0$ such that if $x\in S\cap\overline{W^\epsilon}$ (equivalently, $pers(x)\geq\epsilon$) for $(S, \mu)\in \mathcal{S}$, then $birth(X)\leq C_\epsilon$.
\end{definition}

This is the situation in the top right circle where the $x$ coordinate representing birth has a right hand boundary at $C_\epsilon$.

\begin{definition}
A set  $\mathcal{S}\subset\mathcal{D}$ is {\it uniformly off-diagonal finite (UODF)} if for every $\epsilon>0$ there exists a positive integer$M_\epsilon$ such that $$Mult(D, \overline{W^\epsilon})\leq M_\epsilon$$ for all $D\in\mathcal{S}.$
\end{definition}

This is the situation in the bottom circle. 

Now we list three examples of sets which satisfy only two out of the three conditions. See Figure 5.10.2.

\begin{example}
Let $\mathcal{S}=\{D_n | n\in N\}, D_n=\{(0, 1)\}$ with $\mu_{D_n}(0, 1)=n.$ This is bounded and ODBB but not UODF.
\end{example}

\begin{example}
Let $\mathcal{S}=\{D_n | n\in N\}, D_n=\{(n, n+1)\}$ with $\mu_{D_n}(n, n+1)=1.$ This is bounded and UODF but not ODBB.
\end{example}

\begin{example}
Let $\mathcal{S}=\{D_n | n\in N\}, D_n=\{(0, n)\}$ with $\mu_{D_n}(0, n)=1.$ This is UODF and ODBB but not bounded.
\end{example}

What we are interested in is the topology of continuous functions from $\mathcal{D}$ to $R$. The usual topology on this type of function space is the {\it compact-open topology}.

\begin{definition}
Let $X, Y$ be topological spaces and let $C(X, Y)$ be the set of continuous functions from $X$ to $Y$. Given $K\subset X$ compact and $V\subset Y$ open, let $U(K, V)$ be the set of functions $f\in C(X, Y)$ such that $F(K)\subset V$. The set of $U(K, V)$ for all $K\subset X$ compact and $V\subset Y$ open forms a subbase for a topology on $C(X, Y)$ called the {\it compact-open topology}\index{compact-open topology}. (Recall that this means that open sets are formed by finite intersections and arbitrary unions of the sets $U(K, V)$.)
\end{definition}

\begin{definition}
A topological space is {\it locally compact} if every point has an open neighborhood contained in a compact set. The open set is called a compact neighborhood. 
\end{definition}

Perea, et. al. show the following:

\begin{theorem}
Relatively compact subsets of $(\mathcal{D}, d_B)$ have empty interior. Thus, $(\mathcal{D}, d_B)$ is not relatively compact and no diagram $D\in \mathcal{D}$ has a compact neighborhood. 
\end{theorem}

It turns out that compact subsets of the space of persistence diagrams are nowhere dense (i.e their closure has empty interior.) This means that $\mathcal{D}$ can not be written as a countable union of compact sets. This means that the compact open topology on $C(\mathcal{D}, R)$ is not metrizable (i.e. it can't be made into a metric space.) For more details on why, see Rxample 2.2 , Chapter IV of Conway \cite{Con1}.

Next, Perea et. al. handle the problem of finding compact-open dense subsets of  $C(\mathcal{D}, R)$. In other words, we want a family of functions which will approximate any continuous function $f\in C(\mathcal{D}, R)$. We need the following definition: 

\begin{definition}
A {\it coordinate system} for $\mathcal{D}$ is a collection $\mathcal{F}\subset C(\mathcal{D}, R)$ which separates points. In other words, if $D\neq D'$ are two diagrams in $\mathcal{D}$ then there exists $F\in\mathcal{F}$ for which $F(D)\neq F(D')$.
\end{definition}

Now we want a coordinate system to be small. We could take $\mathcal{F}$ to be the space of all real valued continuous functions on $\mathcal{D}$, but this is too big to be practical. Consider the space $R^n$. We define the cartesian coordinates to be a basis of $n$ continuous linear functions from $R^n$ to $R$. We want to find a continuous embedding (ie, a one-to-one map) of $\mathcal{D}$ into an appropriate topological vector space $V$ which we will need to choose. (Note that a topological vector space is a vector space $V$ over a field $F$ which is a topological space and whose addition and scalar multiplication operations are continuous functions from $V\times V$ to $V$ and $F\times V$ to $V$ respectively. In our case we use the field $R$.) I will skip the rather lengthy derivation and just state the main results. To understand them, we will need a few more definitions. 

\begin{definition}
If $V$ is a topological vector space, its {\it dual} is the vector space $V^*=\{T:V\rightarrow R\}$ such that $T$ is linear and continuous. If the topology of $V$ comes from a norm $||\cdot||_V$, then $V^*$ has the operator norm $$||T||_*=\sup_{||v||_V=1}|T(v)|.$$
\end{definition}

There are several standard topologies on $V^*$ but the one that we will need is the $weak^*$ topology.

\begin{definition}
The $weak^*$ {\it topology} is the smallest topology so that for each $v\in V$, the resulting evaluation function $e_V: V^*\rightarrow R$ defined as $e_V(T)=T(v)$ is continuous. A basis for open neighborhoods of $T\in V^*$ is given by sets of the form $$N(v_1,\cdots,v_k; \epsilon)(T)=\left\{T'\in V^*: \max_{1\leq i\leq k} |T'(v_i)-T(v_i)|<\epsilon\right\}$$ where $v_1, \cdots, v_k\in V$ and $\epsilon>0$.
\end{definition}

We now need a generalization of Banach spaces called {\it locally convex} spaces.

As before let $W$ be the wedge on a persistence diagram and let $K_n$ be a sequence of compact subsets of $W$ with $K_n\subset K_{n+1}$ for all positive integers $n$, and $W=\cup K_n$. Let $C_c(K_n)$ be the set of real valued functions $f$ on $W$ whose support (i.e. $\{x\in W| f(x)\neq 0\}$) is contained in $K_n$. Then $C_c(K_n)$  is a Banach space if endowed with the sup norm $||f||=\sup_{x\in K_n} |f(x)|$. In particular, it is {\it locally convex}.

\begin{definition}
A topological vector space $V$ is {\it locally convex}\index{locally convex} if its topology is generated by a family $\mathcal{P}=\{p_\alpha\}$ of continuous functions $p_\alpha: V\rightarrow [0, \infty)$ so that \begin{enumerate}
\item $p_\alpha(u+v)\leq p_\alpha(u)+p_\alpha(v)$ for all $u, v\in V.$
\item $p_\alpha(\lambda u)=|\lambda| p_\alpha(u)$ for all scalars $\lambda$.
\item If $p_\alpha(u)=0$ for all $\alpha$, then $u=0$.
\item The topology of $V$ is the weakest for which all of the $p_\alpha'$s are continuous.
\end{enumerate}
\end{definition}

In particular, all normed spaces are locally convex. (We just have one $p_\alpha$ with $p_\alpha(v)=||v||$.) Note that each inclusion $C_c(K_n)\subset C_c(K_{n+1})$ is continuous and $C_c(W)=\cup_{n=1}^\infty C_c(K_n).$ 

\begin{definition}
The {\it strict inductive limit topology} is the finest locally convex topology for which each inclusion $C_c(K_n)\subset C_c(K_{n+1})$ is continuous. In this topology, a linear map $T: C_c(W)\rightarrow Y$ where $Y$ is locally convex is continuous if and only if the restriction of $T$ to each $C_c(K_n)$ is continuous. 
\end{definition}

Let $C_c(W)^*$ be the dual of $C_c(W)$ with repect to the strict inductive limit topology and let $C_c(W)^*$ have the weakest topology so that for each $f\in C_c(W)$, the resulting evaluation function $e_f: C_c(W)^*\rightarrow R$ with $e_f(T)=T(f)$ is continuous. This is the corresponding $weak^*$ topology. 

\begin{theorem}
Given a persistence diagram $D=(S, \mu)\in \mathcal{D}$ and a function $f\in C_c(W)$, define $$\nu_D(f)=\sum_{x\in S}\mu(x)f(x).$$ If $C_c(W)$ is endowed with the strict inductive limit topology and $C_c(W)^*$ with the corresponding $weak^*$ topology, then $\nu :  \mathcal{D}\rightarrow C_c(W)^*$ defined by $\nu(D)=\nu_D$ is continuous, injective, and if $A\sqcup B$ denotes the disjoint union of $A$ and $B$, then for $D, D'\in\mathcal{D}$, $\nu(D\sqcup D')=\nu(D)+\nu(D').$
\end{theorem}

So now that we have an embedding $\nu: \mathcal{D}\rightarrow C_c(W)^*$, we an look for coordinate functions. The next theorem characterizes the elements of $C_c(W)^{**}$, the dual of the dual.

\begin{theorem}
Let $V$ be a locally convex space and endow its topological dual $V^*$ with the associated $weak^*$ topology. This is the smallest topology such that all the evaluations $e_V: V^*\rightarrow R$ defined by $e_v(T)=T(v)$ for $v\in V$ are continuous. Then the function $e: V\rightarrow V^{**}$ defined by $e(v)=e_V$ is an isomorphism of locally convex spaces.
\end{theorem}

Applying this theorem to the locally convex space $C_c(W)$ topologized with the strict inductive limit topology implies that the elements of $C_c(W)^{**}$ are evaluations $e_f$ with $f\in C_c(W)$. Composing $e_f$ with $\nu$ yields a continuous function $e_f\nu: \mathcal{D}\rightarrow R$ and $e_f\nu(D)=\nu_D(f)=\nu_f(D),$ where $\nu_D(f)=\nu_f(D)=\sum_{x\in S}\mu(x)f(x).$ The $\nu_f$ are the functions we are looking for.

\begin{definition}
A {\it template system} for $\mathcal{D}$ is a collection $\mathcal{T}\subset C_c(W)$ so that $$\mathcal{F}_\mathcal{T}=\{\nu_f| f\in\mathcal{T}\}$$ is a coordinate system for $\mathcal{D}$. The elements of $\mathcal{T}$ are called {\it template functions}\index{template functions}. 
\end{definition}

Template functions can approximate continuous functions on persistence diagrams in the following sense: 

\begin{theorem}
Let $\mathcal{T}\subset C_c(W)$ be a template system for $\mathcal{D}$, let $\mathcal{C}\subset\mathcal{D}$ be compact, and let $F: \mathcal{C}\rightarrow R$ be continuous. Then for every $\epsilon>0$, there exists a positive integer $N$, a polynomial $p\in R[x_1, \cdots, x_N]$, and template functions $f_1, \cdots, f_N\in \mathcal{T}$ such that $$|p(\nu_D(f_1), \cdots, \nu_D(f_N))-F(D)|<\epsilon$$ for every $D\in\mathcal{C}$. In other words, the collection of functions of the form $D\rightarrow p(\nu_D(f_1), \cdots, \nu_D(f_N))$ is dense in $C(\mathcal{D}, R)$ in the compact-open topology.
\end{theorem}

We now want to see how to construct template functions. The next theorem shows that we can construct a countable template system by translating and rescaling a function $f\in C_c(W)$. 

\begin{theorem}
Let $f\in C_c(W),$ $n$ be a postive integer, and $m\in R^2$ have integer coordinates. Define $$f_{n, m}(x)=f(nx+\frac{m}{n}).$$ If $f\neq 0$, then $$\mathcal{T}=\{f_{n,m}\}\cap{C_c(W)}$$ is a template system for $\mathcal{D}$. Recall that a function $f$ is {\it Lipschitz} if for any $x, y$ in the domain of $f$ there is a constant $C$ independent of $x$ and $y$ such that $|f(x)-f(y)|\leq C|x-y|$. Then if $f$ is Lipschitz, the elements of the coordinate system $$\{\nu_{f_{n, m}}=f_{n, m}\nu| f_{n, m}\in\mathcal{T}\}=\mathcal{F}_\mathcal{T}$$ are Lipschitz on any relatively compact set $\mathcal{S}\subset\mathcal{D}$.
\end{theorem}

Perea, et. al. describe two families of template functions: {\it tent functions} and {\it interpolating polynomials}. I will describe them next. The families will be defined in the birth-lifetime plane instead of the birth-death plane just to keep the definitions a little simpler. The shifted diagrams will be denoted by a tilde. So the wedge $\widetilde{W}$ consists of points $(x, y)\in R^2$ such that $x\geq 0$ and $y>0$, and $\widetilde{W^\epsilon}$ is the subset of $\widetilde{W}$ where $y>\epsilon$. The point $x=(a, b)$ in the birth-death plane is transformed to $\tilde{x}=(a, b-a)$ in the birth-lifetime plane. A diagram $D=(S, \mu)$ is transformed to $\widetilde{D}=(\widetilde{S}, \tilde{\mu})$ where $\tilde{\mu}(\tilde{x})=\mu(x)$. 

The first example is {\it tent functions}\index{tent functions}. Given a point $A=(a, b)\in \widetilde{W}$ and a radius $\delta$ with $0<\delta<\min\{a, b\}$, define the tent function on $\widetilde{W}$ to be $$g_{A, \delta}(x, y)=\left|1-\frac{1}{\delta}\max\{|x-a|, |y-b|\}\right|_+$$ where $|\cdot|_+$ means take the value of the function if it is positive and let it be zero otherwise. Since $\delta<\min\{a, b\}$, the function is zero outside of the compact box $[a-\delta, a+\delta]\times [b-\delta, b+\delta]\subset\widetilde{W}.$

Given a persistence diagram $D=(S, \mu)$, the tent function is $$G_{A, \delta}(D)=\widetilde{G}_{A, \delta}(\widetilde{D})=\sum_{\tilde{x}\in\widetilde{S}}\tilde{\mu}(\tilde{x})g_{A, \delta}(\tilde{x}).$$ To limit the number of these functions, let $\delta>0$ be the partition scale, let $d$ be the number of subdivisions along the diagonal (birth-death) or $y$-axis (birth-lifetime) and let $\epsilon$ be the upward shift. Then use only the tent functions $$\{G_{(\delta i, \delta j+\epsilon), \delta}| 0\leq i\leq d, 1\leq j\leq d\}.$$ An example plot from \cite{PMK} is shown if Figure 5.10.3.

\begin{figure}[ht]
\begin{center}
  \scalebox{0.8}{\includegraphics{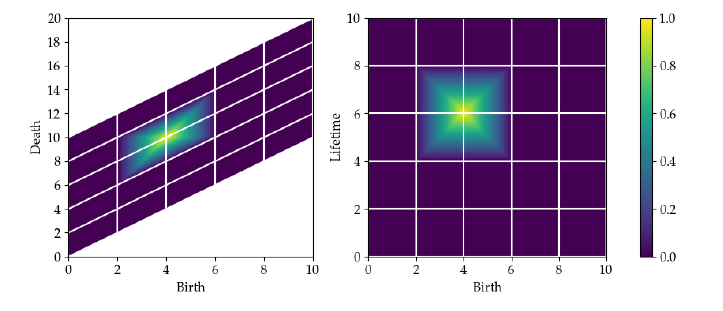}}
\caption{
\rm
Tent function $g_{(4,6), 2}$ drawn in the birth-death plane and the birth-lifetime plane \cite{PMK}. 
}
\end{center}
\end{figure}

Now we have a vector of numbers the output which can be sent to your favorite classifier. Perea, et. al show that tent functions {\it separate points}. In other words, they have different values for different persistence diagrams. 

The other template functions are {\it interpolating polynomials}. They are formed by constructing polynomials with specified values at mesh points of $R^2$ and evaluating them on the points of the persistence diagram. The construction is somewhat more complicated so I will refer you to \cite{PMK} for details. Perea et. al. claimed that these preformed even better than tent functions in their experiments. 

The Teaspoon Package in Python does all of the work for your, constructing both tent functions and interpolating polynomials from your data points. It can be found at https://github.com/lizliz/teaspoon .

\section{Open Source Software and Further Reading}

In this chapter, I showed you what persistence homology is, how it can arise form point clouds, sublevel sets, graphs, and time series, and some methods for pulling out numerical features to feed to a classifier. This should give you lots of things to try on your own problems. 

There is a lot more written on TDA that I haven't covered. Here are some more places to look for further reading. There are five books I can list on TDA that each have a different perspective. The oldest is {\em Computational Homology} by Kaczynski, Mischaikow, and Mrozek \cite{KMM}.This book doesn't cover persistent homology, but covers algorithms for computing homology groups. The books by Zomorodian \cite{Zo1} and Edelsbrunner/Harer \cite{EH} are the earliest books I know of that cover persistent homology. Ghrist's {\em Elementary Applied Topology} \cite{Gh1} gives a taste of a lot of potential applications of topology which you can follow up on if any of the topics interest you. The book by Rabadan and Blumberg \cite{RB} is the most recent and modern, and it is heavily geared towards biomedical applications, especially virus evolution and cancer. The mathematical background material is done very well, and they cover some topics that I skipped. For example, they talk about the use of statistical techniques in topological data analysis. Many of these problems are hard and still unsolved. 

I also, left out the idea of {\it witness complexes} in which a subset of a point cloud called {\it landmark points} is chosen to reduce computational cost. There are extremely fast programs now so this is probably not as much of a problem as it once was.

There are also some good papers that will give you an introduction to TDA. Here are some good papers to get you started. Carlsson \cite{Car1} was one of the earliest introductory papers. Chazal and Michel \cite{CM} is a more modern introduction. Pun, et. al. \cite{PXL} is one of the best organized surveys out there giving you options at every step of the pipeline. They also list software, but they leave out Ripser. Finally, I should mention the introduction by Munch \cite{Munc1} which was used as background for the January 2021 American Mathematical Society short course on Mathematical and Computational Methods for Complex Social Systems. 

One advantage of experimenting with TDA is that there is now a lot of open source software. Ripser is probably the fastest for doing persistent homology, but it doesn't do everything. Javaplex, Hera, and Dionysus are two alternatives. I also mentioned Elizabeth Munch's Teaspoon package which handles template functions. Also check out packages that are available in R. Two comprehensive lists of software and what it can do are found in \cite{PXL} and Appendix A of \cite{RB}.

So now that you know so much about TDA, what comes next. In the next chapter, I will talk about some famous and important algorithms that are topology related but not exactly persistent homology. Then I will briefly discuss a couple of my own ideas that I never got to fully pursue. 

After that, I will take you to the frontier and look at recent developments that could help tap the more advanced techniques of cohomology, Steenrod squares, and homotopy theory. 

\chapter{Related Techniques}

In this chapter, I will cover some applications of algebraic topology which do not involve persistent homology. Section 6.1 will cover Ronald Atkin's "Q-Analysis." Q-Analysis, described in \cite{Atk}, is the earliest attempt I can think of to apply algebraic topology to the real world. In this case, it was to explain the structure and politics of his university, the University of Essex in Colchester, UK. I will spend some time on this subject as nobody else is likely to explain it to you. Q-analysis never really caught on, but I think that Atkin had some interesting things to say about large bureaucracies and how they function. 

In Section 6.2, I will talk about sensor coverage, one of the earliest proposed applications of computational homology theory. Suppose you have an array of sensors that know if their coverage area overlaps with that of neighboring sensors. Can you tell if your coverage region has any holes? 

Section 6.3 will cover Mapper, the main commercial product of Gunnar Carlsson's company, Ayasdi. I will describe the algorithm and its breakthrough application: The discovery of a new class of breast cancer cell. 

In Section 6.4, I will discuss simplicial sets, a generalization of simplical complexes. They are necessary to understand the dimensionality reduction algorithm UMAP, and they are also used in modern computational algebraic topology. Finally, I will discuss the UMAP algorithm itself in Section 6.5. 

\section{Q-Analysis}

In \cite{Atk}, Ronald Atkin looked at the problem of using topology to describe social science data and discover relationships. To do that, he used simplicial complexes to describe binary relations. Remember that an abstract simplical complex consists of a collection of subsets of a vertex set such that if $\sigma$ is in the complex, then so is any face of $\sigma$. The complex then can be used to illustrate relationships. For example, if $X$ is a set of people and $Y$ is a set of committees, we could build a simplicial complex to illustrate the relation "is a member of" by letting $x_i\in \sigma_j$ if person $i$ is a member of committee $j$. We also have the inverse relationship committee $j$ contains person $i$. Q-analysis looks at the geometric structure of these relations at each level of a heirarchy in which the $N+1$ level groups together objects at the $N$ level. For example, if football, baseball, and basketball are objects at the $N$ level they might be grouped together as the object "sports" at the $N+1$ level. At each level we look at how many elements are shared between simplices, {\it patterns} which asscoiate numbers to simplices, and the forces that cause patterns to evolve. 

In Sections 6.1.1, 6.1.2, and 6.1.3, I will make these concepts concrete by discussing topological representation of relations, exterior algebra, and shomotopy respectively. Then I will choose some examples form Atkin's description of the University of Essex to illustrate the concepts. In Section 6.1.4, I will discuss the application to the University's committee structure. Finally, section 6.1.5 will cover rules for success in handling the hierarchy. 

Although there have been a number of papers on Q-analysis since this book was published, you rarely see it in modern discussions of TDA. Still, Atkin had some interesting things to say about the workings of bureaucracies and it would be interesting to apply similar ideas to ohter large organizations such as corporations of Government Agencies. 

\subsection{Topological Representation of Relations}
For the definitions in this section, I will refer to the following example from \cite{Atk}:

Given a relation $\lambda\subset Y\times X$, we have a simplicial complex $KY(X; \lambda)$ where elements of $X$ correspond to vertices and elements of $Y$ correspond to maximal simplices. In what follows the term "simplex" will just be used for the maximal simplices.

\begin{example} 
Suppose we have the following {\it incidence matrix}\index{incidence matrix} $\Lambda$ in which the rows are maximal simplices and the columns are vertices. $\Lambda_{ji}=1$ if $x_i$ is in simplex $Y_j$ and $\Lambda_{ji}=0$ otherwise:
\[\Lambda=\begin{bmatrix}
1 & 1 & 1 & 1 & 0 & 0 & 0 & 0\\
0 & 0 & 1 & 1 & 1 & 0 & 0 & 0\\
0 & 0 & 0 & 0 & 1 & 0 & 0 & 1\\
0 & 0 & 0 & 0 & 0 & 1 & 1 & 1\\
0 & 0 & 1 & 0 & 0 & 0 & 1 & 0\\
0 & 0 & 0 & 1 & 0 & 1 & 0 & 1
\end{bmatrix}\]

In this case, we have 6 maximal simplices: $Y_1=[x_1, x_2, x_3, x_4]$, $Y_2=[x_3, x_4, x_5]$, $Y_3=[x_5, x_8]$, \newline
$Y_4=[x_6, x_7, x_8]$, $Y_5=[x_3, x_7]$, $Y_6=[x_4, x_6, x_8]$.

Note that we can also define the complex $KX(Y; \lambda^{-1}))$ whose incidence matrix is the transpose of $\Lambda$.
\end{example}

Now we look at chains of simplices that have vertices in common.

\begin{definition}
Two simplices $\sigma$ and $\tau$ are $q${\it-connected} if there is a chain $\sigma=\sigma_1, \sigma_2, \cdots, \tau=\sigma_h$ such that each consecutive pair $\{\sigma_i, \sigma_{i+1}\}$ has a face of at least dimension $q$ in common.  We say that this is a chain of $q${\it -connection} of length $h-1$. 
\end{definition}

WARNING: This is not the same as the term $n$-connected used in homotopy theory. 

Also note that I am departing from Atkin's notation where $\sigma_n$ means that $n$ is the dimension of $\sigma.$

A $q$-simplex is always $q$-connected to itself by a chain of length 0. Also, if $\sigma$ and $\tau$ are $q$-connected, they are also $n$-connected for $0\leq n<q$.

\begin{example}
Let $\sigma_1=[x_1, x_2, x_3, x_4, x_5, x_6]$, $\sigma_2=[x_2, x_3, x_4, x_5, x_6, x_9]$, $\sigma_3=[x_4, x_5, x_6, x_7, x_8, x_9]$, and $\sigma_4=[x_6, x_7, x_8, x_9, x_{10}]$. Then $\sigma_1\cap\sigma_2=[x_2, x_3, x_4, x_5, x_6]$, $\sigma_2\cap\sigma_3=[ x_4, x_5, x_6, x_9]$, and $\sigma_3\cap\sigma_4=[ x_6, x_7, x_8, x_9]$. So $\sigma_1$ and $\sigma_2$ share a face of dimension 4, $\sigma_ 2$ and $\sigma_3$ share a face of dimension 3, and $\sigma_ 3$ and $\sigma_4$ share a face of dimension 3. So $\sigma_1$ and $\sigma_4$ are 3-connected through a chain of length 3.
\end{example}

For a simplicial complex $K$ of dimension $n$, fix $q$ with $0\leq q\leq n$. Then the property of being $q$-connected is an equivalence relation on the simplices of $K$. We call the equivalence classes $q${\it-components} and denote their cardinality by $Q_q$. If we analyze $K$ by finding the values of $Q_0, Q_1, \cdots, Q_n$, we say that we have perfomed a {\it Q-analysis} on $K$. 

\begin{example}
Referring to Example 6.1.1, The largest simplex is $Y_1$ of dimension 3. Since there is no other of that dimension, we have $Q_3=1.$ There are 4 simplices of dimension at least 2 (remember we only consider maximal simplices). These are $Y_1$, $Y_2$, $Y_4$ and $Y_6$, and none of these are 2-connected to anything. So $Q_2=4$. We also have $Q_1=4$, since $Y_1$ and $Y_2$ share the 1-face $[x_3, x_4]$, and $Y_4$ and $Y_6$ share the 1-face $[x_6, x_8]$ while $Y_3$ and $Y_5$ are maximal 1-simplices. Finally, $Q_0=1$ as every simplex intersects at leasr one other so that the complex $K$ is connected.
\end{example}

Note that $Q_0$ is the number of connected components, so it always equals the 0-betti number in homology. 

We can do this calculation quickly by mutliplying $\Lambda\Lambda^T$ and subtracting 1 from every entry of the resulting matrix. $\Lambda\Lambda^T$  is a symmetric matrix and $(\Lambda\Lambda^T)_{ij}$ is the number of ones that row $i$ and row $j$ have in common. So subtracting 1 gives the dimension of the largest common face.

\begin{definition}
Given a complex $K$ with dimension $n$ the {\it structure vector} is $Q(K)=[Q_n, \cdots, Q_1, Q_0]$.
\end{definition}

\begin{example}
In our previous example, the structure vector is $Q=[1, 4, 4, 1].$
\end{example}

\begin{definition}
Given a complex $K$ with dimension $n$ and structure vector $Q(K)=[Q_n, \cdots, Q_1, Q_0]$, the {\it obstruction vector} $\hat{Q}$ is defined as $\hat{Q}=[Q_n, Q_{n-1}-1, \cdots, Q_1-1, Q_0-1]$
\end{definition}

The obstruction vector is a property of the complex and acts as an obstruction to the free flow of patterns which will be defined next.

Word reuse warning: This obstruction has nothing to do with the obstruction theory I will define in Chapter 10.

\begin{example}
In our previous example, the obstruction vector is $\hat{Q}=[1, 3, 3, 0].$
\end{example}

\begin{definition}
A {\it pattern} $\pi$ on a complex $K$ is a mapping from the simplices of $K$ to the integers. The restriction of $\pi$ to $t$-simplices is written $\pi^t$, 
\end{definition}

Note that this mapping is not necessarily one-to-one and it is not necessarily a homomorphism.

Atkin describes the underlying geometry as the {\it static backcloth} denoted $S(N)$. An incremental change in pattern $\delta\pi$ on $S(N)$ is the analogue of a force in physics. The {\it intensity} of the force  is on a simplex $\sigma$ is $$\frac{\delta\pi(\sigma)}{\pi(\sigma)}.$$ Changes in patterns in the environment can induce social forces  such as social pressure, organizational pressure, etc. In this case it will be the patterns in the univeristy environment.

\subsection{Exterior Algebra}

In this section, we will discuss the properties of an exterior algebra. You may have seen this if you have ever worked with differential forms. 

Let $V$ be a finite dimensional vector space over a field $F$ with basis $\{v_1,\cdots, v_n\}$. We turn this space into an algebra using the {\it exterior product} or {\it wedge product.}

\begin{definition}
An {\it exterior algebra} $V$ is an  algebra over $F$ whose product is the {\it wedge product} $v\wedge w$ having the properties that $v\wedge w=-w\wedge v$ (antisymmetry) and $v\wedge v=0$ for $v, w\in V$. Using the distributive law, if $a_1, a_2, b_1, b_2\in F$, and $v, w\in V$, then $$(a_1v+a_2w)\wedge (b_1v+b_2w)=(a_1b_2-a_2b_1)(v\wedge w).$$ In addition, the product is associative (i.e. if $v, w, z\in V$, then $(v\wedge w)\wedge z=v\wedge (w\wedge z).$ This makes $V$ into an associative algebra.
\end{definition}

\begin{definition}
Let $V$ be a $n$-dimensional vector space over a field $F$ with basis $\{v_1,\cdots, v_n\}$. Define a new vector space $\wedge^2V$ with basis consisting of all elements of the form $v_i\wedge v_j$ with $v_i$ and $v_j$ basis elements of $V$ and  $i\neq j$. In general can define $\wedge^kV$ for $2\leq k\leq n$ with basis consisting of $v_{i_1}\wedge\cdots\wedge v_{i_k}$ where the product is taken over all $k$-tuples of distinct basis elements of $V$. Let $\wedge^0=F$ and $\wedge^1=V$. Then we define the {\it exterior algebra on} $V$ to be the asscociative algebra $\wedge V$ defined as $$\wedge V=\wedge^0V\oplus\wedge^1V\oplus\wedge^2V\oplus\cdots\oplus\wedge^nV.$$
\end{definition}

Note that if $x\in\wedge^pV$ and $y\in \wedge^qV$, then $(x\wedge y)\in \wedge^{p+q}V$ provided that $p+q\leq n$. Also the dimension of $\wedge^p$ is the nuber of combinations of selecting $p$ objects from $n$ objects (order doesn't matter) or $$\binom{n}{p}=\frac{n!}{(n-p)!p!}.$$

We will now modify this construction for $V$ a module over the integers. (I,e. an abelian group,) Ordering all of the elements of a set $X$, let $V=\{x_1, x_2, \cdots, x_n\}$. We can then associate simplicial complexes with polynomials over $\wedge V$. We will leave out the wedge symbol and just write $v_i\wedge v_j$ as $v_iv_j$. Then a $p$-simplex $[x_{i_0}, x_{i_1}, \cdots, x_{i_p}]$ can be written as the element $x_{i_0}x_{i_1}\cdots x_{i_p}\in\wedge^{p+1}$.

For a complex $K$, and a simplex in $K$, let $\rho(\sigma)$ be the monomial in $\wedge V$ corresponding to $\sigma$. Suppose we have a pattern $\pi$ on $K$. Then the {\it pattern polynomial}, also denoted as $\pi$, is defined to be $$\pi=1+\sum_{\sigma\subset K}\rho(\sigma)\pi(\sigma).$$ 

Now we introduce a function that takes a simplex to the sum of its faces.

\begin{definition}
If $\sigma=[v_0, v_1, \cdots, v_p]$ is a simplex, then we define the face operator $f$ by $$f\sigma=\sum_{i=0}^p[v_, \cdots, \hat{v}_i, \cdots, v_p].$$
\end{definition}

Note that this looks like the boundary but does not involve an alternating sum. Also, we are not working in $Z_2$ so it is not true that applying $f$ twice gives zero.We will let $f^0$ be the identity map and define $f^n=f(f^{n-1})$ for $n>0.$ 

Letting the simplices be monomials in $\wedge V$, $f$ becomes a linear function, i.e. $f(a\sigma)=af(\sigma)$ for $a\in Z$, and $f(\sigma+\tau)=f(\sigma)+f(\tau)$. For a 0-dimensional simpllex $\sigma=[v]$, we define $f(\sigma)=1.$

\begin{example}
In $\wedge V$,  $f(x_1x_2x_3)=x_1x_2+x_1x_3+x_2x_3$. Then $f^2(x_1x_2x_3)=f(x_1x_2+x_1x_3+x_2x_3)=x_1+x_2+x_1+x_3+x_2+x_3=2x_1+2x_2+2x_3\neq 0.$
\end{example}

Now if $\sigma$ and $\tau$ share a vertex, then $\sigma\tau=0$ in $\wedge V$ since the vertex appears twice in the product. If they share two vertices, then we also have $f(\sigma)\tau=\sigma f(\tau).$ But $\sigma f^t(\tau)\neq 0$ for $t>1$. So continuing this argument gives the following:

\begin{theorem}
Two simplices $\sigma$ and $\tau$ in complex $K$ share a $t$-face if and only if their monomial representations in $\wedge V$ satisfy $(f^i\sigma)\tau=0$ or $\sigma(f^i\tau)=0$ for $0\leq i\leq t$ and both terms are nonzero for $i>t$.
\end{theorem}

\subsection{Homotopy Shomotopy}

We now introduce a concept of nearness between simplices.

\begin{definition}
If two simplices share a $q$-face, we say that they are {\it q-near}.
\end{definition}

We can describe a chain of $q$-connectivity in terms of $q$-nearness. In a chain of $q$-connection, successive simplices are $q$-near. In a simplicial complex, we can extend the concept of nearness to chains of connection. We will denotes such a chain from $\sigma$ to $\tau$ as $[\sigma, \tau]$.

\begin{definition}
Let $[\sigma_0, \sigma_n]$ and $[\tau_0, \tau_m]$ be two chains of connection. The intermediate simplices are numbered $\sigma_1, \cdots, \sigma_{n-1}$ and $\tau_1, \cdots, \tau_{m-1}$ repectively.  Then $[\sigma_0, \sigma_n]$ and $[\tau_0, \tau_m]$ are $q$-{\it adjacent} if 
\begin{enumerate}
\item $\sigma_0$ is $q$-near $\tau_0$.
\item $\sigma_n$ is $q$-near $\tau_m$.
\item Given $\sigma_i$ in $[\sigma_0, \sigma_n]$ , there is a $\tau_j$ in $[\tau_0, \tau_m]$ such that $\sigma_i$ is $q$-near $\tau_j$. Let $\tau_{j_1}$ and $\tau_{j_2}$ be the chosen simplices that are $q$-near $\sigma_{i_1}$ and $\sigma_{i_2}$ respectively. Then $i_1<i_2$, implies $j_1<j_2$. 
\item Given $\tau_j$ in $[\tau_0, \tau_m]$ , there is a $\sigma_i$ in $[\sigma_0, \sigma_n]$ such that $\tau_j$ is $q$-near $\sigma_i$. Let $\sigma_{i_1}$ and $\sigma_{i_2}$ be the chosen simplices that are $q$- near$\tau_{j_1}$ and $\tau_{j_2}$  respectively. Then $j_1<j_2$, implies $i_1<i_2$. 
\end{enumerate}
\end{definition}

A $q$-chain of connection is obviously $q$-adjacent to itself. The relation is also symmetric but it is not transitive. 

Atkin now defines a discrete analog of homotopy called a {\it pseudo-homotopy} or {\it shomotopy} for short. 

\begin{definition}
Let $K^+$ be a simplicial complex with a simplex $\sigma_{-1}$ added. $\sigma_{-1}$ is a (-1)-dimensional simplex which will just be something we can use for numbers that are out of range in the function we will now define. Let $c_1$ and $c_2$ be chains of $q$-connection of length $n_1$ and $n_2$ respectively. A (discrete) $q$-{\it shomotopy} is a mapping $Sh: Z\times Z\rightarrow K^+$ such that 
\begin{enumerate}
\item For fixed $x\in Z$, the image of $Sh(x, y)$ is a chain of $q$-connection. As this chain will have finite length $m$ since $K$ is a finite complex, we let $S(x, y)=\sigma_{-1}$ for $y<0$ or $y>m$. 
\item $Sh(0, y)=c_1$ and $Sh(m, y)=c_2$ for some finite $m>0$.
\item $Sh(x, y)$ and $Sh(x+1, y)$ are $q$-adjacent for $0\leq x<m$. 
\end{enumerate}
\end{definition}

Shomotopy is an equivalence relation and the finite set of all chains of connection in $K$ fall into disjoint equivalence classes called $q$-{\it shomotopy classes.}

An analogue of a continuous function (more specifically an isometry) is a {\it face saving map}.

\begin{definition}
Given 2 complexes $K$ and $K'$, let $\beta: K\rightarrow K'$ be an injective map. $\beta$ is a {\it face saving map} if for simplices $\sigma, \tau\in K$, $\sigma$ is $q$-near to $\tau$ in $K$, implies that $\beta(\sigma)$ is $q$-near to $\beta(\tau)$ in $K'$.
\end{definition}

A face saving map preserves $q$-connectivites but need not preserve chain length. This means that these maps also preserve the porperty of being $q$-adjacent and preserve $q$-shomotopies.

Shomotopy gives an analogue of the homotopy groups I will discuss in Chapter 9. As a preview consider the {\it fundamental group} $\pi_1(X, x_0)$. This consists of homotopy classes of maps $f$ from $I=[0,1]$ to $X$ where $f(0)=f(1)=x_0$. These are loops that start and end at the {\it base point} $x_0$. There is a group structure where the multiplication $[x]\dot[y]$ consists of traveling around a representative loop in class $[x]$ followed by one in class $[y]$. The loop where all of $I$ is mapped to the basepoint is the identity and if $[x]$ is represented by $f: I\rightarrow X$, then $[x]^{-1}$ is represented by $g:I\rightarrow X$ where $g(t)=f(1-t)$ so that the loop is traversed in the opposite direction. I will get into much more details in Chapter 9.

For the shomotopy version, let $\sigma_0$ be a fixed simplex in $K$. Then $\sigma_0$ is a {\it base simplex} and a chain of $q$-connection $[\sigma_0, \sigma_0]$ is a $q$-{\it loop} based at $\sigma_0$. Group structures are entirely analogous as the identity is the $q$-loop of length 1 while an inverse of a loop is achieved by listing the simplices in the opposite direction. Loops that are not shomotopic to a loop of length 1 represent holes in the complex.

\subsection{Structures of the University of Essex}

Atkin describes a research study that he did at the University of Essex in Colchester, UK where he was a professor. The idea was to see how the structures described above would apply to various aspects of the University. These categories ranged from phyiscal structure, to administrative organization, committee structures, and social amenities. For each category, there was a set of {\it heirarchical levels} $N-1, N, N+1, \cdots$. Objects at a higher level include objects at the lower levels. Atkin gives an example of social amenities which I have reproduced in Table 6.1.1.

\begin{table} 
\begin{center}
\begin{tabular}{|c|c|}
\hline
Level & Object\\
\hline\hline
N+1 & \{Social Amenity\}\\
\hline
N & \{Catering, Sport\}\\
\hline
N-1 & \{Food, Drink, Groceries, \ldots, Athletics, Ball Games, \ldots\}\\
\hline
N-2 & \{Pies, Chips, \ldots, Beer, Wine, \ldots, Carrots, Milk, Meat, \ldots, Football, Tennis, High-Jump, Billiards, \ldots\}\\
\hline
\end{tabular}
\caption{Social Amenities at the University of Essex\cite{Atk}.}
\end{center}
\end{table}

Note that the levels are not partitions. For example, Drink and Groceries share Milk as a common object. That is where $q$-connectivity plays a role. 

The report covered levels from $N+3$ which was the world outside the University down to the level of the indivisual which was $N-2$, but an individual typically operates at multiple levels depending on the context. 

Within a level $N$, the $N$-level backcloth $S(N)$ is the space in which activity takes place. It is analyzed by looking at the structure vector $Q$ and its associated obstruction vector $\hat{Q}$. This gives a global analysis. A more local one is the study of shomotopy and especially the $q$-holes which the traffic at that level must avoid. These holes look like an object to an individual. There is also the group structure whose role is poorly understood. Anyway, the structure affects how the University community interacts and how decisions are made. 

I will give a couple of examples from the rest of the book. I would have appreciated a lot of them more if I had ever been to the University of Essex as there is a lot of description of the various buildings, dining halls, etc. I imagine it must have changed in some way in the last 50 years, though. I have been to Colchester in 2002 and saw the Colchester Castle, not to be confused with the "citadel" described below.

The first example I will describe is the committee structure. There are 28 of these including the University Council which operates both at level $N+3$ (the interface with the outside world) and at level $N+2$. All of the others operate at level $N+1$.

The maximal simplices are people and the vertices are the committees they are on. The Vice Chancellor is on 20 of these committees representing a simplex of dimension 19. We then have two added components at dimension 9 which are added by two people, one in the History Department and one in the Literature Department. The structure vector is $$Q=\{1, 1, 1, 1, 1, 1, 1, 1, 1, 1, 3, 2, 3, 3, 6. 4. 5. 9. 7. 1\}$$ where the positions go from dimension 19 down to 0 (left to right) and the first 3 is in dimension 9. This gives the obstruction vector $$\hat{Q}=\{2, 1, 2, 2, 5, 3, 4, 8, 6, 0\}.$$ The obstruction vector has its largest value at $q=2$ representing business that must be discussed by 3 committees. As the $q$-connectivity represents committees that are shared by particular individuals, high values in the obstruction vector show more rigidity. A pattern vector may show a ranking each individual assigns to each order of business. To get things done, these patterns need to change and the higher the obstruction value, the more resistance. 

Flipping things around, we get a complex where vertices are people and maximal simplices are committees. The connectivity consist of 2 committees that share people. This is required for good communication and a high obstruction means that there are several groups of committees that don't share members, so they don;t talk to each other. The structure and obstruction vectors in this case also had some large components. 

Atkin proposed that when committees are formed (an annual process), Q-Analysis may help to determine membership that will be more or less resistant to change. 

Another interesting quantity in this context is the {\it eccentricity} of a simplex.

\begin{definition}
Let $\sigma$ be an $r$-dimensional simplex in a complex $K$. Define $\check{q}$ (the bottom-q) to be the least value of $q$ for which a distinct simplex is $q$-connected to $\sigma$. Note that if 2 simplices are $q$-connected, they are also $t-$connected for $t<q$, so for this definition to make sense we use the highest number for which they are $q$-connected. Then define the {\it eccentricity} of $\sigma$ to be $$Ecc(\sigma)=\frac{r-\check{q}}{\check{q}+1}.$$
\end{definition}

Higher eccentricity means a simplex is the most out of step with the others. If $\sigma$ does not share any vertices with another simplex, we say that $\check{q}=-1$ and $Ecc(\sigma)=\infty.$

It turns out that the individual with the highest eccentricity is the Vice Chancellor. Committees with high eccentricity were those who shared the fewest members with any other committee. 

One of the most interesting ideas is that of a hole. In this case, committees that share members form a loop which can't be reconciled. Someone with connections to all of the committees must plug up the hole for business to get done. In this case, the only candidate is the Vice Chancellor, but the rest of the individuals see this as a power grab. It would be interesting to compare this process to that of any large organization.

\subsection{Rules for Success in the Committee Structure}

In general a hole is plugged with a "psuedo-committee" at the next highest level. Here are Atkins two rules:
\begin{enumerate}
\item To {\it succeed} in politics, inject $(N-k)$-level business into $N$-level committees.
\item To {\it thwart} in politics, refer $(N+k)$-level business to $N$-level committees. 
\end{enumerate}

In the latter case, you can purposely create holes. If no such $N$-level committees exist, you can create them and call them "Working Parties." (Working Groups in our terminology.) 

Atkin gives some very specific case studies that are heavily dependent on how his University was structured but he saw a major disconnect between levels $N+2$ (Univeristy Council) and $N+1$ (Committees) as compared with any level below that such as the departments. He refers to levels $N+1$ and $N+2$ as "the citadel". (The other Colchester Castle?) 

Since \cite{Atk} was published, there have been some other papers published on Q-Analysis, but it seems to be mostly missing in modern day TDA discussions. Atkin, himself, didn't seem to publish anything more recently than the early 1980's. Still I think he has some interesting things to say about the workings of a large organization. Case studies of large corporations or Government Agencies might be interesting to see if they behave in similar ways. My guess is that it is very probable.

\section{Sensor Coverage}

I will now move ahead by 30 years and describe one of the first modern applications of topological data analysis. The problem of sensor coverage was described by Vin de Silva and Robert Ghrist, two pioneers in TDA, in \cite{DG}. Suppose we have an array of sensors located at unknown points. (The problem would be trivial if we knew their locations.) Each sensor covers an area which is a closed ball of a fixed radius. Although we don't know the exact locations of the sensors, the sensors can detect neighboring sensors if they are within a certain distance. The question is whether there is complete coverager or is there a hole in the region covered by the sensors. Algebraic topology will come to the rescue and provide a method for solving this problem. In this section, I will summarize \cite{DG} and how it solves this problem. For a more modern approach, see the book by Robinson \cite{Rob}.

We make six assumptions:
\begin{enumerate}
\item Nodes cover a closed ball of covering radius $r_c$.
\item Nodes broadcast unique ID numbers. Each node can detect the identity of any node within radius $r_s$ via a strong signal or within a larger radius $r_w$ via a weak signal.
\item The radii of communication $r_s$ and $r_w$ and the covering radius $r_c$ satisfy $$r_c\geq r_s/\sqrt{2}\mbox{ and }r_w\geq r_s\sqrt{10}.$$
\item Nodes lie in a compact domain $D\subset R^d$. Nodes can detect the presence but not the direction of the boundary $\partial D$ within a fence detection radius $r_f$. 
\item The domain $D-N_{\hat{r}}(\partial D)$ is connected, where $$N_{\hat{r}}(\partial D)=\{x\in D| ||x-\partial D||\leq\hat{r}\}\mbox{ and }\hat{r}=r_f+r_s/\sqrt{2}.$$ Here $||x-\partial D||$ denotes the distance from $x$ to the closest point of $\partial D$.
\item The outer boundary can not exhibit large-scale wrinkling. (The actual statement involves differential geometry terms that I don't want to get into here as that will take us too far afield and are not that relevant in what follows. See \cite{DG} for more details.)
\end{enumerate}

Note that none of these constants mentioned in the assumptions are dependent on the dimension $d$.

Let $X$ be the set of nodes and $X_f$ be the set of {\it fence nodes} which lie within the radius $r_f$ of $\partial D$. We build a collection of graphs: 

$$\begin{tikzpicture}
  \matrix (m) [matrix of math nodes,row sep=3em,column sep=4em,minimum width=2em]
  {
R_s &  R_w\\
F_s &  F_w\\};
  \path[-stealth]

(m-2-1) edge node [left] {$\subset$} (m-1-1)
(m-1-1) edge node [above] {$\subset$} (m-1-2)
(m-2-1) edge node [above] {$\subset$} (m-2-2)
(m-2-2) edge node [right] {$\subset$} (m-1-2)

;

\end{tikzpicture}$$

The graphs $R_s$ and $R_w$ have vertices which are the nodes in $X$ and an edge between 2 nodes if the distance between them is less than $r_s$ and $r_w$ respectively. These are the {\it communication} graphs for the strong and weak signals respectively. The {\it fence subgraphs} $F_s$ and $F_w$ are the subgraphs of $R_s$ and $R_w$ in which we restrict the nodes to the subset $X_f\subset X$. 

For these four graphs we have the corresponding Vietoris-Rips complexes with radius $r_s$ for $R_s$ and $F_s$ and  $r_w$ for $R_w$ and $F_w$. We write corresponding the nested simplicial complexes as:

$$\begin{tikzpicture}
  \matrix (m) [matrix of math nodes,row sep=3em,column sep=4em,minimum width=2em]
  {
\mathcal{R}_s &  \mathcal{R}_w\\
\mathcal{F}_s &  \mathcal{F}_w\\};
  \path[-stealth]

(m-2-1) edge node [left] {$\subset$} (m-1-1)
(m-1-1) edge node [above] {$\subset$} (m-1-2)
(m-2-1) edge node [above] {$\subset$} (m-2-2)
(m-2-2) edge node [right] {$\subset$} (m-1-2)

;

\end{tikzpicture}$$

The sensor cover $\mathcal{U}$ is the union of discs of radius $r_c$ (see assumption 1). Here is the main result of \cite{DG}:

\begin{theorem}
For a fixed set of nodes in a domain $D\subset R^d$ satisfying our six assumptions, the sensor cover $\mathcal{U}$ contains $D-N_{\hat{r}}(\partial D)$ if the homomorphism $$i_*: H_d(\mathcal{R}_s, \mathcal{F}_s)\rightarrow H_d(\mathcal{R}_w, \mathcal{F}_w)$$ induced by inclusion is nonzero.
\end{theorem}

The proof has a lot to do with the geometry of the balls given all the conditions on the radii. I will outline it by mainly focusing on the role of homology. 

Here is where persistence plays a role. The criterion for coverage is equivalent to saying that there is a nonzero generator for the relative homology group $ H_d(\mathcal{R}_s, \mathcal{F}_s)$ that {\it persists} to $H_d(\mathcal{R}_w, \mathcal{F}_w)$. Remember that $r_w>r_s$. 

Assumptions 5 and 6 turn out to insure that $H_d(D, N_{\hat{r}}(\partial D)$ has betti number $\beta=1.$ Also, $H_d(\mathcal{U}\cup N_{\hat{r}}(\partial D),  N_{\hat{r}}(\partial D))$ is nonzero if and only if $\mathcal{U}$ contains $D- N_{\hat{r}}(\partial D)$. Now the complex that would capture the topology of $\mathcal{U}$ exactly would be the \u{C}ech complex. But as we saw in Section 5.1, the \u{C}ech complex is very hard to compute as opposed to the Vietoris-Rips complex. So we would like it to be true that with the Vietoris-Rips complex, $ H_d(\mathcal{R}_s, \mathcal{F}_s)$ is nonzero if and only if $\mathcal{U}$ contains $D- N_{\hat{r}}(\partial D)$.

\begin{figure}[ht]
\begin{center}
  \scalebox{0.4}{\includegraphics{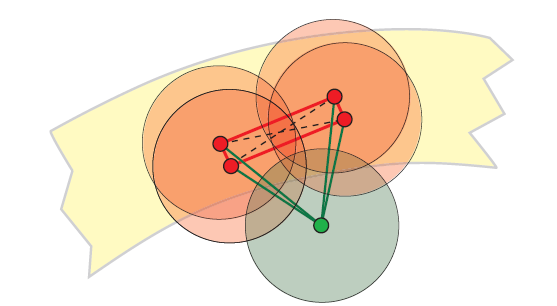}}
\caption{
\rm A generator of $H_2(\mathcal{R}_s, \mathcal{F}_s)$ which is killed by the inclusion $i_*$ into $H_2(\mathcal{R}_w, \mathcal{F}_w)$  The fence nodes are on top and the strip is a border of radius $r_f$. \cite{DG}
}
\end{center}
\end{figure}

Unfortunately, this is not always the case. Consider the situation in Figure 6.2.1. In this case $d=2$, and there is a cycle of points in $\mathcal{F}_s$ which are attached to a single vertex in $\mathcal{R}_s-\mathcal{F}_s$. The cycle has two edges of length $r_s$ and two short edges of length $\epsilon\ll r_s$. Both diagonals of the top rectangle are longer than $r_s$ so they are not present in $\mathcal{F}_s$. This corresponds to a generator in $H_2(\mathcal{R}_s, \mathcal{F}_s)$ which does not imply coverage of the entire domain. When the radius is increased to $r_w$, the diagonals appear and kill the relative 2-cycle by filling in the rectangle. So the image of this fake class is the zero element of $H_2(\mathcal{R}_w, \mathcal{F}_w).$ Now it is possible to produce a new fake class so $H_2(\mathcal{R}_w, \mathcal{F}_w)$ is not necessarily 0, but the original fake class was killed by $i_*$ so we know that there was a gap in coverage. 

The bulk of the work is in analysis of the coverage near the boundary $\partial D$. The proof is lenghty and technical so you can see \cite{DG} for details. The key point, though, is that the homology calculations find holes that the  \u{C}ech complex would pick up but the Vietoris-Rips complex misses. Finally, I will mention that the bounds on the radii can be tightened by introducing a dependence on the dimension $d$ of the domain. Then we have $$r_c\geq r_s\sqrt{\frac{d}{2(d+1)}} \mbox{ and } r_w\geq r_s\sqrt{\frac{7d-5+2\sqrt{2d(d-1)}}{d+1}},$$ and assumption 5 says that the domain $D-N_{\hat{r}}(\partial D)$ is connected, where $$N_{\hat{r}}(\partial D)=\{x\in D| ||x-\partial D||\leq\hat{r}\}\mbox{ and }\hat{r}=r_f+r_s\sqrt{\frac{d-1}{2d}}.$$

\section{Mapper}

Ayasdi is a company specializing in Topological Data Analysis. It was started in 2008 by Gunnar Carlsson, Gurjeet Singh, and Harlan Sexton of Stanford University. One of their biggest commercial products is a visualization algorithm called {\it Mapper} introduced in \cite{SMC}. In this section, I will describe the algorithm and how it works. In Section 6.3.1, I will give the general construction. Section 6.3.2 describes the specific implementation for point cloud data. The construction involves level curves of a height function, so possible height functions are dealt with in Section 6.3.3. Finally, Section 6.3.4 describes an interesting use case relating to cancer research. In \cite{NLC}, Nicolau, Levine, and Carlsson, use Mapper to find a previously undiscovered type of breast cancer cell. Patients with this type of cell exhibited 100\% survival and no metastasis. I will describe how Mapper was used to discover this type of cell.

Mapper software is available in open source versions in both R and Python as well as a version sold by Ayasdi which requires a license.

\subsection{Construction}

\begin{definition}
Let $X$ be a topological space and $\mathcal{U}=\{U_\alpha\}_{\alpha\in A}$ be a cover of $X$ by a finite collection of open sets indexed by a finite indexing set $A$. We define a simplicial complex $N(\mathcal{U})$ called the {\it nerve}\index{nerve} of $\mathcal{U}$ whose vertex set are the elements of $A$ and $\{\alpha_0, \alpha_1, \cdots, \alpha_k\}$ for $\alpha_i\in A$ span a simplex of $N(\mathcal{U})$ if and only if $U_{\alpha_0}\cap U_{\alpha_1}\cap\cdots\cap U_{\alpha_k}\neq\emptyset.$
\end{definition}

Now consider a continuous function $f: X\rightarrow Z$, where $X$ is a topological space and $Z$ is a {\it parameter space.} Usually $Z\subseteq R^n$ for some $n$. Now suppose $\mathcal{V}=\{V_\alpha\}_{\alpha\in A}$ is an open cover of $Z$. Let $f^{-1}(\mathcal{V})=\{f^{-1}(V_\alpha)\}_{\alpha\in A}$. Then since $f$ is continuous, $f^{-1}(\mathcal{V})$ is an open cover of $X$. Given $X$, $Z$, $f$, and $\mathcal{V}$, Mapper returns the nerve $N(f^{-1}(\mathcal{V}))$. We will illustrate this with 2 examples.

\begin{example}
Suppose $X=[-N, N]\subset R$, the parameter space is $Z=[0, +\infty)$, and $f: X\rightarrow Z$ is the absolute value of the quantile function of a Gaussian probability density function with mean 0 and standard deviation $\sigma$. In other words, if $p\in [0, 1]$, $f(p)$ s the absolute value of the number $q$ such that the cumulative distribution function of $q$ has value $p$. (This is a change from \cite{SMC} as just using the probability density function of the Gaussian would not have the stated properties.) The covering $\mathcal{V}$ of $Z$ consists of the four sets $[0, 5)$, $(4, 10)$, $(9, 15)$, and $(14, +\infty)$. Assume that $N$ is large enough so that $f(N)>14$. Then $f^{-1}([0, 5))$ has a single component, but $f^{-1}((4, 10))$, $f^{-1}((9, 15))$, and $f^{-1}((14, +\infty))$ all have 2 components, one on the positive side and one on the negative side. We represent the simplicial complex produced by the nerve of $f^{-1}(\mathcal{V})$ as shown in Figure 6.3.1. The nodes are labled by color and size. The color indicates the value of the function $f$ (red being high and blue being low) at a representative point in the indicated subset of $V\subset Z$. (We could take an average over $V$.) The size represents the number of points in the set corresponding to the node in the case that $X$ is discrete. (As it will be in the next section). In the continuous case, it corresponds to the integral of $f$ over the represented subset of $X$.
\begin{figure}[ht]
\begin{center}
  \scalebox{0.8}{\includegraphics{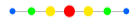}}
\caption{
\rm
Simplicial Complex for Example 6.3.1 \cite{SMC}. 
}
\end{center}
\end{figure}
\end{example}

\begin{example}
Suppose $X=\{(x, y)|x^2+y^2=1\}\subset R^2$ be the unit circle in the plane. Let $Z=[-1, 1]$, and $f(x, y)=x$. The covering $\mathcal{V}$ of $Z$ consists of the three sets $[-1, -\frac{1}{3}]$, $[-\frac{1}{2}, \frac{1}{2}]$, and $[\frac{1}{3}, 1]$.Then $f^{-1}([-1, -\frac{1}{3}])$ and $f^{-1}([\frac{1}{3}, 1])$ have one component each, while $f^{-1}([-\frac{1}{2}, \frac{1}{2}])$ has two components. The corresponding simplicial complex is shown in Figure 6.3.2. (Note that I have modified the example in \cite{SMC} in two ways. First of all, I have replaced $\frac{2}{3}$ in the paper by $\frac{1}{3}$ so that the intervals overlap. Otherwise, the simplicial complex would be four disconnected vertices. Also, I had $f(x, y)=x$ rather than $y$ so the picture would be more intuitive. If $f(x, y)=y$, the figure should be flipped with the two green nodes on the sides and the red node on the top.

\begin{figure}[ht]
\begin{center}
  \scalebox{1.3}{\includegraphics{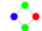}}
\caption{
\rm
Simplicial Complex for Example 6.3.2 \cite{SMC}. 
}
\end{center}
\end{figure}
\end{example}

\subsection{Implementation for Point Cloud Data}

In applying Mapper to point cloud data, we move from what Singh, et. al. call the {\it topological version} to what they call the {\it statistical version}. The main difference is that we use clustering to represent the geometric operation of partitioning a set into its connected components. 

Assume we have a point cloud $X$ with $N$ points. We have a function $f: X\rightarrow R$ on the points called the {\it filter}. Let $Z$ be the range of $f$. Divide this range into a set $S$ of intervals which overlap. Let $L$ be the length of these intervals and $p$ be the fraction of $L$ defining the overlap. For each interval $I_j\in S$, let $X_j=\{x|f(x)\in I_j\}$. Then the sets $X_j$ cover $X$. For each $X_j$, we find clusters $X_{jk}$. If $I_m$ overlaps with $I_j$, we will have overlap between clusters of $X_m$ and clusters of $X_j$. So let each cluster be a vertex in our complex and draw an edge between $X_{jk}$ and $X_{mn}$ if $X_{jk}\cap X_{mn}\neq\emptyset.$

\begin{example}
Suppose $X$ be a set of points sampled from a noisy circle in the plane. Let $Z=[0, 4.2]$, and $f(x, y)=\sqrt{(x-q_x)^2+(y-q_y)^2}$, where $q=(q_x, q_y)$ is the leftmost point in the circle. Let $L=1$ and $p=.2$. Figure 6.3.3 shows the original circle, the unclustered points, and the final complex.

\begin{figure}[ht]
\begin{center}
  \scalebox{1.3}{\includegraphics{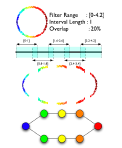}}
\caption{
\rm
Simplicial Complex for Example 6.3.3 \cite{SMC}. 
}
\end{center}
\end{figure}
\end{example}

Any clustering algorithm can be used, but Singh, et. al. use {\it single linkage clustering} \cite{Joh}\cite{JD}. The advantages of this scheme are that its input can be any distance matrix so the points aren't restricted to lie in Euclidean space, and we don't need to know the number of clusters in advance. The algorithm returns a vector $C\in R^{N-1}$, where $N$ is the size of our point cloud. At the start, each point is its own cluster. At each step, we merge the two closest clusters where the distance between clusters $A$ and $B$ is the smallest distance between a point of $A$ and a point of $B$. This distance at each step is held in $C$. We histogram these distances and look for a natural break. Longer distances should separate natural clusters. 

Another issue is that if the filter function is real valued, our complexes are 1-dimensional making them graphs. We may want to look at a higher dimensional simplical complex. This is achieved if we have $R^m$ for $m>1$ as our parameter space. We would then need to have an open cover of our range in $R^m$. For example, we could cover $R^2$ with overlapping rectangles. Then if clusters corresponding to three regions overlapped, we would add a 2-simplex and if there were four regions, we would add a 3-simplex.

\subsection{Height Function Examples}
Singh, et. al. give some examples of filter functions, but any real valued function will work. (Recall that filter and height functions are basically the same thing.)  Here are two interesting examples.

\begin{example}
{\bf Density:} Let $\epsilon>0$, and $$f_\epsilon(x)=C_\epsilon\sum_y \exp\left(\frac{-d(x,y)^2}{\epsilon}\right),$$ where $x, y\in X$, and $C_\epsilon$ is a constant such that $\int f_\epsilon(x)dx=1.$ The function $f_\epsilon$ estimates density and $\epsilon$ controls the smoothness of the estimate. This gives important features of the point cloud. Several other methods of density estimation are explained in \cite{Sil}.
\end{example}

\begin{example}
{\bf Eccentricity:} This is another family of functions that explore the geometry of the set. Eccentricity finds points that are far from the center point and assigns them higher values. On a Gaussian distribution, eccentricity increases when density decreases. Given $p$ with $1\leq p<+\infty$, we set $$E_p(x)=\left(\frac{\sum_{y\in X}d(x,y)^p}{N}\right)^{\frac{1}{p}},$$ and $$E_\infty(x)=\max_{y\in X}d(y, x)$$ where $x\in X$ and $N$ is the number of points in $X$.
\end{example}

\subsection{Application to Breast Cancer Research}

In this section I will summarize the paper of Nicolau, et. al. \cite{NLC} which uses Mapper to identify a previously unknown type of breast cancer cell. This paper was a big historical breakthrough, demonstrating Mapper's usefulness. As the biology is a little outside of my own personal expertise, I will have to use some terms without a rigorous definition and refer you to the paper and its references. What I want you to take from this is the shape of the graph produced by Mapper had a long tail that required an explanation. This tail would not have been visible using standard clustering or common dimensionality reduction techniques such as principal component analysis (PCA) or multidimensional scaling.The procedure outined in \cite{NLC} is called {\it progression analysis of disease} or PAD, and is available as a web tool either with or without Mapper. 

The first step is application of a technique called {\it disease-specific genomic analysis} (DSGA) \cite{NTBJ}, which takes a vector $\overrightarrow{T}$ representing genomic data and represents it as the sum $$\overrightarrow{T}=N_c.\overrightarrow{T}+D_c.\overrightarrow{T},$$ where $N_c.\overrightarrow{T}$ is the {\it normal} component representing healthy tissue and $D_c.\overrightarrow{T}$ is the {\it disease} component. This is done using a method of making high dimensional data less sparse followed by PCA. The reader is referred to \cite{NTBJ} for the details. The idea is to use the disease component representing the deviation from health vector to perform further analysis.

I will now show how to apply Mapper. Note that we could do the same steps on a wide variety of data matrices. It will turn out, though that the Mapper graph produced at the end will have biological significance. Start with a matrix whose columns are patients and whose rows are any genomic variables from a breast cancer microarray gene expression data set. The columns consist of tumor data vectors $\overrightarrow{T_1}, \overrightarrow{T_2}, \cdots, \overrightarrow{T_m}$ and normal tissue vectors $\overrightarrow{N_1}, \overrightarrow{N_2}, \cdots, \overrightarrow{N_k}$. We perform the following steps: \begin{enumerate}
\item Using DSGA, transform all of the data and construct the following two matrices: (i) $DC.mat$, the matrix whose columns $Dc.\overrightarrow{T_1}, Dc.\overrightarrow{T_2}, \cdots, DC.\overrightarrow{T_m}$ are the disease components of the original tumor vectors $\overrightarrow{T_1}, \overrightarrow{T_2}, \cdots, \overrightarrow{T_m}$. (ii) $L1.mat$, the matrix whose columns  $L1.\overrightarrow{N_1}, L1.\overrightarrow{N_2}, \cdots, L1.\overrightarrow{N_k}$ are an estimate of the disease component of normal tissue. These matrices are concatenated to form the matrix $L1Dc.mat$.
\item Threshold data coordinates so that only the genes that show a significant deviation from the healthy state are retained in the data matrix from Step 1. Any appropriate test for significance can be used. 
\item We apply several filter functions to the data points represented by the columns of the matrix. Given a column $v=(v_1,\cdots,v_n)$, the filter function $f_{p,k}(v)$ is defined as $$f_{p,k}(v)=\left(\sum_{i=1}^n|v_i|^p\right)^\frac{k}{p}.$$ Note that for $k=1$, the filter function is just the $L^p$ norm. 
\item Apply Mapper to the columns of the matrix using the filter functions defined in step 3. For the clustering, use the Pearson Correlation Distance $$d_{cor}(x, y)=1-\frac{\sum_{1=1}^n(x_i-\bar{x})(y_i-\bar{y})}{\sqrt{\sum_{i=1}^n(x_i-\bar{x})^2\sum_{i=1}^n(y_i-\bar{y})^2}}.$$
\end{enumerate}
 
Nicolau, et. al. used filters with $p$ ranging from 1 to 5 and $k$ ranging from 1 to 10. Figure 6.3.4 shows the Mapper graph for $f_{p=2, k=4}.$

\begin{figure}[ht]
\begin{center}
  \scalebox{1.0}{\includegraphics{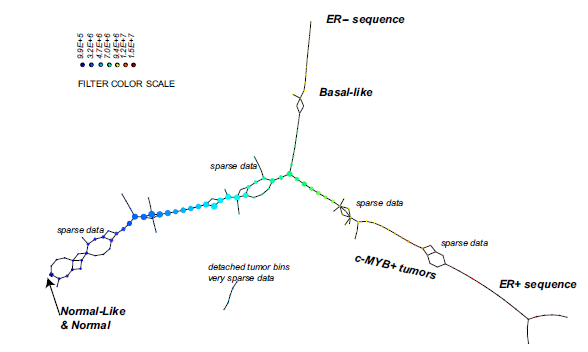}}
\caption{
\rm
Simplicial Complex for Breast Cancer data with filter $p=2$, $k=4$. \cite{NLC}. 
}
\end{center}
\end{figure}

In the figure, the main feature is the long tail at the bottom right. This  corresponds to estrogen receptor positive tumors (ER+) with high levels of the c-MYB gene. This is denoted as c-MYB+ and was represented 7.5\% of the patients. It was newly discovered through this analysis and regular clustering methods would not have discovered it. This group had a 100\% survival rate with no recurrence of disease. See \cite{NTBJ} for more details.

\section{Simplicial Sets}

It will now be necessary to take a long digression in order to understand the UMAP algorithm in the next section. It relies heavily on the idea of {\it simplicial sets}. Simplicial sets, a generalization of simplicial complexes, are a key part of UMAP. In addition, they form the basis for combinatorial homotopy theory which we will need in order to apply obstruction theory to data science.

The idea is as follows: Recall from Chapter 4 that if $K$ and $L$ are simplical complexes, then a {\it simplicial map} $f: K\rightarrow L$ takes vertices $v_i$ in $K$ to vertices $f(v_i)$ in $L$. Then for any point $x\in K$ with barycentric coordinates $\sum t_iv_i$, we have $f(x)=\sum t_if(v_i)$. In other words, a simplicial map is completely determined by its action on the vertices of $K$. Now for distinct $v_i$ we don't make a condition that $f(v_i)$ be distinct. So two vertices in $K$ can be taken to the same vertex in $L$. As an easy example, let $K$ consist of the 2-simplex $[v_0, v_1, v_2]$ and let $L$ be the face $[v_0, v_1]$. Suppose $f(v_0)=v_0$ and $f(v_1)=f(v_2)=v_1.$ This is a map that {\it collapses} $K$ onto one of its faces. We could like to somehow distinguish the simplex $[0, 1]$ as the image of a 2-simplex. We can do this by saying that this simplex is {\it degenerate} and represent it as $[v_0, v_1, v_1]$.

The standard reference for simplicial sets is the comprehensive book by Peter May \cite{May2}. Another reference is \cite{Curt}. These books can be quite difficult for someone new to the subject as they are entirely combinatorial and have little in the way of pictures. To remedy this situation, Greg Friedman, wrote a paper \cite{Fri} with a lot of pictures that is a much gentler introduction to the subject. My previous example came from there, and I will outline his approach  in the remainder of this section. 

\subsection{Ordered Simplicial Complexes}
In what follows, it will be easier to restrict ourselves to {\it ordered simplicial complexes}\index{ordered simplicial complexes} in which the vertices all follow some order. We write a simplex as $[v_{i_0}, \cdots, v_{i_k}]$ if $v_{i_m}<v_{i_n}$  whenever $m<n$. When the entire complex consists of a single simplex, we will talk about the {\it standard (ordered) n-simplex} written as $|\Delta^n|$.\index{standard (ordered) n-simplex}\index{$|\Delta^n|$} In the literature, it is common to write $|\Delta^n|$ as $[0, 1, \cdots, n]$. Then a face of $[0, 1, \cdots, n]$ is of the form $[i_0, i_1, \cdots, i_k]$ where $0\leq i_0<i_1<\cdots<i_k\leq n.$ For each simplex in an ordered simplical complex, there is exactly one way to represent the vertices in this form, and we can think of it is the image of a standard $k$-simplex $|\Delta^k|=[0, 1, \cdots, k].$

Now we would like a way to assign a map to an $n$-simplex that will give us an $(n-1)$-dimensional face.

\begin{definition}
If $[0, 1, \cdots, n]$ is the standard $n-$simplex, then there are $n+1$ {\it face maps}\index{face map} $d_0, d_1, \cdots, d_n$ defined as $d_j([0, 1, \cdots, n])=[0, 1, \cdots, j-1, j+1, \cdots, n]$, i.e. we remove vertex $j$. (Here I am using the notation of \cite{Fri}. In \cite{May2} the more common notation $\partial_j$ is used.)
\end{definition}

\begin{figure}[ht]
\begin{center}
  \scalebox{0.8}{\includegraphics{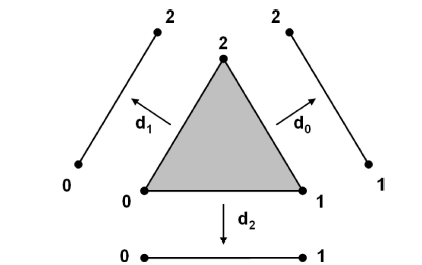}}
\caption{
\rm
Face Maps of $|\Delta^2|$ \cite{Fri}. 
}
\end{center}
\end{figure}

Figure 6.4.1 shows the face maps of $|\Delta^2|$. Note that these maps are assignments rather than continuous maps. In general, we define $d_0, d_1, \cdots, d_n$ on the n simplices of an ordered simplicial complex, where $d_j$ removes the $j$-th vertex in the ordering. We can get more general by allowing composition of face maps to produce a face of dimension less than $n-1$. We just need to have $j$ be less than or equal to $k$ for a $k$-simplex so that the definition still makes sense.

\begin{example}
Let $n=6$ and compute $d_3d_1d_5([0, 1, 2, 3, 4, 5, 6])$. We have $$d_3d_1d_5([0, 1, 2, 3, 4, 5, 6])=d_3d_1([0, 1, 2, 3, 4, 6])=d_3([0, 2, 3, 4, 6])=[0, 2, 3, 6].$$ Note that in the last step we removed vertex 4 instead of 3 since by removing vertex 1, vertex 4 moved into slot 3 (counting from 0).
\end{example}

\begin{example}
Let $n=5$. We have $$d_2d_4([0, 1, 2, 3, 4, 5])=d_2([0, 1, 2, 3, 5])=[0, 1, 3, 5],$$ while $$d_3d_2([0, 1, 2, 3, 4, 5])=d_3([0, 1, 3, 4, 5])=[0, 1, 3, 5].$$ We see that these are equal and in general for $i<j$, removing $i$ first will move $j$ to the $j-1$ spot so that $d_id_j=d_{j-1}d_i$. 
\end{example}

\subsection {Delta Sets}

Before defining simplicial sets, we will look at a concept that is one level of abstraction beyond a simplicial complex.

\begin{definition}
A {\it Delta set}\index{Delta set} (also called a {\it semisimplicial complex}\index{semisimplicial complex}, {\it semisimplicial set}\index{semisimplicial set}, or a $\Delta${\it -set}\index{$\Delta$-set} ) consists of a sequence of sets $X_0, X_1, \cdots$, and for each $n\geq 0$, maps $d_i: X_{n+1}\rightarrow X_n$ for each $i$ with $0\leq i\leq n+1$ such that $d_id_j=d_{j-1}d_i$ whenever $i<j$. (Note that we capitalize Delta to reflect the fact that the Greek letter $\Delta$ is capitalized.)
\end{definition}

The difference between Delta sets and ordered simplicial complexes is that the faces may not be unique. This reflects the fact that we may want to glue faces together. Friedman gives two examples.

\begin{figure}[ht]
\begin{center}
  \scalebox{0.8}{\includegraphics{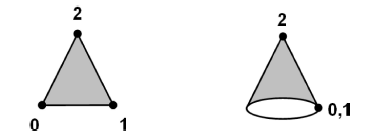}}
\caption{
\rm
$|\Delta^2|$ glued into a cone \cite{Fri}. 
}
\end{center}
\end{figure}

\begin{example}
Figure 6.4.2 shows a 2-simplex $|\Delta^2|$ with the faces $[0, 2]$ and $[1, 2]$ glued together to form a cone. The result is not actually a simplicial complex as the two faces of [0, 1] are not unique. We have $X_0=\{[0], [2]\}$, $X_1=\{[0, 1], [0, 2]\}$ and $X_2=\{[0, 1, 2]\}$. Then the face maps are what we would expect but $d_0([0,1])=[1]=[0]=d_1([0, 1])$. Check that the relations $d_id_j=d_{j-1}d_i$ for $i<j$ still hold in this case. 
\end{example}

\begin{figure}[ht]
\begin{center}
  \scalebox{0.4}{\includegraphics{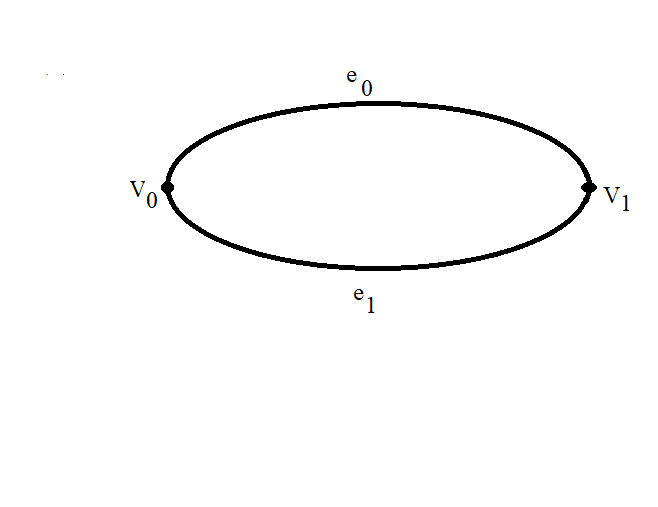}}
\caption{
\rm
A Delta set with two vertices bounding two different edges (based on example from \cite{Fri}). 
}
\end{center}
\end{figure}

\begin{example}
It is possible to have a set of vertices define two different simplices as Figure 6.4.3 shows. Let $e_0$ and $e_1$ be the top and bottom edges respectively. Then $X_0=\{v_0, v_1\}$ and $X_1=\{e_0, e_1\}$. So $d_0(e_0)=d_0(e_1)=v_1$ and  $d_1(e_0)=d_1(e_1)=v_0$. Note that it is the face maps which determine what is glued together.
\end{example}

Now we will give an alternate definition of Delta sets involving category theory. You will need to understand this approach in order to understnad the UMAP paper. First we define a category $\hat{\Delta}$.

\begin{definition}
The category $\hat{\Delta}$ has as objects the sets $[n]=[0, 1, 2, \cdots, n]$, and morphisms the strictly order preserving functions $f: [m]\rightarrow [n]$ for $m\leq n$. (A function $f$ is strictly order preserving if $i<j$ implies $f(i)<f(j)$.)
\end{definition}

Now we look at the dual or opposite category $\hat{\Delta}^{op}$. This category contains the objects $[n]$ but the morphisms go in the opposite direction. Rather than include a face of dimension $[m]$ in a simplex of dimension $[n]$ like the morphisms in $\hat{\Delta}$, we go from a higher dimensional simplex and map it to a particular face. But this is exactly what the face maps in a Delta set do. If $n-m>1$, we just have a composition of face maps. So here is the alternate definition: 

\begin{definition}
A Delta set is a covariant functor $X: \hat{\Delta}^{op}\rightarrow {\bf Set}$, where {\bf Set} is the category of sets and functions. Equivalently, a Delta set is a contravariant functor $\hat{\Delta}\rightarrow {\bf Set}$.
\end{definition}

In this formulation, we think of the functor $X$ of taking the standard $n-$simplex $[n]$ to the collection of all $n-$dimensional simplices of our Delta set. The formula  $d_id_j=d_{j-1}d_i$ for $i<j$ comes from the fact that it holds on $[n]$.

\subsection{Definition of Simplicial Sets}

Recall our cone in Example 6.4.3.  Under the gluing map, the image of $|\Delta^2|$ should still be 2-dimensional (the surface of the cone) but it only has 2 distinct vertices, $[0]=[1]$ and $[2]$. This is an example of a {\it degenerate} simplex. You can think of a degenerate simplex as one without the right number of vertices. But how do we tell it from a 1-dimensional simplex? The answer is that we allow vertices to repeat. So a {\it degenerate simplex}\index{degenerate simplex} is of the form $[v_{i_0}, \cdots, v_{i_n}]$ where the vertices are not all distinct but we still have $i_k\leq i_m$ if $k<m$.  

\begin{example}
How many 1-simplices including degenerate ones are in the 2-simplex $|\Delta^2|=[0, 1, 2]$? There are 6 of them: $[0, 0]$, $[0, 1]$, $[0, 2]$, $[1, 1]$, $[1, 2]$, and $[2, 2]$. In Figure 6.4.4, Friedman illustrates the 1-simplices in $|\Delta^1|$, the 1-simplices in $|\Delta^2|$, and the 2-simplices in $|\Delta^2|$.
\end{example}

There is one more catch. $|\Delta^n|$ can have degenerate simplices of dimension greater than $n$. For example, $|\Delta^0|$ has the degenerate simplex $[0, 0, 0, 0, 0, 0, 0, 0]$. To keep track we need {\it degeneracy maps}.

\begin{definition}
If $[0, 1, \cdots, n]$ is the standard $n-$simplex, then there are $n+1$ {\it degeneracy maps}\index{degeneracy map} $s_0, s_1, \cdots, s_n$ defined as $s_j([0, 1, \cdots, n])=[0, 1, \cdots, j, j, \cdots, n]$, i.e. we repeat vertex $j$. 
\end{definition}

This idea can easily be extended to simplicial complexes and Delta sets. Any degenerate simplex can be obtained from an ordinary simplex by compositions of degeneracy maps. 

\begin{definition}
A {\it simplicial set}\index{simplicial set} consists of a sequence of sets $X_0$, $X_1$, $\cdots$, and for each $n\geq 0$, functions $d_i: X_n\rightarrow X_{n-1}$ and $s_i: X_n\rightarrow X_{n+1}$ for each $i$ with $0\leq i\leq n$ such that: 
\begin{align*}
d_id_j&=d_{j-1}d_i&\mbox{ if } i<j,\\
d_is_j&=s_{j-1}d_i&\mbox{  if } i<j,\\
d_js_j&=d_{j+1}s_j=id,&\\
d_is_j&=s_jd_{i-1}&\mbox{  if } i>j+1,\\
s_is_j&=s_{j+1}s_i&\mbox{ if } i\leq j.
\end{align*}

\end{definition}

Try some examples yourself to convince yourself that these formulas are correct.

Any ordered simplicial complex can be thought of as a simplicial set if we add the degenerate simplices. These are formed by taking any simplex in the complex and repeating any subset of the vertices as many times as we want. Also, a simplicial set is a Delta set if we throw away the degeneracy maps.

\begin{example}
The standard 0-simplex $[0]$ can be thought of a simplical set with an element in each $X_n$ for $n\geq 0$. The element in $X_n$ is $[0, \cdots, 0]$ of length $n+1$.
\end{example}

Recall that we wrote $|\Delta^n|$ for the standard $n$-simplex. We will write $\Delta^n$ with the bars removed to stand for the standard $n$-simplex thought of as a simplicial set.
\newpage
\begin{figure}[ht]
\begin{center}
  \scalebox{0.8}{\includegraphics{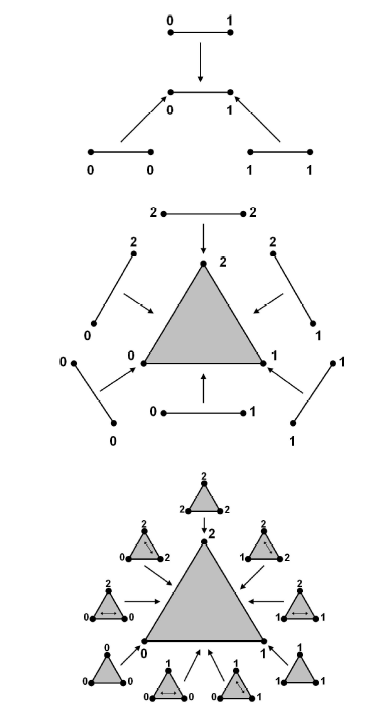}}
\caption{
\rm
The 1-simplices in $|\Delta^1|$, the 1-simplices in $|\Delta^2|$, and the 2-simplices in $|\Delta^2|$ \cite{Fri}. 
}
\end{center}
\end{figure}

\begin{figure}[ht]
\begin{center}
  \scalebox{0.8}{\includegraphics{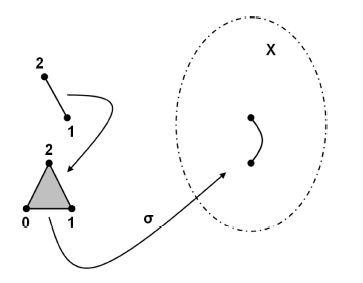}}
\caption{
\rm
Face of a singular simplex \cite{Fri}. 
}
\end{center}
\end{figure}

\begin{figure}[ht]
\begin{center}
  \scalebox{0.8}{\includegraphics{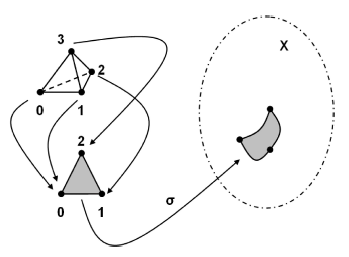}}
\caption{
\rm
Degenerate singular simplex \cite{Fri}. 
}
\end{center}
\end{figure}

\begin{example}
The prototype of simplicial sets comes from singular homology. Recall from Section 4.3.1 that if $X$ is a topological space, then a singular p-simplex is a continuous map $T: \Delta^p\rightarrow X$. We defined the group of singular $p$-chains as the free abelian group generated by these simplices. We got a face map by composing $T$ with a restriction of $\Delta^p$  to one of its faces. (We did this by mapping $\Delta^p$ linearly into $R^\infty$. Also we had followed Munkres' notation and wrote $p$ as a subscript rather than a superscript as we are doing here. The idea is the same though.) For degeneracy operators, include $\Delta^p$ into $\Delta^{p+1}$ and then include in $R^\infty$ and compose $T$ with the result. A degenerate singular simplex is a collapsed version of a higher dimensional simplex. In Figures 6.4.5 and 6.4.6, Friedman \cite{Fri} illustrates a face of a singular simplex and a degenerate singular simplex respectively.
\end{example}

\begin{definition}
A simplex $x\in X_n$ is called {\it nondegenerate} if $x$ can not be written as $s_iy$ for some $i$ and some $y\in X_{n-1}$.
\end{definition}

Every simplex of a simplicial complex or a Delta set is nondegenerate. A degenerate simplex can have a nondegenerate face. For example, any simplex $x$ whether it is degenerate or not can be written in the form $x=d_js_jx.$ It is more surprising that a nondegenerate simplex can have a degenerate face. See the next subsection for an example.

As is the case with Delta sets, simplicial sets also have a categorical definition. We start by defining a category $\Delta$.

\begin{definition}
The category $\Delta$ has objects the finite ordered sets $[n]=\{0, 1, 2,\cdots, n\}$. The morphisms are finite order preserving functions $f: [m]\rightarrow [n].$
\end{definition}

The main difference between the category $\Delta$ and the category $\hat{\Delta}$ is that the morphisms of $\Delta$ are only order preserving rather than strictly order preserving. In other words, we can repeat elements. For example, we can have $f: [3]\rightarrow [5]$ with $f([0, 1, 2, 3])=[0, 2, 2, 4].$

The morphisms in $\Delta$ are generated by $D_i: [n]\rightarrow [n+1]$ and $S_i: [n+1]\rightarrow [n]$ for $0\leq i\leq n.$ These maps are defined as $D_i[0, 1, \cdots, n]=[0, 1, \cdots, \hat{i}, \cdots, n+1]$ and $S_i[0, 1, \cdots, n+1]=[0, 1, \cdots, i, i, \cdots, n]$. 

Now we want our face and degeneracy maps to go in the opposite direction. $D_i$ includes the $i$th face of $[n]$ in $[n+1]$ so the face map $d_i$ should map $[n]$ to this face. The opposite of $S_i$ which collapses $[n+1]$ to $[n]$ by identifying the $i$th and $i+1$th vertex is the map $s_i$ which assigns to an $n$-simplex the degenerate $n+1$ simplex which repeats the $i$th vertex. This leads to the following definition.

\begin{definition}
A {\it simplicial set} is a contravariant functor $X: \Delta\rightarrow {\bf Set}$ or equivalently, a covariant  functor $X: \Delta^{op}\rightarrow {\bf Set}$.
\end{definition}

We can also think of simplicial sets as objects of a category whose morphisms consist of functions $f_n: X_n\rightarrow Y_n$ that commute with face maps and degeneracy maps.

Maps of simplicial sets, unlike the simplicial maps of Chapter 4, are not uniquely defined by what they do to vertices. They are uniquely defined, though, by what they do to nondegenerate simplices.

\subsection{Geometric Realization}

Given a simplicial set, we would like to find an actual geometric object which corresponds to it. This can be done, and the result is called the {\it geometric realization} by May \cite{May2} or just the {\it realization} by Friedman \cite{Fri}.

\begin{definition}
Let $X$ be a simplicial set. Give each set $X_n$ the discrete topology. (Recall that this means that every subset is open.) Let $|\Delta^n|$ be the standard $n$-simplex with the usual topology. The realization $|X|$ is given by $$|X|=\coprod_{n=0}^\infty(X_n\times |\Delta^n|)/\sim,$$ where $\sim$ is the equivalence relation generated by the relations $(x, D_i(p))\sim(d_i(x), p)$ for $x\in X_{n+1}$, $p\in |\Delta^n|$ and the relations $(x, S_i(p))\sim(s_i(x), p)$ for $x\in X_{n-1}$, $p\in |\Delta^n|$, and the symbol $\coprod$ means disjoint union. Here $D_i$ and $S_i$ are the face inclusions and collapses from our discussion of the category $\Delta$.
\end{definition}

This is how the definition works. First of all, we would like a an $n$-simplex for every element of $X_n$, and that is what $X_n\times |\Delta^n|$ provides. The disjoint union means we do this for every $n$. Now look at the relation $(x, D_i(p))\sim(d_i(x), p)$. Then $x$ is an $(n+1)$-simplex of $X$ and $D_i(p)$ is a point on the $ith$ face of a geometric $n+1$-simplex. $(d_i(x), p)$ is the $i$th face of that simplex together with the same point, now in an $n$-simplex. So the relation takes the n simplex corresponding to $d_i(x)$ in $X_n\times |\Delta^n|$  and glues it to the $i$th face of the $(n+1)$-simplex assigned to $x$ in $X_{n+1}\times |\Delta^{n+1}|$. The last step is to get rid of degenerate simplices. The relation $(x, S_i(p))\sim(s_i(x), p)$ takes a degenerate $n$-simplex and a point $p$ in the pre-collapse $n$-simplex $|\Delta^n|$ and glues $p$ to the $(n-1)$-simplex represented by $x$ st the point $S_i(p)$ which is the image of the collapse map. If $x$ happens to be degenerate it is also collapsed.

It turns out that the realization of a simplicial set obtained from an ordered simplicial complex is the original simplicial complex. Let $\partial\Delta^n$ denote the simplicial set obtained from the boundary $\partial|\Delta^n|$ of the ordered simplicial complex $|\Delta^n|$ by adjoining all of the degeneracies. We would like to describe $\partial\Delta^n$ as a simplicial set. One way to do this is to include any $m$-simplex $[i_0,  \cdots, i_m]$ where $0\leq i_1\leq\cdots\leq i_m\leq n$. The one condition is that it is not allowed to include all of the numbers from 0 to $n$ as it would then contain $[0, 1,\cdots, n]$ which would not be included in the boundary as it is all of $|\Delta^n|$. Since the set of simplices of the form $[i_0,  \cdots, i_m]$ whose vertices are nondecreasing and don't include all the numbers from 0 to $n$  forms the simplicial set representing the ordered simplical complex $\partial|\Delta^n|$, we have $|\partial\Delta^n|\cong S^{n-1}$.

A much more efficient way to represent $S^{n-1}$ is to only use the nondegenerate simplices $[0, 1, \cdots, n-1]$ and $[0]$. The only simplex in $X_m$ for $1\leq m<n-1$ is the degenerate simplex $[0, \cdots, 0]$ ($m-1$ times). So all of the faces of $[0, 1, \cdots, n-1]$ are $[0, \cdots, 0]$ and the resulting simplicial set is the equivalent of collapsing the boundary of an $(n-1)$-cell to a point. So this is also a representation of $S^{n-1}$. 

Note that the second construction relied on degenerate sets and could not have been done with Delta sets. It is also an example of a non-degenerate simplex with degenerate faces.

We conclude with some useful facts.The proofs are not hard and can be found in \cite{Fri}.

\begin{theorem}
A degenerate simplex is a degeneracy of a unique nondegenerate simplex. If $z$ is a degenerate simplex, there is a unique non-degenerate simplex $x$ such that $z=s_{i_1}\cdots s_{i_k}x$ for some collection of degeneracy maps $s_{i_1}, \cdots,  s_{i_k}$.
\end{theorem}

\begin{theorem}
If $X$ is a simplicial set then its realization $|X|$ is a CW complex with one $n$-cell for each nondegenerate $n$-simplex of $X$.
\end{theorem}

Our last result will be needed to understand UMAP.

\begin{definition}
Let $C$ and $D$ be categories. An {\it adjunction}\index{adjunction} consists of a pair of functors $F: C\rightarrow D$ and $G: D\rightarrow C$ together with an isomorphism $D(Fc, d)\cong C(c, Gd)$ for each $c\in C$ and $d\in D$ that is natural in both variables. (I will be more explicit below.) $F$ and $G$ are called {\it adjunct functors}\index{adjunct functors}. By {\it isomorphism}, we mean a bijection. We represent the bijection with the pair consisting of the morphism $f^\sharp: Fc\rightarrow d$ in $D$ and the morphism $f^\flat: c\rightarrow Gd$ in $C$. The pair $f^\sharp$ and $f^\flat$ are {\it adjunct} or {\it transposes} of each other. Riehl \cite{Rie}, the source of this form of the definition, shows that the natuarlity condition can be expressed as follows: In the diagram below, the left hand square commutes in category $D$ if and only if the right hand square commutes in category $C$.

$$\begin{tikzpicture}
  \matrix (m) [matrix of math nodes,row sep=3em,column sep=4em,minimum width=2em]
  {
Fc & d & c & Gd\\
Fc' &  d' &c' &  Gd'\\};
  \path[-stealth]

(m-1-1) edge node [above] {$f^\sharp$} (m-1-2)
(m-1-1) edge node [left] {$Fh$} (m-2-1)
(m-1-2) edge node [right] {$k$} (m-2-2)
(m-2-1) edge node [below] {$g^\sharp$} (m-2-2)
(m-1-3) edge node [above] {$f^\flat$} (m-1-4)
(m-1-3) edge node [left] {$h$} (m-2-3)
(m-1-4) edge node [right] {$Gk$} (m-2-4)
(m-2-3) edge node [below] {$g^\flat$} (m-2-4)

;

\end{tikzpicture}$$

\end{definition}

\begin{theorem}
In the above notation, let $C$ be the category of simplicial sets and $D$ be the category of topological spaces. Let $F: C\rightarrow D$ defined by $F(X)=|X|$ be the realization functor and $G: D\rightarrow C$ be the {\it singular set functor}\index{singular set functor} which takes a space $Y$ to the simplicial set $\mathcal{S}(Y)$ where $\mathcal{S}_n(Y)$ consists of the singular $n$-chains, i.e. the continuous maps from $|\Delta^n|$ to  $Y$. Then $F$ and $G$ are adjoint functors.
\end{theorem}

Given $f^\sharp: |X|\rightarrow Y$ and we need to produce $f^\flat: X\rightarrow \mathcal{S}(Y)$. But restricting $f^\sharp$ to a nondegenerate simplex $(x, |\Delta^n|)$ gives a continuous function $|\Delta^n|\rightarrow Y$ which is in $\mathcal{S}(Y)_n$. This gives a map from $X$ to $\mathcal{S}(Y)$ which we can define to be $f^\flat$. If $(x, |\Delta^n|)$ represents a degenerate simplex, then we precompose with the appropriate collapse map of $\Delta^n$ into $|X|$ before applying $f^\sharp$.

If we start with $f^\flat: X\rightarrow \mathcal{S}(Y)$, this assigns to each $n$-simplex $x\in X$ a continuous function $\sigma_x: |\Delta^n|\rightarrow Y$. Let $f^\sharp$ then be the continuous function that acts on the simplex $(x, |\Delta^n|)$ by applying $\sigma_x$ to $|\Delta^n|$.

I will leave it to the reader to check the other parts of the definition of adjoint functors. This result will be crucial to understanding UMAP.

\section{UMAP}

The Uniform Manifold Approximation and Projection for Dimension Reduction or UMAP algorithm of McInnes, Healy, and Melville \cite{MHM} is an algorithm for dimensionality reduction. The description is very long and difficult with a theory that combines topics from algebraic topology, differential geometry, category theory, and fuzzy set theory. In this section, I will outline the algorithm while filling in the mathematical details needed to understand the theory. I will also summarize their discussion on the comparisons to other algorithms and when UMAP  is particularly useful.

In machine learning applications, we tend to deal with high dimensional data. This makes the data difficult to visualize and a large amount of training data is needed to fully train a classifier to deal with new examples. For this reason, dimensionality reduction is needed both to aid visualization and as a preprocessing step for clustering and classification algorithms. Dimension reduction algorithms such as principal component analysis (PCA) \cite{Hot} and multidimensional scaling (MDS) \cite{Kru} preserve pairwise distances among the data points, while algorithms such as t-SNE \cite{vdMH, vdM} favor preservation of local distances over global distance. UMAP falls into the latter category. In this way, t-SNE and UMAP preserve the manifold structure of the data (for example if your data was on the surface of a torus) as opposed to PCA and MDS which squash everything into flat space. 

Following \cite{MHM}, in Section 6.5.1, I will outline the theoretical basis of UMAP. Then 6.5.2 will outline the computational approach. Section 6.5.3 will discuss weaknesses of UMAP and situations in which UMAP should not be used. Finally, Section 6.5.4 will compare UMAP to t-SNE.

Note that this is probably the hardest section in the entire book, so don't be too discouraged if you have trouble on a first reading. The material in later chapters will not be dependent on it.

\subsection{Theoretical Basis}

This section of \cite{MHM} requires an extensive background. I will try to fill in as many gaps as I can, relying on the reader having read the material in this book up to this point. For additional background on the specific constructions used, see \cite{Bar} and \cite{SpiD}. For category theory, \cite{Rie} gives you all the information you need. I will try to give you a feel for the algebraic geometry and differential geometry used but for a full background, you can see \cite{Hart} for algebraic geometry and \cite{SpiM} for differential geometry.

Recall that an $n$-dimensional {\it manifold} is a topological space $M$ such that each point $x\in M$ has an open neighborhood homeomorphic to $R^n$. So a curve is a 1-dimensional manifold and a surface is a 2-dimensional manifold. UMAP assumes that the data is uniformly distributed on some manifold. The algorithm first approximates a manifold on which the data lies and then constructs a {\it fuzzy simplicial set} (to be defined below) representation of the approximated manifold. 

If the data was uniformly distributed on $M$, then away from boundaries, a ball of fixed volume centered at any point of $M$ should contain the same number of data points. As this is not realistic for finite data, we will need some differential geometry. 

\begin{definition}
Let V be a vector space over the field $F$, where $F$ is either the real numbers $R$ or the complex numbers $C$. Then an {\it inner product}\index{inner product} on $V$ is a function $\langle\cdot, \cdot\rangle: V\times V\rightarrow F$ such that for $a, b\in F$ and $x, y, z\in V$: \begin{enumerate}
\item $\langle x, y\rangle=\overline{\langle y, x\rangle}$ (if $F=R$, then $\langle x. y\rangle=\langle y. x\rangle$),
\item $\langle ax+by, z\rangle=a\langle x, z\rangle+b\langle y, z\rangle$
\item $\langle x, x\rangle\geq 0$ and $\langle x, x\rangle=0$ if and only if $x=0$.
\end{enumerate}
The vector space $V$ is called an {\it inner product space}\index{inner product space}.
\end{definition}

If a vector space $V$ has an inner product it automatically has a norm $||x||=\sqrt{\langle x. x\rangle}.$ Then $V$ is also a metric space with $d(x, y)=||x-y||=\sqrt{\langle x-y, x-y\rangle}$. If this metric is complete, i.e. every Cauchy sequence in $V$ converges to a point in $V$, then $V$ is called a {\it Hilbert space}\index{Hilbert space}.

For an $n$-manifold $M$, let $p\in M$. The homeomorphism $\phi: U\rightarrow R^n$ where $U$ is an open subset of $M$ containing $p$ is called a {\it coordinate chart}\index{coordinate chart}. 
\begin{definition}
The tangent space $T_pM$ of an $n$-dimensional manifold $M$ at the point $p$ is defined as follows: Let $\gamma_1, \gamma_2: [-1, 1]\rightarrow M$ be two ciurves on $M$ and suppose a coordinate chart $\phi$ has been chosen. Let $\gamma_1(0)=\gamma_2(0)=p$. Then $\gamma_1$ and $\gamma_2$ are equivalent if the derivatives of $\phi\circ\gamma_1$ and $\phi\circ\gamma_2$ coincide at 0. (These are functions of one real variable so we mean derivative in the ordinary sense.) An equivalence class is called a tangent vector and the set of equivalence classes is called the {\it tangent space}\index{tangent space}. It turns out that the tangent space is independent of the choice of coordinate chart $\phi$. 
\end{definition}

At each point $p$ of an $n$ dimensional manifold $M$, the tangent space $T_pM$ is an n-dimensional vector space. We write $\phi=(x^1, \cdots, x^n): U\rightarrow R^n$. The we write the basis of $T_pM$ as $$\left\{\left(\frac{\partial}{\partial x^1}\right)_p, \cdots, \left(\frac{\partial}{\partial x^n}\right)_p\right\}.$$ Here $\left(\frac{\partial}{\partial x^i}\right)_p$ is the equivalence class of a curve $\gamma$ such that the tangent to $\phi\circ\gamma$ is in the $x^i$ direction. 

Now we let $g$ be a matrix such that $$g_{ij}=\left\langle\left(\frac{\partial}{\partial x^i}\right)_p, \left(\frac{\partial}{\partial x^j}\right)_p\right\rangle.$$ In other words, $g$ determines inner products on $T_pM$. Then $g$ is called a {\it Riemannian metric}\index{Riemannian metric} and $M$ is called a {\it Riemannian manifold}\index{Riemannian manifold}.

Here is the last fact we need before resuming with UMAP.

\begin{theorem}
{\bf Whitney Embedding Theorem:} Let $M$ be an $n$-dimensional manifold that is smooth (i.e if $U$ and $V$ are open subsets of $M$ with corresponding coordinate charts $\phi_U$ and $\phi_V$ then $\phi_U\circ\phi_V^{-1}: R^n\rightarrow R^n$ is differentiable.) Then $M$ can be embedded in a space $R^{2n}$. (I. e. there is a one to one map taking $M$ into $R^{2n}$.)
\end{theorem}

If we embed $M$ in $R^m$ then we call $R^m$ the {\it ambient space}. Now I can state the next result from \cite{MHM}.

\begin{theorem}
Let $(M, g)$ be a Riemannian manifold in an ambient $R^n$, and let $p\in M$ be a point. Let $g$ be constant diagonal matrix in a ball $B\subset U$ centered at $p$ of volume $V_n$ with respect to $g$ where $V_n$ is the volume of a ball in $R^n$ of radius 1. Explicitly, $$V_n=\frac{\pi^{\frac{n}{2}}}{\left(\frac{n}{2}\right)!}\mbox{  for } n \mbox{ even and } V_n=\frac{2^n\pi^{\left(\frac{n-1}{2}\right)}\left(\frac{n-1}{2}\right)}{n!}\mbox{ for }n\mbox{ odd.}$$ Then the {\it geodesic distance} (i.e. the shortest distance along a path in the manifold) between $p$ and another point $q\in B$ is $\frac{1}{r}d(p, q)$ where $r$ is the radius of the ball in the ambient space $R^n$ and $d(p, q)$ is the distance in the ambient space.
\end{theorem}

The point of all of this is that for each point $X_i$ in our data set, we can normalize distances so that a point and its $k$ nearest neighbors are contained in a ball of uniform volume no matter what point we take as our center. So our assumption of the points being uniformly distributed on the manifold will hold locally if we have a good way of pasting each of these local metric spaces together. That is what we will do next.

The next step involves simplicial sets and your author has prepared you for this in the previous section. There are some notational differences I need to point out. The first is the idea of a {\it colimit}. In \cite{MHM}, it is stated that for a simplicial set $X$, we have that $$colim_{x\in X_n}\Delta^n\cong X.$$ This turns out to be equivalent to definition 6.4.10. For a sequence of sets and inclusions the colimit is the disjoint union. (See \cite{Rie} Chapter 3 and especially Example 3.1.26.) We will also use the fact that the singular set functor and realization functor are adjoints. 

Next, we need to talk about {\it fuzzy set theory}\index{fuzzy set theory} \cite{Zad} or {\it fuzzy logic}\index{fuzzy logic} as opposed to faulty logic\index{faulty logic}. Normally, given a set $S$ and an element $x$, either $x\in S$ or $x\notin S$. If $S$ is a fuzzy set, there is a membership function $\mu: S\rightarrow [0, 1]$. So an element $x$ is in $S$ to a varying degree. $\mu(x)$ then is the {\it membership strength} of $x$.

\begin{example}
Consider this scenario. You are getting pretty tired of TDA by now and need to take a break. So you decide to go hiking in the Arizona desert, but unfortunately, you forgot to bring any water with you and you are starting to get pretty thirsty. Then you see two bottles standing up against a cactus. The bottles have unusual labels. The first one says that the probability of the contents being a potable liquid is .9. The second one says that the membership of the contents in the set of potable liquids is .9. What should you do?

Later, someone shows up, apologizes for the confusing messages on the bottles, and tears off the labels. Underneath are a different set of labels. The first bottle turns out to be hydrochloric acid. Somehow, it got mixed up with 9 identical looking water bottles, so the probability of it being water was .9. The second bottle turns out to be a bottle of beer. Assuming the most potable liquid is water, it is similar enough to earn a membership function of .9. There is one other big difference. Now that we have more information, the probability of bottle 1 being potable is now 0, but the fuzzy membership function of bottle number 2 is unchanged despite the additional information.
\end{example} 

We want to represent fuzzy sets in terms of category theory. To do this we need some more category theory. To understand the motivation, I will introduce the concept of a {\it sheaf} from algebraic geometry.

\begin{definition}
Let $X$ be a topological space. A {\it presheaf}\index{presheaf} $\mathcal{P}$ on $X$ consists of the following: \begin{enumerate}
\item For each open set $U\subseteq X$, a set $\mathcal{P}(U)$. The elements in this set are called the {\it sections} of $\mathcal{P}$ over $U$.
\item If $V\subseteq U$ is an inclusion of open sets, we have a corresponding function ${\bf res}_{V, U}: \mathcal{P}(U)\rightarrow\mathcal{P}(V)$ called a {\it restriction function} or {\it restriction morphism.}
\item For every open set $U\subseteq X$, the restriction ${\bf res}_{U, U}: \mathcal{P}(U)\rightarrow\mathcal{P}(U)$ is the identity. 
\item If there are three open sets $W\subseteq V\subseteq U$, then ${\bf res}_{W, V}\circ {\bf res}_{V, U}={\bf res}_{W, U}.$
\end{enumerate}
\end{definition}

For the next definition, if $s\in \mathcal{P}(U)$ is a section for the presheaf $\mathcal{P}$ over the open set $U$ and $V\subseteq U$ is open, then the notation $s|_V$ denotes the element of $\mathcal{P}(V)$ which is equal to ${\bf res}_{V, U}(s)$.

\begin{definition}
A presheaf $\mathcal{P}$ is a {\it sheaf}\index{sheaf} if it has two additional properties:\begin{enumerate}
\item {\bf Locality:} If $\{U_i\}$ is an open cover of an open set $U$, and if $s, t\in\mathcal{P}(U)$ are such that $s|_{U_i}=t|_{U_i}$ for all $i$, then $s=t$.
\item {\bf Gluing:} If $\{U_i\}$ is an open cover of an open set $U$, and if for each $i$, there is a section $s_i\in\mathcal{P}(U_i)$ such that for any pair of open sets $U_i$ and $U_j$ in the cover, $s_i|_{U_i\cap U_j}=s_j|_{U_i\cap U_j}$, then there is a section $s\in \mathcal{P}(U)$ such that $s|_{U_i}=s_i$ for all $i$. 
\end{enumerate}
\end{definition}

\begin{example}
Let $X=R$ and for $U\subset R$, $\mathcal{P}(U)$ is the set of continuous real valued functions on $R$. For $V\subseteq U$, and $f\in \mathcal{P}(U)$, ${\bf res}_{V, U}(f)$ is just the usual restriction of $f$ to $V$.
\end{example}

\begin{example}
Let $X=R$ and for $U\subset R$, $\mathcal{P}(U)$ is the set of {\bf bounded} continuous real valued functions on $R$. Then $\mathcal{P}$ is a presheaf which is not a sheaf. Let $\{U_i\}$ be an open cover of $R$ consisting of the sets $(-i, i)$ Then gluing fails. Letting $f(x)=x$, $f$ is bounded when restricted to each of the $U_i$, but $f$ is unbounded on the real line.
\end{example}

 Now notice that given the category $C$ whose objects are open subsets of $X$ and whose morphisms are inclusions, a presheaf is a functor from $C^{op}$ to {\bf Set}. Now let $I=[0, 1]$ and define a topology whose open sets are intervals of the form $[0, a)$ for $a\in (0, 1]$. Identifying $I$ with the category of subsets of this form and inclusions, we have a presheaf $\mathcal{P}$ which is a functor from $I^{op}$ to {\bf Set}. Define a fuzzy set to be the subcategory whose restrictions from $[0, a)$ to $[0, b)$ are one to one. We can think of a the section $\mathcal{P}([0, a))$ to be the set of all elements whose membership strength is at least $a$.

It turns out that these presheaves are actually sheaves. They form a category whose objects are sheaves and whose morphisms are natural transformations. We are almost ready to define the category of fuzzy sets, but we need a couple more definitions.

\begin{definition}
Suppose $C$ and $D$ are categories. A functor $F: C\rightarrow D$ is {\it full}\index{functor!full} if for each $x, y\in C$, the map $hom_C(x, y)\rightarrow hom_D(Fx, Fy)$ is surjective.
\end{definition}

\begin{definition}
Suppose $C$ and $D$ are categories. A functor $F: C\rightarrow D$ is {\it faithful}\index{functor!faithful} if for each $x, y\in C$, the map $hom_C(x, y)\rightarrow hom_D(Fx, Fy)$ is injective.
\end{definition}

\begin{definition}
A functor $F: C\rightarrow D$ that is both full and faithful is called {\it fully faithful}, If $C$ is a subcategory of $D$ and there is a fully faithful functor  $F: C\rightarrow D$ that is injective on the objects of $C$ then, $C$ is a {\it full subcategory} of $D$.\index{subcategory!full}
\end{definition}

\begin{definition}
The category {\bf Fuzz} of fuzzy sets is the full subcategory of sheaves on $I$ spanned by fuzzy sets.
\end{definition}

Now we can define fuzzy simplicial sets. Recall that a simplicial set was a  functor $X: \Delta^{op}\rightarrow {\bf Set}$. So this is actually a presheaf on $\Delta$ with values in {\bf Set}. So a fuzzy simplicial set will be a presheaf on $\Delta$ with values in {\bf Fuzz}.

\begin{definition}
The category of fuzzy simplicial sets {\bf sFuzz} is the category with objects given by functors from $\Delta^{op}$ to {\bf Fuzz} and morphisms given by natural transformations. 
\end{definition}

A fuzzy simplicial set can also be viewed as a sheaf over $\Delta\times I$ where $\Delta$ has the trivial topology and $\Delta\times I$ has the product topology. Use $\Delta^n_{<a}$ to be the sheaf given by the the representable functor of the object $([n], [0, a))$, i. e. the functor isomorphic to $Hom_(([n], [0, a)), -)$. The importance of the fuzzy version of simplicial sets is their relationship to metric spaces.

\begin{definition}
An {\it extended-pseudo-metric space} $(X, d)$ is a set $X$ and a map $d: X\times X\rightarrow [0, \infty]$ such that for $x, y, z\in X$, 
\begin{enumerate}
\item $d(x, y)\geq 0$ and $d(x, x)=0$,
\item $d(x, y)=d(y, x)$, and
\item $d(x, z)\leq d(x, y)+d(y, z)$ or $d(x, z)=\infty$.
\end{enumerate}
The difference between an extended-pseudo-metric space and a metric space is that an extended-pseudo-metric space can have infinite distances and it is possible for $d(x, y)$ to be 0 even it $x\neq y$.
\end{definition}

The category of extended-pseudo-metric spaces {\bf EPMet} has as objects extended-psuedo-metric spaces. The morphisms are {\it non-expansive maps} which are maps $f: X\rightarrow Y$ with the property that for $x_1, x_2\in X$, $d_Y(f(x_1), f(x_2))\leq d_X(x_1, x_2).$ (Remember that  the map is an {isometry} if $\leq$ is replaced by equality.) The sub-category of {\bf finite} extended-pseudo-metric spaces is denoted {\bf FinEpMet}.

Now construct adjoint functors {\sc Real} and {\sc Sing} between the categories {\bf sFuzz} and {\bf EPMet}. These correspond to the realization and singular functors from classical simplicial set theory. The functor {\sc Real} is defined on the standard fuzzy simplices $\Delta^n_{<a}$ with $$\mbox{{\sc Real}}(\Delta^n_{<a})=\left\{(t_0, \cdots, t_n)\in R^{n+1}|\sum_{i=0}^nt_i=-\log(a), t_i\geq 0 \right\}.$$ The metric on {\sc Real}$(\Delta^n_{<a})$ is inherited from $R^{n+1}$. A morphism $\Delta^n_{<a}\rightarrow\Delta^m_{<b}$ exists only if $a\leq b$ and is determined by a $\Delta$ morphism $\sigma: [n]\rightarrow [m]$. The action of {\sc Real} on such a morphism is given by the map $$(x_0, x_1, \cdots, x_n)\rightarrow \frac{\log(b)}{\log(a)}\left( \sum_{i_0\in\sigma^{-1}(0)} x_{i_0},  \sum_{i_0\in\sigma^{-1}(1)} x_{i_0}, \cdots,  \sum_{i_0\in\sigma^{-1}(m)} x_{i_0}\right).$$ This map is non-expansive since $0\leq a\leq b\leq 1$ implies that $log(b)/log(a)\leq 1.$

Extend this to a general  simplical set via colimits (the disjoint union operation from Definition 6.4.10). So {\sc Real}$(X)=colim_{\Delta^n_{<a}\rightarrow X}${\sc Real}$(\Delta^n_{<a})$.

The analog of the adjoint functor {\sc Sing} is defined for an extended-pseudo-metric space $Y$ by its action on $\Delta\times I$ with $$\mbox{{\sc Sing}}(Y): ([n], [0, a))\rightarrow hom_{\bf EPMet}(\mbox{{\sc Real}}(\Delta^n_{<a}), Y).$$

Since McInnes, et. al. are only interested in finite metric spaces, they consider the subcategory of bounded fuzzy simplicial sets, {\bf Fin-sFuzz}. Define the finite fuzzy realization functor as follows:

\begin{definition}
The functor {\sc FinReal}$: {\bf Fin-sFuzz}\rightarrow{\bf FinEPMet}$ is defined by setting $$\mbox{{\sc FinReal}}(\Delta^n_{<a})=(\{x_1, x_2, \cdots, x_n\}, d_a),$$ where $d_a(x_i, x_j)=-\log(a)$ if $i\neq j$ and $d_a(x_i, x_i)=0$. Then define $$\mbox{{\sc FinReal}}(X)=colim_{\Delta^n_{<a}\rightarrow X}\mbox{{\sc FinReal}}(\Delta^n_{<a}).$$
\end{definition}

The action of  {\sc FinReal} on a map $\Delta^n_{<a}\rightarrow\Delta^m_{<b}$ where $a\leq b$ defined by $\sigma: \Delta^n\rightarrow\Delta^m$ is given by $$(\{x_1, x_2, \cdots, x_n\}, d_a)\rightarrow (\{x_{\sigma(1)}, x_{\sigma(2)}, \cdots, x_{\sigma(n)}\}, d_b),$$ which is a non-expansive map, since if $a\leq b$, then $-\log(a)\geq -\log(b)$ so that $d_a\geq d_b$.

We have a functor {\sc FinSing}$:{\bf FinEPMet}\rightarrow {\bf Fin-sFuzz}$ defined  by $$\mbox{{\sc FinSing}}(Y): ([n], [0, a))\rightarrow hom_{\bf FinEPMet}(\mbox{{FinReal}}(\Delta^n_{<a}), Y).$$ It is then shown that {\sc FinReal} and {\sc FinSing} are adjunct functors.

The point of all of this work is to piece together the metric spaces defined by each $X_i$ in the manifold and its $k$ nearest neighbors. Each of these metric spaces can be transformed into their corresponding fuzzy simplicial set. We then take a {\it fuzzy union} of all of our simplicial sets. This means that the membership function of any element is the maximum of its membership function taken over every set involved in the union. 

\begin{definition}
Let $X=\{X_1, \cdots, X_N\}$ be a dataset in $R^n$. Let $\{(X,  d_i)\}$ for $1\leq i\leq N$ be a family of extended-pseudo-metric spaces with a common underlying set $X$ such that $d_i(X_j, X_k)=d_M(X_j X_k)-\rho$ for $i=j$ or $i=k$ and infinite otherwise where $\rho$ is the distance to the nearest neighbor of $X_i$ in $R^n$ and $d_M$ is the geodesic distance on the manifold either known beforehand or approximated by Theorem 6.5.2. Then the fuzzy topological representation of $X$ is $$\bigcup_{i=1}^n\mbox{{\sc FinSing}}((X, d_i)).$$
\end{definition}

The fuzzy simplicial set union has now merged the different metric spaces and forms a global representation of the manifold. We can now perform dimension reduction by finding low dimensional representations that match the topological structure of the source data.

Let $Y=\{Y_1, \cdots, Y_N\}$ be a set of points corresponding to the data set $X$ but being points in $R^d$ with $d<<n.$ We then compute a fuzzy set representation of $Y$ and compare it to that of $X$. In this case, we usually consider $Y$ as a subset of the manifold $R^d$. 

To compare two fuzzy sets, we need the same reference set. Given a sheaf represntation $\mathcal{P}$ we translate to classical fuzzy sets by setting $A=\cup_{a\in(0, 1]}\mathcal{P}([0, a))$ and membership function $\mu(x)=\sup\{a\in (0, 1]|x\in \mathcal{P}([0, 1))\}.$

\begin{definition}
Consider two fuzzy sets with underlying set $A$ and membership functions $\mu$ and $\nu$ respectively. The {\it cross entropy} $C$ is defined as $$C((A, \mu), (A, \nu))=\sum_{a\in A}\left(\mu(a)\log\left(\frac{\mu(A)}{\nu(A)}\right)+(1-\mu(a))\log\left(\frac{1-\mu(A)}{1-\nu(A)}\right)\right).$$
\end{definition}

The last step is to optimize the embedding $Y$ in $R^d$ with respect to fuzzy cross-entropy using stochastic gradient descent. 

McInnes, et. al. restricted their attention to the 1-skeleton of the fuzzy simplicial sets as this significantly reduced computational costs.

\subsection{Computational View}

McInnes, et. al. \cite{MHM} then provide a computational description of UMAP which I will now summarize. The fuzzy simplicial sets described above are only computationally tractable for the 1-skeleton which can be described as a weighted graph. This makes UMAP a $k$-neighbor based graph learning algorithm. (Recall that we looked at the $k$ nearest neighbors for each data point.) 

The following axioms are assumed to be true:\begin{enumerate}
\item There exists a manifold on which the data is uniformly distributed.
\item The underlying manifold is locally connected, i.e. every point in the manifold has a neighborhood base (see Section 2.2) consisting of open connected sets.
\item The main goal is to preserve the topological structure of this manifold.
\end{enumerate}

Here are the main steps:\begin{enumerate}
\item Graph Construction
\begin{enumerate}
\item Construct a weighted $k$-neighbor graph.
\item Apply a transform on the edges to use the local distance from the ambient space.
\item Deal with the inherent asymmetry of the $k$-neighbor graph.
\end{enumerate}
\item Graph Layout
\begin{enumerate}
\item Define an objective function that preserves the desired characteristics of this $k$-neighbor graph.
\item Find a low dimensional representation which optimizes this objective function.
\end{enumerate}
\end{enumerate}

To perform the graph construction, start with a collection of data points $X=\{x_1,\cdots,x_N\}$ be the input dataset with pairwise distances given by the pseudo-metric $d$. Given an input hyperparameter $k$, for each $x_i$, compute the $k$-nearest neighbors under $d$ labeled $\{x_{i1},\cdots, x_{ik}\}$. UMAP finds these neighbors using the nearest neighbor descent algorithm \cite{DML}. Now define for each $x_i$ two quantities $\rho_i$ and $\sigma_i$. We define $\rho_i$ to be the distance between $x_i$ and the nearest neighbor $x_{i_j}$ such that $1\leq j\leq k$ and $d(x_i, x_{ij})>0$. We set $\sigma_i$ to be the value such that $$\sum_{j=1}^k\exp\left(\frac{-\max(0, d(x_i, x_{ij})-\rho_i)}{\sigma_i}\right)=\log_2(k).$$ Then $\rho_i$ comes from the local connectivity constraint and ensures that $x_i$ connects to at least one other data point with an edge weight of 1. $\sigma_i$ is a normalization factor defining the Riemannian metric local to the point $x_i$.

Now define a directed graph $\overline{G}$ whose vertices are the points  $X=\{x_1,\cdots,x_N\}$. Each point $x_i\in X$ is connected by a directed edge to its $k$ nearest neighbors  $\{x_{i1},\cdots, x_{ik}\}$, where the directed edge from $x_i$ to $x_{ij}$ has a weight of $$\sum_{j=1}^k\exp\left(\frac{-\max(0, d(x_i, x_{ij})-\rho_i)}{\sigma_i}\right).$$ 

The graph $\overline{G}$ is the 1-skeleton of the fuzzy simplicial set associated to the pseudo-metric space associated to $x_i$. The weight associated to the edge is the fuzzy membership strength of the corresponding 1-simplex within the fuzzy simplicial set.

Now let $A$ be the weighted adjacency matrix of $\overline{G}$ and let $$B=A+A^T-A\circ A^T,$$ where $A^T$ is the transpose of $A$ and $X\circ Y$ is the product of two matrices formed by multiplying them componentwise. If $A_{ij}$ is the probability that the directed edge from $x_i$ to $x_j$ exists, then $B_{ij}$ is the probability that at least one of the two directed edges from $x_i$ to $x_j$ or $x_j$ to $x_i$ exists. The UMAP graph $G$ is then the undirected graph whose weighted adjacency matrix is $B$.

Now we need to describe the graph layout step. UMAP uses a {\it force directed graph layout algorithm.} A force directed layout algorithm uses a set of attractive forces along edges and a set of repulsive forces among vertices. The algorithm proceeds by iteratively applying attractive and repulsive forces at each edge or vertex. Slowly decreasing these forces guarantee convergence to a local minimum. 

In UMAP, the attractive force between two vertices $x_i$ and $x_j$ connected by an edge and sitting at coordinates ${\bf y_i}$ and ${\bf y_j}$ is determined by $$\left(\frac{-2ab||{\bf y_i}-{\bf y_j}||_2^{2(b-1)}}{1+||{\bf y_i}-{\bf y_j}||_2^2}\right)w((x_i, x_j))({\bf y_i}-{\bf y_j})$$ where $a$ and $b$ are hyperparameters and $w((x_i, x_j))$ is the weight of the edge between $x_i$ and $x_j$. The repulsive force is given by $$\left(\frac{2b}{\left(\epsilon+||{\bf y_i}-{\bf y_j}||_2^2)\right)\left(1+||{\bf y_i}-{\bf y_j}||_2^{2b})\right)}\right)(1-w((x_i, x_j)))({\bf y_i}-{\bf y_j}),$$ where $\epsilon$ is set to .001 to prevent division by zero. 

These forces are derived from gradients optimizing the edgewise cross-entropy between the weighted graph $G$ and the equivalent weighted graph $H$ constructed from the points $\{y_1,\cdots,y_N\}$.  So the idea is to position the points $y_i$ so that the weighted graph induced by them most closely approximates the graph $G$ where we measure the distance by the total cross-entropy over all edge existence probabilities. Since $G$ matches the topology of the original data, $H$ will be a good low dimensional representation of the topology of this data.

More implementation details and test results are found in \cite{MHM}. In the next two sections, we will briefly look at when you would want to use UMAP as opposed to a competing algorithm.

\subsection{Weaknesses of UMAP}

UMAP is very fast and effective for both visualization and dimesnsion reduction, but there are some situations in which UMAP should not be used.\begin{enumerate}
\item If you need interpretability of the reduced dimensions, UMAP does not give them a specific meaning. Principal Component Analysis (PCA) or the related Non-negative Matrix Factorization (NMF) are much better in this situation. For example, PCA reveals the dimensions which are the directions of greatest variability.
\item Don't use UMAP if you don't expect to have a manifold structure. Otherwise, it could find manifold structure in the noise.
\item UMAP favors local distance over global distance. If you are mainly concerned with global distances, multidimensional scaling (MDS) is better suited to this task.
\item Finally, UMAP makes a number of approximations to improve the computational efficiency of the algorithm. This can impact the results for small data sets with less than 500 points. The authors recommend that UMAP not be used in this case.
\end{enumerate}

\subsection{UMAP vs. t-SNE}

Since I have been mentioning t-SNE throughout this section as UMAP's main competitor, I will now give you an idea of what it actually does. While, \cite{vdMH, vdM} were the original sources, I will base this discussion on the short summary that appears in \cite{RB}. Later, I will compare t-SNE to UMAP and discuss the major differences.

As with UMAP, we face the problem of non-uniform density. The original algorithms was called {\it stochastic neighbor embedding} or SNE \cite{HR}.  

As before, we have a set of data points $\{x_1,\cdots,x_N\}\subset R^n$ and let $\{y_1,\cdots,y_m\}$ be a candidate for the corresponding image points in $R^k$ for $k<n$. Letting $d$ be the distance in $R^n$ and $\overline{d}$ be the distance in $R^k$, let $$p_{j|i}=\frac{e^{-\frac{d(x_i, x_j)^2}{2\sigma^2_i}}}{\sum_{k\neq i}e^{-\frac{d(x_i, x_k)^2}{2\sigma^2_i}}}$$ and  $$q_{j|i}=\frac{e^{-\overline{d}(y_i, y_j)^2}}{\sum_{k\neq i}e^{-\overline{d}(y_i, y_k)^2}}$$ where the variances $\sigma_i$ will be obtained by a process described below, and in the equation for $q_{j|i}$, all of the variances are fixed to be $\frac{\sqrt{2}}{2}$. Set $p_{i|i}=q_{i|i}=0.$

Looking at $p_{j|i}$ and $q_{j|i}$ as probability distributions, we would like them to be as close as possible where the difference will be the {\it Kullback-Leibler divergence}. We will use it in the form of the cost function $$C=\sum_i\sum_jp_{j|i}\log\frac{p_{j|i}}{q_{j|i}}.$$ By minimizing this particular function, we are doing the optimization locally for each $x_i, y_i$ pair taking into account their neighbors, then summing over all of the $x_i$. The density around a point is reflected in this cost function,

The SNE algorithm involves solving for the minimizing points $\{y_i\}$ using gradient descent. Convergence is slow and depends on good choices of $\sigma_i$. To estimate $\sigma_i$, fix a value $$P=2^{-\sum_jp_{j|i}\log_2p_{j|i}}.$$ This value is called the {\it perplexity} and can be thought of as the number of neighbors used. This is a loose estimate of the local dimension. The perplexity is typically between 10 and 100. It is an input parameter, and we solve for values of $\sigma_i$ achieving the desired perplexity. 

SNE has two major problems:\begin{enumerate}
\item The gradient descent procedure is slow and it can be hard to get it to converge. 
\item With a large number of points, visualization can be difficult as the cost function will push most points into the center of mass. 
\end{enumerate}

To address these issues, van der Maaten and Hinton proposed the {\it t-distributed stochastic neighborhood embedding} algorithm or t-SNE. SNE is modified in the following ways: 
\begin{enumerate}
\item We symmetrize $p_{j|i}$ as $$p_{ij}=\frac{p_{j|i}+p_{i|j}}{2}.$$
\item We symmetrize $q_{j|i}$ as $$q_{ij}=\frac{\left(1+\overline{d}(y_i, y_ j)^2\right)^{-1}}{\sum_{k\neq \mathfrak{l}}\left(1+\overline{d}(y_k, y_\mathfrak{l})^2\right)^{-1}}.$$ In the original definition, $q_{j|i}$ was defined using a Gaussian. This definition replaces it with a Student $t$-distribution with one degree of freedom which has more weight in the tails.
\item The cost function is now $$C=\sum_i\sum_jp_{ij}\log\frac{p_{ij}}{q_{ij}}.$$
\end{enumerate}

These modifications lead to heavier tails in the embedding space giving outliers less impact on the overall results and the compression around the center of mass is alleviated. Also, the gradient descent procedure becomes more efficient. 

Now to compare t-SNE and UMAP, we can write the UMAP cross-entropy cost function as $$C_{UMAP}=\sum_{i\neq j} \left(p_{ij}\log\left(\frac{p_{ij}}{q_{ij}}\right)+(1-p_{ij})\log\left(\frac{1-p_{ij}}{1-q_{ij}}\right)\right).$$

Although the cost functions look similar, the definitions for $p_{j|i}$ and $q_{j|i}$ are somewhat different for UMAP. In the high dimensional space, the $p_{j|i}$  are fuzzy membership functions taking the form $$p_{j|i}=e^{\left(\frac{-d(x_i, x_j)-\rho_i}{\sigma_i}\right)}.$$

The values $p_{j|i}$ are only calculated for $k$ nearest neighbors with all other $p_{j|i}=0.$ The distance $d(x_i, x_j)$ in the high dimensional space can be any distance function, $\rho_i$ is the distance to the nearest neighbor of $x_i$, and $\sigma_i$ is the normalizing factor which plays a similar role to the perplexity derived $\sigma_i$ from t-SNE. Symmetrization is carried out by fuzzy set union and can be expressed as $$p_{ij}=(p_{j|i}+p_{i|j})-p_{j|i}p_{i|j}.$$

The low dimensional similarities are given by $$q_{ij}=\left(1+a||y_i-y_j||_2^{2b}\right)^{-1}$$ where $a$ and $b$ are user defined values. UMAP's defaults give $a=1.929$ and $b=0.7915.$ Setting $a=b=1$ results in the Student t-distribution used in t-SNE.

According to $\cite{MHM}$, UMAP is demonstrably faster than t-SNE and provides better scaling. This allows for the generation of high quality embeddings of larger data sets than had ever been possible in the past.

\chapter{Some Unfinished Ideas}

In this very short chapter, I will mention some ideas that have never been fully pursued. They involve the use of algebraic topology for time series of graphs, directed simplicial complexes,  computer intrusion detection, and market basket analysis.

\section{Time Series of Graphs}

In Chapter 5, I discussed time series in general and the basics of graphs. One interesting combination of these two topics is the idea of time series of graphs. There are many problems in which we are interested in {\it network change detection}. This could involve differences in connectivity between nodes or in the weights of the edges between them. Edges can be weighted by attributes of the traffic the network carries. Typically, we want to detect drastic changes or anomalies in the networks.

To convert a series of graphs into a numerical time series, we could use a graph distance measure. Edit distance, discussed in Section 5.7, is especially useful and easy to compute. The book by Wallis, et. al  \cite{WBDK} has several other ideas. We could also get a series by tracking attributes such as number of vertices, number of edges, maximum degree, etc. As I mentioned, a graph can also be converted to a simplicial complex in various ways, so we can track Betti numbers, Euler characteristic, or dimension of the complex. We can also convert any graph into a persistence diagram or landscape and then produce a time series with any of the distance measures we discussed in Chapter 5.  

Given a time series, there is an easy trick to detecting spikes. Slide a window of size $k$ for some small number $k$. Let $\{x_n, x_{n+1}, \cdots, x_{n+k-1}\}$ be the values in this window and let $\mu_n$ and $\sigma_n$ be their mean and standard deviation. Then we can detect a spike at $x_{n+k}$ if $\frac{x_{n+k}-\mu_n}{\sigma_n}$ is above some predetermined threshold. For significant level changes, we can use control chart techniques such as {\it cusum} from statistical quality control. A good book with a lot of examples is \cite{Mon1}. These techniques are very useful if you own a widget factory and a disgruntled employee throws a wrench in the works. For more complex anomalies, you can use SAX as discussed in Chapter 5.

\section{Directed Simplicial Complexes}

Directed simplicial complexes are a  structure that I have just recently learned about. It has a very interesting application in neurology.

The {\it Blue Brain Project}\index{Blue Brain Project} \cite{BB} is a project of EPFL in Lausanne, Switzerland to fully map the brain of a mouse. A mouse's brain has about 100 million neurons and a trillion synapses. The idea is to understand the link between the neural network structure and its function. By constructing graphs that reflect the information flow and using algebraic topology it is possible to find cliques of neurons that respond to specific stimuli and determine how these are organized into cavities. The problem is that algebraic topology mainly deals with data from {\it undirected graphs}. I will now describe the structures needed for this project. The idea of simplicial complexes and persistence related to directed graphs could easily be applied to a variety of problems. In addition, there is an extension of RIPSER called FLAGSER that rapidly performs the required computations. See \cite{Gov1, LGSL, RNSTPCDLHM, Smi1} for more details.

\begin{definition}
An {\it abstract directed simplicial complex}\index{directed simplicial complex} is a collection $S$ of finite ordered sets with the property that if $\sigma$ is in $S$, then so is every subset $\tau$ of $\sigma$ with the ordering inherited from $\sigma$. 
\end{definition}

Directed simplicial complexes are similar to the simplicial complexes that we saw before but the difference is that order of the vertices is now important. So for example, $[v_0, v_1, v_2]$ is not the same as $[v_0, v_2, v_1]$. For this discussion we will call the $i$-th face $\sigma^i$ of $\sigma=[v_0, \cdots, v_n]$ the face determined by removing the vertex $v_{n-i}$

The key structure used in the Blue Brain Project is the {\it directed flag complex}.

\begin{definition}
The {\it directed flag complex}\index{directed flag complex} associated to the the directed graph $G$ is the abstract directed simplicial complex $S(G)$ whose vertices are the vertices of $G$ and whose directed $n$-simplices are $n+1$-tuples $[v_0, \cdots. v_n]$ such that for each pair $i, j$ with $0\leq i<j\leq n$, there is a directed edge from $v_i$ to $v_j$ in the graph $G$. The vertex $v_0$ is called the {\it source} of the simplex, and there is a directed edge from $v_0$ to every other vertex in the simplex. The vertex $v_n$ is called the {\it sink} of the simplex and there is a directed edge from every vertex to $v_n$. 
\end{definition}

Note that the analogue of the directed flag complex for undirected graphs is the complete subgraph or clique complex defined in Section 5.7.

Now we give a definition of directionality in directed graphs.

\begin{definition}
Let $G$ be a directed graph. For each vertex $v$ in $G$, define the {\it signed degree} of $v$ as $$sd(v)=Indegree(v)-Outdegree(v).$$ (Note that for any finite graph, $\sum_{v\in V(G)}sd(v)=0.$) The {\it directionality} of $G$ denoted $Dr(G)$ is defined to be $$Dr(G)=\sum_{v\in V(G)}sd(v)^2.$$
\end{definition}

We then have that a directed $n$-simplex $\sigma$ is a fully connected directed graph on $n+1$ vertices with a unique source and sink. It can be shown if $G$ is any other directed graph on $n+1$ vertices, then $Dr(G)\leq Dr(\sigma).$

To compute homology, we need to define chains and boundaries. We will be working with coefficients in $Z_2$. Then the $n$-chains $C_n$ are the formal sums of the $n$-dimensional directed simplices. We define the boundary $\partial_n: C_n\rightarrow C_{n-1}$ as $$\partial_n(\sigma)=\sigma^0+\sigma^1+ \cdots+\sigma^n,$$ where $\sigma^i$ is the $i$th face of $\sigma$ as defined above. The we define homology in the usual way. 

We can also define persistent homology by taking the sublevel sets of a height function on a directed graph as a filtration. FLAGSER is able to compute betti numbers and persistence diagrams using a variant of RIPSER. One feature is that it can skip matrix columns in the computations resulting in a huge speedup with only small drops in accuracy.

As an example, the Blue Brain project modeled the neocortical column (a cross section of the brain cortex) of a 14 day old rat. The graph has about 30,000 vertices and around 8 million edges. FLAGSER could build a directed flag complex in about 30 seconds on a laptop and compute a fairly accurate estimate of the homology on an HPC node in about 12 hours. See \cite{LGSL} for more details and comparisons.

\section{Computer Intrusion Detection}

How many times has something like this happened to you. You are working on your computer unaware that a fierce lion is watching you closely. After watching you type your password, the lion waits for you to leave and then pounces on your computer, deleting several of your files and writing nasty emails to your boss. (See Figure 7.3.1). How can you convince the system administrator (and your boss) that it wasn't really you.

\begin{figure}[ht]
\begin{center}
  \scalebox{0.8}{\includegraphics{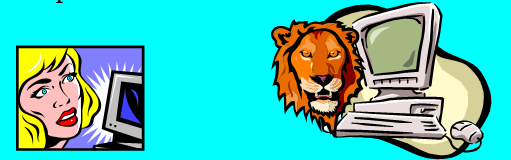}}
\caption{
\rm
Authorized user and lion impostor \cite{PG}. 
}
\end{center}
\end{figure}

Back in 2004, there were already several methods for computer intrusion detection and there are probably even better ones now. Besides, a lion would probably be hitting several keys at once and its spelling and grammar would be terrible. Still, my collaborator Tom Goldring and I thought it would be interesting to try to detect unauthorized users using features from graphs along with the newly developed classification techniques of {\it random forests} \cite{Brei}. Our work was cut short by organizational changes, but I was able to present what we had at the Second Conference on Algebraic Topological Methods in Computer Science in July 2004. There were no proceedings from the conference, but my slides can still be found online  on the conference website \cite{PG}.

The data consisted of 30 sessions each of 10 users on Windows NT machines. The data was represented by 2 directed graphs. {\it Process graphs} have vertices representing the processes called by the user and a directed edge from process 1 to process 2 if process 2 is spawned by process 1. Process graphs are disconnected trees since processes called by the operating system are omitted. {\it Window graphs} have vertices representing windows clicked on by user with a directed edge from a window to the next one clicked on. Window graphs are connected and can have cycles. Both types of graphs have edge weights representing time (in seconds) between the start of processes (process graphs) or the time between first clicks on consecutive windows (window graphs). Figures 7.3.2 and 7.3.3 show examples of a process graph and a window graph respectively. 

\begin{figure}[ht]
\begin{center}
  \scalebox{0.8}{\includegraphics{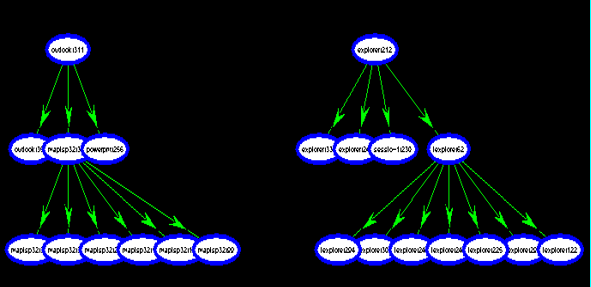}}
\caption{
\rm
Example of a process graph \cite{PG}. 
}
\end{center}
\end{figure}

\begin{figure}[ht]
\begin{center}
  \scalebox{0.8}{\includegraphics{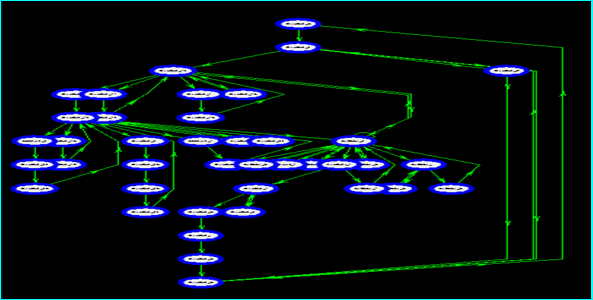}}
\caption{
\rm
Example of a window graph \cite{PG}. 
}
\end{center}
\end{figure}

The data involved 30 sessions each for 10 users. Each session had 56-long feature vectors related to the size and shape of windows and process graphs as well as edge weights representing times between clicks on windows or the start of processes. We then determined the most important features using a random forest. It turned out that window graphs were more important than process graphs. The most important window graph features were sources (representing whether the user returned to the first window used), sinks (representing whether the last window used was one that was used previously), and timing factors such as the total edge weight and the total weighted degree. 

Running the data through a 10-way random forest produced the confusion matrix in Figure 7.3.4. This shows that some users such has users 4, 9, and 10 were pretty easy to distinguish while user 5 was much more difficult. Still, this is a good start for the difficult problem of distinguishing 10 classes. 

\begin{figure}[ht]
\begin{center}
  \scalebox{0.8}{\includegraphics{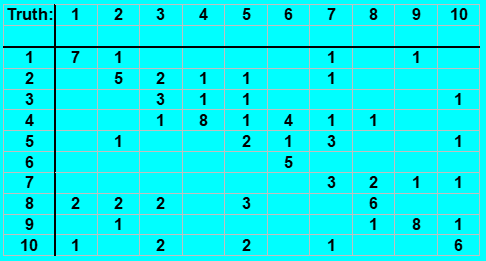}}
\caption{
\rm
Confusion matrix of 10-way random forest for intrusion detection \cite{PG}. 
}
\end{center}
\end{figure}

This was as far as we got, but we had some ideas for possible improvements. The first one was to create a directed graph called an {\it interval graph}. The vertices represent windows and there is a directed edge between 2 windows up at the same time that is directed from window the opened earlier to the one opened later. Interval graphs have the advantage that they distinguish users who open many windows at once from neater users. 

Using algebraic topology, we could add features related to the graph complexes described in Section 5.7. We could also now create the directed flag complex described in Section 7.2. Once the complex is created, we could send a classifier features such as the dimension of the complex, the number of simplices in each dimension, the Euler characteristic, or the Betti numbers of the homology groups.

Although this project was never finished, it may provide some ideas for other types of problems.

\section{Market Basket Analysis}

Back in Section 4.1.1, I mentioned the idea of market basket analysis or association analysis. (For a very understandable and complete description of this topic, see Chapters 6 and 7 of \cite{TSK}.) Suppose a supermarket has a set of items that it sells. Let those be the vertices of a hypergraph. When a customer visits the store, the list of items they buy is called a {\it transaction}. We can think of the set of transactions belonging to a particular customer as a hypergraph with each transaction corresponding to a hyperedge. One interesting problem is to produce a set of {\it association rules}. An example would be that if a customer buys peanut butter and bread, they will also buy jelly. To formulate the problem, we will need some terminology. Let $I=\{i_1, \cdots, i_n\} $ be the set of items the store sells and $T=\{t_1, \cdots, t_k\}$ be the set of transactions of all customers combined. A subset of $I$ is called an {\it itemset}. A transaction $t_i$ contains an itemset $X$ if $X$ is a subset of $t_i$. The number of transactions including the itemset $X$ is called the {\it support count} of $X$ and denoted $\sigma(X)$. 

An {\it association rule} is an expression of the form $X\rightarrow Y$ where $X$ and $Y$ are disjoint itemsets. For example, $\{$peanut butter, bread$\}\rightarrow\{$jelly$\}$ is an association rule. We can always build a rule with any two disjoint sets of items and no rule will always hold. So how do we decide which ones are interesting? We would like to select rules with high values of {\it support} and {\it confidence}. The support of a rule $X\rightarrow Y$ is the fraction of transactions which contain both $X$ and $Y$.  The confidence is the fraction of transactions containing $X$ that also contain $Y$. We would like to find all rules whose support and confidence each exceed a predetermined threshold. But the number of possible rules for a store with $d$ items is $R=3^d-2^{d+1}-1$. (Proof is an exercise in \cite{TSK}.) So we need to do something more efficient. The Apriori Alogrithm \cite{AIS} is designed to find these rules efficiently. Its strategy is to generate {\it frequent itemsets} whose support is above a given threshold and then use that list to find rules of high confidence.

Now suppose we want to use the data to solve a different problem. We would like to classify or cluster the customers based on their transactions. Are they living alone or are they buying for a large family? Can we tell the difference just from what they buy rather than how much?

As the set of each customer's transactions form a hypergraph. As with graphs, there are a number of distance measures to choose from. One example is a generalization of graph edit distance. For hypergraphs $C_1$ and $C_2$, let $\Gamma_i^1$ and $\Gamma_i^2$ be the sets of hyperedges of size $i$ for $C_1$ and $C_2$ receptively. Suppose $m$ is the maximum size of a hyperedge, and $V(C_i)$ is the set of vertices of $C_i$. Then the {\it hypergraph edit distance} is $$d(C_1, C_2)=|V(C_1)|+|V(C_2)|-2|V(C_1)\cap V(C_2)|+\sum_{i=2}^m (|\Gamma_i^1|+|\Gamma_i^2|-2|\Gamma_i^1\cap \Gamma_i^2|).$$ When $m=2$, this reduces to the usual graph edit distance.

Now recall that any hypergraph can be made into a simplicial complex by taking the hyperedges as the maximal simplices and adding all of their subsets to form a simplicial complex. Then we can classify customers by Euler characteristic and betti numbers for each dimension. These techniques are all untried as far as I know and it would be interesting to see what they would produce. 

In a talk at the Workshop on Topological Data Analysis: Theory and Applications sponsored by the Tutte Institute and the University of Western Ontario (May 1-3, 2021), Emilie Purvine of Pacific Northwest National Laboratory argued that the above method of producing a simplicial complex from a hypergraph is inadequate. This is due to the fact that we are adding simplices that aren't really there. For some interesting alternatives, you can watch her talk online \cite{Purv}.

\chapter{Cohomology}

Cohomology is sort of the reverse of homology. Instead of a boundary map that goes down in dimension, we have a coboundary that goes up in dimension. There is a cohomology exact sequence that is like the homology exact sequence only backwards. The Universal Coefficient Theorem (see Section 8.4) shows that if we know  the homology groups of a complex, we can compute the cohomology groups. In fact, two spaces with isomorphic homology groups in every dimension have isomorphic cohomology groups in every dimension. So why do we care?

The answer is that cohomology actually forms a {\it ring}. We have a product called the {\it cup product} (see Section 8.3) which takes a cohomology class of dimension $p$ and a class of dimension $q$ and produces a class of dimension $p+q$. But couldn't we do the same thing in homology? Surprisingly, the answer is no. The argument is subtle and we will see why this is true in Section 8.5. We will also see two spaces with the same homology groups but with a different cohomology ring structure, meaning the spaces are not homotopy equivalent. 

So what does this have to do with data science? Remember that at its heart, algebraic topology is a classifier. The more we know about a pair of objects, the better we can tell them apart. To paraphrase Norman Steenrod \cite{Ste1}, cohomology encodes more of the geometry into the algebra. Steenrod squares, which we will meet in Chapter 11, encode even more of the geometry as they enhance cohomology rings by giving them the richer structure of {\it algebras}. I would conjecture that the richer strucuture will provide even more information, especially if we can adapt the definition of persistence. It turns out that there is already some work in this direction, but there are a lot of open problems. 

One definite practical advantage of cohomology is that the computations are much faster than for persistent homology. Once we do these computations we can find the persistent homology as well. This fact is the main reason for the big speedup of RIPSER over its predecessors. 

In the next section, I will describe four important mathematical tools from the area of homological algebra. They will be necessary to define cohomology and state its most important properties. I will then define cohomology in Section 8.2 and describe the analogs to homology including a new set of Eilenberg-Steenrod axioms. Section 8.3 deals with the cup product which gives cohomology its ring structure. Section 8.4 discusses the two {\it Universal Coefficient Theorems} which allow us to compute cohomology of a space if we know its homology and the homology with coefficients in any abelian group if we know its homology over the integers. Section 8.5 deals with the homology and cohomology of product spaces using the K{\"u}nneth Theorems. In Section 8.6, we will look more closely at the cohomology ring structure of a product space. In Section 8.7, we turn to persistent cohomology, and we conclude the chapter in Section 8.8 by describing the role of persistent cohomology in the RIPSER software.

Unless otherwise stated, the material in Sections 8.1-8.5 comes from Munkres \cite{Mun1}. The other algebraic topology textbooks I mentioned \cite{ES, Hat, Spa} have their own explanations of the same topics. 

\section{Introduction to Homological Algebra}

Homological algebra is an important computational tool in algebraic topology, but it is an entire mathematical field of its own. We will use it to define cohomology and state the Universal Coefficiant and K{\"u}nneth Theorems as well as to perform the related computations. As stated above, the material comes from Munkres \cite{Mun1}, and we will mostly restrict ourselves to abelian groups, but the operations I will define have generalizations to modules. (Recall that an abelian group can be thought of as a module over the ring of integers.) If you want to see it in its more general form, the classic book is the book by Mac Lane with the modest name, "Homology" \cite{MacL2}. Other popular books are the ones by Rotman \cite{Rot} and Hilton and Stammbach \cite{HS}.

\subsection{Hom}

Given two abelian groups, $A$ and $G$, we get a third abelian group $Hom(A, G)$\index{Hom} consisting of the homomorphisms from $A$ to $G$. If $\phi, \psi\in Hom(A, G)$ we define addition by $(\phi+\psi)(a)=\phi(a)+\psi(a)$ for $a\in A$. The reader should check that $(\phi+\psi)$ is a homomorphism and so an element of $Hom(A, G)$. The identity of $Hom(A, G)$ is the function $f$ such that $f(A)=0$ and if $\phi\in Hom(A, G)$ then the inverse $-\phi$ is defined by $(-\phi)(a)=\phi(-a)$. 

\begin{example}
$Hom(Z, G)\cong G$ by the isomorphism $\Phi$ that takes $\phi: Z\rightarrow G$ to $\phi(1)\in G$. To see that this is an isomorphism, suppose $\phi, \psi\in Hom(A, G)$. Then $\Phi(\phi+\psi)=(\phi+\psi)(1)=\phi(1)+\psi(1)=\Phi(\phi)+\Phi(\psi).$ $\Phi$ is surjective since given $g\in G$, we can choose a homomorphism $\phi$ such that $\phi(1)=g$. $\Phi$ is injective since if $\phi: Z\rightarrow G$ is a homomorphism and $\phi(1)=0$, then for any positive integer $n$, $\phi(n)=\phi(1)+\phi(1)+\cdots+\phi(1)=(0+\cdots+0)=0.$ For $n$ negative, $\phi(n)=-\phi(-n)=0$. So the kernel of $\Phi$ is the function which is identically 0. Thus, $\Phi$ is an isomorphism.
\end{example}

The next definition will be key when defining cohomology.

\begin{definition}
A homomorphism $f: A\rightarrow B$ gives rise to a {\it dual homomorphism}\index{dual homomorphism} $$\begin{tikzpicture}
  \matrix (m) [matrix of math nodes,row sep=3em,column sep=4em,minimum width=2em]
  {
Hom(B, G) &  Hom(A, G)\\};
  \path[-stealth]

(m-1-1) edge node [above] {$\tilde{f}$} (m-1-2)

;

\end{tikzpicture}$$ going in the reverse direction. The map $\tilde{f}$ assigns to the homomorphism $\phi: B\rightarrow G$ the composite $$\begin{tikzpicture}
  \matrix (m) [matrix of math nodes,row sep=3em,column sep=4em,minimum width=2em]
  {
A & B & G.\\};
  \path[-stealth]

(m-1-1) edge node [above] {$f$} (m-1-2)
(m-1-2) edge node [above] {$\phi$} (m-1-3)
;

\end{tikzpicture}$$ In other words, $\tilde{f}(\phi)=\phi\circ f$.
\end{definition}

Note that the assignment $A\rightarrow Hom(A, G)$ and $f\rightarrow \tilde{f}$ is a contravariant functor from the category of abelian groups and maps to itself.

\begin{theorem}
If $f$ is a homomorphism and $\tilde{f}$ is its dual homomorphism then:\begin{enumerate}
\item If $f$ is an isomorphism, so is $\tilde{f}$.
\item If $f$ is the zero homomorphism, so is $\tilde{f}$.
\item If $f$ is surjective, then $\tilde{f}$ is injective. 
\end{enumerate}
\end{theorem}

The following result generalizes the last statement. 

\begin{theorem}
If the sequence  $$\begin{tikzpicture}
  \matrix (m) [matrix of math nodes,row sep=3em,column sep=4em,minimum width=2em]
  {
A & B & C & 0\\};
  \path[-stealth]

(m-1-1) edge node [above] {$f$} (m-1-2)
(m-1-2) edge node [above] {$g$} (m-1-3)
(m-1-3) edge  (m-1-4)
;

\end{tikzpicture}$$  is exact, then so is the dual sequence  $$\begin{tikzpicture}
  \matrix (m) [matrix of math nodes,row sep=3em,column sep=4em,minimum width=2em]
  {
Hom(A, G) &  Hom(B, G)  & Hom(C, G) & 0.\\};
  \path[-stealth]

(m-1-2) edge node [above] {$\tilde{f}$} (m-1-1)
(m-1-3) edge node [above] {$\tilde{g}$} (m-1-2)
(m-1-4) edge  (m-1-3)

;

\end{tikzpicture}$$ If $f$ is injective and the first sequence splits, the $\tilde{f}$ is surjective and the second sequence splits.
\end{theorem}

\begin{example}
In general, exactness of a short exact sequence does not imply exactness of the dual sequence. Let $f: Z\rightarrow Z$ be multiplication by 2. Then the sequence $$\begin{tikzpicture}
  \matrix (m) [matrix of math nodes,row sep=3em,column sep=4em,minimum width=2em]
  {
0 & Z & Z & Z_2 & 0\\};
  \path[-stealth]

(m-1-1) edge  (m-1-2)
(m-1-2) edge node [above] {$f$} (m-1-3)
(m-1-3) edge  (m-1-4)
(m-1-4) edge  (m-1-5)
;

\end{tikzpicture}$$  is exact. But $\tilde{f}$ is not surjective. If $\phi\in Hom(Z, Z)$, then $\tilde{f}(\phi)=\phi\circ f$. Now any homomorphism of $Z$ into itself is multiplication by $\alpha$ for some $\alpha\in Z$. So any element of $Hom(Z, Z)$ takes even integers to even integers. So the image of $\tilde{f}$ is not all of $Hom(Z, Z)$ and $\tilde{f}$ is not surjective.
\end{example}

So converting a short exact sequence to Hom reverses the arrows but still needs something else on the left hand side. That will be accomplished using the Ext functor that we will see in Section 8.1.3. We also remark that applying Hom with $G$ in the first coordinate instead of the second produces a covariant functor so the arrows are not reversed.

We conclude this section with a description of the Hom of two finitely generated abelian groups. I will state the theorem in a little less generality than Munkres as I will assume a finite indexing set so we only need to deal with direct sums rather than direct products.

\begin{theorem}
\begin{enumerate}
\item There are isomorphisms $$Hom(\bigoplus_{i=1}^n A_i,G)\cong\bigoplus_{i=1}^nHom(A_i, G)$$ and  $$Hom(A,\bigoplus_{i=1}^m G_i)\cong\bigoplus_{i=1}^m(Hom(A, G_i).$$ 
\item There is an isomorphism of $Hom(Z, G)$ with $G$. If $f: Z\rightarrow Z$ equals multiplication by $m$, then so does $\tilde{f}$. 
\item $Hom(Z_m, G)\cong \ker(G\xrightarrow{m} G)$.
\end{enumerate} 
\end{theorem}

{\bf Proof:} I will give Munkres' proofs of the second half of statement (2) and statement (3). (We have already discussed the first statement of (2) in Example 8.1.1 and (1) is just standard group theory.)

Let $f: Z\rightarrow Z$ be multiplication by $m$. Then for $\phi\in Hom(Z, G)$, and $z\in Z$, $\tilde{f}(\phi)(z)=\phi(f(z))=\phi(mz)=m\phi(z)$. So $\tilde{f}(\phi)=m\phi$. Then $\tilde{f}$ is multiplication by $m$ in $Hom(Z, G)$ and so also in $G$ by way of the isomorphism of $Hom(Z, G)$ with $G$.

To see (3), if $$\begin{tikzpicture}
  \matrix (m) [matrix of math nodes,row sep=3em,column sep=4em,minimum width=2em]
  {
0 & Z & Z & Z_m & 0\\};
  \path[-stealth]

(m-1-1) edge  (m-1-2)
(m-1-2) edge node [above] {$m$} (m-1-3)
(m-1-3) edge  (m-1-4)
(m-1-4) edge  (m-1-5)
;

\end{tikzpicture}$$  is exact, then 
$$\begin{tikzpicture}
  \matrix (m) [matrix of math nodes,row sep=3em,column sep=4em,minimum width=2em]
  {
Hom(Z, G) & Hom(Z, G) & Hom(Z_m, G) & 0\\};
  \path[-stealth]

(m-1-2) edge  node [above] {$m$} (m-1-1)
(m-1-3) edge  (m-1-2)
(m-1-4) edge  (m-1-3)

;

\end{tikzpicture}$$  is exact. Then $\ker(G\xrightarrow{m} G)\cong Hom(Z_m, G)$. ${\blacksquare}$

What if we want to find $Hom(G, H)$ where $G$ and $H$ are finite cyclic groups?

\begin{theorem}
There is an exact sequence $$\begin{tikzpicture}
  \matrix (m) [matrix of math nodes,row sep=3em,column sep=4em,minimum width=2em]
  {
0 & Z_d & Z_n & Z_n & Z_d & 0,\\};
  \path[-stealth]

(m-1-1) edge  (m-1-2)
(m-1-2) edge  (m-1-3)
(m-1-3) edge  node [above] {$m$} (m-1-4)
(m-1-4) edge  (m-1-5)
(m-1-5) edge  (m-1-6)

;

\end{tikzpicture}$$ where $d=gcd(m, n)$, i.e. the greatest common divisor of $m$ and $n$. From this it can be shown that $Hom(Z_m, Z_n)\cong \ker(Z_n\xrightarrow{m} Z_n)\cong Z_d$. 
\end{theorem}

Now we can determine $Hom(A, B)$ when $A$ and $B$ are any finitely generated abelian groups. Reviewing what we know, we have:
$$Hom(Z, G)\cong G$$
$$Hom(Z_m, G)=\ker(G\xrightarrow{m} G)$$
So $$Hom(Z_m, Z)=0$$ and $$Hom(Z_m, Z_n)\cong Z_d,$$ where $d=gcd(m, n).$

We will want to remember these formulas for later. 

\begin{example}
Let $A=Z\oplus Z\oplus Z_{15}\oplus Z_3$, and $B=Z\oplus Z_6\oplus Z_2$. Find $Hom(A, B)$.\begin{align*}
Hom(A, B)&\cong Hom(Z\oplus Z\oplus Z_{15}\oplus Z_3, Z\oplus Z_6\oplus Z_2)\\
&\cong Hom(Z. Z)\oplus Hom(Z, Z_6)\oplus Hom(Z, Z_2) \oplus\\
& Hom(Z. Z)\oplus Hom(Z, Z_6)\oplus Hom(Z, Z_2) \oplus\\
& Hom(Z_{15}, Z)\oplus Hom(Z_{15}, Z_6)\oplus Hom(Z_{15}, Z_2) \oplus\\
& Hom(Z_3. Z)\oplus Hom(Z_3, Z_6)\oplus Hom(Z_3, Z_2)\\
&\cong Z\oplus Z_6\oplus Z_2 \oplus\\
& Z\oplus Z_6\oplus Z_2 \oplus\\
& 0\oplus Z_3\oplus 0 \oplus\\
& 0\oplus Z_3 \oplus 0\\
&\cong Z^2\oplus (Z_6)^2\oplus (Z_3)^2 \oplus (Z_2)^2.\\
\end{align*}

\end{example}

\subsection{Tensor Product}

Just when you thought you could relax, the next topic is guaranteed to make you a lot tenser. (Or maybe tensor?) I will give a description of the important operation, {\it tensor products}.

Before starting, I need to give you a mathematical word reuse warning. Tensor products are not the same thing as the tensors used in physics. Those tensors are sort of like multidimensional matrices. You can think of matrices and vectors as special cases. They do have interesting uses in data science and are an alternate method of analyzing multivariate time series. For anything you would ever want to know about tensors and data science, start with the paper of Tammy Kolda and Brett Bader of Sandia National Laboratory \cite{KB}.

With that out of the way, tensor products are a binary operation on abelian groups, or more generally, on modules. We would like a function $f: A\times B\rightarrow C$, for abelian groups $A, B, C$, that is linear when restricted to either variable. In other words, we want $f$ to be {\it bilinear}. This is what tensor products accomplish.

\begin{definition}
Let $A$ and $B$ be abelian groups. Let $F(A, B)$ be the free abelian group generated be the elements of $A\times B$. Let $R(A, B)$ be the subgroup generated by elements of the form $$(a_1+a_2, b_1)-(a_1, b_1)-(a_2, b_1),$$ $$(a_1, b_1+b_2)-(a_1, b_1)-(a_1, b_2),$$ for $a_1, a_2\in A$ and $b_1, b_2\in B$. Define the {\it tensor product}\index{tensor product} denoted $A\otimes B$ by \index{$A\otimes B$}$$A\otimes B=F(A, B)/R(A, B).$$ The coset of the pair $(a, b)$ is denoted $a\otimes b$.
\end{definition}

Our goal now is to list some properties of tensor product and show how to compute the tensor products of various abelian groups.Then we will compute the tensor product of the two groups from example 8.1.3.

Given a function from $A\times B\rightarrow C$, where $C$ is an abelian group, we have a unique homomorphism of $F(A, B)$ into $C$ since the elements of  $A\times B$ form a basis of $F(A, B)$. This function is bilinear if and only if it maps the subgroup $R(A, B)$ to zero. So every homomorphism of $A\otimes B$ into $C$ determines a bilinear function of $A\times B$ into $C$ and vice versa.

Now any element of $F(A, B)$ is a finite linear combination of pairs $(a, b)$ so any element of $A\otimes B$ is a finite linear combination of elements of the form $a\otimes b$. So $a\otimes b$ is not a typical element of $A\otimes B$, but it is a typical generator.

Now the definition implies that $$(a_1+a_2, b_1)=(a_1, b_1)+(a_2, b_1),$$ and $$(a_1, b_1+b_2)=(a_1, b_1)+(a_1, b_2).$$ So $$a\otimes b=(0+a)\otimes b=0\otimes b+a\otimes b.$$ This gives $0\otimes b=0.$ Similar arguments show $$a\otimes 0=0,$$ $$(-a)\otimes b=-(a\otimes b)=a\otimes (-b),$$ and $$(na)\otimes b=n(a\otimes b)=a\otimes (nb),$$ where $n\in Z.$

\begin{definition}
Let $f: A\rightarrow A'$ and $g: B\rightarrow B'$ be homomorphisms. There is a unique homomorphism $$f\otimes g: A\otimes B\rightarrow A'\otimes B'$$ such that $(f\otimes g)(a\otimes b)=f(a)\otimes g(b)$ for all $a\in A, b\in B$. $f\otimes g$ is called the {\it tensor product} of $f$ and $g$.
\end{definition}

\begin{theorem}
There is an isomorphism $Z\otimes G\cong G$ mapping $n\otimes g$ to $ng$. It is {\it natural} with respect to homomorphisms of $G$. 
\end{theorem}

{\bf Proof:} The function mapping $Z\times G$ to $G$ sending $(n, g)$ to $ng$ is bilinear so it induces a homomorphism $\phi: Z\otimes G\rightarrow G$ sending $n\otimes g$ to $ng$.

Let $\psi: G\rightarrow Z\otimes G$ with $\psi(g)=1\otimes g$. For $g\in G$, $\phi\psi(g)=\phi(1\otimes g)=g$. On a generator $n\otimes g$ of $Z\otimes G$, $\psi(\phi(n\otimes g))=\psi(ng)=1\otimes ng=n\otimes g$. So $\psi$ and $\phi$ are inverses of each other and $Z\otimes G\cong G$.

Naturality comes from the commutativity of the diagram  $$\begin{tikzpicture}
  \matrix (m) [matrix of math nodes,row sep=3em,column sep=4em,minimum width=2em]
  {
Z\otimes G & G&\\
Z\otimes H & H.&\hspace{.5in}\blacksquare\\
};
  \path[-stealth]

(m-1-1) edge node [above] {$\cong$} (m-1-2)
(m-2-1) edge node [above] {$\cong$} (m-2-2)
(m-1-1) edge node [left] {$i_z\otimes f$} (m-2-1)
(m-1-2) edge node [right] {$f$} (m-2-2)
;

\end{tikzpicture}$$

Important note: In this chapter the term {\it natural isomorphism} or {\it natural exact sequence} will mean that maps will satisfy a commutative diagram of the appropriate shape.

\begin{example}
Let $A'$ be a subgroup of $A$ and $B'$ be a subgroup of $B$. It is not necessarily true that $A'\otimes B'$ is a subgroup of $A\otimes B$. The reason is that the tensor product of the inclusion mappings is not necessarily injective. As an example, let $A=Q$, the group of rational numbers under addition. Let $A'=Z$, and $B=B'=Z_2.$ Then $Z\otimes Z_2\cong Z_2$, by our previous theorem, but $Q\otimes Z_2=0$ since in $Q\otimes Z_2$, $a\otimes b=(a/2)\otimes 2b=(a/2)\otimes 0=0.$ This gives $Q\otimes Z_2=0.$ Obviously, $Z_2$ is not a subgroup of 0. This also shows that the tensor product of the inclusion maps, $i: A'\rightarrow A$ and $j: B'\rightarrow B$ are not necessarily injective.
\end{example}

Although the tensor product of injective maps may not be injective, the tensor product of surjective maps is always surjective.

\begin{theorem}
Let $\phi: B\rightarrow C$ and $\phi': B'\rightarrow C'$ are surjective. Then $\phi\otimes\phi': B\otimes B'\rightarrow C\otimes C'$ is surjective, and its kernel is the subgroup of $B\otimes B'$ generated by all elements of the form $b\otimes b'$ where $b\in\ker B$, and $b'\in\ker B'$.
\end{theorem}

{\bf Idea of Proof:} Let $G$ be the subgroup of $B\otimes B'$ generated by elements of the form $b\otimes b'$ described in the statement of the theorem. Then $\phi\otimes\phi'$ maps $G$ to zero, so it induces a homomorphism $\Phi: (B\otimes B')/G\rightarrow C\otimes C'$. We are done if $\Phi$ is an isomorphism, so we build an inverse for it. Let $\psi: C\times C'\rightarrow(B\otimes B')/G$ by $\psi(c, c')=b\otimes b'+G$ where $b, b'$ are chosen so that $\phi(b)=c$ and $\phi(b')=c'$. Check that $\psi$ is well defined and bilinear so that it defines a homomorphism $\Psi: C\otimes C' \rightarrow (B\otimes B')/G.$ Finally show that $\Phi$ and $\Psi$ are inverses.  $\blacksquare$

Now we look at what tensor product does to exact sequences.

\begin{theorem}
 If $$\begin{tikzpicture}
  \matrix (m) [matrix of math nodes,row sep=3em,column sep=4em,minimum width=2em]
  {
A & B & C & 0\\};
  \path[-stealth]

(m-1-1) edge node [above] {$\phi$}  (m-1-2)
(m-1-2) edge node [above] {$\psi$} (m-1-3)
(m-1-3) edge  (m-1-4)

;

\end{tikzpicture}$$  is exact, then 
$$\begin{tikzpicture}
  \matrix (m) [matrix of math nodes,row sep=3em,column sep=4em,minimum width=2em]
  {
A\otimes G & B\otimes G & C\otimes G & 0\\};
  \path[-stealth]

(m-1-1) edge node [above] {$\phi\otimes i_G$}  (m-1-2)
(m-1-2) edge node [above] {$\psi\otimes i_G$} (m-1-3)
(m-1-3) edge  (m-1-4)

;

\end{tikzpicture}$$  is exact. If $\phi$ is injective and the first sequence splits, then $\phi\otimes i_G$ is injective and the second sequence splits.
\end{theorem}

{\bf Proof:} We know from the previous theorem that $\psi\otimes i_G$ is surjective and that the kernel is the subgroup $D$ of $B\otimes G$ generated by elements of the form $b\otimes g$ with $b\in\ker\psi$. The image of $\phi\otimes i_G$ is the subgroup $E$ of $B\otimes G$ generated by elements of the form $\phi(a)\otimes g$. Since $\ker\psi=\mbox{image }\phi$, $D=E$. 

If $\phi$ is injective and the first sequence splits, let $p: B\rightarrow A$ be a homomorphism such that $p\phi=i_A$. (Remember this is the definition of the sequence splitting.) Then $$(p\otimes i_G)\circ(\phi\otimes i_G)=i_A\otimes i_G=i_{A\otimes G}.$$ So $\phi\otimes i_G$ is injective and the second sequence splits via $p\otimes i_G.$ $\blacksquare$

Just as Ext filled out the missing piece of the exact sequence when Hom was applied, Tor does the same thing with tensor product (under additional conditions). I will describe them both in detail in the next two subsections. 

\begin{theorem}
There is a natural isomorphism $$Z_m\otimes G\cong G/mG$$.
\end{theorem}

{\bf Proof:}

Consider the exact sequence $$0\rightarrow Z\xrightarrow{m} Z \rightarrow Z_m\rightarrow 0.$$ (Remember that $m$ above the arrow means multiplication by $m$.) Now tensor it with $G$ to get the exact sequence $$Z\otimes G\xrightarrow{m\otimes i_G} Z\otimes G\rightarrow Z_m\otimes G\rightarrow 0.$$ Since $Z\otimes G\cong G$, this becomes $$G\xrightarrow{m} G\rightarrow Z_m\otimes G\rightarrow 0.$$ This proves the result. $\blacksquare$

The next two results will allow us to compute the tensor product of any two finitely generated abelian groups. The proofs are lengthier, so the reader is referred to Munkres \cite{Mun1}.

\begin{theorem}
We have the following natural isomorphisms:\begin{enumerate}
\item $A\otimes B\cong B\otimes A.$ (Recall that this is not true for Hom. $Hom(A, B)$ is generally not isomorphic to $Hom(B, A)$.)
\item $(\bigoplus_{i=1}^n A_i)\otimes B\cong \bigoplus_{i=1}^n (A_i\otimes B).$ 

$A\otimes(\bigoplus_{i=1}^n B_i) \cong \bigoplus_{i=1}^n (A\otimes B_i).$
\item$A\otimes (B\otimes C)\cong (A\otimes B)\otimes C.$
\end{enumerate}
\end{theorem}

Recall that an additive abelian group is {\it torsion free} if there are no elements of finite order. Note that a torsion free finitely generated abelian group is free.

\begin{theorem}
If $$0\rightarrow A\rightarrow B\rightarrow C\rightarrow 0$$ is exact and $G$ is torsion free, then $$0\rightarrow A\otimes G\rightarrow B\otimes G\rightarrow C\otimes G\rightarrow 0$$ is exact.
\end{theorem}

The proof is easy in the case $G$ is free which will always be the case if $G$ is finitely generated and torsion free, since $D\otimes Z\cong D$ for any abelian group $D$, and direct sums of exact sequences are exact. So exactness of the second sequence follows from exactness of the first sequence. See \cite{Mun1} for the proof in the more general case.

An easy consequence is the following:

\begin{theorem}
If $A$ is free abelian with basis $\{a_i\}$,  and $B$ is free abelian with basis $\{b_j\}$ then $A\otimes B$ is free abelian with basis $\{a_i\otimes b_j\}$.
\end{theorem}

Now we can take the tensor products of any two finitely generated ableian groups. To collect what we know:
$$Z\otimes G\cong G\otimes Z\cong G.$$
$$Z_m\otimes G\cong G\otimes Z_m\cong G/mG.$$ These statements imply that 
$$Z_m\otimes Z\cong Z\otimes Z_m\cong Z/mZ=Z_m,$$ and if $d=gcd(m, n)$, then $$Z_m\otimes Z_n\cong Z_d.$$

\begin{example}

Returning to our running example, 
let $A=Z\oplus Z\oplus Z_{15}\oplus Z_3$, and $B=Z\oplus Z_6\oplus Z_2$. Find $A\otimes B$.\begin{align*}
A\otimes B&\cong (Z\oplus Z\oplus Z_{15}\oplus Z_3)\otimes (Z\oplus Z_6\oplus Z_2)\\
&\cong (Z\otimes Z)\oplus (Z\otimes Z_6)\oplus (Z\otimes Z_2) \oplus\\
& (Z\otimes Z)\oplus (Z\otimes Z_6)\oplus (Z\otimes Z_2) \oplus\\
& (Z_{15}\otimes Z)\oplus (Z_{15}\otimes Z_6)\oplus (Z_{15}\otimes Z_2) \oplus\\
& (Z_3\otimes Z)\oplus (Z_3\otimes Z_6)\oplus (Z_3\otimes Z_2)\\
&\cong Z\oplus Z_6\oplus Z_2 \oplus\\
& Z\oplus Z_6\oplus Z_2 \oplus\\
& Z_{15}\oplus Z_3\oplus 0 \oplus\\
& Z_3\oplus Z_3 \oplus 0\\
&\cong Z^2\oplus Z_{15}\oplus (Z_6)^2\oplus (Z_3)^3 \oplus (Z_2)^2.\\
\end{align*}

\end{example}

Finally, tensor products provide a new definition of homology with coefficients in an arbitrary abelian group $G$. For applications, the most common situation is $G=Z_2$.

\begin{definition}
Let $G$ be an abelian group, and $\mathcal{C}=\{C_p, \partial\}$ be a chain complex. Then the $p$-th homology group of the chain complex $\mathcal{C}\otimes G=\{C_p\otimes G, \partial\otimes i_G\}$ is denoted $H_p(\mathcal{C}; G)$ and called the $p$-th {\it homology group of} $\mathcal{C}$ {\it with coefficients in} $G$.\index{homology with arbitrary coefficients}
\end{definition}

Let's look at the simplicial homology case. How does this definition compare with the one we gave in Section 4.1.3? The group $C_p(K)\otimes G$ is the direct sum of copies of $G$, one for each simplex of $K$. Each element of $C_p(K)\otimes G$ can be represented as a finite sum $\sum\sigma_i\otimes g_i$. Its boundary can be represented as $\sum(\partial\sigma_i)\otimes g_i$. In Section 4.1.3 we represented the chain $c_p$ with coefficients in $G$ by a finite formal sum $c_p=\sum g_i\sigma_i$ with boundary $\partial c_p=\sum g_i(\partial \sigma_i)$. It should be clear that the two chain complexes are isomorphic, but for the theory, $C(K)\otimes G$ is more convenient. We will use this definition for the rest of the chapter.

\subsection{Ext}

This topic is {\bf Ext}remely interesting. In fact it is an example of {\bf Ext}reme algebraic topology. 

Actually, Ext\index{Ext} is an abbreviation of the word extension. It is the missing piece in taking a short exact sequence and applying the Hom functor. Ext will be a crucial piece in converting homology to cohomology using the Universal Coefficient Theorem that I will describe in Section 8.4. 

$Ext(A, B)$ is a functor like $Hom(A, B)$ taking two abelian groups and returning a third abelian group. It is contravariant in the first variable and covariant in the second. This means that given homomorphisms $\gamma: A\rightarrow A'$ and $\delta: B'\rightarrow B$, there is a homomorphism $Ext(\gamma, \delta): Ext(A', B')\rightarrow Ext(A, B).$ 

Considering Ext as a black box for now, the following theorem describes its main properties. We need a definition first.

\begin{definition}
Let $A$ be an abelian group. A {\it free resolution}\index{free resolution} of $A$ is a short exact sequence $$0\rightarrow R\rightarrow F\rightarrow A\rightarrow 0$$ such that $F$ and $R$ are free. Any abelian group $A$ has a free resolution $$0\rightarrow R(A)\rightarrow F(A)\rightarrow A\rightarrow 0,$$ where $F(A)$ is the free abelian group generated by the elements of $A$, and $R(A)$ is the kernel of the projection of $F(A)$ into $A$. (In other words, R(A) is generated by the relations of $A$.) This resolution is called the {\it canonical free resolution}\index{canonical free resolution} of $A$.
\end{definition}

\begin{theorem}
There is a function that assigns to each free resolution  $$0\rightarrow R\xrightarrow{\phi} F\xrightarrow{\psi} A\rightarrow 0$$ of the abelian group $A$, and to each abelian group $B$, an exact sequence $$0\leftarrow Ext(A, B)\xleftarrow{\pi}Hom(R, B)\xleftarrow{\tilde{\phi}}Hom(F, B)\xleftarrow{\tilde{\psi}}Hom(A, B)\leftarrow 0.$$ This function is natural in in the sense that a homomorphism  $$\begin{tikzpicture}
  \matrix (m) [matrix of math nodes,row sep=3em,column sep=4em,minimum width=2em]
  {
0 & R & F & A & 0\\
0 & R' & F' & A' & 0\\};
  \path[-stealth]

(m-1-1) edge  (m-1-2)
(m-1-2) edge  (m-1-3)
(m-1-3) edge  (m-1-4)
(m-1-4) edge  (m-1-5)
(m-2-1) edge  (m-2-2)
(m-2-2) edge  (m-2-3)
(m-2-3) edge  (m-2-4)
(m-2-4) edge  (m-2-5)
(m-1-2) edge node [right] {$\alpha$}  (m-2-2)
(m-1-3) edge node [right] {$\beta$}  (m-2-3)
(m-1-4) edge node [right] {$\gamma$}  (m-2-4)

;

\end{tikzpicture}$$ of free resolutions and a homomorphism $\delta: B'\rightarrow B$ of abelian groups gives rise to a homomorphism of exact sequences:$$\begin{tikzpicture}
  \matrix (m) [matrix of math nodes,row sep=3em,column sep=4em,minimum width=2em]
  {
0 & Ext(A, B) & Hom(R, B) & Hom(F, B) & Hom(A, B) & 0\\
0 & Ext(A', B') & Hom(R', B') & Hom(F', B') & Hom(A', B') & 0.\\};
  \path[-stealth]

(m-1-2) edge  (m-1-1)
(m-1-3) edge  (m-1-2)
(m-1-4) edge  (m-1-3)
(m-1-5) edge  (m-1-4)
(m-1-6) edge  (m-1-5)
(m-2-2) edge  (m-2-1)
(m-2-3) edge  (m-2-2)
(m-2-4) edge  (m-2-3)
(m-2-5) edge  (m-2-4)
(m-2-6) edge  (m-2-5)
(m-2-2) edge node [left] {$Ext(\gamma, \delta)$}  (m-1-2)
(m-2-3) edge node [left] {$Hom(\alpha, \delta)$}  (m-1-3)
(m-2-4) edge node [left] {$Hom(\beta, \delta)$}  (m-1-4)
(m-2-5) edge node [left] {$Hom(\gamma, \delta)$}  (m-1-5)

;

\end{tikzpicture}$$

\end{theorem}

Now the group $Ext(A, B)$ is the cokernel of $\tilde{\phi}$ by exactness. It turns out that the map $Ext(\gamma, \delta)$ does not depend on $\alpha$ or $\beta$ at all. The proof depends on some cohomology theory that I will skip but is explained in Sections 46 and 52 of \cite{Mun1}. It's definition is derived from the commutativity of the two right squares in the diagram at the end of the previous theorem. The two right squares commute by functoriality of Hom, so the left square commutes. We say that $Ext(\gamma, \delta)$ is the homomorphism, induced by $\gamma$ and $\delta$. Munkres also shows that if $\gamma=i_A$ and $\delta=i_B$, then $Ext(i_A, i_B)$ is an isomorphism even if we choose two different free resolutions of $A$.  So let  $$0\rightarrow R(A)\xrightarrow{\phi} F(A)\rightarrow A\rightarrow 0$$ be the canonical free resolution of $A$. Then $$Ext(A, B)=cok(\tilde{\phi})=Hom(R(A), B)/\tilde{\phi}(Hom(F(A), B)).$$ We will use this as our definition.

Now if $\gamma: A\rightarrow A'$ and $\delta: B'\rightarrow B$, then if we extend $\gamma$ to a homomorphism of the canonical free resolution of $A$ to that of $A'$, then we define $Ext(\gamma, \delta): Ext(A', B')\rightarrow Ext(A, B)$ to be the homomorphism induced by $\gamma$ and $\delta$ relative to these free resolutions. So Ext is a functor of two variables which is contravariant in the first variable and covariant in the second. The name $Ext(A, B)$  is short for the {\it group of extensions} of $B$ by $A$.

To compute $Ext(A, B)$ when $A$ and $B$ are finitely generated abelian groups, use the following result,  (Assume only finitely many terms so we can use direct sums.)

\begin{theorem}
\begin{enumerate}
\item There are natural isomorphisms $$Ext(\oplus A_i, B)\cong \oplus Ext(A_i, B),$$ and  $$Ext(A,\oplus B_i)\cong \oplus Ext(A, B_i).$$ 
\item $Ext(A, B)=0$ if $A$ is free.
\item Given $B$, there is an exact sequence $$0\leftarrow Ext(Z_m, B)\leftarrow B\xleftarrow{m} B\leftarrow Hom(Z_m, B)\leftarrow 0.$$
\end{enumerate}
\end{theorem}

{\bf Proof:} Part 1 follows the similar result for Hom. Taking direct sums of the free resolutions result in isomorphisms of the three right hand terms and thus, isomorphisms of the left terms. 

For part 2, a free resolution of $A$ splits if $A$ is free. So taking Hom leaves the sequence exact and $Ext(A, B)=0.$

For part 3, Theorem 8.1.12, start with the free resolution $0\rightarrow Z\xrightarrow{m} Z\rightarrow Z_m\rightarrow 0.$ Applying Theorem 8.1.12 gives the sequence $$0\leftarrow Ext(Z_m, B)\leftarrow Hom(Z, B)\xleftarrow{m}Hom(Z, B)\leftarrow Hom(Z_m, B)\leftarrow 0.$$ Using the fact that $Hom(Z, B)\cong B$, we have the result. $\blacksquare$

The theorem implies that $Ext(Z, G)=0$ and $Ext(Z_m, G)\cong G/mG.$ As a consequence, $Ext(Z_m, Z)\cong Z_m$ and if $d=gcd(m, n)$, then $Ext(Z_m, Z_n)\cong Z_d.$

\begin{example}

Returning to our running example, 
let $A=Z\oplus Z\oplus Z_{15}\oplus Z_3$, and $B=Z\oplus Z_6\oplus Z_2$. Find $A\otimes B$.\begin{align*}
Ext(A, B)&\cong Ext(Z\oplus Z\oplus Z_{15}\oplus Z_3, Z\oplus Z_6\oplus Z_2)\\
&\cong Ext(Z. Z)\oplus Ext(Z, Z_6)\oplus Ext(Z, Z_2) \oplus\\
& Ext(Z. Z)\oplus Ext(Z, Z_6)\oplus Ext(Z, Z_2) \oplus\\
& Ext(Z_{15}. Z)\oplus Ext(Z_{15}, Z_6)\oplus Ext(Z_{15}, Z_2) \oplus\\
& Ext(Z_3. Z)\oplus Ext(Z_3, Z_6)\oplus Ext(Z_3, Z_2)\\
&\cong 0\oplus 0\oplus 0 \oplus\\
& 0\oplus 0\oplus 0 \oplus\\
& Z_{15}\oplus Z_3\oplus 0 \oplus\\
& Z_3\oplus Z_3 \oplus 0\\
&\cong Z_{15}\oplus(Z_3)^3.\\
\end{align*}

\end{example}

\subsection{Tor}

Tor, which stands for "the onion router" is free and open source software for anonymous communication. Traffic is routed through publicly available relays and encrypts the destination IP address. {\it Tor hidden services} allow customers to access a server anonymously and are often used for illegal activities. This subject, while interesting, is NOT the topic of this section. 

The other Tor is short for torsion product and is the last functor we will need. Meanwhile, what do you call a group of mathematicians who give you anonymous help on your homological algebra? A Tor hidden service. 

Tor\index{Tor} will play a similar role to Ext where we will apply tensor product rather than Hom to an exact sequence. Ext and Tor are sometimes called {\it derived functors}\index{derived functor} as Ext is derived from Hom and Tor is derived from tensor product. Tor is a functor that assigns to an ordered pair $A$, $B$ of abelian groups an abelian group $Tor(A, B)$ and to an ordered pair of homomorphisms $\gamma: A\rightarrow A'$ and $\delta: B\rightarrow B'$, a homomorphism $Tor(\gamma, \delta): Tor(A, B)\rightarrow Tor(A', B')$. It is covariant in both variables. 

Note: Munkres uses the notation $A\ast B$ rather than $Tor(A, B)$. I prefer the latter notation which appears in the more general treatment of Mac Lane \cite{MacL2}. 

As the constructions are analogous to those of Ext, I will keep this section brief and refer the reader to Section 54 of \cite{Mun1} for the proofs.

\begin{theorem}
There is a function that assigns to each free resolution  $$0\rightarrow R\xrightarrow{\phi} F\xrightarrow{\psi} A\rightarrow 0$$ of the abelian group $A$, and to each abelian group $B$, an exact sequence $$0\rightarrow Tor(A, B)\xrightarrow{\pi\otimes i_B}R\otimes B\xrightarrow{\phi\otimes i_B}F\otimes B\xrightarrow{\psi\otimes i_B}A\otimes B\rightarrow 0.$$ The function is natural in the sense that a homomorphism of a free resolution of $A$ to a free resolution of $A'$ and a homomorphism of $B$ to $B'$ induce a homomorphism of the tensor product exact sequences. 
\end{theorem}

\begin{definition}
Let $A$ be an abelian group and let $$0\rightarrow R(A)\xrightarrow{\phi} F(A)\rightarrow A\rightarrow 0$$ be the canonical free resolution of $A$. The group $\ker(\phi\otimes i_B)$ is denoted $Tor(A, B)$ and called the {\it torsion product}\index{torsion product} of $A$ and $B$. If $\gamma: A\rightarrow A'$ and $\delta: B\rightarrow B'$ are homomorphisms then extend $\gamma$ to a homomorphism of canonical free resolutions and define $Tor(\gamma, \delta): Tor(A, B)\rightarrow Tor(A', B')$ to be the homomorphism induced by $\gamma$ and $\delta$ relative to these free resolutions. 
\end{definition}

\begin{theorem}
\begin{enumerate}
\item There is a natural isomorphism $Tor(A, B)\cong Tor(B, A)$.
\item There are natural isomorphisms $$Tor(\oplus A_i, B)\cong \oplus Tor(A_i, B),$$ and  $$Tor(A,\oplus B_i)\cong \oplus Tor(A, B_i).$$ 
\item $Tor(A, B)=0$ if $A$ or $B$ is torsion free.
\item Given $B$, there is an exact sequence $$0\rightarrow Tor(Z_m, B)\rightarrow B\xrightarrow{m} B \rightarrow Z_m\otimes B\rightarrow 0.$$
\end{enumerate}
\end{theorem}

This theorem implies that $Tor(Z, G)=0$ and $Tor(Z_m, G)\cong\ker(G\xrightarrow{m} G)$. The latter gives $Tor(Z_m, Z)=0$ and for $d=gcd(m, n)$, $Tor(Z_m, Z_n)\cong Z_d$.

\begin{example}

Returning to our running example, 
let $A=Z\oplus Z\oplus Z_{15}\oplus Z_3$, and $B=Z\oplus Z_6\oplus Z_2$. Find $A\otimes B$.
\begin{align*}
Tor(A, B)&\cong Tor(Z\oplus Z\oplus Z_{15}\oplus Z_3, Z\oplus Z_6\oplus Z_2)\\
&\cong Tor(Z. Z)\oplus Tor(Z, Z_6)\oplus Tor(Z, Z_2) \oplus\\
& Tor(Z. Z)\oplus Tor(Z, Z_6)\oplus Tor(Z, Z_2) \oplus\\
& Tor(Z_{15}. Z)\oplus Tor(Z_{15}, Z_6)\oplus Tor(Z_{15}, Z_2) \oplus\\
& Tor(Z_3. Z)\oplus Tor(Z_3, Z_6)\oplus Tor(Z_3, Z_2)\\
&\cong 0\oplus 0\oplus 0 \oplus\\
& 0\oplus 0\oplus 0 \oplus\\
& 0\oplus Z_3\oplus 0 \oplus\\
& 0\oplus Z_3 \oplus 0\\
&\cong (Z_3)^2.\\
\end{align*}

\end{example}

The following table from \cite{Mun1} summarizes all four functors and allows us to compute them on any pair of finitely generated abelian groups:
\begin{align*}
&Z\otimes G\cong G &Tor(Z, G)=0\\
&Hom(Z, G)\cong G &Ext(Z, G)=0\\
&Z_m\otimes G\cong G/mG &Tor(Z_m, G)\cong\ker(G\xrightarrow{m} G)\\
&Hom(Z_m, G)\cong\ker(G\xrightarrow{m} G) &Ext(Z_m, G)\cong G/mG\\
&Z_m\otimes Z\cong Z_m &Tor(Z_m, Z)=0\\
&Hom(Z_m, Z)=0 &Ext(Z_m, Z)\cong Z_m\\
\end{align*}

If $d=gcd(m, n)$, then $$Z_m\otimes Z_n\cong Tor(Z_m, Z_n)\cong Hom(Z_m, Z_n)\cong Ext(Z_m, Z_n)\cong Z_d.$$

I will conclude this section with some comments on the more general case where $A$ and $B$ are modules over a ring $R$. Let $a\in A$, $b\in B$ and $r\in R$. An element $f\in Hom_R(A, B)$ has the additional condition that $f(ra)=rf(a)$. For tensor products, we write $A\otimes_R B$ and have the additional condition $$r(a\otimes b)=(ra)\otimes b=a\otimes(rb).$$

For Ext and Tor, things become more complicated. Let A be an $R$-module and let $F_0$ be the free module generated by the elements of $A$. Let $\phi: F_0\rightarrow A$ be the natural epimorphism. What changes here is that a submodule of a free module is not necessarily free. So we don't necessarily have an exact sequence $0\rightarrow R\rightarrow F\rightarrow A$ where $R$ and $F$ are both free. Instead, let $F_1$ be the free $R$-module generated by elements of $\ker(\phi)$. This gives an exact sequence $F_1\rightarrow F_0\rightarrow A$ where $F_1$ and $F_0$ are free. Continuing gives a {\it free resolution} of the $R$-module $A$ which is in the form $$\cdots\rightarrow F_k\xrightarrow{\phi_k} F_{k-1}\rightarrow\cdots\xrightarrow{\phi_1}F_0\rightarrow A\rightarrow 0,$$ where each $F_i$ is free. Applying the functor $Hom_R$, we get the sequence $$\cdots\xleftarrow{\tilde{\phi}_2} Hom_R(F_1, G)\xleftarrow{\tilde{\phi}_1}Hom_R(F_0, G)\xleftarrow{\tilde{\phi}}Hom_R(A, G)\leftarrow 0.$$ Exactness only holds for the 2 right hand terms. For $n\geq 1$, $Ext^n_R(A, G)=ker(\tilde{\phi}_{n+1})/im(\tilde{\phi}_n).$ This gives an entire sequence of Ext terms. If $R=Z$, $Ext^1$ is just the Ext we defined earlier, and $Ext^n(A, G)=0$ for $n>1.$ We can define $Tor_R^n(A, G)$ in an analogous way. See \cite{MacL2} for more details.

\section{Definition of Cohomology}

We are now ready to define cohomology. With what you know now. it will be pretty easy. Think of cohomology as a {\it dual} of homology. Boundary maps that decrease the dimension are replaced with coboundary maps that increase the dimension. Here are the specifics: 

\begin{definition}
Let $K$ be a simplicial complex and $G$ an abelian group. The group of {\it p-cochains}\index{cochain} of $K$ with coefficients in $G$ is the group $$C^p(K; G)=Hom(C_p(K), G).$$ Now let $f: C_p(K)\rightarrow G$ be a $p$-cochain. Then we can define a homomorphism $\delta(f): C_{p+1}(K)\rightarrow G$ by $\delta(f)=f\circ \partial$, where $\partial: C_{p+1}(K)\rightarrow C_p(K)$ is the boundary map. So $\delta(f)\in C^{p+1}(K: G)$. We have then defined a map $\delta: C^p(K; G)\rightarrow C^{p+1}(K; G)$ called the {\it coboundary}\index{coboundary}. Then the group $\ker(\delta: C^p(K; G)\rightarrow C^{p+1}(K; G))=Z^p(K; G)$ is called the group of {\it p-cocycles}\index{cocycle}. The group $im(\delta: C^{p-1}(K; G)\rightarrow C^p(K; G))=B^p(K; G))$ is called the group of {\it p-coboundaries}\index{coboundary}. Since $\delta(f)=f\circ\partial$, $\partial^2=0$ implies that $\delta^2=0$, so any coboundariy is a cocycle. Then we define the {p-cohomology}\index{cohomology} group as $$H^p(K; G)=Z^p(K; G)/B^p(K; G).$$
\end{definition}

We will write $H^p(K; G)$ as $H^p(K)$ when $G=Z$. 

It is convenient to use the notation $\langle c^p, c_p\rangle$ to denote the element of $G$ representing the value of the cochain $c^p$ on the chain $c_p$. Then we have $\langle\delta c^p, c_{p+1}\rangle=\langle c_p, \partial c_{p+1}\rangle.$

We still need to define 0-dimensional homology.

\begin{theorem}
Let $K$ be a simplicial complex. Then $H^0(K; G)$ is the group of 0-cochains $c^0$ such that $c^0$ is constant on the set of vertices belonging to a connected component of $|K|$. If $|K|$ is connected, then  $H^0(K; G)\cong G$, and if $G=Z$, the group is generated by the cochain which has the value 1 on each vertex.
\end{theorem}

{\bf Proof:} $H^0(K; G)=Z^0(K; G)$ as there are no coboundaries in dimension 0. If $v$ and $w$ are in the same component of $K$, then there is a 1-chain $c_1$ with $\partial c_1=v-w.$ So for any cocycle $c^0$, $$0=\langle\delta c^0, c_1\rangle=\langle c^0, \partial c_1\rangle=\langle c^0, v\rangle-\langle c^0,w\rangle.$$ On the other hand, if $c^0$ is a cochain such that $\langle c^0, v\rangle-\langle c^0,w\rangle=0$ whenever $v$ and $w$ lie in the same component of $|K|$, then for each oriented 1-simplex $\sigma$ of $K$, $$\langle\delta c^0, \sigma\rangle=\langle c^0, \partial\sigma\rangle=0.$$ Then we have that $\delta c^0=0$, so the $c^0$ is a cocycle and we are done. $\blacksquare$

Next we describe reduced and relative cohomology.

\begin{definition}
Given a complex $K$, let $\epsilon$ be the augmentation map. The diagram $$C_1(K)\xrightarrow{\partial_1}C_0(K)\xrightarrow{\epsilon}Z$$ dualizes to the diagram $$C_1(K; G)\xleftarrow{\delta_1}C_0(K; G)\xleftarrow{\tilde{\epsilon}}Z.$$ The homomorphism $\tilde{\epsilon}$ is called a {\it coaugmentation}\index{coaugmentation}. It is injective and $\delta_1\circ\tilde{\epsilon}=0.$ Define the {\it reduced cohomology}\index{reduced cohomology} of $K$ by setting $\tilde{H}^p(K; G)=H^p(K; G)$ for $p>0$ and $$\tilde{H}^0(K; G)=\ker \delta_1/im\mbox{ } \tilde{\epsilon}.$$
\end{definition}

\begin{theorem}
If $K$ is connected, then  $\tilde{H}^0(K; G)=0$. In general,  $$H^0(K; G)\cong \tilde{H}^0(K; G)\oplus G$$. 
\end{theorem}

\begin{definition}
Let $K$ be a complex and $K_0$ be a subcomplex of $K$. The group of {\it relative cochains}\index{relative cochain} in dimension $p$ is defined as $$C^p(K, K_0; G)=Hom(C_p(K, K_0), G).$$ The relative coboundary operator $\delta$ is the dual of the relative boundary operator. Then the group $\ker(\delta: C^p(K, K_0; G)\rightarrow C^{p+1}(K, K_0; G))=Z^p(K, K_0; G)$ is called the group of {\it relative p-cocycles}\index{relative cocycle}. The group $im(\delta: C^{p-1}(K, K_0; G)\rightarrow C^p(K, K_0; G))=B^p(K, K_0; G))$ is called the group of {\it relative p-coboundaries}\index{relative coboundary}. Then the {\it relative cohomology group}\index{relative cohomology} is $$H^p(K, K_0; G)=Z^p(K, K_0; G)/B^p(K, K_0; G).$$
\end{definition}

The dual of the long exact homology sequence of the pair $(K, K_0)$ leads to a long exact cohomology sequence going the other direction.

\begin{theorem}
Let $K$ be a complex and $K_0$ be a subcomplex. There exists an exact sequence $$\cdots\leftarrow H^p(K_0; G)\leftarrow H^p(K; G)\leftarrow H^p(K, K_0; G)\xleftarrow{\delta} H^{p-1}(K_0; G)\leftarrow\cdots.$$ A similar sequence exists in reduced cohomology if $K_0\neq\emptyset$. A simplical map $f: (K, K_0)\rightarrow (L, L_0)$ induces a homomorphism of long exact cohomology sequences.
\end{theorem}

Like homology theory, cohomology theory has its own set of Eilenberg-Steenrod Axioms. To describe them,.we need the idea of a {\it cochain complex}.

 Let $\mathcal{C}=\{C_p, \partial\}$ be a chain complex and $G$ and abelian group. Then the p-cochain group of $\mathcal{C}$ with coefficients in $G$ is $C^p(\mathcal{C}; G)=Hom(C_p, G)$.  Letting $\delta$ be the dual of $\partial$, we can define the cohomology groups as before. 

\begin{definition}
Let $\mathcal{C}=\{C_p, \partial\}$ and $\mathcal{C'}=\{C'_p, \partial'\}$ be chain complexes. Suppose $\phi: \mathcal{C}\rightarrow\mathcal{C'}$ is a chain map. (Recall that this means that $\partial'\phi=\phi\partial$.) Then the dual homomorphism $$C^p(\mathcal{C}; G)\xleftarrow{\tilde{\phi}}C^p(\mathcal{C'}; G)$$ commutes with $\delta$ and is called a {\it cochain map.} It carries cocycles to cocycles and coboundaries to coboundaries  so it induces a homomorphism $\phi^\ast: H^p(\mathcal{C'}; G)\rightarrow H^p(\mathcal{C}; G)$ on cohomology groups. 
\end{definition}

Now we can state the Eilenberg-Steenrod Axioms for Cohomology. See Section 4.2.3 to compare these to the axioms for homology.

Given an admissible class $\mathcal{A}$ of pairs of spaces $(X, A)$ and an abelian group $G$, a {\it cohomology theory} on $\mathcal{A}$ with coefficients in $G$ consists of the following: \begin{enumerate}
\item A function defined for each integer $p$ and each pair $(X, A)$ in $\mathcal{A}$ whose value is an abelian group $H^p(X, A; G)$.
\item A function that assigns to each continuous map $h: (X, A)\rightarrow (Y, B)$ and each integer $p$ a homomorphism $$H^p(X, A; G)\xleftarrow{h^\ast} H^p(Y, B; G).$$ 
\item A function that assigns to each pair $(X, A)$ in $\mathcal{A}$ and each integer $p$ a homomorphism $$H^p(X, A; G)\xleftarrow{\delta^\ast}H^{p-1}(A; G)$$.
\end{enumerate}

These functions satisfy the following axioms where all pairs of spaces are in $\mathcal{A}$:
\begin{itemize}
\item[$-$]{\bf Axiom 1:} If $i$ is the identity, then $i^*$ is the identity.
\item[$-$]{\bf Axiom 2:} $(kh)^*=h^*k^*.$
\item[$-$]{\bf Axiom 3:} If $f: (X, A)\rightarrow (Y, B)$, then the following diagram commutes:$$\begin{tikzpicture}
  \matrix (m) [matrix of math nodes,row sep=3em,column sep=4em,minimum width=2em]
  {
H^p(X, A; G) & H^p(Y, B; G) \\
H^{p-1}(A; G) & H^{p-1}(B; G) \\};
  \path[-stealth]
    (m-1-2) edge node [above] {$f^*$} (m-1-1)
(m-2-2) edge node [above] {$(f|A)^*$} (m-2-1)
(m-2-1) edge node [right] {$\delta^*$} (m-1-1)
(m-2-2) edge node [right] {$\delta^*$} (m-1-2)
;

\end{tikzpicture}$$
\item[$-$]{\bf Axiom 4:} (Exactness Axiom) The sequence $$\begin{tikzpicture}
  \matrix (m) [matrix of math nodes,row sep=3em,column sep=4em,minimum width=2em]
  {
\cdots & H^p(A; G) & H^p(X; G) &  H^p(X, A; G) & H^{p-1}(A; G) &\cdots \\};
  \path[-stealth]
    (m-1-2) edge (m-1-1)
(m-1-3) edge node [above] {$i^*$} (m-1-2)
(m-1-4) edge node [above] {$\pi^*$} (m-1-3)
(m-1-5) edge node [above] {$\delta^*$} (m-1-4)
  (m-1-6) edge (m-1-5)
;

\end{tikzpicture}$$ is exact where $i: A\rightarrow X$ and $\pi: X\rightarrow (X,A)$ are inclusion maps.
\item[$-$]{\bf Axiom 5:} (Homotopy Axiom) If $h, k: (X, A)\rightarrow(Y, B)$ are homotopic then $h^*=k^*$.
\item[$-$]{\bf Axiom 6:} (Excision Axiom) Given $(X, A)$, let $U$ be an open subset of $X$ such that $\bar{U}\subset A^\circ$ (i.e. the closure of $U$ is contained in the interior of $A)$. Then if $(X-U, A-U)$ is admissible, then inclusion induces an isomorphism $$H^p(X-U, A-U; G)\cong H^p(X, A; G).$$
\item[$-$]{\bf Axiom 7:} (Dimension Axiom) If $P$ is a one point space then $H^p(P; G)=0$ for $p\neq 0$, and $H^0(P; G)\cong G$. 
\end{itemize}

Note that there is no equivalent to the axiom of compact support for cohomology theory. Also, the absence of axiom 7 leads to an {\it extraordinary cohomology theory}.

As with homology, we can define singular and cellular cohomology theories. To define them, we simply modify the definition of cochains in Definition 8.2.1 to have $C_p(K)$ be the group of singular or cellular chains. For a triangulable space, simplicial, singular, and cellular cohomology all coincide.

\begin{theorem}
Let $n>0$. Then \begin{enumerate} 
\item $H^i(S^n; G)\cong G$ for $i=0, n$, 
\item $H^i(B^n, S^{n-1}; G)\cong G$ for $i=n$. 
\end{enumerate}
These groups are zero for all other values of $i$.
\end{theorem}

{\bf Proof:} The first statement easily follows from the cellular chain complex of $S^n$ and the fact that $Hom(Z, G)\cong G$. The second statement comes from the long exact sequence of the pair $(B^n, S^{n-1})$ in cohomology and the fact that $B^n$ is acyclic. $\blacksquare$

Now recall that if $X$ is a CW complex, let $$D_p(X)=H_p(X^p, X^{p-1}).$$ Let $\partial: D_p(x)\rightarrow D_{p-1}(X)$ be the composite $$\begin{tikzpicture}
  \matrix (m) [matrix of math nodes,row sep=3em,column sep=4em,minimum width=2em]
  {
 H_p(X^p, X^{p-1}) & H_{p-1}(X^{p-1}) & H_{p-1}(X^{p-1}, X^{p-2}),\\};
  \path[-stealth]

(m-1-1) edge node [above] {$\partial_*$} (m-1-2)
(m-1-2) edge node [above] {$j_*$} (m-1-3)

;

\end{tikzpicture}$$ where $j$ is inclusion. 

\begin{example}
Let $T$ be the torus and $K$ be the Klein bottle. Find their cohomology with coefficients in the group $G$.

For both of them the cellular chain complex is of the form $$\cdots\rightarrow 0 \rightarrow Z\xrightarrow{\partial_2} Z\oplus Z\xrightarrow{\partial_1} Z\rightarrow 0.$$ Let $\gamma$ generate $D_2(T)$ and $w_1$ and $z_1$ be a basis for $D_1(T)$. For the torus, both $\partial_2$ and $\partial_1$ are zero so from the dual sequence, we get $$H^2(T; G)\cong G,\mbox{        }H^1(T; G)\cong G\oplus G,\mbox{        }H^0(T; G)\cong G.$$

For the Klein bottle $K$, $\partial_1=0$, and we can choose a basis $w_1$ and $z_1$ for $D_1(K)$ such that $\partial_2\gamma=2z_1.$ The dual sequence is of the form $$\cdots\leftarrow 0 \leftarrow G\xleftarrow{\delta_2} G\oplus G\xleftarrow{\delta_1} G\leftarrow 0.$$ Now $Hom(D_1(K), G)\cong G\oplus G$ where the first summand represents homomorphisms $\phi: D_1(K)\rightarrow G$ that vanish on $z_1$ and the second summand represents homomorphisms $\psi: D_1(K)\rightarrow G$ that vanish on $w_1$. Since $\partial_1=0$, its dual is $\delta_1=0$. To compute $\delta_2$, we have that $$\langle\delta_2\phi, \gamma\rangle=\langle\phi, \partial_2\gamma\rangle=\langle\phi, 2z_1\rangle=2\langle\phi, z_1\rangle=0,$$ $$\langle\delta_2\psi, \gamma\rangle=\langle\psi, \partial_2\gamma\rangle=\langle\psi, 2z_1\rangle=2\langle\psi, z_1\rangle.$$ So $\delta_2$ takes the first summand to 0 and equals multiplication by 2 on the second summand. This gives $$H^2(K; G)\cong G/2G,\mbox{        }H^1(K; G)\cong G\oplus \ker(G\xrightarrow{2} G),\mbox{        }H^0(K; G)\cong G.$$ If $G=Z$, this becomes $$H^2(K)\cong Z_2,\mbox{        }H^1(K)\cong Z,\mbox{        }H^0(K)\cong Z.$$

Note that for the torus, the homology and cohomology groups are the same, but for the Klein bottle, they are very different. (Recall that $H_1(K)\cong Z\oplus Z_2$ and $H_2(K)=0$.)
\end{example}

I will conclude this section with a result from Munkres \cite{Mun1} that states that if two spaces have the same homology groups then they have the same cohomology groups. A {\it free} chain complex\index{free chain complex}  is one for which the chain group $C_p$ is free for all $p$. The simplicial, singular, and cellular chain complexes are all free. Also, recall that a chain map is a map that commutes with boundaries. 

\begin{theorem}
Let $\mathcal{C}$ and $\mathcal{D}$ be free chain complexes and let $\phi: \mathcal{C}\rightarrow \mathcal{D}$ be a chain map. If $\phi$ induces isomorphisms in homology in all dimensions than it induces isomorphisms in cohomology in all dimensions. Thus, if $H_p(\mathcal{C})\cong H_p(\mathcal{D})$ for all $p$, then $H^p(\mathcal{C}; G)\cong H^p(\mathcal{D}; G)$ for all $p$ and $G$.
\end{theorem}

This seems like bad news. Cohomology would be no better for classification then homology. As we will see later, the Universal Coeffieicent Theorem allows us to compute cohomology groups from homology groups and vice versa. So why do we bother with cohomology before. As I have already mentioned, cohomology has an additional ring structure that homology lacks. The next section will define {\it cup products} that turn cohomology into a ring. So why isn't homology also a ring? I will answer that question in Section 8.5.

\section{Cup Products}

I will follow the approach in \cite{Mun1} of defining cup products for singular cohomology. Then the isomorphism with simplicial and cellular cohomology will define it for those cases.

In what follows, $R$ will always denote a commutative ring with unit element. Typically, it will be either $Z$ or $Z_n$.

\begin{definition}
Let $X$ be a topological space. Let $S^p(X; R)=Hom(S_p(X), R)$ denote the group of singular $p$-cochains of $X$ with coefficients in $R$. (After all, a ring is always an abelian group under addition.) Define a map $$S^p(X; R)\times  S^q(X; R)\xrightarrow{\cup} S^{p+q}(X; R)$$ as follows: If $T: \Delta_{p+q}\rightarrow X$ is a singular $p+q$-simplex then $$\langle c^p\cup c^q, T\rangle=\langle c^p, T\circ\frak{l}(e_0, \cdots, e_p)\rangle\cdot\langle c^q, T\circ\frak{l}(e_p, \cdots, e_{p+q})\rangle.$$ (If you don't remember what this notation means, see the beginning of Section 4.3.1.) 

The cochain $c^p\cup c^q$ is called the {\it cup product}\index{cup product} of the cochains $c^p$ and $c^q$.
\end{definition}

We can think of $T\circ\frak{l}(e_0, \cdots, e_p)$ as the {\it front p-face} of $\Delta_{p+q}$ and $T\circ\frak{l}(e_p, \cdots, e_{p+q})$ as its {\it back q-face}. The multiplication is the multiplication in $R$.

\begin{theorem}
The cup product of cochains is bilinear and associative. The cochain $z^0$ whose value is 1 on each singular 0-simplex is the unit element. We also have the following coboundary formula: $$\delta(c^p\cup c^q)=(\delta c^p)\cup c^q+(-1)^pc^p\cup (\delta c^q).$$
\end{theorem}

The theorem can be proved directly from the formula for cup product. 

The coboundary formula shows that the cup product of 2 cocycles is itself a cocycle. The cohomology class of $z^p\cup z^q$ for cocycles $z^p$ and $z^q$ depends only on their cohomology classes since $$(z^p+\delta d^{p-1})\cup z^q=z^p\cup z^q +\delta(d^{p-1}\cup z^q)$$ and $$z^p\cup (z^q+\delta d^{q-1})=z^p\cup z^q+(-1)^p\delta(z^p\cup d^{q-1}).$$ This means that we can pass to cohomology.

\begin{theorem}
The cochain cup product induces a cup product $$H^p(X; R)\times H^q(X; R) \xrightarrow{\cup} H^{p+q}(X; R).$$ The cup product on cohomology is bilinear and associative with the cohomology class $\{z^0\}$ acting as the unit element. 
\end{theorem}

\begin{theorem}
If $h: X\rightarrow Y$ is a continuous map then $h^\ast$ preserves cup products. 
\end{theorem}

{\bf Proof:} The theorem follows immediately from the fact that the cochain map $h^\sharp$ preserves cup products of cochains as the value of $h^\sharp(c^p\cup c^q)$ on $T$ equals the value of $c^p\cup c^q$ on $h\circ T$ which is $$\langle c^p, h\circ T\circ\frak{l}(e_0, \cdots, e_p)\rangle\cdot\langle c^q, h\circ T\circ\frak{l}(e_p, \cdots, e_{p+q})\rangle.$$ But this value is the same as $h^\sharp(c^p)\cup h^\sharp(c^q)$. $\blacksquare$

\begin{definition}
Let $H^\ast(X: R)$ denote the direct sum $\oplus_{i=0}^\infty H^i(X; R)$. The cup product makes this group into a ring with unit element. It is called the {\it cohomology ring}\index{cohomology ring} of $X$ with coefficients in $R$. A ring in this form is called a {\it graded ring}\index{graded ring}. 
\end{definition}

Note that a homotopy equivalence induces a ring isomorphism of cohomology rings, so cohomology is a topological invariant.

Is the cohomology ring commutative? It sort of is. If $\alpha^p\in H^p(X; R)$ and $\beta^q\in H^q(X; R)$, then $$\alpha^p \cup \beta^q=(-1)^{pq} \beta^q\cup\alpha^p.$$ This property is called {\it anti-commutativity} but many books just call it {commutativity} despite not being commutative in the usual meaning. I will prove the formula in Section 8.6 with access to more tools. 

We can also define a cup product for simplicial cohomology rather than singular by modifying the cochain formula to $$\langle c^p\cup c^q, [v_0, \cdots, v_{p+q}]\rangle=\langle c^p, [v_0, \cdots, v_p]\rangle\cdot\langle c^q, [v_p, \cdots, v_{p+q}]\rangle.$$ Simplicial cup products retain all of the properties of singular cup products and they form isomorphic rings for triangulable spaces. 

Munkres computes some examples of cohomology rings but there are some difficulties. First of all, we need to find cocycles representing the cohomology groups. Also, isomorphic cell complexes may not always lead to the same cohomology ring. 

We will look at cohomology rings of product spaces in Section 8.5. Munkres discusses cohomology rings of manifolds in connection with Poincar{\'e} duality, a relation between the homology and cohomology groups of a manifold. I won't get too far into it as it is a bit of a diversion, but there is one result about the cohomology ring of real projective space that I will state that will be useful when we talk about Steenrod squares. (I will derive it in a different way in Section 11.3.1.)

\begin{theorem}
If $u$ is the non-zero element of $H^1(P^n; Z_2$) then $u^k$ (i.e. the cup product of $u$ with itself $k$ times) is the nonzero element of  $H^k(P^n; Z_2)$  for $1<k\leq n$. Then $H^\ast(P^n; Z_2)$ is a truncated polynomial ring over $Z_2$ with one generator $u$ in dimension 1 and $u^{n+1}=0$. $H^\ast(P^\infty; Z_2)$ is a polynomial ring over $Z_2$ with a single one dimensional generator.
\end{theorem} 

Finally, I will give an example of two spaces with the same cohomology groups but non-isomorphic cohomology rings. We will need one definition.

\begin{definition}
Let $X$ and $Y$ be topological spaces. The {\it wedge product}\index{wedge product}\index{$X\vee Y$} $X\vee Y$ is the union of $X$ and $Y$ with one point in common.
\end{definition}

\begin{example}
$S^1\vee S^1$ is a figure 8.
\end{example}

\begin{example}
Now consider the following scenario. You want a donut for dessert and you go to your local donut shop. You are a little disappointed when you are handed a sphere with 2 circles attached. (Actually a hollow tetrahedron with two triangles. See Figure 8.3.1.) The manager claims that it is the same thing as all of its homology and cohomology grpups are the same as those of a torus. But something doesn't seem right to you. Who is right?

Let $X=S^1\vee S^1\vee S^2$ and $T=S^1\times S^1$ be the torus. The space $X$ is a CW complex with one cell in dimension 0, two cells in dimension 1, and one cell in dimension 2. Referring to Figure 8.3.1, write the fundamental cycles $$w_1=[a, b]+[b, c]+[c, a],$$ $$z_1=[a, d]+[d, e]+[e, a],$$ $$z_2=\partial[a, f, g, h]$$ for these cells. The boundary operators in the cellular chain complex all vanish. So the cellular chain complex of $X$ is isomorphic to the cellular chain complex of the torus, $T$. Thus, the homology and cohomology groups of $X$ are isomorphic to those of $T$. 

The cohomology rings, however, are not isomorphic. We claim that the cohomology ring of $X$ is {\it trivial}\index{trivial cohomology ring}, i.e. every product of cohomology classes of positive dimension is zero. This is the case because the fundamental cycles in dimension 1 do not carry the boundaries of any two cell so all cup products of one dimensional cocycles are 0. 

Munkres shows that the product of the two generators $\alpha, \beta$ of $H^1(T)$ is $\alpha\cup\beta=\pm\gamma$ where $\gamma$ is a generator of $H^2(T)$. 

Thus $X$ and $T$ do not have isomorphic cohomology rings, which implies that they are not homeomorphic or even homotopy equivalent. So you are well within your rights to have the donut shop give you a refund.

\end{example}

\begin{figure}[ht]
\begin{center}
  \scalebox{0.8}{\includegraphics{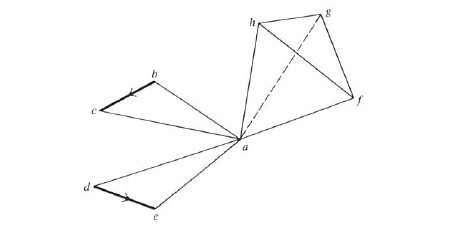}}
\caption{
\rm
$S^1\vee S^1\vee S^2$ \cite{Mun1}. 
}
\end{center}
\end{figure}

In the previous example, I did not duplicate Munkres' computation of the cohomology ring of the torus.(See \cite{Mun1} for more details.)  Finding cocycles generating the cohomology groups is surprisingly hard. This may make the full use of cohomology difficult in data science applications. Still, the rest of the chapter will show that there are more things you can do. 

\section{Universal Coefficient Theorems}

In this section, we will start to make more use of the homological algebra we learned in Section 8.1. We will have two Universal Coefficient Theorems (UCT). The UCT for cohomology allows us to compute cohomology groups given homology groups. The UCT for homology allows us to compute homology groups with coefficients in any other finitely generated abelian groups given the homology groups with integer coefficients.

Letting $\mathcal{C}=\{C_p, \partial\}$ be a chain complex and $G$ be an abelian group, recall that $C^p(\mathcal{C}; G)=Hom(C_p(\mathcal{C}), G)$. Is it true that  $H^p(\mathcal{C}; G)=Hom(H_p(\mathcal{C}), G)$? It would be nice if that were true as it would be pretty easy to compute cohomology groups from homology groups. 

The answer is "almost but not quite". But these will be terms in an exact sequence so we will often have all the information we need.

\begin{definition}
For a chain complex $\mathcal{C}=\{C_p, \partial\}$ and an abelian group $G$, there is a map $$Hom(C_p, G)\times C_p\rightarrow G$$ which takes the pair $(c^p, d_p)$ to the element  $\langle c^p, d_p\rangle\in G$. This is a bilinear map called the {\it evaluation map}. It induces a bilinear map called the {\it Kronecker index}\index{Kronecker index} in which the image of $(\alpha^p, \beta_p)$ is denoted $\langle\alpha^p, \beta_p\rangle$.
\end{definition}

We need to check that this map is well defined. Let $z^p$ be a cocycle and $z_p$ be a cycle. then $$\langle z^p+\delta d^{p-1}, z_p\rangle=\langle z^p, z_p\rangle+\langle d^{p-1}, \partial z_p\rangle=\langle z^p, z_p\rangle+0=\langle z^p, z_p\rangle$$ and $$\langle z^p, z_p+\partial d_{p+1}\rangle=\langle z^p, z_p\rangle+\langle \delta z^p, d_{p+1}\rangle=\langle z^p, z_p\rangle+0=\langle z^p, z_p\rangle.$$

\begin{definition}
The {\it Kronecker map}\index{Kronecker map} is the map $$\kappa : H^p(\mathcal{C}; G)\rightarrow Hom(H_p(\mathcal{C}); G)$$ that sends $\alpha^p$ to the homomorphism $\langle\alpha^p, \rangle$. So $(\kappa(\alpha^p))(\beta_p)=\langle\alpha^p, \beta_p\rangle.$ The map $\kappa$ is a homomorphism because the Kronecker index is linear in the first variable.
\end{definition}

I will now state the UCT for cohomology. See \cite{Mun1} Section 53 for the proof.

\begin{theorem}
{\bf Universal Coefficient Theorem for Cohomology}\index{universal coefficient theorem for cohomology} If $\mathcal{C}$ is a free chain complex and $G$ is an abelian group, then there is an exact sequence $$0\leftarrow Hom(H_p(\mathcal{C}), G)\xleftarrow{\kappa}H^p(\mathcal{C}; G)\leftarrow Ext(H_{p-1}(\mathcal{C}, G)\leftarrow 0$$ that is natural with respect to homomorphisms induced by chain maps. It splits but not naturally. 

If $(X, A)$ is a topological pair, then the exact sequence takes the form $$0\leftarrow Hom(H_p(X, A), G)\xleftarrow{\kappa}H^p(X, A; G)\leftarrow Ext(H_{p-1}(X, A), G)\leftarrow 0$$ and it is natural with respect to homomorphisms induced by continuous maps. ($H_i(X, A)$ is replaced by $H_i(X)$ in the special case where $A=\emptyset$.) This sequence also splits but not naturally. 
\end{theorem}

The next result is stronger than Theorem 8.2.5 and follows from the naturality of the exact sequence in the UCT. 

\begin{theorem}
Let $\mathcal{C}$ and $\mathcal{D}$ be free chain complexes and let $\phi: \mathcal{C}\rightarrow \mathcal{D}$ be a chain map. If $\phi_\ast: H_i(\mathcal{C})\rightarrow H_i(\mathcal{D})$ is an isomorphism for $i=p, p-1$, then  $$H^p(\mathcal{C}; G)\xleftarrow{\phi^\ast} H^p(\mathcal{D}; G)$$ is an isomorphism.
\end{theorem}

\begin{example}
Let $T$ be the torus. From Example 4.3.4, we know that $H_0(T)\cong Z$, $H_1(T)\cong Z\oplus Z$ and $H_2(T)\cong Z$. What are the cohomology groups?

Since, $T$ is connected, we know that $H^0(T)\cong Z$ by Theorem 8.2.1. Now from Section 8.1, $Ext(Z, Z)=0$. So we have from the exact sequence in the UCT that $H^p(T)\cong Hom(H_p(T), Z)$. Since $Hom(Z, Z)\cong Z$, we have $$H^1(T)\cong Hom(H_1(T), Z)\cong Hom(Z\oplus Z, Z)\cong Z\oplus Z$$ and  $$H^2(T)\cong Hom(H_2(T), Z)\cong Hom(Z, Z)\cong Z.$$ This agrees with Example 8.2.1.
\end{example}

\begin{example}
Let $K$ be the Klein bottle.  From Example 4.3.4, we know that $H_0(K)\cong Z$, $H_1(K)\cong Z\oplus Z_2$ and $H_2(K)=0$. What are the cohomology groups?

Again, $K$ is connected, so $H^0(K)\cong Z.$ For $H^1(K)$, we have the exact sequence  $$0\leftarrow Hom(H_1(K), Z)\leftarrow H^1(K)\leftarrow Ext(H_0(K), Z)\leftarrow 0.$$ Now  $Ext(H_0(K), Z)=Ext(Z, Z)=0$, so $$H^1(K)\cong Hom(H_1(K), Z)=Hom(Z\oplus Z_2, Z)=Hom(Z, Z)\oplus Hom(Z_2, Z)=Z\oplus 0=Z.$$ For $H^2(K)$, we have the exact sequence  $$0\leftarrow Hom(H_2(K), Z)\leftarrow H^2(K)\leftarrow Ext(H_1(K), Z)\leftarrow 0.$$ Now $H_2(K)=0$, so $Hom(H_2(K), Z)=0.$ So $$H^2(K)\cong Ext(H_1(K), Z)=Ext(Z\oplus Z_2, Z)=Ext(Z, Z)\oplus Ext (Z_2, Z)=0\oplus Z_2=Z_2.$$ This agrees with Example 8.2.1.
\end{example}

There is an analogous UCT for homology. We want to compute homology with coefficients in some other group given homology with integer coefficients. Tensor product and Tor replace Hom and Ext, and the arrows point in the other direction.

\begin{theorem}
{\bf Universal Coefficient Theorem for Homology}\index{universal coefficient theorem for homology} If $\mathcal{C}$ is a free chain complex and $G$ is an abelian group, then there is an exact sequence $$0\rightarrow H_p(\mathcal{C})\otimes G\rightarrow H_p(\mathcal{C}; G)\rightarrow Tor(H_{p-1}(\mathcal{C}), G)\rightarrow 0$$ that is natural with respect to homomorphisms induced by chain maps. It splits but not naturally. 

If $(X, A)$ is a topological pair, then the exact sequence takes the form $$0\rightarrow H_p(X, A)\otimes G\rightarrow H_p(X, A; G)\rightarrow Tor(H_{p-1}(X, A), G)\rightarrow 0$$ and it is natural with respect to homomorphisms induced by continuous maps. This sequence also splits but not naturally. 
\end{theorem}

\begin{theorem}
Let $\mathcal{C}$ and $\mathcal{D}$ be free chain complexes and let $\phi: \mathcal{C}\rightarrow \mathcal{D}$ be a chain map. If $\phi_\ast: H_i(\mathcal{C})\rightarrow H_i(\mathcal{D})$ is an isomorphism for $i=p, p-1$, then  $$\phi_\ast: H_p(\mathcal{C}; G)\rightarrow H_p(\mathcal{D}; G)$$ is an isomorphism.
\end{theorem}

\begin{example}
Let $T$ be the torus.  Using the fact that $H_0(T)\cong Z$, $H_1(T)\cong Z\oplus Z$ and $H_2(T)\cong Z$, what are the homology groups with coefficients in $Z_2$?

Since, $T$ is connected, we know that $H_0(T; Z_2)\cong Z_2$. Now from Section 8.1, $Tor(Z, Z_2)=0$. So we have from the exact sequence in the UCT for homology that $H_p(T; Z_2)\cong H_p(T)\otimes Z_2$. Since $Z\otimes Z_2\cong Z_2$, we have $$H_1(T; Z_2)\cong H_1(T)\otimes Z_2\cong (Z\oplus Z)\otimes Z_2\cong Z_2\oplus Z_2$$ and  $$H_2(T; Z_2)\cong H_2(T)\otimes Z_2\cong Z\otimes Z_2\cong Z_2.$$
\end{example}

\begin{example}
Let $K$ be the Klein bottle.  Using the fact that $H_0(K)\cong Z$, $H_1(K)\cong Z\oplus Z_2$ and $H_2(K)=0$, what are the homology groups with coefficients in $Z_2$?

Again, $K$ is connected, so $H_0(K; Z_2)\cong Z_2.$ For $p=1$, the sequence becomes $$0\rightarrow H_1(K)\otimes Z_2\rightarrow H_1(K; Z_2)\rightarrow Tor(H_0(K), Z_2)\rightarrow 0.$$ Now $Tor(H_0(K), Z_2)=Tor(Z, Z_2)=0$, so $$H_1(K: Z_2)\cong H_1(K)\otimes Z_2\cong (Z \oplus Z_2)\otimes Z_2\cong Z_2\oplus Z_2.$$ For $p=2$, the sequence becomes $$0\rightarrow H_2(K)\otimes Z_2\rightarrow H_2(K; Z_2)\rightarrow Tor(H_1(K), Z_2)\rightarrow 0.$$ Now $H_2(K)=0$, so $$H_2(K: Z_2)\cong Tor(H_1(K), Z_2)\cong Tor(Z\oplus Z_2, Z_2)\cong Z_2.$$
\end{example}

In $Z_2$ coefficients, homology of the torus and the Klein bottle are the same. In $Z_2$ the "twist" is ignored.

\section{Homology and Cohomology of Product Spaces}

In this section, we will learn how to compute homology and cohomology groups of a product space. We will also see what happens to a cohomology ring. The main result is called the K{\" u}nneth Theorem. There will be a version for both homology and cohomology.

\begin{definition}
Let $\mathcal{C}=\{C_p, \partial\}$ and $\mathcal{C'}=\{C'_p, \partial'\}$ be chain complexes. Their {\it tensor product}\index{tensor product of complexes} is a new complex denoted $\mathcal{C}\otimes\mathcal{C'}$ whose chain group in dimension $m$ is defined by $$(\mathcal{C}\otimes\mathcal{C'})_m=\oplus_{p+q=m}C_p\otimes C'_q,$$ whose boundary $\overline{\partial}$ is defined as $$\overline{\partial}(c_p\otimes c'_q)=\partial c_p\otimes c'_q+(-1)^p c_p\otimes \partial'c'_q.$$ We define $\overline{\partial}$ as a function on $C_p\times C'_q,$ which is bilinear so it induces a homomorphism on tensor product.
\end{definition}

It is easy to see that $\overline{\partial}^2=0$. Also, the tensor product of two chain maps is a chain map on the tensor product complexes. 

The tensor product of two free chain complexes is itself free.

If $K$ and $L$ are simplicial complexes, then $|K|\times |L|$ is a CW complex with cells $\sigma\times\tau$ for $\sigma\in K$ and $\tau \in L$. Letting $\mathcal{D}(K\times L)$ be the corresponding cellular chain complex, the group $D_m(K\times L)$ is free abelian with a basis element for each cell $\sigma^p\times \tau^q$ for which $p+q=m.$ This basis element is a fundamental cycle for $\sigma^p\times\tau^q$.

Orienting the simplices of $K$ and $L$, they form bases for the simplicial chain complexes $\mathcal{C}(K)$ and $\mathcal{C}(L)$. The $m$-dimensional chain group of $(\mathcal{C}(K)\otimes\mathcal{C}(L))$ consists of the elements $\sigma^p\otimes\tau^q$ for $p+q=m.$

From this one to one correspondence, we can choose a group isomorphism that respects boundaries and becomes an isomorphism of chain complexes. This leads to the following result. See \cite{Mun1} for the proof.

\begin{theorem}
If $K$ and $L$ are simplicial complexes and if $|K|\times|L|$ is triangulated so that each cell $\sigma\times\tau$ is the polytope of a subcomplex then $$\mathcal{C}(K)\otimes \mathcal{C}(L)\cong\mathcal{D}(K\times L).$$ So $\mathcal{C}(K)\otimes \mathcal{C}(L)$ can be used to compute the homology of $|K|\times |L|$.
\end{theorem}

Next we get to the K{\" u}nneth Theorem itself. As with the Universal Coefficient Theorems, we will start with chain complexes and then move to topological spaces.

Let $$\Theta: H_p(\mathcal{C})\otimes H_q(\mathcal{C'})\rightarrow H_{p+q}(\mathcal{C}\otimes\mathcal{C'})$$ be defined by $$\Theta(\{z_p\}\otimes\{z'_q\})=\{z_p\otimes z'_q\},$$ where $z_p$ is a $p$-cycle of $\mathcal{C}$, and $z'_q$ is a $q$-cycle of $\mathcal{C'}$. It is easily shown that $\{z_p\otimes z'_q\}$ is a cycle and the map is well defined. Extending $\Theta$ to a direct sum, the map $$\Theta: \oplus_{p+q=m}H_p(\mathcal{C})\otimes H_q(\mathcal{C}')\rightarrow H_m(\mathcal{C}\otimes\mathcal{C'}),$$ can be shown to be a monomorphism, and its image is a direct summand. 

The K{\" u}nneth Theorem identifies the cokernel of $\Theta$. It is our old friend Tor. Again, the reader is referred to \cite{Mun1} for the rather lengthy proof.

\begin{theorem}
{\bf K{\" u}nneth Theorem for Chain Complexes:} Let $\mathcal{C}$ and $\mathcal{C'}$ be free chain complexes. There is an exact sequence $$0\rightarrow \oplus_{p+q=m}H_p(\mathcal{C})\otimes H_q(\mathcal{C}')\xrightarrow{\Theta} H_m(\mathcal{C}\otimes\mathcal{C'})\rightarrow\oplus_{p+q=m}Tor(H_{p-1}(\mathcal{C}), H_q(\mathcal{C'}))\rightarrow 0$$ which is natural with respect to homomorphisms induced by chain maps. The sequence splits but not naturally.
\end{theorem}

\begin{theorem}
If $K$ and $L$ are finite simplicial complexes, the K{\" u}nneth Theorem implies that since  $\mathcal{C}(K)$ and $\mathcal{C'}(L)$ are free, $$H_m(|K|\times |L|)\cong\oplus_{p+q=m}[H_p(K)\otimes H_q(L)\oplus Tor(H_{p-1}(K), H_q(L)].$$
\end{theorem}

\begin{example}
Recall the donut shop from Example 8.3.2. They can't tell a donut form a sphere wedged with 2 circles, so you decide to try a different one. You get a pleasant surprise when a couple of blocks away, you discover a shop that not only sells donuts but the Cartesian product of two spheres of any dimension. Suppose you order a hyperdonut\index{hyperdonut} of the form $S^r\times S^s$ for some $r, s>0$, and you are very curious about its homology. Since $S^r$ has $H_p(S^r)\cong Z$ for $p=0, r$ and is 0 otherwise, the Tor terms are all 0 as $Tor(Z, G)=0$ for any abelian group $G$. So for $p+q=m$, $$H_m(S^r\times S^s)\cong \oplus_{p+q=m} H_p(S^r)\otimes H_q(S^s).$$ So for $r\neq s$, $$H_m(S^r\times S^s)\cong Z$$ if $m=0, r, s, r+s$ and $$H_m(S^r\times S^s)=0$$ otherwise.  

In the case $r=s$, $$H_m(S^r\times S^r)\cong Z$$ for $m=0, 2r$, $$H_m(S^r\times S^r)\cong Z\oplus Z$$ for $m=r,$ and $$H_m(S^r\times S^r)=0$$ otherwise.
\end{example}

\begin{example}
For field coefficients and free chain complexes, the Tor terms also go away. If we have finite simplicial complexes $K$ and $L$, we get a vector space isomorphism $$\oplus_{p+q=m}H_p(K; F)\otimes_F H_q(L; F)\rightarrow H_m(|K|\times |L|; F).$$
\end{example}

There is some background I need to add here. It will be important for applying the K{\" u}nneth Theorem and for looking at the ring structure of product space cohomology in the next section. We will also need it for the construction of Steenrod squares in Chapter 11.

The first step is to establish a relationship between the singular chain complexes of two topological spaces $X$ and $Y$ and the singular chain complex of their product. This accomplished by way of the {\it Eilenberg-Zilber Theorem}. The proof relies on the very abstract idea of {\it acyclic models}.

$C$ will be an arbitrary category but almost always the category of topological spaces or pairs of spaces and the continuous maps between them. $A$ is the category of augmented chain complexes and chain maps. 

\begin{definition}
Let $G$ be a functor from $C$ to $A$. Given an object $X$ of $C$, let $G_p(X)$ be the $p$-dimensional group of the augmented chain complex $G(X)$. Let $\mathcal{M}$ be a collection of objects of $C$ called {\it models}. Then $G$ is {\it acyclic} relative to $\mathcal{M}$ if $G(X)$ is acyclic for each $X\in\mathcal{M}.$\index{acyclic model} We say that $G$ is {\it free} relative to $\mathcal{M}$ if for each $p\geq 0$, there exsts:\begin{enumerate}
\item An index set $J_p$.
\item An indexed family $\{M_\alpha\}_{\alpha\in J_p}$ of objects of $\mathcal{M}$. 
\item An indexed family $\{i_\alpha\}_{\alpha\in J_p}$, where $i_\alpha\in G_p(M_\alpha)$ for each $\alpha$. 
\end{enumerate}
Given $X$, the elements $$G(f)(i_\alpha)\in G_p(X)$$ are distinct and form a basis for $G_p(X)$ as $f$ ranges over $hom(M_\alpha, X)$, and $\alpha$ ranges over $J_p$.
\end{definition}

\begin{example}
Consider the singular chain functor form the category $C$ of topological spaces to $A$. Let $\mathcal{M}$ be the collection $\{\Delta_p\}$ for $p\geq 0$. The functor is acyclic relative to $\mathcal{M}$. It is also free, since for each $p$, let $J_p$ have only one element and the corresponding object of $\mathcal{M}$ be $\Delta_p$. The corresponding element of $S_p(\Delta_p)$ is the identity simplex $i_p$. As $T$ ranges over all continuous maps $\Delta_p\rightarrow X$ the elements $T_\sharp(i_p)=T$ form a basis for $S_p(X)$. 
\end{example}

\begin{example}
Let $G$ be defined on pairs of spaces and $G(X, Y)=\mathcal{S}(X\times Y)$, $G(f, g)=(f\times g)_\sharp$. Let $\mathcal{M}=\{(\Delta_p, \Delta_q)\}$ with $p, q\geq 0$. Then $G$ is acyclic relative to $\mathcal{M}$, since $\Delta_p\times\Delta_q$ is contractible. To show $G$ is free relative to $\mathcal{M}$: For each index $p$, let $J_p$ consist of a single element with the corresponding object of $\mathcal{M}$ being $(\Delta_p, \Delta_p)$ and the corresponding element of $S_p(\Delta_p\times\Delta_p)$ being the diagonal map $d_p(x)=(x, x)$. As $f$ and $g$ range over all maps from $\Delta_p$ to $X$ and $Y$ respectively, $(f\times g)_\sharp(d_p)$ ranges over all maps $\Delta_p\rightarrow X\times Y$ which is a basis for $S_p(X\times Y)$. 
\end{example}

If $G$ is free relative to $\mathcal{M}$, then it is automatically free relative to any larger collection. If it is acyclic relative to $\mathcal{M}$, then it is automatically acyclic relative to any smaller collection. So for $G$ to be both free and acyclic we must choose $\mathcal{M}$ to be exactly the right size.

\begin{theorem}
{\bf Acyclic Model Theorem:} Let $G$ and $G'$ be functors from category $C$ to the category $A$ of augmented chain complexes and chain maps. Let $\mathcal{M}$ be a collection of objects of $C$. If $G$ is free relative to $\mathcal{M}$ and $G'$ is acyclic relative to $\mathcal{M}$, then the following hold:\begin{enumerate}
\item There is a natural transformation $T_X$ of $G$ to $G'$. 
\item Given two natural transformations $T_X$, $T'_X$ of $G$ to $G'$, there is a natural chain homotopy $D_X$ between them. In other words, for each object $X$ of $C$ we have $T_X: G_p(X)\rightarrow G'_p(X)$ and $D_X: G_p(X)\rightarrow G'_{p+1}(X)$, and naturality implies that for each $f\in Hom(X, Y)$, $$G'(f)\circ T_X=T_Y\circ G(f),$$ and $$G'(f)\circ D_X=D_Y\circ G(f).$$
\end{enumerate}
\end{theorem}

The next result will use the Acyclic Model Theorem to show that the chain complex $\mathcal{S}(X)\otimes\mathcal{S}(Y)$ can be used to compute the singular homology of $X\times Y$. 

\begin{theorem}
{\bf Eilenberg-Zilber Theorem:} For every pair $X, Y$ of topological spaces there are chain maps $$\begin{tikzpicture}
  \matrix (m) [matrix of math nodes,row sep=3em,column sep=4em,minimum width=2em]
  {
\mathcal{S}(X)\otimes\mathcal{S}(Y) & \mathcal{S}(X\times Y)\\};
 
\path[->] ($(m-1-1.east)$)+(0, .1) edge node [above] {$\mu$}($ (m-1-2.west)+(0, .1)$);
\path[->] ($(m-1-2.west)$)+(0, -.1)  edge node [below] {$\nu$} ($(m-1-1.east)+(0, -.1)$)

;

\end{tikzpicture}$$ that are chain-homotopy inverse to each other  (i.e. their composition in either direction is chain homotopy equivalent to the identity.) They are natural with respect to chain maps induced by continuous maps.
\end{theorem}

Example 8.5.4 is pretty much the proof. The theorem also justifies the use of topological spaces in the  K{\" u}nneth Theorem. We state it in the following form.

\begin{theorem}
{\bf K{\" u}nneth Theorem for Topological Spaces:} Let $X$ and $Y$ be topological spaces. There is an exact sequence $$0\rightarrow \oplus_{p+q=m}H_p(X)\otimes H_q(Y)\rightarrow  H_m(X\times Y)\rightarrow\oplus_{p+q=m}Tor(H_{p-1}(X), H_q(Y))\rightarrow 0$$ which is natural with respect to homomorphisms induced by continuous maps. The sequence splits but not naturally.
\end{theorem}

The monomorphism $$H_p(X)\otimes H_q(Y)\rightarrow H_m(X\otimes Y)$$ is called the {\it homology cross product}\index{homology cross product}. It is equal to the composite $$H_p(X)\otimes H_q(Y)\xrightarrow{\Theta} H_m(\mathcal{S}(X)\otimes\mathcal{S}(Y))\xrightarrow{\mu_\ast} H_m(X\otimes Y)$$ where $\Theta$ is induced by inclusion and $\mu$ is the chain equivalence from the Eilenberg-Zilber Theorem. We will not use this cross product much but a similar product in cohomology will be much more interesting. 

We will need an explicit formula for the other Eilenberg-Zilber equivalence $\nu$.

\begin{theorem}
Let $\pi_1: X\times Y\rightarrow X$ and $\pi_2: X\times Y\rightarrow Y$ be projections. Let $$\nu: S_m(X\times Y)\rightarrow \oplus_{p+q=m} S_p(X)\otimes S_q(Y)$$ be defined as $$\nu(T)=\sum_{i=0}^m[\pi_1\circ T\circ \frak{l}(e_0,\cdots, e_i)]\otimes [\pi_2\circ T\circ \frak{l}(e_i,\cdots, e_m)].$$ Then $\nu$ is a natural chain map that preserves augmentation. 
\end{theorem}

Using this map, our next goal is to define a cohomology cross product. We will let $R$ be a commutative ring with unit element.

Let $\mathcal{C}$ be a chain complex. This gives rise to a cochain complex $Hom(\mathcal{C}, R)$ whose $p$-dimensional group is $Hom(C_p, R)$. If $\mathcal{C'}$ is another chain complex, we can form the tensor product $Hom(\mathcal{C}, R)\otimes Hom(\mathcal{C'}, R)$ whose $m$-dimensional cochain group is $$\oplus_{p+q=m} Hom(C_p, R)\otimes Hom(C'_q, R)$$ and the coboundary is given by $$\overline{\delta}(\phi^p\otimes\psi^q)=\delta\phi^p\otimes\psi^q+(-1)^p\phi^p\otimes\delta'\psi^q.$$

Instead of taking Hom first and then tensor product, we could perform the two operations in the opposite order and form a different cochain complex. We will define a homomorphism between them. 

\begin{definition}
Let $p+q=m$. We define the homomorphism $$\theta: Hom(C_p, R)\otimes Hom(C'_q, R)\rightarrow Hom((\mathcal{C}\otimes\mathcal{C'})_m, R)$$ by $$\langle\theta(\phi^p\otimes\psi^q), c_r\otimes c'_s\rangle=\langle\phi^p, c_r\rangle\cdot\langle\psi^q, c'_s\rangle.$$ Note that this is zero unless $p=r$ and $q=s$.
\end{definition}

\begin{theorem}
The homomorphism $\theta$ is a natural cochain map.
\end{theorem}

We  have the corresponding cohomology map $$\Theta: H^p(\mathcal{C}; R)\otimes H^q(\mathcal{C'}; R)\rightarrow H^m(\mathcal{C}\otimes\mathcal{C'}); R)$$ for $m=p+q$ defined by $$\Theta(\{\phi\}\otimes\{\psi\})=\{\theta(\phi\otimes\psi)\}.$$

\begin{definition}
Let $m=p+q$. The {\it cohomology cross product}\index{cohomology cross product} is the composite $$H^p(X)\otimes H^q(Y)\xrightarrow{\Theta} H^m(\mathcal{S}(X)\otimes\mathcal{S}(Y))\xrightarrow{\nu^\ast} H^m(X\times Y)$$ where $\nu$ is the Eilenberg-Zilber equivalence. The image of $\alpha^p\otimes\beta^q$ under this homomorphism is denoted $\alpha^p\times\beta^q$.
\end{definition}

Using the homomorphism $\Theta$, we can now state the K{\" u}nneth Theorem for Cohomology.

\begin{theorem}
{\bf K{\" u}nneth Theorem for Cohomology:} Let $\mathcal{C}$ and $\mathcal{C'}$ be free finitely generated chain complexes that vanish below some specified dimension (usually dimension 0). There is an exact sequence $$0\rightarrow \oplus_{p+q=m}H^p(\mathcal{C})\otimes H^q(\mathcal{C}')\xrightarrow{\Theta} H^m(\mathcal{C}\otimes\mathcal{C'})\rightarrow\oplus_{p+q=m}Tor(H^{p+1}(\mathcal{C}), H^q(\mathcal{C'}))\rightarrow 0$$ which is natural with respect to homomorphisms induced by chain maps. The sequence splits but not naturally.
\end{theorem}

Denoting the cohomology cross product by $\times$, we can state the version for topological space. 

\begin{theorem}
Let $X$ and $Y$ be topological spaces and suppose $H_i(X)$ is finitely generated for each $i$. Then there is a natural exact sequence  $$0\rightarrow \oplus_{p+q=m}H^p(X)\otimes H^q(Y)\xrightarrow{\times} H^m(X\times Y)\rightarrow\oplus_{p+q=m}Tor(H^{p+1}(X), H^q(Y))\rightarrow 0.$$ It splits but not naturally if $H_i(Y)$ is finitely generated for each $i$.
\end{theorem}

\section{Ring Structure of the Cohomology of a Product Space}
In this section we will look more closely at the ring structure of a product space in cohomology. We will also settle the question of why cohomology has a ring structure but homology does not.

Again, we let $R$ be a commutative ring with unit element.

The cohomology cross product is the composite of maps $\Theta$ and $\nu^\ast$ which are induced by cochain maps $\theta$ and $\tilde{\nu}$. Defining $c^p\times c^q=\tilde{\nu}\theta(c^p\otimes c^q)$, we call this the {\it cochain cross product.} It is given by the formula $$\langle c^p\times c^q, T\rangle=\langle c^p, \pi_1\circ T\circ \frak{l}(e_0,\cdots, e_p)\rangle\cdot\langle c^q, \pi_2\circ T\circ \frak{l}(e_p,\cdots, e_{p+q})\rangle,$$ where $T: \Delta_{p+q}\rightarrow X\times Y.$ This says that the value of $c^p\times c^q$ on $T$ equals the value of $c^p$ on the front face of the first component of $T$ times the value of $c^q$ on the back face of the second component of $T$.

The basic properties of the cross products are described in the following result. 

\begin{theorem}
\begin{enumerate}
\item If $\lambda: X\times Y\rightarrow Y\times X$ is the map that reverses coordinates, then $$\lambda^\ast(\beta^q\times \alpha^p)=(-1)^{pq}\alpha^p\times \beta^q.$$
\item In $H^\ast(X\times Y\times Z; R)$, we have $(\alpha\times\beta)\times\gamma=\alpha\times(\beta\times\gamma).$
\item Let $\pi_1: X\times Y\rightarrow X$ and $\pi_2: X\times Y\rightarrow Y$ be projections. Let $1_X$ and $1_Y$ be the unit elements of $H^\ast(X; R)$ and $H^\ast(Y; R)$ respectively. Then $$\pi^\ast_1(\alpha^p)=\alpha^p\times 1_Y$$ and $$\pi^\ast_2(\beta^q)=1_X\times\beta^q.$$
\end{enumerate}
\end{theorem}

Now for the punchline: The cup product and cross product have very similar formulas. The next result is the connection between them.

\begin{theorem}
Let $d: X\rightarrow X\times X$ be the diagonal map $d(x)=(x, x)$. Then $$d^\ast(\alpha^p\times\beta^q)=\alpha^p\cup\beta^q.$$
\end{theorem}

{\bf Proof:} Let $z^p$ and $z^q$ be representative cocycles for $\alpha^p$ and $\beta^q$ respectively. If $T: \Delta_{p+q}\rightarrow X$ is a singular simplex, then $$\langle d^\sharp(z^p\times z^q), T\rangle=\langle z^p\times z^q, d\circ T\rangle=\langle z^p, \pi_1\circ(d\circ T)\circ\frak{l}(e_0, \cdots, e_p)\rangle\cdot\langle z^q, \pi_2\circ(d\circ T)\circ\frak{l}(e_p, \cdots, e_{p+q})\rangle$$ $$=\langle z^p, T \circ\frak{l}(e_0, \cdots, e_p)\rangle\cdot\langle z^q, T\circ\frak{l}(e_p, \cdots, e_{p+q})\rangle=\langle c^p\cup c^q\rangle,$$ where we have used the fact that $\pi_1\circ d=\pi_2\circ d=i_x$. $\blacksquare$

So now we have the reason why homology isn't a ring. Both homology and cohomology have cross products but consider the diagram in homology:
$$H_p(X)\times H_q(X)\xrightarrow{\times}H_{p+q}(X\times X)\xleftarrow{d_\ast}H_{p+q}(X).$$ The problem is that for homology, $d_\ast$ goes in the wrong direction. 

Anti-commutativity for the cup product is now easy. 

\begin{theorem}
$$\alpha^p\cup\beta^q=(-1)^{pq}\beta^q\cup\alpha^p.$$
\end{theorem}

{\bf Proof:} Let $\lambda: X\times X\rightarrow X\times X$ be the map that interchanges coordinates. Using the fact that $\lambda\circ d=d$, we have that $$\alpha^p\cup \beta^q=d^\ast(\alpha^p\times\beta^q)=(\lambda\circ d)^\ast(\alpha^p\times\beta^q)=d^\ast((-1)^{pq}\beta_q\times\alpha^p)=(-1)^{pq}\beta^q\cup\alpha^p.$$ $\blacksquare$

I will conclude this section with two more properties of cohomology rings of products.

\begin{theorem}
In the ring $H^\ast(X\times Y; R)$ we have: $$(\alpha\times\beta)\cup(\alpha'\times\beta')=(-1)^{(\mbox{dim } \beta)(\mbox{dim }\alpha')}(\alpha\cup\alpha')\times(\beta\cup\beta').$$
\end{theorem}

\begin{definition}
The tensor product $H^\ast(X; R)\otimes H^\ast(Y; R)$ has the structure of a ring where for $\beta\in H^p(Y; R)$ and $\alpha'\in H^q(X; R)$, $$(\alpha\otimes\beta)\cup(\alpha'\otimes\beta')=(-1)^{pq}(\alpha\cup\alpha')\otimes(\beta\cup\beta').$$
\end{definition}

\begin{theorem}
If $H_i(X)$ is finitely generated for all $i$, then the cross product defines a monomorphism of rings $$H^\ast(X)\otimes H^\ast(Y)\rightarrow H^\ast(X\times Y).$$ If $F$ is a field, it defines an isomorphism of algebras $$H^\ast(X; F)\otimes_F H^\ast(Y; F)\rightarrow H^\ast(X\times Y; F).$$
\end{theorem}

\begin{example}
The hyperdonut shop from Example 8.5.1 was so good, you can't wait to get back there. Now that you have computed the homology groups, you would like to compute the cohomology rings. Consider $S^r\times S^s$ for $r, s>0$. Let $\alpha\in H^r(S^r)$ and $\beta\in H^s(S^s)$ be generators. Then since the torsion terms in the K{\" u}nneth Theorem vanish, $H^\ast (S^r\times S^s)$ is free of rank 4 with basis $1\times 1, \alpha\times 1, 1\times \beta,$ and $\alpha\times\beta.$ The only nonzero product of elements of positive dimension is $(\alpha\times1)\cup(1\times \beta)=(\alpha\times\beta)$.
\end{example}

\begin{example}
Consider the cohomology ring of $P^2\times P^2$ with coefficients in $Z_2.$ Let $u\in H^1(P^2)$ be the nonzero element. We know from Theorem 8.3.4 that $u^2=u\cup u$ is nonzero and $u^n=0$ for $n>2$. By Theorem 8.6.5, $H^\ast (P^2\times P^2; Z_2)$ is a vector space of dimension 9. The basis is:\begin{align*}
\mbox{dim }0\hspace{.3 in} & 1\times 1\\
\mbox{dim }1\hspace{.3 in} & u\times 1, 1\times u\\
\mbox{dim }2\hspace{.3 in} & u^2\times1, u\times u,1\times u^2\\
\mbox{dim }3\hspace{.3 in}& u^2\times u, u\times u^2\\
\mbox{dim }4\hspace{.3 in} & u^2\times u^2.
\end{align*}

Multiplication uses Theorem 8.6.4, and we don't have to worry about signs as we are working in $Z_2$. For example, $$(u^2\times u)\cup (1\times u)=(u^2\cup 1)\times(u\cup u)=u^2\times u^2.$$
\end{example}

\section{Persistent Cohomology}
Although Edelsbrunner and Harer \cite{EH} discuss cohomology, the earliest source I know of is the paper of deSilva et. al \cite{DMV}. The earliest paper that describes an algorithm to compute persistent cohomology is the Ph.D. Thesis of Harer's student at Duke, Aubrey HB \cite{HB}. (I was curious about her unusual last name consisting of 2 initials. According to the biography in her thesis she was born Aubrey Rae Hebda-Bolduc.) In this section, we will look at her description of cohomology and cup product as they relate to persistence. There is a heavy emphasis on manifolds, duality, and characteristic classes. I will briefly discuss the latter in Chapter 11 as they are one of the main applications of Steenrod squares, but I want to mainly focus on other areas. Duality is covered extensively in Munkres \cite{Mun1}, while the classic book on characteristic classes is Milnor and Stasheff \cite{MS}.
 
All of the material in this section is from \cite{HB}

Suppose we have a filtration $$\emptyset=K_0\subset K_1\subset\cdots\subset K_n=K.$$ Suppose we add one simplex at each step so that $K_i=K_{i-1}\cup\sigma_i.$ Also, cohomology will be over $Z_2$.

\begin{definition}
Let $\psi_{j, i}^q: H^q(K_j)\rightarrow H^q(K_i)$ for $i\leq j$ be the map induced by inclusion. (Remember that the arrow goes in the reverse direction.) We say that a $q$-dimensional cohomology class $\Gamma$ is {\it born} at $K_j$ if $\Gamma\in H^q(K_j)$, but $\Gamma\notin Im(\psi_{j+1, j})$. We say that $\Gamma$ {\it dies} entering $K_i$ if $i$ is the largest index such that there exists a $\Lambda\in H^q(K_{j+1})$ with $\psi_{j+1, i}^q(\Lambda)=\psi_{j, i}^q(\Gamma)\in H^q(K_i)$. The {\it copersistence}\index{copersistence} of $\Gamma$ is $j-i$.
\end{definition}

For a simplex $\sigma$ we write $\sigma^\ast$ for the cochain which is one on $\sigma$ and zero everywhere else. Then $\Gamma\in C^q(K^j)$ can be written in the form $\Gamma=\sigma_1^\ast+\cdots+\sigma_n^\ast$ for some $n$.

We call $\sigma$ a {\it coface}\index{coface} of $\tau$ if $\tau$ is a face of $\sigma$. For the computations we will need the {\it q-dimensional coincidence matrix}.\index{coincidence matrix} which we will now define.

\begin{definition}
For a simplicial complex $K$, the {\it q-dimensional coincidence matrix} $CI_q$ is whose rows correspond to the ordered $q$-dimensional simplices and whose columns correspond to the ordered $q-1$-dimensional simplices. If $\sigma^i$ is a  $q$-dimensional simplex, and $\tau_j$ is a $q-1$-dimensional simplex, then $CI_q[i, j]=1$ if $\sigma_i$ is a coface of $\tau_j$ and  $CI_q[i, j]=0$ otherwise.
\end{definition}

The algorithm does a reduction of the {\it total coincidence matrix} $CI$ whose rows and columns are indexed by all simplices. (Not just those of dimension $q$ or $q-1$.) $CI[i, j]=1$ if the simplex $i$ is a coface of $j$ and $j$ represents a simplex of a dimension one lower than $i$. The algorithms in \cite{HB} give generators for $H^p(K)$ at every level of the filtration and determines when these generators die as we move backwards through it. See there for details and numerical examples.

HB's thesis also looks at cup products. In $Z_2$ it is easy to compute cup products. The cup product of two cohomology classes is the pairwise concatenation of the simplices on which the first factor evaluates to one and those on which the second factor does. 

As an example, let $\alpha=[1, 2]^\ast+[1, 4]^\ast+[2, 4]^\ast$ and $\beta=[2, 3]^\ast+[3, 4]^\ast+[2, 4]^\ast$. The only concatenations preserving vertex ordering are $[1, 2, 3]$ and $[1, 2, 4]$. So $\alpha\cup\beta([1, 2, 3])=\alpha([1,2])\cdot\beta([2, 3])=1\cdot 1=1.$ By contrast $\alpha\cup\beta([1, 3, 4])=\alpha([1, 3])\cdot\beta([3, 4])=0\cdot 1=0.$

The focus in \cite{HB} is {\it decomposability}

\begin{definition}
Suppose we have a cohomology class $\beta\in H^{p+q}(K)$. $\beta$ is {\it decomposable} if there exists an $x\in H^p(K)$ and a $y\in H^q(K)$ such that $\beta=x\cup y.$
\end{definition}

The following is the algorithm in \cite{HB} for determining if $\beta$ is decomposable as a cup product with a given $x$.

Let $K$ be a simplicial complex for which we have computed the cohomology groups. Let $[x]$ be a $p$-dimensional class that is born at $K_s$ with $s>n$ and dies at $K_t$ with $t<n-1$. Let $x=\sigma_1^\ast+\cdots+\sigma_n^\ast$ and let $K_n$ be $K_{n-1}$ with a $(p+q)$-dimesnional simplex $\Delta$ added. Also suppose that $\Delta^\ast=\delta(y)$ at the level of the filtration $K_n$. By removing $\Delta$ and moving backwards to $K_{n-1}$, $\delta(y)$ becomes 0, and $y$ is born as a cocycle.

When $\Delta$ is removed, we have $[y]\in H^{p+q-1}(K_{n-1})$. Let $y=\tau_1^\ast+\cdots+\tau_m^\ast.$ Consider the coincidence matrix whose columns are indexed by $\tau_1, \cdots, \tau_m$ and whose rows are codimension one cofaces of the $\tau_j$'s. (In other words, the rows are simplices of dimension one more than those of the columns.) If $\lambda_i$ represents row $i$, then the $i, j$-th entry of the matrix is equal to one if and only if $\tau_j$ is a face of $\lambda_i$. For each $\lambda_i$, if there is a $\sigma_k$ in $x$ such that $\sigma_k$ is a front face of $\lambda_i$, we store $\lambda_i$ in a list denoted $\beta$. After looping through the rows, if the size of $\beta$ is even, then $\beta $ is not a cup product of anything with $x$. (Remember that we are working in $Z_2$.) If the size of the list is odd, then $\beta$ is born as a cup product with $x$. So the birth of this cocycle gives rise to the birth of a cup product.

See \cite{HB} for an explicit numerical example.

\section{Ripser}

As I mentioned earlier, Ripser, which was developed by Ulrich Bauer of the Technical University of Munich is the fastest software I am aware of for computing persistent homology barcodes. The software is available at \cite{Bau1}. A full description of its advantages over earlier algorithms can be found in \cite{Bau2}. Bauer has a nice summary of this paper in a talk he gave at the Alan Turing Institute in 2017 \cite{Bau3}. The material in this section will come from there.

Bauer uses the following example. Take 192 points dispersed over $S^2$. Compute homology barcodes up to dimension 2. There are over 56 million simplices in the 3-skeleton. As an example, javaplex took 53 minutes and 20 seconds and used 12 GB of memory. Dionysus took 8 minutes and 53 seconds and used 3.4 GB. Ripser takes 1.2 seconds and uses 152 MB of memory. 

The improved performance is based on four optimizations:\begin{enumerate}
\item Clearing inessential columns. (I will explain what this means below.)
\item Computing cohomology rather than homology.
\item Implicit matrix reduction. 
\item Apparent and emergent pairs.
\end{enumerate}

It turns out that cohomology yields a considerable speedup but only in conjunction with clearing. 

Consider the standard matrix algorithm. Suppose we have a cloud of points in a finite metric space and we want to compute the perisistent homology over $Z_2$. Start with the {\it filtration boundary matrix} $D$ where the columns and rows are indexed by simplices in the filtration and $D_{i,j}=1$ if $\sigma_i\in\partial\sigma_j$, and $D_{i,j}=0$ otherwise. (See for example \cite{CEM}.)  For a matrix $R$, let $R_i$ be the $i$th column of $R$. The {\it pivot index} of $R_i$ denoted by Pivot $R_i$ is the largest row index of any nonzero entry or 0 if all column entries are 0. Define Pivots $R=\cup_i$ Pivot $R_i\setminus \{0\}$.

A column $R_i$ is {\it reduced} if Pivot $R_i$ cannot be decreased using column additions by scalar multiples of columns $R_j$ with $j<i$. In particular, $R_i$ is reduced of $R_i=0$ or all columns $R_j$ with $j<i$ are reduced and satisfy Pivot $R_j\neq$ Pivot $R_i$. The entire matrix $R$ is {\it reduced} if all of its columns are reduced.

The following proposition is the basis for matrix reduction algorithms to compute persistence homology. 

\begin{theorem}
\cite{CEM}. Let $D$ be a filtration boundary matrix and let $V$ be a full rank upper triangular matrix such that $R=DV$ is reduced. Then the index peristence (i.e. birth-death) pairs are $$\{(i,j)| i=\mbox{Pivot }R_j\neq 0\},$$ and the essential indices (birth times of classes that live forever) are $$\{i | R_i=0, i\notin\mbox{Pivots }R\}.$$
\end{theorem}

\begin{theorem}
The persistent homology is generated by the representative cycles $$\{R_j| j\mbox{ is a death index}\}\cup\{V_i| i\mbox{ is an essential index}\},$$ in the sense that for all filtration indices $k$, $H_\ast(K_k)$ is generated by the cycles $$\{R_j|(i, j)\mbox{ is an index persistence pair with }k\in[i, j)\}\cup\{V_i| i\mbox{ is an essential index with }k\in[i, \infty)\},$$ and for all pairs of indices $k, m$, the image of the map in homology $H_\ast(K_k)\rightarrow H_\ast(K_m)$ induced by inclusion has a basis generated by the cycles $$\{R_j|(i, j)\mbox{ is an index persistence pair with }k, m\in[i, j)\}\cup\{V_i| i\mbox{ is an essential index with }k\in[i, \infty)\}.$$
\end{theorem}

This leads to the first optimization. We don't care about a column with a non-essential birth index. So we set $R_i=0$ if $i=$Pivot $R_j$ for some $j\neq i$. As we still need to compute $V$, we can simply set $V_i=R_j$ as this column is a boundary. Then we still get $R=DV$ reduced and $V$ full rank upper triangular. This is a big savings as reducing birth columns are typically harder than reducing death columns, and we only need to reduce essential birth columns. 

Now for cohomology, we reduce $D^T$ instead of $D$. In low dimensions it turns out that the number of skipped columns is much greater when we use cohomology. The persistence barcodes are easily converted from cohomology to homology by the universal coefficient theorem as we are working over a field.

Returning to our example of 192 points on a sphere, we need to reduce 56,050,096 columns. Using clearing, this reduces to 54,888,816. With cohomology and clearing, the number reduces to 1,161,472 which is less than 1/60 of the work. 

The next optimization is implicit matrix reduction. Ripser only stores the matrix $V$ which is much sparser than $R=DV$. Recompute $D$ and $R$ only when needed. We only need the current column $R_j$ and the pivots of the previous columns to obtain the persistence intervals $[i, j)$. 

The last optimization of apparent and emergent pairs involves ideas from discrete Morse theory where simplices can be collapsed onto faces. See \cite{Bau2} and its references for the details which I will omit in the interest of space.

In conclusion, topological data analysis so far only uses a small fraction of the machinery of algebraic topology. The recent development of the Ripser algorithm demonstrated that cohomology is a very valuable tool that had not previously been explored. Are there more tools from algebraic topology that would be useful in data science? It is hard for me to imagine that this is not the case. In the remainder of this book, I would like to explore this possibility.

\chapter{Homotopy Theory}

Homotopy theory is the branch of algebraic topology dealing with the classification of maps of a sphere into a space. It grew up alongside homology theory and is closely related. Like homology, it is a functor from topological spaces and continuous maps to groups and homomorphisms. If in some dimension, the groups are not isomorphic, the spaces are not homeomorphic or even homotopy equivalent. Also, the theory is somewhat easier to picture for homotopy as opposed to homology.

There are also some major differences. Homology groups as we have seen are computable for a finite simplicial complex, while results for homotopy groups apply under very special cases. While spheres have the simplest homology groups, we don't even know all of the homotopy groups of $S^2$. (And we never will.) In fact homotopy groups of spheres is a special topic with an extensive and very difficult literature. I will briefly touch on it in Chapter 12. The simplest spaces in homotopy which have only one nonzero homotopy group are the Eilenberg-Mac Lane spaces. With the exception of the circle $S^1$, Eilenberg-Mac Lane spaces are generally very complicated and usually infinite dimensional. Sometimes, though, homotopy groups are eaiser to calculate. For example, homotopy groups of product spaces are just direct sums as opposed to homology where we need the more complicated K{\"u}nneth Theorem. And for fiber bundles, a twisted version of a product (eg. a M{\" o}bius strip), homotopy produces a long exact sequence, while for homology, you need to use a spectral sequence, one of the scariest mathematical tools ever dreamed up. I will try to give you a taste of what spectral sequences are at the end of this chapter.

So why am I putting homotopy in a book about data science? The answer is that it is a vital component of {\it obstruction theory}. Suppose $L$ is a subcomplex of a simplicial complex $K$. We have a continuous map $f: L\rightarrow Y$ and we would like to extend it to a map $g: K\rightarrow Y$ such that $g$ is continuous and $g| L=f$. The idea is to extend it over the vertices, edges, triangles, etc. Suppose we have extended it to the $k$-simplices. For the next step, we know that the function has been extended to the boundary of each $(k+1)$-siimplex and we want to extend it to the interior. This turns out to be equivalent to asking if a certain cocycle on $K$ is equal to zero. I will explain the details in Chapter 10.

For now, though, recall that in a filtration, the complex at an earlier filtration index is a subcomplex of the complex at a later one.. For a point cloud, we may ask if a function on the Vietoris-Rps complex at radius $\epsilon_1$ can be extended to a function on the complex at radius $\epsilon_2$, where $\epsilon_2>\epsilon_1$. If not, where in the cloud can it be done, and where are the obstructions. Also, what kind of function and what target space $Y$ would be useful in appliications?

In this chapter, I will restrict myself to giving the basic definitions and results that you need to understand what comes afterwards. I will leave out most of the proofs and refer you to the standard textbooks for more details.

The first book I will recommend is Steenrod's {\em Topology of Fibre Bundles} \cite{Ste3}, the earliest homotopy textbook I know of. (Note: Although "fibre" is often used in the literature, I will use the American spelling "fiber" unless it is in a direct quote.) It is dense and quite dated but very clearly written. My personal favorite is the book by Hu \cite{Hu}. It is also very clear but also a little old. The chapter on spectral sequences, though, is missing some key pieces so I would recommend looking elsewhere on that topic. (See below for some suggestions.) Another classic is the extremely thick and dense book by Whitehead \cite{Whi}. Finally, there is a section on homotopy in Hatcher \cite{Hat}.

Unless otherwise stated, the material in Sections 9.1-9.7 come from Hu \cite{Hu}. Section 9.1 will describe the extension problem, a major theme for the remainder of this book and its dual, the lifting problem. In Section 9.2, I will start with the homotopy group in dimension 1, also known as the fundamental group. Section 9.3 will discuss fiber bundles, and Section 9.4 will discuss an interesting nontrivial example, the Hopf maps. Another example is the spaces of paths and loops which are the subject of Section 9.5. Higher dimensional homotopy groups are the subject of Section 9.6 and Section 9.7 gives some examples of computational tricks. Finally, I will discuss the computational tools of Postnikov systems in Section 9.8 and spectral sequences in Section 9.9.

\section{The Extension and Lifting Problems}

\begin{definition}
Let $f: X\rightarrow Y$ be a continuous map between two topological spaces (generally, the term map will always mean a continuous function), and let $A$ be a subspace of $X$. Then the map $g: A\rightarrow Y$ such that $g(x)=f(x)$ is called the {\it restriction} of $f$ to $A$ and $f$ is called an extension of $g$.
\end{definition}

If $h: A\rightarrow X$ is inclusion, then the {\it extension problem}\index{extension problem} is to determine whether or not a given map $g: A\rightarrow Y$ has an extension over $X$. Equivalently, we want to find $f: X\rightarrow Y$ which makes the following diagram commute (i.e, $g=fh$): 
$$\begin{tikzpicture}
  \matrix (m) [matrix of math nodes,row sep=3em,column sep=4em,minimum width=2em]
  {
A &  & Y\\
 & X & \\};
  \path[-stealth]
(m-1-1) edge node [above] {$g$} (m-1-3)
(m-1-1) edge node [below] {$h$} (m-2-2)
(m-2-2) edge[dashed] node [below] {$f$} (m-1-3)

;

\end{tikzpicture}$$

Using the functoriality properties of homology, if such an extension exists, we can solve find $f_\ast$ so the following diagram commutes:
$$\begin{tikzpicture}
  \matrix (m) [matrix of math nodes,row sep=3em,column sep=4em,minimum width=2em]
  {
H_m(A) &  & H_m(Y)\\
 & H_m(X) & \\};
  \path[-stealth]
(m-1-1) edge node [above] {$g_\ast$} (m-1-3)
(m-1-1) edge node [below] {$h_\ast$} (m-2-2)
(m-2-2) edge[dashed] node [below] {$f_\ast$} (m-1-3)

;

\end{tikzpicture}$$ for each $m$.

Recall that we have looked at this problem for $A=Y=S^n$, $n>0$, and $X=B^{n+1}$. If $g$ is a constant map, let $g(a)=y_0$ for all $a\in A$ and fixed $y_0\in Y$. We can easily extend to all of $X$ by letting $f(x)=y_0$. 

But suppose $g: S^n\rightarrow S^n$ is the identity map. Then $H_n(S^n)\cong Z$ and $g_\ast$ is the identity on $Z$. But $H_n(B^{n+1})=0$ so $h_\ast=0$. Thus, it is impossible to find $f_\ast$ such that $f_\ast h_\ast=g_\ast$ and the identity can not be extended over the entire ball. 

Now the homology problem has a solution even if $g\simeq fh$ meaning that $g$ is homotopic to $fh$ rather than strictly equal to it. We will use this more general problem in what follows. Hu lists a number of related problems. One interesting one is the dual of the extension problem. 

Let $A, X$ and $Y$ be topological spaces. We don't assume anything about inclusion but instead let $f: X\rightarrow A$ be surjective. Given a map $g: Y\rightarrow A$, we would like to find a map $f: Y\rightarrow X$ such that $g=hf$. This is called the {\it lifting problem}.\index{lifting problem} Here is a diagram repositioned to illustrate the name "lifting" but with the spaces and maps having the same names as before. You can see that all of the arrows have been reversed.

$$\begin{tikzpicture}
  \matrix (m) [matrix of math nodes,row sep=3em,column sep=4em,minimum width=2em]
  {
 &  & X\\
Y & & A\\};

  \path[-stealth]
(m-2-1) edge node [below] {$g$} (m-2-3)
(m-1-3) edge node [right] {$h$} (m-2-3)
(m-2-1) edge[dashed] node [above] {$f$} (m-1-3)

;

\end{tikzpicture}$$

Like the extension problem, the lifting problem is also usually generalized to the case of finding a map $f$ such that $g\simeq hf.$

My dissertation \cite{Pos1} actually involved solving a lifting problem involving some strange infinite dimensional function spaces.

\section{The Fundamental Group}
In this section, we will discuss the properties of the one-dimensional homotopy group, also known as the {\it fundamental group}. It classifies maps between $S^1$ and a space $X$. 

First, though, we will look at a simple example of what we will call a fiber bundle in the next section. 

\begin{definition}
The {\it exponential map}\index{exponential map} is the map $p: R\rightarrow S^1$ defined by $$p(x)=e^{2\pi xi},$$ for $x\in R$. 
\end{definition}

The first result is an easy consequence of the definition.

\begin{theorem}
The exponential map $p$ is a homomorphism in the sense that $p(x+y)=p(x)p(y)$ for $x, y\in R$. (Note that multiplying in $S^1$ is actually adding angles.) The kernel $p^{-1}(1)$ of $p$ consists of the integers.
\end{theorem}

\begin{theorem}
For every proper connected subspace $U$ of $S^1$, $p$ takes every component of $p^{-1}(U)$ homeomorphically onto $U$.
\end{theorem}

For the next result, we need to define paths and loops. These will be crucial in what follows.

\begin{definition}
A {\it path}\index{path} in a space $X$ is a map $\sigma: I\rightarrow X$, where $I=[0,1]$. $\sigma(0)$ is called the {\it initial point} of $\sigma$ and $\sigma(1)$ is called the {\it terminal point}. If $\sigma(0)=\sigma(1)=x_0$, then $\sigma$ is called a {\it loop,}\index{loop} and $x_0$ is its {\it base point}.
\end{definition}

\begin{theorem}
{\bf The Covering Path Property:} For eavery path $\sigma: I\rightarrow S^1$ and every point $x_0\in R$ such that $p(x_0)=\sigma(0)$, there exists a unique path $\tau: I\rightarrow R$ such that $\tau(0)=x_0$ and $p\tau=\sigma$.
\end{theorem}

Important: In this subject, it is easy to get lost in all the spaces and maps between them. If you get confused. draw diagrams and make sure you keep track of which spaces are the domain and range of the maps involved.

In the previous result, if $\log$ denotes the natural logarithm, the path $\tau: I\rightarrow R$ is given by $\tau(0)=x_0$, and $$\tau(t)=\frac{1}{2\pi i}\log(\sigma(t))$$ for $t\in I$.

Recall that for $f, g:X\rightarrow Y$, $f$ is homotopic to $g$, denoted $f\simeq g$ if there is a continuous map $H: X\times I\rightarrow Y$ such that  $H(x, 0)=f$ and $H(x, 1)=g$. Following Hu, we will write $f_t: X\rightarrow Y$ and call it a {\it homotopy of f} if $f$ is specified but $g$ is not. 

\begin{theorem}
{\bf The Covering Homotopy Property:} For every map $f: X\rightarrow R$ of a space $X$ into the reals, and every homotopy $h_t: X\rightarrow S^1$, $(0\leq t\leq 1)$ of the map $h=pf$, there exists a unique homotopy $f_t: X\rightarrow R$ of $f$ such that $pf_t=h_t$ for every $t\in I$..
\end{theorem}

Now we define the {\it fundamental group} $\pi_1(X)$. 

Let $X$ be a space and $x_0\in X$. Let $\Omega$ be the set of all maps $f: S^1\rightarrow X$ such that $f(1)=x_0$. If we replace the map $f$ with the map $fp: I\rightarrow X$ where $p$ is the exponential map, we have a loop in $X$ with $x_0$ as its base point. The correspondence is one to one, so identifying $f$ with $fp$, we can consider $\Omega$ to be the set of loops in $X$ with base point $x_0$. 

Letting $$\Omega=\{f: I\rightarrow X| f(0)=x_0=f(1)\},$$ we can define a multiplication in $\Omega$. If $f, g\in \Omega$, the product consists of travelling around loop $f$ and then loop $g$. To be precise, we define the product $f\cdot g$ by $$(f\cdot g)(t)=
\left
\{
\begin{array}{ll}
f(2t) & \mbox{if }0\leq t\leq \frac{1}{2},\\
g(2t-1) & \mbox{if }\frac{1}{2}\leq t\leq 1.
\end{array}
\right.
$$

Now homotopy of maps is an equivalence relation, so we denote the homotopy class of $f$ by $[f]$ and call $f$ a {\it representative} of the class $[f]$. 

Let $d$ denote the degenerate loop $d(I)=x_0$. For a loop $f\in\Omega$, let $f^{-1}$  be the loop $f$ traversed in the opposite direction so that $f^{-1}(t)=f(1-t)$. Then $(f^{-1})^{-1}=f$.

\begin{theorem}
If $f, g, h\in \Omega$ and $d, f^{-1}$ are as defined above, the following hold:\begin{enumerate}
\item If $f\simeq f'$ and $g\simeq g'$, then $f\cdot g\simeq f'\cdot g'$.
\item $(f\cdot g)\cdot  h\simeq f\cdot (g\cdot h).$
\item $f\cdot d\simeq d\cdot f\simeq f$.
\item $f\cdot f^{-1}\simeq f^{-1}\cdot f\simeq d$.
\end{enumerate}
\end{theorem}

By property 1, $[f][g]=[f\cdot g]$. So Theorem 9.2.5 implies that the homotopy classes form a group with identity element $e=[d]$ and the inverse of $[f]$ is $[f^{-1}]$.

\begin{definition}
The homotopy classes of loops $f: I\rightarrow X$ with $f(0)=f(1)=x_0$ form a group with the multiplication defined above called the {\it fundamental group}\index{fundamental group} of $X$ at $x_0$ denoted $\pi_1(X, x_0)$.
\end{definition}

As an example, let $X=S^1$. Remember from Section 4.1.7 that the degree of a map $f: S^n\rightarrow S^n$ determines the homotopy class of $f$. Now $H_n(S^n)\cong Z$, and we defined the degree $deg(f)$ by $f_\ast(z)=deg(f)z$. Theorem 4.1.29 implies that the map $[f]\rightarrow deg(f)$ is a well defined homomorphism from $\pi_1(S^1, 1)$ to $Z$. Hu proves that a map $f: S^1\rightarrow S^1$ of degree $n$ is homotopic to the map $\phi_n: S^1\rightarrow S^1$ defined by $\phi_n(z)=z^n$ and that $\phi_n$ has degree $n$. So the homomorphism is one-to-one and onto and hence $\pi_1(S^1)\cong Z$.

Note that the fundamental group is not computable in general. We can, however, write down generators and relations for the fundamental group of a finite simplicial complex.

Also, unlike homology groups, the fundamental group may not be abelian. For example, let $E$ be a figure 8. Let $x_0$ be the point where the two loops meet, $f$ denote the top loop and $g$ denote the bottom loop. Then $\pi_1(E, x_0)$ is the (non-abelian) free group generated by $f$ and $g$. 

What happens if we pick a different base point $x_1$? If there is a path $\sigma$ from $x_1$ to $x_0$, then if $f$ is a loop with base point $x_1$, and $\tau$ is the reverse of $\sigma$, then $\sigma\cdot f\cdot\tau$ is a loop with base point $x_0$. So the map $[f]\rightarrow [\sigma\cdot f\cdot\tau]$ is an isomorphism $\pi_1(X, x_1)\cong\pi_1(X, x_0)$.

If $X$ is pathwise connected, we can dispense with the base point and simply write $\pi_1(X)$.

Finally, if $\phi: X\rightarrow Y$ and $f: I\rightarrow X$ is a loop on $X$, then $\phi f: I\rightarrow Y$ is a loop on $Y$. So $\phi$ gives rise to a homomorphism $\phi_\ast: \pi_1(X, x_0)\rightarrow \pi_1(Y, y_0)$ for $y_0=\phi(x_0)$ where $\phi_\ast([f])=[\phi f].$

Now we would like to know when $\pi_1(X)=0$ for a pathwise connected space $X$.

\begin{definition}
A pathwise connected space $X$ is {\it simply connected}\index{simply connected} if every pair of paths $\sigma, \tau: I\rightarrow X$ with $\sigma(0)=\tau(0)$ and $\sigma(1)=\tau(1)$ are homotopic by a homotopy that fixes the endpoints. 
\end{definition}

If $f$ is a loop in a simply connected space $X$, and $d$ is the degenerate loop at $f(0)=f(1)=x_0$ then [f]=[d]. So $\pi_1(X)=0$. Hu proves the stronger statement that the two conditions are equivalent.

\begin{theorem}
If $X$ is a nonempty pathwise connected space, then $X$ is simply connected if and only if $\pi_1(X)=0$.
\end{theorem}

\begin{theorem}
$S^n$ is simply connected if and only if $n>1$.
\end{theorem}

{\bf Proof:} $S^0$ consists of two points and is not even pathwise connected, so it is not simply connected. $S^1$ is also not simply connected by the previous theorem since $\pi_1(S^1)\cong Z\neq 0.$

For $n>1$ let $f: I\rightarrow S^n$ be a loop on $S^n$ with base point $x_0$. Since $n>1$, $f(I)$ is a proper subspace of $S^n$ and thus contractible to the point $x_0$. So $[f]=[d]$ and $\pi_1(S^n)=0.$ $\blacksquare$

Since the fundamental group can be non-abelian, how does it compare with the one-dimensional homology group? For a pathwise connected space $X$, there is a simple relationship. 

Let $$h_\ast:\pi_1(X, x_0)\rightarrow H_1(X)$$ be defined as follows. Let $\alpha\in \pi_1(X, x_0).$ choose a representative loop for $\alpha$ called $f: S^1\rightarrow X$ with $f(1)=x_0$. Then $f$ induces a homomorphism $f_\ast: H_1(S^1)\rightarrow H_1(X)$. Let $\iota$ be the generator for the infinite cyclic group $H_1(S^1)$ corresponding to the counterclockwise orientation of $S^1$. Then $f_\ast(\iota)\in H_1(X)$ depends only on the class $\alpha$ so we define $h_\ast$ by $h_\ast(\alpha)=f_\ast(\iota).$ (We will see several variants of this construction in future discussions.) 

Now we need a definition from abstract algebra.

\begin{definition}
Let $G$ be a group and let $g, h\in G$. Then the {\it commutator}\index{commutator} $[g, h]$ of $g$ and $h$ is the element $ghg^{-1}h^{-1}$. The set of commutators is not closed under multiplication, but the subgroup they generate is called the {\it commutator subgroup}.\index{commutator subgroup}
\end{definition}

If $G$ is abelian, any commutator reduces to the identity. If $Comm(G)$ is the commutator subgroup of $G$, then the quotient group $G/Comm(G)$ is abelian and called the {\it abelianization}\index{abelianization} of $G$.

Since $h_\ast$ is a homomorphism and $H_1(X)$ is abelian, $Comm(\pi_1(X, x_0))$ must be contained in the kernel of $h_\ast$. We actually have the following:

\begin{theorem}
If $X$ is pathwise connected then the natural homomorphism $h_\ast$ maps $\pi_1(X, x_0)$ onto $H_1(X)$ with the commutator subgroup $Comm(\pi_1(X, x_0))$ as its kernel. So $H_1(X)$ is isomorphic to the abelianization of $\pi_1(X, x_0)$.
\end{theorem}

Theorem 9.2.8 is the Hurewicz Theorem for dimension one. The Theorem has a nicer form for higher dimensions which we will see in Section 9.7.

\section{Fiber Bundles}

So what is a fiber bundle anyway? The usual definition is a space that is "locally" a cartesian product of two spaces. Another definition is a map that has the {\it covering homotopy property}. 

One thing that will help us a lot is that we are working with finte simplicial complexes in any data science application. The definitions become equivalent if the spaces we are working with are {\it paracompact}.

\begin{definition}
A space $X$ is {\it paracompact}\index{paracompact} if given every open cover $\{U_\alpha\}$ of $X$, has a locally finite refinement. In other words, there is an open cover $\{V_\beta\}$ such that every set $V_\beta$ in the cover is contained in some $U_\alpha$ in the original cover and every point $x\in X$ has a neighborhood which intersects only finitely many of the sets in $\{V_\beta\}$.
\end{definition}

A compact space is easily seen to be paracompact. Any CW complex is paracompact and any finite CW complex is actually compact. Finite simplicial complexes are also compact and those are the only kind of spaces we will deal with in applications.

We will start with the simpler definition found in Hu \cite{Hu}. Then we will compare it to the more elaborate version in Steenrod \cite{Ste3}. Be careful of similar sounding words that may or may not mean the same thing. And again, draw yourself lots of pictures. There are way too many symbols to keep track of otherwise.

Hu starts his discussion with various covering homotopy properties. There are really only three that will concern us: The CHP, the ACHP, and the PCHP. For the definitions refer to the following diagram: 

$$\begin{tikzpicture}
  \matrix (m) [matrix of math nodes,row sep=3em,column sep=4em,minimum width=2em]
  {
 &  & E\\
X & & B\\};

  \path[-stealth]
(m-2-1) edge node [below] {$f$} (m-2-3)
(m-1-3) edge node [right] {$p$} (m-2-3)
(m-2-1) edge node [above] {$f^\ast$} (m-1-3)

;

\end{tikzpicture}$$

\begin{definition}
Let $X$ be a given space and $f: X\rightarrow B$ be a given map, and let $f_t: X\rightarrow B$ for $t\in [0,1]$ be a given homotopy of $f$. A map $f^\ast: X\rightarrow E$ is said to {\it cover} $f$ relative to $p$ if $pf^\ast=f$. A homotopy $f_t^\ast: X\rightarrow E$ for $t\in [0,1]$ of $f^\ast$ is said to cover the homotopy $f_t$ relative to $p$ if $pf_t^\ast=f_t$ for all $t$. Then $f_t^\ast$ is called a {\it covering homotopy} of $f_t$.
\end{definition}

\begin{definition}
The map $p: E\rightarrow B$ is said to have the {\it covering homotopy property}\index{covering homotopy property} or {\it CHP}\index{CHP} for the space $X$ if for every map $f^\ast: X\rightarrow E$ and every homotopy $f_t: X\rightarrow B$ of the map $f=pf^\ast: X\rightarrow B$, there exists a homotopy $f^\ast_t: X\rightarrow E$ of $f^\ast$ which covers the homotopy $f_t$. The map $p: E\rightarrow B$ has the {\it absolute covering homotopy property}\index{absolute covering homotopy property} or {\it ACHP}\index{ACHP} if it has the CHP for every space $X$.The map $p: E\rightarrow B$ has the {\it polyhedral covering homotopy property}\index{polyhedral covering homotopy property} or {\it PCHP}\index{PCHP} if it has the CHP for every {\it triangulable} space $X$.
\end{definition}

\begin{example}
Let $E=B\times D$, and let $p: E\rightarrow B$ be the natural projection. Let $X$ be a given space, $f^\ast: X\rightarrow E$ be a given map, and let $f_t: X\rightarrow B$ a given homotopy of the map $f=pf^\ast$. Then $f_t$ has a covering homotopy $f_t^\ast: X\rightarrow E$ of $f^\ast$ defined by $$f_t^\ast(x)=(f_t(x), qf^\ast(x))$$ for every $x\in X$, $t\in [0,1]$ where $q:E\rightarrow D$ is the other projection. So $p$ has the ACHP.
\end{example}

\begin{definition}
A map $p: E\rightarrow B$ is a {\it fibering}\index{fibering} if it has the PCHP. In this case $E$ is called a {\it fiber space over the base space B with projection} $p:E\rightarrow B$. ($E$ is also called the {\it total space}.) For each point $b\in B$, the subspace $p^{-1}(b)$  of $E$ is called the fiber over $b$.
\end{definition}

For all of our applications a fibering will be equivalent to what Hu calls a {\it bundle space}.

\begin{definition}
A map $p: E\rightarrow B$ has the bundle property if there exists a space $D$ such that for each $b\in B$, there is an open neighborhood $U$ of $b$ in $B$ together with a homeomorphism $$\phi_U: U\times D\rightarrow p^{-1}(U)$$ of $U\times D$ onto $p^{-1}(U)$ satisfying the condition $$p\phi_U(u, d)=u$$ for $u\in U$, $d\in D$. Then $E$ is called the {\it bundle space}\index{bundle space} over the base space $B$ relative to the projection $p$. The space $D$ is called a {\it director space}. The open sets $U$ and the homeomorphisms $\phi_U$ are called the {\it decomposing neighborhoods} and {\it decomposing functions} respectively.
\end{definition}

A bundle space can be thought of as a space that is locally a cartesian product but may have a twist. In this way a mobius strip differs from a cylinder, the latter of which is $S^1\times I$. 

\begin{definition}
Let $E=B\times D$ and $p: E\rightarrow B$ be projection. Then $E$ is a bundle space over $B$. (Just let the neighborhood $U$ be all of $B$.) Then $E$ is called a {\it product bundle}\index{product bundle} or {\it trivial bundle}\index{trivial bundle} over $B$.
\end{definition}

Now you can finally understand the pun from my introduction: What do you get when you cross an elephant and an ant? The trivial elephant bundle over the ant.

The following two results give the connection between a bundle space and the fiber space we defined first. The idea is to use the fact that the projection of a  product space onto one of its factors has the ACHP and thus the PCHP. 

\begin{theorem}
Every bundle space $E$ over $B$ relative to $p: E\rightarrow B$ is a fiber space over $B$ relative to $p$.
\end{theorem}

\begin{theorem}
If a map $p: E\rightarrow B$ has the bundle property, then it has the CHP for every paracompact Hausdorff space and thus any finite CW complex or simplicial complex.
\end{theorem}

Steenrod \cite{Ste3} calls the bundle space we just defined an {\it Ehresmann-Feldbau bundle} or {\it E-F bundle.} It is a special case of what he calls a {\it coordinate bundle} which has a lot more structure. I will now give his definition of a coordinate bundle and the related term {\it fiber bundle.} We first need the notion of a {\it topological group}.

\begin{definition}
A {\it topological group}\index{topological group} $G$ is a set which has a group structure and a topology such that:\begin{enumerate}
\item $g^{-1}$ is continuous for $g\in G$.
\item The map $\mu: G\times G\rightarrow G$ given by $\mu(g_1, g_2)=g_1g_2$ is continuous.
\end{enumerate}
\end{definition}

\begin{definition}
If $G$ is a topological group and $Y$ is a topological space, then we say that $G$ is a {\it topological transformation group}\index{topological transformation group} of $Y$ relative to a map $\eta: G\times Y\rightarrow Y$ if: \begin{enumerate}
\item $\eta$ is continuous.
\item $\eta(e, y)=y$ if $e$ is the identity of $G$.
\item $\eta(g_1g_2, y)=\eta(g_1. \eta(g_2, y))$ for $g_1, g_2\in G$ and $y\in Y$.
\end{enumerate}
We say that $G$ {\it acts on Y}.\index{group acting on a set.} If $Y$ is a set and not necessarily a topological space, a group $G$ (not necessarily topological) acts on $Y$ if statements 2 and 3 hold.
\end{definition}

We usually leave out $\eta$ and write $\eta(g, y)$ as $gy$.

For a fixed $g$, the map $y\rightarrow gy$ is a homeomorphism of $Y$ onto itself since it has a continuous inverse $y\rightarrow g^{-1}y$. So $\eta$ induces a homomorphism from $G$ to the group of homeomorphisms of $Y$. $G$ is {\it effective} if $gy=y$ for all $y\in Y$ implies that $g=e$. Then $G$ is isomorphic to a group of homeomorphisms of $Y$.

In our definition, we will assume that $G$ is effective. I will now define a {\it coordinate bundle}. There are lots of maps here. Draw yourself pictures or I guarantee that you will be totally confused. Also, I will use $E$ for the bundle space and $B$ for the base space which is more conventional, but note that Steenrod uses $B$ for the bundle space and $X$ for the base space.

\begin{definition}
A {\it coordinate bundle}\index{coordinate bundle} $\mathfrak{B}$ consists of the following:\begin{enumerate}
\item A space $E$ called the {\it bundle space}.
\item A space $B$ called the {\it base space}.
\item A map $p: E\rightarrow B$ of $E$ onto $B$ called the {\it projection}.
\item A space $Y$ called the {\it fiber}.
\item An effective topological transformation group $G$ of $Y$ called the {\it group of the bundle}.
\item A family $\{V_j\}$ of open sets covering $B$. These open sets are called {\it coordinate neighborhoods}.
\item For each coordinate neighborhood $V_j\subset B$, we have a homeomorphism $\phi_j: V_j\times Y\rightarrow p^{-1}(V_j)$ called the coordinate function. These functions need to satisfy conditions 8-10.
\item $p\phi_j(b, y)=b.$
\item If the map $\phi_{j, b}: Y\rightarrow p^{-1}(b)$ is defined by setting $$\phi_{j, b}(y)=\phi_j(b, y),$$ then for each pair of indices $i, j$, and each $b\in V_i\cap V_j$, the homeomorphism $$\phi^{-1}_{j, b}\phi_{i, b}: Y \rightarrow Y$$ coincides with the operation of an element of $G$ which is unique since $G$ is effective.
\item For each pair $i, j$ of indices, the map $$g_{ji}: V_i\cap V_j\rightarrow G$$ defined by $$g_{ji}(b)=\phi^{-1}_{j, b}\phi_{i, b}$$ is continuous. These maps are called the {\it coordinate transformations} of the bundle.
\end{enumerate}
\end{definition} 

\begin{definition}
In a coordinate bundle, we denote $p^{-1}(b)=Y_b$ and call it the {\it fiber over b}. 
\end{definition}

\begin{definition}
Two coordinate bundles are equivalent if they have the same bundle space, base space, projection, fiber and group, and their coordinate functions $\{\phi_j\}, \{\phi_k'\}$ satisfy the conditions that $$\overline{g}_{kj}(b)=\phi'^{-1}_{k, b}\phi_{j, b},$$ for $b\in V_j\cap V'_k$ coincides with the operation of an element of $G$ and the map $$\overline{g}_{kj}(b)\rightarrow G$$ is continuous. The definition defines an equivalence relation. An equivalence class of coordinate bundles is called a {\it fiber bundle}\index{fiber bundle}.
\end{definition}

In the bundle space defined in  Hu (our Defintion 9.3.5) the director space is the same as the fiber of a coordinate bundle. In this case, the fibers over each point in the base space are homeomorphic. (See \cite{Ste3} Section 6.) 

For the rest of this chapter and the book in general, I will use the terms fiber space, bundle space, and fiber bundle in accordance with my sources, but remember that in the cases we care about, fiber space and bundle space are pretty much equivalent, and fiber bundles differ with the inclusion of a group of transformations. 

The rest of this section is taken from \cite{Hu}.

\begin{definition} If $E$ is a bundle space over $B$ with projection $p: E\rightarrow B$, a {\it cross section}\index{cross section} in $E$ over a subspace $X$ of $B$ is a map $f: X\rightarrow E$ such that $pf(x)=x$ for $x\in X$.
\end{definition}

If $U\subset B$ is a decomposing neighborhood with decomposing function $\phi_U: U\times D\rightarrow p^{-1}(U)$ and projection $\psi_U: U\times D\rightarrow D$ then for any point $e\in p^{-1}(U)$, there is a cross section $f_e: U\rightarrow E$ given by $f_e(u)=\phi_U(u, d)$ where $d=\psi_U\phi_U^{-1}(e)$ for each $u\in U$. If $u=p(e)$, then $f_e(u)=e$. So in bundle spaces, local cross sections always exist.

For global cross sections in which $f$ is defined on all of $B$, we would have $pf=i_B$ so on homology $p_\ast f_\ast$ is the identity on $H_\ast(B)$. Then $f_\ast: H_m(B)\rightarrow H_m(E)$ is a monomorphism for all $m$, and  $p_\ast: H_m(E)\rightarrow H_m(B)$ is an epimorphism for all $m$. We will see that this condition fails for the Hopf maps from the next section. 

I will conclude this section with a discussion of a useful construction which produces a new fiber space from an old one by way a map of a given space into its base space.

\begin{definition}
Let $p: E\rightarrow B$ and $p': E'\rightarrow B'$ be to fiberings. A map $F: E\rightarrow E'$ is called a {\it fiber map}\index{fiber map} if it carries fibers into fibers, i.e. for every $b\in B$ there is a $b'\in B'$ such that $F$ carries $p^{-1}(b)$ into $p'^{-1}(b')$.  
\end{definition}

In this case, $F$ induces a function $f: B\rightarrow B'$ defined by $$f(b)=p'Fp^{-1}(b)$$ for every $b\in B$. Then $f$ is called the {induced map} of the fiber map $F$. If $E$ is a bundle space over $B$, then $p$ is open so $f$ is continuous. The following diagram is commutative:

$$\begin{tikzpicture}
  \matrix (m) [matrix of math nodes,row sep=3em,column sep=4em,minimum width=2em]
  {
E & E'\\
B & B'\\};

  \path[-stealth]
(m-1-1) edge node [above] {$F$} (m-1-2)
(m-1-1) edge node [left] {$p$} (m-2-1)
(m-1-2) edge node [right] {$p'$} (m-2-2)
(m-2-1) edge node [below] {$f$} (m-2-2)

;

\end{tikzpicture}$$

Suppose we are given a fibering $p': E'\rightarrow B'$ and a map $f: B\rightarrow B'$ of a given space $B'$ into $B$. We can construct a fibering $p: E\rightarrow B$ and fiber map $F: E\rightarrow E'$ that induces $f$. (Hu calls $F$ a {\it lifting} of $f$.) Let $E$ be the subspace of $B\times E'$ given by $$E=\{(b, e')\in B\times E'| f(b)=p'(e')\}$$ and let $p: E\rightarrow B$ denote the projection defined by $p(b, e')=b.$ Let $F: E\rightarrow E'$ be defined by $F(b, e')=e'$. Then by construction $fp=p'F$, so $F$ is a lifting of $f$. 

Hu uses the polyhedral covering homotopy property of $p'$ to show that it also holds for $p$ so $p$ is a fibering called the {\it fibering induced by f.}  

In the case where $B$ is a subspace of $B'$ and $f$ is inclusion, then $E$ can be identified with $p'^{-1}(B)$ and $p$ with $p'|B$. We then call the induced fibering $p: E\rightarrow B$ the {\it restriction} of $p': E'\rightarrow B'$ to $B$. This gives the following result.

\begin{theorem}
If $E$ is a fiber space over a base space $B$ with projection $p:E\rightarrow B$ and if $A$ is a subspace of $B$ then $p^{-1}(A)$ is a fiber space over $A$ with $(p|p^{-1}(A))$ as a projection. 
\end{theorem}

\section{The Hopf Maps}

There are three interesting fiber bundles (or fiber spaces if we ignore their transformation groups) called the {\it Hopf bundles} or {\it Hopf fibrations}. The corresponding projections will be called {\it Hopf maps}. These are the only fiber spaces in which the fiber, bundle space, and base space are all spheres. Proving this is very hard but it is related to the Hopf invariant problem which we will defer to Chapter 11.

The original paper describing these bundles was published by Hopf in 1935 \cite{Hop2}. The three maps are described in more or less detail in Steenrod \cite{Ste3}, Hu \cite{Hu}, and Hatcher \cite{Hat}. I will follow Hatcher's definitions. Then I will talk more about the first bundle which maps $S^3\rightarrow S^2$ following Lyons, who takes a different approach and gives a good intuitive description of where this map comes from. Lastly, I will state some results from Hu on maps which are{\it algebraically trivial}, i.e. trivial on homology for each dimension.

Hopefully it is not too much of a spoiler her to tell you that $\pi_n(X)$ for $n>1$ is the group of homotopy classes of maps from $S^n\rightarrow X$. (The group operation and more details will be explained in their place in section 9.6.) You may hope that like homology and cohomology, the homotopy groups $\pi_m(S^n)=0$ for $m>n.$ Then I could end this book right here. The Hopf maps were the bad news that this wasn't true. As the map from $S^3\rightarrow S^2$ will turn out not to be homotopic to a constant, we know that $\pi_3(S^2)\neq 0$. For obstruction theory applications to data science, this will be one of our simplest examples, and I will discuss the situation in detail in Section 10.2.

Writing our bundles in the form $F\rightarrow E\rightarrow B$ where $F$ is the fiber, $E$ is the bundle space, and $B$ is the base space, the three Hopf bundles are:\begin{enumerate}
\item $S^1\rightarrow S^3\rightarrow S^2.$
\item $S^3\rightarrow S^7\rightarrow S^4.$
\item $S^7\rightarrow S^{15}\rightarrow S^8.$
\end{enumerate}

I will describe these one at a time. The material comes from $\cite{Hat}$.  First, though consider the map $S^n\rightarrow P^n$ which identifies antipodal points on $S^n$. This is a fiber bundle called a {\it covering space} which means it has a discrete fiber. In this case, the fiber consists of two points, so we can think of it as $S^0$. 

Now, look at the complex analogue. We have the fiber space $S^1\rightarrow S^{2n+1}\rightarrow CP^n$. Then $ S^{2n+1}$ is the unit sphere in $C^{n+1}$ and $CP^n$ is the quotient space of $ S^{2n+1}$ under the equivalence relation $(z_0, \cdots, z_n)\sim\lambda(z_0, \cdots, z_n)$ for $\lambda\in S^1$, the unit circle in $C$. The projection $ S^{2n+1}\rightarrow CP^n$ sends $(z_0, \cdots, z_n)$ to its equivalence class $[z_0, \cdots, z_n]$. Since $\lambda$ varies over $S^1$, the fibers are copies of $S^1$. To check that it is actually a bundle space, we let $U_i\subset CP^n$ be the open set of equivalence classes $[z_0, \cdots, z_n]$ such that $z_i\neq 0$. Let $h_i: p^{-1}(U_i)\rightarrow U_i\times S^1$ by $h_i((z_0, \cdots, z_n))=([z_0, \cdots, z_n], z_i/|z_i|).$ We can do this since $z_i\neq 0$. This takes fibers to fibers and is a homeomorphism with inverse $([z_0, \cdots, z_n], \lambda)\rightarrow \lambda z_i^{-1}|z_i|$. So there is a local section.

Now let $n=1$. Recall that $CP^1$ consists of a 0-cell and a 2-cell so it is homeomorphic to $S^2$. So the bundle constructed as above becomes $S^1\rightarrow S^3\rightarrow S^2.$ This is our first Hopf bundle. The projection can be taken to be $(z_0, z_1)\rightarrow z_0/z_1\in C\cup\{\infty\} =S^2.$ The fiber, bundle space, and base space are all spheres. In polar coordinates we have $$p(r_0e^{i\theta_0}, r_1e^{i\theta_1})=\frac{r_0}{r_1}e^{i(\theta_0-\theta_1)}$$ where $r_0^2+r_1^2=1$.

Now replace $C$ by the field $H$ of {\it quaternions}\index{quaternion}. Recall from Example 3.3.3 that the quaternions are of the form $a+bi+cj+dk$ where $a,b, c, $ and $d$ are real numbers and $i^2=j^2=k^2=-1$. We also have $ij=k, jk=i, ki=j,  ji=-k, kj=-i,$ and $ik=-j$. (So the quaternions are not commutative.) Also define the {\it conjugate} $\overline{h}$ for $h\in H$. If $h=a+bi+cj+dk$, then $\overline{h}=a-bi-cj-dk$. Quaternionic projective space $HP^n$ is formed from $S^{4n+3}$ by identifying points which are multiples of each other by unit quaternions analogous to the complex case. It has a CW structure with a cell in dimensions that are a multiple of 4 up to $4n$. So $HP^1$ is $S^4$. There is a fiber bundle $S^3\rightarrow S^{4n+3}\rightarrow HP^n$ . Here $S^3$ is the unit quaternions and $S^{4n+3}$ is the unit sphere in $HP^n$. For $n=1$ this becomes  $S^3\rightarrow S^7\rightarrow S^4$.

We get our last Hopf map by using the octonians.

\begin{definition}
The ring $\mathcal{O}$ of {\it octonions}\index{octonion} or {\it Cayley numbers}\index{Cayley number} consists of pairs $(h_1, h_2)$ of quaternions with multiplication given by $$(a_1, a_2)(b_1, b_2)=(a_1b_1-\overline{b}_2a_2, a_2\overline{b}_1+b_2a_1).$$ Then $\mathcal{O}$ is a ring but it is non-commutative and not even associative. 
\end{definition}

Letting $S^{15}$ be the unit sphere in the 16-dimensional space $\mathcal{O}^2$, the projection map $p: S^{15}\rightarrow S^8=\mathcal{O}\cup\{\infty\}$ is $(z_0, z_1)\rightarrow z_0z_1^{-1}$ where $z_0, z_1\in\mathcal{O}$. This is a fiber bundle with fiber $S^7$, where $S^7$ is the unit octonians, but the proof is complicated by the fact that $\mathcal{O}$ is not associative. See \cite{Hat} for details. It turns out that there is an octonian projective plane $\mathcal{O}P^2$ formed by attaching a 16-cell to $S^8$ via the Hopf map $S^{15}\rightarrow S^8$. But there is no $\mathcal{O}P^n$ for $n>2$ as the associativity of multiplication is needed for the relation $(z_0,\cdots,z_n)\sim\lambda(z_0,\cdots,z_n)$. to be an equivalence relation.

I will now return to the first Hopf bundle $S^1\rightarrow S^3\rightarrow S^2$, and look at it in more detail. In \cite{Lyo}, David Lyons looks at the Hopf map in the context of rotations of 3-space, one of its original motivations. This connection allows for its use in physics areas such as magnetic monopoles \cite{Nak}, rigid body mechanics \cite{MR}, and quantum information theory \cite{MD}. Lyons' paper is aimed at undergraduates and doesn't use much topology but is very interesting for the way it illustrates the map and the spaces involved. I will digress for a while and summarize some of its main points as an aid to working with it in the context of Chapter 10.

Lyons refers to the Hopf bundle as the {\it Hopf fibration}. A fibration is technically a slightly more general term as the fibers are allowed to be homotopy equivalent in a fibration and not necessarily homeomorphic. In the literature, though, fibration and fiber bundle are often used interchangeably, and the Hopf fibration is definitely a fiber bundle in the more strict sense. Still, this allows us to mention the contribution of the Beach Boys to algebraic topology: "I'm picking up good fibrations."

Lyons uses an alternative formula for the projection which will be more useful for us and still produce a $S^1\rightarrow S^3\rightarrow S^2$ fiber bundle. Let $(a, b, c, d)\in S^3$,  so we have $a^2+b^2+c^2+d^2=1.$ Then the Hopf map $p: S^3\rightarrow S^2$ is defined as $$p(a, b, c, d)=(a^2+b^2-c^2-d^2, 2(ad+bc), 2(bd-ac)).$$ The reader should check that the square of the three coordinates on the right sum to $(a^2+b^2+c^2+d^2)^2=1$. 

To understand the properties of the Hopf map, we need to understand more about the quaternions and their relation to rotations. Suppose you want to look at a rotation in $R^3$. We can represent it by choosing an axis of rotation which can be represented by a vector in $R^3$ along with an angle of rotation. So we need a 4-tuple of real numbers. What if you want to look at the composition of two rotations. Given the axis of rotation and angle for each of them can you find these parameters for their composition? William Hamilton invented quaternions in the mid-19th century to handle problems like this. He was inspired by the corresponding problem in $R^2$.

In $R^2$, we can represent rotations of the plane around the origin by unit length complex numbers. If $z_1=e^{i\theta_1}$ and $z_2=e^{i\theta_2}$, then $z_1z_2=e^{i(\theta_1+\theta_2)}$. So if $\theta_1, \theta_2$ represent angles of rotation, then we compose the rotations by multiplying the corresponding complex numbers. We would like quaternions to have a similar property.

Let $H$ be the set of quaternions and $r=a+bi+cj+dk\in H$. Then the {\it conjugate} $\overline{r}=a-bi-cj-dk$ and the {\it norm} $$||r||=\sqrt{a^2+b^2+c^2+d^2}=\sqrt{r\overline{r}}.$$  Now the norm has the property that for $r, s\in H$, $||rs||=||r||||s||,$ so the product of quaternions of unit norm also has unit norm. A unit norm quaternion is a point of $S^3\subset R^4.$ If $r\neq 0$, then the multiplicative inverse of $r$ is $$r^{-1}=\frac{\overline{r}}{||r||^2}.$$ So if $r$ has unit norm, then $r^{-1}=\overline{r}$. Also, multiplication of quaternions is associative but not commutative.

Now we represent rotations in $R^3$ using quaternions as follows. For $p=(x, y, z)\in R^3$, associate a quaternion $p=xi+yj+zk$. (A quaternion with no real part is called {\it pure.}) If $r$ is an arbitrary nonzero quaternion then it turns out that $rpr^{-1}$ is also pure so it is of the form $x'i+y'j+z'k$ and can be associated with the point $(x', y', z')\in R^3$. So $r$ defines a mapping $R_r: R^3\rightarrow R^3$. The following theorem can be proved by direct calculation.

\begin{theorem}
The map $R_r: R^3\rightarrow R^3$ where $r\in H$ has the following properties: \begin{enumerate}
\item $R_r$ is a linear map.
\item If $k$ is a nonzero real number, then $R_{kr}=R_r$.
\item If $r\neq 0$, then $R_r$ is invertible and $(R_r)^{-1}=R_{r^{-1}}$.
\end{enumerate}
\end{theorem} 

Property 2 implies that we can restrict $r$ to have unit norm. 

The next two results show how to use quaternions to define rotations. See \cite{Lyo} for an outline of the proof.

\begin{theorem}
Let $r=a+bi+cj+dk$ be a quaternion of unit length. If $r=\pm 1$, $R_r$ is obviously the identity. Otherwise, $R_r$ is the rotation about the axis determined by the vector $(b, c, d)$ with an angle of rotation $$\theta=2\cos^{-1}(a)=2\sin^{-1}(\sqrt{b^2+c^2+d^2}).$$
\end{theorem}

\begin{theorem}
Let $r$ and $s$ be unit quaternions. Then $R_sR_r=R_{rs}$. So composition of rotations corresponds to multiplication of quaternions. Note: My left side is in reverse order to the statement in \cite{Lyo} to be consistent with my convention of composing functions from right to left. Use $r=i$ and $s=j$ to confirm that this is correct.
\end{theorem}

Now we will redefine the Hopf map in terms of quaternions. Fix a distinguished point $P_0=(1, 0, 0)\in S^2$. Let $r=a+bi+cj+dk$ be a unit quaternion corresponding to the point $(a, b, c, d)\in S^3$. Then the projection $p: S^3\rightarrow S^2$ of the Hopf fibration is defined as $$p(r)=R_r(P_0)=ri\overline{r}.$$ Computing the formula explicitly shows that it corresponds to Lyons' earlier formula. 

Using that formula, what is the inverse image of $(1, 0, 0)\in S^2$? Since, $$p(a, b, c, d)=(a^2+b^2-c^2-d^2, 2(ad+bc), 2(bd-ac))=(1, 0, 0),$$ we see that the set of points $$C=\{(\cos t, \sin t, 0, 0)|t\in R\}$$ maps to $(1, 0, 0)$. Since $$1=a^2+b^2+c^2+d^2=a^2+b^2-c^2-d^2$$ we have that $c, d=0$ so $C$ is the entire inverse image of $(1, 0, 0)$. 

It turns out that the inverse image of every point is a circle, confirming that the fiber of the Hopf filtration is $S^1$. See \cite{Lyo} for more discussion. 

The takeaway for what will come later is that we have some very explicit formulas for a map from a higher to a lower dimensional sphere. The Hopf map turns out to be the generator of $\pi_3(S^2)$, the group of homotopy classes of maps from $S^3$ to $S^2$. I will discuss more implications of this in Chapter 10,

\section{Paths and Loops}

In this section, I will describe an important fiber bundle. The material comes from \cite{Hu}. As we will see, fiber bundles lead to a long exact sequence in homotopy. The {\it path-space fibration} will be our goal. It has the advantage of having a total space which is contractible. 

We will write $Y^X$ for the set of continuous maps from $X$ to $Y$. We will generally consider $Y^X$ to be a topological space having the {\it compact-open topology}. Recall that this means that the set of continuous functions $f: X\rightarrow Y$  such that $f(K)\subset V$ for all $K\subset X$ compact and $V\subset Y$ open forms a subbase for the topology.

Reusing notation, we will use $\Omega=Y^I$ for the space of paths on $Y$ with the compact-open topology. Paths break $Y$ into {\it path components} for which any two points can be connected by a path. If $y\in Y$, then $C_y$ will be the path component containing $y$. 

\begin{definition}
A {\it generalized triad}\index{generalized triad} $(Y; A, B)$ is a space $Y$ together with two subspaces $A$ and $B$. $(Y; A, B)$ is a {\it triad}\index{triad} if $A\cap B\neq \emptyset$. For a generalized triad $(Y; A, B)$ let $[Y; A, B]$ be the subset of $\Omega$ consisting of paths $\sigma$ in $Y$ such that $\sigma(0)\in A$ and $\sigma(1)\in B$.
\end{definition}

Here are some interesting special cases. If $A=Y$  and $B$ consists of a single point $y$, the subspace [Y; Y, y] of $\Omega$ is denoted $\Omega_y$ and called the {\it space of paths that terminate at} $y$. 

Let $\Lambda$ be the space of loops in $Y$. Then $[Y; y, y]$ is the space of paths that start and end at $y$, so by definition it is {\it the space of loops on Y with base point at y}. We write this space as $\Lambda_y$. Let $\delta_y$  be the degenerate loop $\delta_y(I)=y$. Hu uses some facts about the topology of function spaces to prove the following important result:

\begin{theorem}
The space $\Omega_y$ of paths that terminate in $y$, is contractible to the point $\delta_y$.
\end{theorem}

I can now define the {\it path-space fibration}. See \cite{Hu} sections 12 and 13 for the proof that it is actually a fiber space. 

\begin{definition}
The space $\Omega_y$ is a fiber space called the {\it path-space fibration}\index{path-space fibration} over $Y$ relative to the {\it initial projection} defined by $p: \Omega_y\rightarrow Y$ where for $\sigma\in\Omega_y$, $p(\sigma)=\sigma(0)$. Since any path $\sigma\in\Omega_y$ terminates at $y$, if $p(\sigma)=y$, then $\sigma$ also begins at $y$, so $\sigma$ is a loop with base point $y$. So the fiber is $\Lambda_y$. 
\end{definition}

\section{Higher Homotopy Groups}

As I mentioned before, homotopy groups classify maps from $S^n$ to $X$ up to homotopy equivalence. In this section, I will precisely define them and list some of their basic properties.

\subsection{Definition of Higher Homotopy Groups}
 In order to make it easier to define relative homotopy and a group multiplication, I will follow Hu and replace $S^n$ with an $n$-dimensional cube whose boundary is identified to a point. We will use the notation $f:(X, A)\rightarrow (Y, B)$ to mean that $f: X\rightarrow Y$ and $f(A)\subset B$. This can generalize to a triple in the obvious way.

\begin{definition}
Let $n>1$ and $I^n$ be the $n$-dimensional cube that is the product of $I=[0,1]$ with itself $n$ times. Let $\partial I^n$ be the boundary of $I^n$. Let $F^n(X, x_0)$ be the set of maps $f: (I^n, \partial I^n)\rightarrow (X, x_0)$. The maps have an addition $f+g$ defined by $$(f+g)(t)=
\left
\{
\begin{array}{ll}
f(2t_1, t_2, \cdots,t_n) & \mbox{if }0\leq t_1\leq \frac{1}{2},\\
g(2t_1-1,  t_2, \cdots,t_n)  & \mbox{if }\frac{1}{2}\leq t_1\leq 1.
\end{array}\right.
$$for $t=(t_1, \cdots, t_n)\in I^n$. For $f, g\in F^n(X, x_0)$, $f+g\in F^n(X, x_0)$. If $[f]$ is the homotopy class of $f\in F^n(X, x_0)$, then define $[f]+[g]=[f+g]$. The homotopy classes of $F^n(X, x_0)$ form a group under this addition called the {\it n-th homotopy group of X at} $x_0$ denoted $\pi_n(X, x_0)$.\index{homotopy group}\index{$\pi_n(X, x_0)$} The identity is the class [0] of the map $(I^n, \partial I^n)\rightarrow (x_0, x_0)$ and the inverse of $[f]$ is the class $[f\theta]$, where $\theta: I^n\rightarrow I^n$ is defined as $\theta(t)=(1-t_1, t_2, \cdots, t_n)$. 

For $n=1$, $\pi_1(X, x_0)$ is the fundamental group we defined in Section 9.2. For $n=0$, $\pi_0(X, x_0)$ is the set of path components of $X$ and its neutral element is defined to be the component containing $x_0$. We say that $\pi_0(X, x_0)=0$ if $X$ is path connected. Note that $\pi_0(X, x_0)$ is not a group.
\end{definition}

Note that $\pi_n(X, x_0)$ could have been defined if we used maps $f: (S^n, s_0)\rightarrow (X, x_0)$. The two halves of $I^n$ defined by the conditions $t_1\leq \frac{1}{2}$ and $t_1\geq \frac{1}{2}$ correspond to the Southern and Northern hemispheres of $S^n$ respectively. For $n>1$, we can rotate $S^n$ to exchange hemispheres while keeping $s_0$ fixed. This is the intuition for the following important result.

\begin{theorem}
For $n>1$, $\pi_n(X, x_0)$ is an abelian group. 
\end{theorem}

As we saw earlier, $\pi_1(X, x_0)$ need not be abelian. It is abelian, though, in some interesting special cases.

\begin{definition}
Let $X$ be a topological space and $\mu: X\times X\rightarrow X$ define a multiplication on $X$. $X$ is called an {\it H-space}\index{H-space} if the $\mu$ is continuous and there exists a point $e\in X$ such that $\mu(x, e)=\mu(e, x)=x$ for all $x\in X$. 
\end{definition}

The spheres $S^0, S^1, S^3,$ and $S^7$ are H-spaces when considered to be the norm one elements of the reals, complexes, quaternions, and octonions. There is a theorem that these are the only spheres that are H-spaces. $S^7$ is an H-space but not a group as the octonions are not associative.  The space of loops on $X$ with base point $x_0$ is another example of an H-space, where the multiplication consists of travelling around the two loops one after the other. The identity is the degenerate loop consisting of only $x_0$. Hu proves the folowing:

\begin{theorem}
If $X$ is an H-space, then $\pi_1(X, x_0)$ is abelian.
\end{theorem}

Now if $p+q=n$, then $$X^{I^n}=X^{I^p\times I^q}=(X^{I^p})^{I^q},$$ so $\pi_n(X, x_0)=\pi_q(F^p, d_0),$ where $d_0$ is the constant loop. If $p=1$, then $F^p$ is the space of loops in $X$ with base point $x_0$. This proves the following.

\begin{theorem}
$\pi_n(X, x_0)\cong\pi_{n-1}(W, d_0)$ where $W$ is the space of loops on $X$ and $n>0$.
\end{theorem}

Since the space of loops $W$ is an H-space, $\pi_1(W, d_0)$ is abelian. As this must be isomorphic to $\pi_2(X, x_0)$ which is abelian by Theorem 9.6.1, the fact that $\pi_1$ need not always be abelian does not pose a problem in this case. 

Finally, the following is a direct result of the fact that $I^n$ is pathwise connected. 

\begin{theorem}
If $X_0$ denotes the path-component of $X$ containing $x_0$ then for $n>0$, $$\pi_n(X_0, x_0)\cong \pi_n(X, x_0).$$
\end{theorem}

We can already determine the homotopy groups of some spheres. First of all for $n>0$, $S^n$ is connected so $\pi_0(S^n, s_0)=0$. 

Due to the fact that homotopy classes of maps $S^n\rightarrow S^n$  are determined by the degree of the map (See Theorem 4.1.29), we see that $\pi_n(S^n)\cong Z$. Also, it turns out that $\pi_n(S^1)=0$ for $n>1$. This can be proved using covering spaces (covered in books such as \cite{Mun2, Hat, Hu} for example) or the long exact sequence of a fiber space which I will discuss in Section 9.6.4. 

It also turns out that $\pi_m(S^n)=0$ for $0<m<n$. This follows from the fact that if we represent the spheres as CW complexes, $S^n$ consisting of a 0-cell and an $n$-cell. Using the fact that a map between two CW complexes can be approximated by a cellular map, let $f: S^m\rightarrow S^n$ and assume that $f$ is cellular. But then $f$ takes $S^m$ into the $m$-skeleton of $S^n$ which consists of the single 0-cell. So any map $f: S^m\rightarrow S^n$ is homotopic to a constant and  $\pi_m(S^n)=0$ for $m<n$. 

On the other hand, if $m>n$, this is not the case. We will end this section by looking at $\pi_3(S^2)$.

We will start with two theorems proved in \cite{Hu}.

\begin{theorem}
{\bf Hopf Classification Theorem:} Let $X$ be a triangulable space with dimension less than or equal to $n$. Let $\chi$ be a generator of the infinite cyclic group $H^n(S^n)$. The assignment $f\rightarrow f^\ast(\chi)$ sets up a one-to-one correspondence between the homotopy classes of the maps $f: X\rightarrow S^n$ and the elements of the cohomology group $H^n(X)$. 
\end{theorem}

This theorem gives another proof that  $\pi_m(S^n)=0$ for $m<n$. Letting $X=S^m, H^n(S^m)=0$.

\begin{definition}
A map $f: X\rightarrow Y$ is algebraically trivial if for $f_\ast: H_m(X)\rightarrow H_m(Y)$ and $f^\ast: H^m(Y)\rightarrow H^m(X)$, $f_\ast=0$ and $f^\ast=0$ for all $m>0$.
\end{definition}

The Hopf map $p: S^3\rightarrow S^2$ is algebraically trivial as there is no $m>0$ for which $H_m(S^3)$ and $H_m(S^2)$ are both nonzero. The same holds for the cohomology groups. 

Let $X$ be a triangulable space. For an arbitrary map $F: X\rightarrow S^3$, the composed map $f=pF: X\rightarrow S^2$ where $p$ is the Hopf map is also algebraically trivial. 

\begin{theorem}
For any given triangulable space $X$, the assignment $F\rightarrow f=pF$ sets up a one-to-one correspondence between the homotopy classes of the maps $F: X\rightarrow S^3$ and those of the algebraically trivial maps $f: X\rightarrow S^2.$
\end{theorem}

Using Theorems 9.6.5 and 9.6.6 we get the following.

\begin{theorem}
The homotopy classes of the algebraically trivial maps $f: X\rightarrow S^2$ of a 3-dimensional triangulable space $X$ into $S^2$ are in a one-to-one correspondence with the elements of $H^3(X)$. For $\alpha\in H^3(X)$, we associate $\alpha$ with the homotopy class of the map $f=pF: X\rightarrow S^2$ such that $F: X\rightarrow S^3$ is a map with $F^\ast(\chi)=\alpha$ and $\chi$ is as defined in Theorem 9.6.5.
\end{theorem}

Now let $X=S^3$.

\begin{theorem}
The homotopy classes of the maps $f: S^3\rightarrow S^2$ are in one to one correspondence with the integers. For an integer $n$, we associate $n$ with the homotopy class of the map $f=pF: S^3\rightarrow S^2$ such that $F: S^3\rightarrow S^3$ has degree $n$.
\end{theorem}

The immediate consequence of Theorem 9.6.8 is that $\pi_3(S^2)\cong Z$. It turns out that the generator of this group is $[p]$ where $p$ is the Hopf map.

\subsection{Relative Homotopy Groups}

As in homology theory, there are relative homotopy groups and these groups will produce a long exact sequence. 

Let $X$ be a space, $A\subset X$, and $x_0\in A$. We call $(X, A, x_0)$ a {\it triplet}. If $A=x_0$, then we just write the pair $(X, x_0)$. 

\begin{definition}
Consider the $n$-cube $I^n$ for $n>0$. The initial {\it (n-1)-face} of $I^n$ identified with $I^{n-1}$ is defined by $t_n=0$. The union of the remaining $(n-1)$-faces is denoted by $J^{n-1}$. Then we have $$\partial I^n=I^{n-1}\cup J^{n-1}\hspace{.5 in}and\hspace{.5 in}\partial I^{n-1}=I^{n-1}\cap J^{n-1}.$$  (To help visualize this, let $n=2$. Then $I^1$ is the bottom edge and $J^1$ consists of the other three edges.) Let $$f: (I^n, I^{n-1}, J^{n-1})\rightarrow (X, A, x_0).$$ In particular, $f(\partial I^n)\subset A$ and $f(\partial I^{n-1})=x_0.$ Let $F^n=F^n(X, A, x_0)$ be the set of these maps. Then $\pi_n(X, A, x_0)$ is the set of homotopy classes of these maps relative to the system $\{I^{n-1}, A; J^{n-1}, x_0\}$ In other words, if $f$ and $g$ are in the same homotopy class, there is a continuous map $F: I^n\times I\rightarrow X$ such that $F(x, 0)=f(x)$, $F(x, 1)=g(x)$, and for $t\in[0,1]$, $f(t): (I^n, I^{n-1}, J^{n-1})\rightarrow (X, A, x_0)$, The set of these classes is the {\it n-th relative homotopy set} denoted $\pi_n(X, A, x_0)$. If $n>1$, then  $\pi_n(X, A, x_0)$ is a group with addition as defined for $\pi_n(X, x_0)$. This group is called a {\it relative homotopy group}.\index{relative homotopy group} 
\end{definition}

For absolute homotopy groups, $\pi_n(X, x_0)$ is a set for $n=0$, a group for $n>0$ and an abelian group for $n>1$. For relative homotopy groups, it turns out that $\pi_n(X, A, x_0)$  is a set for $n=1$, a group for $n>1$ and an abelian group for $n>2$.

I will now state three analogous properties of relative homotopy groups to those of the absolute homotopy groups. First I need a definition.

\begin{definition}
Let $T=(X, A, x_0)$ be a triplet. Let $X'=[X; X, x_0]$ be the space of paths on $X$ that terminate at $x_0$. Let $p: X'\rightarrow X$ be the initial projection. Let $$A'=p^{-1}(A)=[X; A, x_0]\subset X'.$$ In other words, $A'$ is the set of paths that start in $A$ and terminate at $x_0$. Let $x_0'$ be the degenerate loop $x'_0(I)=x_0$. Then $T'=[X', A', x'_0]$ is called the {\it derived triplet} of $T$.\index{derived triplet} The map $p: (X', A', x'_0)\rightarrow(X, A, x_0)$ is called the {\it derived projection}.
\end{definition}

\begin{theorem}
For every $n>0$, $$\pi_n(X, A, x_0)=\pi_{n-1}(A', x'_0)$$. 
\end{theorem}

This implies that every relative homotopy group can be expressed as an absolute homotopy group.

\begin{theorem}
If $X_0$ denotes the path component of $X$ containing $x_0$ and $A_0$  denotes the path component of $A$ containing $x_0$, then for $n>1$, $$\pi_n(X, A, x_0)=\pi_n(X_0, A_0, x_0).$$
\end{theorem}

\begin{theorem}
If $\alpha\in\pi_n(X, A, x_0),$ is represented by a map $f\in F^n(X, A, x_0)$ such that $f(I^n)\subset A$, then $\alpha=0$.
\end{theorem}

\subsection{Boundary Operator and Induced Homomorphisms}

Now that we have defined relative homotopy groups, we would like to use them in a long exact sequence like we have in homology. One of the components is a boundary operator $$\partial: \pi_n(X, A, x_0)\rightarrow \pi_{n-1}(A, x_0).$$ Unlike in homology, the definition in homotopy is very easy.

\begin{definition}
Let $(X, A, x_0)$ be a triplet. For $n>0$ define  $$\partial: \pi_n(X, A, x_0)\rightarrow \pi_{n-1}(A, x_0)$$ as follows: Let $\alpha\in\pi_n(X, A, x_0)$ be represented by a map $$f:(I^n, I^{n-1}, J^{n-1}).$$ If $n=1$, $f(I^{n-1})$ is a point of $A$ which determines a path component $\beta\in\pi_0(A, x_0)$ of $A$. If $n>1$, then the restriction of $f$ to $I^{n-1}$ is a map of $(I^{n-1}, \partial I^{n-1})$ into $(A, x_0)$, so it represents an element $\beta\in\pi_{n-1}(A, x_0).$ Since $\beta$ does not depend on the map $f$ representing $\alpha$, so define $\partial(\alpha)=\beta$. We call $\partial$ the {\it boundary operator}.
\end{definition}

\begin{theorem}
The boundary operator $\partial$ is a homomorphism for $n>1.$
\end{theorem}

Now consider a map $$f: (X, A, x_0)\rightarrow (Y, B, y_0).$$ We would like homotopy to have the functorial properties that homology has, so we would like to obtain a homomorphism on homotopy groups. This also turns out to be easy.

First of all, since $f$ is continuous, it sends path components of $X$ to path components of $Y$, so it induces a transformation $$f_\ast: \pi_0(X, x_0)\rightarrow\pi_0(Y, y_0)$$ which sends the neutral element of $\pi_0(X, x_0)$ to that of $\pi_0(Y, y_0)$. 

For $n>0$, If $\phi\in F^n(X, A, x_0)$, the composition $f\phi$ is in $F^n(Y, B, y_0)$. The assignment $f\rightarrow f\phi$ defines a map  $$f_\#: F^n(X, A, x_0)\rightarrow F^n(Y, B, y_0).$$  Since $f_\#$ is continuous, it carries the path components of $F^n(X, A, x_0)$ into those of $F^n(Y, B, y_0).$  This induces a transformation $$f_\ast: \pi_n(X, A, x_0)\rightarrow\pi_n(Y, B, y_0)$$ called the {\it induced transformation} which sends the neutral element of $\pi_n(X, A, x_0)$ to that of $\pi_n(Y, B, y_0)$. 

\begin{theorem}
Let $f: (X, A, x_0)\rightarrow (Y, B, y_0).$ If $n=1, f_\ast$ is a transformation which sends the neutral element of $\pi_n(X, A, x_0)$ to that of $\pi_n(Y, B, y_0)$. If $n>1$, then $$f_\ast: \pi_n(X, A, x_0)\rightarrow\pi_n(Y, B, y_0)$$ is a homomorphism.
\end{theorem}

\begin{theorem}
Let $f: (X, x_0)\rightarrow (Y, y_0).$ If $n=0, f_\ast$ is a transformation which sends the neutral element of $\pi_n(X, x_0)$ to that of $\pi_n(Y, y_0)$. If $n>0$, then $$f_\ast: \pi_n(X, x_0)\rightarrow\pi_n(Y, y_0)$$ is a homomorphism.
\end{theorem}

\subsection{Properties of Homotopy Groups}

In this section, I will list some properties of homotopy groups and compare them to those of homology groups. In the next section, I will define a {\it homotopy system}, the analogue of a homology theory. I will state most of these without proof, so see \cite{Hu} for the details.

\begin{theorem}
{\bf Property 1:} If $f: (X, A, x_0)\rightarrow (X, A, x_0)$ is the identity map, then $f_\ast$ is the identity transformation on $\pi_n(X, A, x_0)$ for all $n\geq 0$.
\end{theorem}

\begin{theorem}
{\bf Property 2:} If $f: (X, A, x_0)\rightarrow (Y, B, y_0)$ and $g: (Y, B, y_0)\rightarrow (Z, C, z_0)$ are maps, then $(gf)_\ast=g_\ast f_\ast$ for all $n\geq 0$. So the assignment $(X, A, x_0)\rightarrow \pi_n(X, A, x_0)$ and $f\rightarrow f_\ast$ is a covariant factor from the category of triples to the category of homotopy groups.
\end{theorem}

\begin{theorem}
{\bf Property 3:} If $f: (X, A, x_0)\rightarrow (Y, B, y_0)$ is a map and $g:(A, x_0)\rightarrow (B, y_0)$ is the restriction of $f$ to $A$, then the following rectangle is commutative for all $n>0$:

$$\begin{tikzpicture}
  \matrix (m) [matrix of math nodes,row sep=3em,column sep=4em,minimum width=2em]
  {
\pi_n(X, A, x_0) & \pi_{n-1}(A, x_0)\\
\pi_n(Y, B, y_0) &  \pi_{n-1}(B, y_0)\\};

  \path[-stealth]
(m-1-1) edge node [above] {$\partial$} (m-1-2)
(m-1-1) edge node [left] {$f_\ast$} (m-2-1)
(m-1-2) edge node [right] {$g_\ast$} (m-2-2)
(m-2-1) edge node [below] {$\partial$} (m-2-2)

;

\end{tikzpicture}$$
\end{theorem}

Properties 1-3 correspond to those of homology groups. There is also a long exact sequence for a triplet which differs from the long exact sequence of a pair in homology only by the inclusion of the base point. For the pieces that are sets and not groups, the term {\it kernel} will refer to the inverse image of the neutral element which must then coincide with the image of the previous map to extend the definition of exactness.

\begin{theorem}
{\bf Property 4 (Exactness Property):} If $(X, A, x_0)$ is a triplet, let $i: (A, x_0)\subset (X, x_0)$ and $j: (X, x_0)=(X, x_0, x_0)\subset (X, A, x_0)$ be inclusions, and $i_\ast, j_\ast$ be the induced transformations. Let $\partial$ be the boundary operator defined above. Then the following {\it homotopy sequence of the triplet} is exact: $$\cdots\xrightarrow{j_\ast}\pi_{n+1}(X, A, x_0)\xrightarrow{\partial}\pi_n(A, x_0)\xrightarrow{i_\ast}\pi_n(X, x_0)\xrightarrow{j_\ast}\pi_n(X, A, x_0)\xrightarrow{\partial}\cdots$$ $$\cdots\xrightarrow{j_\ast}\pi_1(X, A, x_0)\xrightarrow{\partial}\pi_0(A, x_0)\xrightarrow{i_\ast}\pi_0(X, x_0)$$
\end{theorem}

\begin{theorem}
{\bf Property 5 (Homotopy Property):}  If $f, g: (X, A, x_0)\rightarrow (Y, B, x_0)$ are homotopic, then $f_\ast, g_\ast: \pi_n(X, A, x_0)\rightarrow \pi_n(Y, B, x_0)$  are equal for every $n$.
\end{theorem}

Now topological invariance of homotopy groups are an immediate consequence of properties 1, 2, and 5.

\begin{theorem}
If $f: (X, A, x_0)\rightarrow (Y, B, x_0)$ is a homotopy equivalence, then the induced transformation $f_\ast: \pi_n(X, A, x_0)\rightarrow \pi_n(Y, B, x_0)$ is an isomorphism for all $n$. (One to one and onto for $n=1$ where $f$ is a map of sets rather than groups.)
\end{theorem}

Now you will see the special place that fiber spaces have in homotopy theory. They have a long exact sequence of their own which is very useful for computation. For homology, though things aren't so nice. That is where spectral sequences, which will be discussed in Section 9.8, come in.

\begin{theorem}
{\bf Property 6 (Fibering Property):}  If $p: (E, C, x_0)\rightarrow (B, D, y_0)$ is a fibering, and $C=p^{-1}(D)$, then $$p_\ast: \pi_n(E, C, x_0)\rightarrow \pi_n(B, D, y_0)$$ is an isomorphism for all $n$. (One to one and onto for $n=1$ where $f$ is a map of sets rather than groups.)
\end{theorem}

A very important consequence of this property is the following. Let $D=y_0$ and $C=p^{-1}(y_0)=F$ where $F$ is the fiber of $p$ by the definition of the fiber. So our fibering becomes  $p: (E, F, x_0)\rightarrow (B, y_0)$ . By Property 6, $p_\ast$ is an isomorphism on homotopy groups, so replacing $\pi_n(E, F, x_0)$ by $\pi_n(B, y_0)$ in the long exact sequence of the triplet $(E, F, x_0)$ gives the long exact sequence of the fiber space.

\begin{theorem}
If $p: E\rightarrow B$ is a fibering with fiber $F$, and suppose all $E$,$B$, and $F$ are all path connected (so we can ignore the base point), then we have a long exact sequence: $$\cdots\xrightarrow{j_\ast}\pi_{n+1}(B)\xrightarrow{\partial}\pi_n(F)\xrightarrow{i_\ast}\pi_n(E)\xrightarrow{j_\ast}\pi_n(B)\xrightarrow{\partial}\cdots$$ $$\cdots\xrightarrow{j_\ast}\pi_1(B)\xrightarrow{\partial}\pi_0(F)\xrightarrow{i_\ast}\pi_0(E)$$
\end{theorem}

\begin{example}
For the Hopf map, $E=S^3, B=S^2$, and $F=S^1$. Since $\pi_n(S^1)=0$ for $n>1$,we get from the long exact sequence that for $n\geq 3$, $\pi_n(S^3)\cong\pi_n(S^2)$. (Write out some terms and check this for yourself.) 
\end{example}

In homology theory, the fibering property is generally false. It is replaced by the excision property which is false in homotopy theory.

For our final property, let $X$ consist of one point $x_0$. Then the only map from $I^n$ into $X$ is the constant map. This gives the following.

\begin{theorem}
{\bf Property 7 (Triviality Property):}  If $X$ consists of the single point $x_0$, then $\pi_n(X, x_0)=0$ for all $n$.
\end{theorem}

\subsection{Homotopy Systems vs. Eilenberg-Steenrod Axioms}

Properties 1-7 describe a {\it homotopy system}. In this section, I will give a precise definition. Compare it to the definition of a homology theory we saw in Section 4.2.3.

\begin{definition} 
A {\it homotopy system}\index{homotopy system} $$H=\{\pi, \partial, \ast\}$$ consists of three functions $\pi$, $\partial$, $\ast$. The function $\pi$ assigns to each triplet $(X, A, x_0)$ and each integer $n>0$, the set $\pi_n(X, A, x_0)$. The function $\partial$ assigns to each triplet $(X, A, x_0)$ and each integer $n>0$, a transformation $$\partial: \pi_n(X, A, x_0)\rightarrow \pi_{n-1}(A, x_0),$$ where if $n=1$, $\pi_0(A, x_0)$ denotes the set of all path components of $A$, and for $n>0$, $\pi_n(A, x_0)$ is defined to be $\pi_n(A, x_0, x_0)$. The function $\ast$ assigns to each map $$f: (X, A, x_0)\rightarrow (Y, B, y_0)$$ and each integer $n>0$, a transformation $$f_\ast: \pi_n(X, A, x_0)\rightarrow \pi_n(Y, B, y_0).$$

The system $H$ must satisfy 7 axioms. For $1\leq i\leq 7$, Axiom $i=$Property $i$ from above with 2 exceptions. 

{\bf Axiom 4:} Similar to Property 4 but the homotopy sequence of a pair need only be {\it weakly exact}. This means that if $\pi_n(X, x_0)=0$ for all $n$, then $\partial$ is one to one and onto.

{\bf Axiom 6:} If $p: (X', A', x'_0)\rightarrow (X, A, x_0)$ is the derived projection then $p_\ast$ is a one to one and onto map $$p: \pi_n(X', A', x'_0)\rightarrow\pi_n (X, A, x_0)$$ for all $n>0$.

The derived projection is the same as the path-space fibration, so it is actually weaker than Property 6.
\end{definition}

\begin{theorem}
$H=\{\pi, \partial, \ast\}$ as we have defined the three functions in this chapter is an example of a homotopy system.
\end{theorem}

Hu proves that any homotopy system is in a sense equivalent to ours. See \cite{Hu} for details.

\subsection{Operation of the Fundamental Group on the Higher Homotopy Groups}

The goal of the remainder of this section is to determine if we really need the base point and what it does. When can we dispose of it altogether and just talk about the homotopygroups $\pi_n(X)$?

We need a couple of definitions.

\begin{definition}
An {\it automorphism}\index{automorphism} of a group is an isomorphism of the group to itself.
\end{definition}

\begin{definition}
The {\it action} of a group $G$ on a set $S$ is a function $f: G\times S\rightarrow S$ where we write $f(g, s)$ as $gs$ and such that for $g, h\in G$ and $s\in S$, \begin{enumerate}
\item $1_Gs=s$.
\item $g(hs)=(gh)s.$
\end{enumerate}
\end{definition}

Let $X$ be a space and let $x_0, x_1\in X$. Let $\sigma: I\rightarrow X$ be a path connecting them so that $\sigma(0)=x_0$ and $\sigma(1)=x_1$. Then $x_0, x_1$ lie in the same path component of $X$ so $\pi_0(X, x_0)=\pi_0(X, x_1)$. For $n>0$, Hu proves the following result.

\begin{theorem}
For each $n>0$, every path $\sigma: I\rightarrow X$ with $x_0=\sigma(0)$ and $x_1=\sigma(1)$ gives an isomorphism $$\sigma_n: \pi_n(X, x_1)\cong\pi_n(X, x_0)$$ which depends only on the homotopy class of $\sigma$ relative to the endpoints. (I. e. the endpoints are kept fixed.) If $\sigma$ is the degenerate path $\sigma(I)=x_0$ then $\sigma_n$ is the identity automorphism. If $\sigma, \tau$ are paths with $\tau(0)=\sigma(1)$, then $(\sigma\tau)_n=\sigma_n\tau_n$. For every path $\sigma: I\rightarrow X$ and map $f: X\rightarrow Y$, let $\tau=f\sigma$ be the corresponding path in $Y$ with $y_0=f(x_0)$ and $y_1=f(x_1)$. Then the following commutes:
$$\begin{tikzpicture}
  \matrix (m) [matrix of math nodes,row sep=3em,column sep=4em,minimum width=2em]
  {
\pi_n(X, x_1) & \pi_n(X, x_0)\\
\pi_n(Y, y_1) &  \pi_n(Y, y_0)\\};

  \path[-stealth]
(m-1-1) edge node [above] {$\sigma_n$} (m-1-2)
(m-1-1) edge node [left] {$f_\ast$} (m-2-1)
(m-1-2) edge node [right] {$f_\ast$} (m-2-2)
(m-2-1) edge node [below] {$\tau_n$} (m-2-2)

;

\end{tikzpicture}$$

\end{theorem}

As a loop is really a path, the theorem immediately implies the following.

\begin{theorem}
The fundamental group $\pi_1(X, x_0)$ acts on $\pi_n(X, x_0)$ for $n\geq 1$ as a group of automorphisms. 
\end{theorem}

Theorem 9.6.25 implies that for a pathwise connected space, the homotopy groups are isomorphic for any choice of a basepoint and we can write $\pi_n(X)$ for these groups. To give these groups a geometric meaning, though, we need one more condition on $X$. We need it to be {\it n-simple}.

\begin{definition}
A group $G$ acts {\it simply} on a set $S$ if $gs=s$ for all $g\in G$, $s\in S$. A space $X$ is {\it n-simple}\index{n-simple} if $\pi_1(X, x_0)$ acts simply on $\pi_n(X, x_0)$ for all $x_0\in X$.
\end{definition}

The following properties are immediate. 

\begin{theorem}
A pathwise connected space is $n$-simple if there exists an $x_0\in X$ such that $\pi_1(X, x_0)$ acts simply on $\pi_n(X, x_0)$.
\end{theorem}

\begin{theorem}
A simply connected space is $n$-simple for every $n>0$.
\end{theorem}

\begin{theorem}
A pathwise connected space is $n$-simple if $\pi_n(X)=0.$.
\end{theorem}

It can be shown that $\pi_1(X, x_0)$ acts on itself by {\it conjugation}. In other words, for $g, h$ in $\pi_1(X, x_0)$, $hg=h^{-1}gh$. Then we have the following.

\begin{theorem}
A pathwise connected space $X$ is 1-simple if and only if $\pi_1(X)$ is abelian.
\end{theorem}

\begin{example}
The sphere $S^m$ is $n$-simple for $m, n>0$. This holds for $m>1$ and any $n$ by Theorem 9.6.28. By Theorem 9.6.29, $S^1$ is $n$-simple for $n>1$, and we know that $S^1$ is 1-simple by Theorem 9.6.30.
\end{example}

To get a better feel for the geometric meaning of $n$-simplicity we have the following:

\begin{theorem}
A space $X$ is $n$-simple if and only if for every point $x_0\in X$ and maps $f, g: S^n\rightarrow X$ with $f(s_0)=x_0=g(s_0)$, $f\simeq g$ implies $f\simeq g$ rel $s_0$.
\end{theorem}

{\bf Proof:} Let $X$ be $n$-simple. Then $\pi_1(X, x_0)$ acts simply on $\pi_n(X, x_0)$. Since $f\simeq g$, there exists a homotopy $h_t$ with $h_0=f$ and $h_1=g$. Let $f$ and $g$ represent the elements $\alpha$ and $\beta$ of $\pi_n(X, x_0)$ respectively. Let $\sigma: I\rightarrow X$ be a path defined by $\sigma(t)=h_t(s_0)$ for $t\in I$. Since $\sigma(0)=x_0=\sigma(1)$, $\sigma$ represents and element $w\in \pi_1(X, x_0)$ and it can be shown that $\alpha=w\beta$. Since $\pi_1(X, x_0)$ acts simply on $\pi_n(X, x_0)$, $\alpha=\beta$. So by definition of $\pi_n(X, x_0)$, $f\simeq g$ rel $s_0$.

For the other direction, let $w\in \pi_1(X, x_0)$ and let $\sigma$ be a loop representing $w$. If $\alpha\in\pi_n(X, x_0)$ is represented by a map $f: S^n\rightarrow X$ with $f(s_0)=x_0$, then the element $w\alpha$ of $\pi_n(X, x_0)$ is represented by a map $g:S^n\rightarrow X$ with $g(s_0)=x_0$ and satisfying $f\simeq g$. This implies that $f\simeq g$ rel $s_0$, so $w\alpha=\alpha$ and $X$ is $n$-simple. $\blacksquare$

We can use this result to prove the following.

\begin{theorem}
A pathwise connected topological group $X$ is $n$-simple for every $n>0$.
\end{theorem}

{\bf Proof:} Let $x_0$ be the identity element of $X$, and $f, g: S^n\rightarrow X$ be two homotopic maps with $f(s_0)=g(s_0)=x_0$. Then there is a homotopy $h_t$ with $h_0=f$, and $h_1=g$. Define a homotopy $k_t: S^n\rightarrow X$, by taking $$k_t(s)=[h_t(s_0)]^{-1}[h_t(s)],$$ for $s\in S^n$ and $t\in I$. Then $k_0=f$ and $k_1=g$ since $f(s_0)=g(s_0)=x_0$, and $k_t(s_0)=x_0$ for all $t\in I$. Then $f\simeq g$ rel $s_0$. By the previous theorem, $X$ is $n$-simple. $\blacksquare$

If we consider $\pi_n(X)$ to be the set of homotopy classes of maps $S^n\rightarrow X$ and $\pi_n(X, x_0)$ to be the set of homotopy classes of maps $(S^n, s_0)\rightarrow (X, x_0)$. Obviously $\pi_n(X, x_0)$ is a subgroup of  $\pi_n(X)$. Hu shows the following.

\begin{theorem}
If $X$ is pathwise connected and $n$-simple, then the inclusion  $\pi_n(X, x_0)\rightarrow\pi_n(X)$ is one to one and onto.
\end{theorem}

In this case, we can safely drop the basepoint and talk about the group $\pi_n(X)$.

\section{Calculation of Homotopy Groups}

So how do you calculate homotopy groups? It's really hard. There is no general algorithm like there is for homology. Also, even for the sphere $S^2$ we don't know alll of the homotopy groups, and we never will. Still, there are some useful results that can help in a lot of cases. I will mostly list some of the well known theorems in this section. You should think about what is easy in homotopy and hard in homology or vice versa. In this section, the material comes from Hu \cite{Hu} unless otherwise stated.

\subsection{Products and One Point Union of Spaces}

The product of two spaces is a good example of a case that is easier for homotopy than for homology. Recall that for homology and cohomology, we had the K\"{u}nneth Theorems (Theorems 8.5.6 and 8.5.10) which give the product as a term in an exact sequence involving the factors. For homotopy, the formula takes a much simpler form.

\begin{theorem}
Let $X$ and $Y$ be spaces with $x_0\in X$ and $y_0\in Y$. Let $Z=X\times Y$ and $z_0=(x_0, y_0)\in Z$. Then for every $n>0$, $$\pi_n(Z, z_0)\cong\pi_n(X, x_0)\oplus\pi_n(Y, y_0).$$ (Here we are using direct additive notation despite the fact that $\pi_1$ may be non-abelian.)
\end{theorem}

{\bf Proof:} Let $p: (Z, z_0)\rightarrow (X, x_0)$ and $q: (Z, z_0)\rightarrow (Y, y_0)$ be projections and $i: (X, x_0)\rightarrow (Z, z_0)$ and $j: (Y, y_0)\rightarrow (Z, z_0)$ be inclusions. passing to the corresponding homotopy groups gives $$p_\ast i_\ast=1,\hspace{1 in} q_\ast j_\ast=1, \hspace{1 in}p_\ast j_\ast=0, \hspace{1 in}q_\ast i_\ast=0.$$

Let $$h: \pi_n(Z, z_0)\rightarrow\pi_n(X, x_0)\oplus\pi_n(Y, y_0)$$ be a homomorphism defined by $h(\alpha)=(p_\ast(\alpha), q_\ast(\alpha))$ for $\alpha\in \pi_n(Z, z_0)$. We need to show that $h$ is an isomorphism. Letting $\alpha\in\pi_n(X, x_0)$ and $\beta\in\pi_n(Y, y_0)$, let $\gamma=i_\ast(\alpha)+j_\ast(\beta)\in \pi_n(Z, z_0)$. Then $$h(\gamma)=(p_\ast i_\ast\alpha+p_\ast j_\ast\beta, q_\ast i_\ast\alpha+q_\ast j_\ast\beta)=(\alpha, \beta),$$ so $h$ is an epimorphism.

Now let $\delta\in\pi_n(Z, Z_0)$ and $h(\delta)=0.$ Then we have $p_\ast\delta=0$ and $q_\ast\delta=0$. Let $f: (I^n, \partial I^n)\rightarrow(Z, z_0)$ represent $\delta$. Then since $p_\ast\delta=0$ and $q_\ast\delta=0$, there are homotopies $g_t: (I^n, \partial I^n)\rightarrow(X, x_0)$ and $h_t: (I^n, \partial I^n)\rightarrow(Y, y_0)$ with $g_0=pf$, $h_0=qf$ $g_1(I^n)\rightarrow x_0$ and $h_1(I^n)\rightarrow y_0$.  Define a homotopy $f_t: I^n\rightarrow Z$ by letting $f_t(s)=(g_t(s), h_t(s)).$ Then $f_0=f$, $f_1(I^n)=z_0$, and $f_t(\partial I^n)=z_0$ for all $t\in I$. Thus $\delta=0$ and $h$ is a monomorphism. Thus $h$ is an isomorphism. $\blacksquare$

Note that the inverse of $h$ is $$h^{-1}: \pi_n(X, x_0)\oplus\pi_n(Y, y_0)\rightarrow\pi_n(Z, z_0)$$ is given by $h^{-1}(\alpha, \beta)=i_\ast(\alpha)+j_\ast(\beta).$

\begin{example}
Lets return to our favorite hyperdonut shop. First we will be boring and get a regular donut $T=S^1\times S^1$. Then $\pi_1(T)=\pi_1(S_1)\oplus \pi_1(S_1)\cong Z\oplus Z$, and $\pi_n(T)=0$ for $n>1$. The formula extends to any number of factors so $\pi_1(S^1\times\cdots\times S^1)\cong Z\oplus\cdots\oplus Z$ where the free abelian group on the right has rank equal to the number of factors in the product and $\pi_n(S^1\times\cdots\times S^1)=0$. For a hyperdonut $X=S^m\times S^n$ for $m, n>0$, $\pi_i(X)\cong\pi_i(S^m)\oplus\pi_i(S^n)$. 
\end{example}

Our next example is a little more complicated for homotopy than for homology. Let $x_0\in X$ and $y_0\in Y$. Identify $x_0$ with $y_0$ and let $X\vee Y$ be the result.  Then we call $X\vee Y$ the {\it one point union} of $X$ and $Y$. It can be imbedded in $X\times Y$  by a map $k$ taking $X\vee Y$ to the subset of $X\times Y$ of the form $(x, y_0)\cup(x_0, y)$. For homology we have that for $n>0$, $H_n(X\vee Y)\cong H_n(X)\oplus H_n(Y)$. For example, we can use a Mayer-Vietoris sequence as $(x, y_0)$ and $(x_0, y)$ are homotopy equivalent to $X$ and $Y$ respectively and their intersection is a single point. In homotopy we have a slightly more complicated formula.

\begin{theorem}
For every $n>1$, we have $$\pi_n(X\vee Y, (x_0\sim y_0))\cong \pi_n(X, x_0)\oplus\pi_n(Y, Y_0)\oplus\pi_{n+1}(X\times Y, X\vee Y, (x_0, y_0)).$$
\end{theorem}

See Hu for the details of the proof but I will comment that we need $n>1$ so that all of the groups are abelian. The extra term $\pi_{n+1}(X\times Y, X\vee Y, (x_0, y_0))$ on the right comes from the long exact sequence of the pair $(X\times Y, X\vee Y).$

Hu shows that in the special case of spheres we have the following. (This should make you hungry for some more hyperdonuts.

\begin{theorem}
For every $p, q>0$ and $n<p+q-1$, we have $$\pi_n(S^p\vee S^q)\cong \pi_n(S^p)\oplus\pi_n(S^q).$$
\end{theorem}

\subsection{Hurewicz Theorem}

I will now outline the famous Hurewicz Theorem which discusses the relation between homotopy and homology groups of a space. I will refer you to Hu for the full proof but it is important for the subsequent material that you understand the form of the map from the first nonzero homotopy group to the homology group of the same dimenstion.

Let $\alpha$ be an element of the relative homology group $\pi_n(X, A, x_0)$ for $n>0.$ (Remember that for $n=1$, this is a set and not a group.) Let $$\phi: (E^n, S^{n-1}, s_0)\rightarrow (X, A, x_0)$$ represent $\alpha$, and let $E^n$ be a unit $n$-ball in $R^n$, $S^{n-1}$ be the unit $(n-1)$-sphere bounding $E^n$ and $s_0=(1,0,\cdots, 0).$ The coordinate system in $R^n$ determines an orientation and thus a generator $\xi_n$ of the group $H_n(E^n, S^{n-1})\cong Z$. 

Since $\phi$ maps $(E^n, S^{n-1})$ into $(X, A)$ it induces a map on homology $\phi_\ast: H_n(E^n, S^{n-1})\rightarrow H_n(X, A)$ where $H_n(X, A)$ denotes the singular homology group with integral coefficients. $\phi_\ast$ depends only on $\alpha\in\pi_n(X, A, x_0)$ so we have a function $$\chi_n: \pi_n(X, A, x_0)\rightarrow H_n(X, A)$$ where $\chi_n(\alpha)=\phi_\ast(\xi_n)$.

\begin{theorem}
If either $n>1$ or $A=x_0$, then $\chi_n$ is a homomorphism which will be called the {\it natural homomorphism} of $\pi_n(X, A, x_0)$ into $H_n(X, A)$.
\end{theorem}

\begin{theorem}
For any map $f: (X, A, x_0)\rightarrow (Y, B, y_0)$, the following is commutative:

$$\begin{tikzpicture}
  \matrix (m) [matrix of math nodes,row sep=3em,column sep=4em,minimum width=2em]
  {
\pi_n(X, A, x_0) & \pi_n(Y, B, y_0)\\
H_n(X, A) &  H_n(Y, B)\\};

  \path[-stealth]
(m-1-1) edge node [above] {$f_\ast$} (m-1-2)
(m-1-1) edge node [left] {$\chi_ n$} (m-2-1)
(m-1-2) edge node [right] {$\chi_n$} (m-2-2)
(m-2-1) edge node [below] {$f_\ast$} (m-2-2)

;

\end{tikzpicture}$$
 
\end{theorem}

For the case $A=x_0$, there is an isomorphism $$j_\ast; H_n(X)\cong H_n(X, x_0).$$ Let $$h_n=j^{-1}_\ast\chi_n: \pi_n(X, x_0)\rightarrow H_n(X).$$ This is called the {\it natural homomorphism} or {\it Hurewicz homomorphism}\index{Hurewicz homomorphism} of $\pi_n(X, x_0)$ into $H_n(X).$  For $n=1$, $h_1$ corresponds to the homomorphism $h_\ast$ in the proof of 
Theorem 9.2.7.

\begin{definition}
A space $X$ for $n\geq 0$ is {\it n-connected}\index{n-connected} if it is pathwise connected and $\pi_m(X)=0$ for all $0<m\leq n$. So a 1-connected space is simply connected.
\end{definition}

We can now state the Hurewicz Theorem. See Hu \cite{Hu} for the proof.

\begin{theorem}
{\bf Hurewicz Theorem:} If $X$ is an $(n-1)$-connected finite simplicial complex with $n>1$, then the reduced homology groups $\tilde{H}_i(X)=0$ for $0\leq i<n$, and the natural homomorphism $h_n$ is an isomorphism. 
\end{theorem}

So the first nonzero homotopy group is isomorphic to the homology group of the same dimension for a simply connected finite simplicial complex.

Recall that the case $n=1$ was handled in Theorem 9.2.8. As $\pi_1(X)$ may be non-abelian, we have to mod out by the commutator subgroup and make it abelian to produce $H_1(X)$. 

The map $h_n$ will be important for some of our later constructions.

\subsection{Freudenthal's Suspension Theorem}

In this section,we will follow Hatcher \cite{Hat}.

As stated above, excision generally doesn't hold for homotopy. A version of it does hold, though, in a range of dimensions called the {\it stable range}. 

\begin{theorem}
Let $X$ be a CW complex, and let $A, B\subset X$ be subcomplexes of $X$ such that $X=A\cup B$, and $C=A\cap B\neq \emptyset.$ If $(A, C)$ is $m$-connected and $(B, C)$ is $n$-connected with $m, n\geq 0$, then the map $\pi_i(A, C)\rightarrow \pi_i(X, B)$ induced by inclusion is an isomorphism for $i<m+n$ and an epimorphism for $i=m+n$.
\end{theorem}

This theorem is sometimes called the {\it homotopy excision theorem} as the rather lengthy proof makes use of an excision argument that holds in the specified dimension range.

Now recall that for a complex $X$, the cone $CX$ is the complex obtained by taking a point $x_0$ not in $X$ and connecting every point in $X$ with a line segment to $x_0$. As an example, it turns a circle into the usual meaning of a cone. Also remember that cones plug up holes and that they are acyclic, so that all reduced homology groups and homotopy groups are zero. The suspension $SX$ involves taking two points $x_0$ and $x_1$ and joining each to every point in $X$. The suspension of a circle is $S^2$ as you can see by holding two ice cream cones together along their wide ends. In fact we have that $S(S^n)\sim S^{n+1}$ for $n\geq 0$. 

\begin{definition}
The {\it suspension map}\index{suspension map} $\pi_i(X)\rightarrow \pi_{i+1}(SX)$ is defined as follows: Let $SX=C_+X\cup C_-X$ where $C_+X$ and $C_-X$ are two cones over $X$. The suspension map is the map $$\pi_i(X)\approx\pi_{i+1}(C_+X, X)\rightarrow \pi_{i+1}(SX, C_-X)\approx\pi_{i+1}(SX),$$ where the isomorphisms on the two ends come from the long exact homotopy sequences of the pairs $(C_+X, X)$ and $(SX, C_-X)$ respectively and the middle map is induced by inclusion.
\end{definition}

We can now state the important {\it Freundenthal Suspension Theorem}.

\begin{theorem}
{\bf Freudenthal Suspension Theorem:} The suspension map $\pi_i(S^n)\rightarrow \pi_{i+1}(S^{n+1})$ for $i>0$ is an isomorphism for $i<2n-1$ and an epimorphism for $i=2n-1$. This holds more generally for the suspension $\pi_i(X)\rightarrow \pi_{i+1}(SX)$ whenever $X$ is an $(n-1)$-connected CW complex. 
\end{theorem}

{\bf Proof:} From the long exact sequence of the pair $(C_+X, X)$, we see that $\pi_{i+1}(C_+X, X)\cong \pi_i(X)$ so that $(C_+X, X)$ is $n$-connected if $X$ is $(n-1)$-connected. The same holds for the pair $(C_-X, X)$. Replacing $i$ by $i+1$, $m$ by $n$, $A$ by $C_+X$, $B$ by $C_-X$, $X$ by $SX$ and $C$ by $X$ in Theorem 9.7.7, gives that the middle map in the suspension is an isomorphism for $i+1<n+n=2n$ and an epimorphism for $i+1=2n$. $\blacksquare$

\begin{definition}
The range of dimensions where the suspension $\pi_i(S^n)\rightarrow \pi_{i+1}(S^{n+1})$ is an isomorphism is called the {\it stable range}.\index{stable range} The study of homotopy groups of spheres in the stable range is called {\it stable homotopy theory}.\index{stable homotopy theory}
\end{definition}

\begin{theorem}
$\pi_n(S^n)\cong Z$ generated by the identity map for all $n>0$. The degree map $\pi_n(S^n)\rightarrow Z$ is an isomorphism. 
\end{theorem}

We have already seen this result, but it is also an immediate consequnce of the Freudenthal suspension theorem since for $i>1$, $\pi_i(S^i)\cong \pi_{i+1}(S^{i+1})$. Although we would only know that $\pi_1(S^1)\rightarrow \pi_2(S^2)$ is an epimorphism, we could use the Hopf fibration $S^1\rightarrow S^3\rightarrow S^2$ to show it is actually an isomorphism.

The Freudenthal Suspension Theorem is the most basic tool for computing homotopy groups of spheres. We will take another look at it in Chapter 12.

\subsection{Whitehead's Theorem}

It turns out as in the case of homology that two spaces that have the same homotopy groups do not have to be homotopy equivalent (let alone homeomorphic). We do have the following useful result for CW complexes. 

\begin{theorem}
{\bf Whitehead's Theorem:} Let $f: X\rightarrow Y$ induce isomorphisms $f_\ast: \pi_n(X)\rightarrow \pi_n(Y)$ for all $n\geq 0$. Then $f$ is a homotopy equivalence.
\end{theorem}

See Hatcher \cite{Hat} for a proof..

\section{Eilenberg-MacLane Spaces and Postnikov Systems}

In homology, spheres are the simplest type of space. In reduced homology, $\tilde{H}_n(S^n)\cong Z$ and $\tilde{H}_i(S^n)=0$ for $i\neq n$. In homotopy, this only holds for $S^1$, but for $S^n$ with $n>1$, we don't know all of the homotopy groups and never will. I will talk more about that in Chapter 12.

So what are the simplest spaces when it comes to homotopy? 

\begin{definition}
Let $n>1$ and $\pi$ be an abelian group. Then the {\it Eilenberg-MacLane space}\index{Eilenberg-MacLane space} denoted $K(\pi, n)$ \index{$K(\pi, n)$} is a space whose only nonzero homotopy group is $\pi_n(K(\pi, n))=\pi$. If $n=1$, we have a similar definition but $\pi$ need not be abelian.
\end{definition}

\begin{example}
We have already shown that $S^1$ is a $K(Z, 1)$ space.
\end{example}

\begin{example}
Consider the sphere $S^\infty$, which is the limit of including $S^n$ in $S^{n+1}$ as $n\rightarrow \infty$. This space is contractible (we will see why in Chapter 11), and the quotient map $p: S^\infty\rightarrow P^\infty$ which glues together antipodal points is a covering space, i.e. it is a fiber bundle with a discrete fiber. In this case, the fiber $F$  consists of two points, and $\pi_0(F)=Z_2$ while $\pi_n(F)=0$ for $n>0$.  Restricting to $S^n\rightarrow P^n$, we have the fiber bundle $F\rightarrow S^n\rightarrow P^n$ and the long exact sequence of the fiber bundle shows that $\pi_1(P^n)=Z_2$ and $\pi_i(P^n)=0$ for $1<i<n$. Letting $n\rightarrow \infty$ shows that $P^\infty$ is a $K(Z_2, 1)$ space. 
\end{example}

\begin{example}
From the bundle $S^1\rightarrow S^\infty\rightarrow CP^\infty$, the long exact sequence shows that $\pi_i(CP^\infty)\cong\pi_{i-1}(S^1)$ for all $i$ since $S^\infty$ is contractible. This makes $CP^\infty$ a $K(Z, 2)$ space.
\end{example}

So is there a $K(\pi, n)$ for any abelian $\pi$ and $n>1$? The answer is yes as we will now show. Note that most of these spaces are infinite dimensional and you will only visualize them in your worst nightmares. But building a CW complex of that form is surprisingly easy. I will present the argument found in the book by Mosher and Tangora \cite{MT}. We will make a lot of use of their book in Chapter 11.

\begin{theorem}
If $n>1$ and $\pi$ is an abelian group, then there exists a CW complex with the homotopy type of a $K(\pi, n)$ space.
\end{theorem}

{\bf Proof:} Let $$0\rightarrow R\rightarrow F\rightarrow \pi\rightarrow 0$$ be a free resolution of $\pi$ and let $\{a_i\}_{i\in I}$ and $\{b_j\}_{j\in J}$ be bases of $R$ and $F$ respectively where $I$ and $J$ are index sets. Let $K$ be the wedge (one point union) of spheres $S^n_j$ for $j \in J$. A more poetic and often used term is that $K$ is a {\it bouquet of spheres}\index{bouquet of spheres}. For each $a_i$, which represents a relation in $\pi$, take an $(n+1)$-cell $e_i$ and attach it to $K$ by a map $f_i: \dot{e}_i\rightarrow K$ where $[f_i]=a_i\in \pi_n(K)=F$, and $\dot{e}_i$ is the boundary of $e_i$. Let $X$ be the space formed from $K$ by attaching the cells $e_i$ as described. The dimension of $X$ is at most $n+1$ and $\pi_n(X)=\pi$. It is also $(n-1)$-connected by the Hurewicz Theorem as there are no cells in dimension less than $n$ so the reduced homology groups are $\tilde{H}_i(X)=0$ for $0\leq i<n$.

We now need to kill off the higher homotopy groups. That will be done with the following result.

\begin{theorem}
Let $Z$ be the complex formed by attaching an $(m+1)$-cell to a CW complex $Y$ by a map $f: S^m\rightarrow Y$. Then the inclusion $j: Y\rightarrow Z$ induces isomorphisms $j_\ast: \pi_i(Y)\equiv\pi_i(Z)$ for $i<m$.  In dimension $m$, $j_\ast$ is an epimorphism which takes the subgroup generated by $[f]$ to zero.
\end{theorem}

{\bf Proof:} The isomorphism for $i<m$ follows from the fact that $j$ can be approximated by a cellular map and that $Y$ and $Z$ have the same $m$-skeleton.  If $i, j$, and $k$ are inclusions the following diagram commutes. 

$$\begin{tikzpicture}
  \matrix (m) [matrix of math nodes,row sep=3em,column sep=4em,minimum width=2em]
  {
S^m=\dot{e}^{m+1} & e^{m+1}\\
Y &  Z=Y\cup_fe^{m+1}\\};

\path[-stealth]
(m-1-1) edge node [above] {$i$} (m-1-2)
(m-1-1) edge node [left] {$f$} (m-2-1)
(m-1-2) edge node [right] {$k$} (m-2-2)
(m-2-1) edge node [below] {$j$} (m-2-2)

;

\end{tikzpicture}$$

In homotopy we get $j_\ast f_\ast=k_\ast i_\ast=0$, since $e^{m+1}$ acyclic implies that $i_\ast=0$. If $[1]$ is the generator of $\pi_m(S^m)\cong Z$, then $j_\ast[f]=j_\ast f_\ast=0$. Since $j_\ast$ is a homomorphism, it takes the subgroup generated by $[f]$ to zero. Finally, the fact that $j_\ast$ is onto $\pi_m(Z)$ is true by cellular approximation since any map of $S^m$ into $Z$ can be deformed into the $m$-skeleton $Z^m=Y$. This proves Theorem 9.8.2. So the attaching maps are all killed off and the higher homotopy groups are zero, proving Theorem 9.8.1. $\blacksquare$

We will see that any two $K(\pi, n)$ spaces with $n>1$ have the same homotopy type. This will come from an important theorem about cohomology operations, the subject of Chapter 11.

The analog of an Eilenberg-MacLane space for homology is a {\it Moore Space}.

\begin{definition}
Let $n\geq 1$ and $G$ be an abelian group. Then the {\it Moore space}\index{Moore space} denoted $M(G, n)$ \index{$M(G, n)$} is a space whose only nonzero reduced homology group is $H_n(M(G, n))=G$. For $n>1$ we will need the space to be simply connected. 
\end{definition}

\begin{example}
Any $S^n$ is a $M(Z, n)$ space. 
\end{example}

\begin{example}
Let $X$ be $S^n$ with a cell $e^{n+1}$ attached by a map $S^n\rightarrow S^n$ having degree $m$. Then $X$ is a $M(Z_m, n)$ space. We can then form $M(G, n)$ with $G$ any finitely generated abelian group by taking the wedge of spheres and spaces of the form $M(Z_m, n)$. 
\end{example}

For a general abelian group $G$, there is a construction involving the free resolution of $G$ which is similar ot the construction of Eilenberg-MacLane spaces above. See Hatcher \cite{Hat} for details.

The final topic in this section is {\it Postnikov systems}\index{Postnikov system} also known as {\it Postnikov towers}\index{Postnikov tower}. They were invented in 1951 by Mikhail Postnikov and are a means to decompose a space with more than one nontrivial homotopy groups into simpler spaces. (I would have invented them myself if the "nikov" had been replaced by an "ol".) Postnikov's early papers are mainly in Russian, but the American Mathematical Society published a translation in 1957 \cite{Postn}. My advisor, Donald Kahn, wrote an early paper in 1963 \cite{Kahn} which shows how a map $f: X\rightarrow Y$ gives rise to a map of the Postnikov system for $X$ into the Postnikov system for $Y$. Postnikov systems are covered in a number of books (eg. \cite{Hat, Whi}), but I will follow the description in Mosher and Tangora \cite{MT}. Postnikov systems will be needed for the calculations in the next chapter that are involved in applying obstruction theory to data science.

We would like to find invariants which represent the homotopy type of a space through dimension $n$. Unfortunately, the $n$-skeleton  does not uniquely determine homotopy type through dimension $n$ as the following example shows. 

\begin{example}
Let $K$ denote $S^n$ as a CW complex with one $n$-cell and one 0-cell, where the $n$-cell is attached to the 0-cell along its boundary. Let $L$ be $S^n$ as a CW-complex built in a different way. $S^n$ consists of $S^{n-1}$ with two $n$-cells attached. $S^{n-1}$ is the equator and the two $n$-cells consist of the northern and southern hemispheres. Then $K$ and $L$ are homeomorphic, but the $(n-1)$-skeleton $K^{n-1}$ of $K$ is a single point, but $L^{n-1}$ is $S^{n-1}$. 
\end{example}

To get around this problem, we need Whitehead's theory of $n$-types \cite{JHCW}. First we state his {\it cellular approximation theorem} which we have informally used in our discussion of Eilenberg-MacLane spaces. 

\begin{definition}
A map $f: K\rightarrow L$ is {\it cellular} if $f(K^n)\subset L^n$. A homotopy $F: K\times I\rightarrow L$ between cellular maps is {\it cellular} if $F(K^n \times I)\subset L^{n+1}$ for every $n$.
\end{definition}

\begin{theorem}
{\bf Cellular Approximation Theorem: } Let $K_0$ be a subcomplex of $K$ and let $f: K\rightarrow L$ be a map such that $f|K_0$ is cellular. Then there exists a cellular map $g: K\rightarrow L$ such that $g\sim f$ rel $K_0$. 
\end{theorem}

\begin{theorem}
If $f: K\rightarrow L$, there exists a cellular map $g: K\rightarrow L$ which is homotopic to $f$. If two cellular maps $f$ and $g$ are homotopic, there exists a cellular homotopy between them.
\end{theorem}

Next we will define the ideas of $n$-homotopy type and $n$-type. 

\begin{definition}
Two maps $f, g: X\rightarrow Y$ are {\it n-homotopic}\index{n-homotopic} if for every CW complex $K$ of dimension at most $n$ and for every map $\phi: K\rightarrow X$, we have that $f\phi, g\phi:K\rightarrow Y$ are homotopic.
\end{definition}

The property of being $n$-homotopic is an equivalence relation. If $f$ and $g$ are $n$-homotopic maps from a complex $K$ into a space $X$, their restrictions to $K^n$ are homotopic by the definition. Conversely, if their restrictions to $K^n$ are homotopic, then they are $n$-homotopic by the Cellular Approximation Theorem.

\begin{definition}
Two spaces $X$ and $Y$ have the same {\it n-homotopy type}\index{n-homotopy type} if there exists maps $f: X\rightarrow Y$ and $g: Y\rightarrow X$ such that $fg$ is $n$-homotopic to $id_Y$ and $gf$ is $n$-homotopic to $id_X$. Then we say that $(f, g)$ is a $n$-homotopy equivalence and $g$ is an $n$-homotopy inverse of $f$. 
\end{definition}

The next definition will be more important. 

\begin{definition}
Two CW complexes $K$ and $L$ have the same {\it n-type}\index{n-type} if their $n$-skeletons $K^n$ and $L^n$ have the same $(n-1)$-homotopy type. In other words, $K$ and $L$ have the same $n$-type if there are maps between $K^n$ and $L^n$ whose compositions are $(n-1)$-homotopic to the identity maps of $K^n$ and $L^n$. 
\end{definition}

WARNING: $n$-type and $n$-homotopy type do not mean the same thing. Be careful not to confuse them.

Note that when $n=\infty$, $\infty$-homotopic means homotopic and $K^\infty=K$. 

\begin{theorem}
If $K$ and $L$ have the same $n$-type, then they have the same $m$-type for $m<n$. This statement holds even for $n=\infty$. If  $K$ and $L$ have the same $\infty$-type, then they have the same $m$-type for all finite $m$. So $n$-type is a homotopy invariant of a complex.
\end{theorem}

\begin{example}
Let $K$ and $L$ be the two complexes representing $S^n$ given in Example 9.8.6. Let $n=3$. Then $K^2=e^0$ (i. e. a single point) and $L^2=S^2$. Then since any map of a 1-complex into $S^2$ is nullhomotopic (i.e. homotopic to a constant map), $K^2$ and $L^2$ have the same 1-homotopy type, so $K$ and $L$ have the same 2-type. 
\end{example}

It turns out that the 1-type of a complex is a measure of the number of connected components, and Whitehead showed that two complexes with the same 2-type have the same fundamental group. Postnikov systems will characterize $n$-type for higher values of $n$.

Recall that $[X, Y]$ denotes homotopy classes of maps from $X$ to $Y$.

\begin{theorem} 
If $L_1$ and $L_2$ have the same $n$-type and the dimension of $K$ is less than or equal to $n-1$ then there is a one to one correspondence between the sets $[K, L_1]$ and $[K, L_2].$
\end{theorem}

Letting $K=S^i$ for $i<n$ gives the following.

\begin{theorem} 
If $L_1$ and $L_2$ have the same $n$-type then $\pi_i(L_1)\cong\pi_i(L_2)$ for all $i<n$.
\end{theorem}

The converse is false, but Whitehead \cite{JHCW} proves the following.

\begin{theorem} 
Suppose that there is a map $F: K^n\rightarrow L^n$ which induces isomorphisms on the homotopy groups $\pi_i(K^n)\cong\pi_i(L^n)$ for all $i<n$. Then $K$ and $L$ have the same $n$-type.
\end{theorem}

Now we are almost ready to define Postnikov systems. If $L$ is an $(n-1)$ connected complex and $\pi=\pi_n(L)$ then $L$ and $K(\pi, n)$ have the same homotopy groups through dimension $n$.

\begin{theorem} 
$L$ and $K(\pi, n)$ have the same (n+1)-type.
\end{theorem}

We need a map between the two spaces that induces the isomorphism on $\pi$ so that we can use Theorem 9.8.8. This comes from a correspondence between $[L, K(\pi, n)]$ and $H^n(L; \pi)$. This is an important result in the theory of cohomology operations. Mosher and Tangora cover it earlier in their book, but we will accept it for now and talk about it more in Chapter 11.

We would like to modify $K(\pi, n)$ in such a way as to build a space with the same $(n+2)$-type as the complex $L$. This is done through a {\it Postnikov system}.

\begin{definition}
Let $X$ be an $(n-1)$-connected CW complex with $n\geq 2$ so that $X$ is simply connected. A diagram of the form shown is a {\it Postnikov system}\index{Postnikov system} if it satisfies the following conditions:

$$\begin{tikzpicture}
  \matrix (m) [matrix of math nodes,row sep=3em,column sep=4em,minimum width=2em]
  {
& \vdots &\\
& X_{m+1} &\\
X & X_m & K(\pi_{m+1}, m+2)\\
& \vdots &\\
& X_n=K(\pi_n, n) &\\
& (\ast) &\\};

\path[-stealth]
(m-1-2) edge  (m-2-2)
(m-2-2) edge  (m-3-2)
(m-3-2) edge  (m-4-2)
(m-4-2) edge  (m-5-2)
(m-5-2) edge  (m-6-2)
(m-3-1) edge node [left] {$\rho_{m+1}$} (m-2-2)
(m-3-1) edge node [below] {$\rho_m$} (m-3-2)
(m-3-2) edge node [above] {$k_m(X)$} (m-3-3)

;

\end{tikzpicture}$$
\begin{enumerate}
\item Each $X_m$ for $m\geq n$ has the same $(m+1)$-type as $X$ and there is a map $\rho_m: X\rightarrow X_m$ inducing the isomorphisms $\pi_i(X)\cong\pi_i(X_m)$ for all $i\leq m$. 
\item $\pi_i(X_m)=0$ for $i>m$.
\item $X_{m+1}$ is the fiber space over $X_m$ induced by $k_m(X)$ from the path-space fibration over $K(\pi_{m+1}, m+2)$.
\item The diagram commutes up to homotopy.
\end{enumerate}
In the diagram, $\pi_i$ denotes $\pi_i(X)$, $k_m$ denotes a homotopy class of maps $X_m\rightarrow K(\pi_{m+1}, m+2)$ and $(\ast)$ is a one point space.
\end{definition}

The tools that Mosher and Tangora use to construct a Postnikov system are beyond what you have. (See the first 12 chapters of \cite{MT}.) In chapter 10, I will outline a method that is claimed to be computable in polynomial time using a very different technique.

\section{Spectral Sequences}

Spectral sequences are the most ugly, hideous algebraic objects ever invented. Imagine a two dimensional array whose grid points contain abelian groups or more generally modules. These objects are also constantly changing over time as we stand by helplessly hoping they converge to something we can deal with. But in algebraic topology, we need to face up to them eventually. For example, suppose we have a fiber space $F\rightarrow E\rightarrow B$. In homotopy theory, we have an exact sequence. But what about homology? If you know the homology of $F$ and $B$ can you find the homology of $E$. As we will see in Chapter 11, finding cohomology of Eilenberg MacLane spaces involves problems of this type. The Leray-Serre spectral sequence comes to the rescue here. Another spectral sequence, the Adams spectral sequence, is used in computing homotopy groups of spheres. 

The main book on spectral sequences that will tell you everything you ever wanted to know abot the topic is the book by McCleary \cite{McCl}. Mac Lane's book \cite{MacL2} has a more compact description and there is also a description in Mosher and Tangora \cite{MT} which is more focused on cohomology operations. The description in Hu \cite{Hu} seems to be missing some key pieces so I don't recommend it for this topic. Also, note that Edelsbrunner and Harer \cite{EH} have a very brief description in their TDA book, but they don't really explore how spectral sequences might be applied.

I am presenting a brief description as it is important in understanding topics such as cohomology operations and homotopy groups of speheres. Also, one of the main ways that spectral sequences arise is through filtrations so it may connect to TDA in that way. In what follows, I will mainly follow Mac Lane with help from other sources as needed. I will start with some definitions, then describe the Leray-Serre spectral sequence. Then I will discuss the two ways a spectral sequence can arise: through a filtration or through an exact couple. Finally, I will mention the Kenzo software of Dousson, et. al. \cite{Kenz}, which can automatically do computations in spectral sequences and Postnikov towers in a number of cases. I have not had the chance to experiment with this software but it could be the basis of future efforts.

\subsection{Basic Definitions}

\begin{definition} 
A {\it differential bigraded module}\index{differential bigraded module} over a ring $R$ is family of modules $E=\{E_{p, q}\}$ for integers $p$ and $q$ together with a {\it differential} which is a family of homomorphisms $d: E\rightarrow E$ of bidegree $(-r, r-1)$ and $d\cdot d=0$. The bidegree means that for all $p, q\in Z$, $$d: E_{p, q}\rightarrow E_{p-r, q+r-1}.$$ 
\end{definition} 

The condition $d\cdot d=0$ means that we can consider $E$ to be a chain complex with $d$ as a boundary map and define homology as $$H_{p, q}(E)=\frac{\ker(d: E_{p, q}\rightarrow E_{p-r, q+r-1})}{im(d: E_{p+r, q-r+1}\rightarrow E_{p, q})}.$$  

We can make $E$ into a singly graded module by letting $E=\{E_n\}$ with total degree $n$ by letting $$E_n=\sum_{p+q=n}E_{p, q}.$$ Then the differential $d$ becomes a differential $d: E_n\rightarrow E_{n-1}$ and $$H_n(E)=\sum_{p+q=n}H_{p, q}(E).$$

\begin{definition}
A {\it spectral sequence}\index{spectral sequence} $E=\{E^r, d^r\}$ is a sequence $E^2, E^3, \cdots$ of bigraded modules each with a differential $$d^r: E^r_{p, q}\rightarrow E^r_{p-r, q+r-1}$$ for $r=2, 3, \cdots$ of bidegree $(-r, r-1)$ and $E^{r+1}=H(E^r, d^r)$. In other words, $E^{r+1}$ is the bigraded homology module of the preceeding $(E^r, d^r)$. $E^r$ and $d^r$ determine $E^{r+1}$ but not necessarily $d^{r+1}$. $E^2$ is called the {\it initial term} of the spectral sequence. (Note that it is conventional for a spectral sequence to start at $r=2$ but there is no reason it couldn't start at $r=1$ instead.)
\end{definition}

In what follows, I will use $d^2$ for the differential when $r=2$ and $d\cdot d$ for $d$ composed with itself.

\begin{definition}
If $E_1, E_2$ are two spectral sequences then a {\it homomorphism} $f: E_1\rightarrow E_2$ is a family of homomorphisms of bigraded modules $$f^r: {E^r_1}_{(p, q)}\rightarrow{E^r_2}_{(p, q)}$$ such that $d^rf^r=f^rd^r$ and each $f^{r+1}$ is the map induced by $f^r$ on homology.
\end{definition}

Now we will look at the sequence in terms of submodules of $E^2$. Recall that $E^{r+1}=H(E^r, d^r)$. So $E^3=H(E^2, d^2)=C^2/B^2$ where $C^2$ is the kernel of $d^2$ and $B^2$ is the image of $d^2$. $C^2, B^2$ are then submodules of $E^2$. Similarly $E^4=H(E^3, d^3)=C^3/B^3$ where $C_3/B^2=\ker d^3$ and $B^3/B^2=im\hspace{.02 in}d^3$. Here $B^3\subset C^3$. Iterating gives a chain $$0=B^1\subset B^2\subset B^3\subset\cdots\subset C^3\subset C^2\subset C^1=E^2$$ of bigraded submodules of $E^2$ with $E^{r+1}=C^r/B^r$ and $$d^r: C^{r-1}/B^{r-1}\rightarrow C^{r-1}/B^{r-1}$$ has kernel $C^r/B^{r-1}$ and image $B^r/B^{r-1}$.

We can think of $C^{r-1}$ as the module of elements that live until stage $r$ and $B^{r-1}$ as the module of elements that bound by stage $r$. Letting $r\rightarrow \infty$, $C^\infty=\cap_{r=2}^\infty C^r$ is the module of elements that live forever and $B^\infty=\cup_{r=2}^\infty B^r$ is the module of elements that eventually bound.

Now $B^\infty\subset C^\infty$ so the spectral sequence determines a bigraded module $$E^\infty=\{E^\infty_{p, q}=C^\infty_{p, q}/B^\infty_{p, q}\}.$$ Generally $E^\infty$ will be what you are trying to find. Also note that a homomorphism $f: E_1\rightarrow E_2$ will induce homomorphisms $f^r: E_1^r\rightarrow E_2^r$ for $2\leq r\leq\infty$. 

In homology calculations we often use a {\it first quadrant spectral sequence} in which $E_{p, q}^r=0$ if $p<0$ or $q<0$. The modules $E_{p, q}$ are displayed at the grid points $(p, q)$ of the first quadrant of the $p-q$ plane as shown in Figure 9.9.1.

\begin{figure}[ht]
\begin{center}
  \scalebox{0.8}{\includegraphics{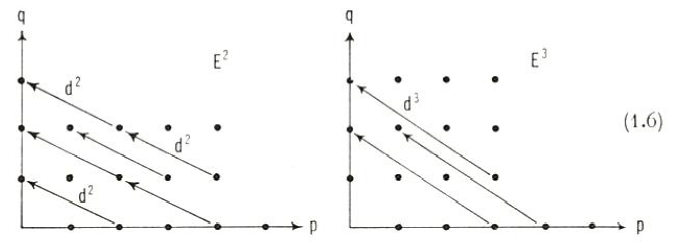}}
\caption{
\rm
First Quadrant Spectral Sequence \cite{MacL2}. 
}
\end{center}
\end{figure}

The differential $d^r$ is marked by an arrow. The terms of total degree $n$ lie on a line of slope -1, and the differentials go from the line to one point below it. At the grid point $(p, q)$, the value of $E^{r+1}_{p, q}$ is the kernel of the arrow starting at $E^r_{p, q}$ modulo the image of the arrow that ends there. There is no change if the out going arrow ends outside the first quadrant (if $r>p$) or the incoming arrow starts outside of it ($r>q+1)$. So for $\max(p, q+1)<r<\infty$, $E^{r+1}_{p, q}=E^r_{p, q}$. So for fixed, $p$ and $q$, $E^r_{p, q}$ is eventually constant for large enough $r$ as $r$ increases.

The terms $E_{p, 0}$ on the $p$-axis are called the {\it base} terms. Each arrow $d^r$ ending on the base comes from below so it is zero. So each $E^{r+1}_{p, 0}$ is a submodule of $E^r_{p, 0}$ equal to the kernel of $d^r: E^r_{p, 0}\rightarrow E^r_{p-r, r-1}.$ This gives a sequence of monomorphisms $$E_{p, 0}^\infty=E_{p, 0}^{p+1}\rightarrow\cdots\rightarrow E_{p, 0}^4\rightarrow E_{p, 0}^3\rightarrow E_{p, 0}^2.$$

The terms $E_{0, q}$ on the $q$-axis are called the {\it fiber} terms. Each arrow from a fiber term ends at zero. So the kernel of $d^r$ us all of $E^r_{0, q}$, and $E^{r+1}_{0, q}$ is the quotient of $E^r_{0, q}$ by the image of $d^r$. This gives a sequence of epimorphisms $$E^2_{0, q}\rightarrow E^3_{0, q}\rightarrow E^4_{0, q}\rightarrow\cdots E^{q+2}_{0, q}=E^\infty_{0, q}.$$ The two sequences are known as {\it edge homomorphisms}. The terms {\it base} and {\it fiber} should sound familiar. You will see where they come from in the next subsection.

\begin{definition}
A spectral sequence is {\it bounded below} if for each degree $n$, there is a number $s$ dependent on $n$ such that $E^2_{p, q}=0$ when $p<s$ and $p+q=n$. This means that on a line of constant degree,the terms are eventually zero as $p$ decreases. This is certainly true in a first quadrant sequence.
\end{definition}

\begin{theorem}
{\bf Mapping Theorem:} If $f: E_1\rightarrow E_2$ is a homomorphism of spectral sequences and if $f^t: E_1^t\rightarrow E_2^t$ is an isomorphism for some $t$, then $f^r: E_1^r\rightarrow E_2^r$ is an isomorphism for $r\geq t$. If in addition, $E_1$ and $E_2$ are bounded below, $f^\infty: E_1^\infty\rightarrow E_2^\infty$ is an isomorphism.
\end{theorem}

\subsection{Leray-Serre Spectral Sequence}

There are a number of famous spectral sequences (see \cite{McCl} for several examples), but I will focus on one of them. The Leray-Serre spectral sequence is the main tool for computing homology of fiber spaces. If you can read French, J-P Serre's Ph.D. thesis \cite{Ser} has a great explanation and the proofs of the theorems I will state. I could have only dreamed of having a thesis that influential. It was considered so importnat that it was published in its entirety. 

Consider a fiber space $F\rightarrow E\rightarrow B$. We will assume that $B$ is pathwise and simply connected. Then any fiber has the same homology groups, so we can talk about the groups $H_p(B; H_q(F))$. This is the (singular) homology of $B$ with coefficients in the group $H_q(F)$. If $p=0$, then $H_0(B; H_q(F))\cong H_q(F).$ 

\begin{theorem}
{\bf Leray-Serre:} If $f: E\rightarrow B$ is a fiber map with base $B$ pathwise and simply connected, and fiber $F$ pathwise connected, then for each $n$, there is a nested family of subgroups of the singular homology group $H_n(E)$, $$0=H_{-1, n+1}\subset H_{0, n}\subset H_{1, n-1}\subset \cdots\subset H_{n-1, 1}\subset H_{n, 0}=H_n(E),$$ and a first quadrant spectral sequence such that $$E^2_{p, q}\cong H_p(B; H_q(F)),$$ and $$E^\infty_{p, q}\cong H_{p, q}/H_{p-1, q+1}.$$ 

If $e_B$ is the composite edge homomorphism on the base (defined in Section 9.9.1), the composite  $$H_p(E)=H_{p, 0}\rightarrow H_{p, 0}/H_{p-1, 1}\cong E^\infty_{p, 0}\xrightarrow{e_B}E^2_{p, 0}\cong H_p(B; H_0(F))\cong H_p(B)$$ is the homomorphism induced on homology by the fiber map $f: E\rightarrow B$. 

If $e_F$ is the composite edge homomorphism on the fiber, the composite $$H_q(F)\cong H_0(B; H_q(F))\cong E^2_{0, q}\xrightarrow{e_F} E^\infty_{0, q}\rightarrow H_q(E)$$ is the homomorphism induced on homology by the inclusion $F\subset E$.
\end{theorem}

The sequence relates the homology of the base and the fiber to that of the total space $E$. $E^\infty$ gives successive factor groups in the filtration of the homology of $E$.

Now recall the Universal Coefficient Theorem for Homology (Theorem 8.4.3) states that the sequence $$0\rightarrow H_p(B)\otimes H_q(F)\rightarrow E^2_{p, q}=H_p(B; H_q(F))\rightarrow Tor(H_{p-1}(B), H_q(F))\rightarrow 0$$ is exact. If the $H_{p-1}(B)$ are all torsion free then $E^2_{p, q}\cong H_p(B)\otimes H_q(F).$

Mac Lane gives three examples. I will describe them in detail. They are instructive in seeing how a spectral sequence works.

\begin{theorem}
{\bf The Wang sequence:\index{Wang sequence}} If $f: E\rightarrow S^k$ is a fiber space for $k\geq 2$, and the fiber $F$ is pathwise connected, then there is an exact sequence $$\cdots\rightarrow H_n(E)\rightarrow H_{n-k}(F)\xrightarrow{d^k} H_{n-1}(F)\rightarrow H_{n-1}(E)\rightarrow\cdots .$$
\end{theorem}

{\bf Proof:} The base space $S^k$ is simply connected and has homology $H_p(S^k)\cong Z$ if $p=0, k$ and $H_p(S^k)=0$ otherwise. Then we have $E^2_{p, q}=H_p(S^k; H_q(F))\cong H_p(B)\otimes H_q(F),$ since all of the homology groups of $S^k$ are torsion free. Also, since $Z\otimes G\cong G$ for any abelian group $G$, we have $E^2_{p, q}\cong H_q(F)$ for $p=0, k$ and $E^2_{p, q}=0$ otherwise. The nonzero terms of $E^2_{p, q}$ are all on the vertical lines $p=0$ and $p=k$, so the only nonzero differential $d^r$ for $r\geq 2$ is $d^k$. (Recall that $d^k$ lowers $p$ by $k$ and raises $q$ by $k-1$.) This means that $E^2=E^3=\cdots=E^k$ and $E^{k+1}=E^{k+2}=\cdots=E^\infty$. Since $E^{k+1}=E^\infty$ is the homology of $(E^k, d^k)$ we have the exact sequence $$0\rightarrow E^\infty_{k, q}\rightarrow E^k_{k, q}\xrightarrow{d^k} E^k_{0, q+k-1}\rightarrow  E^\infty_{0, q+k-1}\rightarrow 0.$$ 

Now the filtration for $H_n(E)$ has only two nonzero quotient modules so it collapses to $$0\subset H_{0, n}=H_{k-1, n-k+1}\subset H_{k, n-k}=H_n.$$ Since by definition $E^\infty_{p, q}\cong H_{p, q}/H_{p-1, q+1},$ we get a short exact sequence $$0\rightarrow E^\infty_{0, n}\rightarrow H_n(E)\rightarrow E^\infty_{k, n-k}\rightarrow 0.$$ Now splice this exact sequence and the previous one together substituting $q=n-k$, $E^k_{k, q}=E^2_{k, q}\cong H_q(F)$, and $E^k_{0, q}=E^2_{0, q}\cong H_q(F)$. We then get: 

$$\begin{tikzpicture}
  \matrix (m) [matrix of math nodes,row sep=3em,column sep=4em,minimum width=2em]
  {
& H_n(E) & & & 0 &\\
0 & E^\infty_{k, n-k} & H_{n-k}(F) & H_{n-1}(F) & E^\infty_{0, n-1} & 0\\
& 0 & & & H_{n-1}(E) &\\};
  \path[-stealth]
(m-2-1) edge  (m-2-2)
(m-2-2) edge  (m-2-3)
(m-2-3) edge  node[above] {$d^k$}(m-2-4)
(m-2-4) edge  (m-2-5)
(m-2-5) edge  (m-2-6)

(m-1-2) edge  (m-2-2)
(m-2-2) edge  (m-3-2)
(m-1-5) edge  (m-2-5)
(m-2-5) edge  (m-3-5)
(m-1-2) edge[dashed] (m-2-3)
(m-2-4) edge[dashed] (m-3-5)
;

\end{tikzpicture}$$

Following the dashed lines with $d^k$ in the middle gives the desired result. $\blacksquare$

For the next example, we will compute the homology of the loop space of a sphere.

Notation warning: Recall that Hu \cite{Hu} uses the symbol $\Omega$ for a path space and $\Lambda$ for a loop space. Mac Lane \cite{MacL2} and most other books use $\Omega$ for a loop space. I will switch notation here to agree with him. From now on, I will specify what the symbol represents. In this section $\Omega X$ will mean the space of loops on $X$. Mac Lane calls the space of paths on $X$, $L(X)$. 

Now if $B$ is connected and simply connected, we have a fiber space $p: L(B)\rightarrow B$. Let $b_0$ be a base point and $L(B)$ the set of paths that end at $b_0$. Then $p$ takes a path in $B$ to it's initial point and $p^{-1}(b_0)$ is the space of loops starting and ending at $b_0$. We let $B$ be pathwise and simply connected so we can ignore the base point. Then the fiber is $F=\Omega B$. Also recall that $L(B)$ is contractible. 

We will look at the homology of the space of loops of the sphere $S^k$. You may find this useful if $k=2$ and you find yourself going around in circles.

\begin{theorem}
The loop space $\Omega S^k$ of the sphere $S^k$ with $k>1$ has homology $H_n(\Omega S^k)\cong Z$ if $n\equiv 0$ (mod $k-1$) and $H_n(\Omega S^k)=0$ otherwise for $n\geq 0$.
\end{theorem}

{\bf Proof:} $S^k$ for $k>1$ is simply connected so each loop can be contracted to the zero loop. Thus $\Omega S^k$ is path connected so that $H_0(\Omega S^k)\cong Z$. Now $E=L(B)$ is contractible so $E$ is acyclic. So every third term in the Wang sequence is zero except $H_0(E)\cong Z$. So the sequence gives isomorphisms $H_{n-k}(\Omega S^k)\cong H_{n-1}(\Omega S^k)$ and the result follows from the initial value at $H_0(\Omega S^k).$ $\blacksquare$

\begin{theorem}
{\bf: The Gysin sequence:\index{Gysin sequence}} If $p: E\rightarrow B$ is a fiber space with simply connected base $B$ and with fiber $F=S^k$ with $k\geq 1$, then there is an exact sequence $$\cdots\rightarrow H_n(E)\xrightarrow{p_\ast}H_n(B)\xrightarrow{d^{k+1}}H_{n-k-1}(B)\rightarrow H_{n-1}(E)\rightarrow\cdots .$$
\end{theorem}

{\bf Proof:} Since $H_q(F)=H_q(S^k)$ is only nonzero if $q=0, k$, we have $E^2_{p, q}\cong H_p(B)$ if $q=0$ or $q=k$, and $E^2_{p, q}=0$ otherwise. The spectral sequence then lies on the two horizontal lines $q=0$ and $q=k$, and the only nonzero differential is $d^{k+1}$. We then get the exact sequences $$0\rightarrow E^\infty_{n, 0}\rightarrow   E^2_{n, 0}\xrightarrow{d^{k+1}}  E^2_{n-k-1, k}\rightarrow  E^\infty_{n-k-1, k}\rightarrow  0$$ and $$0\rightarrow  E^\infty_{n-k-1, k+1}\rightarrow H_n(E)\rightarrow E^\infty_{n, 0}\rightarrow 0,$$ and we splice them together as in Theorem 9.9.3 to get our result. $\blacksquare$

\subsection{Filtrations and Exact Couples}

I will now discuss two situations that can give rise to spectral sequences. 

In the first situation we have a {\it filtration}. A filtration $F$ of a module $A$ is a family of submodules ${F_pA}$ for $p\in Z$ such that $$\cdots\subset F_{p-1}A\subset F_pA\subset F_{p+1}A\subset\cdots.$$

Each filtration $F$ of $A$ determines a {\it associated graded module} $G^FA=\{(G^FA)_p=F_pA/F_{p-1}A\}$, consisting of the successive factor modules in the chain. 

If $F$ and $F'$ are filtrations of $A$ and $A'$ respectively, then a homomorphism $F: A\rightarrow A'$ of filtered modules is a module homomorphism with $f(F_pA)\subset F'_pA'$. The filtration of a graded module $A$ with a differential and homomorphisms preserving the grading induces a filtration of the homology module $H(A)$ with $F_p(H(A))$ defined as the image of $H(F_pA)$ under the injection $F_pA\rightarrow A$. Since $A$ itself is graded, the filtration $F$ of $A$ determines a filtration $F_pA_n$ of each $A_n$ and the differential of $A$ induces homomorphisms $\partial  F_pA_n\rightarrow F_pA_{n-1}$ for each $p$ and each $n$. The family $\{F_pA_n\}$ is a bigraded module. Letting $q=n-p$ we call $p$ the {\it filtration degree} and $q$ the {\it complementary degree}. The bigraded module then takes the  form $\{F_pA_{p+q}\}$. We call this a $FDG_Z$ module as an abbreviation of a {\it filtered differential Z-graded module.} 

A filtration $F$ of a differential graded module $A$ is said to be {\it bounded} if there are integers $s, t$ dependent on $n$ such that $s<t$, $F_sA_n=0$, and $F_tA_n=A_n$. This makes the filtration of finite length and it takes the form $$0=F_sA_n\subset F_{s+1}A_n\subset\cdots\subset F_tA_n=A_n.$$

A spectral sequence $\{E^r_p, d^r\}$ is said to {\it converge} to a graded module $H$ if there is a filtration $F$ of $H$ and isomorphisms $E^\infty\cong F_pH/F_{p-1}H$ of graded modules for each $p$. To put the spectral sequence in the more familiar bigraded form, let $E^r_p$ be the graded module $\{E^r_{p, q}\}$ graded by the complementary degree $q$. 

\begin{theorem}
Each filtration $F$ of a differential graded module $A$ determines a spectral sequence $(E^r, d^r)$ for $r=1, 2, \cdots$ which is a covariant functor of $(F, A)$ together with the isomorphisms $E^1_{p, q}\cong H_{p+q}(F_pA/F_{p-1}A).$ If $F$ is bounded, then the sequence converges to $H(A)$, i.e. $E^\infty_{p, q}\cong F_p(H_{p+q}A)/F_{p-1}(H_{p+q} A).$
\end{theorem}

Question to think about: Persistent homology starts with a filtration of a differential graded module. (See Sections 5.1 and the more general version in Section 5.2.) Would building the associated spectral sequence give you interesting information about your data? As I mentioned above, Edelsbrunner and Harer mention spectral sequences in \cite{EH} but to the best of my knowledge, the meaning of a spectral sequence derived from a point cloud has not been fully explored.

See \cite{MacL2} for the rather lengthy proof of Theorem 9.9.6 and \cite{McCl} for a more complete discussion of the relationship between a filtration and the associated spectral sequence.

An alternate way of deriving spectral sequences is through {\it exact couples}. 

\begin{definition}
An {\it exact couple}\index{exact couple} $\mathfrak{C}=\{D, E, i, j, k\}$ consists of two modules $D$ and $E$ together with homomorphisms $i$, $j$, and $k$ such that the triangle below is exact:$$\begin{tikzpicture}
  \matrix (m) [matrix of math nodes,row sep=3em,column sep=4em,minimum width=2em]
  {
D &  & D\\
& E &\\};
  \path[-stealth]
(m-1-1) edge  node[above] {$i$}(m-1-3)
(m-2-2) edge  node[above] {$k$}(m-1-1)
(m-1-3) edge  node[above] {$j$}(m-2-2)

;

\end{tikzpicture}$$

The modules $D$ and $E$ may be graded or bigraded.
\end{definition}

Since the triangle is exact, the composite $jk: E\rightarrow E$ has square $(jk)(jk)=j(kj)k=j0k=0$, since $kj=0$ by exactness. So $jk$ is a differntial on $E$. Let $E'=H(E, jk)$ and $D'=i(D)$.  Let $i'=i|D'$, $j'(id)=j(d)+jkE$, and $k'(e+jkE)=ke$ for $d\in D$ and $e\in E$ where $jk(e)=0$. 

Now $id=0$ implies $d\in kE$ by exactnesss, so $jd\in jkE$ and $j'$ is well defined. Also, $jke=0$ implies that $ke\in iD$ so $k'$ is well defined.

It is easy to show that that $\mathfrak{C}'=\{D', E', i', j', k'\}$ forms and exact triangle. (Try it.) $\mathfrak{C'}$ is then an exact couple called the {\it derived couple} of $\mathfrak{C}$. Iterating the process produces the {\it rth-derived couple} $\mathfrak{C}^r=\{D^r, E^r, i^r, j^r, k^r\}$.

 Let $\mathfrak{C}^1=\mathfrak{C}$, $\mathfrak{C}^2=\mathfrak{C}'$, and $\mathfrak{C}^{r+1}=(\mathfrak{C}^r)'$. Then we have the following:

\begin{theorem}
An exact couple of bigraded modules $D$ and $E$ with maps of bidegrees deg$(i)=(1, -1)$, deg$(j)=(0, 0)$, and deg$(k)=(-1, 0)$ determines a spectral sequence $(E^r, d^r)$ with $d^r=j^rk^r$.
\end{theorem}

For the exact couple $\mathfrak(C)^r$, Mac Lane shows that we have bidegrees deg$(i^r)=(1, -1)$, deg$(j^r)=(-r+1, r-1)$, and deg$(k^r)=(-1, 0)$, so that deg$(j^rk^r)=(-r, r-1)$. Thus, each $E^{r+1}$ is the homology of $E^r$  with the appropriate bidegree for a spectral sequence. The first stage of the sequence can be expressed as shown: 

$$\begin{tikzpicture}
  \matrix (m) [matrix of math nodes,row sep=3em,column sep=4em,minimum width=2em]
  {
& & \vdots & & \vdots &\\
\cdots & E_{p, q+1} & D_{p-1, q+1} & E_{p-1, q+1} & D_{p-2, q+1} & \cdots\\
\cdots & E_{p+1, q} & D_{p, q} & E_{p, q} & D_{p-1, q} & \cdots\\
\cdots & E_{p+2, q-1} & D_{p+1, q-1} & E_{p+1, q-1} & D_{p, q-1} & \cdots\\
& & \vdots & & \vdots &\\};
  \path[-stealth]
(m-2-1) edge  node[above] {$j$} (m-2-2)
(m-2-2) edge  node[above] {$k$} (m-2-3)
(m-2-3) edge  node[above] {$j$}(m-2-4)
(m-2-4) edge  node[above] {$k$} (m-2-5)
(m-2-5) edge  node[above] {$j$}(m-2-6)
(m-3-1) edge  node[above] {$j$} (m-3-2)
(m-3-2) edge  node[above] {$k$} (m-3-3)
(m-3-3) edge  node[above] {$j$}(m-3-4)
(m-3-4) edge  node[above] {$k$} (m-3-5)
(m-3-5) edge  node[above] {$j$}(m-3-6)
(m-4-1) edge  node[above] {$j$} (m-4-2)
(m-4-2) edge  node[above] {$k$} (m-4-3)
(m-4-3) edge  node[above] {$j$}(m-4-4)
(m-4-4) edge  node[above] {$k$} (m-4-5)
(m-4-5) edge  node[above] {$j$}(m-4-6)

(m-1-3) edge  node[right] {$i$} (m-2-3)
(m-2-3) edge  node[right] {$i$} (m-3-3)
(m-1-5) edge  node[right] {$i$} (m-2-5)
(m-2-5) edge  node[right] {$i$} (m-3-5)
(m-3-3) edge  node[right] {$i$} (m-4-3)
(m-4-3) edge  node[right] {$i$} (m-5-3)
(m-3-5) edge  node[right] {$i$} (m-4-5)
(m-4-5) edge  node[right] {$i$} (m-5-5)

;

\end{tikzpicture}$$

A filtration of a graded differential module $A$ determines an exact couple as follows. The short exact sequence of complexes $$0\rightarrow F_{p-1}A\rightarrow F_pA\rightarrow F_pA/F_{p-1}A\rightarrow 0$$ yields the exact homology sequence $$\cdots \rightarrow H_n(F_{p-1}A)\xrightarrow{i}H_n(F_pA)\xrightarrow{j}H_n(F_pA/F_{p-1}A)\xrightarrow{k}H_{n-1}(F_{p-1}A)\rightarrow \cdots,$$ where $i$ is injection, $j$ is projection, and $k$ is the homology boundary homomorphism. These sequences give an exact couple with $D_{p, q}=H_{p+q}(F_pA)$ and $E_{p, q}=H_{p+q}(F_pA/F_{p-1}A)$ and the degrees of $i$, $j$, and $k$ as in Theorem 9.9.7. This is called the {\it exact couple of the filtration} $F$.

\begin{theorem}
The spectral sequence of $F$ is isomorphic to that of the exact couple of $F$. 
\end{theorem}

See Mac Lane \cite{MacL2} for the proof. Also, filtrations give rise to exact couples but exact couples do not have to come from a filtration. See Mac Lane for an example.

\subsection{Kenzo}

Kenzo is a software package for computational homological algebra and algebraic topology developed by Francis Sergeraert of Grenoble Alpes University in Grenoble, France. The package was written in LISP and originally called EAT (effective algebraic topology) and then CAT (constructive algebraic topology.) The latter acronym inspired Sergeraert to name the package after his cat, Kenzo. (See Figure 9.9.2). Kenzo was written in LISP and has an extensive website \cite{Kenz} with examples and lots of documentation. 

\begin{figure}[ht]
\begin{center}
  \scalebox{0.8}{\includegraphics{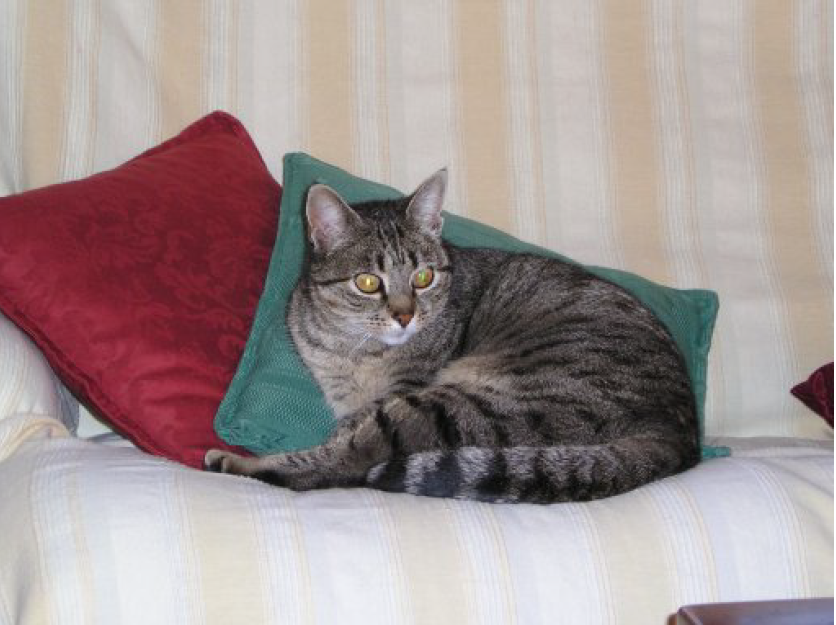}}
\caption{
\rm
Kenzo the Cat \cite{Kenz}. 
}
\end{center}
\end{figure}

Kenzo is not necessary for computing homology and cohomology groups, even though you can use it that way. Rather, it is your best hope for computations of objects such as spectral sequences and Postnikov systems. It is also the basis for the obstruction theory algorithms I will mention in the next chapter. The paper that brought Kenzo to my attention was the introduction to the 2006 Summer School in Genova (Genoa) Italy \cite{RS} written by Sergeraert and Julio Rubio of University of La Rioja, Spain. Rubio's Ph.D student, Ana Romero, is a joint author on a paper that specifically deals with what type of spectral sequences can be computed with Kenzo. See \cite{RRS} for details.

I will spend some time on some key concepts. Although you can prove some pretty cool results with exact sequences and spectral sequences, neither one is actually an algorithm. There are plenty of times where they introduce a lot of uncertainty and the fact that you can set up your problem that way does not guarantee a solution. We want an algorithm that is {\it constructive}, and that is the type of problem that Kenzo can handle.

Rubio and Sergeraert \cite{RS} discuss the notion of constructive mathematics. There are proofs that an object exists that don't actually produce the object. Here is one of my favorites. 

\begin{definition}
If $X$ is a topological space then a {\it sequence}\index{sequence} in $X$ is a function from the nonnegative integers to $X$. If we relax this condition to allow for the domain to be any linearly ordered set, the resulting function is called a {\it net}\index{net}. An {\it ultranet}\index{ultranet} $U$  is a net such that for any subset $A\subset X$, $U$ is eventually in $A$ or $U$ is eventually not in $A$. (By this we mean that if the ultranet is a map from an indexing set $J$ to $X$, and $x_\alpha$ is the point in $X$ that is the image of $\alpha\in J$, then there exists $\beta\in J$, such that $x_\gamma\in A$ for all $\gamma>\beta$ or $x_\gamma\notin A$ for all $\gamma>\beta$. Which of these conditions holds depends on the particular subset $A$.)
\end{definition}

\begin{example}
Let a net (or even a sequence) have the property that for all $\gamma>\beta$, $x_\gamma=x_0$ where $x_0$ is a fixed point in $X$ and $\beta\in J$ is fixed. Then for a given choice of  $A\subset X$,  the net is eventually in $A$ if $x_0\in A$ and eventually not in $A$ if $x_0\notin A.$ This type of ultranet is called a {\it trivial ultranet}.\index{trivial ultranet} 
\end{example}

\begin{theorem}
There exists a non-trivial ultranet.
\end{theorem}

I will refer you to Willard \cite{Wil} for example if you want to see the proof. He uses a one to one correspondence between nets and a special collection of subsets of a space called a {\it filter}. The counterpart to ultranets, ultrafilters, are a lot easier to construct. These imply the existence of nontrivial ultranets, but to my knowledge, nobody has ever constructed one. 

Part of the idea of constructive mathematics, then, is that that there is a difference between saying the converse of a proposition is false and that the proposition is true. Here is an example from \cite{RS}.

\begin{example}
Consider these two statements. Do they mean the same thing?\begin{enumerate}
\item It is false that there is no book about constructive analysis in the library.
\item The upper shelf to the left of the east window at the second floor of the library has a book about constructive analysis.
\end{enumerate}
\end{example}

Both statements imply that such a book exists, but only the second one enables it to be found. Romero and Sergeraert suggest the book by Bishop and Bridges \cite{BiBr} as a good example. 

Here is another example from \cite{RS}.

\begin{example}
Are there two irrational numbers $a$ and $b$ such that $a^b$ is rational? 

{\bf Solution 1:} Let $c=\sqrt{2}^{\sqrt{2}}$.  If $c$ is rational then $a=b=\sqrt{2}$ is a solution. If $c$ is irrational, then let $a=c$ and $b=\sqrt{2}$. Then $$a^b=(\sqrt{2}^{\sqrt{2}})^{\sqrt{2}}=\sqrt{2}^{(\sqrt{2}\sqrt{2})}=\sqrt{2}^2=2.$$ 

So either $(a, b)=(\sqrt{2}, \sqrt{2})$ is a solution or $(a, b)=(\sqrt{2}^{\sqrt{2}}, \sqrt{2})$ is a solution.

{\bf Solution 2:} Let $a=\sqrt{2}$ and $b=2\log_23.$ Both $a$ and $b$ are known to be irrational and $$a^b=\sqrt{2}^{2\log_23}=2^{\frac{1}{2}\cdot 2\log_23}=2^{\log_23}=3.$$ So $(a, b)=(\sqrt{2}, 2\log_23)$ is a solution.

Solution 1 is not a constructive solution. It gives two answers and one of them must be right but we don't know which one. Solution 2 gives an actual answer.
\end{example}

So what is a constructive algorithm in homological algebra? An example of a failure is J.P. Serre's attempt to determine $\pi_6(S^3)$ through building fibrations and using his spectral sequence. I will omit the details but in the end he got the exact sequence $$0\leftarrow Z_6\leftarrow \pi_6(S^3)\leftarrow Z_2\leftarrow 0.$$ This shows that $\pi_6(S^3)$ has 12 elements but he couldn't decide between $Z_{12}$ (later proved to be the correct answer) or $Z_6\oplus Z_2$. The problem was that the homology group from the spectral sequence corresponding to $Z_6$ did not have an explicit isomorphism exhibited between corresponding cycles and elements of $Z_6$. This would have solved the problem.

The constructiveness requirement in homological algorithm is stated in \cite{RS} as follows:

\begin{definition} 
Let $R$ be a ground ring and $C$ a chain complex of $R$-modules. A {\it solution} $S$ of the homological problem for $C$ is a set $S=\{\sigma_i\}_{1\leq i\leq 5}$ of five algorithms: \begin{enumerate}
\item $\sigma_1: C\rightarrow\{True, False\}$ is a predicate for deciding for every $n\in Z$ and every $n$-chain $c\in C_n$ whether $c$ is a cycle. 
\item $\sigma_2: Z\rightarrow\{R$-modules$\}$ associates to every integer $n$ an $R$-module $\sigma_2(n)$ which is isomorphic to $H_n(C)$.
\item The algorithm $\sigma_3$ is indexed by $Z$. For each $n\in Z$, $\sigma_{3, n}: \sigma_2(n)\rightarrow Z_n(C)$ associates to every $n$-homology class $h\in \sigma_2(n)$ a cycle $\sigma_{3, n}(h)$ representing this homology class.
\item The algorithm $\sigma_4$ is indexed by $Z$. For each $n\in Z$, $\sigma_{4, n}: C_n\supset Z_n(C)\rightarrow \sigma_2(n)$ associates to every $n$-cycle $z\in Z_n(C)$ the homology class of $z$ coded as an element of $\sigma_2(n)$.
\item The algorithm $\sigma_5$ is indexed by $Z$. For each $n\in Z$, $\sigma_{5, n}: Z_n(C)\rightarrow C_{n+1}$ associates to every $n$-cycle $z\in Z_n(C)$ that $\sigma_4$ determined to be a boundary, a chain $c\in C_{n+1}$ such that $dc=z$. 
\end{enumerate}
\end{definition} 

Another concept is a {\it locally effective} versus an {\it effective} object. For a locally effective object, we can perform calculations without storing the entire object in memory, just as a pocket calculator can add integers without storing all of $Z$ in memory. For some, algorithms, though, we do need to know everything about the object. If we want to compute $H_n(C)$, we need a full knowledge of $H_k(C)$ for $c=n+1, n, n-1$ and full knowledge of the differentials $d_k$ and $d_{k+1}$ so that we can compute $\ker d_k$ and image $d_{k+1}$.

The last concept I will define is a {\it reduction}. It allows us to perform calculations in a simpler complex that is equivalent in some way to a more complicated one. 

\begin{definition}
A reduction $\rho: \hat{C}\Rightarrow C$ is a diagram: 
$$\begin{tikzpicture}
  \matrix (m) [matrix of math nodes,row sep=3em,column sep=4em,minimum width=2em]
  {
\hat{C} & C\\};
  \path[-stealth]
(m-1-2) edge [bend right]  node[above] {$g$} (m-1-1)
(m-1-1) edge [bend right] node[below] {$f$} (m-1-2)
(m-1-1) edge [loop left]  node[left] {$h$} (m-1-1)

;

\end{tikzpicture}$$
where: \begin{enumerate}
\item $\hat{C}$ and $C$ are chain complexes.
\item $f$ and $g$ are chain-complex morphisms.
\item $h$ is a chain homotopy (degree $+1).$
\item These relations are satisfied:
\begin{enumerate}
\item $fg=id_c.$
\item $gf+dh+hd=id_{\hat{C}}.$
\item $fh=hg=hh=0.$
\end{enumerate}
\end{enumerate}
\end{definition}

\begin{theorem}
Let $\rho: \hat{C}\Rightarrow C$ be a reduction. Then $\rho$ is equivalent to a decomposition: $$\hat{C}=A\oplus B\oplus C'$$ such that:
\begin{enumerate}
\item $ C'=\mbox{im }g$ is a subcomplex of $\hat{C}$.
\item $A\oplus B=\ker f$ is a subcomplex of $\hat{C}$.
\item $A=\ker f\cap\ker h$ is not in general a subcomplex of $\hat{C}$.
\item $B=\ker f\cap\ker d$ is a a subcomplex of $\hat{C}$ with null differentials. 
\item The chain morphisms $f$ and $g$ are inverse isomorphisms between $C'$ and $C$. 
\item The arrows $d$ and $h$ are module homomorphisms of degrees -1 and 1 between $A$ and $B$.
\end{enumerate}
\end{theorem}

Figure 9.9.3 shows all of the reduction maps in a single picture. 

\begin{figure}[ht]
\begin{center}
  \scalebox{0.8}{\includegraphics{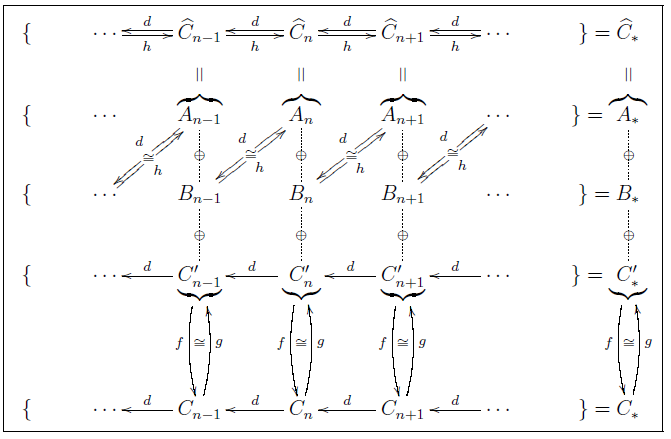}}
\caption{
\rm
Maps Involved in a Reduction \cite{RS}. 
}
\end{center}
\end{figure}

In the next chapter we will need to determine if maps between certain simplicial complexes and spheres when restricted to the boundaries of simplices are zero elements in the resulting  homotopy groups. Building on the theory of constructive algebraic topology, this seems to be possible in a lot of cases. I will make this explicit next.

\chapter{Obstruction Theory}

Back in 1989, one of the first chapters I wrote on my dissertation heavily relied on obstruction theory. My advisor said that it was a good start but I would have to add to it or people would call me an "obstructionist". Thirty years later, I was telling anyone who would listen that I thought obstruction theory would have an interesting application in data science. TDA deals with filtrations and these are really simplicial complexes that are growing. What if there was an interesting function from these complexes to some space $Y$? At each stage in adding a new simplex, we could ask if this function could be extended continuously to the larger complex. Would that tell you something interesting about your data that you didn't know before? Maybe you could find something interesting about the point cloud itself. All you would have to do is compute cohomology groups of the complex with coefficients in the homotopy groups of the target space. Most people could only nod their heads at this one.

Then with COVID and the lockdown, I suddenly had a lot more free time. I decided I would review my homotopy theory (pretty rusty after 30 years) and attack this problem once and for all. As of the time of this writing (July 2021) I don't have proof of an application that will change the world. I can argue, though, that this problem is possible due to recent advances in computational topology such as the ones I mentioned at the end of the last chapter. If you have read the first 9 chapters of this book, I can now explain to you what I mean.

So the idea is this: Let $X$ be a point cloud. and $VR_\epsilon(X)$ be the Vietrois-Rips complex associated with $X$ at radius $\epsilon$. For $\delta>\epsilon$, $VR_\epsilon(X)\subset VR_\delta(X)$. Let $f: X\rightarrow Y$ be an "interesting" function and let $Y$ be an "interesting" space. As we will see, we need to use homotopy groups of $Y$, so we would like $Y$ to be a low dimensional sphere $S^k$ where $k\geq 2$. (Remember that life could get unpleasant if $Y$ is not simply connected.) Suppose we had a continuous extension $g: VR_\epsilon(X)\rightarrow Y$ that agreed with $f$ on the vertices of $VR_\epsilon(X)$ , i.e. on the original point cloud. Could we extend $g$ to a continuous function $\tilde{g}: VR_\delta(X)\rightarrow Y$ such that $\tilde{g}=g$ on $VR_\epsilon(X)$? 

You should recognize this as the {\it extension problem}. To attack it, consider the two prototype examples. Recall that if $f: A\rightarrow Y$ is a constant function of the form $f(A)=y_0\in Y$, we can extend $f$ continuously to $g; X\rightarrow y_0\in Y$ for any $X\supset A$. This also works if $f$ is only homotopic to a constant map. If $A=S^n$, then this means that $f$ represents the zero element of $\pi_n(Y)$. 

At the other extreme, we showed in Theorem 4.1.29 that if $f: S^n\rightarrow S^n$ has a nonzero degree, then it does not extend to a continuous map $h: B^{n+1}\rightarrow S^n$. This holds in particular if $f$ is the identity on $S^n$. Since everything holds up to homotopy we can say that $f$ extends to $h$ if and only if $[f]=0$ in $\pi_n(S^n)$.

The idea then behind obstruction theory is this: Suppose $K$ is a complex and $L$ is a subcomplex. Let $f: L\rightarrow Y$. Extend $f$ over the vertices, then the edges, then the 2-simplices of $K\setminus L$. Suppose we have extended $f$ to the $n$-simplices of $K\setminus L$. We want to extend it over the $(n+1)$-simplices. If $\sigma$ is an $(n+1)$-simplex, then $f$ is already defined on the boundary $\partial\sigma\sim S^n$. So the restriction of $f$ to the boundary of $\sigma$ is map $S^n\rightarrow Y$. We can extend $f$ to $\sigma$ if and only if $f|\partial\sigma$ represents the zero element of $\pi_n(Y)$. Since we now have an element of $\pi_n(Y)$ for every $(n+1)$-simplex, the result is a cohain in $H^{n+1}(K, L; \pi_n(Y))$ known as the obstruction cochain. We can extend $f$ to the $(n+1)$-simplices if and only if this cochain is zero. 

We can tell if there is an extension to each individual simplex $\sigma$ if we can decide whether $f|\partial\sigma$ represents the zero element of $\pi_n(Y)$. How could we do this and would we live long enough to determine the answer. Francis Sergeraert comes to our rescue here along with a group of colleagues from Masaryk and Charles Universities in the Czech Republic. The paper \cite{CKMSVW} describes a polynomial time algorithm to compute all maps into a sphere from a space $X$.This set is presented as a group and a given map can be associated with a group element. The paper \cite{CKMVW2} provides a polynomial time computation of homotopy groups and Postnikov systems, and also discusses the extension problem, which can be solved under certain dimensionality conditions. The companion paper \cite{CKMVW1} shows that these conditions can not be relaxed. Finally, \cite{FV} gives an algorithm for showing if two given maps are homotopic. Especially the first three papers are very long and complicated, so I will state the claims they make that apply to the problem at hand and outline their approach. 

In the first section, I will give a more rigorous description of obstruction theory as it relates to the extension problem. In Section 10.2, I will present a possible data science application of obstruction theory and the Hopf map $f: S^3\rightarrow S^2$. Finally in Section 10.3, I will discuss some additional simplical set concepts needed to understand the computational papers I mentioned above. The papers themselves are the subject of Section 10.4.

I should mention that obstruction theory was used by Steenrod for the problem of constructing a cross section of a fiber bundle. See Steenrod \cite{Ste3}. Smith, Bendich, and Harer \cite{SBH} applied obstructions to constructions of cross sections for a data merging application. In what follows I will restrict myself to the extension problem.

\section{Classical Obstruction Theory and the Extension Problem}

This section is taken from Hu \cite{Hu}. Since we are interested in data science applications, we will assume our complexes are finite and triangulable. We will work with cell complexes here and simplicial sets in Section 10.4, but in our case they will be equivalent. 

Let $K$ be a finite cell complex, $L$ a subcomplex of $K$, and $Y$ a given pathwise connected space with a fixed point $y_0\in Y$. Let $K^n$ be the $n$-skeleton of $K$. The restricted extension problem is whether a map $f: L\rightarrow Y$ can be extended continuously over $K$ (as opposed to the more general problem of a map which is homotopic to an extension). In any case we would meet in data science, extendability is actually equivalent to a weaker condition. The extension problem depends only on the homotopy class of $f$. 

\begin{definition}
Let $f: X\rightarrow Y$ and $A\subset X$. Then a homotopy $h_t; A\rightarrow Y$ is a {\it partial homotopy}\index{partial homotopy} of $f$ if $f|A=h_0$. 
\end{definition}

\begin{definition}
A subspace $A$ of a space $X$ has the {\it homotopy extension property} or {\it HEP}\index{homotopy extension property}\index{HEP} in $X$ with respect to a space $Y$ if every partial homotopy $h_t: A\rightarrow Y$ of a map $f: X\rightarrow Y$ has an extension $g_t: X\rightarrow Y$ such that $g_0=f$. $A$ has the {\it absolute homotopy extension property} or {\it AHEP}\index{absolute homotopy extension property}\index{AHEP} in $X$ if it has the AHEP in $X$ with respect to every space $Y$. 
\end{definition}

For our applications, the AHEP will always hold as the following theorem shows. (See I.9. of \cite{Hu} for the proof.)

\begin{theorem}
Let $(X, A)$ be a finitely trianguable pair. Then $A$ has the AHEP in $X$. 
\end{theorem}

As we are dealing with finite simplicial complexes, A will have the AHEP in $X$ with respect to any $Y$. In this case, we can loosen our problem to ask whether there is a map $f: X\rightarrow Y$ such that $fh$ is {\it homotopic} to a given $g: A\rightarrow Y$ if $h: A\rightarrow X$ is inclusion.

The method will be to extend our given map step by step over the subcomplexes $\bar{K}^n=L\cup K^n$ until we hit an obstruction. It may also be possible to go back a step and try a different extension that will avoid the obstruction. The rest of this section will illustrate the process.

\subsection{Extension Index}
Again, let $K$ be a cell complex, $L$ a subcomplex of $K$, and $\bar{K}^n=L\cup K^n$. Also, again assume that $Y$ is pathwise connected.

\begin{definition}
A map $f: L\rightarrow Y$ is {\it n-extensible}\index{n-extensible} over $K$ for a given $n\geq 0$ if it has an extension over the subcomplex $\bar{K}^n$ of $K$. The {\it extension index}\index{extension index} of $f$ over $K$ is the least upper bound of integers $n$ such that $f$ is $n$-extensible over $K$.
\end{definition}

Since $Y$ is pathwise connected, we can map the vertices of $K\setminus L$ to arbitrary points of $Y$ and the one cells to paths connecting their endpoints. So we have the following.

\begin{theorem}
Every map $f: L\rightarrow Y$ is 1-extensible over $K$.
\end{theorem}

\begin{theorem}
Homotopic maps have the same extension index.
\end{theorem}

This holds since $L$ has the AHEP in $\bar{K}^n$.

Hu shows that nothing changes if there is a further subdivision or a different triangulation is used so we can assume that all maps are simplicial.

\subsection{The Obstruction Cochain}

Now fix an integer $n$. We will assume that $Y$ is $n$-simple. Recall that this means that the action of $\pi_1(Y, y_0)$ on $\pi_n(Y, y_0)$ is simple. (See Section 9.6.6). We can ensure this condition by letting $Y$ be simply connected. In applications, we will have $Y=S^m$ where $m\geq 2$.

Now let $g: \bar{K}^n\rightarrow Y$ be given. The map $g$ determines an $(n+1)$-cochain $c^{n+1}(g)$ with coefficients in $\pi_n(Y)$ as follows:

Let $\sigma$ be any $(n+1)$-cell in $K$. The boundary $\partial\sigma$ is an oriented $n$-sphere. Since $\partial\sigma\subset K^n$, the map $g_\sigma=g|\partial\sigma$ is already defined and as a map from $S^n\rightarrow Y$, determines an element $[g_\sigma]\in \pi_n(Y)$. Define the {\it obstruction cochain}\index{obstruction cochain} $c^{n+1}(g)$ of the map $g$ by taking $$[c^{n+1}(g)](\sigma)=[g_\sigma]\in\pi_n(Y)$$ for every $(n+1)$-cell $\sigma$ of $K$. 

It turns out that this cochain is actually a cocycle as the next theorem shows.

\begin{theorem}
The obstruction $c^{n+1}(g)$ is a relative $(n+1)$-cocycle of $K$ modulo $L$, i.e. $$c^{n+1}(g)\in Z^{n+1}(K, L; \pi_n(Y)).$$
\end{theorem}

{\bf Proof:} We first need to show that $$c^{n+1}(g)\in C^{n+1}(K, L; \pi_n(Y)).$$ We have already shown that it is a $(n+1)$-cochain of $K$ with coefficients in $\pi_n(Y)$. So we just need to show it vanishes on $L$. Let $\sigma$ be an $(n+1)$-cell of $L$. Since $g$ is defined on all of $L$, it has an extension over the closure of $\sigma$ so $[g_\sigma]=0$ in $\pi_n(Y)$.

Now we need to prove that $c^{n+1}(g)$ is a cocycle. Given an $(n+2)$-cell $\sigma$ in $K$, we need to show that $[\delta c^{n+1}(g)(\sigma)]=0$. Let $B$ denote the subcomplex $\partial\sigma$ of $K$ and $B^n$ denote its $n$-skeleton. Consider the homomorphisms in the following diagram.

$$C_{n+1}(B)\xrightarrow{\partial} Z_n(B)=Z_n(B^n)=H_n(B^n)\xleftarrow{h}\pi_n(B^n)\xrightarrow{k_\ast}\pi_n(Y),$$ where $h$ denotes the natural homomorphism between homotopy and homology (see Section 9.7.2.) and $k_\ast$ is induced by the map $k=g|B^n$. 

If $n>1$, then $B^n$ is $(n-1)$-connected so $h$ is an isomorphism by the Hurewicz Theorem. If $n=1$, then $h$ is an epimorphism and the kernel of $h$ is contained in the kernel of $k_\ast$ since $\pi_n(y)$ is abelian (Theorem 9.6.30). So in either case we get a well defined homomorphism $$\phi=k_\ast h^{-1}: Z_n(B)\rightarrow \pi_n(Y).$$ Since $C_n(B)$ is an abelian group, the kernel $Z_n(B)$ of $\partial: C_n(B)\rightarrow C_{n-1}(B)$ is a direct summand of $C_n(B)$. So, $\phi$ has an extension $$d: C_n(B)\rightarrow \pi_n(Y).$$ 

Now let $\tau$ be an $(n+1)$-cell in $B$. Then $[c^{n+1}(g)](\tau)$ is represented by a map $k|\partial\tau.$ We then get $$[c^{n+1}(g)](\tau)=k_\ast h^{-1}(\partial\tau)=d(\partial\tau).$$ This implies that $$[\delta c^{n+1}(g)](\sigma)=[c^{n+1}(g)](\partial\sigma)=d(\partial\partial\sigma)=0.$$ So $c^{n+1}(g)\in Z^{n+1}(K, L; \pi_n(Y)).$ $\blacksquare$

The next two theorems are easy consequences.

\begin{theorem}
The map $g: \bar{K}^n\rightarrow Y$ has an extension over $\bar{K}^{n+1}$ if and only if $c^{n+1}(g)=0$.
\end{theorem}

\begin{theorem}
If $g_0, g_1: \bar{K}^n\rightarrow Y$ are homotopic maps, then $c^{n+1}(g_0)=c^{n+1}(g_1)$.
\end{theorem}

\subsection{The Difference Cochain}

Now consider two maps $g_0, g_1: \bar{K}^n\rightarrow Y$ which are homotopic on $\bar{K}^{n-1}$. The difference of the obstruction cocycles of $g_0$ and $g_1$ is a coboundary.

To see this, consider a homotopy $h_t: \bar{K}^{n-1}\rightarrow Y$ such that $h_0=g_0|\bar{K}^{n-1}$ and $h_1=g_1|\bar{K}^{n-1}.$

Think of the unit interval $I=[0, 1]$ as a cell complex with a 1-cell $I$ and two 0-cells, 0 and 1. We have $\delta 0=-I$ and $\delta 1=I$. Then the product $J=K\times I$ is also a cell complex. Let $J^n$ be its $n$-skeleton. Let $$\bar{J}^n=(L\times I)\cup J^n=(\bar{K}^n\times 0)\cup (\bar{K}^{n-1}\times I)\cup (\bar{K}^n\times 1).$$

Let $F: J^n\rightarrow Y$ be defined by $F(x, 0)=g_0(x)$ and $F(x, 1)=g_1(x)$ for $x\in \bar{K}^n$ and $F(x, t)=h_t(x)$ for $x\in \bar{K}^{n-1}$. 

Then $F$ determines an obstruction cocycle $c^{n+1}(F)$ of the complex $J$ modulo $L\times I$ with coefficients in $\pi_n(Y)$. Then $c^{n+1}(F)$ coincides with $c^{n+1}(g_0)\times 0$ on $K\times 0$ and with $c^{n+1}(g_1)\times 1$ on $K\times 1$. 

Let $M$ denote the subcomplex $(K\times 0)\cup (L\times I)\cup (K\times 1)$ of $J=K\times I$. Then it follows that $$c^{n+1}(F)-c^{n+1}(g_0)\times 0-c^{n+1}(g_1)\times 1$$ is a cochain of $J$ modulo $M$ with coefficients in $\pi_n(Y)$. Since $\sigma\rightarrow\sigma\times I$ is a one-to-one correspondence between the $n$-cells of $K\setminus L$ and the $(n+1)$-cells of $J\setminus M$, it defines the isomorphism $$k: C^n(K, L, \pi_n(Y))\approx C^{n+1}(J, M; \pi_n(Y)).$$

So there is a unique cochain $d^n(g_0, g_1; h_t)\in C^n(K, L, \pi_n(Y))$ called the {\it deformation cochain}\index{deformation cochain} such that $$kd^n(g_0, g_1; h_t)=(-1)^{n+1}\{c^{n+1}(F)-c^{n+1}(g_0)\times 0-c^{n+1}(g_1)\times 1\}.$$  If we actually have $g_0|\bar{K}^{n-1}=g_1|\bar{K}^{n-1}$ and $h_t(x)=g_0(x)=g_1(x)$ for every $x\in \bar{K}^{n-1}$ and every $t\in [0, 1]$, we write $d^n(g_0, g_1; h_t)$ as $d^n(g_0, g_1)$ and call it the {\it difference cochain}\index{difference cochain}.

So the difference cochain has the property that $g_0$ and $g_1$ are actually equal on $\bar{K}^{n-1}$ rather than just homotopic.

\begin{theorem}
In the situation above, the homotopy $h_t: \bar{K}^{n-1}\rightarrow Y$ has an extension $H_t: \bar{K}^n\rightarrow Y$ such that $H_0=g_0$ and $H_1=g_1$ if and only if $d^n(g_0, g_1; h_t)=0.$
\end{theorem}

The deformation cochain has a nice coboundary formula. 

\begin{theorem}
$$\delta d^n(g_0, g_1; h_t)=c^{n+1}(g_0)-c^{n+1}(g_1).$$
\end{theorem}

Since $\delta I=0$, the isomorphism $k$ taking $\sigma$ to $\sigma\times I$ defined above commutes with $\delta$. Applying $\delta$ to the formula for $k$ and using the facts that $\delta 0=-I$, $\delta 1=I$ and that $c^{n+1}(F), c^{n+1}(g_0)$, and $c^{n+1}(g_1)$ are all cocyles  proves the theorem.

The  next result shows the importance of the deformation (difference) cochain. See \cite{Hu} for the proof.

\begin{theorem}
Let $g_0: \bar{K}^n\rightarrow Y$ and suppose there is a homotopy $h_t: \bar{K}^{n-1}\rightarrow Y$ with $h_0=g_0|\bar{K}^{n-1}$. For every cochain $c\in C^n(K, L; \pi_n(Y))$, there exists a map $g_1: \bar{K}^n\rightarrow Y$ such that $g_1|\bar{K}^{n-1}=h_1$ and $d^n(g_0, g_1; h_t)=c.$ If $z$ is a cocycle in $Z^{n+1}(K, L; \pi_n(Y))$ such that $z$ is in the same cohomology class as $c^{n+1}(g_0)$ mod $L$, there exists a map $g_1: \bar{K}^n\rightarrow Y$ such that $g_1|\bar{K}^{n-1}=h_1$ and $c^{n+1}(g_1)=z.$ 
\end{theorem}

In particular, if $h_t=g_0|\bar{K}^{n-1}$ for all $t\in [0,1]$, then there exists a map $g_1:  \bar{K}^n\rightarrow Y$ such that $g_0|\bar{K}^{n-1}=g_1|\bar{K}^{n-1}$ and $c^{n+1}(g_1)=z.$

\subsection{Eilenberg's Extension Theorem}

Let $g: \bar{K}^n\rightarrow Y$. Then $g$ determines an obstruction cocycle $c^{n+1}(g)\in Z^{n+1}(K, L; \pi_n(Y)).$ Since  $c^{n+1}(g)$ is a cocycle, it represents a cohomology class $\gamma^{n+1}(g)\in H^{n+1}(K, L; \pi_n(Y))$ called the {\it obstruction cohomology class}\index{obstruction cohomology class}. The next theorem shows its significance.

\begin{theorem}
{\bf Eilenberg's Extension Theorem:} $\gamma^{n+1}(g)=0$ if and only if there exists a map $H: \bar{K}^{n+1}\rightarrow Y$ such that $H|\bar{K}^{n-1}=g|\bar{K}^{n-1}.$
\end{theorem}

{\bf Proof:} Assume that such an $H$ exists and let $h=H|\bar{K}^n$. Since $h$ extends to $\bar{K}^{n+1}$, $c^{n+1}(h)=0.$ Since  $g|\bar{K}^{n-1}=h|\bar{K}^{n-1},$ the difference cochain $d^n(g, h)$ is defined. By Theorem 10.1.8, $c^{n+1}(h)=0$ implies that $c^{n+1}(g)$ is the coboundary of $d^n(g, h)$. Hence, $\gamma^{n+1}(g)=0$.

Now assume that $\gamma^{n+1}(g)=0$. Then $c^{n+1}(g)$ is in the same cohomology class as 0 mod $L$.  By Theorem 10.1.9, there exists a map $h: \bar{K}^n\rightarrow Y$ such that $g|\bar{K}^{n-1}=h|\bar{K}^{n-1}$, and $c_{n+1}(h)=0$. Then $h$ has an extension $H: \bar{K}^{n+1}\rightarrow Y$. $\blacksquare$

Now suppose that $g$ is an extension of $f: L\rightarrow Y$. If $c^{n+1}(g)\neq 0$, then $g$ can not be extended over $ \bar{K}^{n+1}$. But if $c^{n+1}(g)$ is in the zero cohomology class, then the obstruction is removable by modifying the values of $g$ only on the $n$-cells of $K\setminus L$.

\subsection{Obstruction Sets}

I will conclude this section with one last definition.

\begin{definition}
Let $f:L\rightarrow Y$. The {\it (n+1)-dimensional obstruction set}\index{obstruction set} $$O^{n+1}(f)\subset H^{n+1}(K, L; \pi_n(Y))$$ is defined as follows: If $f$ is not $n$-extensible over $K$ then $O^{n+1}(f)=\emptyset.$ So suppose instead that there exists an extension $g: \bar{K}^n\rightarrow Y$ of $f$. The cohomology class $\gamma^{n+1}(g)$ is called an {\it (n+1)-dimensional obstruction element of f}. Then $O^{n+1}(g)$ is the set of all $(n+1)$-dimensional obstruction elements of $f$.
\end{definition}

So we have a $(n+1)$-dimensional obstruction element of $f$ for each homotopy class of extensions $\bar{K}^n\rightarrow Y$ of $f$.

\begin{theorem}
Homotopic maps have the same $(n+1)$-dimesnional obstruction set.
\end{theorem}

Now let $(K', L')$ be another cellular pair and $\phi: (K', L')\rightarrow (K, L)$ be a proper cellular map. Then $\phi$ induces a homomorphism $$\phi^\ast: H^{n+1}(K, L; \pi_n(Y))\rightarrow  H^{n+1}(K', L'; \pi_n(Y)).$$ 

\begin{theorem}
If $f'=f\phi: L'\rightarrow Y$, then $\phi^\ast$ sends $O^{n+1}(f)$ into $O^{n+1}(f')$.
\end{theorem}

If $(K, L)$ is a subdivision of $(K', L')$ and $\phi$ is the identity, then $\phi^\ast$ is an isomorphism that sends $O^{n+1}(f)$ onto a subset of $O^{n+1}(f')$ in a one to one fashion. If $(K, L)$ and $(K', L')$ are simplicial, then the identity $\phi^{-1}$ is homotopic to a simplicial map and $\phi^\ast$ is an isomorphism that sends $O^{n+1}(f)$ onto $O^{n+1}(f')$. So the $(n+1)$-dimensional obstruction set is smaller on a subdivision and smallest on a simplicial pair. 

\begin{theorem}
If $(K, L)$ is a simplicial  pair and $f: L\rightarrow Y$ is a map, then $O^{n+1}(f)$ is a topological invariant and does not depend on the triangulation of $(K, L)$.
\end{theorem}

Now we let $(K, L)$ be a cellular pair, but in applications it will actually be simplicial. The following theorems come from the definition of obstruction sets and the Eilenberg Extension Theorem. They apply to both the cellular and simplicial case.

\begin{theorem}
The map $f: L\rightarrow Y$ is $n$-extensible over $K$ if and only if $O^{n+1}(f)\neq\emptyset.$
\end{theorem}

\begin{theorem}
The map $f: L\rightarrow Y$ is $(n+1)$-extensible over $K$ if and only if $O^{n+1}(f)$ contains the zero element of $H^{n+1}(K, L; \pi_n(Y))$.
\end{theorem}

\begin{theorem}
If $Y$ is $r$-simple and $H^{r+1}(K, L; \pi_r(Y))=0$ for every $r$ with $n\leq r<m$ then the $n$-extensibility of $f: L\rightarrow Y$ over $K$ implies its $m$-extensibility over $K$.
\end{theorem}

If $K\setminus L$ is of dimension not exceeding $m$ then the hypothesis of Theorem 10.1.16 implies that $f: L\rightarrow Y$ has an extension over $K$ if an only if it is $n$-extensible over $K$. 

\begin{theorem}
If $Y$ is $r$-simple and $H^{r+1}(K, L; \pi_r(Y))=0$ for every $r\geq 1$, then every map $f: L\rightarrow Y$ has an extension over $K$.
\end{theorem}

\subsection{Application to Data Science}

Let's look at what we know in relation to data science. Consider a point cloud $X$ again and let $VR_\epsilon(X)$ be the Vietoris-Rips complex at radius $\epsilon$. Let $Y$ be a simply connected space. It is then $r$-simple for every $r$. The most tractable situation is to let $Y=S^m$, where $m\geq 2$. As we will see in Chapter 12, finding homotopy groups of spheres is an extremely difficult problem involving some very sophisticated algerbraic topology tools. On the other hand, a huge number of these groups are known and low dimensional groups of low dimensional spheres  can be read from a table. (I will give you one in Chapter 12.) 

Now suppose $\delta>\epsilon$. Then $VR_\epsilon(X)\subset VR_\delta(X)$. Suppose $f$ is an interesting function going from $VR_\epsilon(X)\rightarrow Y$. Can $f$ be extended to $g: VR_\delta(X)\rightarrow Y$. Looking at the obstruction cochain $c^{n+1}(f)$ we can discover potential problem simplices. Let $L=VR_\epsilon(X)$ and $K=VR_\delta(X)$. If $Y=S^m$, we need to look at $H^{n+1}(K, L, \pi_n(S^m))$.

Now if $Y=S^m$, then $\pi_n(S^m)=0$ for $n<m$, so we know that we can extend $f$ to $\bar{K}^m$. For $n\geq m$, $\pi_n(S^m)$ is generally nonzero. So there are two things to look at. First of all is $H^{n+1}(K, L, \pi_n(S^m))=0$? Given our knowledge of $\pi_n(S^m)$  we can determine if the group is zero by computing cohomolgy with the help of the long exact sequence of the pair $(K, L)$ and the Universal Coefficient Theorem. If there is no simplex in $K\setminus L$ of dimension more than $k$, then $H^{k+1}(K, L, \pi_k(S^m))=0.$

Two more problems remain. \begin{enumerate}
\item Didn't we say that we need the coefficients to be in a field? $\pi_n(S^m)$ is not necessarily going to be a field. For example, $\pi_m(S^m)\cong Z$ for all $m>0$. This doesn't actually matter. Although we can't use the persistent homology algorithms, we still can compute the homology or cohomology for $\epsilon$ and $\delta$ fixed. Our problem may suggest interesting $(\epsilon, \delta)$ pairs or we can use every pair at which a simplex is added to the Vietoris-Rips complex at $\delta$. After all, this will only be finitely many times. In a different step, we can compute the persistent homology of our cloud with $Z_2$ coefficients and use features from both persistent homology and obstruction theory.

\item What if $H^{n+1}(K, L, \pi_n(S^m))\neq 0$? We can still extend $f$ to an $(n+1)$-simplex $\sigma$ if $f|\partial\sigma$ is the zero element of $\pi_n(S^m)$. How can we tell if that is the case. Using the papers \cite{CKMSVW, CKMVW1,  CKMVW2} of Cadek, et. al, we can determine whether a map is the zero element of a homotopy group provided certain restrictions on dimensionality are met. This leaves us with a lot of cases in which this is possible. These papers will be the subject of Section 10.4.\end{enumerate}

 So I will leave you with some questions to think about. \begin{enumerate}
\item For what data science problems would it be interesting to pick out a particular portion of a point cloud? (Or a graph, time series, etc?)
\item Is there a meaningful function from these points to a sphere for which obstructions in part of your cloud would have an interesting meaning?
\item Could $n$-extensibility of such a function be a feature that would help solve a difficult classification problem?
\end{enumerate}

I will give a rough idea of a possible application in the next section.

\section{Possible Application to Image Classification}

I will give a disclaimer at the beginning that this idea is not completely thought out. I know that there are a lot of great ways to segment and classify images, but the hope is that we will find something interesting that may have not shown up with other methods. Even if this doesn't accomplish that goal, it might inspire other situations for which obstruction theory is more suited.

Consider an image having $k\times n$ pixels. Consider each one to be a point on a grid. As a start, we could take the distance between them to be Euclidean distance where we consider neighboring pixels to be a distance of one unit apart. Pixels adjacent along a diagonal would be $\sqrt{2}$ units apart. To make things more interesting, we could vary inter-pixel distances in some way. Even more interesting, we could keep distances the same but stack the frames to represent a moving picture. Distances to the pixels on different levels could represent units of time. 

Now each pixel represents three eight byte numbers corresponding to $Y$ (luminance or greyscale) and two chrominance values $C_r$ and $C_b$. Each of these is a number between 0 and 255. Rescale the numbers by subtracting 128 and adding some small number like .001 to avoid all three numbers being 0. Let $v=[\bar{Y}, \bar{C_r}, \bar{C_b}]$ be the vector of rescaled values. Then $v$ corresponds to the point $v/||v||$ on $S^2$. So we have a map $f$ from the vertices of $VR_\epsilon(X)$ to $S^2$, and we can extend it continuously to the one and two dimensional simplices since $S^2$ is path connected and simply connected. The interesting case is when we try to add 3-simplices and 4-simplices. Recall that $\pi_2(S^2)\cong Z$ is generated by the $id_{S^2}$ and $f|\partial\sigma_3$  is the zero element if and only if it has degree 0. This will determine if the map can be extended to a given 3-simplex $\sigma_3$. For the 4-simplex $\sigma_4$, $\pi_3(S^2)\cong Z$ is generated by the Hopf map. We can use the method of the Cadek, et. al. papers to determine if $f|\partial\sigma_4$ is the zero element of this group. Due to the dimensionality restrictions in their results, this is as far as you can go.

Which simplices would have non-zero obstructions? Would anything interesting be going on there? Especially in the moving picture case, it would make for some interesting experiments.

\section{More on Simplicial Sets}

As computational homotopy theory can't really handle topological spaces themselves, it has to work with simplicial sets. Peter May's book \cite{May2} has a very well developed homotopy theory in the simplicial set context. It is probably a good idea to work through it if you want a full understanding of the algorithms in Section 10.4. For this book, I will just state the main theorems and outline some of the algorithms which get you there. Still we will need to know a little more about simplicial sets. I will again get the material from the more brief and easier to understand discussion in Friedman \cite{Fri}. While nowhere near as complete as \cite{May2}, it gives an idea of the major concepts. I will focus on Friedman's description of product sets and Kan complexes. You may find it helpful to review Section 6.4 at this point.

Two concepts from there will be especially important to remember: \begin{enumerate}
\item According to Theorem 6.4.1, any degenerate simplex is formed by applying a sequence of degeneracy maps to a unique non-degenerate simplex.
\item According to Theorem 6.4.2, a simplicial set $X$ can be realized as a CW complex with one $n$-cell for each nondegenerate $n$-simplex of $X$.
\end{enumerate}

\subsection{Products}

Recall from Section 8.5 that we had to work pretty hard to work with products of simplicial complexes and chain complexes. Simplicial sets are much easier by comparison. 

\begin{definition}
Let $X$ and $Y$ be simplicial sets. Their product $X\times Y$ is defined by its $n$-simplices, face maps, and degeneracy maps as follows: \begin{enumerate}
\item $(X\times Y)_n=X_n\times Y_n=\{(x, y)| x\in X_n, y\in Y_n\},$
\item If $(x, y)\in (X\times Y)_n$, then $d_i(x, y)=(d_ix, d_iy)$.
\item If $(x, y)\in (X\times Y)_n$, then $s_i(x, y)=(s_ix, s_iy)$. 
\end{enumerate}
\end{definition}

\begin{theorem}
Let $X$ and $Y$ be simplicial sets and $|X|$ and $|Y|$ their geometric realizations. If $|X\times Y|$ is a CW-complex, then $|X\times Y|$ and $|X|\times |Y|$ are homeomorphic. 
\end{theorem}

To see what this looks like, I will take Friedman's example of a square.

\begin{example}
Consider $\Delta^1\times\Delta^1$. We claim that $|\Delta^1\times\Delta^1|$ is the square $|\Delta^1|\times|\Delta^1|$. We need to find the nondegenerate simplices of $\Delta^1\times\Delta^1$.

Figure 10.3.1 \cite{Fri} illustrates the situation.

\begin{figure}[ht]
\begin{center}
  \scalebox{0.8}{\includegraphics{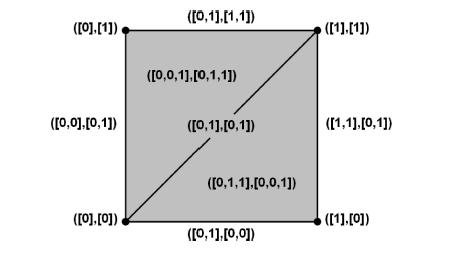}}
\caption{
\rm
Realization of $\Delta^1\times\Delta^1$ \cite{Fri}. 
}
\end{center}
\end{figure}

$\Delta^1$ has the two 0-simplices $[0]$ and $[1]$, so the product 0-simplices are $$X_0=\{([0], [0]), ([0], [1]), ([1], [0]), ([1], [1])\},$$ which are the four vertices of the square.

For dimension 1, the one simplices of $\Delta^1$ are $[0, 0]$, $[0, 1]$, and $[1, 1]$. So the product has nine 1-simplices. The only one made up of nondegenerate simplices is $([0, 1], [0, 1])$. Since $d_0([0, 1], [0, 1])=([1], [1])$ and $d_1([0, 1], [0, 1])=([0], [0])$, we know that $([0, 1], [0, 1])$ must be the diagonal. Those with one degenerate and one nondegenerate 1-simplex are $([0, 0], [0, 1]), ([0, 1], [0, 0]), ([1, 1],  [0, 1]),$ and $([0, 1], [1, 1])$, which by checking $d_0$ and $d_1$, we can see are the left, bottom, right, and top of the square respectively. The other four 1- simplices are degeneracies of the vertices. For example, $([0, 0], [1, 1])=s_0([0], [1]).$

Next, there are four 2-simplices of $\Delta^1: [0, 0, 0], [0, 0, 1], [0, 1, 1]$, and $[1, 1, 1]$. So there are sixteen 2-simplices of $\Delta^1\times\Delta^1$. There are also two degeneracy maps $s_0$ an $s_1$ from $(\Delta^1\times\Delta^1)_1\rightarrow(\Delta^1\times\Delta^1)_2.$ These act on the nine 1-simplices, but we reduce the number from 18. since $s_0s_0=s_1s_0$ and there are four degenerate 1-simplices $s_0v_i$ of $\Delta^1\times\Delta^1$ corresponding to the degeneracies of the four vertices. This reduces the number to 14 degenerate 2-simplices. There are two nondegnerate two simplices $([0, 0, 1], [0, 1, 1])$ and $([0, 1, 1], [0, 0, 1])$ which are the two triangles. We can see this by computing face maps. 

Now we need to show that all higher dimensional simplices in $\Delta^1\times\Delta^1$ are degenerate. Since there are no nondegenreate simplices of $\Delta^1$ of dimension greater than 1, each 3-simplex must be a double degeneracy of a 1-simplex. But there are only six options for a 1-simplex $e$. They are $s_0s_0e, s_0s_1e, s_1s_0e, s_1s_1e, s_2s_0e$, and $s_2s_1e$. But the simplicial set axioms reduce these to $s_1s_0e, s_2s_0e,$ and $s_2s_1e$. But the axioms then show that for 1-simplices $e$ and $f$, \begin{align*}
(s_1s_0e, s_1s_0f)&=s_1(s_0e, s_0f)\\
(s_1s_0e, s_2s_0f)&=(s_0s_0e, s_0s_1f)=s_0(s_0e, s_1f)\\
(s_1s_0e, s_2s_1f)&=(s_1s_0e, s_1s_1f)=s_1(s_0e, s_1f)\\
(s_2s_0e, s_1s_0f)&=(s_0s_1e, s_0s_0f)=s_0(s_1e, s_0f)\\
(s_2s_0e, s_2s_0f)&=s_2(s_0e, s_0f)\\
(s_2s_0e, s_2s_1f)&=s_2(s_0e, s_1f)\\
(s_2s_1e, s_1s_0f)&=(s_1s_1e, s_1s_0f)=s_1(s_1e, s_0f)\\
(s_2s_1e, s_2s_0f)&=s_2(s_1e, s_0f)\\
(s_2s_1e, s_2s_1f)&=s_2(s_1e, s_1f)
\end{align*}

So all 3-simplices of $\Delta^1\times\Delta^1$ are degenerate. So all higher dimensional simplices are also degenerate.

\end{example}

Friedman also looks at $\Delta^p\times \Delta^q$, See \cite{Fri} for the details.

\subsection{Kan Complexes}

In order to do homotopy theory with simplicial sets we need an equivalent condition to the homotopy extension property for topological spaces. We saw in Theorem 10.1.1 that the AHEP holds for a finitely triangulable pair. To apply this to simplicial sets, we need some definitions.

\begin{definition}
As a simplicial complex, the {\it kth horn} $|\Lambda^n_k|$ on the $n$-simplex $|\Delta^n|$ is the subcomplex of $|\Delta^n|$ obtained by removing the interior of $|\Delta^n|$ and the interior of the face $d_k\Delta^n$.  Let $\Lambda^n_k$ be the associated simplicial set. This set consists of simplices $[i_0,\cdots, i_m]$ with $0\leq i_0\leq\cdots\leq i_m\leq n$ such that not all numbers from 1 to $n$ are included (this would be the top face or a degeneracy of it), and we never have all numbers other than $k$ represented as this would be the $(n-1)$-face that was removed or a degeneracy of it. Figure 10.3.2 from \cite{Fri} shows an example.
\end{definition}

\begin{figure}[ht]
\begin{center}
  \scalebox{0.8}{\includegraphics{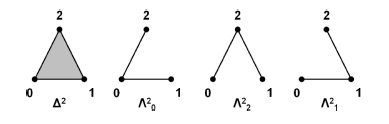}}
\caption{
\rm
Horns on $|\Delta_2|$ \cite{Fri}. 
}
\end{center}
\end{figure}

\begin{definition}
The simplicial set $X$ satisfies the {\it extension condition}\index{extension condition} or {\it Kan condition}\index{Kan condition} if any morphism of simplicial sets $\Lambda^n_k\rightarrow X$ can be extended to a simplicial morphism $\Delta^n\rightarrow X$. $X$ is then called a {\it Kan complex}\index{Kan complex} or a {\it fibrant}\index{fibrant}.
\end{definition}

Note that Kan complexes are named after Daniel Kan \cite{Kan} as opposed to my advisor, Donald Kahn despite having similar names.

Here is an alternate definition of the Kan condition.

\begin{definition}
The simplicial set $X$ satisfies the {\it Kan condition} for any collection of $(n-1)$ simplices $$x_0, x_1, \cdots, x_{k-1}, x_{k+1}, \cdots, x_n$$ in $X$ such that $d_ix_j=d_{j-1}x_i$ for any $i<j$ with $i\neq k$ and $j\neq k$ there is an $n$-simplex $x$ in $X$ such that $d_ix=x_i$ for all $i\neq k$.
\end{definition}

The alternative definition's condition on the simplices glues them together to form the horn $\Lambda^n_k$ (possibly with degenerate faces) and the definition says we can extend the horn to an $n$-simplex in $X$ which is possibly degenerate.

\begin{example}
The standard simplices $\Delta^n$ for $n>0$ do not satisfy the Kan condition. Let $\Delta^1=[0, 1]$ and $\Lambda^2_0$ be the horn depicted in Figure 10.3.2 consisting of $[0, 1], [0, 2]$ and their degeneracies. Let $f: \Lambda^2_0\rightarrow \Delta^1$ be a simplicial morphism with $f([0, 2])=[0, 0]$ and $f([0, 1])=[0, 1]$. Then $f$ is unique as it is specified on the nondegenerate simplices of $ \Lambda^2_0$.  It is also a simplicial map since all functions on all simplices are order preserving. But this can not be extended to a simplicial map $\Delta^2\rightarrow \Delta^1$ since $f(0)=0, f(1)=1,$ and $f(2)=0$. so this map is not order preserving on $\Delta^2$. 
\end{example}

Recall that the singular set functor, briefly mentioned in Theorem 6.4.3 is defined as follows:

\begin{definition}
If $Y$ is a topological space, let $\mathcal{S}(Y)$ be the simplicial set where $\mathcal{S}_n(Y)$ consists of the singular $n$-chains, (i,e. the continuous maps from $|\Delta^n|\rightarrow Y$). Then $\mathcal{S}$ is called the {\it singular set functor}\index{singular set functor}.
\end{definition}

Friedman states that this example is critical.

\begin{example}
If $Y$ is a topological space, then $\mathcal{S}(Y)$ satisfies the Kan condition. To see this let $f: \Lambda^n_k\rightarrow \mathcal{S}(Y)$. This map assigns for each $n-1$ face $d_i\Delta^n$ with $i\neq k$  of $\Delta ^n$ a singular simplex $\sigma_i: |\Delta^{n-1}|\rightarrow Y$. Every other simplex of $\Lambda^n_k$ is a face or a degeneracy of a face of one of these $(n-1)$ simplices, so the map $f$ is completely determined. Also the simplicial set axioms ensure that the $\sigma_i$ glued together yield a continuous function $f: |\Lambda^n_k|\rightarrow Y$. We can now extend this function to all of $|\Delta^n|$ by letting $r: |\Delta^n|\rightarrow |\Lambda^n_k|$ be a continuous retraction. Such a retraction exists since $(|\Delta^n|, |\Lambda^n_k|)$ is homeomorphic to $(I^{n-1}\times I, I^{n-1}\times 0)$. Then let $\sigma=fr: |\Delta^n|\rightarrow Y$. This is a singular $n$-simplex with $d_if=\sigma_i$ for all $i\neq k$. This then determines the extension.
\end{example}

Figure 10.3.3 from $\cite{Fri}$ illustrates the extension.

\begin{figure}[ht]
\begin{center}
  \scalebox{0.8}{\includegraphics{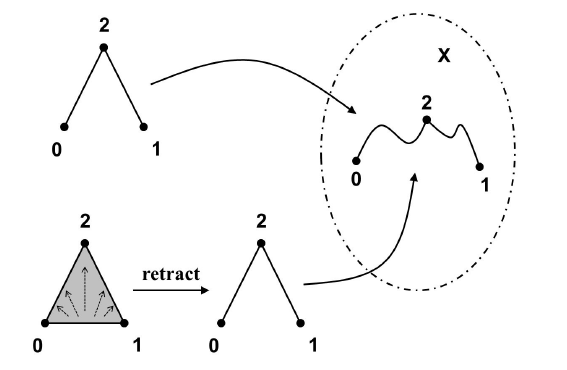}}
\caption{
\rm
A singular set satisfies the Kan condition \cite{Fri}. 
}
\end{center}
\end{figure}

\section{Computability of the Extension Problem}

Author's Note: This section is mostly still accurate with the exception of the paper \cite{CKMVW2}, which has some unaddressed issues. These were pointed out to me in a communication from Francis Sergeraert \cite{Serg1} shortly after the originial version of this book was published. As detailed in \cite{Serg2}, he objected that the authors had not addressed the role of functional programming that is the basis of Kenzo. A true cost computation would need to depend on programming closures which define a local scope. The polynomial time proof is then faulty. Unfortunately, the director of the Czech Team, Jiri Matousek, died around the time \cite{Serg2} was written, so the issues were never addressed. This makes the tractability of the extension problem still an open issue. My impression is that Kenzo could be used as a tool for this type of computation but it is still unknown under what conditions it would be practical. I think this would be a very interesting research problem for someone to try. As the issues are a little outside of my expertise, I will not change the original text of Section 10.4, but leave it to the reader to look at \cite{Serg2} and compare it with what I have written below.

In this section, I will discuss the possibility of deciding the extension problem on a computer. Although obstruction theory and the extension problem were one of Steenrod's main interests, he never found a systematic procedure. See \cite{Ste1} for an interesting description, although it is rather dated now. I will have a lot more to say about that paper in Chapter 11.

The only early paper addressing an algorithmic approach was the paper of E. H. Brown \cite{Bro} from 1957. Brown showed that $[X, Y]$ is computable if $Y$ is simply connected and all higher homotopy groups of $Y$ are finite. This is pretty limiting as it rules out $Y$ being a sphere, since $\pi_n(S^n)\cong Z.$ He also came up with a method to compute higher homotopy groups of $Y$ which drops the finiteness condition but is hard to generalize. We also know that it is undecidable to determine if $[S^1, Y]=\pi_1(Y)$ is trivial. 

The two papers of \u{C}adek et. al. \cite{CKMSVW, CKMVW2} give an algorithm for deciding if a function can be extended and the related question of characterizing maps into a sphere. As I mentioned above, the main problem to decide if $f: \partial\sigma\rightarrow Y$ can be extended to all of $\sigma$ is to determine if $[f]$ is the zero element of $\pi_n(Y)$ assuming that $\sigma$ is an $(n+1)$ simplex.

This problem is a special case of determining the group of maps from a space to a sphere $Y=S^d$. In \cite{CKMSVW} an even more general result is proved.

\begin{theorem}
Let $d\geq 2$. There is an algorithm that given finite simplicial complexes (or simplicial sets) $X$ and $Y$, where dim $X\leq 2d-2$ and $Y$ is $(d-1)$-connected, computes the abelian group $[X, Y]$ of homotopy classes of simplicial maps from $X$ to $Y$. Since this is a finitely generated abelian group, our result expresses it as a direct product of cyclic groups. Not only that, but given a simplicial map $f: X\rightarrow Y$, the element of $[X, Y]$ corresponding to $[f]$ can be determined. If $d$ is fixed, the algorithm runs in time polynomial to the number of simplices (or nondegenerate simplices in the simplicial set case) of $X$ and $Y$.
\end{theorem}

Based on this result, \cite{CKMVW2} explicitly applies it to the extension problem.

\begin{theorem}
Let $d\geq 2$. There is a polynomial time algorithm that given finite simplicial complexes $X$ and $Y$ and a subcomplex $A\subset X$, where dim $X\leq 2d-1$ and $Y$ is $(d-1)$-connected, determines whether a simplicial map $f: A\rightarrow Y$ admits an extension to a (not necessarily simplicial) map $f: X\rightarrow Y$. If the extension exists, dim $X\leq 2d-2$,  and $[X, Y]_f\subset[X, Y]$ is the subset of $[X, Y]$ consisting of the maps $X\rightarrow Y$ extending $f$, then the algorithm computes the structure of $[X, Y]_f$ as a coset in the abelian group $[X, Y]$. More generally $X$ can be a finite simplicial set (as opposed to a simplicial complex), and $Y$ can be a simplicial set whose homology can be computed in polynomial time.
\end{theorem}

Another \u{C}adek et. al. paper \cite{CKMVW1} shows the importance of the bounds on dimensionality in the previous two results. 

\begin{theorem}
The extendability problem is undecidable for finite simplicial complexes $A\subset X$ and $Y$ and a simplicial map $f: A\rightarrow Y$ if dim $X=2d$ and $Y$ is $(d-1)$-connected. In addition, for every $d\geq 2$, there is a fixed $(d-1)$-connected $Y=Y_d$ such that for maps into $Y_d$ with $A, X$, and $f$ as input, and dim $X=2d$, the extension problem is undecidable. 
\end{theorem}

In particular, for every {\it even} $d\geq 2$, the extension problem is undecidable for $X$ of dimension $2d$ and $Y=S^d$. For odd $d$, however, it was shown in \cite{Vok} that if $Y=S^d$, then the extension problem is decidable without any restriction on the dimension of $X$.

With the background I have given you so far, you should be able to read the papers, but they are quite long and involved. Rather than give a complete proof of the results I have just stated, I will go over some of the ideas. I will mainly focus on  \cite{CKMSVW} but I will also mention the extension proof in \cite{CKMVW2}.

To start, although the results would presumably address the data science applications I have in mind, there is no mention of actual code or any experiments that I could find. It is true that Francis Sergeraert, the inventor of Kenzo, is a coauthor on \cite{CKMSVW}, so it is possible that Kenzo could implement the algorithm, but I have not looked closely at it. I also have not found any mention of using these algorithms for data science applications. It seems to be an open area for future research and the results give a lot of hope that obstruction theory for many special cases in data science could be tractable.

I will now describe some of the ideas that lie behind the algorithm described in Theorem 10.4.1 and give the brief outline found in \cite{CKMSVW}

First of all, the step by step method of obstruction theory is not really a practical algorithm. The extensions at each step are not unique and we could end up searching an infinite tree of extensions as each step depends on previous choices. This is the reason that Brown \cite{Bro} needed all of the higher homotopy groups of $Y$ to be finite. To keep this process under control,  \cite{CKMSVW} uses the group structure on $[X, Y]$.

In representing abelian groups, we need to know what information is really available to us. {\it Semi-effective} abelian groups have a representation of their elements that can be stored on a computer and we can compute the group operations of addition and inverse on the level of representatives of their homotopy classes. A {\it fully effective} representation means that we have a list of generators and relations, and we can express any given element in terms of the generators. A homomorphism $f$ between two semi-effective abelian groups is {\it locally effective} if there is an algorithm that takes a representative of an element $a$ in the domain of $f$ and returns a representative of $f(a)$. 

Our geometric objects are represented as simplicial sets. A finite simplical set is just encoded as a list of its nondegenerate simplices and how they are glued together. We also need a way to represent an infinite simplicial set such as an Eilenberg-Mac Lane space. To handle this, we use a framework developed by Sergeraert (see \cite{RS}) in which the simplicial set is treated as a black box. This means we can encode the set and operations like the face operators. The simplicial set is called {\it locally effective.} A simplicial map between locally effective simplicial sets is {\it locally effective} if there is an algorithm that evaluates it on any given simplex of the domain simplicial set. 

One more step before outlining the algorithm is to talk about Postnikov systems. To convert to the notation used by \cite{CKMSVW}, I will modify the diagram from Definition 9.8.7 that we previously used for a Postnikov system.

$$\begin{tikzpicture}
  \matrix (m) [matrix of math nodes,row sep=3em,column sep=4em,minimum width=2em]
  {
& \vdots &\\
& P_{m+1} &\\
Y & P_m & K_{m+2}=K(\pi_{m+1}, m+2)\\
& \vdots &\\
& P_d=K(\pi_d, d) &\\
& (\ast) &\\};

\path[-stealth]
(m-1-2) edge  (m-2-2)
(m-2-2) edge  (m-3-2)
(m-3-2) edge  (m-4-2)
(m-4-2) edge  (m-5-2)
(m-5-2) edge  (m-6-2)
(m-3-1) edge node [left] {$\phi_{m+1}$} (m-2-2)
(m-3-1) edge node [below] {$\phi_m$} (m-3-2)
(m-3-2) edge node [above] {$k_m(Y)$} (m-3-3)

;

\end{tikzpicture}$$

The maps $\phi_i: Y\rightarrow P_i$ induce bijections $[X, Y]\rightarrow [X, P_i]$ if dim $X\leq i$. If $Y$ is $(d-1)$-connected, then for $i\leq 2d-2$, the Postnikov stage $P_i$ is an H-space or a group up to homotopy. This in turn induces an abelian group structure on $[X, P_i]$ for every space $X$ regardless of the dimension of $X$. But if dim $X\leq 2d-2$, then we have the bijection $[X, Y]\rightarrow [X, P_{2d-2}]$ which then gives the abelian group structure for $[X, Y]$. To prove Theorem 10.4.1, we compute the stages $P_0, \cdots,  P_{2d-2}$ of a Postnikov system of $Y$ and then by induction on $i$, determine $[X, P_i]$ for $i\leq 2d-2$. We then return the description of $[X, P_{2d-2}]$ as an abelian group.

For the inductive computation of $[X, P_i]$ we can ignore the dimension restriction on $X$ as we will also need to compute $[SX, P_{i-1}]$ where $SX$ is the suspension of $X$ and is of dimension one greater than that of $X$.

Finally, we need the $P_i$ to be Kan simplicial sets.  This insures that every continuous map $X\rightarrow P_i$ is homotopic to a simplicial map. Simplicial maps are discrete and finitely representable objects.

Here is an outline of the algorithm described in Theorem 10.4.1.

{\bf Main Algorithm:}

\begin{enumerate}
\item Using the algorithm in Section 4 of \cite{CKMVW2}, find a suitable representation of the first $2d-2$ stages of a Postnikov system for $Y$. This algorithm provides us with the isomorphism types of the first $2d-2$ homotopy groups $\pi_i=\pi_i(Y)$ along with the Postnikov stages $P_i$, and the Eilenberg-Mac Lane spaces $K_{i+1}=K(\pi_i, i+1)$ and $L_i=K(\pi_i, i)$ for $i\leq 2d-2$ as locally effective simplicial sets. It also outputs the maps between these spaces such as the Postnikov classes $k_{i-1}: P_{i-1}\rightarrow K_{i+1}$ as locally effective simplicial maps. 
\item Given a finite simplicial set $X$, the algorithm computes $[X, P_i]$ as a fully effective abelian group by induction on $i$ for $i\leq 2d-2$. The final output is then $[X, P_{2d-2}]$. This step involves the following:
\begin{enumerate}
\item Construct locally effective simplicial maps $\boxplus_i: P_i\times P_i\rightarrow P_i$ and $\boxminus_i: P_i\rightarrow P_i$ representing group addition and inverse respectively. This does imply that the spaces $P_i$ have some sort of group structure. \u{C}adek, et. al. state that the $P_i$ have an H-space structure (i. e. it is group up to homotopy) since $P_i$ has nonzero homotopy groups in the range from $d$ to $2d-1$ where $Y$ is $(d-1)$-connected. They outline the proof in Section 5 of \cite{CKMSVW}, but refer the reader to Whitehead \cite{Whi} for a full proof.
\item Using the group structure on the $P_i$, induce group operations $\boxplus_{i\ast}$ and $\boxminus_{i\ast}$ on the set $SMap(X, P_i)$ of the simplicial maps from $X$ to $P_i$. These correspond to the group operations in $[X, P_i]$, the group of homotopy classes of $SMap(X, P_i)$.The result is a semi-effective representation for $[X, P_i]$. 
\item The last step is to convert this semi-effective representation of $[X, P_i]$ into a fully effective one. For this step, $[X,  L_i]$ and $[X, K_{i+1}]$ are easy to compute as fully effective abelian groups, since the groups of homotopy classes of maps from $X$ into Eilenberg-Mac Lane spaces are isomorphic to certain cohomology groups of $X$. (I will explain this in Chapter 11.) In addition to $[X,  L_i]$ and $[X, K_{i+1}]$, assume that we have already computed $[SX, P_{i-1}]$ and $[X, P_{i-1}]$ in the induction step. (Here $SX$ is the suspension of $X$.)
\item The four groups from the previous step fit into an exact sequence of abelian groups as follows:
$$\begin{tikzpicture}
  \matrix (m) [matrix of math nodes,row sep=3em,column sep=4em,minimum width=2em]
  {
\left[SX, P_{i-1}\right] & \left[X, L_i\right] & \left[X, P_i\right] &\\
 & &\left[X, P_{i-1}\right] &\left[X, K_{i+1}\right]\\};

\path[-stealth]
(m-1-1) edge  (m-1-2)
(m-1-2) edge  (m-1-3)
(m-1-3) edge  (m-2-3)
(m-2-3) edge  (m-2-4)

;

\end{tikzpicture}$$ For the specific maps between them and how you extract the value of $[X, P_i]$, I will refer you to Section 6 of \cite{CKMSVW} as some detailed explanation is required. The idea, though is that we filter out the maps $X\rightarrow P_{i-1}$ that can be lifted to maps $X\rightarrow P_i$. (This requires the computation of a different type of obstruction. See for example, \cite{Hu}.) For each map that can be lifted, we determine all possible liftings and check on which lifts are homotopic. Since there are infinitely many homotopy classes of maps involved in these operations, we need to work with generators and relations in the appropriate abelian groups of homotopy classes. In the end, we have a fully effective representation of $[X, P_i]$.
\end{enumerate}
\end{enumerate}

If $Y$ is fixed, we can evaluate the Postnikov classes $k_i$ for $i\leq 2d-2$ as one time work. Note that if $Y=S^d$, then $k_d$ corresponds to the Steenrod square $Sq^2$. You will learn all about Steenrod squares in Chapter 11.

The reader is encouraged to look at Sections 3-6 of \cite{CKMSVW} for explicit details of the algorithm we have just outlined.

Finally, we will see how Theorem 10.4.1 can be used to derive the extension algorithm in Theorem 10.4.2. The proof comes from \cite{CKMVW2}.

{\bf Proof of Theorem 10.4.2:}

Suppose $A\subset X$ and $Y$ are simplicial sets and $f: A\rightarrow Y$ is a simplicial map. Let $X$ be finite with dim $X\leq 2k-1$ and let $Y$ be $(k-1)$-connected. By Theorem 7.6.22 of \cite{Spa}, a continuous extension to $X$ exists if and only if $\phi_{2k-2}f: A\rightarrow P_{2k-2}$ can be continuously extended to $X$, where $\phi_{2k-2}: Y\rightarrow P_{2k-2}$ is the map in the Postnikov system for $Y$. By the homotopy extension property, this happens when there is a map $X\rightarrow P_{2k-2}$ whose restriction to $A$ is homotopic to $\phi_{2k-2}f$. In terms of homotopy classes, we need $[\phi_{2k-2}f]$ to be in the image of the restriction map $\rho:[X, P_{2k-2}]\rightarrow [A, P_{2k-2}].$

Now to compute $[X, Y]$ we compute the group $[X, P_{2k-2}]$.  If dim $X=2k-1$ the two groups are not isomorphic but the computation of $[X, P_{2k-2}]$ still works. By the methods of the algorithm for Theorem 10.4.1, this group is fully effective and the generators are specified as simplicial maps. We also use the same methods to compute $[A, P_{2k-2}]$. The group operation in  $[X, P_{2k-2}]$ is induced by an operation $\boxplus$ on $SMap(X, P_{2k-2})$ which is defined on each simplex, so the restriction map $\rho$ is a group homomorphism. 

If $[g]\in [X, P_{2k-2}]$ is represented by $g$, then the restriction $g|A$ is a representative of an element of $[A, P_{2k-2}]$ which we can express in terms of the generators of $[A, P_{2k-2}]$. So $\rho$ is polynomial time computable and we can compute its image as a subgroup of $[A, P_{2k-2}]$ by computing the images of the generators of $[X, P_{2k-2}]$ . Given a simplicial map $f: A\rightarrow Y$, we compute the corresponding element $[\phi_{2k-2}f]\in [A, P_{2k-2}]$ and test if it lies in the image of $\rho$. This is the required polynomial time algorithm.

If we have the additional restriction that dim $X\leq 2k-2$, then $[X, Y]\cong [X, P_{2k-2}]$ and $[A, Y]\cong [A, P_{2k-2}]$, so if $x$ is in the image of $\rho$, we can compute the preimage $\rho^{-1}(x)$ as a coset of $[X, P_{2k-2}]$. (Here $\rho$ is represented by a matrix.) This coset is isomorphic to the group $[X, Y]_f$ of maps $X\rightarrow Y$ which are extensions of $f$. This proves Theorem 10.4.2.  $\blacksquare$

This leaves us with some questions:\begin{enumerate}
\item Would there be an interesting map $f$ from a Vietoris-Rips complex to a space $Y$ which may not be a sphere but is $(k-1)$-connected?
\item The Hopf map is within range of Theorem 10..4.2. ($X=S^3$ has dimension 3, and $S^2$ is 1-connected. Letting k=2, we have dim $X=2(2)-1.$) What would the algorithm do in the case of classifying still or moving images. How long would it take?
\item Many homotopy groups of spheres are tabulated. Would knowing these speed up the algorithm? 
\end{enumerate}

This area is wide open. Assuming the algorithms from this section, it would be very interesting to see what situations are tractable and what we could learn from our data.

\chapter{Steenrod Squares and Reduced Powers}

With this chapter, we are definitely at the frontier of where algebraic topology and data science intersect. I have seen very little discussion anywhere of the use of Steenrod squares for a data science application. So why am I including it in this book?

\begin{enumerate}
\item Steenrod squares provide an even richer algebraic structure than cohomology. Using $Z_2$ coefficients (as we must with Steenrod squares), cup product provides a map $$\cup: H^p(X; Z_2)\times H^q(X; Z_2)\rightarrow H^{p+q}(X; Z_2).$$ Steenrod squares provide maps $$Sq^i: H^p(X; Z_2)\rightarrow H^{p+i}(X; Z_2),$$ where $0\leq i\leq p$. It turns out that for $a\in H^p(X; Z_2),$ we have $Sq^0(a)=a$ and $Sq^p(a)=a\cup a$. (Hence the term Steenrod {\it squares}.) This also seems relevant as we normally work in $Z_2$ for persistent homology. But what if you hate $Z_2$ and really want to do all your work in $Z_{683}$? (Yes, 683 is a prime number.) You can always use reduced powers. They are of the form $$\mathcal{P}^i: H^q(X; Z_p)\rightarrow H^{q+2i(p-1)}(X; Z_p),$$ where $p$ is an odd prime. I would probably try to talk you out of it in a practical situation, especially $p=683$. 

\item There are explicit formulas you can program onto a computer for finding Steenrod squares of simplicial sets. The paper of Gonz\'{a}lez-D\'{\i}az and Real \cite{GDR1} provides such a formula. The bad news is that the complexity is exponential. If $c\in C^{i+k}(X; Z_2)$ then the number of face operators taking part in the formula for $Sq^i(c)$ is $O(i^{k+1})$. The good news is that these numbers are very manageable if you keep the dimensions low. Also, Aubrey HB \cite{HB} came up with a simplified formula for $Sq^1$. Just computing a small number of the squares could make the difference in a hard classification problem.

\item Steenrod squares are used to construct {\it Stiefel-Whitney classes}, which are cohomology classes used to classify {\it vector bundles}. Vector bundles are fiber bundles whose fiber is a vector space. The point is that Steenrod squares are already used in a classification problem. 

\item Steenrod squares are historically tied to obstruction theory and the extension problem. Not only was it a motivation for their discovery, but the main theorem about the correspondence between cohomology operations and cohomology classes of Eilenberg-Mac Lane spaces relies on obstruction theory for its proof. (This will be discussed in Section 11.2.)
\end{enumerate}

In Section 11.1, I will discuss two historical papers by Steenrod that give some insight into his original motivation for Steenrod squares. Then I will go through the modern theory mainly using the book by Mosher and Tangora \cite{MT}. This will include the correspondence between cohomology operations (of which Steenrod squares are a special case) and the cohomology of Eilenberg-Mac Lane spaces, construction and properties of Steenrod squares, the Hopf invariant, the structure of the Steenrod algebra, and some approaches to the very hard problem of computing the cohomology of Eilenberg-Mac Lane spaces. Then I will briefly discuss reduced powers which involve cohomology with coefficients in $Z_p$, where $p$ is an odd prime. Mosher and Tangora don't discuss this case, so we will need to refer to the book by Steenrod and Epstein \cite{SE}. Then I will briefly describe vector bundles and Stiefel-Whitney classes. The classic book on this subject is Milnor and Stasheff \cite{MS}, but Aubrey HB deals with them as well in \cite{HB}. I will conclude with the formulas found in \cite{GDR1} and \cite{HB} that allow for computations on a computer.

\section{Steenrod's Original Motivation}

Steenrod squares first appeared in the paper \cite{Ste2}, written in 1946 and published in 1947. At that time, homotopy theory was in its very early stages and there was a lot of interest in describing groups of homotopy classes of maps between one type of complex and another. The first such result was classfiying maps from $S^n$ to itself in terms of degree. In the 1930's, Hopf discovered that $\pi_3(S^2)\cong Z$ generated by the Hopf map. Then Freudenthal and Pontrjagin discovered that $\pi_{n+1}(S^n)\cong Z_2$ for $n>2$. 

The problem Steenrod wanted to generalize was first solved by Hopf in 1933. Rather than map $S^n$ to itself, map an $n$-dimensional complex $K^n$ to $S^n$. Hopf showed that two maps $f, g: K^n\rightarrow S^n$ are homotopic (written $f\sim g$) if and only if $f^\ast(s^n)=g^\ast(s^n)$ where $s^n$ is the generator of $H^n(S^n)\cong Z$. Any $n$-cocycle in $K^n$ is $f^\ast(s^n)$ for some $f$. So the homotopy classes are in one to one correspondence with $H^n(K^n)$. 

In 1941, Pontrjagin \cite{Pont1} enumerated all of the homotopy classes of maps of a 3-complex $K^3$ to a 2-sphere $S^2$. If $f\sim g$, then $f^\ast(s^2)=g^\ast(s^2)$. In this case, $g\sim g'$, where $g'$ is a function that coincides with $f$ on the 2-skeleton $K^2$. For a 3-cell $\sigma$, $f$ and $g'$ define an element $d^3(f, g', \sigma)\in \pi_3(S^2)$. Then Pontrjagin showed that $f\sim g'$ and therefore to $g$ if there exists a 1-cocycle $e^1$ in $K^3$ such that $d^3(f, g')$ is cohomologous to $2e^1\cup f^\ast (s^2)$, where $d^3(f, g')$ is a 3-cocycle in $Z^3(K^3; \pi_3(S^2))=Z^3(K^3; Z)$. He also showed that any pair consisting of a 2-cocycle and a 3-cocycle on $K^3$ can be realized as $f^\ast (s^2)$ and $d^3(f, g)$ for a suitable $f$ and $g$.

Steenrod's goal was to generalize this result to the case of maps from an (n+1)-complex $K$ to $S^n$. To do this, the cup product is generalized to the {\it cup-i} product $\cup_i$. For $i\geq 0$, and cochains $u$ and $v$ of dimensions $p$ and $q$ respectively, $u\cup_i v$ has dimension $p+q-i$. If $i=0$, $u\cup_0 v=u\cup v$. 

In the original formula, the cup-$i$ product was similar to our original definition of simplicial cup product. (See the formula immediately after Definition 8.3.2.) If $u$ and $v$ are as above, then for a $(p+q)$-simplex, apply $u$ to the $p$-simplex spanned by the first $p+1$ vertices and $v$ to the $q$-simplex spanned by the last $q+1$ vertices. Then multiply the result. For the cup-$i$ product, we have some overlap in vertices and need to reorder them in a specific way. The result can then be positive or negative depending on whether the reordering involves an even or odd permutation. (Remember this means flipping two elements an even or odd number of times.) 

The coboundary formula is also somewhat more complicated. Recall from Theorem 8.31 that $$\delta(c^p\cup c^q)=(\delta c^p)\cup c^q+(-1)^pc^p\cup (\delta c^q).$$ This means the cup product of 2 cocycles is itself a cocycle.

Here is the corresponding formula for cup-$i$ products. 

\begin{theorem}
If $u$ and $v$ are cochains of dimensions $p$ and $q$, then $$\delta(u\cup_i v)=(-1)^{p+q-i} u \cup_{i-1} v+(-1)^{pq+p+q}v\cup_{i-1} u +\delta u\cup_i v +(-1)^pu \cup_i \delta v.$$
\end{theorem}

The proof is long, ugly, and painful. See \cite{Ste2} if you are curious. 

The point, though, is that the cup-i product of 2 cocycles is no longer a cocycle. But all is not lost. If we take $u\cup_i u$ and $u$ is a cocycle, the last 2 terms go away. If $u$ has dimension $p$, then $$\delta(u\cup_i u)=(-1)^{2p-i} u \cup_{i-1} u+(-1)^{p^2+2p}u\cup_{i-1} u=((-1)^{2p-i}+(-1)^{p^2+p}) u \cup_{i-1} u.$$ This is equal to $2( u \cup_{i-1} u)$ if $p-i$ is even or 0 if $p-i$ is odd. Working mod 2, $u\cup_i u$ is a cocycle if $u$ is a cocycle. We replace $(u\cup_i u)$ with the homomorphism $$Sq_i: H^p(K; Z_2)\rightarrow H^{2p-i}(K; Z_2).$$ Later we will switch to the more modern notation $Sq^i$ which is equal to $Sq_{p-i}$.

The Pontrjagin formula for $n>2$ becomes $d^{n+1}(f, g)$ is cohomlogous to $e^{n-1}\cup_{n-3}e^{n-1}$. 

The homotopy classification theorem is obtained from an extension theorem. It states that a map $f$ from the $n$-skeleton $K^n$ of an $(n+2)$-complex $K$ to $S^n$  can be extended to $K$ if and only if $f^\ast(s^n)$ is a cocycle in $K$ and $f^\ast(s^n)\cup_{n-2} f^\ast(s^n)=0.$ See \cite{Ste2} for the proof and all of the details. 

The other interesting historical paper is \cite{Ste1}. Steenrod died in 1971, but he had given permission to the journal "Advanced Mathematics" to publish the paper which was based on the Colloquium Lectures given at the Annual Meeting of the American Mathematical Society and Pennsylvania State University in August 1957. It was published posthumously in 1972. It is much more polished than his original paper and a great introduction, but a little dated. What I want to give you is a taste of the connection between the development of cohomology operations and the extension problem. I will summarize the first part of \cite{Ste1} and use Mosher and Tangora \cite{MT} for the actual theory.

Steenrod's introduction to the subject focused on the extension problem. I will restate it in the way it appeared at the beginning of Chapter 9. Let $X$ and $Y$ be topological spaces. Let $A$ be a closed subset of $X$ and let $h: A\rightarrow Y$ be a mapping. A mapping $f: X\rightarrow Y$ is called an {\it extension} of $h$ if $f(x)=h(x)$ for each $x\in A$. Letting $g: A\rightarrow X$ be inclusion we want to find a map $f$ such that $fg=h$. In terms of diagrams, we have the following: 

$$\begin{tikzpicture}
  \matrix (m) [matrix of math nodes,row sep=3em,column sep=4em,minimum width=2em]
  {
& X &\\
A & & Y\\};

\path[-stealth]

(m-2-1) edge node [left] {$g$} (m-1-2)
(m-2-1) edge node [below] {$h$} (m-2-3)
(m-1-2) edge [dashed] node [right] {$f$} (m-2-3)

;

\end{tikzpicture}$$

We would like to determine if such an extension $f$ exists. 

We use algebraic topology to attack this problem by turning a geometry problem into an algebra problem. Using homology, we want to find $f_\ast$ for each $q$  in the following diagram: 

$$\begin{tikzpicture}
  \matrix (m) [matrix of math nodes,row sep=3em,column sep=4em,minimum width=2em]
  {
& H_q(X) &\\
H_q(A) & & H_q(Y)\\};

\path[-stealth]

(m-2-1) edge node [left] {$g_\ast$} (m-1-2)
(m-2-1) edge node [below] {$h_\ast$} (m-2-3)
(m-1-2) edge [dashed] node [right] {$f_\ast$} (m-2-3)

;

\end{tikzpicture}$$

Here all the maps are group homomorphisms. A necessary condition for an extension to exist is that there is a group homomorphism $f_\ast$ such that $f_\ast g_\ast=h_\ast$. Unfortunately, this condition is not sufficient. As Steenrod points out, too much information is lost in the translation from spaces to groups. We would like to encode as much algebra in the geometry as possible.

We can do better by using cohomology. Now we have the diagram 

$$\begin{tikzpicture}
  \matrix (m) [matrix of math nodes,row sep=3em,column sep=4em,minimum width=2em]
  {
& H^q(X) &\\
H^q(A) & & H^q(Y)\\};

\path[-stealth]

(m-1-2) edge node [left] {$g^\ast$} (m-2-1)
(m-2-3) edge node [below] {$h^\ast$} (m-2-1)
(m-2-3) edge [dashed] node [right] {$f^\ast$} (m-1-2)

;

\end{tikzpicture}$$

Now we want to find an $f^\ast$ such that $h^\ast=g^\ast f^\ast$. Another key difference, form the previous case is that our maps are now {\it ring} homomorphisms. This means that they not only respect addition, but also ring multiplication. In other words $$f^\ast(u\cup v)=f^\ast(u)\cup f^\ast(v).$$  Below I will give some examples on how this approach is an improvement. 

Meanwhile, we can do even better. A {\it cohomology operation} $T$ relative to dimensions $q$ and $r$ is a collection of functions $\{T_X\}$, one for each space $X$, such that $$T_X: H^q(X)\rightarrow H^r(X)$$ and for each mapping $f: X\rightarrow Y,$ $$f^\ast T_Y u=T_X f^\ast u$$ for all $u\in H^q(Y)$. (I will state a more general definition in the next section.)

$Sq^i$ and $\mathcal{P}^i$ are two examples of cohomology operations. 

If $f$ is the solution to an extension problem then $f^\ast$ must be a ring homomorphism with $f^\ast T_y=T_x f^\ast$ for every cohomology operation $T$. 

Some examples will illustrate these ideas.

\begin{example}
We have already seen that the identity on $S^n$ can not be extended to the ball $B^{n+1}$ of which $S^n$ forms the boundary. I will review this case here for comparison. Converting to homology, we have the diagram

$$\begin{tikzpicture}
  \matrix (m) [matrix of math nodes,row sep=3em,column sep=4em,minimum width=2em]
  {
& H_q(B^{n+1}) &\\
H_q(S^n) & & H_q(S^n)\\};

\path[-stealth]

(m-2-1) edge node [left] {$g_\ast$} (m-1-2)
(m-2-1) edge node [below] {$h_\ast$} (m-2-3)
(m-1-2) edge [dashed] node [right] {$f_\ast$} (m-2-3)

;

\end{tikzpicture}$$

Let $q=n$. Then $H_n(S^n)\cong Z$. $B^{n+1}$ is acyclic, so in particular $H_n(B^{n+1})=0$. Since $h$ is the identity, $h_\ast$ is the identity homomorphism on $Z$. But then $g_\ast=0$ so  for any $f_\ast$, we have $f_\ast g_\ast=0\neq h_\ast=id_Z$.
\end{example}

Next we have a case where we need to use cohomology to make a decision.

\begin{example}
Let $X=CP^2$, the complex projective plane. (See Definition 4.3.9.) By Theorem 4.3.8, $X$ is a CW complex with one cell each in dimensions 0, 2, and 4. Let $A$ be $CP^1$. Then $A$ is the complex projective line. It has one 0-cell and one 2-cell so it is homeomorphic to $S^2$. I claim that $A$ is not a retract of $X$. Again, this means that the identity on $A$ can not be extended to all of $X$. 

Suppose this was not the case. If $f: X\rightarrow A$ is a retract, then if $g: A\rightarrow X$ is inclusion, we would have $g^\ast f^\ast=identity$,  This means that $H^q(X)=$im $f^\ast \oplus\ker g^\ast$ for any $q$. We can abbreviate this to $$H^\ast(X)=\mbox{im }f^\ast \oplus\ker g^\ast.$$ In addition, $g^\ast: $im $f^\ast\rightarrow H^\ast(A)$ is an isomorpism. Using the fact that $f^\ast$ and $g^\ast$ are ring homomorphisms, we have that im $f^\ast$ is a subring and $\ker g^\ast$ is an ideal. 

In our case, we need $\ker g^\ast$ to be a direct summand. The table shows the cohomology of $X$ and $A$ with integer coefficients.

\kbordermatrix{ & H^0 &  H^1 &  H^2 &  H^3 &  H^4\\
A & Z & 0 & Z & 0 & 0\\
X & Z & 0 & Z & 0 & Z} \qquad

Now $g^\ast: H^q(X)\rightarrow H^q(Z)$ is an isomorphism in dimensions 0 and 2 and $g^\ast=0$ in dimension 4. So the necessary direct sum composition exists and is unique. In dimensions 0 and 2, $\ker g^\ast=0$ and im $f^\ast=H^q(X)$. In dimension 4, $g^\ast=0$, so $\ker g^\ast=H^4(X)$ and  im $f^\ast=0$. So using only the group structure, we might think that the retraction exists.

But looking at the ring structure, we see that this is not the case. This is because im $f^\ast$ is not a subring. To see this, note that in dimension 2, $H^2(X)\cong Z$, so let $u$ be a generator. Then $u$ is in  im $f^\ast$. Since $X$ is a manifold, we can use {\it Poincar\'{e} duality} to show that $u\cup u$ generates $H^4(X)$. (I haven't covered Poincar\'{e} duality, which relates homology and cohomology of manifolds as we won't make much use of it. If you are curious, see Munkres \cite{Mun1} Sections 65 and 68.)  Anyway, since we said that im $f^\ast=0$ in dimension 4,  we see that $u\in$im $f^\ast$,  but $u\cup u\notin$im $f^\ast$. So $f^\ast$ is not a subring and we have a contradiction. 

So we know there is no retraction from $CP^2$ to $CP^1$ but we needed the ring structure to see that.
\end{example}

The final example is a retraction problem that can not be decided by the cohomology ring but can be determined using the Steenrod squaring operations.

\begin{example}
Let $P^n$ be $n$-dimensional projective space. Let $X$ be $P^5$ with the subspace $P^2$ collapsed to a point. Let $A\subset X$ be the image of $P^4$ under the collapsing map $P^5\rightarrow X$. We claim that $A$ is not a retract of $X$.

Again, suppose $A$ is a retract of $X$. Then as before, $\ker g^\ast$ is a direct summand of $H^\ast(X; Z_2)$. The cohomology of $A$ and $X$ with $Z_2$ coefficients is shown in the following table: 

\kbordermatrix{ & H^0 &  H^1 &  H^2 &  H^3 &  H^4 & H^5\\
A & Z_2 & 0 & 0 & Z_2 & Z_2 & 0\\
X & Z_2 & 0 & 0 & Z_2 & Z_2 & Z_2} \qquad

Now $g^\ast: H^\ast(X; Z_2)\rightarrow H^\ast(A; Z_2)$ is an isomorphism in dimensions less than 5. So im $f^\ast=H^q(X; Z_2)$ in dimensions less than 5. and  im $f^\ast=0$ in dimension 5. In this case, though, im $f^\ast$ actually is a subring. To see this, the smallest positive dimension with $H^q(X; Z_2)\neq 0$ is dimension 3. But if $u$ is the nonzero element of $H^q(X; Z_2)$ for $q\geq 3$, then $u\cup u$ has dimension at least 6. Since $H^q(X; Z_2)=0$ for $q\geq 6$, $u\cup u=0$. Thus im $f^\ast$ is a subring. So it is still possible for $A$ to be a retract of $X$.

Now consider the cohomology operation $Sq^2: H^3(X; Z_2)\rightarrow H^5(X; Z_2).$ Let $u$ be the nonzero element of $H^3(X; Z_2)$. It turns out that $Sq^2 u$ is the nonzero element of $H^5(X; Z_2)$. Now im $f^\ast$ contains $u$, but is zero in dimension 5. So if $v\in H^3(A; Z_2)$ is the nonzero element in dimension 3, then $f^\ast(v)=u$. Then $f^\ast Sq^2 v=f^\ast(0)=0$, but $Sq^2 f^\ast(v)=Sq^2(u)\neq 0$. So $f^\ast Sq^2\neq Sq^2 f^\ast$ and the retraction can not exist.
\end{example}

The last example is an argument for trying Steenrod squares in data science applications. They are able to settle the existence of extensions that even cohomology can't handle. We will need to get a better idea of what they actually are, what properties they have, and how to calculate them. The rest of this chapter will answer those questions. 

\section{Cohomology Operations}

The material from here until the end of Section 11.7 is taken from Mosher and Tangora \cite{MT}. I will begin with a formal definition of a {\it cohomology operation}.

\begin{definition}
A {\it cohomology operation}\index{cohomology operation} of type $(\pi, n; G, m)$ is a family of functions $\theta_X: H^n(X; \pi)\rightarrow H^m(X; G)$, one for each space $X$ such that for any map $f: X\rightarrow Y$, $f^\ast\theta_Y=\theta_Xf^\ast.$ Cohomology operations are not necessarily homomorphisms. We denote the set of all cohomology operations of type $(\pi, n; G, m)$ as $\mathcal{O}(\pi, n; G, m)$
\end{definition}

\begin{example}
Let $R$ be a ring. Then the square $u\rightarrow u^2=u\cup u$ gives for each $n$ an operation $$H^n(X; R)\rightarrow H^{2n}(X; R).$$
\end{example}

\begin{definition}
For a space $X$ and $i>0$, define the {\it Hurewicz homomorphism}\index{Hurewicz homomorphism} $h: \pi_i(X)\rightarrow H_i(X)$ by choosing a generator $u$ of $H_i(S^i)\cong Z$ and letting $h([f])=f_\ast(u)$ where $f: S^i\rightarrow X$.
\end{definition}

Recall that the Hurewicz Theorem states that if $X$ is $(n-1)$-connected for $n>1$ then $h$ is an isomorphism in dimensions less than or equal to $n$ and $h$ is an epimorphism in dimension $n+1$. 

Now the Universal Coefficient Theorem for Cohomology gives an exact sequence $$0\rightarrow Ext(H_{n-1}(X), \pi)\rightarrow H^n(X; \pi)\rightarrow Hom(H_n(X), \pi)\rightarrow 0$$ for any space $X$. If $X$ is $(n-1)$-connected, then $H_{n-1}(X)=0$, so the Ext term goes away and we have $$H^n(X; \pi)\cong Hom(H_n(X), \pi).$$

In this case if $\pi=\pi_n(X)$, then the group $Hom(H_n(X), \pi_n(X))$ contains the inverse $h^{-1}$ of the Hurewicz homomorphism. The inverse exists in dimension $n$ since $X$ is $(n-1)$-connected.

\begin{definition}
Let $X$ be $(n-1)$-connected. The {\it fundamental class}\index{fundamental cohomology class} of $X$ is the cohomology class $\imath\in H^n(X; \pi_n(X))$ corresponding to $h^{-1}$ under the above isomorphism. The class may also be denoted by $\imath_X$ or $\imath_n$.
\end{definition}

In particular, the Eilenberg-Mac Lane space $K(\pi, n)$ has a fundamental class $\imath_n\in H^n(K(\pi, n); \pi)$.

Here is the main theorem.

\begin{theorem}
There is a one to one correspondence $[X, K(\pi, n)]\leftrightarrow H^n(X; \pi)$ given by $[f]\leftrightarrow f^\ast(\imath_n)$. (The square brackets represent homotopy classes.)
\end{theorem}

I will refer you to \cite{MT} for the proof, but I will point out that it makes heavy use of obstruction theory, both for extending maps and extending homotopies. So obstruction theory provided both the motivation and the tools for studying cohomology operations.

I will devote the rest of this section to some consequences of Theorem 11.2.1.

\begin{theorem}
There is a one to one correspondence $$[K(\pi, n), K(\pi', n)]\leftrightarrow Hom(\pi, \pi').$$
\end{theorem}

{\bf Proof:} By Theorem 11.2.1, $[K(\pi, n), K(\pi', n)]\leftrightarrow H^n(K(\pi, n); \pi')$. By the Universal Coefficient Theorem, $H^n(K(\pi, n); \pi')\cong Hom(H_n(K(\pi, n)), \pi')=Hom(\pi, \pi')$, since $H_n(K(\pi, n))=\pi$ by the Hurewicz Theorem. $\blacksquare$

We have seen earlier that the Eilenberg-Mac Lane spaces have the homotopy type of a CW complex which is unique up to homotopy.We will always assume we are working with CW complexes. We will abbreviate $H^m(K(\pi, n); G)$ as $H^m(\pi, n; G)$. 

Let $\theta$ be a cohomology operation of type $(\pi, n; G, m)$. Since $\imath_n\in H^n(\pi, n; \pi)$, we have that $\theta(\imath_n)\in H^m(\pi, n; G)$. 

\begin{theorem}
There is a one to one correspondence $$\mathcal{O}(\pi, n; G, m)\leftrightarrow H^m(\pi, n; G)$$ given by $\theta\leftrightarrow \theta(\imath_n).$
\end{theorem}

{\bf Proof:} The idea is to show that the function  $\theta\leftrightarrow \theta(\imath_n)$ has a two sided inverse.

Let $\phi\in H^m(\pi, n; G)$. We will use $\phi$ to denote an operation of type $(\pi, n; G, m)$ defined as follows. If $u\in H^n(X; \pi)$, let $\phi(u)=f^\ast(\phi)\in H^m(X; G)$, where $f: X\rightarrow K(\pi, n)$, and $[f]$ corresponds to $u$ in Theorem 11.2.1. 

We now have functions in both directions between $\mathcal{O}(\pi, n; G, m)$ and $H^m(\pi, n; G)$ and we need to see that their composition in either order is the identity. If $X=K(\pi, n)$ and $u=\imath_n$, then $f$ must be homotopic to the identity, so $\phi(\imath_n)=f^\ast(\phi)=\phi$. Going in the other direction, if $\phi=\theta(\imath_n)$, then the operation $\phi$ is equal to $\theta$ since $\phi(u)=f^\ast(\phi)=f^\ast(\theta(\imath_n))=\theta(f^\ast(\imath_n))=\theta(u)$ and we are done. $\blacksquare$

We will often use a symbol for both and operation and its corresponding cohomology class.

The following theorem is a direct result of Theorems 11.2.1 and 11.2.3. 

\begin{theorem}
There is a one to one correspondence $$\mathcal{O}(\pi, n; G, m)\leftrightarrow [K(\pi, n), K(G, m)].$$ 
\end{theorem}

By Theorem 11.2.3, the problem of finding all cohomology operations comes down to finding the cohmology groups of the appropriate Eilenberg-Mac Lane space. This should be easy as we know there are efficient ways to compute cohomology groups. The bad news is that our methods will not usually work. $K(Z, 1)$ is a circle, but Eilenberg-Mac Lane spaces tend to be infinite dimensional, so they are not representable as finite complexes, and our current methods would literaly take forever.

We will see how to calculate cohomology groups of some specific Eilenberg-Mac Lane spaces in Section 11.7. Until then we will focus on a specific type of cohomology operation that we have already mentioned: Steenrod squares.

\section{Construction of Steenrod Squares}

In this section, I will outline the construction of Steenrod squares found in Mosher and Tangora. These are cohomology operations of type $(Z_2, n; Z_2, n+i)$. 

\subsection{Cohomology of $K(Z_2, 1)$}

Recall that we saw in Example 9.8.2 that the infinite dimensional projective space $P^\infty$ is a $K(Z_2, 1)$ space. We will look at the homology and cohomology of this space. 

We start by looking at the cell structure. Letting $S^\infty=\cup_{n=0}^\infty S^n$, we can give $S^\infty$ the structure of a CW-complex with two cells in each dimension $i$ for $i\geq 0$ denoted $d_i$ and $Td_i$. These represent the Northern and Southern hemispheres with $T$ being the map that interchanges them. The boundary for homology is given by $$\partial d_i=d_{i-1}+(-1)^iTd_{i-1},$$ with $\partial T=T\partial$ and $T^2=1$. (Think of $S^2$ where $d_{i-1}$ and $Td_{i-1}$ are the two pieces of the equator.) 

Then $S^\infty$ is acyclic. This is because in even dimensions, the nonzero cycles are generated by $d_{2j}-Td_{2j}=\partial d_{2j+1}$, so cycles coincide with boundaries. In odd dimensions, the cycles are generated by $d_{2j-1}+Td_{2j-1}=\partial d_{2j}$. 

Now consider the homology of $P^\infty$. We obtain $P^\infty$ from $S^\infty$ by identifying $d_i$ and $Td_i$ so that there is now one cell in every dimension. We will call that cell $e_i$. Then the boundary formula becomes $\partial e_i=e_{i-1}+(-1)^i e_{i-1}$. So $\partial e_i=2e_{i-1}$ for $i$ even and $\partial e_i=0$ for $i$ odd. So the reduced homology $\tilde{H}_i(P^\infty)$ is $Z_2$ for $i$ odd and 0 for $i$ even. 

By the Universal Coefficient Theorems, since $H_\ast(P^\infty; Z)$ is $Z_2$ in odd dimensions, we have that $H^\ast(P^\infty; Z)$ is $Z_2$ in positive even dimensions. Then $H_\ast(P^\infty; Z_2)$ and $H^\ast(P^\infty; Z_2)$ are $Z_2$ in every dimension.

We next look at the ring structure. Let $W$ be the cellular chain complex of $S^\infty$. Then $W$ is a $Z_2$-free acyclic chain complex with two generators in each dimension $i$ for $i\geq 0$.

To get the ring structure of $P^\infty$ we need a diagonal map for $W$. (You may want to review Section 8.6 to see the connection.)  Let $T$ act on $W\otimes W$ by $T(u\otimes v)=T(u)\otimes T(v)$, where $T$ is the action that flips hemispheres described above.

Let $r: W\rightarrow W\otimes W$ be defined by $$r(d_i)=\sum_{j=0}^i(-1)^{j(i-j)}d_j\otimes T^jd_{i-j}$$ and $$r(Td_i)=T(r(d_i)).$$ Here $T^j=T$ if $j$ is odd and $T^j=1$ if $j$ is even, since $T^2=1$. Then $r$ is a chain map with respect to the boundary in $W\otimes W$. (Recall that $\partial(u\otimes v)=\partial u\otimes v +(-1)^{dim (u)} u\otimes \partial v.$)

Let $h$ denote the diagonal map of $Z_2$, so that $h(0)=(0, 0)$ and $h(1)=(1, 1)$. Then $r$ is {\it h-equivariant}, i.e. $r(gw)=h(g)r(w)$ for $g\in Z_2$ and $w\in W$. So $r$ induces a chain map $S: W/T\rightarrow W/T\otimes W/T$, where $W/T$ is the chain complex of $P^\infty$. This means that $$s(e_i)=\sum_{j=0}^i(-1)^{j(i-j)}e_j\otimes e_{i-j}.$$  Then $s$ is a chain approximation to the diagonal map $\Delta$ of $P^\infty$, and we can use it to find cup products in $H^\ast(P^\infty; Z_2)$. Let $\alpha_i$ be the nonzero element of $H^i(P^\infty; Z_2)\cong Z_2$.  The summation for $s(e_{j+k})$ contains a term $e_j\otimes e_k$ with coefficient 1 mod 2. Then \begin{align*}
\langle\alpha_j\cup \alpha_k, e_{j+k}\rangle &=\langle\Delta^\ast(\alpha_j\times \alpha_k), e_{j+k}\rangle\\
&=\langle\alpha_j\times \alpha_k, s(e_{j+k})\rangle\\
& =1 \mod 2.
\end{align*} Thus, $\alpha_j\cup\alpha_k=\alpha_{j+k}$. This proves Theorem 8.3.4 which states that $H^\ast(P^\infty; Z_2)$ is a polynomial ring generated by the nonzero element $u\in H^1(P^\infty; Z_2)$.

\subsection{Acyclic Carriers}

For the construction of Steenrod squares, we will need a theorem on {\it acyclic carriers}\index{acyclic carrier}. First we will need some definitions. In what follows, $\pi$ and $G$ are groups.

\begin{definition}
The {\it group ring}\index{group ring} $Z[\pi]$ consists of formal sums $\sum z_i g_i$ where $z_i\in Z$ and $g_i\in \pi$. Addition is carried out component wise and multiplication using the distributive law, the usual multiplication in $Z$, and the group multiplication for $\pi$. For example,  $$(z_1 g_1 + z_2 g_2)(z_3 g_3 + z_4 g_4)=z_1 z_3 g_1 g_3 + z_1 z_4 g_1 g_4 + z_2 z_3 g_2 g_3 + z_2 z_4 g_2 g_4.$$
\end{definition}

Let $K$ be a $\pi$-free chain complex with a $Z[\pi]$-basis $B$ of elements called {\it cells}. For two cells, $\sigma$ and $\tau$, let $[\tau: \sigma]\in Z[\pi]$ be the coefficient of $\sigma$ in $\partial \tau$. Let $L$ be a chain complex acted on by $G$, and let $h: \pi\rightarrow G$ be a homomorphism.

\begin{definition}
An {\it h-equivariant carrier}\index{equivariant carrier} $\mathcal{C}$ from $K$ to $L$ is a function $\mathcal{C}$ from the basis $B$ of $K$ to the subcomplexes of $L$ such that: \begin{enumerate}
\item If $[\tau: \sigma]\neq 0$, then $\mathcal{C}(\sigma)\subset\mathcal{C}(\tau).$
\item If $x\in\pi$ and $\sigma\in B$, then $h(x)\mathcal{C}(\sigma)\subset\mathcal{C}(\sigma).$
\end{enumerate}

$\mathcal{C}$ is an {\it acyclic carrier}\index{acyclic carrier} if the subcomplex $\mathcal{C}(\sigma)$ is acyclic for every cell $\sigma\in B$. The $h$-equivariant chain map $f: K\rightarrow L$ is {\it carried} by $\mathcal{C}$ if $f(\sigma)\in \mathcal{C}(\sigma)$ for every $\sigma\in B$.
\end{definition}

\begin{theorem}
{\bf (Equivariant) Acyclic Carrier Theorem:} Let $\mathcal{C}$ be an acyclic carrier from $K$ to $L$. Let $K'$ be a subcomplex of $K$ which is a $Z[\pi]$-free complex on a subset of $B$. Let $f: K'\rightarrow L$ be an $h$-equivariant chain map carried by $\mathcal{C}$. Then $f$ extends over all of $K$ to an $h$-equivariant chain map carried by $\mathcal{C}$. This extension is unique up to an $h$-equivariant chain homotopy carried by $\mathcal{C}$. 
\end{theorem}

{\bf Proof:} We proceed by induction on dimension. Suppose that $f$ has been extended over $K^q$. Let $\tau\in B$ be a $(q+1)$-cell. Then $\partial\tau=\sum a_i\sigma_i$ where $a_i=[\tau: \sigma_i]\in Z(\pi)$. Then $f(\partial\tau)=\sum f(a_i\sigma_i)=\sum h(a_i)f(\sigma_i)$,which is in $\mathcal{C}(\tau)$ by the definition of an $h$-equivariant carrier. Since $f$ is a chain map, $f(\partial\tau)=\partial f(\tau)$ is a cycle, but it is also a boundary since $\mathcal{C}(\tau)$ is acyclic. So there exists $x\in \mathcal{C}(\tau)$ such that $\partial x=f(\partial \tau)$. Choose any such $x$ and let $f(\tau)=x$. Then $f$ is extended over $K^{q+1}$ by ensuring it is $h$-equivariant. Uniqueness is proved by applying the construction to the complex $K\times I$ and the subcomplex $K'\times I\cup K\times \dot{I}$. $\blacksquare$

\subsection{The Cup-i Products}

We will next see how the cup-i products mentioned in Section 11.1 are derived using the Acyclic Carrier Theorem. We will still have a little way to go as this will involve cohomology with integer coefficients.

Let $K$ be the chain complex of a simplicial complex and let $W$ be the $Z_2$-free complex of $S^\infty$ described earlier. Let $T$ generate the action of $Z_2$ on $W\otimes K$ by $T(w\otimes k)=T(w)\otimes k$, and on $K\otimes K$ by $T(x\otimes y)=(-1)^{dim(x)dim(y)}(y\otimes x)$. Let $K=C(X)$ be the chain complex of the simplicial complex $X$. By the Eilenberg-Zilber Theorem (Theorem 8.5.5), there is a chain-homotopy equivalence $\Psi: C(X\times X)\rightarrow C(X)\times C(X).$ For a generator $\sigma$ of $K$ (i.e a simplex of $X$),$\Psi(C(\sigma\times\sigma))$ is a subcomplex of $C(X)\otimes C(X)$. Then let the carrier $\mathcal{C}: W\otimes K\rightarrow K\otimes K$ be defined as $\mathcal{C}(d_i\otimes\sigma)=\Psi(C(\sigma\times\sigma))$. $\mathcal{C}$ is  acyclic and $h$-equivariant, where $h$ is the identity map of $Z_2$. Then there exists an $h$-equivariant chain map $\phi: W\otimes K\rightarrow K\otimes K$ carried by $\mathcal{C}$.

The map $\phi$ is what we need for the cup-i products. If we restrict $\phi$ to $d_0\otimes K$, we can view it as a map $\phi_0: K\rightarrow K\otimes K$ which is carried by the diagonal carrier. So it is chain homotopic to the diagonal map and can be used to compute cup products in $K$. (Review Section 8.6.) The same can be said of $T\phi_0: \sigma\rightarrow \phi(Td_0\otimes\sigma)$. So $\phi_0$ and $T\phi_0$ are both carried by $\mathcal{C}$, so they are equivariantly homotopic. A chain homotopy is given by $\phi_1: K\rightarrow K\otimes K$ where $\phi_1(\sigma)=\phi(d_1\otimes \sigma)$. Then $\phi_1$ and $T\phi_1$ are equivariantly homotopic homotopies and a chain homotopy is given by an analogously defined $\phi_2$.

\begin{definition}
For each integer $i\geq 0$, define a {\it cup-i product}\index{cup-i product} $$C^p(K)\otimes C^q(K)\rightarrow C^{p+q-i}(K)\mbox{ with }(u, v)\rightarrow u\cup_i v$$ by the formula $$(u\cup_i v)(c)=(u\otimes v)\phi(d_i\otimes c)$$ for $c\in C_{p+q-i}(K).$ The definition depends on a specific choice of $\phi$ but it is independent of that choice as we will soon see below.
\end{definition}

Mosher and Tangora give another proof of the coboundary formula (Theorem 11.1.1) with this definition. It is not hard but messy so I will refer you to \cite{MT} for details.

\subsection{Steenrod Squares}

As we mentioned in Section 11.1, if $u\in C^p(K)$ is a cocycle mod 2, (i.e. $\delta u=2a$ for some $a\in C^{p+1}(K)$), then $u \cup_i u$ is also a cocycle mod 2. As before, we get a map $Sq_i: Z^p(K; Z_2)\rightarrow Z^{2p-i}(K ; Z_2)$ taking $u$ to $u\cup_i u$. This operation also takes coboundaries to coboundaries, so it passes to a (group) homomorphism: $Sq_i: H^p(K; Z_2)\rightarrow H^{2p-i}(K ; Z_2)$.

\begin{theorem}
Let $f: K\rightarrow L$ be a continuous map. Then the following diagram commutes:

$$\begin{tikzpicture}
  \matrix (m) [matrix of math nodes,row sep=3em,column sep=4em,minimum width=2em]
  {
H^p(L; Z_2) & H^{2p-i}(L; Z_2)\\
H^p(K; Z_2) &  H^{2p-i}(K; Z_2)\\};

\path[-stealth]
(m-1-1) edge node [above] {$Sq_i$} (m-1-2)
(m-1-1) edge node [left] {$f^\ast$} (m-2-1)
(m-1-2) edge node [right] {$f^\ast$} (m-2-2)
(m-2-1) edge node [below] {$Sq_i$} (m-2-2)

;

\end{tikzpicture}$$
\end{theorem}

{\bf Proof:} We can assume $f$ is simplicial as otherwise we could substitute it with a simplicial approximation. Let $u\in C^p(L)$, and $c\in C_{2p-i}(K)$. Then we have $$f^\ast(Sq_i(u))(c)=(u\otimes u)\phi_L(d_i\otimes f(c))=(u\otimes u)\phi_L(1\otimes f)(d_i\otimes c)$$ $$Sq_i(f^\ast(u))(c)=(f^\ast u\otimes f^\ast u)\phi_K(d_i\otimes c)=(u\otimes u)(f\otimes f)\phi_K(d_i\otimes c).$$ Now the two chain maps $\phi_L(1\otimes f)$ and $(f\otimes f)\phi_K$ are both carried by the acyclic carrier $\mathcal{C}$ from $W\otimes K$ to $L\otimes L$ defined by $\mathcal{C}(d_i\otimes \sigma)=C(f(\sigma)\otimes f(\sigma)).$ So they are equivariantly chain homotopic and the images of the two functions above are cohomologous. $\blacksquare$

\begin{theorem}
The operation $Sq_i$ is independent of the choice of $\phi$.
\end{theorem}

{\bf Proof:} In the previous theorem, let $K=L$ and $f$ be the identity. Letting $\phi_L$ and $\phi_K$ be two different choices of $\phi$ implies the result. $\blacksquare$

Recall that in defining cup-i products, we said that $(u\cup_0 u)(c)=(u\otimes u)\phi_0(c)$ and $\phi_0$ is chain homotopic to the diagonal map and thus suitable for computing cup products. This shows the following.

\begin{theorem}
If $u\in C^p(X; Z_2)$ then $Sq_0(u)=u^2\equiv  u\cup u$.
\end{theorem}

Recall that $Sq_i$ was the notation used in Steenrod's original paper \cite{Ste2}. From now on, we will use the more modern notation $Sq^i$ defined below.

\begin{definition}
We define $Sq^i: H^p(X; Z_2)\rightarrow H^{p+i}(X; Z_2)$ for $0\leq i\leq p$ by $Sq^i=Sq_{p-i}$. For values of i outside this range, define $Sq^i=0$. We refer to $Sq^i$ as a {\it Steenrod square}\index{Steenrod square}\index{$Sq^i$}.
\end{definition}

$Sq^i$ raises the cohomology dimension by $i$. If $\{u\}\in H^p(X; Z_2)$, then $Sq^p(u)=Sq_{p-p}(u)=Sq_0(u)=u^2$ by Theorem 11.3.4. The name "Steenrod square" is derived from this fact.

We can also define Steenrod squares for relative cohomology. Consider the exact sequence $$0\rightarrow C^\ast(K, L)\xrightarrow{q^\ast}C^\ast(K)\xrightarrow{j^\ast}C^\ast(L)\rightarrow 0,$$ where $j$ and $q$ are inclusions. We can let $\phi_L=\phi_K | W\otimes L,$ since $\phi_k(d_i\otimes \sigma)\in C(\sigma\times\sigma)\subset L\otimes L$ for $\sigma\in L$. Then for $u, v\in C^\ast(K), j^\ast(u\cup_i v)=j^\ast(u)\cup_i j^\ast(v)$. To define a relative cup-i product, let $u, v\in C^\ast(K, L)$. Then $ j^\ast(q^\ast u\cup_i q^\ast v)=j^\ast q^\ast(u)\cup_i j^\ast q^\ast(v)=0$ since $j^\ast q^\ast=0$. But by exactness, we have that $(q^\ast u\cup_i q^\ast v)$ is in the image of $q^\ast$. Since $q^\ast$ is a monomorphism, we can define $u\cup_i v$ the unique cochain in $C^\ast(K, L)$ such that $q^\ast(u\cup_i v)=q^\ast u\cup_i q^\ast v$. The coboundary formula stays the same in this case, so we get a homomorphism $$Sq^i: H^p(K, L; Z_2)\rightarrow H^{p+i}(K, L; Z_2)$$ and $q^\ast Sq^i=Sq^i q^\ast$.

Letting $\delta^\ast$ be the coboundary homomorphism  $\delta^\ast: H^p(L)\rightarrow H^{p+1}(K, L)$ in the exact cohomology sequence of the pair $(K, L)$, Mosher and Tangora show the following.

\begin{theorem}
$Sq^i$ commutes with $\delta^\ast$ as in the diagram below.

$$\begin{tikzpicture}
  \matrix (m) [matrix of math nodes,row sep=3em,column sep=4em,minimum width=2em]
  {
H^p(L; Z_2) & H^{p+i}(L; Z_2)\\
H^{p+1}(K, L; Z_2) &  H^{p+i+1}(K, L; Z_2)\\};

\path[-stealth]
(m-1-1) edge node [above] {$Sq^i$} (m-1-2)
(m-1-1) edge node [left] {$\delta^\ast$} (m-2-1)
(m-1-2) edge node [right] {$\delta^\ast$} (m-2-2)
(m-2-1) edge node [below] {$Sq^i$} (m-2-2)

;

\end{tikzpicture}$$
\end{theorem}

We conclude this section with one more property that will provide an interesting example. Recall that the {\it suspension} of $X$ is the result of taking $X\times I$ and collapsing each of $X\times \{0\}$ and $X\times \{1\}$ to single points. By an argument similar to the proof of Theorem 4.2.1, there is an isomorphism $S^\ast: \tilde{H}^p(X)\approx \tilde{H}^{p+1}(SX)$. 

\begin{theorem}
The Steenrod squares commute with suspension, i.e. $$Sq^i S^\ast=S^\ast Sq^i: \tilde{H}^p(X)\rightarrow \tilde{H}^{p+i+1}(SX).$$
\end{theorem}

\begin{example}
This example shows a nontrivial Steenrod square which is not a cup product. Let $P^2$ be the projective plane. $H^\ast(P^2; Z_2)$ is the truncated polynomial ring generated by the nonzero class $u\in H^1(P^2; Z_2)\cong Z_2$. We have that $u^3=0$. Now $Sq^1(u)=Sq_0(u)=u^2$, so $S^\ast Sq^1(u)\neq 0$. By Theorem 11.3.6, $Sq^1 S^\ast(u)$ is also nonzero, so the operation $Sq^1$ is nonzero. Thus, the operation $Sq^1$ is nontrivial in $H^2(SP^2; Z_2)$.
\end{example}

\section{Basic Properties}

\begin{theorem}
The Steenrod squares $Sq^i$ for $i\geq 0$ have the following properties:\begin{enumerate}
\item $Sq^i$ is a natural homomorphism $H^p(K, L; Z_2)\rightarrow H^{p+i}(K, L; Z_2).$
\item If $i>p$, $Sq^i(x)=0$ for $x\in H^p(K, L; Z_2)$.
\item If $x\in H^p(K, L; Z_2)$, then $Sq^p(x)=x^2=x\cup x$.
\item $Sq^0$ is the identity homomorphism.
\item $Sq^1$ is the {\it Bockstein homomorphism} defined below.
\item $\delta^\ast Sq^i=Sq^i\delta^\ast$ where $\delta^\ast$ is the coboundary homomorphism $\delta^\ast: H^p(L; Z_2)\rightarrow H^{p+1}(K, L; Z_2)$.
\item {\bf Cartan Formula:} Writing $xy$ for $x\cup y$ we have $$Sq^i(xy)=\sum_j(Sq^jx)(Sq^{i-j}y).$$
\item {\bf Adem Relations:} For $a<2b$, $$Sq^aSq^b=\sum_c\dbinom{b-c-1}{a-2c}Sq^{a+b-c}Sq^c,$$ where the binomial coefficient is taken mod 2.
\end{enumerate}
\end{theorem}

We have already seen properties 1, 2, 3, and 6. Chapter 3 of Mosher and Tangora \cite{MT} is devoted to proving the other properties given our definitions. Rather than give all of the details which involve some lengthy calculations, I will outline some of the interesting points. See \cite{MT} for more details. 

Note that these properties can be taken as axioms that completely characterize the squares. This is done in Steenrod and Epstein \cite{SE}.

\subsection{$Sq^1$ and $Sq^0$}

First of all, I will define $Sq^1$ as promised. Consider the exact sequence $$0\rightarrow Z\xrightarrow{m} Z\rightarrow Z_2\rightarrow 0,$$ where $m$ is multiplication by 2. Define the {\it Bockstein homomorphism}\index{Bockstein homomorphism} $\beta: H^p(K, L; Z_2)\rightarrow H^{p+1}(K, L; Z)$ as follows. Let $x\in H^p(K, L; Z_2)$. Represent $x$ by a cocycle $c$ and choose an integral cochain $c'$ which maps to $c$ under reduction mod 2. Then $\delta c'$ is a multiple of 2, and $\frac{1}{2}(\delta c')$ represents $\beta x$.  The composition of $\beta$ and the reduction homomorphism  $H^{p+1}(K, L; Z)\rightarrow H^{p+1}(K, L; Z_2)$ gives a homomorphism $$\delta_2:  H^p(K, L; Z_2)\rightarrow H^{p+1}(K, L; Z_2)$$ which is also called the {\it Bockstein homomorphism}. I claim that $Sq^1=\delta_2$.

To outline the proof, we start with the following. It is also a special case of the Adem Relations.

\begin{theorem}
$\delta_2 Sq^j=0$ if $j$ is odd, and $\delta^2 Sq^j=Sq^{j+1}$ if $j$ is even.
\end{theorem}

An immediate consequence is that Property 3 implies Property 4. See \cite{MT} for the proof of Theorem 11.4.2. I will give the proof of Property 3.

First of all, let $K=P^2$. Then if $u$ is the nonzero element of $H^1(P^2, Z_2)$, $\delta^2 Sq^0(u)=Sq^1(u)=u^2\neq 0.$ So $Sq^0(u)\neq 0$ and $Sq^0(u)=u$, since $u$ is the only nonzero element of $H^1(P^2, Z_2)$. So Property 3 holds in this case. 

If $K=S^1$, then letting $f: S^1\rightarrow P^2$ be such that $f^\ast(u)=\sigma$ where $\sigma$ generates $H^1(S^1; Z_2)$, we have $$Sq^0(\sigma)=Sq^0(f^\ast(u))=f^\ast(Sq^0(u))=f^\ast(u)=\sigma.$$ Using suspensions and Theorem 11.3.6, we have that Property 3 holds for any sphere $S^n$. We extend this to any complex $K$ of dimension $n$ by mapping $K$ to $S^n$ so that if $\sigma$ generates $H^n(S^n; Z_2)$ we let $f^\ast(\sigma)$ map to a given $n$-dimensional cohomology class of $K$. We then generalize to an infinite dimensional complex $K$ through inclusion of $K^n$ into $K$ and finally to pairs via an excision argument. (See \cite{MT} for the details in these cases.)

\subsection{The Cartan Formula}

I won't prove the Cartan Formula, but an interesting observation is that while the squares $Sq^i$ are homomorphisms with respect to group addition, they are obviously not ring homomorphisms. Otherwise we would have $Sq^i(xy)=Sq^i(x)Sq^i(y)$, which contradicts the Cartan Formula. On the other hand, there is something we can do. 

\begin{definition}
Define $Sq: H^\ast(K; Z_2)\rightarrow H^\ast(K; Z_2)$ by $$Sq(u)=\sum_iSq^iu.$$  This sum is finite and the image $Sq(u)$ is not homogeneous. In other words it is the sum of elements of differing dimensions. 
\end{definition}

\begin{theorem}
$Sq$ is a ring homomorphism.
\end{theorem}

{\bf Proof:} $Sq(u)\cup Sq(v)=(\sum_i Sq^iu)\cup (\sum_j Sq^j v)$ has $Sq^i(u\cup v)$ as its $(p+q+i)$-dimensional term by the Cartan formula. (Here $u$ has dimension $p$ and $v$ has dimension $q$.) Then $Sq(u)\cup Sq(v)=Sq(u\cup v)$. $\blacksquare$

\begin{theorem}
For $u\in H^1(K; Z_2)$, we have $Sq^i(u^j)=\binom{j}{i}u^{j+1}.$
\end{theorem}

{\bf Proof:} $Sq(u)=Sq^0(u)+Sq^1(u)=u+u^2.$ Since $Sq$ is a ring homomorphism, $Sq(u^j)=(u+u^2)^j=u^j\sum_k\binom{j}{k}u^k$. Comparing coefficients proves the result. $\blacksquare$

In particular, this gives the action of all $Sq^i$ on $H^\ast(P^\infty; Z_2).$

\subsection{Cartesian Products of $P^\infty$}

This subsection is necessary if you want to understand more of the theory from Mosher and Tangora in later chapters. It is especially needed to understand the structure of the Steenrod algebra that I will discuss in Section 11.6.

Let $K=P^\infty=K(Z_2, 1).$ We will write $K_n=K\times\cdots\times K$ ($n$ times). Since $H^\ast (K; Z_2)$ is the polynomial ring on the one dimensional class, the K\"{u}nneth Theorem implies that $H^\ast(K_n; Z_2)$ is the polynomial ring over $Z_2$ generated by $x_1, \cdots, x_n$, where $x_i$ is the nonzero one dimensional class in the $i$-th copy of $K$. (Recall that for coefficients in a field, the Tor terms go away. See Example 8.5.2.) 

\begin{definition}
The {\it elementary symmetric polynomials}\index{elementary symmetric polynomials} in $n$ variables $x_1, \cdots, x_n$ are defined as follows:\begin{align*}
\sigma_0&=1\\
\sigma_1&=\sum_{1\leq i\leq n} x_i\\
\sigma_2&=\sum_{1\leq i<j\leq n} x_i x_j\\
\sigma_3&=\sum_{1\leq i<j<k\leq n} x_i x_j x_k\\
\cdots&\cdots\cdots\\
\sigma_n&=x_1 x_2\cdots x_n\\
\end{align*}
\end{definition}

\begin{definition}
A {\it symmetric polynomial}\index{symmetric polynomial} in $n$ variables is a polynomial where interchanging two variables produces the same polynomial.
\end{definition}

 Obviously the elementary symmetric polynomial are all symmetric polynomials. Let $S$ be the subring of $Z_2[x_1, \cdots, x_n]$ consisting of the symmetric polynomials.

\begin{theorem}
{\bf Fundamental Theorem of Symmetric Algebra:} Let $R$ be a commutative ring with unit element and $S$ the subring of the polynomial ring $R[x_1, \cdots, x_n]$ consisting of the symmetric polynomials. Then $S$ is equal to the polynomial ring $R[\sigma_1, \cdots, \sigma_n]$ where $\sigma_i$ is the $i$th elementary symmetric polynomial in $n$ variables. (Note that $\sigma_0=1$ is excluded.)
\end{theorem}

So in our case, $S=Z_2[\sigma_1, \cdots, \sigma_n]$. First, we will see what the squares do to $\sigma_i$.

\begin{theorem}
In $H^\ast(K_n; Z_2)$, $Sq^i(\sigma_n)=\sigma_n\sigma_i$ for $0\leq i\leq n$.
\end{theorem}

{\bf Proof:} $$Sq(\sigma_n)=Sq(x_1 x_2\cdots x_n)=Sq(x_1) Sq(x_2)\cdots Sq(x_n)=\prod_{i=1}^n(x_i+x_i^2)=\sigma_n(\prod_{i=1}^n(1+x_i))=\sigma_n\sum_{i=0}^n\sigma_i.$$ The term of degree $n+i$ is then $\sigma_n\sigma_i$. $\blacksquare$

As in Section 11,2, we will abbreviate $H^\ast(K(\pi, n); G)$ as $H^\ast(\pi, n; G)$. Let $\imath_n\in H^n(Z_2, n; Z_2)$, be the fundamental class. (See Definition 11.2.3.)

\begin{theorem}
In $H^\ast(Z_2, n; Z_2)$, $Sq^i(\imath_n)\neq 0$ for $0\leq i\leq n$. 
\end{theorem}

{\bf Proof:} By Theorem 11.2.1, there is a map $f: K_n\rightarrow K(Z_2, n)$ such that $f^\ast(\imath_n)=\sigma_n$. Then $$f^\ast(Sq^i(\imath_n))=Sq^i f^\ast(\imath_n)=Sq^i(\sigma_n)=\sigma_n\sigma_i\neq 0,$$ and we are done.  $\blacksquare$

We can find more linearly independent elements of $H^\ast(Z_2, n; Z_2)$ by using compositions of the squares.

Let $I=\{i_1, \cdots, i_r\}$ be a sequence of positive integers. Then $Sq^I$ denotes the composifion $Sq^{i_1} Sq^{i_2} \cdots Sq^{i_r}$. If $I=\{\}$ is the empty sequence, we let $Sq^I=Sq^0$. 

The following definition will be very important in Section 11.6 when we talk about the Steenrod algebra.

\begin{definition}
A sequence $I$ as described above is {\it admissible}\index{admissible sequence} if $i_j\geq 2(i_{j+1})$ for every $j<r$. This is automatically satisfied if $I$ is empty or $I=\{i_1\}$. If $I$ is an admissible sequence, we also refer to $Sq^I$ as admissible. The {\it length} of $I$ is the number $r$ of terms. The degree $d(I)$ is the sum of the terms $\sum_{j=1}^r i_j$. So $Sq^I$ raises the dimension by $d(I)$. For an admissible sequence $I$, the {\it excess} is $e(I)=2i_1-d(I)$.
\end{definition}

\begin{theorem}
If $e(I)$ is the excess, then $$e(I)=2i_1-d(I)=i_1-i_2-\cdots-i_r=(i_1-2i_2)+(i_2-2i_3)+\cdots+i_r.$$
\end{theorem}

Again, let $S$ be the symmetric polynomial subring of the polynomial ring $H^\ast(K_n; Z_2)=Z_2[x_1,\cdots, x_n]$. Define an ordering of the monomials of $S$ as follows. Given any such monomial, write it as $$m=\sigma_{j_1}^{e_1}\sigma_{j_2}^{e_2}\cdots \sigma_{j_s}^{e_s}$$ such that $j_1>j_2>\cdots>j_s$. Then put $m<m'$ if $j_1<j_1'$ or if $j_1=j_1'$ and $(m/\sigma_{j_1})<(m'/\sigma_{j_1}).$

\begin{theorem}
If $d(I)\leq n$, then $Sq^I(\sigma_n)$ can be written as $\sigma_n Q_I$ where $$Q_I=\sigma_{i_1}\cdots\sigma_{i_r}+\mbox{ (a sum of monomials of lower order)}.$$
\end{theorem}

The result is proved using induction on the length of $I$, the Cartan Formula, and Theorem 11.4.6.

As $I$ runs through all the admissible sequences of degree $\leq n$, the monomials $\sigma_I=\sigma_{i_1}\sigma_{i_2}\cdots\sigma_{i_r}$ are linearly independent in $S$, and hence in $H^\ast(K_n; Z_2)$. The theorem shows that the $Sq^I(\sigma_n)$ are also linearly independent. Choosing a map $f: K^n\rightarrow K(Z_2, n)$ such that $f^\ast(\imath_n)=\sigma_n.$ shows the following.

\begin{theorem}
As $I$ runs through all the admissible sequences of degree $\leq n$, the elements $Sq^I(\imath_n)$ are linearly independent in $H^\ast(Z_2, n; Z_2)$.
\end{theorem}

\begin{theorem}
If $u\in H^n(K; Z_2)$ for some space $K$ and $I$ has excess $e(I)>n$, then $Sq^Iu=0$. If $e(I)=n$, then $Sq^Iu=(Sq^Ju)^2$, where $J$ denores the sequence obtained from $I$ by dropping $i_1$.
\end{theorem}

{\bf Proof:} If $e(I)=i_1-i_2-\cdots-i_r>n$, then $i_1>n+i_2+\cdots+i_r=\mbox{dim }(Sq^J)(u)$. So by Property 2 of Theorem 11.4.1, $Sq^Iu=0$. 

If $e(I)=n$, then we replace the greater than sign with an equal sign and $i_1=\mbox{dim }(Sq^J)(u)$.  By Property 3 of Theorem 11.4.1,  $Sq^Iu=(Sq^Ju)^2$. $\blacksquare$

These results are included in a theorem of Serre which states that $H^\ast(Z_2, n; Z_2)$ is the polynomial ring over $Z_2$ with generators $\{Sq^I(\imath_n)\}$ as $I$ runs through all admissible sequence of excess less than $n$. Mosher and Tangora use this result to prove the Adem relations but don't prove it until later using the machinery of spectral sequences. For completeness. I will discuss it and some related results in Section 11.7.

\subsection{Adem Relations}

The Adem Relations (Property 8 of Theorem 11.4.1) relate to what happens when you switch the order of the compostion of 2 Steenrod squares. Rather than copy Mosher and Tangora's very long proof, I will focus on showing you how to use them.

Letting $\lfloor x\rfloor$ denote the floor of $x$ or the largest integer $z\leq x$ an Adem relation has the form $$R=Sq^aSq^b+\sum_{c=0}^{\lfloor \frac{a}{2}\rfloor}\binom{b-c-1}{a-2c}Sq^{a+b-c}Sq^c\equiv 0\mod 2,$$ where $a<2b$. We use the convention that the binomial coefficient $\binom{x}{y}=0$ if $y<0$ or $x<y$.

Let's try some examples.

\begin{example}
$$Sq^2 Sq^6=Sq^7 Sq^1$$

{\bf Proof:} Let $a=2, b=6$. Then since $2<2(6)=12$, we can use the Adem relation \begin{align*}
Sq^2Sq^6&=\sum_{c=0}^1 \binom{6-c-1}{2-2c}Sq^{8-c}Sq^c\\
&=\binom{6-0-1}{2-2(0)}Sq^8 Sq^0+\binom{6-1-1}{2-2(1)} Sq^7 Sq^1\\
&=\binom{5}{2}Sq^8 Sq^0+ \binom{4}{0} Sq^7 Sq^1\\
&=10 Sq^8 Sq^0+ Sq^7 Sq^1\\
&=Sq^7 Sq^1,
\end{align*} since we are working modulo 2. $\blacksquare$
\end{example}

\begin{example}
$$Sq^2 Sq^4=Sq^6+Sq^5 Sq^1$$ 

{\bf Proof:} Let $a=2, b=4$. Then since $2<2(4)=8$, we can use the Adem relation \begin{align*}
Sq^2Sq^4&=\sum_{c=0}^1 \binom{4-c-1}{2-2c}Sq^{6-c}Sq^c\\
&=\binom{4-0-1}{2-2(0)}Sq^6 Sq^0+\binom{4-1-1}{2-2(1)} Sq^5 Sq^1\\
&=\binom{3}{2}Sq^6 Sq^0+ \binom{2}{0} Sq^5 Sq^1\\
&=3 Sq^6 Sq^0+ Sq^5 Sq^1\\
&=Sq^6+Sq^5 Sq^1,
\end{align*} since we are working modulo 2 and $Sq^0$ is the identity so it can be dropped. $\blacksquare$
\end{example}

\begin{example}
$$Sq^{2n-1} Sq^{n}=0$$ 

{\bf Proof:} Let $a=2n-1, b=n$. Then since $2n-1<2n$, we can use the Adem relation 
$$Sq^{2n-1} Sq^{n}=\sum_{c=0}^{\lfloor\frac{2n-1}{2}\rfloor} \binom{n-c-1}{2n-1-2c}Sq^{3n-1-c}Sq^c.$$

Now for the binomial coefficient to be nonzero, we need $2n-1-2c>0$, so $c<\lfloor\frac{2n-1}{2}\rfloor=n-1$. This condition is already taken care of in the summation. But we also need \begin{align*}
2n-1-2c&<n-c-1\\
n-2c&<-c\\
n<c.
\end{align*} But this is impossible since $c<n-1$. Thus $Sq^{2n-1} Sq^{n}=0$.
$\blacksquare$
\end{example}

\section{The Hopf Invariant}

The next topic has little to do with data science and is a little tangential but a rather interesting application of Steenrod squares. I will outline some of the arguments and refer the reader to Chapter 4 of Mosher and Tangora \cite{MT} for the details. 

Consider the sphere $S^n$ where $n\geq 2$. Let $f: S^{2n-1}\rightarrow S^n$ be given. Let $S^{2n-1}$ be the boundary of an oriented $2n$-cell and form the cell complex $K=S^n\cup_f e^{2n}$. Recall that this means we added a $2n$-cell to $K$ with the boundary attached to $S^ n$ by $f$. The integral cohomology of $K$ is zero except in dimensions $0, n$, and $2n$, and $H^i(K; Z)\cong Z$ in those dimensions. Let $\sigma$ and $\tau$ generate $H^n(K; Z)$ and $H^{2n}(K; Z)$ respectively. Then $\sigma^2=\sigma\cup\sigma$ is an integral multiple of $\tau$.

\begin{definition}
The {\it Hopf invariant}\index{Hopf invariant} of $f$ is the integer $H(f)$ such that $\sigma^2=H(f)\tau$.
\end{definition}

{\bf Question:} Under what circumstances does $H(f)=1$?

The homotopy type of $K$ depends only on the homotopy class of $f$. So the Hopf invariant defines a map $H: \pi_{2n-1}(S^n)\rightarrow Z$ where the class $\alpha\in\pi_{2n-1}(S^n)$ is taken to $H(f)$ where $f$ is any map representing $\alpha$.

If $n$ is odd, then by anti-commutativity of cup product (Theorem 8.6.3), $$\sigma\cup\sigma=(-1)^{n^2}\sigma\cup\sigma=-\sigma\cup\sigma,$$ so $\sigma^2=0$. So for $n$ odd, $H(f)=0$.

If $n=2, 4,$ or $8$, we have the Hopf maps $S^3\rightarrow S^2$, $S^7\rightarrow S^4$, and $S^{15}\rightarrow S^8$. (See Section 9.4.) These maps are known to have Hopf invariant one as they turn out to be attaching maps of $2n$ cells in the complex, quaternionic, and Cayley projective planes respectively and these spaces have cohomology rings that are truncated polynomial rings so that powers of generators are also generators. So $\sigma^2=\tau.$

\begin{theorem}
If $n$ is even, there exists a map $f: S^{2n-1}\rightarrow S^n$ with Hopf invariant 2.
\end{theorem}

\begin{figure}[ht]
\begin{center}
  \scalebox{0.8}{\includegraphics{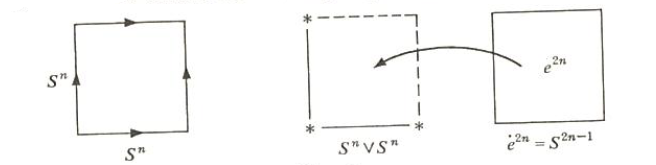}}
\caption{
\rm
Attaching a $2n$-cell to $S^n\vee S^n$ to form $S^n\times S^n$. \cite{MT}. 
}
\end{center}
\end{figure}

It takes some work to prove this, but the map itself is the map $Fg: S^{2n-1}\rightarrow S^n$ formed as follows: Consider the product $S^n\times S^n$ as a cell complex formed by attaching a $2n$-cell to $S^n\vee S^n$. Let $g: S^{2n-1}\rightarrow S^n\vee S^n$ be the attaching map shown in Figure 11.5.1 \cite{MT}. Let $F: S^n\vee S^n\rightarrow S^n$ be the folding map. Then Mosher and Tangora show that $Fg$ has Hopf invariant 2.

\begin{theorem}
The transformation $H: \pi_{2n-1}(S^n)\rightarrow Z$ is a homomorphism, i.e. if $f, g$ represent elements of $\pi_{2n-1}(S^n)$, then $H(f+g)=H(f)+H(g)$.
\end{theorem}

Again, I will refer you to \cite{MT} for the proof, but we immediately have the following.

\begin{theorem}
If $n$ is even, then $\pi_{2n-1}(S^n)$ contains $Z$ as a direct summand.
\end{theorem}

Now we will work towards showing that there does not exist a map $f:  S^{2n-1}\rightarrow S^n$ with Hopf invariant 1 if $n$ is not a power of 2. The key is the fact that $f$ has an odd Hopf invariant, so $\sigma^2$ is an odd multiple of $\tau$ in $K=S^n\cup_fe^{2n}$. But then, since $\sigma^2=Sq^n\sigma$, we have $Sq^n\sigma=\tau$ in the mod 2 cohomology of $K$.

\begin{definition}
$Sq^i$ is said to be {\it decomposable}\index{decomposable Steenrod square} if $$Sq^i=\sum_{t<i}a_t Sq^t,$$ where each $a^t$ is a sequence of squaring operations. If $Sq^i$ is not decomposable, then it is {\t indecomposable}\index{indecomposable Steenrod square}.
\end{definition}

\begin{example}
$Sq^1$ is indecomposable. Also $Sq^2$ is indecomposable since $Sq^1Sq^1=0$.
\end{example}

\begin{example}
$Sq^3=Sq^1Sq^2$ by the Adem relations, so it is decomposable. Also, we showed in Example 11.4.2 that $Sq^2 Sq^4=Sq^6+Sq^5 Sq^1$, so $Sq^6$ is decomposable.
\end{example}

\begin{theorem}
$Sq^i$ is indecomposable if and only if $i$ is a power of 2.
\end{theorem}

{\bf Proof:} Let $i$ be a power of 2 and let $u\in H^1(P^\infty; Z_2)$ be the generator. Using the ring homomorphism $Sq$, $$Sq(u^i)=(Sq (u))^i=(u+u^2)^i=u^i+u^{2i} \mod 2.$$ So $Sq^t(u^i)=0$ unless $t=0$ or $t=i$, and $Sq^i(u^i)=u^{2i}$. The fact that $u^{2i}\neq 0$ shows that $Sq^i$ is indecomposable, since otherwise, $$u^{2i}=Sq^i(u^i)=\sum_{t<i} a_t Sq^t(u^i)=0.$$

For the converse, let $i=a+2^k$, where $0<a<2^k$. Writing $b$ for $2^k$, the Adem relations give $$Sq^a Sq^b=\binom{b-1}{a}Sq^{a+b}+\sum_{c>0}\binom{b-c-1}{a-2c} Sq^{a+b-c} Sq^c.$$ Since $b=2^k$ is a power of 2, the next theorem will show that $\binom{b-1}{a}=1 \mod 2.$ So $Sq^i=Sq^{a+b}$ is decomposable. $\blacksquare$

\begin{theorem}
Let $p$ be a prime and let $a$ and $b$ have expansions $a=\sum_{i=0}^m a_i p^i$, and $b=\sum_{i=0}^m b_i p^i$. Then $$\binom{b}{a}=\prod_{i=0}^m \binom{b_i}{a_i} \mod p.$$
\end{theorem}

{\bf Proof:} In the ring $Z_p[x]$, $(1+x)^p=1+x^p$. (This makes high school algebra students especially happy.) Thus $(1+x)^b=\prod(1+x)^{b_i p^i}=\prod (1+x^{p^i})^{b_i}.$ Then $\binom{b}{a}$ is the coefficient of $x^a$ in this expansion as seen from the first expression, but the last expression shows that this coefficient is $\prod_{i=0}^m \binom{b_i}{a_i}$. $\blacksquare$

I can now state Mosher and Tangora's main result.

\begin{theorem}
If there exists a map $f: S^{2n-1}\rightarrow S^n$ of Hopf invariant 1. Then $n$ is a power of 2.
\end{theorem}

{\bf Proof:} Recall that if $f$ exists, then $Sq^n \sigma=\tau$ in the complex $K=S^n\cup_f e^{2n}$, where $\sigma, \tau$ are the generators of $H^\ast(K; Z_2)$ in dimensions $n$ and $2n$ respectively. If $n$ is not a power of 2, then $Sq^n$ is decomposable, but in $K$, $Sq^i=0$ for $0<i<n$. This is a contradiction, so $n$ must be a power of 2. $\blacksquare$

J. Frank Adams proved an even stronger result. 

\begin{theorem}
If there exists a map $f: S^{2n-1}\rightarrow S^n$ of Hopf invariant 1. Then $n=2, 4,$ or $8$. So the only maps of Hopf invariant one are the three Hopf maps.
\end{theorem}

The proof is the subject of a very dense 85 page paper \cite{JFA1}. You should try it if you are looking for a challenge.

\section{The Steenrod Algebra}

Let's recall the sructures of abstract algebra. A {\it group} has a multiplication or an addition but not necessarily both. A {\it ring} has addition and multiplication tied together by a distributive law. A {\it vector space} has addition and scalar multiplication, and an {\it algebra} has addition, multiplication, and scalar multiplication. The typical first example of an algebra is the set of square $n\times n$ matrices over a field. 

Just when you thought things couldn't get more complicated, I will introduce the {\it Steenrod algebra}. It involves our friends the Steenrod squares. The Steenrod algebra is an example of a {\it Hopf algebra}. And cohomology rings are also graded modules over the Steenrod algebra. 

So what is our interest? As I mentioned earlier, the more algebraic structure that encodes our data, the more information we have for classification. For now, other than use specific Steenrod squares $Sq^i$ as features, I don't know what else to do with the extra structure. But I will leave it to you to think about.

Meanwhile, Frank Adams used the Steenrod algebra to build a new spectral sequence to study homotopy groups of spheres. The Adams spectral sequence provides information about the 2-component of $\pi_m(S^n)$, where $m$ is in the {\it stable range} (i.e. $m\leq 2n-2$) and the 2-component is the set of elements whose orders are powers of 2. The $E_2$ term which starts the sequence involves $Ext_\mathcal{A}(M, N)$ where $M$ and $N$ are modules over the Steenrod algebra $\mathcal{A}$ that we will define in this chapter. While this topic is somewhat more advanced than I want to cover in this book, the idea is that even this much algebraic structure leads to meaningful results in geometry. I believe that the potential for data science goes well beyond persistent homology, and it would be very interesting to see where it could lead.

In the next chapter, I will focus on some more basic approaches to computing homotopy groups of spheres. If you are intrigued by the discussion of the Adams spectral sequence, you can learn more in Adams' original papers \cite{JFA1, JFA2}, the last chapter of Mosher and Tangora \cite{MT} or the  book by McLeary \cite{McCl}.

The rest of this section is taken from \cite{MT}. I will discuss some of the major properties of the Steenrod algebra and state some results mostly without the proofs.

The Steenrod algebra is an example of a {\it graded algebra}. I will start by building up to the definition of a graded algebra. Then I will define the Steenrod algebra as an example.

Recall that we defined a module in Definition 3.4.8. It is basically a generalization of a vector space where the scalars can be an arbitrary ring as opposed to a field. If the ring is non-commutative, we have to distinguish between left and right modules depending on which side the scalars are multiplied. In our case, the scalars will always be a {\it commutative} ring $R$ with unit element. In this case, we don't run into this problem and generally multiply with the scalars on the left. So $r\in R$ and $m\in M$ implies $rm\in M$.

\begin{definition}
Let $M$ be a module over a commutative ring with unit element $R$. Then $M$ is {\it unitary}\index{unitary module} or {\it unital}\index{unital module} if $1m=m$, where $1$ is the unit element of $R$.
\end{definition}

\begin{definition}
$M$ is a {\it graded R-module}\index{graded R-module} if $M=\cup_{i=0}^\infty M_i$ where the $M_i$ are unitary $R$-modules. A {\it homomorphism} $f: M\rightarrow N$ is a sequence $f_i$ of $R$-homomorphisms $f_i: M_i\rightarrow N_i$. (I.e. if $r\in R$ and $m_1, m_2\in M_i$ then $f_i(rm_1)=rf_i(m_1)$ and $f_i(m_1+m_2)=f_i(m_1)+f_i(m_2).$) The tensor product of $M\otimes N$ of two graded $R$-modules is defined by letting $$(M\otimes N)_t=\sum_{i=0}^\infty M_i\otimes N_{t-i}.$$
\end{definition}

\begin{definition}
$A$ is a {\it graded R-algebra}\index{graded R-algebra} if $A$ is a graded $R$-module with a multiplication $\phi: A\otimes A\rightarrow A$, where $\phi$ is a homomorphism of graded $R$-modules and there exists an element $e\in A$ such that for $a\in A$, $\phi(e\otimes a)=\phi(a\otimes e)=a.$ We will write $\phi(m\otimes n)$ as $mn.$ We call $e$ the unit of $A$. A homomorphism of graded $R$-algebras is a graded $R$-module homomorphism which respects multiplication and units. 
\end{definition}

\begin{definition}
A graded $R$-algebra $A$ is {\it associative}\index{associative graded R-algebra} if $\phi(\phi\otimes 1_A)=\phi(1_A\otimes \phi)$. If $a, b, c\in A$, then this means that $(ab)c=a(bc)$, corresponding to our usual definition.
\end{definition}

\begin{definition}
A graded $R$-algebra $A$ is {\it commutative}\index{commutative graded R-algebra} if $\phi(T)=\phi: A\otimes A\rightarrow A$, where $T: A\otimes A\rightarrow A\otimes A$ is defined as $$T(m\otimes n)=(-1)^{deg(m)deg(n)}(n\otimes m),$$ where $m, n\in A$. This means that $mn=(-1)^{deg(m)deg(n)}nm.$
\end{definition}

\begin{definition}
A graded $R$-algebra $A$ is {\it augmented}\index{augmented graded R-algebra} if there is a graded algebra homomorphism $\epsilon: A\rightarrow R$, where $R$ is considered as a graded $R$-algebra such that $R_0=R$, and $R_i=0$ for $i>0$. This means that since $\epsilon$ respects the grading, we have that  $\epsilon: A_0\rightarrow R$. An augmented graded $R$-algebra is {\it connected}\index{connected augmented graded R-algebra} if $\epsilon$ is an isomophism.
\end{definition}

Let $A$ and $B$ be graded $R$-algebras. Then $A\otimes B$ can be thought of as the tensor product of graded $R$-modules. We can give it an algebra structure by defining $$(a_1\otimes b_1)(a_2\otimes b_2)=(-1)^{deg(a_2)deg(b_1)}(a_1 a_2)\otimes (b_1 b_2),$$ for $a_1, a_2\in A$ and $b_1, b_2\in B$.

\begin{definition}
Let $M$ be an $R$-module. The {\it tensor algebra}\index{tensor algebra} $\Gamma(M)$ is defined as follows: Let $M^0=R$, $M^1=M$, $M^2=M\otimes M$, and $M^n=M\otimes \cdots\otimes M$ ($n$ times).Then $\Gamma(M)$ is the graded $R$-algebra defined by $\Gamma(M)_r=M^r$, and the product is given by the isomorphism $M^s\otimes M^t\cong M^{s+t}$.
\end{definition}

The tensor algebra $\Gamma(M)$ is associative but not commutative. 

Now we are finally ready to define the Steenrod algebra. Let $R=Z_2$ and $M$ be the graded $Z_2$-module such that $M_i=Z_2$ generated by the symbol $Sq^i$ for $i\geq 0$. Then the tensor algebra $\Gamma(M)$ is bigraded where for example $Sq^p\otimes Sq^q$ is in $\Gamma(M)_{2, p+q}$ and represents the composition $Sq^p Sq^q$. For each pair of integers $(a, b)$ such that $0<a<2b$ let $$R(a, b)=Sq^a\otimes Sq^b+\sum_c\binom{b-c-1}{a-2c} Sq^{a+b-c}\otimes Sq^c.$$ (By now you should recognize these as the Adem relations.)

\begin{definition}
The {\it Steenrod algebra}\index{Steenrod algebra} $\mathcal{A}$ is the quotient algebra $\Gamma(M)/Q$, where $Q$ is the ideal of $\Gamma(M)$ generated by the $R(a, b)$ defined above along with $1+Sq^0$.
\end{definition}

The elements of $\mathcal{A}$ are polynomials in $Sq^i$ with coefficients in $Z_2$ subject to the Adem relations.

\begin{theorem}
The monomials $Sq^I$ where $I$ runs through the admissible sequences form a basis for $\mathcal{A}$ as a $Z_2$ module.
\end{theorem}

{\bf Idea of proof:} Linear independence follows from Theorem 11.4.10. To show that these elements span $\mathcal{A}$, define the {\it moment} $m(I)$ by the formula $$m(I)=m(\{i_1,\cdots, i_r\})=\sum_{s=1}^r i_s s.$$  If $I$ is not admissible, there is a pair $i_s, i_{s+1}$ with $i_s<2i_{s+1}$. Starting at the right and applying the Adem relations to the first such pair, leads to a sum of monomials with strictly lower moment than $I$. Since the moment function is bounded below by 0, the process terminates and the admissible $Sq^I$ span $\mathcal{A}$. $\blacksquare$

\begin{example}
Find a basis for $\mathcal{A}_7$. We need to find sequences $\{i_1,\cdots, i_r\}$ that sum to 7 where  $i_s\geq 2i_{s+1}$ for every consecutive pair of indices $i_s, i_{s+1}$.Starting at 7 and working downwards gives that $\mathcal{A}_7$ is generated by $Sq^7, Sq^6 Sq^1, Sq^5 Sq^2,$ and  $Sq^4 Sq^2 Sq^1$.
\end{example}

\begin{example}
Find a basis for $\mathcal{A}_9$. Starting at 9 and working downwards gives that $\mathcal{A}_9$ is generated by $Sq^9, Sq^8 Sq^1, Sq^7 Sq^2,$ and $Sq^6 Sq^3$.
\end{example}

Now we will give $\mathcal{A}$ even more structure. We will produce an algebra homomorphism $\psi: \mathcal{A}\rightarrow \mathcal{A}\otimes \mathcal{A}$ called the {\it diagonal map} of $\mathcal{A}$.

Let $M$ be the graded $Z_2$ module generated in each $i\geq 0$ by $Sq^i$. Define $\psi: \Gamma(M)\rightarrow \Gamma(M)\otimes \Gamma(M)$ by the formula $$\psi(Sq^i)=\sum_j Sq^j\otimes Sq^{i-j}$$ and the requirement that $\psi$ be an algebra homomorphism so that $$\psi(Sq^r\otimes Sq^s)=\psi(Sq^r)\otimes \psi(Sq^s)=\sum_a (Sq^a\otimes Sq^{r-a})\otimes \sum_b (Sq^b\otimes Sq^{s-b}).$$

\begin{theorem}
The map $\psi$ induces an algebra homomorphism  $\psi: \mathcal{A}\rightarrow \mathcal{A}\otimes \mathcal{A}$.
\end{theorem}

The proof takes some work but the basic idea is to show that the kernel of the natural projection $$p: \Gamma(M)\rightarrow \Gamma(M)/Q=\mathcal{A}$$ is contained in the kernel of $\psi$ defined above.

The diagonal map turns $\mathcal{A}$ into a {\it Hopf algebra} which I will now define.

Let $A$ be a connected graded $R$-algebra with unit. The existence of the unit is expressed by the fact that both compositions in this diagram are the identity, where $1$ represents the identity map on $A$, and $\eta$ is the {\it coaugmentation}\index{coaugmentation}, which is the inverse of the isomorphism $\epsilon: A_0\rightarrow R$. Also, $\phi$ is the multiplication in $A$. So we have a unit element $1\in R$ with $1a=a1=a$ for $a\in A$.

$$\begin{tikzpicture}
  \matrix (m) [matrix of math nodes,row sep=3em,column sep=4em,minimum width=2em]
  {
& A\otimes R & & \\
A & & A\otimes A & A\\
& R\otimes A & & \\};

\path[-stealth]
(m-2-1) edge node [above] {$\cong$} (m-1-2)
(m-2-1) edge node [below] {$\cong$} (m-3-2)
(m-1-2) edge node [above] {$1\otimes\eta$} (m-2-3)
(m-3-2) edge node [below] {$\eta\otimes 1$} (m-2-3)
(m-2-3) edge node [above] {$\phi$} (m-2-4)
;

\end{tikzpicture}$$

Now let $A$ be a graded $R$-module with a given augmentation, which is an $R$-homomorphism $\epsilon: A\rightarrow R$. We say that $A$ is a {\it co-algebra}\index{co-algebra} with {\it co-unit}\index{co-unit} if we are given an $R$-homomorphism $\psi: A\rightarrow A\otimes A$ such that both compositions are the identity in the dual diagram (i.e. arrows turned backwards) below.

$$\begin{tikzpicture}
  \matrix (m) [matrix of math nodes,row sep=3em,column sep=4em,minimum width=2em]
  {
& A\otimes R & & \\
A & & A\otimes A & A\\
& R\otimes A & & \\};

\path[-stealth]
(m-1-2) edge node [above] {$\cong$} (m-2-1)
(m-3-2) edge node [below] {$\cong$} (m-2-1)
(m-2-3) edge node [above] {$1\otimes\epsilon$}(m-1-2)
(m-2-3) edge node [below] {$\epsilon\otimes 1$} (m-3-2) 
(m-2-4) edge node [above] {$\psi$} (m-2-3)
;

\end{tikzpicture}$$

Then $\psi$ is called the diagonal map or the {\it co-multiplication}.

\begin{definition}
Let $A$ be a connected $R$-algebra with augmentation $\epsilon$. Let $A$ also have a co-algebra structure with co-unit and that the diagonal map $\psi: A\rightarrow A\otimes A$ is a homomorphism of $R$-algebras. Then $A$ is a {\it (connected) Hopf algebra}\index{Hopf algebra}.
\end{definition}

The above discussion shows the following: 

\begin{theorem}
The Steenrod algebra $\mathcal{A}$ is a Hopf algebra.
\end{theorem}

Now let $R$ be a field (think of $Z_2$) and let $A$ be a connected Hopf algebra over $R$. Since $R$ is a field, each $A_i$ is a vector space. Suppose each $A_i$ is finite dimensional over $R$. We define the dual Hopf algebra $A^\ast$ by setting $(A^\ast)_i=(A_i)^\ast$, the vector space dual. (Recall from linear algebra that if $V$ is a vector space over a field $F$ then the {\it dual space}\index{dual space} $V^\ast$ is the vector space of homomorphisms (i.e. linear transformations) from $V$ to $F$.)

The multiplication $\phi$ in $A$ gives rise to the diagonal map $\phi^\ast$ in $A^\ast$ and the diagonal map $\psi$ in $A$ gives rise to the multiplication map $\psi^\ast$ in $A^\ast$. Then $A^\ast $ is itself a Hopf algebra and $\phi^\ast $ is associative (commutative) if and only if $\phi$ is associative (commutative) and the same holds for $\psi$ and $\psi^\ast$. $A$ and $A^\ast$ are isomorphic as $R$-modules but not in general as algebras. 

Now let $\mathcal{A}^\ast$ be the dual of the Steenrod algebra. Then $\mathcal{A}^\ast$ is a Hopf algebra with an associative diagonal map $\phi^\ast$ and with associative and commutative multiplication $\psi^\ast.$ It turns out that $\mathcal{A}^\ast$ is a polynomial ring.

\begin{definition}
Let $u$ be the generator of $H^1(P^\infty; Z_2)$. For each $i\geq 0$, let $\xi_i$ be the element of $\mathcal{A}_{2^i-1}^\ast$ such that $$\xi_i(\theta)(u)^{2^i}=\theta(u)\in H^{2^i}(P^\infty; Z_2)$$ for all $\theta\in \mathcal{A}_{2^i-1}$. We set $\xi_0$ to be the unit in $\mathcal{A}^\ast$.
\end{definition}

To break this down, $\xi_i(\theta)$ is an element of $Z_2$. If this element is $z$, we want $z(u^{2^i})=\theta(u)$. So $\xi_i(\theta)$ is the coefficient of $u^{2^i}$ when $\theta$ is applied to $u$.

Here is the main structure theorem for $\mathcal{A}^\ast$. See \cite{MT} for the proof.

\begin{theorem}
As an algebra, $\mathcal{A}^\ast$ is the polynomial ring over $Z_2$ generated by the $\{\xi_i\}$ for $i\geq 1.$
\end{theorem}

Now the Steenrod algebra actually acts on the cohomology ring over a space. This makes the cohomology ring an algebra over an algebra.

To put this more precisely, let $A$ be a graded $R$-algebra and $M$ be a graded $R$-module, where $R$ is a commutative ring with unit element. 

\begin{definition}
$M$ is a {\it graded A-module} if there is an $R$-module homomorphism $\mu: A\otimes M\rightarrow M$ such that $\mu(1\otimes m)=m$ and $$\mu(\phi_A\otimes 1)=\mu(1\otimes \mu): A\otimes A\otimes M\rightarrow M,$$ where $\phi_A$ is the multiplication in $A$. 
\end{definition}

Now suppose that $A$ is a Hopf algebra and that $M$ is an $A$-module which is also an $R$-algebra. Then $M\otimes M$ is an $A$-module under the composition $$\mu': A\otimes M\otimes M \xrightarrow{\psi\otimes 1\otimes 1} A\otimes A\otimes M\otimes M\xrightarrow{T} A\otimes M\otimes A\otimes M\xrightarrow{\mu\otimes\mu} M\otimes M,$$ where $\psi$ is the diagonal map of $A$, and $T$ interchanges the second and third terms of the tensor product. 

\begin{definition}
$M$ is an algebra over the Hopf algebra $A$ if the multiplication $\phi_M: M\otimes M\rightarrow M$ is a homomorphism of modules. This means that the following diagram commutes:
$$\begin{tikzpicture}
  \matrix (m) [matrix of math nodes,row sep=3em,column sep=4em,minimum width=2em]
  {
A\otimes M\otimes M & M\otimes M\\
A\otimes M &  M\\};

\path[-stealth]
(m-1-1) edge node [above] {$\mu'$} (m-1-2)
(m-1-1) edge node [left] {$1\otimes \phi_M$} (m-2-1)
(m-1-2) edge node [right] {$\phi_M$} (m-2-2)
(m-2-1) edge node [below] {$\mu$} (m-2-2)

;

\end{tikzpicture}$$
\end{definition}

In the case of interest to us, let $R=Z_2, A=\mathcal{A},$ and $M=H^\ast(X; Z_2)$, where $X$ is a topological space. Recall that the diagonal map $\psi$ of $\mathcal{A}$ is defined by  $$\psi(Sq^i)=\sum_j Sq^j\otimes Sq^{i-j}$$ and we have $$\psi(Sq^r\otimes Sq^s)=\psi(Sq^r)\otimes \psi(Sq^s)=\sum_a (Sq^a\otimes Sq^{r-a})\otimes \sum_b (Sq^b\otimes Sq^{s-b}).$$ If $\theta\in\mathcal{A}$ and $x, y\in H^\ast(X; Z_2)$, then the Cartan formula implies that $$\theta(\phi_M(x\otimes y))=\theta(x\cup y)=\phi_M(\psi(\theta)(x\otimes y)).$$ So the cohomology ring with $Z_2$ coefficients is actually an algebra over the Hopf algebra $\mathcal{A}$.

Mosher and Tangora also derive a formula for the  diagonal map of the dual $\mathcal{A}^\ast$ which I will state for completeness.

The diagonal map $\phi^\ast$ of $\mathcal{A}^\ast$ is given by the formula $$\phi^\ast(\xi_k)=\sum_{i=0}^k(\xi_{k-i})^{2^i}\otimes\xi_i.$$

An interesting question is what the extra structure of Steenrod squares could tell you about your data. Could the extra structure help with hard classification problems now that more of the geometry is included? The good news is that in low dimensional cases, there are algorithms that can compute them in a reasonable time period. These issues will be the subject of the last section of this chapter.

\section{Cohomology of Eilenberg-Mac Lane Spaces}

Recall way back in Section 11.2 that we said that the cohomology operations of the form $H^n(X; \pi)\rightarrow H^m(X; G)$ are in one to one correspondence with the cohomology group $H^m(K(\pi, n); G)$. (We abbreviate the latter as $H^m(\pi, n; G)$.) To recall how the correspondence worked, we let $\theta$ be such an operation. By the Universal Coefficient Theorem, $$H^n(X; \pi_n(X))\cong Hom(H_n(X); \pi_n(X))$$ if $X$ is $(n-1)$-connected. $K(\pi, n)$ definitely meets this criterion, so we let $\imath_n$ be the element of $H^n(\pi, n; \pi)$ corresponding to the inverse $h^{-1}$ of the Hurewicz isomorphism $h: \pi_n(X)\rightarrow H_n(X)$ where $X=K(\pi, n)$. Then $\imath_n \in H^n(\pi, n, \pi)$ is called the {\it fundamental class} and the cohomology operation $\theta$ is paired with $\theta(\imath_n)\in H^m(\pi, n; G)$.

Now look at the cohomology operation $Sq^i$. It is of the form $H^n(X; Z_2)\rightarrow H^{n+i}(X; Z_2)$. Here $G=\pi=Z_2$, and $m=n+i$. So $Sq^i$ is paired with an element of $H^{n+i}(Z_2, n, Z_2)$. So we would expect a close connection between the structure of the cohomology ring $H^\ast(Z_2, n, Z_2)$ and the Steenrod squares.. We will see that this is a polynomial ring whose generators are certain Steenrod squares. Proving this fact is not easy, though. It will involve some sophisticated machinery, especially spectral sequences.

The calculations are not as much hard as messy. To understand them, though, I will outline some key points. Then you can look at Mosher and Tangora if you are interested in the details.

At this point, you should review Section 9.9. I will assume you know the basic definitions of spectral sequences, the use of the Leray-Serre spectral sequence in computing homology of fiber spaces, and how to get a spectral sequence from an exact couple. 

\subsection{Bockstein Exact Couple}

The Bockstein exact couple will be needed for Section 11.7.8 to describe the structure of $H^\ast (Z_{2^m}, 1; Z_2)$. 

\begin{definition}
The {\it Bockstein exact couple}\index{Bockstein exact couple} is of the form:
$$\begin{tikzpicture}
  \matrix (m) [matrix of math nodes,row sep=3em,column sep=4em,minimum width=2em]
  {
D^1=H^\ast( ; Z) & & H^\ast( ; Z)\\
& E^1=H^\ast( ; Z_2)&\\};

\path[-stealth]
(m-1-1) edge node [above] {$i^1$} (m-1-3)
(m-2-2) edge node [below] {$k^1$} (m-1-1)
(m-1-3) edge node [below] {$j^1$} (m-2-2)

;

\end{tikzpicture}$$

Here $i^1$ is induced by multiplication by 2 in $Z$; $j^1$ is reduced by the reduction mod 2 homomorphism $\rho: Z\rightarrow Z_2$, and $k^1$ is the Bockstein homomorphism $\beta$ defined at the start of Section 11.4.1. The differential $d=j^1 k^1$ is the Bockstein homomorphism $\delta_2$ (Section 11.4.1). 
\end{definition}

Mosher and Tangora use subscripts rather than superscripts for successive differentials. So $d_1=d, d_2=j^2 k^2: E^2\rightarrow E^2$, etc.

The operation $d_r$ acts as follows. Take a cocycle in $Z_2$ coefficients, represent it by an integral cocycle, take its coboundary, divide by $2^r$, and reduce the coefficients mod 2. (The division is possible since $d_r$ is defined only on the kernel of $d_{r-1}$. Every $d_r$ raises the dimension by 1 in $H^\ast( ; Z_2).$

The Bockstein differentials act on $Z_2$ cohomology, but they can give information on integral cohomology as well. The Universal coefficeint theorems can be used to prove the following.

\begin{theorem}
Elements of $H^\ast(X; Z_2)$ which come from free integral classes lie in $\ker d_r$ for every $r$ and not in im $d_r$ for any $r$. If $z$ generates a summand $Z_{2^r}$ in $H^{n+1}(X; Z)$ then there exist corresponding summands $Z_2$ in $H^n(X; Z_2)$ and  $H^{n+1}(X; Z_2)$. If we call their generators $z'$ and $z''$ respectively, then $d_i(z')$ and $d_i(z'')$ are zero for $i<r$ and $d_r(z')=z''$. This means that $z'$ and $z''$ are not in im $d_i$ for $i<r$. We say that the image of $\rho$ on the free subgroup of $H^\ast(X; Z)$ "persists to $E^\infty$", and that $z'$ and $z''$ "persist to $E^r$ but not to $E^{r+1}$."
\end{theorem}

Our last result will be useful in calculating homotopy groups of spheres.

\begin{theorem}
Suppose that $H^i(X; Z_2)=0$ for $i<n$, and $H^n(X; Z_2)=Z_2$ with generator $z$. Then the part of $H^n(X; Z)$ not involving odd primes is $Z$ if $d_r(z)=0$ for all $r$ and $Z_{2^n}$ if $d_i(z)=0$ for $i<n$, and $d_n(z)\neq 0$.
\end{theorem}

\subsection{Serre's Exact Sequence for a Fiber Space}

Recall from Theorem 9.6.22 that a fiber space has a long exact sequence in homotopy. Is there something similar in homology? Not always, but we can show the following using the Leray-Serre spectral sequence (Section 9.9.2).

\begin{theorem}
Let $F\xrightarrow{i} E\xrightarrow{p} B$ be a fiber space with $B$ simply connected. Suppose that $H_i(B)=0$ for $0<i<p$ and $H_j(F)=0$ for $0<j<q$. Then there is an exact sequence $$H_{p+q-1}(F)\xrightarrow{i_\ast} H_{p+q-1}(E)\xrightarrow{p_\ast} H_{p+q-1}(B)\xrightarrow{\tau}H_{p+q-2}(F)\rightarrow\cdots\rightarrow H_1(E)\rightarrow 0.$$
\end{theorem}

{\bf Proof:} Since $E_{i, j}^2=H_p(B; H_q(F))$ in the Leray-Serre spectral sequence, $E_{i, j}^2=0$, when either $0<i<p$ or $0<j<q$. Looking at the $E_{i, j}^\infty$ terms where $i+j=n$ gives the exact sequence $$0\rightarrow E_{0, n}^\infty\rightarrow H_n(E)\rightarrow E_{n, 0}^\infty\rightarrow 0.$$ In general, we have the exact sequence $$0\rightarrow E_{n, 0}^\infty\rightarrow E_{n, 0}^n\xrightarrow{d^n} E_{0, n-1}^n\rightarrow E_{0, n-1}^\infty\rightarrow 0.$$ If $n<p+q$, $E_{n, 0}^n\cong H_n(B)$ and $E_{0, n-1}^n\cong H_{n-1}(F).$ Substituting into the previous sequence in splicing it into the one above for all $n<p+q$ proves the theorem. $\blacksquare$

The sequence of the theorem is called {\it Serre's exact sequence for a fiber space.}\index{Serre's exact sequence for a fiber space}

\subsection{Transgression}

The map $\tau: H_n(B)\rightarrow H_{n-1}(F)$ of Theorem 11.7.3 corresponds to $d^n_{n, 0}$ and is called the {\it transgression}\index{transgression}. This map was only defined if $n<p+q$. and $F$ and $B$ are $(q-1)$-connected and $(p-1)$-connected respectively. 

To define it more generally as $d^n_{n, 0}$, $\tau$ has a subgroup of $H_n(B)$ as its domain and a quotient of $H_{n-1}(F)$ as its range.

Here is a different but equivalent description which will be more useful. In the fiber space $F\rightarrow E\xrightarrow{p} B$, let $p_0:(E, F)\rightarrow (B, \ast)$ be the map of pairs induced by $p$. (Here $\ast$ denotes a single point of $B$.) For any $n$, let $\overline{\tau}$ be a map from im $(p_0)_\ast$ which is a subgroup of $H_n(B, \ast)$ to a quotient group of $H_{n-1}(F)$ defined as follows: If $x\in \mbox{im }(p_0)_\ast\subset H_n(B, \ast)$, choose $y\in H_n(E, F)$ such that $(p_0)_\ast(y)=x$, and take $\partial y$ as a representative of $\overline{\tau}(x)$. Then $\overline{\tau}(x)$ is in a quotient group because of the indeterminacy in choosing $y$.

\begin{definition}
We say that $x\in H_n(B)$ is {\it transgressive}\index{transgressive} if $\tau(x)$ is defined. Since $H_n(B)=E_{n, 0}^2$, this is the same as saying that $d^i(x)=0$ for all $i<n$.
\end{definition}

The above is equivalent to the condition that $x\in \mbox{im }(p_0)_\ast$ and that if $x=(p_0)_\ast(y)$ then $\tau{x}$ is the homology class of $\partial y$. So $\tau=\overline{\tau}$. The proof is in Serre's thesis \cite{Ser}.

\subsection{Cohomology Version of Leray-Serre Spectral Sequence}

As we want to interact with Steenrod squares, we really would like a spectral sequence for cohomology of a fiber space as opposed to homology. Fortunately, Serre supplied us with one of those as well. We will rely on it heavily for computing cohomology of some Eilenberg-Mac Lane spaces.

We will assume $F\rightarrow E\xrightarrow{p} B$ is a fiber space and make life easier by letting $B$ be simply connected. We will get a spectral sequence $(E_r, d_r)$ of which $E_2^{p, q}=H^p(B; H^q(F))$. Note that we distinguish the cohomology spectral sequence from the homology spectral sequence by turning subscripts into superscripts and vice versa. We have $$d_r^{p, q}: E_r^{p, q}\rightarrow E_r^{p+r, q-r+1},$$ so that the bidegree of $d_r$ is $(r, -r+1)$ which is the opposite direction of the homology differential $d^r$. The sequence converges to $H^\ast(E)$.

\begin{theorem}
If $R$ is a commutative ring with unit element, then there is a spectral sequence with $E_2^{p, q}=H^p(B; H^q(F; R))$ converging to $H^\ast(E; R)$ such that:\begin{enumerate}
\item For each $r$, $E_r$ is a bigraded ring. If $\mu$ is the ring multiplication, then $$\mu: E_r^{p, q}\otimes  E_r^{p', q'}\rightarrow  E_r^{p+p', q+q'}.$$
\item In $E_r$, $d_r$ is an {\it anti-derivation}, which means that $$d_r(ab)=d_r(a)b+(-1)^kad_r(b),$$ where $k$ is the total degree of $a$. The formula holds wherever it makes sense.
\item The product in the ring $E_{r+1}$ is induced by the product in $E_r$, and the product in $E_\infty$ is the cup product in $H^\ast (E; R)$. 
\end{enumerate}
\end{theorem}

Note that if $R$ is a field, the K\"{u}nneth Theorem implies that $E_2=H^\ast(B; R)\otimes H^\ast(F; R)$.

\begin{theorem}
{\bf: Serre's Cohomology Exact Sequence:} If $B$ is $(p-1)$-connected $F$ is $(q-1)$-connected, then there is an exact sequence analogous to the homology version which terminates as: $$\cdots\rightarrow H^{p+q-2}(F)\xrightarrow{\tau} H^{p+q-1}(B)\xrightarrow{p^\ast} H^{p+q-1}(E)\xrightarrow{i^\ast} H^{p+q-1}(F).$$
\end{theorem}

The cohomology transgression is $$\tau=d_n^{0, n-1}: E_n^{0, n-1}\rightarrow E_n^{n, 0}.$$ 

\begin{definition}
$x\in H^{n-1}(F)$ is {\it transgressive} if $\tau(x)$ is defined. This is true if $d_i(x)=0$ for all $i<n$, or equivalently, $\delta x$ lies in im $p^\ast\subset H^n(E, F).$ If $\delta x=p^\ast(y)$, then $\tau(x)$ contains $y$. 
\end{definition}

\begin{theorem}
If $x\in H^{n-1}(F)$  is transgressive, then so is $Sq^i(x)$ and if $y\in \tau(x)$, then $Sq^i(y)\in \tau(Sq^i(x)).$
\end{theorem}

{\bf Proof:} If $y\in \tau(x)$, then $p^\ast y=\delta x$, so $Sq^i(p^\ast y)=Sq^i(\delta x)$. By naturality, $p^\ast(Sq^i(y))=\delta(Sq^i(x))$ which means that $Sq^i(y)=\tau(Sq^i(x))$. $\blacksquare$

\subsection{$H^\ast(Z, 2; Z)$}

For the remainder of this section, I will state some results about cohomology of Eilenberg-Mac Lane spaces. I will limit myself to stating results with some short comments. See Mosher and Tangora for detailed calculations, all of which involve cohomology spectral sequences.

Two types of fiber spaces will be useful here. The first is the path-space fibration of Definition 9.5.2. We will write it as $\Omega B\rightarrow E\xrightarrow{p} B$. (I am following the notation of Mosher and Tangora that $\Omega B$ is the loops on $B$ which is more common than the notation $\Lambda B$ found in Hu.) Recall that $E$ is the space of paths on $B$ and is contractible. If $B=K(\pi, n)$ then the homotopy exact sequence shows that the fiber is $F=\Omega B=K(\pi, n-1)$.

The other type comes from a short exact sequence of abelian groups $$0\rightarrow A\rightarrow B\rightarrow C\rightarrow 0.$$ Mosher and Tangora show that we can use this sequence to construct a fiber space $$K(A, n)\xrightarrow{i} K(B, n)\xrightarrow{p} K(C, n).$$

The work will involve setting up the appropriate fiber space and calculating the cohomology using the Serre spectral sequence.

\begin{theorem}
$H^\ast(Z, 2; Z)$ is the polynomial ring $Z[\imath_2]$ where $\imath_2$ is of degree 2 and $(\imath_2)^n$ generates $H^{2n}(Z, 2; Z)$. 
\end{theorem}

The theorem is proved using the spectral sequence of the fiber space $S^1=K(Z, 1)\rightarrow E\rightarrow K(Z, 2)$ where $E$ is contractible and $R=Z$. 

\subsection{$H^\ast(Z_2, n; Z_2)$}

The first approach in Mosher and Tangora is to calculate the cohomology of $H^\ast(Z_2, 2; Z_2)$ in low dimensions using the fiber space $$F=K(Z_2, 1)=P^\infty\rightarrow E\rightarrow B=K(Z_2, 2).$$ We already know that the cohomology of $F$ is a polynomial ring on one generator and the results about transgression are used to relate it to the generators of the cohomology of $B$. To find the ring structure in general, we need another theorem.

\begin{definition}
A graded ring $R$ over $Z_2$ has the ordered set $x_1, x_2, \cdots$ as a {\it simple system of generators} if the monomials $\{x_{i_1}x_{i_2}\cdots x_{i_r}|i_1<i_2<\cdots<i_r\}$ form a $Z_2$-basis for $R$ and if for each $n$, only finitely many $x_i$ have degree $n$.
\end{definition}

\begin{theorem}
{\bf Borel's Theorem:} Let $F\rightarrow E\xrightarrow{p} B$ be a fiber space with $E$ acyclic, and suppose that $H^\ast(F; Z_2)$ has a simple system $\{x_i\}$ of transgressive generators. Then $H^\ast(B; Z_2)$ is a polynomial ring generated by $\{\tau(x_i)\}$. 
\end{theorem}

\begin{theorem}
$H^\ast(Z_2, n; Z_2)$ is a polynomial ring over $Z_2$ with generators $\{Sq^I(\imath_n)\}$ where $I$ runs through all admissible sequences with excess less than $n$. 
\end{theorem}

We know that the theorem is true for $n=1$ so we proceed by induction on $n$. We suppose it is true for $n$ and use the fiber space  $$F=K(Z_2, n)\rightarrow E\rightarrow B=K(Z_2, n+1)$$ with $E$ acyclic. We apply Borel's Theorem and go through a long calculation with admissible sequences. See \cite{MT} for the details. 

\subsection{Additional Examples}

Here are four more results from \cite{MT} on cohomology of Eilenberg-Mac Lane spaces. 

\begin{theorem}
$H^\ast(Z, 2; Z_2)$ is the polynomial ring over $Z_2$ generated by $\imath_2\in H^2(Z, 2; Z_2)$.
\end{theorem}

\begin{theorem}
$H^\ast(Z, n; Z_2)$ is the polynomial ring over $Z_2$ generated by the elements $\{Sq^I(\imath_n)\}$ where $I$ runs through all admissible sequences of excess $e(I)<n$ , and where the last entry $i_r$ of $I$ is not equal to one.
\end{theorem}

For the next theorem. we need to define an {\it exterior algebra}.

\begin{definition}
Let $V$ be a vector space over a field $F$. (More generally, we could start with a module over a commutative ring with unit element, but we will be using the field $Z_2$ in our case.) The {\it exterior algebra}\index{exterior algebra} $E(V)$ is the quotient $\Gamma(V)/J$ where $\Gamma(V)$ is the tensor algebra over $V$ (Definition 11.6.7), and $J$ is generated by the elements of the form $v\otimes v$ for $v\in V$. The exterior algebra $E(v_1,\cdots,v_r)$ for $v_1,\cdots,v_r$ is the exterior algebra over the subspace of $V$ generated by $v_1,\cdots,v_r$.
\end{definition}

If $f: \Gamma(V)\rightarrow \Gamma(V)/J=E(V)$ is the projection, then we write $v\wedge w$ for $f(v\otimes w)$. So we have that $v\wedge v=0$ and $$(v+ w)\wedge (v+w)=(v\wedge v)+(v\wedge w)+(w\wedge v)+ (w\wedge w)=(v\wedge w)+(w\wedge v)=0,$$ so $(v\wedge w)=-(w\wedge v)$.

\begin{theorem}
$H^\ast(Z_{2^m}, 1; Z_2)=P[d_m(\imath_1)]\otimes E(\imath_1)$ for $m\geq 2$. Here $d_m$ denotes the $m$th Bockstein homomorphism (See Sections 11.7.1-11.7.4.), $P[d_m(\imath_1)]$ is the polynomial ring over $Z_2$ generated by $d_m(\imath_1)\in H^2(Z_{2^m}, 1; Z_2)$ and  $E(\imath_1)$ is the exterior algebra generated by $\imath_1\in H^1(Z_{2^m}, 1; Z_2).$ (Recall that the Bockstein homomorphisms raise the dimension by one. This happens in both the homology and the cohomology versions.)
\end{theorem}

\begin{theorem}
$H^\ast(Z_{2^m}, n; Z_2)$ is the polynomial ring over $Z_2$ generated by elements $\{Sq^{I_m}(\imath_n)\}$ where we define $Sq^{I_m}=Sq^I$ if $I$ terminates with $i_r>1$, and we replace $Sq^{i_r}$ with $d_m$ in $Sq^I$ if $i_r=1$. In this case as well, $I$ runs through all admissible sequences of excess $e(I)<n$.
\end{theorem}

One further question. We know that we could not have computed these cohomology rings with software like Ripser, since Eilenberg-Mac Lane spaces are infinite dimensional and therefore have infinitely many simplices. Given that Kenzo can compute spectral sequences in some cases, could it have found these results? That could be the subject of some interesting experiments.

\section{Reduced Powers}

It is well known that if Superman is exposed to kryptonite, he becomes less powerful and thus has {\it reduced powers}. That, however is not the subject of this section.

Steenrod reduced powers are an additional family of cohomology operations. Suppose that 2 is not your favorite prime number, and you really prefer to work over $Z_{41}$. Reduced powers handle the case of cohomology over $Z_p$ where $p$ is an odd prime number. 

The construction is discussed in Steenrod's lectures \cite{Ste1} but presented in much more detail in Steenrod and Epstein \cite{SE}. Since Mosher and Tangora leave out reduced powers entirely, the material for this section comes from \cite{SE}. I will refer you there for the full constructions, but I will mention that the construction in Section 11.3 is modified to look at the action of a subgroup $\pi$ of the permutation group $S_n$ on a $W\otimes K\otimes\cdots\otimes K=W\otimes K^{\otimes n}$, where $W$ is now a $\pi$-free acyclic complex, and $K^{\otimes n}$ is the tensor product of $K$ with itself $n$ times. 

The result is the {\it (cyclic) reduced powers}\index{reduced powers}\index{Steenrod reduced powers} taking the form $$\mathcal{P}^i: H^q(K; Z_p)\rightarrow H^{q+2i(p-1)}(K; Z_p).$$

The noation leaves ambiguous which value of $p$ is being used. I used to imagine notation like $\mathfrak{3}^i$ for $p=3$, but that could obvously get confusing. In \cite{Ste1}, Steenrod sometimes uses the notation $\mathcal{P}_p^i$, but I have never seen that anywhere else. So for now, we will use the conventional notation, $\mathcal{P}^i$ and understand which value of $p$ we are using by context.

Our first task is to write down the properties corresponding to those listed in Theorem 11.4.1. We start with a Bockstein homomorphism analogous to the one representing $Sq^1$.

Let $p$ be an odd prime and let $\beta: H^q(X; Z_p)\rightarrow H^{q+1}(X; Z_p)$ be the Bockstein operator associated with the exact coefficient sequence $$0\rightarrow Z\xrightarrow{m} Z\rightarrow Z_p\rightarrow 0,$$ where $m$ is multiplication by $p$. To define it, we start with an homomorphism $\overline{\beta}: H^q(X; Z_p)\rightarrow H^{q+1}(X; Z)$ defined in an analogous way to the mod 2 case. Let $x\in H^q(X; Z_p)$. Represent $x$ by a cocycle $c$ and choose an integral cochain $c'$ which maps to $c$ under reduction mod $p$. Then $\delta c'$ is a multiple of p, and $\frac{1}{p}(\delta c')$ represents $\overline{\beta} x$.  The composition of $\overline{\beta}$ and the reduction homomorphism  $H^{q+1}(X; Z)\rightarrow H^{q+1}(X, Z_p)$ gives a homomorphism $$\beta:  H^q(X; Z_p)\rightarrow H^{q+1}(X; Z_p).$$

We have that $\beta^2=0$ and writing $xy$ for $x\cup y$, $$\beta(xy)=\beta(x)y+(-1)^{dim(X)}x\beta(y).$$

We now give the main axioms for $\mathcal{P}^i$. Compare them to the axioms for $Sq^i$ in Theorem 11.4.1.

\begin{theorem}
The Steenrod reduced powers $\mathcal{P}^i$ for $i\geq 0$ have the following properties:
\begin{enumerate}
\item  $\mathcal{P}^i$ is a natural homomorphism $H^q(K, L; Z_p)\rightarrow H^{q+2i(p-1)}(K, L; Z_p).$
\item If $2k>q$, then $\mathcal{P}^k(x)=0$ for $x\in H^q(K; L; Z_p)$.
\item If $x\in H^{2k}(K; L; Z_p)$, then $\mathcal{P}^k(x)=x^p=x\cup\cdots\cup x$ ($p$ times).
\item $\mathcal{P}^0$ is the identity homomorphism.
\item $\delta^\ast\mathcal{P}^i=\mathcal{P}^i \delta^\ast$, where $\delta^\ast$ is the coboundary homomorphism $\delta^\ast: H^q(L; Z_p)\rightarrow H^{q+1}(K, L; Z_p)$.
\item {\bf Cartan Formula:} Writing $xy$ for $x\cup y$ we have $$\mathcal{P}^k(xy)=\sum_i(\mathcal{P}^ix)(\mathcal{P}^{k-i}y).$$
\item {\bf Adem Relations:} For $a<pb$, $$\mathcal{P}^a \mathcal{P}^b=\sum_{t=0}^{\lfloor \frac{a}{p}\rfloor}(-1)^{a+t}\dbinom{(p-1)(b-t)-1}{a-pt}\mathcal{P}^{a+b-t}\mathcal{P}^t,$$ where the binomial coefficient is taken mod p. 

If $a\leq b$ then \begin{align*}
\mathcal{P}^a \beta\mathcal{P}^b&=\sum_{t=0}^{\lfloor \frac{a}{p}\rfloor}(-1)^{a+t}\dbinom{(p-1)(b-t)}{a-pt}\beta\mathcal{P}^{a+b-t}\mathcal{P}^t\\
&+\sum_{t=0}^{\lfloor \frac{a-1}{p}\rfloor}(-1)^{a+t-1}\dbinom{(p-1)(b-t)-1}{a-pt-1}\mathcal{P}^{a+b-t}\beta\mathcal{P}^t.
\end{align*}
\end{enumerate}
\end{theorem}

It turns out that like $Sq^i$, $\mathcal{P}^i$ also commutes with suspension. The Bockstein homomorphism $\beta$ also commutes with both coboundary and suspension.

Also, like $Sq^i$, the reduced powers form an algebra.

\begin{definition}
Define the Steenrod algebra $\mathcal{A}(p)$ to be the graded associative algebra generated by elements $\mathcal{P}^i$ of degree $2i(p-1)$ and $\beta$ of degree 1, subject to $\beta^2=0$, the Adem Relations, and $\mathcal{P}^0=1$.
\end{definition}

A monomial in $\mathcal{A}(p)$ can be written in the form $$\beta^{\epsilon_0}\mathcal{P}^{s_1}\beta^{\epsilon_1}\cdots\mathcal{P}^{s_k}\beta^{\epsilon_k},$$ where $\epsilon_i\in\{0, 1\}$ and $s_i$ is a positive integer. We denote this monomial by $\mathcal{P}^I$, where $$I=\{\epsilon_0, s_1, \epsilon_1, \cdots, s_k, \epsilon_k, 0, 0, \cdots\}.$$ A sequence $I$ is {\it admissible} if $s_i>ps_{i+1}+\epsilon_i$ for each $i\geq 1$. The corresponding $\mathcal{P}^I$ as well as $\mathcal{P}^0$ are called {\it admissible monomials.} Define the moment of $I$ to be $\sum_{i=1}^k i(s_i+\epsilon_i)$. The {\it degree} of $I$ is the degree of $\mathcal{P}^I$. 

Steenrod and Epstein prove the following.

\begin{theorem}
The admissible monomials form a basis for $\mathcal{A}(p)$.
\end{theorem}

Steenrod and Epstein have a similar discussion to that of Steenrod squares of the structure of $\mathcal{A}(p)$ and its dual $\mathcal{A}(p)^\ast$. For example, $\mathcal{A}(p)$ is a Hopf algebra, so reduced powers make the cohomology ring $H^\ast(X; Z_p)$ an algebra over a Hopf algebra. See there for more details. 

Obviously, Steenrod squares and reduced powers do not represent all possible cohomology operations. For example, Steenrod \cite{Ste1} mentions the Pontrjagin $p$th power $$\mathfrak{P}: H^{2q}(K; Z_{p^k})\rightarrow H^{2pq}(K; Z_{p^{k+1}}).$$ Pontrjagin \cite{Pont2} discovered these operations for $p=2$ and Thomas \cite{Thom1, Thom2} discovered them for primes $p>2$.

There are also {\it higher order} cohomology operations derived from the Steenrod squares that have proved useful in the computation of homotopy groups of spheres. I will not discuss them in this book but there is a good description in Chapter 16 of Mosher and Tangora \cite{MT}. 

Finally, there is some work out there on computer calculation of Steenrod reduced powers. See for example the papers of Gonzalez-Diaz and Real \cite{GDR2, GDR3} on this subject. Still, the idea of using cohomology operations for data science is so new, I will stick with the simplest ones, the Steenrod squares, for the remainder of this chapter.

\section{Vector Bundles and Stiefel-Whitney Classes}

In this section, I will describe a classical example of Steenrod squares being used in a classification problem. Algebraic topology has often been focused on solving problems in differential geometry. As an example, recall the problem of the existence of a nonzero tangent vector field for $S^n$ that we addressed in Theorem 4.1.35. This section deals with {\it vector bundles} which lie on the boundary between these two fields.

Roughly speaking, vector bundles are fiber spaces whose fiber is a vector space. An example from differential geometry is the {\it tangent bundle} on a manifold. The base space is the manifold and the fiber at a point on the manifold is the tangent space at that point. Vector bundles can be classified by certain cohomology classes of the base space known as {\it characteristic classes.} I will focus on a specific type called {\it Stiefel-Whitney classes}. These are classes in the cohomology over $Z_2$, and it turns out they can be represented using Steenrod squares. Two vector bundles that are homeomorphic have the same Stiefel-Whitney classes, but the converse is not always true. 

Question: If our data lies on a manifold, could we solve a classification problem by calculating the Stiefel-Whitney classes of some appropriate bundle?

I won't address computation issues, but Aubrey HB touches on this area in her thesis \cite{HB}. See there for details.

In what follows, I will give the definition of vector bundles and Stiefel-Whitney classes along with stating some of their properties. I will describe a notable example, the Grassmann manifold, and describe its cohomology. Finally, I will describe how the Stiefel-Whitney classes are represented by Steenrod squares. 

The classic book on characteristic classes if the book by Milnor and Stasheff \cite{MS}. All of the material in this section comes from there.

\subsection{Vector Bundles}

I will start by defining a vector bundle. Compare this to the definition of a fiber space from Section 9.3.

\begin{definition}
A real {\it vector bundle}\index{vector bundle} $\xi$ over a topological space $B$ (called the base space) consists of the following:
\begin{enumerate}
\item A topological space $E=E(\xi)$ called the {\it total space}.
\item A map $\pi: E\rightarrow B$ called the {\it projection map}.
\item For each $b\in B$, $\pi^{-1}(b)$ is a vector space over the reals called the {\it fiber} over $b$ and denoted by $F_b$.
 \end{enumerate}
In addition, the bundle must satisfy the condition of {\it local triviality}\index{local triviality}: For each $b\in B$, there exists a neighborhood $U\subset B$, an integer $n\geq 0$, and a homeomorphism $$h: U\times R^n\rightarrow \pi^{-1}(U)$$ so that for each $b\in U$, the correspondence $x\rightarrow h(b, x)$ defines a vector space isomorphism between $R^n$ and $\pi^{-1}(b).$
\end{definition}

The pair $(U, h)$ is called a {\it local coordinate system for} $\xi$ about $b$. If $U$ can be chosen to be all of $B$, then $\xi$ is called a {\it trivial bundle}.

The bundle $\xi$ is called an {\it n-plane bundle} if the dimension of $F_b$ is $n$ for all $b\in B$.

\begin{definition}
Two bundles $\xi$ and $\eta$ over the same base space $B$ are isomorphic, written $\xi\cong\eta$  if there exists a homeomorphism $f: E(\xi)\rightarrow E(\eta)$ which maps each vector space $F_b(\xi)$ isomorphically onto the corresponding vector space $F_b(\eta)$.
\end{definition}

\begin{example}
The trivial bundle with total space $B\times R^n$ and projection $\pi(b, x)=b$ is a vector bundle denoted $\epsilon^n_B$. 
\end{example}

\begin{example}
The {\it tangent bundle}\index{tangent bundle} $\tau_M$ of a smooth manifold $M$ (the coordinate charts are infinitely differentiable) has a total space $DM$ which consists of all pairs $(x, v)$ with $x\in M$ and $v$ in the tangent space $T_xM$. (See Definition 6.5.2.) The projection map is $\pi: DM\rightarrow M$ with $\pi(x, v)=x$. If $\tau_M$ is trivial, then $M$ is called {\it parallelizable}.\index{parallelizable manifold}
\end{example}

\begin{example}
The {\it normal bundle}\index{normal bundle} $\nu$ of a smooth manifold $M\in R^n$ has a total space $E\subset M\times R^n$ which consists of all pairs $(x, v)$ such that $v$ is orthogonal to the tangent space $T_xM$.
\end{example}

\begin{definition}
The {\it canonical line bundle}\index{canonical line bundle} $\gamma_n^1$ over the projective space $P^n$ is the bundle with total space $E(\gamma_n^1)\subset P^n\times R^{n+1}$ consisting of all pairs $(\{\pm x\}, v)$ such that $v$ is a multiple of $x$. Each fiber $\pi^{-1}(\{\pm x\})$ can be identified with the line through $x$ and $-x$ in $R^{n+1}$. 
\end{definition}

\begin{theorem}
The bundle $\gamma_n^1$ over $P^n$ is not trivial for $n\geq 1$.
\end{theorem}

\begin{definition}
A {\it cross-section} of a vector bundle $\xi$ with base space $B$ is a continuous function $s: B\rightarrow E(\xi)$ which takes $b\in B$ into the corresponding fiber $F_b(\xi)$. The cross section is nowhere zero if $s(b)$ is a nonzero vector in $F_b(\xi)$ for each $b$.
\end{definition}

A vector field is a cross-section of the tangent bundle of a manifold. 

A trivial $R^1$ bundle has a nowhere zero cross-section. The bundle $\gamma_n^1$ turns out to have no such cross section. Milnor and Stasheff prove a more general result.

\begin{theorem}
An $R^n$ bundle $\xi$ is trivial if and only if $\xi$ admits $n$ cross-sections $s_1, \cdots, s_n$ which are nowhere dependent, i.e, the vectors $s_1(b), \cdots, s_n(b)$ are linearly independent for all $b\in B.$
\end{theorem}

\subsection{New Bundles from Old Ones}

Now that we know what a vector bundle is, here are some constructions of new ones from old ones. They will be necessary to understand the axioms for Stiefel-Whitney classes.

\begin{definition}
Let $\xi$ be a vector bundle with projection $\pi: E\rightarrow B$, and let $B'\subset B$. Then letting $E'=\pi^{-1}(B')$ and $\pi': E'\rightarrow B'$ be the restriction of $\pi'$ to $E'$, we get a vector bundle $\xi|B'$ called the {\it restriction} of $\xi$ to $B'$. Each fiber $F_b(\xi|B')$ is equal to the corresponding fiber $F_b(\xi)$.
\end{definition}

The next example is entirely analogous to induced fibering construction at the end of Section 9.3.

\begin{definition}
Let $\xi$ be a vector bundle with projection $\pi: E\rightarrow B$, and let $B_1$ be a topological space. Given a map $f: B_1\rightarrow B$, we get the {\it induced bundle} $f^\ast \xi$ over $B_1$. The total space $E_1$ of $f^\ast \xi$ is the subset of $B_1\times E$ consisting of pairs $(b, e)$ with $f(b)=\pi(e)$. 
\end{definition}

\begin{definition}
Let $\xi_1$ and $\xi_2$ be vector bundles with projections $\pi_i: E_i\rightarrow B_i$, for $i=1, 2$. Then the {\it Cartesian product} $\xi_1\times\xi_2$ is the bundle with projection map $$\pi_1\times \pi_2: E_1\times E_2\rightarrow B_1\times B_2$$ with fibers $$(\pi_1\times\pi_2)^{-1}(b_1, b_2)=F_{b_1}(\xi_1)\times F_{b_2}(\xi_2).$$
\end{definition}

\begin{definition}
Let $\xi_1, \xi_2$ be bundles over the same base space $B$. Let $d: B\rightarrow B\times B$ be the diagonal map $d(b)=(b, b)$. The induced bundle $d^\ast(\xi_1\times \xi_2)$ over $B$ is called the {\it Whitney sum}\index{Whitney sum} of $\xi_1$ and $\xi_2$ and is denoted $\xi_1\oplus\xi_2$. Each fiber $F_b(\xi_1\oplus\xi_2)$ is isomorphic to the direct sum $F_b(\xi_1)\oplus F_b(\xi_2).$
\end{definition}

\begin{definition}
Let $\xi, \eta$ be bundles over the same base space $B$, and suppose that $E(\xi)\subset E(\eta)$. The bundle $\xi$ is a {\it sub-bundle} of $\eta$ if each fiber $F_b(\xi)$ is a vector subspace of the corresponding fiber $F_b(\eta)$.
\end{definition}

\begin{theorem}
Let $\xi_1, \xi_2$ be sub-bundles of $\eta$ such that each fiber $F_b(\eta)$ is the vector space direct sum $F_b(\xi_1)\oplus F_b(\xi_2).$ Then $\eta$ is isomorphic to the Whitney sum $\xi_1\oplus\xi_2$.
\end{theorem}

Recall that a {\it Hilbert space} is a vector space $V$ with an inner product. Assuming the scalars are the reals, we have an inner product $v.w\in R$ for $v, w\in V$. Milnor and Stasheff call this a {\it Euclidean vector space}. There is a function $\mu: V\rightarrow R$ with $\mu(v)=v.v$. In general $$v.w=\frac{1}{2}(\mu(u+v)-\mu(v)-\mu(w)).$$

The function $\mu$ is {\it quadratic}, which means that $\mu(v)=\sum_i \mathfrak{l}_i(v)\mathfrak{l}'_i(v)$ where $\mathfrak{l}$ and $\mathfrak{l}'$ are linear, and $\mu$ is {\it positive definite}, which means that $\mu(v)>0$ if $v\neq 0$.

\begin{definition}
A {\it Euclidean vector bundle} is a real vector bundle $\xi$ together with a continuous function $\mu: V\rightarrow R$ such that the restriction of $\mu$ to each fiber of $\xi$ is quadratic and positive definite. $\mu$ is called a {\it Euclidean metric} on $\xi$.
\end{definition}

Now let $\eta$ be a Euclidean vector bundle and $\xi$ a sub-bundle of $\eta$. Let $F_b(\xi^\perp)$ consist of all vectors $v\in F_b(\eta)$ such that $v.w=0$ for all $w\in F_b(\xi)$. Let $E(\xi^\perp)$ be the union of the $F_b(\xi^\perp).$ Then $E(\xi^\perp)$ is the total space of a sub-bundle $\xi^\perp$ of $\eta$, and $\eta$ is isomorphic to the Whitney sum $\xi\oplus\xi^\perp.$ We call $\xi^\perp$ the {\it orthogonal complement} of $\xi$.

As an example, the orthogonal complement of the tangent bundle of a manifold is the normal bundle.

\subsection{Stiefel-Whitney Classes}

We can now define Stiefel-Whitney classes. They will be elements of the cohomology ring over $Z_2$ of the base space of a vector bundle. Milnor and Stasheff also discuss other types of characteristic classes such as the Euler classes for integral cohomology and Chern classes for vector bundles whose fibers are vector spaces over the complex numbers. See \cite{MS} for more details.

As was the case for Steenrod squares, Stiefel-Whitney classes can be defined axiomatically without worrying about whether they actually exist. This is the approach of Milnor and Stasheff. After the axioms and some properties are described, \cite{MS} describes the construction in terms of Steenrod squares. It is then shown that classes constructed in this way satisfy the required axioms. I will follow their approach including a short discussion of an interesting example known as a Grassmann manifold.

\begin{definition}
{\bf Axiomatic definition:} To each real vector bundle $\xi$ there is a sequence of cohomology classes $$w_i(\xi)\in H^i(B(\xi); Z_2), i=0, 1, 2, \cdots,$$ called the {\it Stiefel-Whitney classes}\index{Stiefel-Whitney classes} of $\xi$. They satisfy the following four axioms:\begin{enumerate}
\item The class $w_0(\xi)$ is equal to the unit element $1\in H^0(B(\xi); Z_2)$, and $w_i(\xi)=0$ for $i>n$ if $\xi$ is an $n$-plane bundle.
\item {\bf Naturality:} If $f: B(\xi)\rightarrow B(\eta)$ is covered by a bundle map from $\xi$ to $\eta$ (i.e. there is a map $E(\xi)\rightarrow E(\eta)$ that preserves fibers.) Then $$w_i(\xi)=f^\ast w_i(\eta).$$
\item {\bf The Whitney Product Theorem:} If $\xi$ and $\eta$ are vector bundles over the same base space, then $$w_k(\xi\oplus \eta)=\sum_{i=0}^k w_i(\xi)\cup w_{k-i}(\eta).$$
\item For the line bundle $\gamma_1^1$ over the circle $P^1$, $w_1(\gamma_1^1)\neq 0.$
\end{enumerate}
\end{definition}

Assuming such classes exist, here are some consequences.

\begin{theorem}
If $\xi$ is isomorphic to $\eta$, then $w_i(\xi)=w_i(\eta)$ for all $i$. 
\end{theorem}

This shows that Stiefel-Whitney classes classify vector bundles just as homology, cohomology, and homotopy classify topological spaces. A similarity is that isomorphic vector bundles must have the same Stiefel-Whitney classes, but having the same classes does not guarantee that two bundles are isomorphic. 

\begin{theorem}
If $\epsilon$ is a trivial vector bundle then $w_i(\epsilon)=0$ for $i>0$.
\end{theorem}

This is true since for a trivial vector bundle $\epsilon$, there is a bundle map from $\epsilon$ to a bundle over a point. This theorem combined with the Whitney product theorem gives the following.

\begin{theorem}
If $\epsilon$ is trivial, and $\eta$ is another vector bundle over $B$, then $w_i(\epsilon\oplus\eta)=w_i(\eta).$
\end{theorem}

\begin{theorem}
If $\xi$ is a $n$-plane bundle with a Euclidean metric which possesses a nowhere zero cross-section, then $w_n(\xi)=0$. If $\xi$ has $k$ cross-sections that are nowhere linearly dependent then $$w_{n-k+1}(\xi)=w_{n-k+2}(\xi)=\cdots=w_n(\xi)=0.$$
\end{theorem}

\subsection{Grassmann Manifolds}

Grassmann manifolds are an interesting class of manifolds whose cohomology ring can be described in terms of Stiefel-Whitney classes. I will describe what they are and state some of their properties.

\begin{definition} 
A {\it Grassmann manifold}\index{Grassmann manifold} $G_n(R^{n+k})$ is the set of all $n$-dimensional planes through the origin of $R^{n+k}$.
\end{definition}

\begin{definition} 
An {\it n-frame} in $R^{n+k}$ is an $n$-tuple of linearly independent vectors in $R^{n+k}$ The collection of all $n$-frames in $R^{n+k}$ is an open subset of the $n$-fold Cartesian product $R^{n+k}\times\cdots\times R^{n+k}$ called a {\it Stiefel manifold}\index{Stiefel manifold} and denoted $V_n(R^{n+k})$.
\end{definition}

There is a function $$q: V_n(R^{n+k})\rightarrow G_n(R^{n+k})$$ which maps an $n$-frame to the $n$-plane it spans. Milnor and Stasheff define the topology on  $G_n(R^{n+k})$ as the quotient topology determined by this map and use it to prove the following.

\begin{theorem}
The Grassmann manifold  $G_n(R^{n+k})$ is a compact manifold of dimension $nk$. The correspondence $X\rightarrow X^\perp$ which assigns each $n$-plane to its orthogonal $k$-plane defines a homeomorphism between $G_n(R^{n+k})$ and $G_k(R^{n+k}).$
\end{theorem}

The proof of this statement is lengthy, but I will comment on the dimension $nk$. The dimension of a manifold is $d$ if every open neighborhood of a point is homeomorphic to $R^d$. Now if $X_0\in G_n(R^{n+k})$ then $X_0$ is an $n$-plane in $R^{n+k}$, so let $X_0^\perp$ be its orthogonal $k$-plane, so that $R^{n+k}=X_0\oplus X_0^\perp$. Then let $U$ be the open subset of $G_n(R^{n+k})$ consisting of all $n$-planes $Y$ such that the projection $p: X_0\oplus X_0^\perp\rightarrow X_0$ maps $Y$ onto $X_0$. In other words, $Y\cap X_0^\perp=0.$ So each $Y\in U$ can be considered to be the graph of a linear transformation $T(Y): X_0\rightarrow X_0^\perp$. So there is a one to one correspondence $T$ between $U$ and $Hom(X_0, X_0^\perp)$, and the latter has dimension $nk$. It takes some work, but we are done once we prove that $T$ is a homeomorphism.(See \cite{MS}.)

A canonical vector bundle $\gamma_n(R^{n+k})$ can be constructed over $G_n(R^{n+k})$ by taking as the total space $E(\gamma_n(R^{n+k}))$ the set of pairs $$(n\mbox{-plane in }R^{n+k}, \mbox{ vector in that plane}).$$

We can generalize this construction to $R^\infty$ and call the bundle $\gamma^n$.

Let $G_n=G_n(R^\infty)$ be the set of $n$-planes in $R^\infty$. Milnor and Stasheff prove the following.

\begin{theorem}
The cohomology ring $H^\ast(G_n; Z_2)$ is a polynomial ring over $Z_2$ freely generated (ie. there are no poynomial relations) by $w_1(\gamma_n),\cdots, w_n(\gamma_n)$.
\end{theorem}

\subsection{Representation of Stiefel-Whitney Classes as Steenrod Squares}

We are finally ready to see the connection between Stiefel-Whitney classes and Steenrod squares. 

Let $\xi$ be an $n$-plane vector bundle with total space $E$, base space $B$, and projection map $\pi$. Let $E_0$ be the nonzero elements of $E$. (Remember that an element of $E$ is an element of some fiber $F_b=\pi^{-1}(b)$. So any element of $E$ is an element of some vector space and so it makes sense to talk about nonzero elements. If $F$ is a typical fiber, let $F_0=F\cap E_0$ denote the nonzero elements of $F$.

The following two theorems are  proved in a later chapter of \cite{MS} with more advanced tools that we have not covered. See there (Chapter 10) if you are curious. 

In this section, we will always work with $Z_2$ coefficients, so $H^i(X)$ will be understood to be $H^i(X; Z_2)$.

\begin{theorem}
$H^i(F, F_0)=0$ for $i\neq n$ and $H^n(F, F_0)=Z_2.$

$H^i(E, E_0)=0$ for $i<n$ and $H^i(E, E_0)=H^{i-n}(B)$ for $i\geq n$.
\end{theorem}

\begin{theorem}
$H^i(E, E_0)=0$ for $i<n$, and $H^n(E, E_0)$ contains a unique class $u$ called the {\it fundamental cohomology class} such that for each fiber $F$, the restriction $u|(F, F_0)\in H^n(F, F_0)$ is the unique nonzero class in $H^n(F, F_0)$. In addition the correspondence $x\rightarrow x\cup u$ defines an isomorphism $H^k(E)\rightarrow H^{k+n}(E, E_0)$ for every $k$.
\end{theorem}

Now a vector bundle has a zero cross-section which embeds $B$ as a deformation retract of $E$, so that the projection $\pi: E\rightarrow B$ induces an isomorphism $\pi^\ast: H^k(B)\rightarrow H^k(E)$. 

\begin{definition}
The {\it Thom isomorphism}\index{Thom isomorphism} $\phi: H^k(B)\rightarrow H^{k+n}(E, E_0)$ is defined to be the composition of the two isomorphisms $$H^k(B)\xrightarrow{\pi^\ast} H^k(E)\xrightarrow{\cup u} H^{k+n}(E, E_0).$$
\end{definition}

Using $\phi$ we define the Stiefel-Whitney class $w_i(\xi)\in H^i(B)$ by $$w_i(\xi)=\phi^{-1} Sq^i \phi(1).$$ Now $\phi(1)=1\cup u=u$, so $Sq^i(\phi(1))=Sq^i(u)$, and $w_i(\xi)$ is the unique cohomology class in $H^i(B)$ such that $\phi(w_i(\xi))=\pi^\ast w_i(\xi)\cup u$ is equal to $Sq^i(u)$.

It remains to check that this definition satisfies the four axioms of Stiefel-Whitney classes.

The main takeaway from this is that Steenrod squares have appeared in a classical classification problem from geometry. Although there is still the weakness that two non-isomorphic bundles could have the same Stiefel-Whitney classes, they may still be useful given the right type of data. 

For any practical problem involving Steenrod squares in data science, we will need to know how to compute them. That is the topic of the next section.

\section{ Computer Computation of Steenrod Squares}

Author’s Note: Anibal Medina Mardones has made me aware of some work he has done on computational Steenrod squares that I was not aware of at the time I wrote this book. In particular, the paper \cite{Lupo} written with Umberto Lupo and Guillaume Tauzin describes an efficient method for computing barcodes representing the persistence of Steenrod squares with an outer loop that runs in quadratic time. In addition, they have developed a Python package “steenroder” that runs this algorithm. It is available at https://github.com/Steenroder/steenroder. This work has the advantage of directly fitting into a TDA workflow and is faster than the methods I mention below,

In this section, I will describe in detail the formula for computing Steenrod squares derived by Gonz\'{a}lez-D\'{\i}az and Real. All of the material comes from \cite{GDR1}, except for a simplified formula for $Sq^1$ from Aubrey HB's thesis \cite{HB}.

The computations will be performed on our old friends, the simplicial sets. We also use a different construction of Steenrod squares than the one presented in Section 11.2. We will use the construction of Steenrod and Epstein \cite{SE} in which they construct a series of morphisms $\{D_i\}$ called {\it higher diagonal approximations} which are proved to exist by the acyclic model theorem. Then letting $C_\ast^N(X)$ denote the {\it normalized simplicial set} consisting of only non-degenerate simplices of $X$, if $c\in  Hom(C_j^N(X), Z_2),$ and $x\in C_{i+j}^N(X)$ then $$Sq^i(c)(x)=\mu(\langle c\otimes c, D_{j-i}(x)\rangle), \hspace{.5in} i\leq j,$$ and $Sq^i(c)(x)=0$ if $i>j$. Here $\mu$ is the homomorphism induced by multiplication in $Z_2$.

So we are left with the problem of writing down a formula for $D_i$. Now Real derived such a formula in \cite{Real2} in terms of the component morphisms of the {\it Elilenberg-Zilber contraction} from $C_\ast^N(X\times X)$ onto $C_\ast^N(X)\otimes C_\ast^N(X)).$ Now the formula for $D_i$ involves an explicit formula for the {\it homotopy operator}. (All of these terms will be defined soon.) But the homotopy operator is defined in terms of partitions or {\it shuffles} of the face and degeneracy operators of the simplicial sets. This would give complexity of order $2^n$ to evaluate an element of dimension $n$. The main contribution of \cite{GDR1} is to drastically reduce this complexity by realizing that the face and degeneracy operators of $X$ can always be written as $$s_{j_t}\cdots s_{j_1}\partial_{i_1}\cdots\partial{i_s},$$ where $j_t>\cdots>j_1\geq 0$, and $i_s>\cdots>i_1\geq 0$. Since the image of the morphisms of the Eilenbeg-Zilber contraction is always normalized, we can throw away any terms with a degeneracy operator in its expression. This produces a much simpler formula for the $D_i$.

I won't repeat the definition of a simplicial set, but we will use $\partial_i$ for face operators and $s_i$ for degeneracy operators. A simplex $x$ is {\it degenerate} if $x=s_iy$ for some simplex $y$ and degeneracy operator $s_i$. A simplex which is not degenerate is {\it non-degenerate}. The proofs in the paper make heavy use of the relation $\partial_i s_j=s_{j-1}\partial_i$ for $i<j$. This puts a composition of face and degeneracy operators in the standard form and determines which terms are degenerate.

For an $n$ simplex, we use $d_n=\sum_{i=0}^n (-1)^i\partial_i$ as the homology boundary and $\delta^n$ as the corresponding coboundary. Letting $X$ be a simplicial set and $C_\ast(X)$ be the chain complex $\{C_n(X), d_n\}$ define $s(C_\ast(X))$ as the graded $R$-module generated by all of the degenerate simplices. In $C_\ast(X)$, $d_n(s(C_{n-1}(X)))\subset s(C_{n-2}(X))$, so $C_\ast^N(X)=\{C_n(X)/s(C_{n-1}(X)), d_n\}$ is a chain complex called the {\it normalized chain complex associated to X}. If $G$ is a ring, let $C^\ast(X: G)$ be the cochain complex associated to $C_\ast^N(X)$.

At the end of Chapter 9, I spoke about {\it reductions}. Following Eilenberg and Mac Lane \cite{EM1}, I will use the term {\it contraction} and repeat the definition with the notation used in \cite{GDR1}.

\begin{definition}
Let $N$ and $M$ be chain complexes. A {\it contraction}\index{contraction} from $N$ onto $M$ is a triple $(f, g, \phi)$, where $f: N\rightarrow M$ is called {\it projection}, $g: M\rightarrow N$ is called {\it inclusion}, and $\phi: N\rightarrow N$ is called {\it homotopy} and raises the degree by 1. (Note that the analogy is with {\it chain homotopy}.) These maps satisfy the following relations:\begin{enumerate}
\item $fg=1_M.$
\item $\phi d+d\phi+gf=1_N$.
\item $\phi g=0$.
\item $f\phi=0.$
\item $\phi\phi=0.$
\end{enumerate}
\end{definition}

This definition implies an equivalence between big complex $N$ and small complex $M$. 

\begin{definition}
Let $(f_1, g_1, \phi_1)$ and $(f_2, g_2, \phi_2)$ be two contractions. We can construct two additional contractions:\begin{enumerate}
\item The {\it tensor product contraction} $(f_1\otimes f_2, g_1\otimes g_2, \phi_1\otimes g_2 f_2+1_{N_1}\otimes\phi_2)$ from $N_1\otimes N_2$ to $M_1\otimes M_2$.
\item If $N_2=M_1$, the {\it composition contraction} $(f_2 f_1, g_1 g_2, \phi_1+g_1\phi_2 f_1)$ from $N_1$ to $M_2$.
\end{enumerate}
\end{definition}

The formula for Steenrod squares comes from an explicit formula for a particular contraction from $C_\ast^N(X\times Y)$ to $C_\ast^N(X)\otimes C_\ast^N(Y)$  called the {\it Eilenberg-Zilber contraction} \cite{EZ}. Before I define it, I will have to teach you a new dance, the {\it (p, q)-shuffle.}

\begin{definition}
If $p$ and $q$ are non-negative integers, a {\it (p, q)-shuffle}\index{(p, q)-shuffle} $(\alpha, \beta)$ is a partition of the set $$\{0, 1, \cdots, p+q-1\}$$ of integers into disjoint subsets $\alpha_1<\cdots<\alpha_p$ and $\beta_1<\cdots<\beta_q$ of $p$ and $q$ integers respectively. The {\it signature} of the shuffle $(\alpha, \beta)$ is defined by $$sig(\alpha, \beta)=\sum_{i=1}^p \alpha_i-(i-1).$$
\end{definition} 

\begin{example}
Let $p=4, q=6, \alpha=\{3, 5, 6, 8\}$, and $\beta=\{0, 1, 2, 4, 7, 9\}$. Then $$sig(\alpha, \beta)=(3-0)+(5-1)+(6-2)+(8-3)=3+4+4+5=16.$$
\end{example}

\begin{definition}
The {\it Eilenberg-Zilber contraction}\index{Eilenberg-Zilber contraction} from $C_\ast^N(X\times Y)$ to $C_\ast^N(X)\otimes C_\ast^N(Y)$, where $X$ and $Y$ are simplicial sets, consists of the triple (AW, EML, SHI) defined by:\begin{enumerate}
\item The {\it Alexander-Whitney operator}\index{AW} $$AW: C_\ast^N(X\times Y)\rightarrow C_\ast^N(X)\otimes C_\ast^N(Y)$$ defined as $$AW(a_m\times b_m)=\sum_{i=0}^m \partial_{i+1} \cdots\partial_m a_m\otimes \partial_0\cdots\partial_{i-1} b_m.$$ If $X=Y$, AW is a simplicial approximation to the diagonal and this operator allows for the construction of cup products. Interchanging $a_m$ and $b_m$ produces a different approximation and comparing them leads to the Steenrod squares.
\item The {\it Eilenberg-Mac Lane operator}\index{EML} $$EML: C_\ast^N(X)\otimes C_\ast^N(Y)\rightarrow C_\ast^N(X\times Y)$$ defined by $$EML(a_p\otimes b_q)=\sum_{(\alpha, \beta)\in \{(p, q)\mbox{-shuffles}\}} (-1)^{sig(\alpha, \beta)}s_{\beta_q}\cdots s_{\beta_1} a_p\times s_{\alpha_p}\cdots s_{\alpha_1}  b_q.$$ This operator is a process of "triangulation" in $X\times Y$.
\item The {\it Shih operator}\index{SHI} $$SHI: C_\ast^N(X\times Y)\rightarrow C_{\ast+1}^N(X\times Y)$$ is defined by $$SHI(a_0\times b_0)=0;$$ \begin{align*}
SHI(a_m\times b_m)=\sum_{\substack{0\leq q\leq m-1\\0\leq p\leq m-q-1\\(\alpha, \beta)\in \{(p+1, q)\mbox{-shuffles}\}}}(-1)^{\overline{m}+sig(\alpha, \beta)+1} &s_{\beta_q+\overline{m}}\cdots s_{\beta_1+\overline{m}}s_{\overline{m}-1}\partial_{m-q+1}\cdots\partial_m a_m\\
&\times s_{\alpha_{p+1}+\overline{m}}\cdots s_{\alpha_1+\overline{m}}\partial_{\overline{m}}\cdots\partial_{m-q-1} b_m;\end{align*} where $\overline{m}=m-p-q$ and $sig(\alpha, \beta)=\sum_{i=1}^{p+1} \alpha_i-(i-1).$
\end{enumerate}
\end{definition}

Eilenberg and Mac Lane \cite{EM2} derived a recursive formula for SHI, but the one here was stated by Rubio in \cite{Rub} and proved by Morace in the appendix to \cite{Real1}.

We will need to define some more maps.

The diagonal map $$\Delta: C_\ast^N(X)\rightarrow  C_\ast^N(X\times X)$$ is defined by $\Delta(x)=(x, x).$

The automorphism $$t:  C_\ast^N(X\times X)\rightarrow  C_\ast^N(X\times X)$$ is defined by $t(x_1, x_2)=(x_2, x_1).$

The automorphism $$T:  C_\ast^N(X)\otimes  C_\ast^N(X)\rightarrow  C_\ast^N(X)\otimes  C_\ast^N(X)$$ is defined by $T(x_1\otimes x_2)=(-1)^{dim(x_1)dim(x_2)}x_2\otimes x_1$.

Now the AW operator is not commutative. Being an approximation to the diagonal, though, it can be used to compute cup products. Let $R$ be a ring, $c\in C^i(X; R), c'\in C^j(X; R)$, and $x\in C_{i+j}^N(X)$, the cup product of $c$ and $c'$ is \begin{align*}
(c\cup c')(x) &=\mu(\langle c\otimes c', AW\Delta(x)\rangle)\\
&=\mu(\langle c, \partial_{i+1}\cdots\partial_{i+j}x\rangle\otimes\langle c', \partial_0\cdots\partial_{i-1}x\rangle),
\end{align*} where $\mu$ is induced by multiplication in $R$.

In \cite{Ste4}, Steenrod builds a sequence of maps $\{D_i\}$ called {\it higher diagonal approximations} where $$D_i: C_\ast^N(X)\rightarrow  C_\ast^N(X)\otimes C_\ast^N(X)$$ defined by \begin{align*}
D_0 &=AW \Delta\\
d_\otimes D_{i+1}+(-1)^i D_{i+1} d &=T D_i +(-1)^{i+1}D_i,
\end{align*} where $d$ and $d_\otimes$ are the homology differentials (boundary maps) of $C_\ast^N(X)$ and $C_\ast^N(X)\otimes C_\ast^N(X)$ respectively. It turns out that the $D_i$ can be expressed in the form $D_i=h_i \Delta$ where $$h_i: C_\ast^N(X\times X)\rightarrow C_\ast^N(X)\otimes C_\ast^N(X)$$ is a homomorphism of degree $i$  and whose existence can be proved using the acyclic model theorem. This gives a recursive formula for the $\{D_i\}$.

Instead of following this approach, Gonz\'{a}lez-D\'{\i}az and Real make use of the explixit formula for the Eilenberg-Zilber contraction, It is shown in \cite{Real2} that $h_i=AW(t(SHI)^i)$ for all $i$. The formula for the Steenrod square $$Sq^i: H^j(X; Z_2)\rightarrow H^{j+i}(X; Z_2)$$ is now $$Sq^i(c)(x)=\mu(\langle c\otimes c, AW(t(SHI)^{j-i}(x, x))\rangle)$$ for $i\leq j$, and $Sq^i(c)(x)=0$ for $i>j$.

Now the image of $h_i$ lies in $C_\ast^N(X)\otimes C_\ast^N(X)$, so expressing the compositions of face and degeneracy operators of the summands of the formulas in the standard form $$s_{j_t}\cdots s_{j_1}\partial_{i_1}\cdots\partial{i_s},$$ determines which terms can be eliminated. We only keep terms with no degeneracy operators. In this way, Gonz\'{a}lez-D\'{\i}az and Real reduce the amount of work required to compute Steenrod squares and make these computations tractable in a lot of cases.

Now I will state some of the explicit formulas from \cite{GDR1}. As scary as they may look at first, think of how easy they would be to implement on a computer. The proofs involve some combinatorics and a lot of messy arithmetic, but follow directly from the formulas we have given. See \cite{GDR1} for the details.

\begin{theorem}
Let $R$ be a ring and $X$ a simplicial set. Let $(AW, EML, SHI)$ be the Eilenberg-Zilber contraction from $C_\ast^N(X\times X)$ onto $C_\ast^N(X)\otimes C_\ast^N(X)$. Then the morphism $$h_n=AW((t(SHI))^n): C_m^N(X\times X)\rightarrow (C_\ast^N(X)\otimes C_\ast^N(X))_{m+n}$$ can be expressed in the form:\begin{itemize}
\item \begin{align*}
AW((t(SHI))^n)=&\sum_{i_n=n}^m \hspace{.1 in}\sum_{i_{n-1}=n-1}^{i_n-1}\cdots\sum_{i_0=0}^{i_1-1}(-1)^{A(n)+B(n, m, \overline{\imath})+C(n, \overline{\imath})+D(n, m, \overline{\imath})}\\
&\partial_{i_0+1}\cdots\partial_{i_1-1}\partial_{i_2+1}\cdots\partial_{i_{n-1}-1}\partial_{i_n+1}\cdots\partial_{m}\\
&\otimes\partial_{0}\cdots\partial_{i_0-1}\partial_{i_1+1}\cdots\partial_{i_{n-2}-1}\partial_{i_{n-1}+1}\cdots\partial_{i_n-1}
\end{align*} if $n$ is even or 
\item \begin{align*}
AW((t(SHI))^n)=&\sum_{i_n=n}^m \hspace{.1 in}\sum_{i_{n-1}=n-1}^{i_n-1}\cdots\sum_{i_0=0}^{i_1-1}(-1)^{A(n)+B(n, m, \overline{\imath})+C(n, \overline{\imath})+D(n, m, \overline{\imath})}\\
&\partial_{i_0+1}\cdots\partial_{i_1-1}\partial_{i_2+1}\cdots\partial_{i_{n-2}-1}\partial_{i_{n-1}+1}\cdots\partial_{i_n-1}\\
&\otimes\partial_{0}\cdots\partial_{i_0-1}\partial_{i_1+1}\cdots\partial_{i_{n-1}-1}\partial_{i_n+1}\cdots\partial_{m}
\end{align*} if $n$ is odd, \end{itemize}

where 
$$A(n)=
\left
\{
\begin{array}{ll}
1 & \mbox{if }n\equiv3, 4, 5, 6\mod 8\\
0 & \mbox{otherwise.}
\end{array}
\right.
$$

$$B(n, m, \overline{\imath})=
\left
\{
\begin{array}{ll}
\sum_{j=0}^{\lfloor\frac{n}{2}\rfloor} i_{2j} & \mbox{if }n\equiv 1, 2 \mod 4\\
\\
\sum_{j=0}^{\lfloor\frac{n-1}{2}\rfloor} i_{2j+1}+nm & \mbox{if }n\equiv 0, 3 \mod 4
\end{array}
\right.
$$

$$C(n, \overline{\imath})=\sum_{j=0}^{\lfloor\frac{n}{2}\rfloor}( i_{2j}+i_{2j-1})(i_{2j-1}+\cdots+i_0)$$

and

$$D(n, m, \overline{\imath})=
\left
\{
\begin{array}{ll}
(m+i_n)(i_n+\cdots+i_0) & \mbox{if }n\mbox{ is even},\\
0 & \mbox{if }n\mbox{ is odd},
\end{array}
\right.
$$
where $\overline{\imath}=(i_0, i_1, \cdots, i_n).$
\end{theorem}

The theorem gives an explicit formula for the cup-$i$ product since for $c\in C^p(X; R)$ and $c'\in C^q(X: R)$, we have $$(c\cup_i c')(x)=\mu(\langle c\otimes c', D_i(x)\rangle),$$ for all $i$ and $x\in C_{p+q-i}^N(X).$

In a similar way, we get a formula for Steenrod squares. Since we said that $$Sq^i(c)(x)=\mu(\langle c\otimes c, AW(t(SHI)^{j-i}(x, x))\rangle),$$ we get the following.

\begin{theorem}
Let $X$ be a simplicial set. Then if $c\in C^j(X; Z_2)$ and $x\in C_{i+j}^N(X)$ then $$Sq^i: H^j(X; Z_2)\rightarrow H^{j+i}(X: Z_2)$$ is defined by: 
\begin{itemize}
\item If $i\leq j$ and $i+j$ is even, then: \begin{align*}
Sq^i(c)(x)=&\sum_{i_n=S(n)}^m \hspace{.1 in}\sum_{i_{n-1}=S(n-1)}^{i_n-1}\cdots\sum_{i_1=S(1)}^{i_2-1}\\
&\mu(\langle c, \partial_{i_0+1}\cdots\partial_{i_1-1}\partial_{i_2+1}\cdots\partial_{i_{n-1}-1}\partial_{i_n+1}\cdots\partial_{m}x\rangle\\
&\otimes\langle c, \partial_{0}\cdots\partial_{i_0-1}\partial_{i_1+1}\cdots\partial_{i_{n-2}-1}\partial_{i_{n-1}+1}\cdots\partial_{i_n-1}x\rangle).
\end{align*}
\item If $i\leq j$ and $i+j$ is odd, then: \begin{align*}
Sq^i(c)(x)=&\sum_{i_n=S(n)}^m \hspace{.1 in}\sum_{i_{n-1}=S(n-1)}^{i_n-1}\cdots\sum_{i_1=S(1)}^{i_2-1}\\
&\mu(\langle c, \partial_{i_0+1}\cdots\partial_{i_1-1}\partial_{i_2+1}\cdots\partial_{i_{n-2}-1}\partial_{i_{n-1}+1}\cdots\partial_{i_n-1}x\rangle\\
&\otimes\langle c, \partial_{0}\cdots\partial_{i_0-1}\partial_{i_1+1}\cdots\partial_{i_{n-1}-1}\partial_{i_n+1}\cdots\partial_{m}x\rangle).
\end{align*}
\item If $i>j$, then $Sq^i(c)(x)=0.$
\end{itemize}
Here $n=j-i$, $m=i+j$, $$S(k)=i_{k+1}-i_{k+2}+\cdots+(-1)^{k+n-1}i_n+(-1)^{k+n}\lfloor\frac{m+1}{2}\rfloor+\lfloor\frac{k}{2}\rfloor,$$ for all $0\leq k\leq n$ and $i_0=S(0)$.
\end{theorem}

Finally, \cite{GDR1} evaluates the complexity. Here is the main result.

\begin{theorem}
Let $X$ be a simplicial set and $k$ be a non-negative integer. If $c\in C^{i+k}(X; Z_2)$, then the number of face operators taking part in the formula for $Sq^i(c)$ is $O(i^{k+1}).$
\end{theorem}

This looks like bad news. The complexity is exponential. Is it so bad that you won't live long enough to see the computation finish? Let's take a closer look.We start with the example in \cite{GDR1}.

\begin{example}
Suppose there are $O(k^2)$ non-degenerate simplices in each $X_k$. The number of face operators in $Sq^i(c)$ for $c\in C^{i+2}(X; Z_2)$ is $O(i^5)$. This is because the number of face operators in the formula is $O(i^3)$ and we multiply this by the number of simplices in $X_{2i+2}$ as we need to evaluate $Sq^i(c)$ on each of them. 
\end{example}

Here are some reasons why this is not so bad. \begin{enumerate}
\item In TDA, we normally work in low dimensions. If $c$ is of dimension 3 for example, the worst complexity is $3^3$ or 27 face operators per simplex.
\item Remember what a face operator is. We have an ordered list of vertices and we remove one corresponding to a particular index. That should be easy and fast to implement.
\item The lower the dimension you start with, the fewer nonzero Steenrod squares. 
\item If you want to just use $Sq^1$ as a feature, you can implement HB's simplified formula I will state below.
\end{enumerate}

What I have not seen is any timing experiments. It might be worth a try before giving up on Steenrod squares. I have not done any myself, but I think you will be pleasantly surprised.

To conclude, here is the promised simplified formula for $Sq^1$ from Aubrey HB's thesis \cite{HB}. See there for the proof.

\begin{theorem}
Let $c\in C^p(X; Z_2)$. Then \begin{align*}
Sq^1(c)(\langle v_0, v_1, \cdots, v_{p+1}\rangle)=&\sum_{i=0}^{\lfloor\frac{p}{2}\rfloor}\sum_{j=0}^{i-1} c(\langle v_0, \cdots, \widehat{v_{2i}}, \cdots, v_{p+1}\rangle)\cdot c(\langle v_0, \cdots, \widehat{v_{2j}}, \cdots, v_{p+1}\rangle)\\
&+\sum_{i=0}^{\lfloor\frac{p}{2}\rfloor}\sum_{j=0}^{i-1} c(\langle v_0, \cdots, \widehat{v_{2i-1}}, \cdots, v_{p+1}\rangle)\cdot c(\langle v_0, \cdots, \widehat{v_{2j-1}}, \cdots, v_{p+1}\rangle),
\end{align*} where the notation $\widehat{v_k}$ means that vertex $v_k$ is left out, and all addition and multiplication is carried out in $Z_2.$
\end{theorem}

\chapter{Homotopy Groups of Spheres}

As homotopy groups of spheres is the last chapter in Hu's book, it seems a good place to conclude. A lot has happened since 1959, and this topic is now the subject of a number of entire books, so I am not going to even try to cover all of it. Instead, I will give you a small taste of it and provide you with a table which you can use for obstruction theory or just satisfy your curiosity.

I will start with a short description of {\it stable homotopy groups} of spheres that appears in Hatcher \cite{Hat}. These are sequences of homotopy groups you can get using the Freudenthal suspension theorem. Although these groups seem chaotic, there are some interesting patterns when you look at $p$-components of these groups. (I.e. quotients where elements of order prime to $p$ are eliminated.) Then I will talk about {\it classes} of abelian groups in which notions such as isomorphism are generalized. This is a useful tool for proving theorems about $p$-components of abelian groups. Then I will outline some of the techniques and examples from the last chapter of Hu \cite{Hu}. I would encourage you to read it if you ever need to deal with these groups as it is one of the only places I know of that includes explicit generators. Then I will briefly mention some of the more advanced techniques such as the Adams spectral sequence. I will also mention a very recent paper that claims to have found a huge number of new groups with machine computations. Finally, I will give you a table of $\pi_i(S^n)$ for $1\leq i\leq 15$ and $1\leq n\leq 10$ courtesy of  H. Toda \cite{Toda}. 

This material will stretch the limit of my own knowledge, so I can only give you an idea of a lot of it, but it is something I would like to look into further. Maybe I will have more to say in a sequel.

\section{Stable Homotopy Groups}

The material in this section is taken from Hatcher \cite{Hat}.

The Freudenthal suspension theorem (Theorem 9.7.8) states that the suspension map $\pi_i(X)\rightarrow \pi_{i+1}(SX)$ for an $(n-1)$-connected CW complex $X$ is an isomorphism for $i<2n-1$. We will raise the dimensions by 1 so for an $n$-connected complex, the map is an isomorphism for $i<2(n+1)-1=2n+1.$ This means that it is an isomorphism for $i\leq n$, so the suspension $SX$ is $(n+1)$-connected.

\begin{definition}
The sequence of iterated suspensions $$\pi_i(X)\rightarrow \pi_{i+1}(SX)\rightarrow\pi_{i+2}(S^2X)\rightarrow\cdots$$ eventually has all maps become isomorphisms. The group that results is the {\it stable homotopy group}\index{stable homotopy group} $\pi_i^S(X)$\index{$\pi_i^S(X)$}. The range of dimensions where the suspension is an isomorphism is called the {\it stable range}\index{stable range} and its study is called {\it stable homotopy theory.}\index{stable homotopy theory} The terms {\it unstable range} and {\it unstable homotopy theory} are often used for the opposites.
\end{definition}

\begin{definition}
The group $\pi_i^S(S^0)$ equals $\pi_{i+n}(S^n)$ for $n>i+1$. We write $\pi_i^S$ for $\pi_i^S(S^0)$ and call it the {\it stable i-stem}\index{stable i-stem}. 
\end{definition}

It turns out that $\pi_i^S$ is finite for $i>0$. At the time that Hatcher was writing his book, these were known for $i\leq 61$. The recent paper by Isaksen, Wang, and Xu \cite{IWX} has increased this to 90 with only a couple of exceptions. I will talk more about that later. Table 12.1.1 was borrowed by Hatcher from \cite{Toda}. It shows the groups $\pi_i^S$ for $i\leq 19$. Note that I will use the notation $(Z_p)^n$ for $Z_p\oplus \cdots \oplus Z_p$ ($n$ times).

\begin{table} 
\begin{center}
\begin{tabular}{|c|c|c|c|c|c|c|c|c|c|c|}
\hline
$i$ & 0 & 1 & 2 & 3 & 4 & 5 & 6 & 7 & 8 & 9\\
\hline
$\pi_i^s$ & $Z$ & $Z_2$ & $Z_2$ & $Z_{24}$ & 0 & 0 & $Z_2$ & $Z_{240}$ & $(Z_2)^2$ & $(Z_2)^3$\\ 
\hline
$i$ & 10 & 11 & 12 & 13 & 14 & 15 & 16 & 17 & 18 & 19\\
\hline
$\pi_i^s$ & $Z_6$ & $Z_{504}$ & $0$ & $Z_3$ & $(Z_2)^2$ & $Z_{480}\oplus Z_2$ & $(Z_2)^2$ & $(Z_2)^4$ & $Z_8\oplus Z_2$ & $Z_{264}\oplus Z_2$\\
\hline
\end{tabular}
\caption{ Stable Homotopy Groups $\pi_i^S$ for $i\leq 19$. \cite{Hat, Toda}.}
\end{center}
\end{table}

There are some interesting patterns when looking at {\it p-components} of these groups. 

\begin{definition}
Let $G$ be a group. Then the {\it p-component}\index{p-component} of $G$ is the quotient of $G$ obtained by factoring out all elements of $G$ whose orders are relatively prime to $p$.
\end{definition}

\begin{example}
Let $G=Z\oplus (Z_2)^3\oplus (Z_5)^2$. Then the 2-component of $G$ is $Z\oplus (Z_2)^3$ and the 5-component is  $Z\oplus (Z_5)^2$. The 3-component is $Z$. 
\end{example}

\begin{example}
For a finite abelian group, the p-component consists of all elements whose order is a power of $p$. Letting $H=(Z_2)^3\oplus (Z_5)^2$, the 2-component of $H$ is $(Z_2)^3$, the 5-component is  $(Z_5)^2$, and the 3-component is 0. 
\end{example}

\begin{figure}[ht]
\begin{center}
  \scalebox{0.6}{\includegraphics{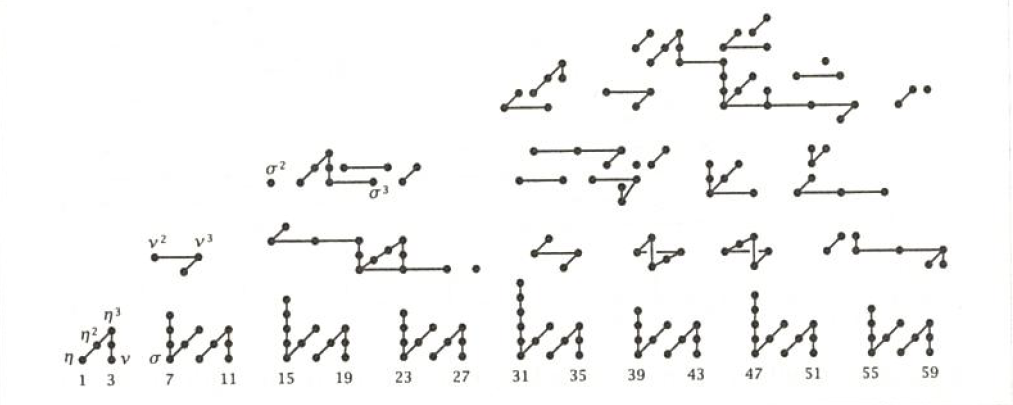}}
\caption{
\rm
2-components of $\pi_i^S$ for $i\leq 60$ \cite{Hat}. 
}
\end{center}
\end{figure}

Figure 12.1.1 from \cite{Hat} is a diagram of the 2-components of $\pi_i^S$ for $i\leq 60$. A vertical chain of $n$ dots in column $i$ represents a $Z_{2^n}$ summand of $\pi_i^S$. The bottom dot is a generator and moving up vertically represents multiplication by 2. The three generators $\eta$, $\nu$, and $\sigma$ for $i=1, 3, 7$ are the Hopf maps $S^3\rightarrow S^2$, $S^7\rightarrow S^4$, and $S^{15}\rightarrow S^8$ respectively. The horizontal and diagonal lines provide information about compositions of maps between spheres. We have products $\pi_i^S\times \pi_j^S\rightarrow \pi_{i+j}^S$ defined by compositions $S^{i+j+k}\rightarrow S^{j+k}\rightarrow S^k$. Hatcher proves that this product produces a graded ring structure on stable homotopy groups. 

\begin{theorem}
The composition products $\pi_i^S\times \pi_j^S\rightarrow \pi_{i+j}^S$  induce a graded ring structure on $\pi_\ast^S=\oplus_i\pi_i^S$ such that for $\alpha\in \pi_i^S$ and $\beta\in\pi_j^S$, $\alpha\beta=(-1)^{ij}\beta\alpha$.
\end{theorem}

It then follows that the $p$-components $_p\pi_\ast^S=\oplus{_i} (_p\pi_i^S)$ also form a graded ring with the same property. In $_2\pi_i^S$ many of the compositions with suspensions of the Hopf maps $\eta$ and $\nu$ are nontrivial. These are indicated in the diagram by segments extending one or two units to the right, diagonally for $\eta$ and horizontally for $\nu$. In the far left corner, you can see the relation $\eta^3=4\nu$ (each vertical step is multiplication by 2) in $_2\pi_3^S$. Now Table 12.1.1 shows $\pi_3^S\cong Z_{24}$ so the 2-component is $Z_8$. So in  $\pi_3^S$ the relation is $\eta^3=12\nu$ since $2\eta=0$ (there is no dot directly above $\eta$), so $2\eta^3=0$ and $\eta^3$ must be the unique element of order 2 in $Z_{24}.$ 

At the bottom we see a repeating pattern with spikes in dimensions $8k-1$. The spike in dimension $2^m(2n+1)-1$ has height $m+1$. (Try some examples.)

\begin{figure}[ht]
\begin{center}
  \scalebox{0.6}{\includegraphics{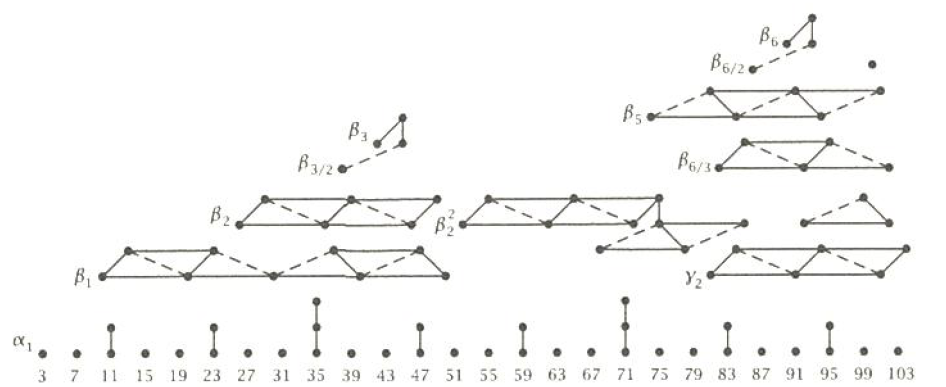}}
\caption{
\rm
3-components of $\pi_i^S$ for $i\leq 100$ \cite{Hat}. 
}
\end{center}
\end{figure}

Figure 12.1.2 shows the 3-components of $\pi_i^S$ for $i\leq 100$. Now vertical segments are multiplication by 3, the solid diagonal and horizontal segments denote compostion with $\alpha_1\in$ $_3\pi_3^S$ and $\beta_1\in$ $_3\pi_{10}^S$ respectively. The dashed lines involve a more complicated composition known as a {\it Toda bracket}. They are described in Hatcher briefly and Mosher and Tangora \cite{MT} in more detail. For now, I would call your attention to the regularity of the bottom of the diagram where there is a spike of height $m+1$ in dimension $4k-1$ where $m$ has the property that $3^m$ is the highest power of 3 dividing $4k$. 

\begin{figure}[ht]
\begin{center}
  \scalebox{0.6}{\includegraphics{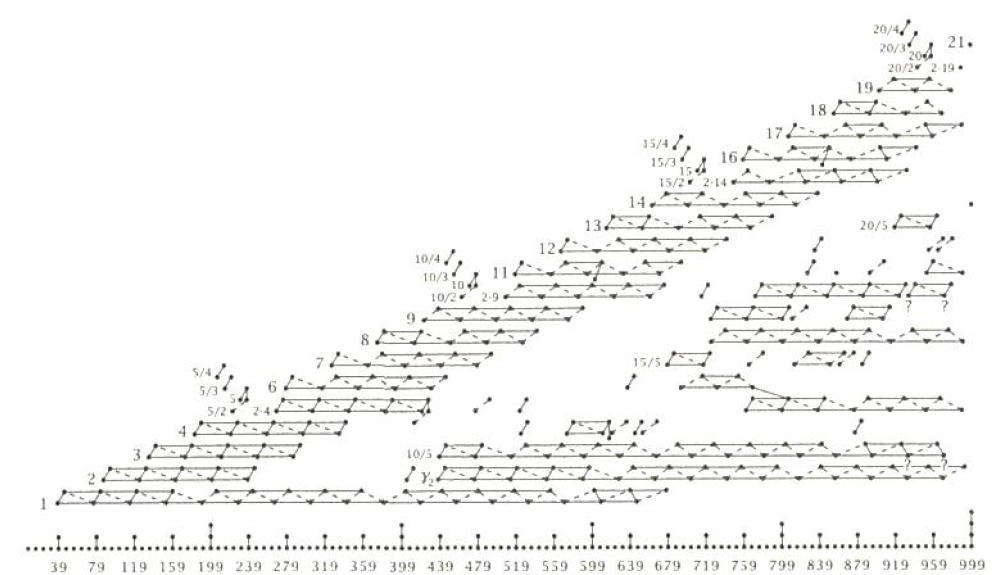}}
\caption{
\rm
5-components of $\pi_i^S$ for $i\leq 1000$ \cite{Hat}. 
}
\end{center}
\end{figure}

Even more regularity appears with higher primes. Figure 12.1.3 shows the 5-components for $i\leq 1000$. 

Hatcher produced these diagrams from tables published in Kochman \cite{Koc} and Kochman and Mahowald \cite{KM} for $p=2$. The computations for $p=3, 5$ are from Ravenel \cite{Rav}.

I haven't told you how any of those groups were computed. I will give some small examples. But first we will need to talk about {\it classes} of abelian groups. 

\section{Classes of Abelian Groups}

Although this subject is also covered in Hu \cite{Hu}, I will take the material in this section from Mosher and Tangora \cite{MT}.

Our goal will be to show that if $X$ and $Y$ are simply connected and all of their homology groups are finitely generated (always true for a finite complex), then if $f: X\rightarrow Y$ induces an isomorphism on homology with $Z_p$ coefficients, the $p$-components of the homotopy groups are also isomorphic. This is a valuable tool in computing homotopy groups.

\begin{definition}
A {\it class of abelian groups}\index{class of abelian groups} is a collection $\mathcal{C}$ of abelian groups satisfying the following axiom: 

{\bf Axiom 1:} If $0\rightarrow A'\rightarrow A \rightarrow A''\rightarrow 0$ is a short exact sequence, then $A$ is in $\mathcal{C}$ if and only if $A'$ and $A''$ are in $\mathcal{C}$.
\end{definition}

This means that classes of abelian groups are closed under formation of subgroups, quotient groups, and group extensions. (By definition $A$ is an {\it extension} of $A''$ by $A'$ as it is the middle term in the short exact sequence. This is a direct product if the sequence splits.)

We will always assume that $\mathcal{C}$ is nonempty.

\begin{definition}
A homomorphism $f: A\rightarrow B$ is a $\mathcal{C}${\it -monomorphism} if $\ker f\in \mathcal{C}$. It is a $\mathcal{C}${\it -epimorphism} if {\it coker }$f\in \mathcal{C}$. It is a $\mathcal{C}${\it -isomorphism} if it is both a $\mathcal{C}$-monomorphism and a $\mathcal{C}$-epimorphism. \end{definition}

Note that a $\mathcal{C}$-isomorphism does not necessarily have an inverse. For this reason, the relation $A\sim B$ if there is a $\mathcal{C}$-isomorphism from $A$ to $B$ is reflexive but not symmetric. We say $A$ and $B$ are $\mathcal{C}$-isomorphic if they are equivalent by the smallest reflexive, symmetric, and transitive relation containing $\sim$. Equivalently, there is a finite sequence $\{A=A_0, A_1, A_2, \cdots, A_n=B\}$ such that for each $i$ with $0\leq i<n$, there is a $\mathcal{C}$-isomorphism between $A_i$ and $A_{i+1}$ in one direction or the other.

We will also make use of the following axioms:

{\bf Axiom 2A:} If $A$ and $B$ are in $\mathcal{C}$, then so are $A\otimes B$ and $Tor(A, B)$.

{\bf Axiom 2B:} If $A$ is in $\mathcal{C}$, then so is $A\otimes B$ for every abelian group $B$. 

{\bf Axiom 3:}  If $A$ is in $\mathcal{C}$, then so is $H_n(A, 1; Z)$ for every $n>0$.

Axiom 2B implies 2A since $Tor(A, B)$ is a subgroup of $A\otimes R$ for a certain group $R$. 

\begin{example}
$\mathcal{C}_0$ contains only the group $G=\{0\}$ having one element. Then $\mathcal{C}_0$-monomorphism, $\mathcal{C}_0$-epimorphism, and $\mathcal{C}_0$-isomorphism reduce to the usual definitions. This satisfies axioms 1, 2B, and 3. Axiom 2A doesn't apply since $\mathcal{C}_0$ does not contain two distinct groups. 
\end{example}

\begin{example}
$\mathcal{C}_{FG}$ is the class of finitely generated abelian groups. It satisfies axioms 1, 2A, and 3. It can fail for 2B if $B$ is not finitely generated. 

For axiom 3, we know that $K(Z, 1)=S^1$ and $K(Z_2, 1)=P^\infty$. They have finitely many cells in each dimension so their homology is finitely generated. A similar structure can be built for $K(Z_p, 1),$ where $p$ is an odd prime. So axiom 3 holds.
\end{example}

\begin{example}
$\mathcal{C}_p$ where $p$ is a prime is the class of abelian torsion groups of finite exponent (i.e the least common multiple of the orders of the group elements is finite), and the order of every element is relatively prime to $p$. This class satisfies axioms 1, 2B, and 3, but the latter is not obvious. 
\end{example}

We will generalize three theorems we saw earlier to make use of classes of abelian groups. See \cite{MT} for the proofs.

\begin{theorem}
{\bf Hurewicz Theorem mod $\mathcal{C}$:} If $X$ is simply connected and $\mathcal{C}$ is a class satisfying axioms 1, 2A, and 3, then if $\pi_i(X)$ is in $\mathcal{C}$ for all $i<n$, we have that $H_i(X)$ is in $\mathcal{C}$ for all $i<n$, and the Hurewicz homomorphism $h: \pi_n(X)\rightarrow H_n(X)$ is a $\mathcal{C}$-isomorphism.
\end{theorem}

\begin{theorem}
{\bf Relative Hurewicz Theorem mod $\mathcal{C}$:} Let $A$ be a subspace of $X$ and let $X$ and $A$ be simply connected. Also suppose that $i_\#: \pi_2(A)\rightarrow \pi_2(X)$ is an epimorphism where $i: A\rightarrow X$ is inclusion. Suppose $\mathcal{C}$ is a class satisfying axioms 1, 2B, and 3. Then if $\pi_i(X, A)$ is in $\mathcal{C}$ for all $i<n$, we have that $H_i(X, A)$ is in $\mathcal{C}$ for all $i<n$, and $h: \pi_n(X, A)\rightarrow H_n(X, A)$ is a $\mathcal{C}$-isomorphism.
\end{theorem}

\begin{theorem}
{\bf Whitehead's Theorem mod $\mathcal{C}$:} Let $f: A\rightarrow X$ where $X$ and $A$ are simply connected and suppose $f_\#: \pi_2(A)\rightarrow \pi_2(X)$ is an isomorphism. Suppose $\mathcal{C}$ is a class satisfying axioms 1, 2B, and 3. Then $f_\#: \pi_i(A)\rightarrow \pi_i(X)$ is a $\mathcal{C}$-isomorphism for all $i<n$ and $f_\#: \pi_n(A)\rightarrow \pi_n(X)$ is a $\mathcal{C}$-epimorphism if and only if the same statements hold for $f_\ast$ on $H_\ast.$
\end{theorem}

It is an easy fact that for $\mathcal {C}$ satisfying Axiom 1, if $A\in\mathcal {C}$ and $A$ and $B$ are $\mathcal {C}$-isomorphic, then $B\in\mathcal {C}$. Then Theorem 12.2.1 implies the following.

\begin{theorem}
All homotopy groups of a finite simply connected complex are finitely generated.
\end{theorem}

We now state the main theorem for this statement. The proof is a pretty easy consequence of our three generalized theorems. (Their proofs are more complicated. See \cite{MT}.) For the final part we need one more Theorem proved in \cite{MT}.

\begin{theorem}
Let $f: A_1\rightarrow A_2$ be a homomorphism of finitely generated abelian groups. Suppose $f$ is a $\mathcal{C}_p$-isomorphism. Then $A_1$ and $A_2$ have isomorphic $p$-components. 
\end{theorem}

\begin{theorem}
{\bf $\mathcal{C}_p$ Approximation Theorem:} Let $X$ and $A$ be simply connected and {\it nice}. (CW-complexes are fine.) Suppose $H_i(A)$ and $H_i(X)$ are finitely generated for every $i$. Let $f: A\rightarrow X$ be such that $f_\#; \pi_2(A)\rightarrow \pi_2(X)$ is an epimorphism. (The map $f$ is then homotopic to an inclusion map, so we can let it be inclusion without loss of generality.) Then conditions 1-6 below are equivalent and imply condition 7.\begin{enumerate}
\item $f^\ast: H^i(X; Z_p)\rightarrow H^i(A; Z_p)$ is isomorphic for $i<n$ and monomorphic for $i=n$.
\item $f_\ast: H_i(A; Z_p)\rightarrow H_i(X; Z_p)$ is isomorphic for $i<n$ and epimorphic for $i=n$.
\item $H_i(X, A; Z_p)=0$ for $i\leq n$.
\item $H_i(X, A; Z)\in \mathcal{C}_p$ for $i\leq n$.
\item $\pi_i(X, A)\in \mathcal{C}_p$ for $i\leq n$.
\item $f_\#: \pi_i(A)\rightarrow \pi_i(X)$ is a $\mathcal{C}_p$-isomorphism for $i<n$ and a $\mathcal{C}_p$-epimorphism for $i=n$.
\item $\pi_i(A)$ and $\pi_i(X)$ have isomorphic $p$-components for $i<n$.
\end{enumerate}
\end{theorem}

So the theorem implies that if we want to find the $p$-component of $\pi_i(X)$, find a space $A$ with the same cohomology in dimension $i$ with $Z_p$ coefficients and a map of $A$ into $X$ inducing isomorphisms in $Z_p$ coefficients.

{\bf Proof:} Conditions 1 and 2 are equivalent by vector space duality, since $Z_p$ is a field and all $H_i$ are finitely generated.

Conditions 2 and 3 are equivalent by the exact homology sequence of the pair $(X, A)$.

Condition 3 implies condition 4, since the universal coefficient theorem gives an exact sequence $$0\rightarrow H_i(X, A)\otimes Z_p\rightarrow H_i(X, A; Z_p)\rightarrow Tor(H_{i-1}(X, A), Z_p)\rightarrow 0$$ and the middle group is zero by condition 3. But this implies that the left group is zero since it maps into zero by a monomorphism. Since $H_i(X, A)$ is finitely generated, $H_i(X, A)\otimes Z_p=0$ implies that $H_i(X, A)$ is the direct sum of finite groups of order prime to $p$, so $H_i(X, A)\in\mathcal{C}_p$, proving Condition 4. 

Condition 4 applied to $H_i$ and $H_{i-1}$ combined with the fact that $\mathcal{C}_p$ satisfies axiom 2B and thus axiom 2A, shows that the tensor product or Tor of these groups with $Z_p$ is in $\mathcal{C}_p$. Then the exact sequence shows that Condition 3 holds.

Conditions 4 and 5 are equivalent by the relative Hurewicz therorem. 

Conditions 5 and 6 are equivalent by the exact homotopy sequence mod $\mathcal{C}_p$ of the pair $(X, A)$.

Then Condition 6 implies condition 7 by the Hurewicz theorem mod $\mathcal{C}_{FG}$, which impliies that the groups $\pi_i(A)$ and $\pi_i(X)$ are finitely generated. Condition 7 now follows by Theorem 12.2.5. $\blacksquare$

\section{Some Techniques from the Early Years of the Subject}

In this section, I will describe some techniques and results from Hu's book \cite{Hu}. The book was written in 1959, so the Serre spectral sequence was already discovered, but not the Adams spectral sequence which I will describe in the next section. Steenrod squares had already been developed, but Hu doesn't cover them. Still, he has some interesting results derived using tools I have already covered, and it is instructive to see some examples. I will briefly outline what has happened in the last 60 years in the next section, and I will cite some references you can read if the subject interests you.

Hu describes the techniques involving the Freudenthal suspension theorem we have already covered in the first part of the chapter. I will now give some additional examples, all from his book.

\subsection{Finiteness of Homotopy Groups of Odd Dimensional Spheres}

We already know the homotopy groups of $S^1$ so assume that we are looking at $S^n$ for $n$ odd and $n\geq 3$.  The main theorem is as follows.

\begin{theorem}
If $S^n$ is an odd dimensional sphere and $m>n$ then $\pi_m(S^n)$ is finite.
\end{theorem}

We will need a definition.

\begin{definition}
A fiber space $p:E\rightarrow B$ where $B$ is pathwise connected is called {\it n-connective}\index{n-connective fiber space} if $E$ is $n$-connected and $$p_\ast: \pi_m(E, e_0)\rightarrow \pi_m(B, b_0)$$ is an isomorphism for $m>n$ where $e_0\in E$ and $p(e_0)=b_0$ and hence $\pi_{n+1}(B, b_0)$ is isomorphic to $H_{n+1}(E)$.
\end{definition}

Let $X$ be an $n$-connective fiber space over $S^n$. By the long exact sequence of a fiber space, $F$ is a $K(Z, n-1)$ space. Now by Theorem 12.2.4, all homotopy groups of $S^n$ are finitely generated. So by Theorem 12.2.1, all homology groups of $X$ are finitely generated. 

Using the Leray-Serre spectral sequence, it can be shown that $H^\ast(F; Z)$ is a polynomial algebra generated by an element of degree $n-1$. (Similar to the proof of Theorem 11.7.7 which is the special case $n=3$.) Hu then uses the Wang exact sequence for cohomology $$\cdots\rightarrow H^m(X; Z)\rightarrow H^m(F; Z)\xrightarrow{\rho^\ast}H^{m-n+1}(F; Z)\rightarrow H^{m+1}(X: Z)\rightarrow\cdots$$ to show that $$\rho^\ast: H^{p(n-1)}(F; Z)\rightarrow H^{(p-1)(n-1)}(F; Z)$$ is an isomorphism for every positive integer $p$. Then exactness shows that $H^m(X; Z)=0$ for every $m>0$. Since $H_m(X)$ is finitely generated, the universal coefficient theorem shows that $H_m(X)$ is finite for $m>0$. But then the Hurewicz Theorem mod $\mathcal{C}_{FG}$ shows that $\pi_m(X)$ is finite. So the theorem holds since $X$ is $n$-connective.

\subsection{Iterated Suspension}

For this chapter, I will break with Hu's notation and use the more common $\Omega X$ rather than $\Lambda X$ for the space of loops on $X$.

Pick a point $s_0\in S^n$ to be the base point. Let $S^n$ be the equator of the sphere $S^{n+1}$. Let $W=\Omega S^{n+1}$ which is the space of loops in $S^{n+1}$ starting and ending at $s_0$. Let $i: S^n\rightarrow W$ be defined as follows: Let $u$ and $v$ be the North and South poles of $S^{n+1}$ respectively. Define $i(x)$ for $x\neq s_0$ to be loop formed by joining $s_0$ to $u$ to $x$ to $v$ and then back to $s_0$ all by the shortest geodesic arcs. The loop $i(s_0)$ is the loop from $s_0$ to $u$ to $s_0$ to $v$ and back to $s_0$. The map $i$ is a homeomorphism of $S^n$ into $W$.

We have that $i(s_0)$ is homotopic to the degenerate loop $w_0\in W$, so $i(s_0)$ can be joined to $w_0$ by a path $\sigma$. For $m>0$ we have a homomorphism $$i_\ast: \pi_m(S^n, s_0)\rightarrow \pi_m(W, s_0)$$ where we identify $x$ with $i(x)$ so we write $s_0$ for $i(s_0)$. The path $\sigma$ induces an isomorphism $$\sigma_\ast: \pi_m(W, s_0)\approx \pi_m(W, w_0).$$ We also have that since $W$ is the space of loops on $S^{n+1}$ then there is an isomorphism $$h_\ast: \pi_m(W, w_0)\approx\pi_{m+1}(S^{n+1}, s_0).$$ Composing these three maps gives the map $$\Sigma=h_\ast\sigma_\ast i_\ast: \pi_m(S^n, s_0)\rightarrow \pi_{m+1}(S^{n+1}, s_0).$$ This map is called the {\it suspension map}\index{suspension map} and is equivalent to the suspension map defined in Definition 9.7.2. The Freudenthal suspension theorem (Theorem 9.7.8) also applies with this definition. 

If we repeat the process with $$j: \Omega(S^{n+1})\rightarrow \Omega^2(S^{n+2}),$$ we get an embedding $$k=ji: S^n\rightarrow\Omega^2(S^{n+2}).$$ For each $m$, we have that $k$ induces an isomorphism $$k_\ast: \pi_m(S^n, s_0)\rightarrow \pi_m(\Omega^2(S^{n+2}), s_0).$$ Similar to before, we also have an isomorphism $$\mathfrak{l}_\ast: \pi_m(\Omega^2(S^{n+2}), s_0)\rightarrow \pi_{m+2}(S^{n+2}, s_0).$$

\begin{theorem}
The map $\mathfrak{l}_\ast k_\ast$ is equal to the iterated suspension $\Sigma^2$. In addition $k_\ast$ is an isomorphism for $m<2n-1$ and an epimorphism for $m=2n-1$.
\end{theorem}

The next result uses $\Sigma^2$ to obtain information on $p$-primary components of the homotopy groups.

\begin{theorem}
Let $n\geq 3$ be an odd integer and $p$ a prime. Then the iterated suspension $$\Sigma^2: \pi_m(S^n)\rightarrow\pi_{m+2}(S^{n+2})$$ is a $\mathcal{C}_p$ isomorphism if $m<p(n+1)-3$ and is a $\mathcal{C}_p$ epimorphism if $m=p(n+1)-3.$
\end{theorem}

\begin{theorem}
If $n\geq 3$ is an odd integer, $p$ is prime, and $m<n+4p-6$,  then the $p$-primary components of $\pi_m(S^n)$ and $\pi_{m-n+3}(S^3)$ are isomorphic.
\end{theorem}

{\bf Proof:} We will proceed by induction on $n$. When $n=3$, this is obvious. Assume $q\geq 5$ is an odd integer and the theorem is true for every odd integer $n$ with $3\leq n<q$.

By the previous theorem, the $p$-primary components of $\pi_m(S^q)$ and $\pi_{m-2}(S^{q-2})$ are isomorphic if $m-2<p(q-1)-3.$ Since $q\geq 5$, we have $(p-1)(q-5)\geq 0$ so $$q+4p-8\leq p(q-1)-3.$$  So $m-2<p(q-1)-3$ whenever $m<q+4p-6.$

By our induction hypothesis, the $p$ primary components of $\pi_{m-2}(S^{q-2})$ and $\pi_{m-q+3}(S^3)$ are isomorphic if $m<q+4p-6$. So the result is true if $n=q$. $\blacksquare$

\begin{example}
Let $p=7$. Then the $7$-primary components of $\pi_m(S^n)$ and $\pi_{m-n+3}(S^3)$ are isomorphic if $m<n+28-6=n+22$ where $n$ is odd. So for example, the 7-components of $\pi_m(S^{11})$ and $\pi_{m-8}(S^3)$ are isomorphic if $8<m<33$. 
\end{example}

\subsection{The $p$-primary components of $\pi_m(S^3)$}
Let $X$ be a 3-connective fiber space over $S^3$. We know the fiber $F$ is a $K(Z, 2)$ space. 

\begin{theorem}
$H_m(X)=0$ if $m$ is odd and $H_{2n}(X)\cong Z_n$ for every $n>0$.
\end{theorem}

This is also derived from the Wang exact sequence for cohomology and the fact that $H^\ast(F)$ is a polynomial algebra with generator of dimension 2. 

\begin{theorem}
If $p$ is prime then the $p$-primary component of $\pi_m(S^3)$ is 0 if $m<2p$ and $Z_p$ is $m=2p$.
\end{theorem}

{\bf Proof:} Letting $\mathcal{C}_p$ be that class of abelian groups of order prime to $p$, we have by Therorem 12.3.5 that $H_m(X)\in \mathcal{C}_p$ if $0<m<2p$.  By the Hurewicz theorem mod $\mathcal{C}_p$, we have that $\pi_m(X)\in \mathcal{C}_p$ for $0<m<2p$, and $\pi_{2p}(X)$ is $\mathcal{C}_p$-isomorphic to $Z_p$. Since $X$ is a 3-connective fiber space over $S^3$, we have that $\pi_m(S^n)\cong \pi_m(X)$ for each $m>3$. $\blacksquare$

\begin{theorem}
If $n\geq 3$ is an odd integer and $p$ is prime then the $p$-primary component of $\pi_m(S^n)$ is 0 if $m<n+2p-3$ and is $Z_p$ if $m=n+2p-3.$
\end{theorem}

{\bf Proof:} 
Since $p\geq 2$, we have $$n+2p-3<n+4p-6.$$ Then by Theorem 12.3.4, the $p$-primary components of  $\pi_m(S^n)$ and $\pi_{m-n+3}(S^3)$ are isomorphic.$\blacksquare$

\begin{example}
If $p=7$ and $n=11$ then the $7$-primary component of $\pi_m(S^{11})$ is 0 if $m<11+14-3=22$ and is $Z_7$ if $m=22$.
\end{example}

\subsection{Pseudo-projective Spaces}

\begin{definition}
A space $P=P_h^{n+1}$ is a {\it pseudo-projective space}\index{pseudo-projective space} if it is formed by attaching an $(n+1)$-cell $E^{n+1}$ to $S^n$ by a map $\phi: \partial E^{n+1}=S^n\rightarrow S^n$ of degree $h>0$. 
\end{definition}

The homology groups are $H_0(P)\cong Z$, $H_n(P)\cong Z_h$ and $H_i(P)=0$ for $i\neq 0, n$. 

\begin{theorem}
For every $m<2n-1$, there is an exact sequence $$0\rightarrow \pi_m(S^n)\otimes Z_h\rightarrow\pi_m(P)\rightarrow Tor(\pi_{m-1}(S^n), Z_h)\rightarrow 0.$$
\end{theorem}

Let $X$ be a 3-connective fiber space over $S^3$ with projection $\omega: X\rightarrow S^3.$ Since the $p$-primary component of $\pi_{2p}(X)$ is $Z_p$ , there is a map $f: S^{2p}\rightarrow X$ representing a generator $[f]$ of this component. 

Now let $P=P_p^{2p+1}.$ Then $S^{2p}\subset P$, and $p[f]=0$ implies that $f$ can be extended to a map $g: P\rightarrow X$. Compose with the projection $\omega: X\rightarrow S^3$ to get $\chi=\omega g: P\rightarrow S^3$. 

\begin{theorem}
$\chi_\ast: \pi_m(P)\rightarrow \pi_m(S^3)$ is a monomorphism for $m<4p-1$. It sends the $p$-primary component of $\pi_m(P)$ onto that of $\pi_m(S^3)$ for $m\leq 4p-1$.
\end{theorem}

\begin{theorem}
If $p$ is a prime number and $m\leq 4p-2$, then \begin{align*}
\pi_{2p}(P)&\cong Z_p,\\
\pi_{4p-3}(P)&\cong Z_p,\\
\pi_{4p-2}(P)&\cong Z_p\hspace{.4 in}\mbox{if }p>2,\\
\pi_m(P)&=0\hspace{.4 in}\mbox{otherwise.}
\end{align*}
\end{theorem}

{\bf Proof:} By Theorem 12.3.5 and Hurewicz's Theorem, $\pi_m(P)=0$ for $m<2p$ and $\pi_{2p}(P)=Z_p$. Applying Theorem 12.3.8 with $n=2p$, $h=p$, and $m\leq 4p-3$, we get an exact sequence $$0\rightarrow \pi_m(S^{2p})\otimes Z_p\rightarrow\pi_m(P)\rightarrow Tor(\pi_{m-1}(S^{2p}), Z_p)\rightarrow 0.$$  By the Freudenthal suspension theorem, $\pi_m(S^{2p})\cong\pi_{m-1}(S^{2p-1})$ for every $m\leq 4p-3$. Since $2p-1\geq 3$, Theorem 12.3.7 implies that the $p$-primary component of $\pi_{m-1}(S^{2p-1})$ is 0 if $m<4p-3$ and $Z_p$ if $m=4p-3$. This gives $\pi_m(P)=0$ for $2p<m<4p-3$ and $\pi_{4p-3}(P)\cong Z_p$.

By Theorem 12.3.9, the $p$-primary component of $\pi_m(S^3)$ is 0 whenever $2p<m<4p-3,$ and is $Z_p$ if $m=4p-3$. Then by Theorem 12.3.4, if $n\geq 3$ is odd, the $p$-primary component of $\pi_m(S^n)$ is 0 for $n+2p-3<m<n+4p-6$. Now $\pi_m(S^{2p})\cong \pi_{m+1}(S^{2p+1})$ for every $m\leq 4p-2$. So the $p$-primary component of $\pi_{m+1}(S^{2p+1})$ is 0 if $4p-3<m<6p-6.$

If $p>2$ then $4p-2<6p-6$ so the $p$-primary component of $\pi_{4p-2}(S^{2p})$ is 0. By Theorem 12.3.8 with $n=2p$, $h=p$, and $m= 4p-2$, we get $\pi_{4p-2}(P)\cong Z_p$. $\blacksquare$

Now we have two facts about the homotopy groups of $S^3$.

\begin{theorem}
If $p$ is a prime number, then the $p$-primary component of $\pi_m(S^3)$ is 0 if $2p<m<4p-3$ and is $Z_p$ if $m=4p-3$. If $p>2$, then the $p$-primary component of $\pi_{4p-2}(S^3)$ is $Z_p$.
\end{theorem}

\begin{theorem}
If $n\geq 3$, is an odd integer and $p$ prime then the $p$-primary component of $\pi_m(S^n)$ is 0 if $n+2p-3<m<n+4p-6$ and that of $\pi_{n+4p-6}(S^n)$ is 0 or $Z_p$.
\end{theorem}

\subsection{The Hopf Invariant Revisited}

Recall the definition of the Hopf invariant from Section 11.5. We have a map $f: S^{2n-1}\rightarrow S^n$. Letting $K=S^n\cup_f e^{2n}$, we have that if $\sigma$ and $\tau$ are generators of $H^n(K; Z)$ and $H^{2n}(K; Z)$ respectively, then the Hopf invariant $H(f)$ of $f$ is defined by $\sigma\cup\sigma=H(f)\tau$. It turns out that $H(f)=0$ if $n$ odd. Hu proves the following: 

\begin{theorem}
Let $n$ be even and $f: S^{2n-1}\rightarrow S^n$ be a map with Hopf invariant $H(f)=k\neq 0$. Let $\mathcal{C}$ denote the class of all finite abelian groups of order dividing some power of $k$. If $$\chi_f: \pi_{m-1}(S^{n-1})\oplus\pi_m(S^{2n-1})\rightarrow \pi_m(S^n)$$ is the homomorphism defined by $$\chi_f(\alpha, \beta)=\Sigma(\alpha)\oplus f_\ast(\beta)$$ where $\alpha\in\pi_{m-1}(S^{n-1})$, $\beta\in \pi_m(S^{2n-1})$, and $\Sigma$ is the suspension map, then $\chi_f$ is a $\mathcal{C}$-isomorphism for every $m>1$.
\end{theorem}

\begin{theorem}
If $H(f)=\pm 1$, then $\chi_f$ is an isomorphism and the suspension $\Sigma: \pi_{m-1}(S^{n-1})\rightarrow \pi_m(S^n)$ is a monomorphism for every $m>1.$
\end{theorem}

Theorem 11.5.1 states that there exists a map of Hopf invariant 2, so we have the following.

\begin{theorem}
Let $n$ be even and let $\mathcal{C}$ denote the class of all finite abelian groups of order dividing some power of $2$. Then the homotopy group $\pi_m(S^n)$ is $\mathcal{C}$-isomorphic to $\pi_{m-1}(S^{n-1})\oplus\pi_m(S^{2n-1})$.
\end{theorem}

\subsection{$\pi_{n+1}(S^n)$ and $\pi_{n+2}(S^n)$}

Since $\pi_m(S^1)=0$ for $m>1$ and $\pi_m(S^2)\cong\pi_m(S^3)$ for $m>2$, assume $n\geq 3.$

\begin{theorem}
$\pi_{n+1}(S^n)\cong Z_2$ for $n\geq 3$. 
\end{theorem}

{\bf Proof:} Let $X$ be a 3-connective fiber space over $S^3$. Then $$\pi_4(S^3)\cong\pi_4(X)\cong H_4(X)\cong Z_2,$$ where the last equality follows by Theorem 12.3.5. By the Freudenthal suspension theorem we have $$\pi_4(S^3)\cong\pi_5(S^4)\cong\cdots\pi_{n+1}(S^n)\cong\cdots,$$ so $\pi_{n+1}(S^n)\cong Z_2$ for $n\geq 3$. $\blacksquare$

Now recall that the generator of $\pi_3(S^2)$ is the Hopf map $p$. The suspension $\Sigma: \pi_3(S^2)\rightarrow \pi_4(S^3)$ is an epimorphism, so $\Sigma p: S^4\rightarrow S^3$ is the generator of $\pi_4(S^3)$. In general, the generator of $\pi_{n+1}(S^n)$ is the $(n-2)$-fold suspension $\Sigma^{n-2}p$ of the Hopf map $p$.

\begin{theorem}
$\pi_{n+2}(S^n)\cong Z_2$ for $n\geq 3$. 
\end{theorem}

{\bf Proof:} Apply Theorem 12.3.8 with $n=4$, $h=2$, and $m=5$. Then we get an exact sequence $$0\rightarrow\pi_5(S^4)\otimes Z_2 \rightarrow \pi_5(P_2^5)\rightarrow Tor(Z, Z_2)\rightarrow 0.$$ Now $\pi_5(S^4)\cong Z_2$ by the previous result and $Z_2\otimes Z_2=Z_2$. Also $Tor(Z, Z_2)=0$, so $\pi_5(P_2^5)\cong Z_2$. Then Theorem 12.3.9 implies that the 2-primary component of $\pi_5(S^3)$ is isomorphic to $Z_2$. By Theorem 12.3.3, The $p$-primary component of $\pi_5(S^3)$ is 0 for $p>2$, so $\pi_5(S^3)\cong Z_2$.

Since the Hopf map $S^7\rightarrow S^4$ has Hopf invariant 1, Apply Theorem 12.3.14 with $n=4$ and $m=6$ so that $$\pi_6(S^4)\cong \pi_5(S^3)\oplus \pi_6(S^7)\cong Z_2\oplus 0=Z_2.$$ By the Freudenthal suspension theorem we have $$\pi_6(S^4)\cong\pi_7(S^5)\cong\cdots\pi_{n+2}(S^n)\cong\cdots,$$ so $\pi_{n+2}(S^n)\cong Z_2$ for $n\geq 3$. $\blacksquare$

In general, let $p: S^3\rightarrow S^2$ be the Hopf map. Then the suspension is $\Sigma p: S^4\rightarrow S^3$ and the 2-fold suspension is $\Sigma^2 p: S^5\rightarrow S^4$. The generator of $\pi_5(S^3)$ is represented by $q=\Sigma p\Sigma^2 p: S^5\rightarrow S^3$. In general, $\pi_{n+2}(S^n)$ is generated by the $(n-3)$-fold suspension $\Sigma^{n-3} q$.

\subsection{$\pi_{n+r}(S^n)$ for $3\leq r\leq 8$}

Hu lists some results for $\pi_{n+r}(S^n)$ for $3\leq r\leq 8$. For $r=3, 4$, he computes them with methods similar to those we have just used, and he refers readers to \cite{Toda} for the rest. Here is the rest of his list:

$r=3$
\begin{align*}
\pi_6(S^3)&\cong Z_{12}\\
\pi_7(S^4)&\cong Z\oplus Z_{12}\\
\pi_{n+3}(S^n)&\cong Z_{24}\hspace{.4 in}\mbox{if }n\geq 5.\\
\end{align*}

$r=4$
\begin{align*}
\pi_7(S^3)&\cong Z_2\\
\pi_8(S^4)&\cong Z_2\oplus Z_2\\
\pi_9(S^5)&\cong Z_2\\
\pi_{n+4}(S^n)&=0\hspace{.4 in}\mbox{if }n\geq 6.\\
\end{align*}

$r=5$
\begin{align*}
\pi_8(S^3)&\cong Z_2\\
\pi_9(S^4)&\cong Z_2\oplus Z_2\\
\pi_{10}(S^5)&\cong Z_2\\
\pi_{11}(S^6)&\cong Z\\
\pi_{n+5}(S^n)&=0\hspace{.4 in}\mbox{if }n\geq 7.\\
\end{align*}

$r=6$
\begin{align*}
\pi_9(S^3)&\cong Z_3\\
\pi_{10}(S^4)&\cong Z_{24}\oplus Z_3\\
\pi_{n+6}(S^n)&\cong Z_2\hspace{.4 in}\mbox{if }n\geq 5.\\
\end{align*}

$r=7$
\begin{align*}
\pi_{10}(S^3)&\cong Z_{15}\\
\pi_{11}(S^4)&\cong Z_{15}\\
\pi_{12}(S^5)&\cong Z_{30}\\
\pi_{13}(S^6)&\cong Z_{60}\\
\pi_{14}(S^7)&\cong Z_{120}\\
\pi_{15}(S^8)&\cong Z\oplus Z_{120}\\
\pi_{n+7}(S^n)&\cong Z_{240}\hspace{.4 in}\mbox{if }n\geq 9.\\
\end{align*}

$r=8$
\begin{align*}
\pi_{11}(S^3)&\cong Z_2\\
\pi_{12}(S^4)&\cong Z_2\\
\pi_{13}(S^5)&\cong Z_2\\
\pi_{14}(S^6)&\cong Z_{24}\oplus Z_2\\
\pi_{15}(S^7)&\cong Z_2\oplus Z_2\oplus Z_2\\
\pi_{16}(S^8)&\cong Z_2\oplus Z_2\oplus Z_2\oplus Z_2\\
\pi_{17}(S^9)&\cong Z_2\oplus Z_2\oplus Z_2\\
\pi_{n+8}(S^n)&\cong Z_2\oplus Z_2\hspace{.4 in}\mbox{if }n\geq 10.\\
\end{align*}

One final note: Hu gives $\pi_{10}(S^4)\cong Z_{24}\oplus Z_2$, but after checking \cite{Toda}, I believe this to be a typo and that $\pi_{10}(S^4)\cong Z_{24}\oplus Z_3$ is actually correct.

\section{More Modern Techniques and Further Reading}

The entire second half of Mosher and Tangora \cite{MT} is devoted to computing homotopy groups of spheres. 

Their first approach to build an approximation to $S^n$ using a Postnikov system. As a first approximation, $K(Z, n)$ has the same cohomology and homology groups as $S^n$ up to dimension $n$. Unfortuantely, $K(Z, n)$ has nonzero cohomology in higher dimensions. We would like to have a better approximation. Here we look at 2-components and use mod 2 coefficients. Now $H^{n+1}(Z, n; Z_2)=0$ and $H^{n+2}(Z, n; Z_2)\cong Z_2$ generated by $Sq^2(\imath_n)$ where $\imath_n$ is the fundamental class. (See Theorem 11.7.11) This class defines a map $K(Z, n)\rightarrow K(Z_2, n+2)$ by Theorem 11.2.1. From the standard contractible fibre space (i.e the path-space fibration) over $K(Z_2, n+2)$, the map induces a fiber space $X_1$ over $K(Z, n)$ with fiber $\Omega K(Z_2, n+2)=K(Z_2, n+1)$. It turns out that $X_1$ is a better approximation to $S^n$ than $K(Z, n)$ . This is because the class $Sq^2(\imath_n)\in H^{n+2}(Z, n; Z_2)$ has been killed so that $H^{n+2}(X_1; Z_2)=H^{n+2}(S^n; Z_2)=0$. Using this construction along with the $\mathcal{C}_p$ Approximation Theorem (Theorem 12.2.6), Mosher and Tangora are able to compute the 2-components of stable homotopy groups $\pi_{n+k}(S^n)$ for $k\leq 7$. Note that I have left out a lot of details, so see the description in \cite{MT}.

Mosher and Tangora use the Postnikov decomposition of $X$ to construct a spectral sequence for $[Y, X]$, the homotopy classes of maps from $Y$ to $X$. If $Y=S^m$ and $X=S^n$ where $m\geq n$, we get $\pi_m(S^n)$ as a special case. The problem is that this spectral sequence starts at $H^\ast(Y; \pi_\ast X)$ (See \cite{MT}, Chapter 14.) so it requires us to know the homotopy groups of $X$ already. This lead to the big breakthrough that happened just as Hu's book was being wirtten. This was the discovery of the Adams spectral sequence \cite{JFA2}.

The Adams spectral sequence gives information about the 2-component of $[Y, X]$. It is valid in the stable range, so if $X$ is $(n-1)$-connected, we need $dim(Y)\leq 2n-2.$ Finally, the $E^2$ term is the module $Ext_\mathcal{A}(H^\ast(X),H^\ast(Y))$ where we are talking about the cohomology of $X$ and $Y$ as modules over the Steenrod algebra $\mathcal{A}$. It turns out that $Ext_\mathcal{A}(H^\ast(X),H^\ast(Y))$ also has a multiplicative structure. Although \cite{MT} doesn't get into that, there is more in McLeary's book on spectral sequences \cite{McCl} as well as Adams' original paper \cite{JFA2}. 

To get the full homotopy groups, we need the $p$-primary components where $p>2$ is an odd prime. For this, we need Steenrod reduced powers and the corresponding mod $p$ Steenrod algebras $\mathcal{A}_p$. For more on $p$-components of homotopy groups of spheres, see Steenrod and Epstein \cite{SE} and McLeary \cite{McCl}.

Meanwhile, Toda had computed several homotopy groups of spheres using the {\it EHP spectral sequence} related to the fibration $$S^n\rightarrow\Omega S^{n+1}\rightarrow \Omega S^{2n+1}$$ over $Z_2$. The first term is $E_1^{k, n}=\pi_{k+n}(S^{2n-1})$ where the 2-primary components are implied. He used slightly different fibrations for the $p$-primary components. His methods were used to make the table in Hu's book. See \cite{Toda} for details.

It turns out that calculating the Ext term for the Adams Spectral sequence is very hard. One approach is to use the May spectral sequence \cite{May3, May4}. The Adams spectral sequence can be strengthened by replacing cohomology mod $p$ by a generalized cohomology theory such as {\it complex cobordism}. This was the innovation of Novikov \cite{Nov} in 1967. For a complete description of cobordism theory including complex cobordism, see \cite{Stong}. A good description of how all of this fits together to compute stable homotopy groups of spheres, see Ravenel \cite{Rav}. The book assumes a strong background, so I would recommend reading Mosher and Tangora and the relevant parts of McLeary first.

Finally, Isaksen, Wang, and Xu \cite{IWX} made the most recent advance in the subject in 2020. The {\it motivic homotopy theory} of Morel and Voevodsky \cite{MV} was originally used to apply homotopy methods to algebraic geometry. Using these methods, Isaksen, et. al. have computed the stable homotopy groups $\pi_i^S$ up to $i=90$ beating out the old result of $i\leq 61$. This has been made possible by computer calculations of Ext groups involved in the Adams spectral sequence. This has produced data up to dimension 200 which has not yet been fully interpreted. The relevant algorithms were developed by Bruner \cite{Bru1, Bru2, Bru3}, Nassau \cite{Nas}, and Wang \cite{Wang}. See \cite{IWX}  and its references for the details.

\section{Table of Some Homotopy Groups of Spheres}

I will conclude this section with the promised table. Recall that we want $\pi_i(S^n)$ for $1\leq n\leq 10$ and $1\leq i\leq 15.$

First of all, we know that $\pi_1(S^1)\cong Z$ and $\pi_i(S^1)=0$ for $i>1$. 

We also know that $\pi_n(S^n)\cong Z$ and $\pi_i(S^n)=0$ for $i<n$.

Theorem 9.6.8 says that $\pi_3(S^2)\cong Z$ and Example 9.6.1 says that $\pi_n(S^2)\cong\pi_n(S^3)$ for $n\geq 3$.

As an exercise, step through the results in Section 12.3 and see what you can fill in. Then use the stable stems that were listed but not proved in Section 12.1. There will still be a few groups that you will not find, but you will get a lot of them.

Table 12.5.1 is the result. I will write $Z_p^n$ for the direct sum of $Z_p$ with itself $n$-times.

One last comment. The one place where I can see homoopy groups of spheres showing up in data science would be in obstruction theory. Recall that obstructions to the extension problem where we want to extend $f: A\rightarrow Y$ to $f: X\rightarrow Y$ lie in cohomology groups whose coefficients are homotopy groups of $Y$. If $Y$ is a sphere, it would be useful to know some of its low dimensional homotopy groups. 

\renewcommand\arraystretch{1.4}
\begin{table} 
\begin{center}
\begin{tabular}{|c|c|c|c|c|c|c|c|c|c|c|}
\hline
$S^n$ & $S^1$ & $S^2$ & $S^3$ & $S^4$ & $S^5$ & $S^6$ & $S^7$ & $S^8$ & $S^9$ & $S^{10}$\\
\hline
$\pi_1(S^n)$ & $Z$ & $0$ & $0$ & $0$ & $0$ & $0$ & $0$ & $0$ & $0$ & $0$\\ 
\hline
$\pi_2(S^n)$ & $0$ & $Z$ & $0$ & $0$ & $0$ & $0$ & $0$ & $0$ & $0$ & $0$\\ 
\hline
$\pi_3(S^n)$ & $0$ & $Z$ & $Z$ & $0$ & $0$ & $0$ & $0$ & $0$ & $0$ & $0$\\ 
\hline
$\pi_4(S^n)$ & $0$ & $Z_2$ & $Z_2$ & $Z$ & $0$ & $0$ & $0$ & $0$ & $0$ & $0$\\ 
\hline
$\pi_5(S^n)$ & $0$ & $Z_2$ & $Z_2$ & $Z_2$ & $Z$ & $0$ & $0$ & $0$ & $0$ & $0$\\ 
\hline
$\pi_6(S^n)$ & $0$ & $Z_{12}$ & $Z_{12}$ & $Z_2$ & $Z_2$ & $Z$ & $0$ & $0$ & $0$ & $0$\\ 
\hline
$\pi_7(S^n)$ & $0$ & $Z_2$ & $Z_2$ & $Z\oplus Z_{12}$ & $Z_2$ & $Z_2$ & $Z$ & $0$ & $0$ & $0$\\ 
\hline
$\pi_8(S^n)$ & $0$ & $Z_2$ & $Z_2$ & $Z_2^2$ & $Z_{24}$ & $Z_2$ & $Z_2$ & $Z$ & $0$ & $0$\\ 
\hline
$\pi_9(S^n)$ & $0$ & $Z_3$ & $Z_3$ & $Z_2^2$ & $Z_2$ & $Z_{24}$ & $Z_2$ & $Z_2$ & $Z$ & $0$\\ 
\hline
$\pi_{10}(S^n)$ & $0$ & $Z_{15}$ & $Z_{15}$ & $Z_{24}\oplus Z_3$ & $Z_2$ & $0$ & $Z_{24}$ & $Z_2$ & $Z_2$ & $Z$\\ 
\hline
$\pi_{11}(S^n)$ & $0$ & $Z_2$ & $Z_2$ & $Z_{15}$ & $Z_2$ & $Z$ & $0$ & $Z_{24}$ & $Z_2$ & $Z_2$\\ 
\hline
$\pi_{12}(S^n)$ & $0$ & $Z_2^2$ & $Z_2^2$ & $Z_2$ & $Z_{30}$ & $Z_2$ & $0$ & $0$ & $Z_{24}$ & $Z_2$\\ 
\hline
$\pi_{13}(S^n)$ & $0$ & $Z_{12}\oplus Z_2$ & $Z_{12}\oplus Z_2$ & $Z_2^3$ & $Z_2$ & $Z_{60}$ & $Z_2$ & $0$ & $0$ & $Z_{24}$\\ 
\hline
$\pi_{14}(S^n)$ & $0$ & $Z_{84}\oplus Z_2^2$ & $Z_{84}\oplus Z_2^2$ & $Z_{120}\oplus Z_{12}\oplus Z_2$ & $Z_2^3$ & $Z_{24}\oplus Z_2$ & $Z_{120}$ & $Z_2$ & $0$ & $0$\\ 
\hline
$\pi_{15}(S^n)$ & $0$ & $Z_2^2$ & $Z_2^2$ & $Z_{84}\oplus Z_2^5$ & $Z_{72}\oplus Z_2$ & $Z_2^3$ & $Z_2^3$ & $Z\oplus Z_{120}$ & $Z_2$ & $0$\\ 
\hline

\end{tabular}
\caption{ Homotopy Groups $\pi_i(S^n)$ for  $1\leq n\leq 10$ and $1\leq i\leq 15$ \cite{Toda}.}
\end{center}
\end{table}

\chapter{Conclusion}

This book serves three functions. Chapters 2-4 give you the math background in point set topology, abstract algebra, and homology to enable you to understand most current papers on topological data analysis. Chapters 5-7 give you a taste of what is going on currently. It should give you an idea of what to do with data in the form of a point cloud, graph, image, or time series. It also includes a list of open source software that should enable you to easily perform experiments. Finally in Chapters 8-12, I give you a taste of some more advanced topics in algebraic topology such as cohomology, homotopy, obstruction theory, and Steenrod cohomology operations. Although there are questions of complexity, I have pointed you to several papers that have addressed algorithmic issues and a number of open source software projects. Although persistent homology works very well in a lot of practical cases, it would be interesting to experiment with some of the other techniques and see if they can help in solving hard classification problems. Their application in data science is an area that is currently wide open.

As I write this book, topological data analysis is evolving rapidly. The goal of this book is to enable you to understand the range of topics and enable you to experiment and contribute to this exciting field.

% For those of you who use bibliography files, the lines commented out
% below are for you to use (and comment out the lines that make use of
% the ``thebibliography'' environment below).  Be sure to choose the
% style you want (if other than `plain') and insert the filename of
% your ``.bib'' file as the argument to the `\bibliography' command.

%\bibliographystyle{plain}
%\bibliography{YourBibFile}
\pagebreak
% The standard ``thebibliography'' environment, adjusted slightly in
% the `mypaper.cls' file.

\printindex

\end{document}